\documentclass[11pt]{amsart}

\usepackage[english]{babel}
\usepackage[utf8]{inputenc}
\usepackage[OT1]{fontenc}
\usepackage{lmodern}
\usepackage{mathrsfs}
\usepackage{bbm}
\usepackage{amsmath}
\usepackage{amssymb}
\usepackage{mathtools}
\usepackage{color}
\usepackage{hyperref}
\usepackage{tikz-cd}
\usepackage{nicefrac}
\usepackage{enumitem}
\usepackage{xcolor}
\usepackage{graphicx}
\usepackage{adjustbox}
\usepackage{multicol}
\usepackage{csquotes}
\usepackage{caption}
\usepackage{subcaption}
\usepackage{titlesec}
\usepackage{xpatch}

\usepackage{geometry}
\setlist[itemize]{
    label=\textbullet,
    leftmargin=3em,
    itemsep=0em,
    topsep=0pt,
}

\setlist[enumerate]{
    leftmargin=3em,
    itemsep=0pt,
    topsep=0pt,
}

\usepackage[
	backend=biber,
	style=alphabetic,
	maxbibnames=99,
	ibidtracker=true,
	sorting=nyt,
	doi=false,
    	isbn=false,
    	url=false]{biblatex}
 	
\renewbibmacro{in:}{}

\newcommand{\C}{\mathbb{C}}
\newcommand{\Z}{\mathbb{Z}}
\newcommand{\N}{\mathbb{N}}
\newcommand{\PP}{\mathbb{P}}
\newcommand{\kk}{\mathbbm{k}}
\newcommand\bbS{\mathbb{S}}

\newcommand\bi{\mathbf i}
\newcommand\bj{\mathbf j}
\newcommand\bR{\mathbf R}
\newcommand\bS{\mathbf S}
\newcommand\bT{\mathbf T}
\newcommand{\bfU}{\mathbf{U}}

\newcommand\cB{\mathcal{B}}
\newcommand\cE{\mathcal{E}}
\newcommand\cZ{\mathcal{Z}}
\newcommand\cR{\mathcal{R}}
\newcommand\cP{\mathcal{P}}
\newcommand\cI{\mathcal{I}}
\newcommand\cJ{\mathcal{J}}
\newcommand\cF{\mathcal{F}}
\newcommand\cG{\mathcal{G}}
\renewcommand\O{\mathcal{O}}
\newcommand\cL{\mathcal{L}}

\newcommand\fg{\mathfrak g}

\newcommand\fb{\mathfrak b}
\newcommand\fh{\mathfrak h}
\newcommand\fn{\mathfrak n}
\newcommand\fr{\mathfrak r}
\newcommand\fsl{\mathfrak{sl}}
\newcommand{\fgl}{\mathfrak{gl}}

\newcommand{\sfPsi}{\mathsf{\Psi}}
\newcommand{\sfA}{\mathsf{A}}
\newcommand{\sfY}{\mathsf{Y}}

\newcommand\scrC{\mathscr{C}}

\newcommand{\la}{\lambda}

\newcommand{\op}[1]{\operatorname{#1}}
\DeclareMathOperator{\wt}{wt}
\DeclareMathOperator{\awt}{awt}
\newcommand{\aT}{A}
\newcommand{\ul}[1]{\underline{#1}}
\newcommand{\Gr}{\op{Gr}}
\newcommand{\Pic}{\op{Pic}}
\newcommand{\Fl}{\op{Fl}}
\newcommand{\Pro}{\op{Pro}}
\newcommand{\Var}{\op{Var}}
\newcommand{\Hom}{\op{Hom}}
\newcommand{\Oh}{\op{Coh}}
\newcommand{\Spec}{\op{Spec}}
\DeclareMathOperator{\Ext}{Ext}
\DeclareMathOperator{\Ker}{Ker}
\DeclareMathOperator{\spa}{span}

\DeclareMathOperator{\GKdim}{GKdim}
\DeclareMathOperator{\head}{top}
\DeclareMathOperator{\soc}{soc}
\DeclareMathOperator{\htt}{ht}
\DeclareMathOperator{\res}{res}
\DeclareMathOperator{\infl}{infl}
\newcommand{\Qq}[1]{Q_{#1}}
\newcommand{\coord}{S}
\newcommand{\fact}{F}
\newcommand{\topRM}{\mathrm{top}}
\newcommand{\EGT}{\cE_{\mathrm{GT}}}
\newcommand{\chiGT}{\chi_{\mathrm{GT}}}
\newcommand{\rht}{\mathrm{ht}}
\newcommand{\End}{\mathrm{End}}
\newcommand{\Seq}{\mathrm{Seq}}
\newcommand{\fmod}{\text{-}\mathrm{fmod}}
\newcommand{\Res}{\mathrm{Res}}
\newcommand{\ch}{\mathrm{ch}}
\newcommand{\KLRchar}{\mathcal{S}}
\usepackage{shuffle}
\newcommand{\trou}{\underline{\hspace{0.7em}}}
\newcommand{\ulambda}{{\underline{\lambda}}}
\newcommand{\bbslash}{\backslash\hspace{-0.3em}\backslash}
\newcommand{\umu}{\underline{\smash{\mu}}}

\newcommand{\redbold}[1]{{\color{red}\bf #1}}
\newcommand{\redboldell}{{\ensuremath{\color{red}\boldsymbol\ell}}}
\newcommand{\SL}{\mathrm{SL}}
\newcommand{\invomega}{\omega}
\newcommand{\low}{\mathrm{low}}
\newcommand{\smallheart}{\,{\scriptstyle\heartsuit}\,}
\newcommand{\timesop}{\times^{\mathrm{op}}}
\newcommand{\inverseiota}{\pi}
\newcommand{\AGHL}{\mathscr{A}}
\newcommand{\AGHLcomp}{\hat{\mathscr{A}}}

\newcommand{\Oshtot}{\hat{\O}_{sh}}

\newcommand{\dsop}{\mathbin{\diamondsuit}}
\newcommand{\hsop}{\mathbin{\heartsuit}}

\renewcommand{\simeq}{\cong}
\renewcommand{\subset}{\subseteq}

\newlength{\globalspace}
\newlength{\globalspaceeq}
\newlength{\globalspacethm}
\newlength{\globalspacesec}

\newtheoremstyle
	{main_plain}					%
	{\globalspacethm}				%
	{0pt}						%
	{\itshape}					%
	{}							%
	{\bfseries}					%
	{.}							%
	{\globalspacethm}				%
	{}							%
	
\newtheoremstyle
	{main_definition}			%
	{\globalspacethm}			%
	{0pt}						%
	{\normalfont}				%
	{}							%
	{\bfseries}					%
	{.}							%
	{\globalspacethm}			%
	{}	

\theoremstyle{main_plain}
\newtheorem{Theorem}{Theorem}[section]
\newtheorem{Corollary}[Theorem]{Corollary}
\newtheorem{Conjecture}[Theorem]{Conjecture}
\newtheorem{Question}[Theorem]{Question}
\newtheorem{Proposition}[Theorem]{Proposition}
\newtheorem{Lemma}[Theorem]{Lemma}

\theoremstyle{main_definition}
\newtheorem{Def}[Theorem]{Definition}
\newtheorem{Rem}[Theorem]{Remark}
\newtheorem{Example}[Theorem]{Example}

\newtheoremstyle
  {intro}						%
  {\globalspacethm}					%
  {0pt}							%
  {\normalfont}					%
  {}								%
  {\bfseries}					%
  {.}							%
  {\globalspacethm} 				%
  {}								%

\newtheoremstyle
  {intro_it}						%
  {\globalspacethm}					%
  {0pt}							%
  {\itshape}						%
  {}								%
  {\bfseries}					%
  {.}							%
  {\globalspacethm} 				%
  {}								%

\newcounter{intro_counter}
\newcounter{abc}

\theoremstyle{intro}
\newtheorem{intro_con}[abc]{Conjecture}
\newtheorem{intro_que}[intro_counter]{Question}
\newtheorem*{intro_rem}{Remark}
\newtheorem*{intro_ex}{Example}
\theoremstyle{intro_it}
\newtheorem{intro_thm}[intro_counter]{Theorem}
\newtheorem{intro_cor}[intro_counter]{Corollary}

\makeatletter

\renewcommand{\section}{%
  \@startsection{section}{1}			%
  {\z@}								%
  {\globalspacesec}					%
  {2pt}								%
  {\normalfont\scshape\centering}}	%

\renewcommand{\subsection}{			%
  \@startsection{subsection}{2}		%
  {\z@}								%
  {0pt}								%
  {-.5em}							%
  {\normalfont\bfseries}}			%

\renewcommand{\subsubsection}{%
  \@startsection{subsubsection}{3}	%
  {\z@}								%
  {0pt}								%
  {-.5em}							%
  {\normalfont\itshape}}				%

\makeatother

\makeatletter
\def\thm@space@setup{%
  \thm@preskip=0pt \thm@postskip=0pt
}
\makeatother

\makeatletter
	\xpatchcmd{\proof}{\topsep6\p@\@plus6\p@\relax}{}{}{}
\makeatother

\usetikzlibrary{decorations.markings}
\tikzset{double line with arrow/.style args={#1,#2}{decorate,decoration={markings,%
mark=at position 0 with {\coordinate (ta-base-1) at (0,1pt);
\coordinate (ta-base-2) at (0,-1pt);},
mark=at position 1 with {\draw[#1] (ta-base-1) -- (0,1pt);
\draw[#2] (ta-base-2) -- (0,-1pt);
}}}}

\author[J. Kamnitzer]{Joel Kamnitzer}
\address{J. Kamnitzer \\McGill University, Montréal, Québec, Canada}
\email{joel.kamnitzer@mcgill.ca}

\author[A. Labelle]{Antoine Labelle}
\address{A. Labelle \\Universität Bonn, Bonn, Nordrhein-Westfalen, Deutschland}
\email{labelle@mpim-bonn.mpg.de}

\author[A. Leroux-Lapierre]{Alexis Leroux-Lapierre}
\address{A. Leroux-Lapierre \\McGill University, Montréal, Québec, Canada}
\email{alexis.leroux-lapierre@mail.mcgill.ca}

\author[T. Pinet]{Théo Pinet}
\address{T. Pinet \\McGill University, Montréal, Québec, Canada}
\email{theo.pinet@mcgill.ca}

\author[A. Weekes]{Alex Weekes}
\address{A. Weekes\\Université de Sherbrooke, Sherbrooke, Québec, Canada}
\email{alex.weekes@usherbrooke.ca}

\begin{document}
\setlength{\medskipamount}{0pt}
\setlength{\smallskipamount}{0pt}

\setlength{\abovedisplayskip}{\globalspaceeq}
\setlength{\belowdisplayskip}{\globalspaceeq}
\setlength{\abovedisplayshortskip}{\globalspaceeq}
\setlength{\belowdisplayshortskip}{\globalspaceeq}

\setlength{\belowcaptionskip}{\globalspacethm} 			
\setlength{\textfloatsep}{\globalspacethm}

\begin{abstract}
In this paper, we study the category $\mathcal{O}$ of representations of shifted Yangians associated to a simply-laced simple Lie algebra $\mathfrak{g}$ over $\mathbb{C}$. In particular, we prove that the (complexified) Grothendieck ring of this category is isomorphic to the Cox ring of the~\textit{open bi-infinite Bott--Samelson variety}, which is a pro-variety we construct from Bott--Samelson varieties for alternating heaps. Combining this with work of Francone--Leclerc, we prove a conjecture of Hernandez--Zhang by identifying the above Grothendieck ring with a cluster algebra introduced by Geiss--Hernandez--Leclerc. Our methods also yield an action of the Langlands dual group $G^{\vee}$ on this Grothendieck ring, and show that the shifted coproducts defined in work of the first and fifth authors with collaborators give rise to coproducts for truncated shifted Yangians. This machinery then allows us to prove further conjectures of Frenkel--Hernandez and Geiss--Hernandez--Leclerc on extended $QQ$-systems, and to obtain a generalization of a duality defined by Hernandez--Leclerc.
\end{abstract}

\title[Truncated shifted Yangians and the bi-infinite Bott-Samelson variety]{Category $\O$ for truncated shifted Yangians \\ and the
bi-infinite Bott-Samelson variety}
\maketitle

\setcounter{tocdepth}{1}						
\tableofcontents

\setlength{\parskip}{\globalspace}
\section{Introduction}\label{sec:Intro}

\subsection{Shifted quantum groups and Grothendieck rings}\label{sec:introPart1}

Fix a simply-laced simple Lie algebra $\fg$ over $\C$ with Dynkin diagram $I$. In his seminal paper \cite{drinfeld1985hopf}, Drinfeld introduced the \textit{Yangian} $Y(\fg)$, an associative algebra deforming the universal enveloping algebra of the current algebra $\fg[t]$. A key property of $Y(\fg)$ is the existence of a coproduct
$$\Delta : Y(\fg)\to Y(\fg)\otimes Y(\fg),$$
which endows the category $\mathscr{C}$ of finite-dimensional $Y(\fg)$-modules with a monoidal structure, and thus makes the associated Grothendieck group $K_0(\mathscr{C})$ into a ring. The ring $K_0(\mathscr{C})$,~and its trigonometric analogue (coming from representations of quantum affine algebras),~have been extensively studied in recent years (see,~e.g.,~\cite{chari1990yangians,hernandez2021quantum,
kashiwara2024monoidal}).~In~particular,~a major breakthrough %
occurred when Hernandez--Leclerc \cite{hernandez2010cluster} proposed a geometric~interpretation of the intricate relations governing the \textit{Grothendieck ring} $K_0(\mathscr{C})$~in~terms~of~Lie-theoretic varieties and their cluster structures. In this paper, we extend this interpretation by relating the \textit{category $\mathcal{O}$} of representations of \textit{shifted Yangians} with a pro-variety we call the \textit{bi-infinite Bott--Samelson variety} and a cluster algebra introduced %
in \cite{geiss2024representations}. %
\medskip\par

The shifted Yangians $Y_{\mu}(\fg)$ are %
associative algebras that depend on a coweight $\mu\in P^\vee$ of $\fg$%
. They were defined for $\fg=\fgl_n$ and $\mu$ dominant %
by Brundan--Kleshchev \cite{brundan2006shifted},~and~were generalized to arbitrary $\fg$~and~$\mu$ by the first and fifth~authors~and~collaborators~in~\cite{braverman2016coulomb,finkelberg2018comultiplication}.  Moreover, as shown~in \cite{finkelberg2018comultiplication}, there is, for any $\mu_1,\mu_2\in P^{\vee}$, a \textit{shifted coproduct}
\begin{equation}\label{eq:intro_def_of_shifted_coprod}
\Delta_{\mu_1,\mu_2} : Y_{\mu_1 + \mu_2}(\fg) \rightarrow Y_{\mu_1}(\fg) \otimes Y_{\mu_2}(\fg)
\end{equation}
that recovers the coproduct for $Y(\fg)=Y_0(\fg)$ when $\mu_1=\mu_2=0$.\medskip\par
In \cite{hernandez2024shifted}, Hernandez--Zhang began a systematic study of the category $\O$ associated to $Y_{\mu}(\fg)$, which we denote %
by $\mathcal{O}_{\mu}$ (see Definition \ref{def:Omu} for details). By their work, simple~objects of $\O_\mu $ are characterized using their \textit{highest $ \ell$-weights}, which are $I$-tuples of rational functions $\psi = (\psi_i(u))_{i \in I}$ whose expansion at $u=\infty$ describes the action of the \textit{loop-Cartan subalgebra} of $Y_{\mu}(\fg)$ on a highest $\ell$-weight vector. Moreover, Hernandez--Zhang defined a generalization of Knight's \textit{$\ell$-character map} \cite{knight1995spectra}
\begin{equation}\label{eq:intro_ell_character}
\chi_\ell:K_0(\O_{\mu})\to \cE_\ell
\end{equation}
which sends the class of a module $V$ of $\mathcal{O}_{\mu}$ to a generating series that records the dimensions of the weight spaces of $V$ for the action of the loop Cartan subalgebra~of~$Y_{\mu}(\fg)$.
By \cite[Theorem 3.14]{hernandez2024shifted}, this map gives an injective morphism of rings
\begin{equation*}
\textstyle \chi_\ell:K_0(\bigoplus_{\mu\in P^\vee} \O_\mu)\to \cE_\ell,
\end{equation*}
where the ring structure for the left-hand side comes from the shifted coproducts \eqref{eq:intro_def_of_shifted_coprod} (and where the multiplication for $\mathcal{E}_{\ell}$ is similar to the usual product of generating series).\medskip\par
We call the category
\begin{equation*}
\O_{sh}:=\textstyle \bigoplus_{\mu\in P^\vee} \O_\mu
\end{equation*}
the \textit{category $\mathcal{O}$ for shifted Yangians}. By the above, its Grothendieck group $K_0(\O_{sh})$ has~the structure of a commutative and associative ring. 
\begin{intro_rem} 
The above discussion has a trigonometric analogue in terms of %
\textit{shifted quantum affine algebras} \cite{finkelberg2019multiplicative,hernandez2023representations}. 
Although formulated differently, this alternative framework is (at least for our purposes) interchangeable with the present setting via the equivalences~of ``truncated categories'' given in \cite[Corollary 1.2.1]{varagnolo2025representations} (see also \cite{dumanski2025k} and \cite{gautam2016yangians}).~More precisely, combining these equivalences with \cite[Theorem 5.4]{hernandez2025jordan} and \cite[Corollary 5.22]{kamnitzer2019category} %
gives a ring isomorphism between $K_0(\mathcal{O}_{sh})$ and its trigonometric counterpart (studied for example in \cite{hernandez2023representations,geiss2024representations}) that is compatible with $\ell$-characters on both sides and maps simple classes to simple classes. Accordingly, we freely pass between the trigonometric and rational (i.e.~shifted Yangian) frameworks throughout the paper.
\end{intro_rem}

We will be interested in a particular subcategory of $\O_{sh}$. Fix a bipartition on the Dynkin diagram $I=I_0\sqcup I_1$ into even and odd vertices and let
$$I\times_2\Z = \{(i,a)\in I\times \Z\,|\,i \text{ and } a\text{ have the same parity}\}. $$
We say that a highest $\ell$-weight $\psi$ is \emph{integral} if the zeros and poles of each rational function $\psi_i(u)$ are integers of the same parity as~$i$. Denote by $\O_{sh}^\Z$ the full subcategory of $\O_{sh}$~whose objects have finite length and for which each simple constituent $L(\psi)$ has an integral highest $\ell$-weight $\psi$. In Section \ref{sec:genericsimplicity}, we prove that this subcategory is closed under tensor products and that the study of the ring $K_0(\O_{sh})$ essentially reduces to the study of $K_0(\O_{sh}^\Z)$. \medskip\par

The main motivating question of our current work is:
\begin{intro_que}\label{que:intro_what_is_this}
What is the scheme $\Spec K_\C(\O_{sh}^\Z)$?
\end{intro_que}

Major progress regarding Question \ref{que:intro_what_is_this} was recently obtained by Leclerc and collaborators. A first milestone was crossed in \cite{geiss2024representations} where Geiss--Hernandez--Leclerc defined a cluster algebra $\mathcal A \subset \cE_\ell$ whose initial cluster variables are Frenkel--Hernandez's \textit{$Q$-variables} \cite{frenkel2024extended}. These elements, denoted by $Q_{w\varpi_i^\vee,a}$ with $(i,a)\in I\times_2 \Z$ and $w\in W$, are %
formal solutions~in $\cE_{\ell}$ of the \textit{extended $QQ$-system}, which is the system of equations given by
\begin{equation}\label{eq:intro_the_extended_QQ_system}
\textstyle Q_{w \varpi_i^\vee, a} Q_{w s_i \varpi_i^\vee, a+2} - Q_{w s_i \varpi_i^\vee, a} Q_{w \varpi_i^\vee, a+2} = \prod_{j \sim i} Q_{w \varpi_j^\vee, a+1},
\end{equation}
for $(i,a)\in I\times_2\Z$ and $w\in W$ such that $\ell(ws_i)>\ell(w)$. Furthermore, the properties of the $Q$-variables allowed Geiss--Hernandez--Leclerc to construct an isomorphism
\begin{equation}\label{eq:intro_isoCompleted}
\smash{\mathcal{E}_{\ell}\simeq \AGHLcomp}
\end{equation}
with $\smash{\AGHLcomp}$ a completion of the cluster algebra $\AGHL$. This led to the following conjecture, which we reformulate slightly as was done in the refinement \cite[Conjecture 6.5]{hernandez2025jordan}:
\begin{intro_con}[\cite{geiss2024representations,hernandez2025jordan}]\label{intro_con_GHL}
The isomorphism \eqref{eq:intro_isoCompleted} restricts to an isomorphism
$$ \smash{K_0(\mathcal{O}_{sh}^{\Z})\simeq \AGHL.}$$
\end{intro_con}\vspace*{-1.25mm}
Alongside Conjecture \ref{intro_con_GHL}, Geiss--Hernandez--Leclerc recalled the following conjecture, first proposed in \cite[Conjectures 6.8 and 6.11]{frenkel2024extended}:
\begin{intro_con}[\cite{frenkel2024extended,geiss2024representations}]\label{intro_con_FH}
There are explicit simple objects $L_{w\varpi_i^{\vee},a}$ such that
\begin{equation*}
\smash{Q_{w\varpi_i^{\vee},a}=\chi_{\ell}(L_{w\varpi_i^{\vee},a}).}
\end{equation*}
In particular, the classes in $K_0(\mathcal{O}_{sh}^{\Z})$ of the $L_{w\varpi_i^{\vee},a}$'s satisfy the extended $QQ$-system \eqref{eq:intro_the_extended_QQ_system}.
\end{intro_con}
\begin{intro_rem} Conjecture \ref{intro_con_FH} is shown
in \cite[Section 7]{frenkel2024extended} for rank 2 Lie algebras and~in~\cite[Sections 9.3.1--9.3.2]{geiss2024representations} assuming $\ell(w)\in\{0,1,\ell(w_0)\}$. Moreover, Conjectures \ref{intro_con_GHL}--\ref{intro_con_FH} are~part~of a broader conjecture of \cite{geiss2024representations} (known to hold for $\mathfrak{g}=\mathfrak{sl}_2$ by \cite[Section 9.4]{geiss2024representations}) which states that %
$\mathcal{O}_{sh}^{\Z}$ is a \textit{monoidal categorification} of the cluster algebra $\AGHL$. With this in mind,\newpage\noindent Conjecture \ref{intro_con_FH} can be seen as proposing a categorification for the initial cluster variables~and some initial mutation relations (namely, the extended $QQ$-system).
\end{intro_rem}

Coming back to Question \ref{que:intro_what_is_this}, another substantial contribution appeared in \cite{francone2025cluster}, %
where Francone--Leclerc constructed, for each Coxeter element $c\in W$, an algebra isomorphism
\begin{equation}\label{eq:intro_isoFL}
 \C\otimes_\Z \AGHL\simeq \C[B(G^{\vee},c)]
\end{equation}
where $B(G^\vee,c)$ is the \emph{scheme of $(G^{\vee},c)$-bands}. (Here $G^{\vee}$ is the connected simply-connected group whose Lie algebra is the Langlands dual Lie algebra $\fg^{\vee}$%
.) Points of Francone--Leclerc's scheme $B(G^\vee,c)$ are bi-infinite sequences $(g_s)_{s\in \Z}\in (G^\vee)^\Z$ such that
\begin{equation*}
g_sg_{s+1}^{-1}\in N^\vee \dot{c} \cap \dot{c}N^\vee_-
\end{equation*}
with $N^\vee,N^\vee_-$ the unipotent radicals of fixed opposite Borel subgroups $B^\vee,B_-^\vee\subseteq G^\vee$. Note that $B(G^{\vee},c)$ depends intrinsically on the choice of Coxeter element $c\in W$.

Thus, combining Conjecture \ref{intro_con_GHL} with \eqref{eq:intro_isoFL} shows that the scheme of bands $B(G^{\vee},c)$ gives~a conjectural answer to Question  \ref{que:intro_what_is_this}. We show here that this conjectural answer is correct, but argue that, while the scheme $B(G^{\vee},c)$ is well-adapted to the cluster perspective, it~is~{not~so natural from the point of view of representation theory} (more precisely from the~viewpoint of \textit{truncated shifted Yangians}). Hence, we define a new pro-variety $\mathcal{Z}_{\infty}^{\circ}$, the \textit{open bi-infinite Bott--Samelson variety},
and show that the~Cox~ring $\cR$ of this pro-variety satisfies
\begin{equation}\label{eq:intro_isoRK0}
K_{\C}(\mathcal{O}_{sh}^{\Z})\simeq \cR.
\end{equation}
Moreover, we show that any choice of Coxeter element $c$ induces an isomorphism
\begin{equation*}
\Spec \cR \simeq B(G^{\vee},c)
\end{equation*}
which allows us to think of $\Spec \cR$ as a Coxeter-independent version of Francone--Leclerc's scheme of bands. Finally, blending together \eqref{eq:intro_isoFL}--\eqref{eq:intro_isoRK0}, we obtain an algebra isomorphism
\begin{equation*}
\C\otimes_{\Z}\AGHL \simeq K_{\C}(\mathcal{O}_{sh}^{\Z}),
\end{equation*}
which we prove is the restriction of \eqref{eq:intro_isoCompleted}. This leads to the proof of Conjecture \ref{intro_con_GHL} (Theorem \ref{thm:conjecture_HZ}). Along the way, we also prove Conjecture \ref{intro_con_FH} (Theorem \ref{thm:ConjFH}) using the tools described below (which also lead naturally to a more precise definition of $\cZ_{\infty}^{\circ}$).

\subsection{Truncated shifted Yangians}\label{sec:intro_TSY}

The first and fifth authors and collaborators defined in \cite{kamnitzer2014yangians,braverman2016coulomb} distinguished quotients of the shifted Yangians $Y_{\mu}:=Y_{\mu}(\fg)$. These quotients, called \textit{truncated shifted Yangians} and denoted $Y_{\mu}^{\la}(\bR)$, %
depend on the additional data of a dominant coweight $ \lambda = \sum_{i\in I} \lambda_i \varpi_i \in P_+^{\vee}$ and a \textit{set of parameters} $\bR=(R_i)_{i\in I}\in \Z^{\la}$
(which~is~simply an $I$-tuple of multisets $R_i$ of size $\la_i$ with integer entries of the same parity as $i$). They are constructed using an algebra morphism called the \textit{GKLO map} (see Section \ref{section: tsy} for details) and are quantized Coulomb branch algebras by the results of \cite{braverman2016coulomb}.

The \textit{truncations} $Y^\la_\mu(\bR)$ %
quantize slices in the affine Grassmannian of $G$ (where $G$ denotes the adjoint group with Lie algebra $ \fg$). More precisely, results of \cite{kamnitzer2014yangians,braverman2016coulomb,weekes2019} %
show that the associated graded algebra $\operatorname{gr} Y^\la_\mu(\bR)$ %
satisfies
$$\operatorname{gr} Y^\la_\mu(\bR) \cong \C[\overline{\mathcal{W}}{}^\la_\mu],$$
where $\overline{\mathcal{W}}{}^\la_\mu$ is the generalized affine Grassmannian slice attached to the pair $(\la,\mu)$ (which~can also be seen as the Coulomb branch of the corresponding quiver gauge theory).\par 
In addition, %
one can define a \textit{characteristic cycle map} \cite{Alexis}
\begin{equation}\label{eq:intro_CC}
\op{CC}:K_0(\O^\la_\mu(\bR)) \rightarrow \op{H}_{\text{top}}(\overline{\Gr}{}^\la \cap S^\mu),
\end{equation}
where $ \Gr^\lambda $ is a spherical Schubert cell in the affine Grassmannian, $ S^\mu $ is a semi-infinite~orbit, and $\mathcal{O}_{\mu}^{\la}(\bR)\subseteq \mathcal{O}_{\mu}$ is the full subcategory of objects annihilated by the kernel of the defining map $Y_{\mu}\twoheadrightarrow Y_{\mu}^{\la}(\bR)$. The codomain of \eqref{eq:intro_CC} obtains a representation theoretic interpretation under the celebrated \textit{geometric Satake correspondence}, which provides an isomorphism
$$\C\otimes_\Z\op{H}_{\text{top}}(\overline{\Gr}{}^\la \cap S^\mu)\simeq V(\la)_\mu,$$
where $ V(\la) $ is the irreducible module of highest~weight $\la$ for the Langlands dual group $G^\vee$. Letting $ \O^\la_{sh}(\bR) := \bigoplus_\mu \O^\la_\mu(\bR)$, the characteristic cycle gives a map
$$ K_{\C}(\smash{\mathcal{O}_{sh}^{\la}(\bR)})\to V(\la),$$
and %
it is therefore natural to expect a categorical $\fg^\vee $-action on $\O^\la_{sh}(\bR)$. Such an action was constructed in \cite{kamnitzer2019category} using an equivalence relating the \textit{truncated categories} $\mathcal{O}_{\mu}^{\la}(\bR)$ to categories of finite-dimensional modules over \textit{parity KLRW algebras}. This makes
$$ \mathsf{V}(\lambda,\bR):=\smash{K_\C(\O^\la_{sh}(\bR))}$$
into a $G^{\vee}$-module, which is known to satisfy
\begin{equation}\label{eq:intro_VlaRincl}
\textstyle V(\la)\subseteq \mathsf{V}(\la,\bR)\subseteq \bigotimes_{i\in I}V(\varpi_i^{\vee})^{\otimes \la_i}.
\end{equation}
Furthermore, the results of \cite{kamnitzer2019category} show that the highest $\ell$-weights of simple objects of $\O^\la_{sh}(\bR)$ are given by the \textit{product monomial crystals} of Section \ref{sec:intro_DecompCrystal}. This~shows in particular that each simple object of $\O_{sh}^\Z$ lies in one of the %
subcategories $\O^\la_{sh}(\bR)$. Th%
us, the inclusions
\begin{equation}\label{eq:intro_VlambdaR_in_KOsh}
\mathsf{V}(\lambda,\bR)\subset K_\C(\O_{sh}^\Z)%
\end{equation}
give a filtration of %
$K_\C(\O_{sh}^\Z)$.  Our scheme $\cZ^\circ_{\infty}$ is constructed so that the subspaces $\mathsf{V}(\la,\bR)$ come from sections of line bundles on a compactification $ \cZ_\infty $.

\subsection{Combining structures}

Summarizing parts of the two previous sections, we see that %
there are two structures at play at the level of Grothendieck rings:
\begin{enumerate}[label=(\roman*)]
\item\label{item:intro_ring_structure_from_coprod} a ring structure on $K_0(\O_{sh}) $ coming from the shifted coproducts $\Delta_{\mu_1, \mu_2}$, and
\item\label{item:intro_Gvee_action_from_trunc} a collection of subspaces $K_\C(\O^\la_{sh}(\bR)) \subset K_\C(\O_{sh}^\Z) $, each carrying a $G^\vee$-action.
\end{enumerate}
We show that \ref{item:intro_ring_structure_from_coprod} ``localizes'' while \ref{item:intro_Gvee_action_from_trunc} ``globalizes'' as follows.
\subsubsection{Truncated shifted coproducts}

For \ref{item:intro_ring_structure_from_coprod}, we establish that shifted coproducts restrict to truncations. More precisely, we show the result below, which solves a conjecture~%
of~\cite{finkelberg2018comultiplication}.

\begin{intro_thm}\label{thm:truncoprod}
Let $\la_1, \la_2\in P_+^{\vee}$ and fix $\mu_1, \mu_2\in P^{\vee}$ such that both $V(\la_1)_{\mu_1}$ and $V(\la_2)_{\mu_2}$ are non-zero. Then, for any $(\bR_1, \bR_2)$, there exists a unique algebra morphism
\begin{equation}\label{eq:intro_TSC}
Y^{\la_1+\la_2}_{\mu_1+\mu_2}(\bR_1 \cup \bR_2) \longrightarrow Y^{\la_1}_{\mu_1}(\bR_1) \otimes  Y^{\la_2}_{\mu_2}(\bR_2)
\end{equation}
making the diagram\vspace*{-0.75mm}
$$\adjustbox{scale=0.9,center}{
\begin{tikzcd}[row sep =1em]
Y_{\mu_1+\mu_2} \ar[r,"\Delta_{\mu_1,\mu_2}"] \ar[d,two heads] & Y_{\mu_1}\otimes Y_{\mu_2} \ar[d,two heads] \\
Y^{\la_1+\la_2}_{\mu_1+\mu_2}(\bR_1 \cup \bR_2) \ar[r,dashed] & Y_{\mu_1}^{\lambda_1}(\bR_1) \otimes Y_{\mu_2}^{\lambda_2}(\bR_2)
\end{tikzcd}
}$$
commutative (where the vertical arrows are the defining projections).
\end{intro_thm}\newpage
The proof appears in Section \ref{sec:TSC}. First, we show that it suffices to prove the result for a Zariski dense set of pairs $(\mathbf{R}_1, \mathbf{R}_2)$, utilizing the fact that the algebras $Y_\mu^\la(\mathbf{R})$~are~fibres~of~a well-behaved family $Y_\mu^\la$
of algebras over the base $\C^\la=\prod_{i\in I} \C^{\lambda_i}/\Sigma_{\lambda_i}
 $. The main idea of the proof is then to describe the defining ideal $\Ker(Y_\mu \twoheadrightarrow Y_\mu^\la(\mathbf{R}))$ in terms of modules~in~$\mathcal{O}_\mu^\la(\mathbf{R})$ and then to exploit the generic simplicity of products $L(\psi_1) \otimes L(\psi_2)$
established in Section \ref{sec:genericsimplicity} (via new $R$-matrices in the category $\mathcal{O}_{sh}$).
We expect that the hypotheses about $V(\la_1)_{\mu_1}$ and
$V(\la_2)_{\mu_2}$ can be dropped, see Remark \ref{rem:truncoprodA}. We also conclude in Remark \ref{rem:injective} that the maps \eqref{eq:intro_TSC} are injective and that they quantize the multiplication maps
of generalized affine Grassmannian slices defined in \cite[Section 2(vii)]{braverman2016coulomb}.\medskip\par
A direct consequence of Theorem \ref{thm:truncoprod} is that tensor product gives a multiplication map
\begin{equation}\label{eq:intro_mult_of_trunc}
\mathsf{V}(\la_1,\bR_1)\otimes \mathsf{V}(\la_2,\bR_2)\to \mathsf{V}(\la_1+\la_2,\bR_1\cup\bR_2).
\end{equation}
\subsubsection{$G^\vee$-action on $K_{\C}(\O_{sh}^{\Z})$}

For the structure \ref{item:intro_Gvee_action_from_trunc}, we prove that the different $G^\vee$-actions on the subspaces $\mathsf{V}(\la,\bR)\subseteq K_{\C}(\mathcal{O}_{sh}^{\Z})$ glue into a well-defined global action:
\begin{intro_thm}[Theorems \ref{thm:construction_of_gaction_on_KOsh} and \ref{thm:mult_in_KOsh_is_G_equiv}]\label{thm:intro_glueing}
There exists a unique $G^\vee$-action on $K_\C(\O_{sh}^\Z)$ making $G^{\vee}$-equivariant each inclusion
\begin{equation*}
\mathsf{V}(\la,\bR)\subseteq K_{\C}(\mathcal{O}_{sh}^{\Z}).
\end{equation*}
Moreover, the multiplication in $ K_\C(\O_{sh}^\Z) $ is $G^\vee$-equivariant.
\end{intro_thm}
Our proof of this result appears in Section \ref{sec:glueing_of_g_action} and is somewhat indirect. Indeed, instead of studying the gluing of the different categorical $G^{\vee}$-actions on the categories $\O^\la_{sh}(\bR)$,~we relate $\ell$-characters to a notion of characters for modules over parity KLRW algebras. %

\subsection{The bi-infinite Bott--Samelson variety}\label{sec:intro_biinfinite}
Alongside Theorem \ref{thm:intro_glueing}, we prove in Section \ref{sec:glueing_of_g_action} that $K_\C(\O_{sh}^\Z)$ is generated by the fundamental subspaces
\begin{equation*}
\{K_\C(\O^{\varpi_i^{\vee}}_{sh}(a))=\mathsf{V}(\varpi_i^{\vee},a)\simeq V(\varpi_i^{\vee})\}_{(i,a)\in I\times_2\Z}.
\end{equation*}
Hence, as the multi-homogeneous coordinate ring of the Grassmannian $ \Gr(i):=P_i^{-} \backslash G^\vee$ (with $P_i^{-}$ the maximal parabolic subgroup associated to $i\in I$) is $\bigoplus_{k\geq 0} V(k\varpi_i^\vee)$, we deduce
that $\Spec K_\C(\O_{sh}^\Z) $ embeds into the infinite product of cones (see Section \ref{subsec:cox_ring_of_Zinftycirc})
\begin{equation}\label{eq:intro_prodconesGr}
\textstyle \prod_{(i,a)\in I\times_2\Z} \widehat{\Gr}(i).
\end{equation}
We introduce the \emph{bi-infinite Bott--Samelson variety} to describe the image of $\Spec  K_\C(\O_{sh}^\Z)$ in the product \eqref{eq:intro_prodconesGr}. This pro-variety, denoted $\mathcal{Z}_{\infty}$,  parametrizes collections $(x_{i,a})_{(i,a)\in I\times_2\Z}$ satisfying $x_{i,a}\in \Gr(i)$ for all $(i,a)$ as well as an incidence condition for each pair $(x_{i,a}, x_{j, a+1})$ with $i\sim j$. An important feature of $\mathcal{Z}_{\infty}$ is that, for every choice of \textit{height function} $\xi:I\to \Z$ (i.e.~$\xi$ satisfies $|\xi_i-\xi_j|=1$ whenever $i\sim j$ and $i$ always has the same parity as $\xi_i$), there is a projection onto the flag variety
\begin{equation}\label{eq:intro_proj_onto_Fl}
\textstyle \cZ_\infty\twoheadrightarrow \Fl=B^\vee_-\backslash G^\vee\hookrightarrow \prod_{i\in I}\Gr(i)
\end{equation}
given by selecting elements indexed by $\{(i,\xi_i)\}_{i\in I}$ in $\prod_{(i,a)\in I\times_2\Z}\Gr(i)$.

\begin{intro_ex}
For $G^\vee=\SL_n$, the pro-variety $\mathcal{Z}_{\infty}$ parametrizes collections of subspaces of $\C^n$ $(W_{i,a})_{(i,a)\in I\times_2\Z}$ such that $W_{i,a}\subset W_{i+1,a\pm 1}$ for all $(i,a)\in I\times_2\Z$ (see Figure \ref{fig:dessin_bella_et_diag_biinfini}).
\begin{figure}[ht]
\centering
\begin{subfigure}[c]{0.5\textwidth}
  \centering
  \includegraphics[width=0.7\linewidth]{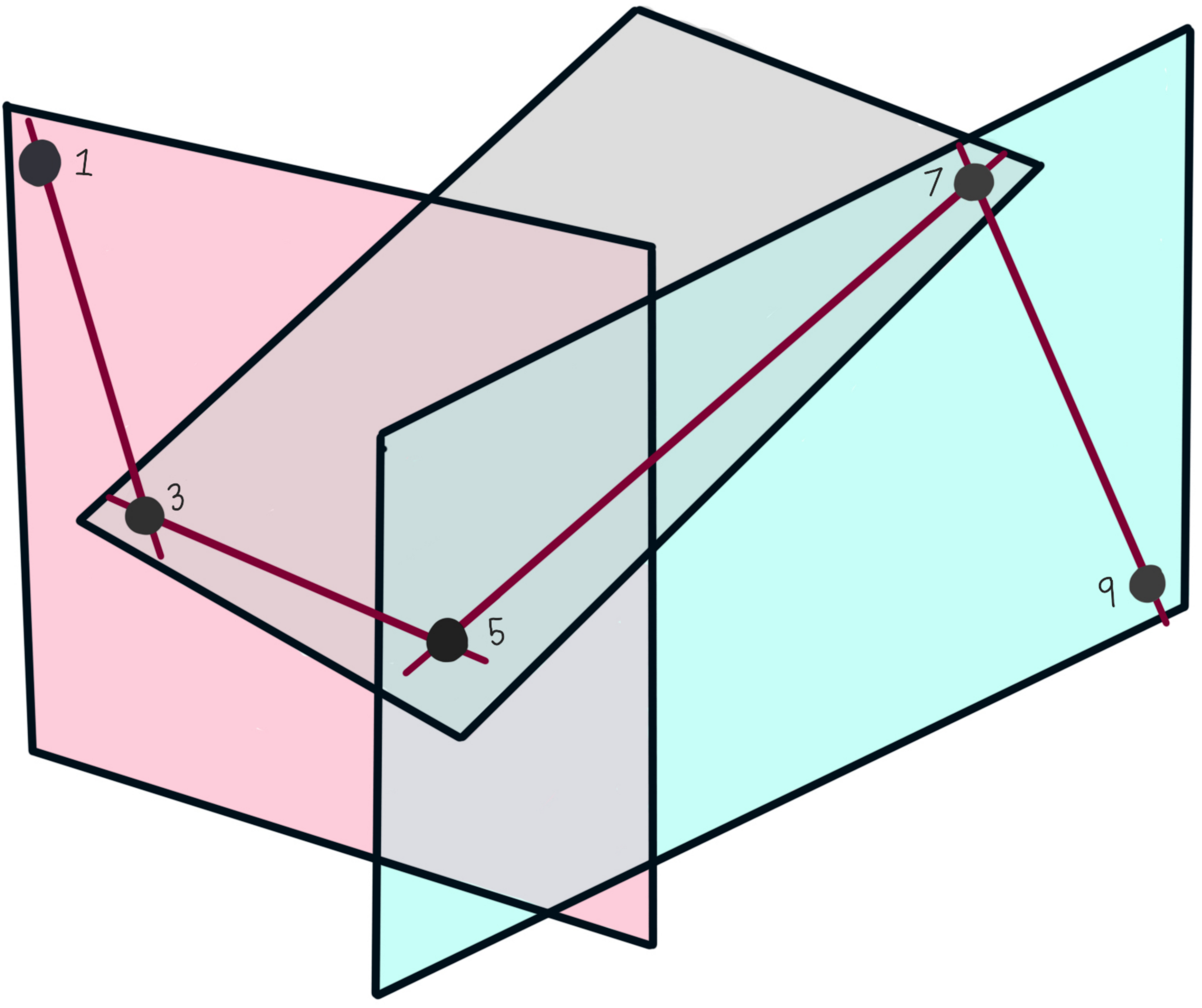}
\end{subfigure}%
\begin{subfigure}[c]{.5\textwidth}
  \centering
  \begin{tikzpicture}[scale=0.65,every node/.style={scale=0.95}]
\def\rad{0.5}
\def\xnudge{0.8}
\def\ynudge{0.8}
\def\nbVertices{5}
\foreach \i in {2,4} {
	\foreach \x in {-2,0,2,4} {
			\ifnum\x<4
				\draw (\i,\x) -- (\i+1,\x+1);	
			\else
			\fi
			\draw (\i,\x) -- (\i+1,\x-1);
	}
}
\foreach \i in {1,3} {
	\foreach \x in {-3,-1,1,3} {
		\ifnum \x>-3
			\draw (\i,\x) -- (\i+1,\x-1);
		\else
		\fi
		\draw (\i,\x) -- (\i+1,\x+1);

	}
}

\foreach \i in {1,3,5} {
	\foreach \x in {-1,1,3} {
		\fill[fill=white] (\i,\x) circle (\rad);
		\ifnum\x<5
			\node at (\i,\x) {$W_{\i,\x}$};
		\else
		\fi
	}
	\fill[fill=white] (\i,-3) circle (\rad);
	\node at (\i,-3) {$\vdots$};
}

\foreach \i in {2,4} {
	\foreach \x in {-2,0,2} {
		\fill[fill=white] (\i,\x) circle (\rad);
		\ifnum\x>-4
			\node at (\i,\x) {$W_{\i,\x}$};
		\else
		\fi
	}
	\fill[fill=white] (\i,4) circle (\rad);
	\node at (\i,4) {$\vdots$};
}

\def\hI{-4}
\def\sizeV{0.05}
\def\xZ{6.5}
\def\ymaxZ{4}
\def\yminZ{-3}

\draw (1,\hI) -- (\nbVertices,\hI);
\foreach \i in {1,...,\nbVertices} {
	\fill (\i,\hI) circle (\sizeV) node[below] {\footnotesize $\i$};
}

\node at (1-0.75,\hI) {$I$};

\draw[->] (\xZ,\yminZ)--(\xZ,\ymaxZ) node [above right] {$\Z$};

\foreach \n in {-2,...,3} {
	\draw (\xZ-0.1,\n)-- (\xZ+0.1,\n);
}
\end{tikzpicture}
\end{subfigure}
\caption{Part of a projectivized point of $\cZ_\infty$ for $G^\vee=\SL_4$ (left) and the bi-infinite diagram parametrizing points of $\cZ_{\infty}$ for $G^\vee=\SL_6$ (right).}
\label{fig:dessin_bella_et_diag_biinfini}
\end{figure}
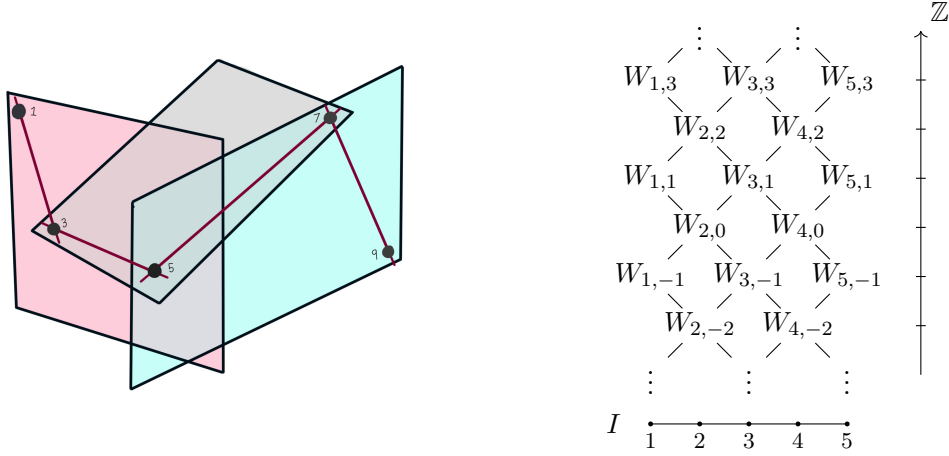
\end{intro_ex}\vspace*{-3mm}

For every $(i,a)\in I\times_2\Z$, the pro-variety $\cZ_\infty$ carries a $G^\vee$-equivariant line bundle $\O_{i,a}(1)$ obtained by pulling back the antitautological line bundle $\O(1)$ on $ \Gr(i)$. Thus, given any set of parameters $\bR$ of size $\lambda$, we can define a $G^\vee$-equivariant line bundle $\cL_{\la,\bR}$ by tensoring the line bundles $\O_{i,a}(1)$ as $i$ ranges over $I$ and $a$ ranges over the elements of $R_i$.
Combining Gibson's~Demazure~character formula \cite{gibson2021demazure} with the theory of Bott--Samelson varieties, we prove:

\begin{intro_thm}[Theorem \ref{th:charcR}]\label{thm:intro_char_of_rep_from_sections}
For any set of parameters $\bR$, the $G^\vee$-module $\op{H}^0(\cZ_\infty, \cL_{\la,\bR})$ is isomorphic to $\mathsf{V}(\lambda,\bR)=K_{\C}(\mathcal{O}_{sh}^{\la}(\bR))$.
\end{intro_thm}
\begin{intro_rem}
Theorem \ref{thm:intro_char_of_rep_from_sections} can be understood as a Borel--Weil description of the module $\mathsf{V}(\lambda,\bR)$ since specializing to line bundles pulled back along \eqref{eq:intro_proj_onto_Fl} recovers the Borel--Weil construction of simple $G^{\vee}$-modules as spaces of sections of line bundles over the flag variety $\Fl$.
\end{intro_rem}

We define the \textit{open bi-infinite Bott--Samelson variety} %
as the open subvariety~of~$\mathcal{Z}_{\infty}$~given by the condition $ x_{i,a} \ne x_{i,a+2} $ for every $ (i,a) \in I \times_2 \Z $.
Restricting the line bundles $\cL_{\la,\bR}$~%
to $ \cZ_\infty^\circ $ %
gives injective $G^\vee$-equivariant maps
\begin{equation}\label{eq:intro_injection_given_by_restriction}
\op{H}^0(\cZ_\infty, \cL_{\la,\bR}) \hookrightarrow \op{H}^0(\cZ_\infty^\circ,\cL_{\la,\bR}\vert_{\cZ_\infty^\circ}).
\end{equation}
However, many of the line bundles $\cL_{\la,\bR}$ become isomorphic under \eqref{eq:intro_injection_given_by_restriction}, as the Picard group $\cP:=\Pic(\cZ_\infty^\circ)$ is a free abelian group of finite rank. We are therefore led to the study~of~the Cox ring of the open bi-infinite Bott--Samelson variety
\begin{equation}\label{eq:intro_def_of_cox}
\cR:=\textstyle\bigoplus_{\tau\in \cP} \op{H}^0(\cZ_\infty^\circ,\cL_{\tau})
\end{equation}
and of the scheme
\begin{equation*}
\widehat{\cZ}_\infty^\circ:=\Spec\cR,
\end{equation*}
which we
characterize as the universal principal torus bundle $\widehat{\cZ}_\infty^\circ\to \cZ_\infty^\circ$ for the torus $A$ whose character lattice is $\cP$. This is precisely the scheme mentioned in Section~\ref{sec:introPart1}~that~gives a Coxeter-independent presentation of Francone--Leclerc's scheme of bands.

Unwrapping the construction, we see that $\widehat{\cZ}_\infty^\circ$ is the scheme which parametrizes arrays
$(x_{i,a})$ in the infinite product \eqref{eq:intro_prodconesGr} satisfying incidence conditions with some explicit normalization constraints (see Definition \ref{def:Zhat} for details). Moreover, the variety $\widehat{\cZ}_\infty^\circ$ is built~in such a way that, for each choice of height function $\xi:I\to \Z$, there is a map $\widehat{\cZ}_\infty^\circ\twoheadrightarrow N^\vee_-\backslash G^\vee$ which fits into a Cartesian square
\begin{equation}\label{eq:intro_cart_square_base_affine}
\adjustbox{scale=0.9}{
\begin{tikzcd}[row sep =1.5em]
\widehat{\cZ}_\infty^\circ \arrow[d,two heads] \arrow[r]& \cZ_\infty^\circ\arrow[d, two heads]\\
N^\vee_-\backslash G^\vee \arrow[r] & B^\vee_-\backslash G^\vee
\end{tikzcd}}
\end{equation}
where the left surjection is also given by selecting elements of the array $(x_{i,a})$ according to $\xi$. Notice that this choice of height function identifies the tori $A$ and $T^\vee$.

From the definition of $\cZ_\infty^\circ$, we are able to deduce:
\begin{intro_thm}[Corollary \ref{cor:presentation-of-R}]\label{thm:intro_presentation_of_biinfinite_bott}
The Cox ring $\cR$ is the quotient of the polynomial ring
\begin{equation*}
\bigotimes_{(i,a)\in I\times_2\Z} \op{Sym}^\bullet V(\varpi_i^{\vee},a),
\end{equation*}
(where $V(\varpi_i^{\vee},a)$ is a copy of $V(\varpi_i^{\vee})$ indexed by $a$)
by an ideal generated by four families~of explicit relations.
\end{intro_thm}

The above theorem suggests that identifying the $G^{\vee}$-modules $V(\varpi_i^{\vee},a)$ and
\begin{equation*}
\mathsf{V}(\varpi_i^{\vee},a)=K_{\C}(\mathcal{O}^{\varpi_i^{\vee}}_{sh}(a)) \subseteq K_{\C}(\mathcal{O}_{sh}^{\Z})
\end{equation*}
should extend to an algebra morphism $\Omega:\cR\to K_{\C}(\mathcal{O}_{sh}^{\Z})$. To prove that this is indeed the case, we need to check that the four families of relations hinted at in Theorem \ref{thm:intro_presentation_of_biinfinite_bott} hold in $K_{\C}(\mathcal{O}_{sh}^{\Z})$. This leads~to~the study of extended $QQ$-systems and chamber modules.

\subsection{Chamber modules and extended QQ-system}
For $(i,a)\in I\times_2\Z$ and $\gamma\in W\varpi_i^{\vee}$, the category $\O_{\gamma}^{\varpi_i^\vee}(a)$ contains a unique simple object (up to isomorphism), which we denote $L_{\gamma,a}$ and call \textit{chamber module} (since it depends on a chamber weight $\gamma$).
\begin{intro_rem} The family of chamber modules contains all positive and negative prefundamental representations of \cite{hernandez2012asymptotic,hernandez2024shifted}. Also, if $\gamma=w\varpi^\vee_i$ for $w\in W$, then $L_{w\varpi^\vee_i,a}$ is precisely the module appearing in Conjecture \ref{intro_con_FH} (see Lemma \ref{lem:highest_ell_weight_chamber_module_braid}).
\end{intro_rem}
To study the relations tying together the classes of chamber modules in $K_0(\O_{sh})$,~we~use the identity
\begin{equation*}
\textstyle\varpi_i^\vee+s_i\varpi_i^\vee=2\varpi_i^\vee-\alpha_i^\vee=\sum_{j\sim i}\varpi_j^\vee
\end{equation*}
which holds since we work with a simply-laced Dynkin diagram. An immediate consequence of this is that the simple $G^\vee$-representation $V(2\varpi_i^\vee-\alpha_i^\vee)$ embeds
in the two tensor products $\bigotimes_{j\sim i} V(\varpi_j^\vee)$ and $V(\varpi_i^\vee)\otimes V(\varpi_i^\vee)$ via
\begin{equation*}
\textstyle v_{2\varpi_i^\vee-\alpha_i^\vee}\mapsto \bigotimes_{j\sim i} v_{\varpi_j^\vee}%
\end{equation*}
and 
\begin{equation*}
v_{2\varpi_i^\vee-\alpha_i^\vee}\mapsto v_{\varpi_i^\vee}\otimes v_{s_i\varpi_i^\vee}-v_{s_i\varpi_i^\vee}\otimes v_{\varpi_i^\vee},
\end{equation*}
respectively. Using this, we show (see Section \ref{subsec:cat_O_of_QQ} for details on the notation and the maps):

\newpage
\begin{intro_thm}[Theorem \ref{thm:the_heptagon_commutes}]\label{thm:intro_heptagon}
The heptagon of $G^{\vee}$-equivariant maps
\begin{equation*}
\adjustbox{scale=0.9,center}{
\begin{tikzpicture}[
  baseline=(current bounding box.center),
  every node/.style={inner sep=2pt},
  arr/.style={->, shorten >=4pt, shorten <=4pt}
]
  \node (A1) at (0,2.3) {$V(2\varpi_i^{\vee}-\alpha_i^{\vee})$};

  \node (A2) at (-3.6,1.5) {$ \bigotimes_{j\sim i} V(\varpi_j^{\vee})$};
  \node (A7) at (3.6,1.5) {$V(\varpi_i^{\vee})\otimes V(\varpi_i^{\vee})$};

  \node (A3) at (-3.6,0.1) {$ \bigotimes_{j\sim i} K_\C(\O_{sh}^{\varpi_j^{\vee}}(a+1))$};
  \node (A6) at (3.6,0.1) {$K_\C(\O_{sh}^{\varpi_i^{\vee}}(a))\otimes K_\C(\O_{sh}^{\varpi_i^{\vee}}(a+2))$};

  \node (A4) at (-2.4,-1.15)
  {$K_\C(\O^{2\varpi_i^{\vee}-\alpha_i^{\vee}}_{sh}(a+1))$};

\node (A5) at ( 2.4,-1.15)
  {$K_\C(\O^{2\varpi_i^{\vee}}_{sh}(a,a+2))$};

  \draw[arr] (A1) -- (A2);
  \draw[arr] (A1) -- (A7);

  \draw[arr] (A2) -- (A3);
  \draw[arr] (A7) -- (A6);

  \draw[arr] (A3) -- (A4);
  \draw[arr] (A6) -- (A5);

  \draw[arr] (A4) -- (A5);
\end{tikzpicture}
}
\end{equation*}
commutes.
\end{intro_thm}
In particular, each extremal vector in $V(2\varpi_i^{\vee}-\alpha_i^{\vee})$ gives a relation in $K_{\C}(\O_{sh}^{2\varpi_i^{\vee}}(a,a+2))$ by travelling left and right in the above commutative diagram. This gives:

\begin{intro_thm}[Theorem \ref{thm:extended_QQ_holds}]
For $(i,a)\in I\times_2\Z$ with $w\in W$ such that $\ell(ws_i)>\ell(w)$,
\begin{equation*}
\textstyle [L_{w \varpi_i^\vee, a}][L_{w s_i \varpi_i^\vee, a+2}] - [L_{w s_i \varpi_i^\vee, a}][L_{w \varpi_i^\vee, a+2}] = \prod_{j \sim i} [L_{w \varpi_j^\vee, a+1}]
\end{equation*}
in $K_0(\O_{sh})$. Hence chamber modules solve the extended $QQ$-system.
\end{intro_thm}
Moreover, via uniqueness results for solutions of the extended $QQ$-system inspired from \cite[Section 2]{frenkel2022weyl}, we deduce a proof of Conjecture \ref{intro_con_FH}, that is:
\begin{intro_thm}[Theorem \ref{thm:ConjFH}]\label{thm:intro_conFH} For $(i,a)\in I\times_2\Z$ with $w\in W$, we have, in $\cE_{\ell}$,
$$\chi_{\ell}(L_{w\varpi^\vee_i,a})=Q_{w\varpi^\vee_i,a}.$$
\end{intro_thm}
The relations coming from the diagram of Theorem~\ref{thm:intro_heptagon} correspond to the image in $K_{\C}(\mathcal{O}_{sh}^{\Z})$ of one of the families of defining relations %
of Theorem~\ref{thm:intro_presentation_of_biinfinite_bott}. For the other relations, we need:

\subsection{Height functions, the base affine space, and the dual canonical basis}\label{sec:intro_BaseAff}

Two of the other families of relations alluded to in Theorem \ref{thm:intro_presentation_of_biinfinite_bott} arise naturally from the geometry of the base affine space. As explained in \eqref{eq:intro_cart_square_base_affine}, each choice of height function gives a surjection $\widehat{\cZ}_\infty^\circ\twoheadrightarrow N^\vee_-\backslash G^\vee$. On the level of functions, this induces an inclusion
\begin{equation}\label{eq:intro_incBaseAff}
\textstyle\bigoplus_{\lambda\in P_+^\vee} V(\lambda)\simeq\C[N^\vee_-\backslash G^\vee]\hookrightarrow \C[\widehat{\cZ}_\infty^\circ]%
.
\end{equation}
We give a categorical counterpart of the above map. Let $\O^\lambda_{sh}(\xi)$ be the truncated category corresponding to the set of parameters
$\bR$ whose $i^{\text{th}}$-component $R_i$ consists of $\lambda_i$ copies of %
$\xi_i$ for all $i\in I$. Using the representation theory of cyclotomic KLR algebras, we deduce:

\begin{intro_thm}[Theorem \ref{thm:subalgebra_iso_base_affine_space}]\label{thm:intro_baseAff}
For each choice of a height function $\xi$, the full subcategory
$$\textstyle \O_{sh}(\xi):=\bigoplus_{\lambda\in P_+^\vee} \O^\lambda_{sh}(\xi)$$
is closed under tensor products.
Moreover, there is a 
$G^{\vee}$-equivariant inclusion
$$ \C[N_-^{\vee}\backslash G_-^\vee]\simeq K_{\C}(\mathcal{O}_{sh}(\xi))\hookrightarrow K_{\C}(\mathcal{O}_{sh}^{\Z})$$
that recovers \eqref{eq:intro_incBaseAff} via the map $\Omega$ of Section \ref{sec:intro_biinfinite}.
\end{intro_thm}

From this result, we show Corollary \ref{co:Rel12} which leads to the proof that two of the families of relations appearing in Theorem \ref{thm:intro_presentation_of_biinfinite_bott} hold in $K_{\C}(\mathcal{O}_{sh}^{\Z})$. We also use the relation between~the subcategory $\mathcal{O}_{sh}(\xi)\subseteq \mathcal{O}_{sh}$ and cyclotomic KLR algebras to show:

\begin{intro_cor}
For each choice of a height function $\xi$, the image of Lusztig's dual canonical basis under the composite
\begin{equation*}
V(\lambda)\simeq K_\C(\O^\lambda_{sh}(\xi))\hookrightarrow K_\C(\O_{sh}^\Z)
\end{equation*}
belongs to the basis of simple modules.
\end{intro_cor}

Combining this with our main results, the above corollary proves part of a conjecture of Francone--Leclerc concerning the dual canonical basis of $\C[G^\vee]$ (see \cite[Section 8.5]{francone2025cluster}).

\begin{intro_rem} Theorem \ref{thm:intro_baseAff} above generalizes. Indeed, we introduce in Section \ref{subsec:subalgebras_of_KOsh} subcategories $\mathcal{O}_{sh}(\xi,\xi')\subseteq \mathcal{O}_{sh}^{\Z}$ that are indexed by ordered pairs of height functions $(\xi,\xi')$, and are closed under tensor products. Moreover, we identify the geometric counterpart of these categories via the morphism $\Omega$. Notable cases are $\mathcal{O}_{sh}(\xi,\xi)=\mathcal{O}_{sh}(\xi)$ (which categorifies $\C[N_-^{\vee}\backslash G^{\vee}]$ by Theorem \ref{thm:intro_baseAff}) and $\mathcal{O}_{sh}(\xi,\xi^*+h)$ (cf.~Section \ref{sec:Notation} for notation), which categorifies $\C[G^{\vee}]$.
\end{intro_rem}

\subsection{The isomorphism $\Omega$} Using the results of the above sections, it is not hard to deduce that the map $\Omega:\cR\to K_{\C}(\mathcal{O}_{sh}^{\Z})$ indeed defines an algebra morphism. Moreover, as $K_0(\mathcal{O}_{sh}^{\Z})$ is generated by its submodules $\{\mathsf{V}(\varpi_i^{\vee},a)\}_{(i,a)\in I\times_2\Z}$ (see the beginning of Section \ref{sec:intro_biinfinite}), this morphism is surjective. We show in Section \ref{sec:iso_between_KOsh_and_R} that it is also injective.~More~precisely, we give two separate proofs of this injectivity: one via Geiss--Hernandez--Leclerc's~isomorphism \eqref{eq:intro_isoCompleted} and the other relying on Theorem \ref{thm:intro_char_of_rep_from_sections}. This finishes the proof of %
Conjecture \ref{intro_con_GHL} and answers Question \ref{que:intro_what_is_this}. \medskip\par Summarizing everything, the main result of this paper can be stated as:

\begin{intro_thm}[Theorems \ref{thm:Omega_is_surj} and \ref{thm:injectivity-geo}]\label{intro_Omega_iso}
The morphism $\Omega:\cR\to K_\C(\O_{sh}^\Z)$ of Section \ref{sec:intro_biinfinite}~is a $G^\vee$-equivariant isomorphism of $\C$-algebras. Furthermore, it restricts to isomorphisms of $G^\vee$-representations 
$$\op{H}^0(\cZ_\infty, \cL_{\la,\bR}) \simeq \mathsf{V}(\la,\bR)=K_{\C}(\mathcal{O}_{sh}^{\la}(\bR))$$ 
for each $\lambda$ and $\bR$.
\end{intro_thm}

We now turn to a brief description of some byproducts of our work.
\subsection{Inclusions of crystals and decomposition of $\mathcal{O}_{sh}^{\Z}$}\label{sec:intro_DecompCrystal}
As hinted at in Section \ref{sec:intro_TSY},~the highest $\ell$-weights $\cB(\lambda,\bR)$ of simple objects of the truncated category $\mathcal{O}_{sh}^{\la}(\bR)$ form a normal $\fg^{\vee}$-crystal called the \textit{product monomial crystal}. Moreover, the corresponding $G^{\vee}$-module is precisely $\mathsf{V}(\la,\bR)$, and \eqref{eq:intro_VlaRincl} is equivalent to
\begin{equation*}
B(\la) \subseteq \cB(\la, \bR) \subseteq \textstyle\bigotimes_{i \in I} B(\varpi_i^\vee)^{\otimes \lambda_i},
\end{equation*}
where $B(\la)$ is the crystal of the irreducible representation $V(\la)$. In Theorem \ref{thm:charaterization_max_sing_crystals}, we show that the first inclusion is an equality if and only if %
$\bR$ is chosen according to~a~height~function (as in Section \ref{sec:intro_BaseAff}). \medskip\par
The product monomial crystal is a subcrystal of \textit{Nakajima's monomial crystal} $\cB$, which, as a set, is the collection  of Laurent monomials in variables $\{y_{i,a}\}_{(i,a)\in I\times_2\Z}$. Clearly, $\cB$ is a group. We let $\Gamma$ (resp.~$\Gamma_+$) be the subgroup (resp.~submonoid) of $\cB$ generated by the $$z_{i,a}=\dfrac{y_{i,a}y_{i,a+2}}{\textstyle\prod_{j\sim i} y_{j,a+1}}$$
for $(i,a)\in I\times_2\Z$. 

The group %
$ \cB $ and its quotient $ \cB/\Gamma$ have interesting geometric and representation theoretic incarnations in our work.  Geometrically, the map $ y_{i,a} \mapsto \O_{i,a} $ gives an isomorphism between $ \cB$ and the Picard group $ \Pic(\cZ_\infty)$ of the bi-infinite Bott--Samelson variety.  This descends~to an isomorphism
\begin{equation*}
\cB/\Gamma \cong \cP=\Pic(\cZ_{\infty}^{\circ})
\end{equation*}
Algebraically, we identify $\cB$ with the group $\fr$ of highest $\ell$-weights~of simple objects of $\mathcal{O}_{sh}^{\Z}$.
Using \textit{Nakajima's partial order} $\preceq$ on $\mathfrak{r}$, we prove the result below, which, after passing to Grothendieck groups, recovers \eqref{eq:intro_def_of_cox}.
\begin{intro_thm}[Theorem \ref{thm:BlockDec}]
The category $\O_{sh}^\Z$ decomposes as a direct sum 
\begin{equation*}
\O_{sh}^\Z=\textstyle\bigoplus_{\tau\in \cP} {}_\tau\O_{sh}^\Z,
\end{equation*}
where ${}_\tau\O_{sh}^\Z\subseteq \mathcal{O}_{sh}$ is the Serre subcategory generated by the simple objects $L(\psi)$ with highest $\ell$-weights in the coset $\tau\in \cP\simeq \cB/\Gamma$.
\end{intro_thm}
\begin{intro_rem}
This decomposition generalizes one found for finite-dimensional representations of quantum affine algebras by Chari--Moura \cite{chari2005characters} using work of Etingof--Moura \cite{etingof2003elliptic}.
\end{intro_rem}
Studying separately the geometric and categorical sides shows that, for each~$\tau$, both $\op{H}^0(\cZ_\infty^\circ, \cL_{\tau})$ and ${}_\tau\O_{sh}^\Z$ admit a filtration indexed by sets of parameters $\bR$. We identify these two filtrations in:
\begin{intro_thm}[Theorem \ref{th:TFAE}, Theorem \ref{thm:TFAEsurjection} and Proposition \ref{th:cRgamma}]
For sets of parameters $\bR_1=(R_{1,i})_{i\in I}\in \Z^{\la_1}$ and $\bR_2=(R_{2,i})_{i\in I}\in \Z^{\la_2}$, the following are equivalent:
\begin{enumerate}
\item There is an inclusion of sets $\cB(\lambda_1,\bR_1)\subset\cB(\lambda_2,\bR_2)$.
\item The monomials $y_{\bR_1}=\prod_{i\in I}\prod_{a\in R_{1,i}}y_{i,a}$ and $y_{\bR_2}=\prod_{j\in I}\prod_{b\in R_{2,j}}y_{j,b}$ are related by
\begin{equation*}
y_{\bR_2}\in y_{\bR_1}\Gamma_+.
\end{equation*}
\item Upon restriction to $\cZ_\infty^\circ$ as in \eqref{eq:intro_injection_given_by_restriction}, there is an inclusion of $G^\vee$-representations
\begin{equation*}
\op{H}^0(\cZ_\infty, \cL_{\la_1,\bR_1})\subset \op{H}^0(\cZ_\infty, \cL_{\la_2,\bR_2}).
\end{equation*}
\item There is a commutative diagram
\begin{equation*}
\adjustbox{scale = 0.9}{
\begin{tikzcd}[scale = 0.35,column sep=0.5em, row sep=1.5em]
& Y_{\mu} \ar[dl,swap,two heads] \ar[dr,two heads] &\\
Y^{\la_2}_{\mu}(\bR_2) \ar[rr, two heads, dashed] &  & Y_{\mu}^{\lambda_1}(\bR_1)
\end{tikzcd}}
\end{equation*}
where the downward arrows are the defining projections. Upon pullback, this gives~an inclusion of categories $\O^{\lambda_1}_{sh}(\bR_1)\subset\O^{\lambda_2}_{sh}(\bR_2)$ as subcategories of $\O_{sh}^\Z$.
\end{enumerate}
\end{intro_thm}

\subsection{Hernandez--Leclerc's duality} 
Another consequence of our results concerns a morphism of Grothendieck rings defined in \cite{hernandez2016cluster} in the context of monoidal categorifications of cluster algebras arising from the representation theory of Borel subalgebras of quantum affine algebras. In Section \ref{sec:duality}, we define an involutive automorphism of the scheme $\widehat{\mathcal{Z}}_{\infty}^{\circ}$ which gives rise to an algebra involution of $K_{\C}(\mathcal{O}_{sh})$ using the map~$\Omega$~of Theorem \ref{intro_Omega_iso}. We then show the following theorem using results on dual canonical bases of tensor products \cite[Chapter 24]{lusztig1993introduction}:

\begin{intro_thm} 
The morphism $D$ sends classes of simple objects to classes of simple objects. In particular, it sends the class of the positive prefundamental module $L(\sfPsi_{i,a})$ to the class of the negative prefundamental module $L(\sfPsi_{i,-a}^{-1})$ for all $(i,a)\in I\times_2\Z$. Moreover, the map $D$ restricts~to Hernandez--Leclerc's morphism (after using isomorphisms of \cite{hernandez2023representations,varagnolo2025representations}).
\end{intro_thm}

\subsection{Further research directions}
Throughout this paper, we make progress towards Geiss--Hernandez--Leclerc's conjecture that $\mathcal{O}_{sh}^{\Z}$ is a monoidal categorification of the cluster algebra $\AGHL$. Indeed, in Theorem \ref{thm:chamber_modules_are_prime}, we prove that all chamber modules, that is all the modules corresponding to cluster variables in the initial seeds described in \cite{geiss2024representations},~are real (i.e.~of simple tensor square) and prime (i.e.~they cannot be factorized non-trivially as tensor products in $\mathcal{O}_{sh}^{\Z}$). Furthermore, we highlight in Section \ref{subsec:KLRW_and_TSY} that the graded structure on KLRW algebras should give rise to a ``natural'' $t$-deformation of $K_0(\O_{sh}^\Z)$. We argue~here that this deformation, along with the $R$-matrices we construct in Section \ref{sec:TensO} (generalizing results of \cite{hernandez2024shifted}), should help prove Geiss--Hernandez--Leclerc's conjecture using a method analogous to that used in \cite{kang2018monoidal,cautis2019cluster}. We also think that the deformation of $K_0(\mathcal{O}_{sh}^{\Z})$ coming from KLRW algebras should be compatible with  the limit of the standard Poisson structure on Bott--Samelson varieties studied in \cite{elek2021bott}, as well as with Paganelli's quantum cluster algebra $\AGHL_{t,w_0}$ \cite{paganelli2025quantum}, and Fujita--Hernandez's monoidal Jantzen filtrations \cite{fujita2024monoidal}.

By carefully tracing through the different definitions, it is not hard to see that both $\cZ_{\infty}$ and the Cox ring $ \cR $  are defined over any field $k$ of characteristic $0$. Our proofs actually show that $\Omega$ is an isomorphism of $k$-algebras from $ \cR$ to $k\otimes_\Z K_0(\O_{sh})$. It would be interesting to determine if this works over fields of positive characteristic as well.

\subsection{Acknowledgements}

We are indebted to H.~Murata who carefully explained to~us the shuffle product for KLRW algebras.
We thank B.~Webster for many interesting discussions which had an influence on content appearing in Section \ref{sec:KLR}.
We also thank P.~Baumann for carefully reading Section \ref{subsec:effect_of_inv_on_basis} and for sharing great ideas concerning dual canonical bases of tensor products.
Our ongoing collaboration with A.~Kalmykov \cite{otherpaper} had a great impact on this project and we are grateful for all the ideas he shared with us.
We thank H.~Zhang for sharing his progress on Conjecture \ref{conj:Associators}, as well as
L.~Francone, D.~Hernandez, V.~Krylov, B.~Leclerc, S.~Gautam and L.~Rybnikov for interesting discussions. The authors used \textit{ChatGPT-Plus-Sol-5.6} to correct typos. Much of the research leading to this paper was conducted at \textsc{lacim} and we are grateful for their hospitality.

During this project, 
the first author held grants from the \textit{Fonds de recherche du Québec --- Nature et Technologies} (\textsc{frqnt}) and the
\textit{Natural Sciences and Engineering Research Council of Canada} (\textsc{nserc}).
The second author was supported by the master's graduate scholarships of \textsc{nserc} and \textsc{frqnt}.
The third author was supported by \textsc{frqnt}'s doctoral scholarship and the \textit{Institut des Sciences Mathématiques}.
The fourth author was supported by postdoctoral and doctoral fellowships from \textsc{frqnt} and \textsc{nserc}.
The fifth author held a grant from \textsc{nserc}.
This support is gratefully acknowledged.

\subsection{Outline of the paper}

{\it Section \ref{sec:Crystal}} introduces Nakajima's monomial crystal and the product monomial crystal. We establish combinatorial criteria that will be used throughout the paper in Theorem \ref{th:TFAE} and we characterize maximally singular sets of parameters in Theorem \ref{thm:charaterization_max_sing_crystals}. A brief discussion concerning the non-integral case appears in Section \ref{sec:nonInt}.
{\it Section \ref{sec:TSY}} defines shifted Yangians, shifted coproducts, truncated shifted Yangians together with all the necessary technical tools related to those algebras. {\it Section \ref{sec:rep_of_shifted_yangians}} studies modules over shifted Yangians. There, we recall the definition of the category $\O_{sh}$ with the truncated categories $\O_\mu^\lambda(\bR)$. We give a direct sum decomposition for $\O_{sh}$ in Theorem \ref{thm:BlockDec} and define chamber modules in Section \ref{sec:Chambermodules}. We also introduce GT-characters in Section \ref{subsec:GT_weights} which are one of the main tools used to study the $G^\vee$-action on $K_\C(\O_{sh})$. {\it Section \ref{sec:TensorResults}} starts by giving combinatorial criteria for identifying highest $\ell$-weight objects which, when combined with Theorem \ref{thm:trTensor}, produce families of $R$-matrices. We then prove generic simplicity results for tensor products of simple objects in Theorem \ref{thm:TensSimpInt} and show in Section \ref{sec:TSC} that coproducts descend to truncations. \medskip\par

{\it The \nameref{sec:interlude}} marks an important separation in the paper, where we switch the roles of $G$ and $G^\vee$ and restrict ourselves to the integral part of $\mathcal{O}_{sh}$. \medskip\par

{\it Section \ref{sec:KLR}} recalls KLR algebras, KLRW algebras and parity KLRW algebras. A yoga of shuffle products for parity KLRW algebra modules is developed in Section \ref{subsec:char_parity_KLRW} and Theorem \ref{thm:the_big_pentagone_commutes} compares it with the product of GT-characters.
{\it Section \ref{sec:glueing_of_g_action}} explains how to glue the $G$-action on each $K_\C(\O_{sh}^\lambda(\bR))$ into a $G$-action on $K_\C(\O_{sh})$.
{\it Section \ref{sec:description_of_mult}} applies previous results to the proof of certain relations in $K_0(\O_{sh})$. Theorem \ref{thm:subalgebra_iso_base_affine_space} studies relations coming from the base affine space while Theorem \ref{thm:extended_QQ_holds} shows that isoclasses of chamber modules satisfy the extended $QQ$-system. A discussion about real and prime objects can be found in Section \ref{sec:application_real_and_prime}.
{\it Section \ref{sec:bott-samselson-background}} switches gears and tackles the geometry of Bott--Samelson varieties. {\it Section \ref{sec:bi-infinite-bott-samelson}} introduces the bi-infinite Bott--Samelson variety $\mathcal{Z}_\infty$ and the Cox ring $\cR$ of the open cell $\mathcal{Z}_\infty^\circ\subset \mathcal{Z}_\infty$, for which a presentation is given in
Corollary \ref{cor:presentation-of-R}.
Finally, Section \ref{subsec:filt_of_R} studies spaces of sections of the line bundles $\cL_{\la,\bR}$ and proves they have~the~same~characters as the product monomial crystals $\cB(\lambda,\bR)$.
{\it Section \ref{sec:iso_between_KOsh_and_R}} relates the geometry to the category $\mathcal{O}_{sh}$. This culminates in a proof that $K_\C(\O_{sh})$ is isomorphic to $\cR$, which can be found in Theorems \ref{thm:Omega_is_surj}, \ref{thm:omega_is_inj_cluster_side} and \ref{thm:injectivity-geo}. After this, Section \ref{subsec:subalgebras_of_KOsh} identifies notable full subcategories $\O_{sh}(\xi,\xi')\subset\O_{sh}$ categorifying the coordinate rings $\C[\widehat{\cZ}_{\xi,\xi'}^\circ]$. Special cases recover $\C[N_{-}\backslash G]$ and $\C[G]$. Corollary \ref{cor:alternate_pres_type_A} gives an alternative presentation of $K_0(\O_{sh})$ in type A. {\it Section \ref{sec:duality}} studies an involution of the scheme $\widehat{\cZ}_\infty^\circ$ (cf.~Theorem \ref{thm:F_is_involution_of_bi_infinite_bs}) and transports it to $K_\C(\O_{sh})$ via the isomorphism $\Omega$. We show in Theorem  \ref{thm:DLpsi_is_still_the_class_of_a_simple} that the resulting involution of $K_{\C}(\mathcal{O}_{sh})$ preserves classes of simple objects and prove in Section \ref{sec:HLD} that this involution restricts~to a duality studied previously by Hernandez--Leclerc.
\medskip\par
\
The article concludes with two appendices. {\it Appendix \ref{app:pro_varieties}} gives the necessary tool-kit for dealing with the pro-varieties of Section \ref{sec:bi-infinite-bott-samelson}. {\it Appendix \ref{sec:Inflations}} studies~how~truncated~categories
behave under inclusions of Dynkin diagrams. This leads to a notion of inflation for modules over shifted Yangians and proves a conjecture made by the fourth author in \cite{pinetinflations}. Results of this appendix are also used in Section \ref{subsec:extendQQ_and_principal_block} in order to give a conjectural categorification of the extended $QQ$-system via an exact sequence happening in the principal block for $\fsl_2$.

\subsection{Notation}\label{sec:Notation}
The symmetric group on $m$ letters is denoted by $\Sigma_m$. The collection of multisets of complex numbers of size $m$ is denoted by $\C^m/\Sigma_m$. For a multiset $A\in \C^m/\Sigma_m$, we let $ p_A(u) = \prod_{a \in A} (u - a) \in \C[u]$.  This realizes a bijection between multisets and monic polynomials.  Note that $ p_{A \cup B}(u) = p_A(u) p_B(u) $, where $ \cup $ denotes union (with multiplicities) of multisets.\medskip\par

Let $ X $ be a set. We will write $\Z^X$ for the collection of functions from $X$ to $\Z$. We will represent such functions as formal linear combinations of the form $ \sum_{x \in X} a_x [x] $ for $ a_x \in \Z $. Here, $[x]\in \Z^X$ denotes the delta function at $x\in X$, i.e.~$[x](x)=1$ and $[x](y)=0$ for all $y\neq x$.
This function will also be denoted by $\delta_x$.\medskip\par

Throughout, $I$ will denote a simply-laced Dynkin diagram. For $i,j\in I$, write $i\sim j$ if $i$ is connected to $j$. Fix a decomposition $I=I_{0}\sqcup I_{1}$ into even and odd vertices which makes $ I$ into a bipartite graph (i.e. no edge has vertices of the same parity). We will denote by $\overline{i}\in \Z/2\Z$ the parity of $i\in I$. We will often use the two sets
\begin{align*}
I \times_2 \Z := \{ (i,a) \in I \times \Z : a-i\in 2\Z \},\hspace{1em} I \times_2^{\mathrm{op}} \Z := \{ (i,b) \in I \times \Z : b-i\in 2\Z+1\}.
\end{align*}
Accordingly, we have $\overline{i}+2\Z=\{a\in \Z\,|\,(i,a)\in I\times_2\Z\}$. We endow $I \times_2 \Z$ with the partial order $\le$ generated by the relations $(i,a)\le (i,a+2)$ and $(i,a)\le (j,a+1)$ whenever $i\sim j$. Equivalently, we have $(i,a) \le (j,b)$ if and only if $a+d(i,j)\le b$, where $d(i,j)$ denotes the distance in $I$ between the vertices $i$ and $j$. We also endow $I \times_2^{\mathrm{op}} \Z$ with the partial order defined in the same way. In Section \ref{subsec:parity_KLRW_alg}, we will be using a total order on $I \times_2 \Z$ which should not be confused with the previously defined partial order. It will depend on the extra data of a total order on $I$ such that even vertices are smaller than odd vertices.

\subsubsection{The Grothendieck group}
Let $\mathcal{C}$ be an abelian category. We denote its Grothendieck group by $K_0(\mathcal{C})$. When $\mathcal{C}$ is of finite-length, the abelian group $K_0(\mathcal{C})$ has a basis given by isoclasses of simple objects. If $\mathcal{C}$ also has a tensor product $\otimes:\mathcal{C}\times \mathcal{C}\to \mathcal{C}$, the Grothendieck group $K_0(\mathcal{C})$ becomes a ring by defining $[M_1]\cdot [M_2]:=[M_1\otimes M_2]$%
. The complexification of $K_0(\mathcal{C})$ is denoted by $K_\C(\mathcal{C}):=\C\otimes_\Z K_0(\mathcal{C})$.

\subsubsection{Lie theoretic notation}\label{subsubsec:lie_theory_notation}
Let $\fg$ denote the finite-dimensional Lie algebra over $\C$ associated to the simply-laced Dynkin diagram $I$. Given a complex algebraic group $G$ whose Lie algebra is $\fg$, we fix $B$ a choice of Borel subgroup of $G$ together with $T\subset B$ a maximal torus. We denote by $N$ the unipotent radical of $B$. We let $B_-$ be the opposite Borel subgroup with respect to $B$ and $N_-$ be its unipotent radical. The complex Lie algebras of the groups $T$, $N$, $N_-$, $B$, $B_-$ are denoted by $\fh$, $\fn$, $\fn_-$, $\fb$, $\fb_-$ respectively.
\medskip\par
Let $\Delta$ be the set of roots of $\fg$. The choice of Borel subgroup $B$ determines a subset of positive roots $\Delta_+\subset\Delta$ along with a set of simple roots $\{\alpha_i\}_{i\in I}$. The root lattice is denoted by $Q:=\bigoplus_{i\in I} \Z\alpha_i$ and the positive root cone is defined as $Q_+:=\bigoplus_{i\in I} \Z_{\geq 0}\alpha_i\subset Q$. The simple coroots for the dual root system are denoted by $\{\alpha_i^\vee\}_{i\in I}$. We fix fundamental weights $\{\varpi_i\}_{i\in I}\subset \fh^\ast$ defined by $\langle \varpi_i,\alpha_j^\vee\rangle=\delta_{ij}$. Similarly, we have fundamental coweights $\{\varpi_i^\vee\}_{i\in I}\subset \fh$. The Cartan matrix of $\fg$ is denoted by $C=(c_{ij})_{i,j\in I}$, where $c_{ij}=\langle\alpha_j,\alpha_i^\vee\rangle$. The weight lattice of $\fg$ is denoted by $P=\bigoplus_{i\in I} \Z\varpi_i$ and the cone of dominant weights is labelled $P_+:=\bigoplus_{i\in I}\Z_{\geq 0}\varpi_i$.

For $\lambda\in P_+$, let $V(\lambda)$ be the irreducible $\fg$-representation of highest weight $\lambda$. We~fix~once and for all a choice of highest weight vector $v_\lambda\in V(\lambda)$ for each such representation. We also endow each $V(\lambda)$ with its dual canonical basis $B(\lambda)$, where $v_\lambda\in B(\lambda)$.

The Weyl group of $\fg$ will be denoted by $W$ and its canonical generators, the simple reflections, will be denoted by $\{s_i\}_{i\in I}$. The length function will be denoted by $\ell(\trou)$. The longest element of $W$ will be denoted $w_0$.
We let $(\trou)^\ast:I\to I$ be the involution satisfying $w_0\alpha_i=-\alpha_{i^\ast}$.
A Coxeter element of $W$ is a product of all the simple reflections, each appearing exactly once in the product. The order of such an element is independent of the choice in which the simple reflections appear and it is called the Coxeter number of $\fg$. We denote this non-negative integer by $h$.

We also fix a Serre presentation for $\fg$, where the Chevalley generators are labelled for each $i\in I$ by $e_i\in \fg_{\alpha_i}$, $f_i\in \fg_{-\alpha_i}$ and $h_i:=[e_i,f_i]$. 

We label data concerning the Langlands dual group using the $(\trou)^\vee$ notation. 

\subsubsection{Height functions}\label{subsubsec:height_functions}

We will make extensive use of \textit{height functions}, which are sequences $ \xi = (\xi_i)_{i \in I } \in \prod_{i\in I}(\overline{i}+2\Z)$ such that
\begin{equation*}
|\xi_i - \xi_j | = 1 \text{ whenever } i \sim j.
\end{equation*}
A choice of height function $\xi$ determines an orientation for the Dynkin diagram $I$ %
(where we put a directed arrow $i\rightarrow j$ whenever $i\sim j$ and $\xi_i-\xi_j=1$). Such a choice also specifies a Coxeter element $c_{\xi}\in W$ 
via
\begin{equation*}
c_\xi := \prod^{\rightarrow}_{k\in \Z}\prod_{\substack{i\in I\\ \xi_i=k}} s_i.
\end{equation*}
The height function $\xi$ also specifies a dual height function $\xi^\ast$ defined by $(\xi^\ast)_i:=\xi_{i^\ast}$ for all $i\in I$. It is a well-known combinatorial fact that the integers $m_i=\tfrac{1}{2}(\xi_i^*-\xi_i+h)\in \Z$ are the smallest non-negative integers such that 
\begin{equation*}
(c_\xi)^{m_i}\varpi_i=-\varpi_{i^\ast}.
\end{equation*}
The collection of all height functions forms a partially ordered set by declaring $\xi\leq \xi'$ if $\xi_i\leq \xi_{i}'$ for all $i\in I$.

\section{Product monomial crystals}\label{sec:Crystal}
This section recalls facts about the \textit{product monomial crystals} of \cite{kamnitzer2019highest} and introduces a second notion of weight for %
these crystals. It also includes a characterization of all \textit{maximally singular} product monomial crystals (%
which answers a question of \cite{kamnitzer2019highest}).
\subsection{Crystals}
Denote by $\cB$ the \textit{monomial crystal} of \cite{nakajima2003tanalogues}, that is the set of (all) Laurent monomials in the infinitely-many variables $$\{y_{i,a}\,|\,(i,a)\in I\times_2\mathbb{Z}\}.$$
Clearly, $\cB$ is an abelian group under multiplication%
. It is moreover a normal $\fg^\vee$-crystal. To describe the corresponding crystal structure, let
\begin{equation}\label{eq:zia}
z_{i,a}=\frac{y_{i,a}y_{i,a+2}}{\prod_{j\sim i} y_{j,a+1}}
\end{equation}
for $(i,a) \in I\times_2\mathbb{Z}$. Let also $\wt: \cB \rightarrow P^\vee$ be the group morphism, called \textit{weight map}, given by $\wt(y_{i,a})=\varpi_i^{\vee}$. %
Consider the maps $\{\varepsilon_i,\varphi_i:\cB\rightarrow \Z\}_{i\in I}$ given on a monomial~$m\in \cB$ by
$$  \varepsilon_i(m)=\max_{\substack{a\in \overline{i}+2\Z}} \varepsilon_i^{a}(m)\ \,\text{and}\ \, \varphi_i(m)=\max_{a\in \overline{i}+2\Z} \varphi^{a}_i(m) $$
where, for $a\in \overline{i}+2\Z$ and $m=\prod_{(i,a)\in I\times_2\Z}y_{i,a}^{n_{i,a}}$, %
$$ \textstyle \varepsilon_i^a(m) = -\sum_{b\leq a}n_{i,b}\ \, \text{and}\ \, \varphi_i^a(m)=\sum_{b\geq a}n_{i,b}.$$
Finally, denote by $\{\tilde{e}_i,\tilde{f}_i:\cB\rightarrow\cB\cup\{0\}\}_{i\in I}$ the \textit{crystal operators} given by
\begin{align*}
\tilde{e}_i(m)&=\begin{cases}
0 & {\text{if }\varepsilon_i(m)=0,}\\
z_{i,a} m & {\text{otherwise, where $a\in \overline{i}+2\Z$ is minimal such that $\varepsilon_i(m)=\varepsilon_i^a(m)$;}}
\end{cases}
\end{align*}
and
\begin{align*}
\tilde{f}_i(m)&=\begin{cases}
0 & {\text{if }\varphi_i(m)=0,}\\
z_{i,a-2}^{-1} m & {\text{otherwise, where $a\in \overline{i}+2\Z$ is maximal such that  $\varphi_i(m)=\varphi_i^a(m)$.}}
\end{cases}
\end{align*}
Kashiwara's result can be written more explicitly as:
\begin{Theorem}[{\cite[Proposition~3.1]{kashiwara2003realizations}}]\label{thm:MonCrys}
The set $\cB$ together with the maps $\wt,\varepsilon_i,\varphi_i,\tilde{e}_i,\tilde{f}_i$ above defines a normal $\fg^\vee$-crystal.
\end{Theorem}
For $\lambda\in P^{\vee}_+$, %
let $B(\lambda)$ be the crystal of the finite-dimensional simple $\fg^{\vee}$-module of highest weight $\lambda$. By \cite[Theorem 4.3]{kashiwara2003realizations}, there are infinitely-many distinct copies of $B(\lambda)$ in %
$\cB$. In particular, for $(i,a)\in I\times_2\Z$, %
the subcrystal $\cB(\varpi_i^{\vee},a)$ of $\cB$ generated by the variable~$y_{i,a}$ satisfies%
, as abstract crystals,
$$ \cB(\varpi_i^{\vee},a)\simeq B(\varpi_i^{\vee}).$$
Unfortunately, describing all %
elements in $ \cB(\varpi_i^{\vee}, a)$ is a tedious task, especially when~$i\in I$ is not minuscule. 
Nevertheless, we can combine \cite[Proposition 3.2]{pressley1991fundamental} with the upcoming Theorem \ref{thm:NakCrys} to get:
\begin{Lemma}\label{th:lowest}
Let $ (i,a) \in I \times_2 \Z $. The lowest weight element of $\cB(\varpi_i^{\vee}, a)$ is $y_{i^*, a-h}^{-1}$.
\end{Lemma}
The rough idea of product monomial crystals \cite{kamnitzer2019highest} is to multiply together different \textit{fundamental crystals} $\cB(\varpi_i^{\vee},a)$ to generate bigger parts of the monomial crystal $\cB$.

Let $\lambda\in P^{\vee}_+$ be a dominant coweight and write $\lambda = \sum_{i\in I} \lambda_i \varpi_i^\vee$
with $\{\lambda_i\}_{i\in I}\subseteq \mathbb{Z}_{\geq 0}$. Set
$$\textstyle \C^\lambda = \prod_{i\in I} \C^{\lambda_i} / \Sigma_{\lambda_i},$$
where $\Sigma_k$ is the symmetric group on a set of size $k$. %
We will call a point $\bR = (R_i)_{i\in I}$ in $\C^{\lambda}$ a \emph{set of parameters of coweight} $\lambda$. Note that such a point is simply a collection of multisets $R_i$ of size $\lambda_i$ for each $i$. For example,
\begin{equation}\label{eq:bR}
(\{-1,1,1\},\{0,0\}) \in \C^{3\varpi_1^{\vee}+2\varpi_2^{\vee}}.
\end{equation}
An element $\bR\in \mathbb{C}^{\lambda}$ is %
 \textit{integral} if $R_i\subseteq \overline{i}+2\mathbb{Z}$ for each $i\in I$. Note that this holds for \eqref{eq:bR} (if we choose $\overline{1}=1$ and $\overline{2}=0$).
 
\begin{Def}[\cite{kamnitzer2019highest}] 
Fix $\lambda\in P_+^{\vee}$ and let $\bR=(R_i)_{i\in I}\in \mathbb{C}^{\lambda}$ be integral. Then~the \textit{product monomial crystal} is the set
\begin{equation}\label{eq:MonCrys}
\cB(\lambda,\bR)=\prod_{i\in I}\prod_{a\in R_i} \cB(\varpi_i^{\vee},a)
\end{equation}
where the products denote multiplication of monomials in $\cB$. Thus $\cB(\lambda,\bR)$ is the subset of $\cB$ formed by multiplying monomials coming from the \textit{fundamental crystals} $\cB(\varpi_i^{\vee},a)$
(and then forgetting multiplicities of the resulting elements of $\cB$).
\end{Def}

As the name suggests, product monomial crystals $\cB(\lambda,\bR)$ are %
(normal) subcrystals of~$\mathcal{B}$.%

\begin{Theorem}[{\cite[Theorem 2.2]{kamnitzer2019highest}}]\label{thm:ProdMonCrys}
For $\lambda\in P_{+}^{\vee}$ and integral $\bR\in \mathbb{C}^{\lambda}$, the product monomial crystal $\cB(\lambda,\bR)$ is a normal subcrystal of $\cB$. Also, there exist crystal inclusions
\begin{equation*}
\textstyle B(\lambda)\subseteq \cB(\lambda,\bR)\subseteq \bigotimes_{i\in I} B(\varpi_i^{\vee})^{\otimes \lambda_i}.
\end{equation*}
\end{Theorem}

Now, for $\mu\in P^{\vee}$, let
$$\cB(\la,\bR)_{\mu} = \{p\in \cB(\la,\bR)\,|\,\wt(p)=\mu\}\ \text{ and }\ B(\lambda)_{\mu}=\{x\in B(\la)\,|\,\wt(x)=\mu\}$$
be the weight-$\mu$ components of the crystals $\cB(\la,\bR)$ and $B(\la)$ (respectively). We will need:
\begin{Lemma}\label{lem:ConvexHullWts} Fix $\la=\sum_{i\in I}\la_i\varpi_i^{\vee}\in P_+^{\vee}$ and suppose that $\bR\in \C^{\la}$ is integral. Then
$$\{\mu\in P^{\vee}\,|\,\cB(\la,\bR)_{\mu}\neq \emptyset\}=\{\mu\in P^{\vee}\,|\,B(\la)_{\mu}\neq \emptyset\},$$
that is the crystals $\cB(\la,\bR)_{\mu}$ and $B(\la)_{\mu}$ have the same weights.
\end{Lemma}

\begin{proof}
The inclusion $\supseteq$ follows from Theorem \ref{thm:ProdMonCrys}. For the other one, recall that
the set of weights of the irreducible representation $V(\la)$ is equal to $(\la-Q^{\vee}_+)\cap \text{Conv}(W\la)$ by \cite[Chapitre 8\S 7]{bourbaki1975lie} (where $\text{Conv}(W\la)$ is the convex hull of $W\la$). In particular, the set of weights of $V=\bigotimes_{i\in I} V(\varpi_i^{\vee})^{\otimes \la_i}$ is the intersection of $\la-Q_+^{\vee}$ with
$$ \textstyle \sum_{i\in I}\lambda_i\text{Conv}(W\varpi_i^{\vee})=\text{Conv}(W\la),$$
where the last equality is well-known (and follows for example from \cite[Lemma 6.1]{kamnitzer2005mirkovic}). Hence, weights of $V$ are all weights of $V(\la)$ and Theorem \ref{thm:ProdMonCrys} implies the result.
\end{proof}

\begin{Rem}\label{rem:nonIntCryst}%
Product monomial crystals $\cB(\lambda,\bR)$
can also be constructed for non-integral sets of parameters $\bR$ \cite{kamnitzer2019highest}. More precisely, for $r\in \mathbb{C}$ with $i\in I$, fix $a\in \overline{i}+2\Z$ and let $\cB(\varpi_i^{\vee},r)$ be the crystal which has the same crystal operators as $\cB(\varpi_i^{\vee},a)$, but for which the monomials are translated according to $y_{j,b}\mapsto y_{j,b+r-a}$ (for every $(j,b)\in I\times_2\Z$). One can define a set $\cB(\lambda,\bR)$ using \eqref{eq:MonCrys} for every (i.e.~even non-integral) set of parameters $\bR \in \C^{\lambda}$, but the resulting set is no longer a subset of the monomial crystal $\cB$. Moreover, in this more general situation, the set $\cB(\lambda,\bR)$ is no longer a $\fg^{\vee}$-crystal, but rather a crystal for the Lie algebra $(\fg^{\vee})^{\oplus k}$ where $k$ is the number of distinct integrality classes in $\bR$ (cf.~\cite{kamnitzer2019highest}).
\end{Rem}

\begin{Def} A set of parameters $\bR\in \C^{\lambda}$ is \textit{maximally singular} if $\cB(\lambda,\bR)\simeq B(\lambda)$.
\end{Def}

\begin{Example}\label{ex:sl3_prod_mon_crystal} 
For $\fg=\mathfrak{sl}_3$ and $a\in \C$, the fundamental crystal $\cB(\varpi_1^{\vee},a)$ is given by
\begin{center}
\begin{tikzpicture}[yscale=1.35]
\def\radR{1.25}
\node (A) at (60:\radR) {$y_{1,a}$};
\node (B) at (180:\radR) {$\frac{y_{2,a-1}}{y_{1,a-2}}$};
\node (C) at (-60:\radR) {$\frac{1}{y_{1,a-2}y_{2,a-3}}$};
\draw[-stealth,very thick,orange] (A) -- (B);
\draw[-stealth,very thick,blue] (B) -- (C);
\end{tikzpicture}
\end{center}
where {\color{orange}\bf orange} (resp.~{\color{blue}\bf blue}) arrows correspond to the action of 
$\tilde{f}_1$ (resp.~$\tilde{f}_2$). Now, consider the set of parameters 
$$\bR_k=(\{a\},\{a+2k+1\})\in \C^{\lambda}$$ 
for $\lambda=\varpi_1^{\vee}+\varpi_2^{\vee}$ and $k\in \Z$. If $k=0$, the product monomial crystal $\cB(\lambda,\bR_0)$ is
\begin{center}
\begin{tikzpicture}[yscale=0.85]
\def\radR{3}
\def\xshift{1}
\def\yshift{0}
\node (A) at (90:\radR) {$y_{1,a}y_{2,a+1}$};
\node (B) at (30:\radR) {$\frac{y_{1,a}^2}{y_{2,a-1}}$};
\node (C) at (150:\radR) {$\frac{y_{2,a-1}y_{2,a+1}}{y_{1,a-2}}$};

\node (D) at (-\xshift,-\yshift) {$\frac{y_{1,a}}{y_{1,a-2}}$};
\node (E) at (\xshift,\yshift) {$\frac{y_{2,a+1}}{y_{2,a-3}}$};

\node (F) at (-30:\radR) {$\frac{y_{1,a}}{y_{2,a-1}y_{2,a-3}}$};
\node (G) at (210:\radR) {$\frac{y_{2,a-1}}{y_{1,a-2}^2}$};
\node (H) at (270:\radR) {$\frac{1}{y_{1,a-2}y_{2,a-3}}$};

\draw[-stealth,very thick,orange] (A) -- (C);
\draw[-stealth,very thick,blue] (A) -- (B);

\draw[-stealth,very thick,blue] (C) [out=-30, in=135] to (E);
\draw[-stealth,very thick,orange] (B) [out=210, in=45] to (D);

\draw[-stealth,very thick,blue] (E) -- (F);
\draw[-stealth,very thick,orange] (D) -- (G);

\draw[-stealth,very thick,blue] (G) -- (H);
\draw[-stealth,very thick,orange] (F) -- (H);
\end{tikzpicture}
\end{center}
where we use the same {\color{orange}\bf orange} and {\color{blue}\bf blue} convention as above. Clearly, the monomial
$$ \textstyle \frac{y_{1,a}}{y_{1,a-2}}=y_{1,a} \frac{1}{y_{1,a-2}} =  \frac{y_{2,a-1}}{y_{1,a-2}} \frac{y_{1,a}}{y_{2,a-1}}\in \cB(\lambda,\bR_0) $$
can be expressed in two distinct ways as a product of a monomial in $\cB(\varpi_1^{\vee},a)$ with one in $\cB(\varpi_2^{\vee},a+1)$. This explains why $\bR_0$ is maximally singular, that is why, as $\fg^{\vee}$-crystals,
$$\cB(\varpi_1^{\vee}+\varpi_2^{\vee},\bR_0)\simeq B(\varpi_1^{\vee}+\varpi_2^{\vee}).$$
On the other hand, for $k\neq 0$, the crystal $\cB(\lambda,\bR_k)$ is easily seen to be isomorphic to $$ B(\varpi_1^{\vee}+\varpi_2^{\vee})\oplus B(0) \simeq B(\varpi_1^{\vee})\otimes B(\varpi_2^{\vee})$$
and thus contains 9 (distinct) elements.
\end{Example}

\subsection{A second notion of weight}\label{sec:awt}
Write $ \Gamma \subseteq \cB $ for the subgroup of all Laurent monomials in the elements $\{z_{i,a}\}_{(i,a)\in I\times_2\Z}$ of \eqref{eq:zia} %
and let $$\awt : \cB \twoheadrightarrow \cB/\Gamma$$ be the natural projection. Then $\cB/\Gamma$ is (non-canonically) isomorphic to the coweight lattice $P^{\vee}$ and is hence a free abelian group. Fix a height function $\xi$ and let $c=c_\xi\in W$ be the corresponding Coxeter element (see Section \ref{subsubsec:height_functions}). 
The following easily proven result uses the element $c=c_{\xi}$ to give an isomorphism $\cB/\Gamma\simeq P^{\vee}$ and thus makes it reasonable to call the map $\awt:\cB\twoheadrightarrow \cB/\Gamma$ a ``second weight map''.

\begin{Lemma}\label{lem:height_function_identifies_groups}
Fix a height function $ \xi $ with associated Coxeter element $ c=c_{\xi}$.
\begin{enumerate}
\item The group $ \cB/\Gamma$ is a free abelian group with generators $ \{\awt(y_{i, \xi_i}) : i \in I \} $.
\item The map
$ \cB/\Gamma \rightarrow P^\vee $
given by $ \awt(y_{i,\xi_i + 2s}) \mapsto c^s \varpi_i^\vee $ %
gives a group isomorphism.
\end{enumerate}
\end{Lemma}
Consider the submonoids $\cB_+\subseteq \cB$ and $\Gamma_+\subseteq \Gamma$ of actual (i.e.~not Laurent) monomials in (respectively) the $y_{i,a}$'s and $z_{i,a}$'s with $(i,a)\in I\times_2\Z$. Note that $\Gamma_+ $ is not contained~in~$\cB_+$. Remark also that, if $ m \in \cB $ and $ \tilde{e}_i(m) \ne 0 $, then there exists $ a\in \overline{i}+2\Z $ with $ z_{i,a} m = \tilde e_i(m) $. This observation leads us to the results below. (Recall that $(i^*)^*=i$ for all $i\in I$.)
\begin{Lemma} \label{th:qm}
For all $ m \in \cB $, there exists $ q \in \Gamma_+ $ such that $ qm \in \cB_+$.
\end{Lemma}
\begin{proof}
By multiplicativity, we can assume $m=y_{i,a}^{-1}$ where $(i,a)\in I\times_2 \Z$. Thus Lemma~\ref{th:lowest} gives %
$m\in \cB(\varpi_{i^*},a+h)$ and the observation above implies that $y_{i^*,a+h} \in m\Gamma_+\cap \cB_+$.
\end{proof}
Given $\lambda \in P_+^{\vee}$ and $\bR=(R_i)_{i\in I} \in \C^{\lambda}$ an integral set of parameters, we define
\begin{equation}\label{eq:yR}
y_{\bR} = \prod_{i\in I}\prod_{a\in R_i}y_{i,a} \in \cB_+.
\end{equation}
The result below follows directly from \eqref{eq:MonCrys} and Lemma \ref{th:qm}.
\begin{Lemma} \label{th:qm2}
Fix $\lambda, \bR$ as above with $ p \in \cB(\la, \bR) $. Then $ qp = y_{\bR}$ for some $ q \in \Gamma_+ $.
\end{Lemma}

We will be interested in studying containments of product monomial crystals, that is in studying inclusions $\cB(\la_1,\bR_1)\subseteq\cB(\la_2,\bR_2)$ with $\lambda_1,\lambda_2\in P_+^{\vee}$, $\bR_1\in \C^{\lambda_1}$ and $\bR_2\in \C^{\lambda_2}$.
For this purpose, remark that the map $ \awt:\cB\twoheadrightarrow \cB/\Gamma $ is constant on each fundamental crystal $\cB(\varpi_i^{\vee}, a) $ and thus on each product monomial crystal $\cB(\la, \bR) $ by multiplicativity. In fact:%

\begin{Theorem}\label{th:TFAE}
Fix $\la_1,\la_2\in P_{+}^{\vee}$ with $\bR_1\in \C^{\la_1}$ and $ \bR_2\in \C^{\la_2}$ integral. Then the following statements are equivalent:
\begin{enumerate}
\item $\cB(\la_1, \bR_1) \subseteq \cB(\la_2, \bR_2)$,
\item $ \emptyset \neq \cB(\lambda_1, \bR_1)_\mu \subseteq \cB(\lambda_2, \bR_2)_\mu$ for some $\mu\in P^{\vee}$,
\item $ y_{\bR_1} \in \cB(\lambda_2, \bR_2)$, and
\item $ y_{\bR_2} \in y_{\bR_1} \Gamma_+$.
\end{enumerate}
\end{Theorem}

\begin{proof} The implication $(1) \Longrightarrow(2)$ is clear and $(3)\Longrightarrow (4)$ follows from Lemma \ref{th:qm2}. %
For $(2)\Longrightarrow (3)$, use Lemma \ref{lem:ConvexHullWts} to fix $m\in \cB(\la_1,\bR_1)_{\mu}$ in the connected subcrystal generated by $y_{\bR_1}$. Then $m\in \cB(\lambda_2,\bR_2)_{\mu}$ %
and applying $\tilde{e}_i$'s gives $y_{\bR_1}\in \cB(\la_2,\bR_2)$. Let~us~finally~show $(4)\Longrightarrow(1)$. For this, assume $y_{\bR_2}\in y_{\bR_1}\Gamma_+$ and fix $p\in \cB(\la_1,\bR_1)$. Use the~equivalence~class\footnote{Note that the equivalence classes defined in \cite{kamnitzer2019highest} for ``\textit{$\bR$-monomial data}'' do not depend on $\bR$.}
$$[p]\subseteq\cB(\la_1,\bR_1)$$
defined in \cite[Proof of Theorem 6.4]{kamnitzer2019highest}. By this same reference, there is a monomial $m\in \cB_+\cap [p]%
$. By hypothesis and since $m\in \cB(\la_1,\bR_1)$,
$$ y_{\bR_2}\in y_{\bR_1}\Gamma_+\subseteq m \Gamma_+ $$
and \cite[Lemma 6.7]{kamnitzer2019highest} gives $m\in \cB(\la_2,\bR_2)$. It then follows (by \cite[Proof of Theorem 6.4]{kamnitzer2019highest} again) that $[p]\subseteq \cB(\la_2,\bR_2)$. This ends the proof.
\end{proof}

An interesting corollary of the above result is that product monomial crystals with fixed value under the map $\awt:\cB\twoheadrightarrow\cB/\Gamma$ form a directed system under inclusion.
\begin{Corollary} \label{co:directed}
Fix $\lambda_1,\lambda_2,\bR_1,\bR_2$ as in Theorem \ref{th:TFAE} and suppose $\awt(y_{\bR_1})=\awt(y_{\bR_2})$. Then there exists $\lambda\in P_+^{\vee}$ and $\bR\in \C^{\la}$ integral with $\cB(\la_1,\bR_1)\cup\cB(\la_2,\bR_2)\subseteq \cB(\la,\bR)$.
\end{Corollary}
\begin{proof}
Since $\awt(y_{\bR_1})=\awt(y_{\bR_2})$, we can find $q_1,q_2\in \Gamma_+$ such that $y_{\bR_1}q_1=y_{\bR_2}q_2$. Also, by Lemma \ref{th:qm}, there exists $q\in \Gamma_+$ such that $y_{\bR_1}q_1q \in \cB_+$. Clearly, one can write
\begin{equation}\label{eq:directed}
y_{\bR_2}q_2q=y_{\bR_1}q_1q = \prod_{i\in I}\prod_{a\in R_i} y_{i,a} = y_{\bR}
\end{equation}
for some integral set $\bR=(R_i)_{i\in I}$ of parameters. Let $\lambda = \sum_{i\in I}\lambda_i \varpi_i^{\vee}$ with $\lambda_i$ the cardinality of $R_i$ (for $i\in I$). Then $\bR\in \C^{\la}$, and \eqref{eq:directed} with Theorem \ref{th:TFAE} imply the result.
\end{proof}

Another nice corollary is that each element of $\cB$ is contained in a unique minimal product monomial crystal (which is equivalent to saying that any intersection of product monomial crystals is a union of product monomial crystals).

\begin{Corollary} \label{cor:minimal-pmc}
For $p\in \cB$, there exist unique $\lambda_1\in P^{\vee}_+$ and $\bR_1\in \C^{\la_1}$ such that,~for every $\lambda_2\in P^\vee_+$ and $\bR_2\in \C^{\lambda_2}$,
\begin{equation*}
p\in \cB(\la_2,\bR_2) \iff \cB(\la_1,\bR_1)\subseteq \cB(\la_2,\bR_2)
\end{equation*}
\end{Corollary}

\begin{proof}
Fix, as in the proof of Theorem \ref{th:TFAE}, $m\in \cB_+\cap [p]$. Then $m=y_{\bR_1}$ for some $\la_1\in P^{\vee}_+$ and $\bR_1\in \C^{\la_1}$. Also, for $\la_2\in P^{\vee}_+$ and $\bR_2\in \C^{\la_2}$, by \cite[Proof of Theorem 6.4]{kamnitzer2019highest},
$$ p\in \cB(\la_2,\bR_2) \iff [p]\subseteq \cB(\la_2,\bR_2) \iff m=y_{\bR_1}\in \cB(\la_2,\bR_2)$$
with the latter condition equivalent to the containment $\cB(\la_1,\bR_1)\subseteq \cB(\la_2,\bR_2)$ by Theorem \ref{th:TFAE}(3). This ends the proof as the pair $(\la_1,\bR_1)$ is clearly unique.
\end{proof}

\subsection{Maximally singular crystals}\label{subsec:maximally_sing_crystals}
Fix $\la%
\in P_+^{\vee}$ with a height function $\xi=(\xi_i)_{i\in I}$. Then a set of parameters $\bR=(R_i)_{i\in I}\in\C^{\la}$ is said to be \textit{given by the height function $\xi$} if, for all $i\in I$, the multiset $R_i$ is formed by $\langle \la,\alpha_i\rangle$ copies of the integer $\xi_i$ (and nothing else), i.e.
$$R_i=\{%
\xi_i,\dots,\xi_i%
\}$$
where $|R_i|=\langle\la,\alpha_i\rangle$.
The theorem below answers a question asked in \cite[p.10]{kamnitzer2019highest}.~For $i,j\in I$, let $d(i,j)$ be the length of the shortest path joining $i$ and $j$ in the Dynkin diagram of $\fg$. Note that $d(i,i)=0$ and $d(i,j)=d(j,i)$ for all $i,j\in I$.

\begin{Theorem}\label{thm:charaterization_max_sing_crystals}
Fix $\la\in P_+^{\vee}$ and $\bR=(R_i)_{i\in I}\in \C^{\la}$ integral. Then the following statements are equivalent:
\begin{itemize}
\item[(1)] $\bR$ is maximally singular,
\item[(2)] $\bR$ is given by some height function $\xi$,
\item[(3)] $|a-b|\leq d(i,j)$ for all $i,j\in I$, $a\in R_i$ and $b\in R_j$.
\end{itemize}
\end{Theorem}

\begin{proof}[Partial proof]
The equivalence between (2) and (3) is clear. For $(1)\Longrightarrow (3)$, suppose that $|a-b|>d(i,j)$ for some $i,j\in I$, $a\in R_i$ and $b\in R_j$. Assume also, without loss of generality, that $a\leq b$ and fix a minimal path
$ j = k_1\rightarrow k_2 \rightarrow \dots \rightarrow k_{d+1}=i$
with $d=d(i,j)$. Consider %
the monomial
$$ x=\tilde{f}_{k_{d+1}}\dots \tilde{f}_{k_2}\tilde{f}_{k_1}y_{j,b}. $$
Then%
, as is easily shown,
$$ \textstyle x = y_{j,b}z_{k_1,b-2}^{-1}z_{k_2,b-3}^{-1}\dots z_{k_{d+1},b-d-2}^{-1} \in y_{i,b-d-2}^{-1}\cB_+$$
and it follows in particular that $x$ belongs to the crystal $\cB(\varpi_j^{\vee},b)$ (since it is non-zero).~Thus, letting $y=y_{\bR}y_{j,b}^{-1}y_{i,a}^{-1}\in \cB_+$, we get that $xyy_{i,a} \in \cB(\la,\bR)$ with
\begin{equation}\label{eq:MaxSing1}
\textstyle xyy_{i,a} \in \frac{y_{i,a}}{y_{i,b-d-2}}\cB_+.
\end{equation}
Clearly $\tilde{e}_{\ell}(xyy_{i,a})=0$ if $\ell \neq i$. Moreover, $\tilde{e}_i(\frac{y_{i,a}}{y_{i,b-d-2}})=0$ (as $b-a> d$ by assumption) and it follows from \eqref{eq:MaxSing1} and the definition of the crystal operator $\tilde{e}_i$ that $\tilde{e}_i(xyy_{i,a})=0$. Hence, the monomial $xyy_{i,a}$ generates a component of $\cB(\la,\bR)$ that is disjoint from~the~one~generated by the dominant monomial $y_{\bR}$. In particular, %
$\bR$ is not maximally singular.
This ends the proof of $(1)\Longrightarrow (3)$. The proof of %
$(2)\Longrightarrow (1)$, which essentially uses \cite[Corollary~5.14]{gibson2021demazure}, is postponed to the end of Section \ref{subsec:filt_of_R}.
\end{proof}
\begin{Rem} As mentioned in \cite{kamnitzer2019highest}, it would be interesting to characterize all~sets of parameters $\bR\in \C^{\la}$ for which $\cB(\la,\bR)\simeq \bigotimes_{i\in I}B(\varpi_i^{\vee})^{\otimes \langle \la,\alpha_i\rangle}$. Another fascinating open problem related to %
product monomial crystals is to find what interpolating representations $V(\la,\bR)$ of $\fg^{\vee}$ (between $V(\la)$ and $\bigotimes_{i\in I}V(\varpi_i^{\vee})^{\otimes \langle \la,\alpha_i\rangle}$) can give rise to a crystal $\cB(\la,\bR)$~(see \cite{gibson2021demazure} and the upcoming Theorem \ref{th:charcR} for partial results in this direction%
). %
\end{Rem}
\subsection{Non-integral case}\label{sec:nonInt}
The goal of this subsection is to prove a %
 version of Corollary \ref{co:directed} where the sets %
$\bR_1,\bR_2$ are not assumed to be integral. For this, denote by $\cB[r]$ (resp.~$\Gamma[r]$) the set of Laurent monomials in the $y_{i,a+r}$'s (resp.~$z_{i,a+r}$) where $(i,a)\in I\times_2\Z$ and $r\in \C$. Let also $\cB_{+}[r]$ and $\Gamma_{+}[r]$ be the associated sets of monomials and define, for a subset $S\subseteq \C$,
$$ \textstyle
\cB_S = \prod_{r\in S}\cB[r],\quad \Gamma_S = \prod_{r\in S}\Gamma[r], \quad \cB_{S,+} = \prod_{r\in S} \cB_+[r]\  \text{ and }\ \Gamma_{S,+} = \prod_{r\in S}\Gamma_+[r].$$
Then $\cB_{\C}$ contains all the (possibly non-integral) product monomial crystals of Remark \ref{rem:nonIntCryst} and $\cB_{\emptyset}=\{1\}$ by convention. Also, $\cB_{S,+}=\cB_S\cap \cB_{\C,+}$ and $\Gamma_{S,+}=\Gamma_S\cap\Gamma_{\C,+}$ by definition. %
\begin{Lemma}\label{lem:GammaS} Fix $S\subseteq \C$. Then $\Gamma_S=\cB_S\cap \Gamma_{\C}$%
.
\end{Lemma}
\begin{proof} The inclusion $\subseteq$ is clear. For the other one, fix $x\in \Gamma_{\C}$ and suppose $x\not\in \Gamma_S$. Then
\begin{equation}\label{eq:xLemNonIntGamma}
x=z_{i_1,a_1}^{r_1}\dots z_{i_s,a_s}^{r_s}
\end{equation}
for some distinct pairs $(i_1,a_1),\dots,(i_s,a_s) \in I\times \C$ and non-zero $r_1,\dots r_s\in \Z$. Also, there~is $1\leq k\leq s$ for which $a_k\not\in \overline{i_k}+S+2\Z$. Let $1\leq K\leq s$ be such that $a_K$ is minimal\footnote{Here ``minimal'' means that the real part of $a_K$ is minimal among the real parts of the elements of \eqref{eq:setLemSC}.} in%
\begin{equation}\label{eq:setLemSC}
\{a\in \C\,|\, a = a_k\text{ for some }1\leq k\leq s \text{ and }a\not\in \overline{i_k}+S+2\Z\}\neq \emptyset. 
\end{equation}
Then the variable $y_{i_K,a_K}$ appearing (with exponent $r_K$) in $z_{i_K,a_K}^{r_K}$ is not cancelled by another variable %
in \eqref{eq:xLemNonIntGamma}. Thus $x\not\in \cB_S$ and the result follows.
\end{proof}

Take $y_{\bR_1},y_{\bR_2}\in \cB_{S,+}$ for some $S\subseteq \C$ and suppose $y_{\bR_1}\Qq{1}=y_{\bR_2}\Qq{2}$ where $\Qq{1},\Qq{2}\in \Gamma_{\C,+}$. By Lemma \ref{lem:GammaS},
$$\textstyle \frac{\Qq{1}}{\Qq{2}} =\frac{y_{\bR_2}}{y_{\bR_1}} \in \cB_S\cap \Gamma_{\C}=\Gamma_{S}$$
so that %
$\Qq{1},\Qq{2}\in \Gamma_{S,+}$ without loss of generality. Choose $a\in S$ and write $S_a=S\backslash(a+2\Z)$. Then $\cB_{S,+}=\cB_{+}[a]\cB_{S_a,+}$ and there are thus factorizations
$$ y_{\bR_1} = y_1y_1' \text{ and } y_{\bR_2}=y_2y_2'$$
where $y_1,y_2\in \cB_+[a]$ and $y_1',y_2'\in \cB_{S_a,+}$. Similarly, $\Gamma_{S,+}=\Gamma_+[a]\Gamma_{S_a,+}$ and
$$ \Qq{1} = q_1q_1' \text{ with } \Qq{2} = q_2q_2'$$
for some $q_1,q_2\in \Gamma_+[a]$ and $q_1',q_2'\in \Gamma_{S_a,+}$. Hence, using $y_{\bR_1}\Qq{1}=y_{\bR_2}\Qq{2}$, we get
$$ \textstyle \frac{y_1'q_1'}{y_2'q_2'}=\frac{y_2q_2}{y_1q_1}\in \cB[a]\cap\cB_{S_a}=\{1\} $$
so that $y_1q_1=y_2q_2$ with $y_1'q_1'=y_2'q_2'$%
. %
This leads to our generalization of Corollary \ref{co:directed}:

\begin{Lemma}\label{lem:InductiveSystemC} 
Choose $\la,\la'\in P_+^{\vee}$ with $\bR\in \C^{\la}$ and $\bR'\in \C^{\la'}$ such that $y_{\bR}\in y_{\bR'}\Gamma_{\C}$. Then, there exists $\la_{tot}\in P_+^{\vee}$ and $\bR_{tot}\in \C^{\la_{tot}}$ satisfying $\cB(\la,\bR)\cup\cB(\la',\bR')\subseteq \cB(\la_{tot},\bR_{tot})$.
\end{Lemma}

\begin{proof}
Fix $q,q'\in \Gamma_{\C,+}$ for which $y_{\bR}q=y_{\bR'}q'$. Take also $a_1,\dots,a_r\in\C$ with $a_k-a_{\ell}\not\in 2\Z$ if $k\neq \ell$ and such that $y_{\bR},y_{\bR'}\in \cB_{S,+}$ for $S=\{a_1,\dots,a_r\}$. Then
$$y_{\bR}=y_1\dots y_r\ \text{ and }\ y_{\bR'}=y_1'\dots y_r'$$
for some $y_1,\dots,y_r,y_1',\dots,y_r' \in \cB_{S,+}$ with $y_k,y_k'\in \cB_+[a_k]$ for all $1\leq k\leq r$.\medskip\par
Fix $k\in \{1,\dots,r\}$ and write
$$ \textstyle
y_k =\prod_{i\in I}\prod_{b\in R_{k,i}}y_{i,b} = y_{\bR_{k}}\text{ with }y_k' =\prod_{i\in I}\prod_{b\in R_{k,i}'}y_{i,b} = y_{\bR_{k}'}$$
for sets $\bR_k=(R_{k,i})_{i\in I}$ and $\bR_k'=(R_{k,i}')_{i\in I}$ of parameters satisfying $R_{k,i}\subseteq a_k+\overline{i}+2\Z$ for every $i\in I$. Set also
$$ \textstyle \lambda_k = \sum_{i\in I}\lambda_{k,i}\varpi_i^{\vee}\text{ and }\lambda_k'=\sum_{i\in I}\lambda_{k,i}'\varpi_i^{\vee}$$
with $\lambda_{k,i}$ (resp.~$\lambda_{k,i}'$) the cardinality of $R_{k,i}$ (resp.~$R_{k,i}'$) for all $i\in I$. Then the definition of the product monomial crystals $\cB(\la,\bR)$ and $\cB(\la',\bR')$ gives
$$ \textstyle \cB(\la,\bR)=\prod_{k=1}^r \cB(\lambda_k,\bR_k)\text{ and }\cB(\la',\bR')=\prod_{k=1}^r \cB(\lambda_k',\bR_k').$$
On the other hand, %
as before the lemma, we can find $q_1,\dots,q_r,q_1',\dots,q_r'\in \Gamma_{S,+}$ satisfying
$$q=q_1\dots q_r\text{ with }q'=q_1'\dots q_r',$$
as well as $q_k,q_k'\in \Gamma_{+}[a_k]$ and $y_{\bR_k}q_k=y_{\bR_k'}q_k'$ for all $1\leq k\leq r$. Hence, fixing $k\in\{1,\dots,r\}$, by (the ``$a_k$-translate'' of) Lemma \ref{th:qm}, there is $p_k\in \Gamma_+[a_k]$ with $y_{\bR_k}q_kp_k \in \cB_+[a_k]$. Write
$$ \textstyle
y_{\bR_k}q_kp_k = y_{\bR_k'}q_k'p_k = \prod_{i\in I}\prod_{b\in R^{(k)}_i}y_{i,b}=y_{\bR^{(k)}}$$
for some set $\bR^{(k)}=(R_i^{(k)})_{i\in I}$ of parameters with $R_i^{(k)}\subseteq a_k+\overline{i}+2\Z$ for all $i\in I$. Set also
$$\textstyle\lambda^{(k)}=\sum_{i\in I}\lambda^{(k)}_i\varpi_i^\vee$$
where $\lambda_i^{(k)}$ is the cardinality of $R^{(k)}_i$. %
Then (the ``$a_k$-translate'' of) Theorem \ref{th:TFAE} gives
$$\cB(\la_k,\bR_k)\cup\cB(\la_k',\bR_k')\subseteq \cB(\lambda^{(k)},\bR^{(k)})$$
for all $1\leq k\leq r$. Thus, using again \eqref{eq:MonCrys} with $\lambda_{tot}=\sum_{k=1}^r\lambda^{(k)}\in P_+^{\vee}$ and $\textstyle \bR_{tot}=%
\bigcup_{k=1}^r \bR^{(k)}%
$
(the %
union of %
multisets%
), we get that
$$\textstyle \cB(\la,\bR)\cup\cB(\la',\bR')\subseteq \prod_{k=1}^r(\cB(\la_k,\bR_k)\cup\cB(\la'_k,\bR'_k))\subseteq \prod_{k=1}^r\cB(\lambda^{(k)},\bR^{(k)}) =\cB(\la_{tot},\bR_{tot})$$
as claimed. This ends the proof.
\end{proof} 
 			
\section{Truncated shifted Yangians}\label{sec:TSY}
This section recalls the definition of truncated shifted Yangians as images of remarkable homomorphisms \cite{gerasimov2005class} and gives the important properties of the shifted coproducts~$\Delta_{\mu,\nu}$ of \cite{finkelberg2018comultiplication}. Our exposition %
closely follows that of \cite{finkelberg2018comultiplication,hernandez2024shifted,kamnitzer2022lie}.
\subsection{Shifted Yangians}\label{sec:Algebras} %
As mentioned in the introduction, dominantly shifted Yangians of type A were first introduced in \cite{brundan2006shifted} as subalgebras of Drinfeld's Yangian $Y(\mathfrak{gl}_n)$. This definition was then extended to all simple types in \cite{kamnitzer2014yangians}, and, using ``\textit{Cartan doubled Yangians}'', to arbitrary (i.e.~not necessarily dominant) shifts in \cite{braverman2016coulomb,finkelberg2018comultiplication}.
\begin{Def}[\cite{braverman2016coulomb,finkelberg2018comultiplication}]\label{cartandoubleyangian}
The {\em Cartan doubled Yangian} $Y_\infty %
$ (of $\fg$) is the algebra generated by elements $ e_{i,q}, f_{i,q}, h_{i,p} $ for $ i\in I$, $ q \in \mathbb{Z}_{>0} $ and $ p \in \mathbb{Z} $, with defining relations\footnote{Recall that we only work with simply-laced Lie algebras. Hence $\alpha_i\cdot\alpha_j=c_{ij}$ in the notation of \cite{braverman2016coulomb}. Moreover, we adopt here the convention that $\hbar=2$. This gives a slightly different presentation than the usual one (which can be recovered using $h_{i,p}\mapsto 2^{\langle \mu,\alpha_i\rangle+p}h_{i,p}$, %
$e_{i,q}\mapsto 2^{\langle \mu,\alpha_i\rangle+q}e_{i,q}$ and $f_{i,q}\mapsto 2^{q}f_{i,q}$%
).
}
\begin{align}
\label{H,H} [h_{i,p}, h_{j,p'}] &= 0,  \\
\label{E,F} [e_{i,q}, f_{j,q'}] &= %
2\delta_{ij}h_{i,q+q'-1},\\
\label{H,E} [h_{i,p+1},e_{j,q}] - [h_{i,p}, e_{j,q+1}] &= %
c_{ij} (h_{i,p} e_{j,q} + e_{j,q} h_{i,p}), \\
\label{H,F} [h_{i,p+1},f_{j,q}] - [h_{i,p}, f_{j,q+1}] &= %
-c_{ij} (h_{i,p} f_{j,q} + f_{j,q} h_{i,p}) , \\
\label{E,E} [e_{i,q+1}, e_{j,q'}] - [e_{i,q}, e_{j,q'+1}] %
&= c_{ij} (e_{i,q} e_{j,q'} + e_{j,q'} e_{i,q}), \\
\label{F,F} [f_{i,q+1}, f_{j,q'}] - [f_{i,q}, f_{j,q'+1}] %
&= -c_{ij} (f_{i,q} f_{j,q'} + f_{j,q'} f_{i,q}),\\
i \neq j,\,N = 1 - c_{ij} \Longrightarrow
\label{symE} \operatorname{sym} &[e_{i^,q_1}, [e_{i,q_2}, \cdots [e_{i,q_N}, e_{j,q'}]\cdots]] = 0, \\
i \neq j,\,N = 1 - c_{ij} \Longrightarrow
\label{symF} \operatorname{sym} &[f_{i,q_1}, [f_{i,q_2}, \cdots [f_{i,q_N}, f_{j,q'}]\cdots]] = 0.
\end{align}
\end{Def}
The algebra $Y_{\infty}$ admits all shifted Yangians (associated to $\mathfrak{g}$) as quotients.
\begin{Def}[\cite{braverman2016coulomb,finkelberg2018comultiplication}]\label{shiftedyangian}
Fix $\mu\in P^{\vee}$ arbitrary. The shifted Yangian $Y_\mu$ (of $\fg$) is the quotient of $Y_\infty$ by the relations $h_{i,-\langle \mu,\alpha_i\rangle}=1$ and $h_{i,p}=0$ for $p<-\langle \mu,\alpha_i\rangle$.
\end{Def}
\begin{Rem}
When $\mu=0$, the algebra $Y=Y_0$ is the usual Yangian of $\mathfrak{g}$ (with the above generators and relations recovering its Drinfeld ``new presentation''). %
\end{Rem}
Fix $\mu\in P^{\vee}$. We will use the series $e_i(u) = \sum_{q\geq 1}e_{i,q}u^{-q}$, $f_i(u)=\sum_{q\geq 1} f_{i,q}u^{-q}$ and
$$\textstyle h_i(u)= \sum_{p\in \mathbb{Z}} h_{i,p}u^{-p} =u^{\langle \mu,\alpha_i\rangle} + \sum_{p\geq 1} h_{i, -\langle \mu,\alpha_i\rangle + p} u^{\langle \mu,\alpha_i\rangle-p}$$
for $i\in I$ as well as the triangular decomposition \cite{finkelberg2018comultiplication}
\begin{equation}\label{eq:TriangularDec}
Y_{\mu} \simeq Y_{\mu}^-\otimes Y_{\mu}^0\otimes Y_{\mu}^+
\end{equation}
with $Y_{\mu}^-$, $Y_{\mu}^0$ and $Y_{\mu}^+$ respectively the subalgebras of $Y_{\mu}$ generated by the $f_{i,q}$'s, the $h_{i,p}$'s and the $e_{i,q}$'s (with $i\in I$, $q\in \mathbb{Z}_{>0}$ and $p>-\langle \mu,\alpha_i\rangle$). We will also need the shift morphisms
$$\iota_{\mu,\mu_1,\mu_2}: Y_\mu \rightarrow Y_{\mu+\mu_1+\mu_2}$$ defined in \cite{finkelberg2018comultiplication} for all antidominant coweights $\mu_1,\mu_2\in-P^{\vee}_+$ by
$$ e_{i,q} \mapsto e_{i,q-\langle\mu_1,\alpha_i\rangle},\quad f_{i,q} \mapsto f_{i,q-\langle\mu_2,\alpha_i\rangle} \ \text{ and }\ h_{i,p}\mapsto h_{i,p-\langle\mu_1+\mu_2,\alpha_i\rangle}.$$
\begin{Lemma}[\cite{finkelberg2018comultiplication}] Shift morphisms are injective. Also, for $\mu\in P^{\vee}_+$, the map $\iota_{\mu,0,-\mu}$ recovers the %
definition of dominantly shifted Yangians as subalgebras of the %
Yangian $Y$.
\end{Lemma}

\subsection{Truncations}
\label{section: tsy}
Fix $\lambda\in P^{\vee}_+$ and $\bR=(R_i)_{i\in I}\in \C^{\lambda}$ as in Section \ref{sec:Crystal}. Define polynomials
\begin{equation}
\label{eq:polyp}
\textstyle p_{R_i}(u) = \prod_{c \in R_i} (u-c)
\end{equation}
and fix $\mu\in P^{\vee}$ such that $\lambda-\mu=\sum_{i\in I} m_i\alpha_i^{\vee}$ for some $(m_i)_{i\in I}\in \mathbb{N}^I$ (i.e. $\mu\leq \lambda$). Following \cite[Lemma 2.1]{gerasimov2005class}, there exist in $Y_{\mu}^0$ unique elements $(a_{i,r})_{i\in I,r\in\Z_{>0}}$ such~that, for $i\in I$, %
\begin{equation}
\label{eq: def of A gens}
h_i(u) = p_{R_i}(u) \frac{\prod_{j \sim i} (u-1)^{m_j}}{u^{m_i} (u-2)^{m_i}} \frac{\prod_{j \sim i} a_j(u-1 )}{a_i(u) a_i(u-2)},
\end{equation}
with $a_i(u) = 1 + \sum_{r \geq 1 } a_{i,r} u^{-r}$. These elements, which depend intrinsically on $\lambda$, $\mu$ and $\bR$, are one of the two key ingredients needed to define truncated shifted Yangians. The second one, also due to \cite{gerasimov2005class}, is the algebra of difference operators $ \widetilde{\mathscr A}_{\lambda-\mu} $ defined below.
\begin{Def}[\cite{gerasimov2005class}]\label{def:DiffOpAlg} Let $\mathscr{A}_{\lambda-\mu}$ be the algebra generated by elements
$$\{w_{i,r},\beta_{i,r}^{\pm 1}\}_{i\in I,1\leq r\leq m_i},$$ with all commutators trivial except
$$
[\beta_{i,r}^{\pm 1}, w_{i,r} ] = \pm 2\beta_{i,r}^{\pm 1}
$$
for $i\in I$ and $1\leq r\leq m_i$. Let $\widetilde{\mathscr{A}}_{\lambda-\mu}$ be the Ore localization of $\mathscr{A}$, i.e.~the ring obtained by formally introducing inverses
$(w_{i,r} - w_{i,s}+ k )^{-1}$
for all $i \in I$, $1\leq r \neq s \leq m_i$ and $k \in 2\mathbb{Z}$.
\end{Def}
The algebra of difference operators 
$\widetilde{\mathscr{A}}_{\lambda-\mu}$ can be related to the algebra $Y_{\mu}$ by a remarkable morphism %
which was first introduced (for $\mu=0$) in \cite{gerasimov2005class}, further studied in \cite{kamnitzer2014yangians}, %
and connected with the theory of Coulomb branches in \cite[Appendix~B]{braverman2016coulomb}. To recall~the definition of this morphism, set, for $i\in I$ and $1\leq r\leq m_i$,
$$\textstyle W_i (u) = \prod_{1\leq s\leq m_i}(u - w_{i,s}) \ \text{ and } \ W_{i,r}(u) = \prod_{\substack{1\leq s\leq m_i \\ s\neq r}} (u - w_{i,s}).$$

\begin{Theorem}[\mbox{\cite[Theorem B.15]{braverman2016coulomb}}]\label{GKLO homomorphism}
Fix any orientation of the Dynkin diagram of~$\fg$. Then there is a homomorphism of algebras
$$\Phi_\mu^\lambda(\bR): Y_\mu \longrightarrow \widetilde{\mathscr A}_{\lambda-\mu},$$
uniquely determined by the assignments $a_i(u) \mapsto u^{-m_i} W_i(u)$,
\begin{align*}
e_i(u) & \mapsto -\sum_{r=1}^{m_i}  \frac{\prod_{j\rightarrow i} W_j(w_{i,r} - 1)}{(u-w_{i,r}) W_{i,r}(w_{i,r})} \beta_{i,r}^{-1} \text{ and} \\
f_i(u) & \mapsto \sum_{r=1}^{m_i} \frac{p_{R_i}(w_{i,r}+2) \prod_{j\leftarrow i} W_j(w_{i,r} + 1 )}{(u-w_{i,r} - 2) W_{i,r}(w_{i,r})} \beta_{i,r}.
\end{align*}
\end{Theorem}

\begin{Def}
The \emph{truncated shifted Yangian} $Y_\mu^\lambda(\bR)$ is defined to be the image of $\Phi_\mu^\lambda(\bR)$.
\end{Def}
\begin{Rem}
\label{rem:tsyideal}
Some articles (e.g.~\cite{kamnitzer2019highest,hernandez2024shifted}) use \textit{truncated shifted Yangian} for the quotient %
of $Y_{\mu}$ by the 2-sided ideal generated by
$ \{a_{i,r}\,|\,i\in I, r>m_i\}$.~It~is easy to see that there is a surjection of this quotient onto $Y_{\mu}^{\lambda}(\mathbf{R})=\text{Im\,} \Phi_{\mu}^{\lambda}(\mathbf{R})$, but results of \cite{Besson2026} %
show that it is not an isomorphism in general (though it is so in type A, as shown in \cite{kamnitzer2019highest}).
\end{Rem}
\subsection{Shifted coproducts}\label{sec:Coproduct}
As is well-known since the pioneering work of Drinfeld \cite{drinfeld1987new}, the Yangian $Y=Y_0$ admits a coproduct $\Delta:Y\rightarrow Y\otimes Y$. To describe it, recall the inclusion $U(\mathfrak{g})\subseteq Y$ defined on Chevalley generators by $e_i\mapsto \frac{1}{2}e_{i,1}$, $f_i\mapsto  \frac{1}{2}f_{i,1}$ and $h_i\mapsto \frac{1}{2} h_{i,1}$. Then, $\Delta$ is uniquely determined by
\begin{equation}\label{eq:Delta1}
\Delta(x) = x\otimes 1+1\otimes x
\end{equation}
for $x\in U(\mathfrak{g})$ and
\begin{equation}\label{eq:Delta2}
\textstyle \Delta(h_{i,2}) = h_{i,2}\otimes 1+1\otimes h_{i,2}+h_{i,1}\otimes h_{i,1}+\sum_{\gamma\in \Delta_+}\langle \alpha_i,\gamma\rangle f_{\gamma}\otimes e_{\gamma}
\end{equation}
for $i\in I$, where $\Delta_+$ is the set of positive roots of $\mathfrak{g}$ and where $e_{\gamma}\in \mathfrak{g}_{\gamma}$ and $f_{\gamma}\in \mathfrak{g}_{-\gamma}$ suitably normalized root vectors (cf.~\cite[Section 3.12]{finkelberg2018comultiplication} or \cite[Section 4.2]{guay2018coproduct} for details). \medskip\par
The following result introduces analogues of the coproduct $\Delta$, called \textit{shifted coproducts}, for other shifted Yangians. These analogues will play a crucial role in the rest of the article.
\begin{Theorem}[\cite{finkelberg2018comultiplication}]\label{thm:CoproductUnicity} There exists a unique family $\{\Delta_{\mu,\nu}: Y_{\mu+\nu}\rightarrow Y_{\mu}\otimes Y_{\nu}\}_{\mu,\nu\in P^{\vee}}$ of algebra morphisms such that 
\begin{enumerate}[label=(\roman*)]
\item $\Delta_{0,0}$ is the usual coproduct for the Yangian $Y=Y_0$,
\item $\Delta_{\mu,\nu}(e_{i,n})=e_{i,n}\otimes 1$ for $n<-\langle\mu,\alpha_i\rangle$ and $\mu,\nu$ antidominant,
\item $\Delta_{\mu,\nu}(f_{i,m})=1\otimes f_{i,m}$ for $m<-\langle\nu,\alpha_i\rangle$ and $\mu,\nu$ antidominant,
\item the diagram 
{
$$
\adjustbox{scale=0.88,center}{
\begin{tikzcd}[column sep = 4.5em]
Y_{\mu+\nu} \ar[r,"\Delta_{\mu,\nu}"] \ar[d,"\iota_{\mu+\nu,\zeta,\eta}",swap] & Y_{\mu}\otimes Y_{\nu}\ar[d,"\iota_{\mu,\zeta,0}\otimes\iota_{\nu,0,\eta}"] \\
Y_{\mu+\nu+\zeta+\eta} \ar[r,"\Delta_{\mu+\zeta,\nu+\eta}"] & Y_{\mu+\zeta} \otimes Y_{\nu+\eta}
\end{tikzcd}}
$$}\noindent
commutes for $\zeta,\eta$ antidominant (where the $\iota$'s are the shift morphisms), and
\item the diagram %
{
$$
\adjustbox{scale=0.88,center}{
\begin{tikzcd}[column sep = 4.5em]
Y_{\mu+\rho+\nu} \ar[r,"\Delta_{\mu+\rho,\nu}"] \ar[d,"\Delta_{\mu,\rho+\nu}",swap] & Y_{\mu+\rho}\otimes Y_{\nu}\ar[d,"\Delta_{\mu,\rho}\otimes 1"] \\
Y_{\mu}\otimes Y_{\rho+\nu} \ar[r,"1\otimes \Delta_{\rho,\nu}"] & Y_{\mu}\otimes Y_{\rho}\otimes Y_{\nu}
\end{tikzcd}}
$$}\noindent
commutes for $\rho$ antidominant.
\end{enumerate}
\end{Theorem}
	
\section{Representations of shifted Yangians}\label{sec:rep_of_shifted_yangians}
This section recalls facts about the representation theory of (truncated) shifted Yangians and defines the important categories for the sections to come. In particular, we recall from \cite{hernandez2024shifted} the definition of the category $\mathcal{O}_{sh}$ mentioned in Section \ref{sec:Intro} and study ``\textit{truncations}'' $\mathcal{O}_{sh}^{\lambda}(\mathbf{R})$ of this category using results of \cite{kamnitzer2019category}. We also give a new \textit{block decomposition} for $\mathcal{O}_{sh}$ using \textit{$\ell$-characters} and study the notion of \textit{GT-weights}\footnote{Here ``\textit{GT}'' stands for ``\textit{Gelfand-Tsetlin}''.} (as defined in \cite{kamnitzer2019category}). %
\subsection{The category $\mathcal{O}_{sh}$}\label{sec:Osh} There are three notions of \textit{weight spaces} that one may encounter while studying modules over shifted Yangians. More precisely, fixing $\mu\in P^{\vee}$ and a representation $V$ of $Y_{\mu}$, one can consider either:
\begin{enumerate}
\item simultaneous \textbf{generalized}\footnote{The fact that we use generalized eigenspaces here instead of genuine eigenspaces is slightly non-standard (see upcoming Remark \ref{rem:2categoriesO}), but is necessary if we want to use the results of \cite{kamnitzer2019category} on KLRW algebras%
.} eigenspaces for the finitely-many
commuting generators $h_{i,1-\langle \mu,\alpha_i\rangle}$ (with $i\in I$), that is spaces%
\begin{equation*}
V_{\omega} = \{v\in V\,|\, \exists p\in \mathbb{N} \text{ such that } (h_{i,1-\langle \mu,\alpha_i\rangle}-2\langle\alpha_i^{\vee},\omega\rangle)^pv=0 \text{ for all } i\in I\}
\end{equation*}
defined from an element $\omega \in \mathfrak{h}^*$, or%
\item simultaneous generalized eigenspaces for the infinitely-many commuting generators $h_{i,r}$ (with $i\in I$ and $r>-\langle \mu,\alpha_i\rangle$), that is spaces
$$ V_{\psi} = \{v\in V\,|\, \exists p\in \mathbb{N}\text{ such that } (h_{i,r}-\psi_{i,r})^p v=0 \text{ for all } i\in I \text{ and }r>-\langle\mu,\alpha_i\rangle\}$$
defined from a sequence $\psi =(\psi_{i,r})_{i\in I,r>-\langle \mu,\alpha_i\rangle}$ of complex numbers, or, lastly,
\item simultaneous generalized eigenspaces for the %
finitely-many
commuting\footnote{Clearly, $[a_{i,r},a_{j,s}]=0$ for $i,j\in I$ and $r,s\in \mathbb{Z}_{>0}$ as all $a_{i,r}$'s lie in the commutative subalgebra $Y_{\mu}^0\subseteq Y_{\mu}$.
Moreover, as we will see in Section \ref{subsec:GT_weights}, using $a$'s instead of $h$'s makes essentially no difference for modules~$V$ over the truncation $Y_{\mu}^{\la}(\bR)$ (which are $Y_{\mu}$-modules by pullback via the defining projection $Y_{\mu}\twoheadrightarrow Y_{\mu}^{\la}(\bR)$).
}
 elements~$a_{i,r}$
 (with $i\in I$ and $1\leq r\leq m_i$)
 defined in \eqref{eq: def of A gens} for a given choice of
 dominant
 coweight $\lambda\in P^{\vee}_+$ satisfying $\lambda-\mu=\sum_{i\in I}m_i\alpha_i^{\vee}\in Q^{\vee}_+$
 and set of parameters $\mathbf{R}\in \C^{\lambda}$
 .
\end{enumerate}
A \textit{weight of $V$} %
is an element $\omega\in \mathfrak{h}^*$ for which the generalized eigenspace $V_{\omega}$ (called \textit{weight space}) is non-zero. Analogously, an \textit{$\ell$-weight of $V$} is a sequence $\psi$ of complex numbers~as in (2) %
for which the generalized eigenspace $V_{\psi}$ (referred to as an \textit{$\ell$-weight space}) is non-zero. These notions (with the concept of \textit{GT-weight} that will be defined in Section \ref{subsec:GT_weights} using the last type of eigenspaces above) are critical for the study of the representation theory of the shifted Yangian $Y_{\mu}$. In particular, the first notion allows us to define a \textit{BGG category~$\mathcal{O}$}~in the category $Y_{\mu}$--Mod of all (left) $Y_{\mu}$-modules.
\begin{Def}[{\cite[Section 5]{kamnitzer2022lie}}]\label{def:Omu} The category $\mathcal{O}_{\mu}$ is the full subcategory of $Y_{\mu}$--Mod containing the objects $V$ satisfying the conditions:
\begin{itemize}
\item[($\mathcal{O}1$)] $V=\bigoplus_{\omega\in \mathfrak{h}^*} V_{\omega}$ as vector spaces,
\item[($\mathcal{O}2$)] $\dim V_{\omega} < \infty$ for all $\omega\in \mathfrak{h}^*$, and
\item[($\mathcal{O}3$)] there are $\omega_1,\dots,\omega_N\in \mathfrak{h}^*$ such that, for all $\omega\in \mathfrak{h}^*$,
$$V_{\omega}\neq 0 \Longrightarrow %
\omega\leq \omega_k
\text{ for some }1\leq k\leq N.$$
\end{itemize}
\end{Def}
\begin{Rem}\label{rem:2categoriesO} In \cite[Section 3.3]{hernandez2024shifted}, the authors give another definition for the category $\mathcal{O}_{\mu}$ where the weight spaces $V_{\omega}$ are replaced with actual eigenspaces, i.e.~spaces of the form
$$ \{v\in V\,|\, h_{i,1-\langle \mu,\alpha_i\rangle}v=\omega_i v\text{ for all }i\in I\}.$$
The resulting category $\mathcal{O}_{\mu}^{\text{HZ}}$ --- which is properly contained in the one of Definition \ref{def:Omu}~--- is, as will be detailed in Remark \ref{rem:OHZ_is_too_small}, too small for our purposes. Nevertheless, as we will explain below, both %
$\mathcal{O}_{\mu}$ and $\mathcal{O}_{\mu}^{\text{HZ}}$ lead to the same collection of simple objects (and hence to the same Grothendieck group), allowing us to apply directly most of Hernandez--Zhang's results to our (slightly more general) context.
\end{Rem}
By \cite[Lemma 2.5]{hernandez2024shifted} and Definition \ref{def:Omu}, the shifted coproducts $\Delta_{\mu,\nu}:Y_{\mu+\nu}\rightarrow Y_{\mu}\otimes Y_{\nu}$ of Section \ref{sec:Coproduct} naturally endow the
direct sum of abelian categories
$$\textstyle \mathcal{O}_{sh}=\bigoplus_{\mu\in P^{\vee}} \mathcal{O}_{\mu}$$
with a tensor product $\otimes$ (defined by pullback) and provide in particular the group $K_0(\mathcal{O}_{sh})$ with a ring structure. Now, as mentioned in Section \ref{sec:Intro}, the primordial goal of this paper is to show that the ring $K_0(\mathcal{O}_{sh})$ (or rather the complexified Grothendieck ring of a ``monoidal skeleton'' of $\mathcal{O}_{sh}$) is isomorphic to the coordinate ring of the \textit{scheme of bands} of  \cite{francone2025cluster}. To do this, however, we first need a better understanding of the simple objects in the category $\mathcal{O}_{sh}$ and of the tensor structure on this category. 
\begin{Rem} The coproducts $\Delta_{\mu,\nu}$ are not co-associative (i.e.~the diagram in property~(v) of Theorem \ref{thm:CoproductUnicity} does not commute for arbitrary $\rho$) and it is not known if the category $\mathcal{O}_{sh}$ is monoidal%
. Nevertheless,~we~will slightly abuse terminology in this article and use~the~term ``monoidal'' when talking about categorifications of cluster algebras or other related notions.
\end{Rem}\newpage
We follow \cite{hernandez2024shifted} and package, for $\mu\in P^{\vee}$, sequences $\psi=(\psi_{i,r})_{i\in I,r>-\langle\mu,\alpha_i\rangle}$ of complex numbers into $I$-tuples of Laurent series $(\psi_i(u))_{i\in I}$, where
\begin{equation}\label{eq:lwt}
\textstyle  \psi_i(u) = u^{\langle \mu,\alpha_i\rangle}+ \sum_{r>-\langle \mu,\alpha_i\rangle } \psi_{i,r}u^{-r}%
\end{equation}
for $i\in I$. Using this, the $\ell$-weight space $V_{\psi}$ of a $Y_{\mu}$-module $V$ can be alternatively defined as a simultaneous generalized eigenspace for the (finitely-many) currents $h_i(u)$ (with $i\in I$), i.e.~as the space
\begin{equation*}
\smash{V_{\psi} = \{v\in V\,|\,\exists p\text{ such that }(h_i(u)-\psi_i(u))^pv=0 \text{ for all } i\in I\}.}
\end{equation*}
Denote by $\mathfrak{r}_{\mu}$ the set of Laurent series $\psi=(\psi_{i}(u))_{i\in I}$ as in \eqref{eq:lwt} that are expansions~at~$\infty$ of rational functions in $u$, and write $$\textstyle\smash{\mathfrak{r}=\bigsqcup_{\mu\in P^{\vee}}\mathfrak{r}_{\mu}.}$$ 
For $\psi=(\psi_i(u))_{i\in I}\in \C(u)^I$,
$$\textstyle \psi\in \mathfrak{r}\iff \text{all }\psi_i(u)\text{'s can be written as }\psi_i(u)=\frac{p_i(u)}{q_i(u)} \text{ for monic }p_i(u),q_i(u)\in \C[u].$$
Moreover, $\mathfrak{r}_{\mu}$ is simply the subset of $\mathfrak{r}$ where the polynomials $(p_i(u),q_i(u))_{i\in I}$ on the right-hand side of the above equivalence satisfy
$\deg p_i(u)-\deg q_i(u)=\langle \mu,\alpha_i\rangle$
for all $i\in I$.
\begin{Def}%
\label{def:highestlwtmod}
Fix a $Y_{\mu}$-module $V$ and let $\psi=(\psi_i(u))_{i\in I}\in \mathfrak{r}_{\mu}$. Then $V$~is~said to be
of \textit{highest $\ell$-weight $\psi$} if there exists $v\in V$ such that $V=Y_{\mu}^-v$ with
$$ h_i(u)v = \psi_i(u)v\ \text{ and }\ e_i(u)v=0$$
for all $i\in I$. The vector $v\in V$ is then called a \textit{highest $\ell$-weight vector} for $V$.
\end{Def}
One can also define, as in \cite[Section 3.1]{hernandez2024shifted}, highest $\ell$-weight objects associated to more general sequences of Laurent series $\psi=(\psi_i(u))_{i\in I}$ (i.e.~sequences that do not belong to $\mathfrak{r}_{\mu}$). The next result makes however such a definition obsolete for the study of the~category~$\mathcal{O}_{\mu}$.
\begin{Proposition}[{\cite[Lemma 3.11]{hernandez2024shifted}}]\label{prop:lwtsOmu} Fix $V$ in $\mathcal{O}_{\mu}$. %
Then all $\ell$-weights of $V$ lie in $\mathfrak{r}_{\mu}$.
\end{Proposition}
\begin{Rem} The above result is stated in \cite{hernandez2024shifted} for the category $\mathcal{O}_{\mu}^{\text{HZ}}\subseteq \mathcal{O}_{\mu}$ mentioned in Remark \ref{rem:2categoriesO}, but the proof (which in this case can also be adapted from that of \cite[Proposition 3.6]{gautam2016yangians}) also works directly in our more general setting. This comment applies to other results in this section %
and we stop making it to avoid overcomplicating our exposition.
\end{Rem}
Fix a simple object $V$ of $\mathcal{O}_{\mu}$. By standard arguments (see,~e.g.,~\cite[p.9]{pinet2024functor}), the object~$V$ must be of highest $\ell$-weight $\psi$ for some sequence of Laurent series $\psi=(\psi_i(u))_{i\in I}$. Moreover, by Proposition \ref{prop:lwtsOmu}, this sequence $\psi$ must belong to $\mathfrak{r}_{\mu}$. On the other hand, \cite[Theorem 3.12]{hernandez2024shifted} shows that, for every $\psi \in \mathfrak{r}_{\mu}$, there exists, up to isomorphism,  a unique simple object $L(\psi)$ in $\mathcal{O}_{\mu}^{\text{HZ}}\subseteq \mathcal{O}_{\mu}$ with highest $\ell$-weight $\psi$. This leads to the following result.
\begin{Theorem}[{\cite[Theorem 3.12]{hernandez2024shifted}}]\label{thm:simplesOmu} The set $\{L(\psi)\}_{\psi\in \mathfrak{r}_{\mu}}$ gives a complete collection of mutually non-isomorphic simple objects for $\mathcal{O}_{\mu}$ and $\mathcal{O}_{\mu}^{\text{\normalfont{HZ}}}$.
\end{Theorem}
\begin{Rem}[{\cite[Example 3.3]{hernandez2024shifted}}]\label{rem:prefundMonGen} Fix $\psi,\xi\in \mathfrak{r}$ and let $v$ and $w$ be the respective highest $\ell$-weight vectors of $L(\psi)$ and $L(\xi)$. Then, by \cite[Lemma 2.5]{hernandez2024shifted}, for all $i\in I$,
$$ e_i(u)\cdot(v \otimes w) = 0\ \text{ and }\ h_i(u)\cdot(v\otimes w)=\psi_i(u)\xi_i(u)v\otimes w$$
so that $v\otimes w$ generates a submodule of highest $\ell$-weight $\psi\xi$ in the tensor product $L(\psi)\otimes L(\xi)$. In particular, the simple module $L(\psi\xi)$ is necessarily a composition factor of $L(\psi)\otimes L(\xi)$.
\end{Rem}
We finish this subsection by giving explicit examples of simple modules in $\mathcal{O}_{sh}$. We will need,  for $i\in I$ and $a\in \C$, the \textit{prefundamental $\ell$-weight} $\mathsf{\Psi}_{i,a}\in \mathfrak{r}_{\varpi_i^{\vee}}$ given by
$$ (\mathsf{\Psi}_{i,a})_j(u)=\left\{ \begin{array}{ll} u-a & \text{if }i=j, \\ 1 &\text{else}
\end{array}\right.$$
and the \textit{fundamental $\ell$-weight} $\mathsf{Y}_{i,a}\in \mathfrak{r}_0$ defined by
$$ \mathsf{Y}_{i,a} = \mathsf{\Psi}_{i,a-2}\mathsf{\Psi}_{i,a}^{-1}$$
using the natural group structure on $\mathfrak{r}=\bigsqcup_{\mu\in P^{\vee}}\mathfrak{r}_{\mu}$ (where the operation is given by pointwise multiplication of rational functions).
\begin{Example}[\cite{hernandez2024shifted,zhang2020yangians}]\label{ex:Prefund} Fix $i\in I$ and $a\in \mathbb{C}$. The simple module $L(\mathsf{\Psi}_{i,a})$ of $\mathcal{O}_{\varpi_i^{\vee}}$ is called \textit{positive prefundamental representation}. It has dimension 1 and satisfies
$$ e_j(u)v=f_j(u)v=0\ \text{ and }\ h_j(u)v=(\mathsf{\Psi}_{i,a})_j(u)v$$
for all $j\in I$ and $v\in L(\mathsf{\Psi}_{i,a})$. Similarly, the simple object $L(\mathsf{\Psi}_{i,a}^{-1})$  of $\mathcal{O}_{-\varpi_i^{\vee}}$ is called \textit{negative prefundamental representation}. It is infinite-dimensional and can be realized explicitly, for $\mathfrak{g}=\mathfrak{sl}_2$ (and $i=1$), as the module obtained by endowing the vector space defined over the basis $\{v_n\}_{n\geq 0}$ with the $Y_{-\varpi_1^{\vee}}$-action given by\footnote{Recall that we work with the convention that $\hbar=2$.}
$$  \textstyle e_1(u)v_n = \frac{1-\delta_{n,0}}{u-a+2(n-1)}v_{n-1},\ \  f_1(u)v_n=\frac{2(n+1)}{u-a+2n}v_{n+1}\ \text{\,and\,}\ h_1(u)v_n=\frac{(u-a-2)}{(u-a+2(n-1))(u-a+2n)}v_n$$
for all $n\geq 0$. No such explicit realization for $L(\mathsf{\Psi}_{i,a}^{-1})$ is known in general\footnote{Some realizations, via quantum unipotent rings and Lax matrices, are known under restrictive conditions (see, e.g.,~\cite{bazhanov2011baxter,jang2025unipotent} and the references therein).}.
\end{Example} %
\begin{Rem}\label{rem:PrefundImp} Positive and negative prefundamental representations are of critical importance for the study of the representation theory of shifted and unshifted Yangians. Indeed, they are closely related to the famous Baxter's $Q$-operators of quantum integrable systems %
\cite{frenkel2015baxter,zhang2020yangians}, give rise to cluster variables in various monoidal categorifications of cluster algebras \cite{hernandez2016cluster,geiss2024representations}, and are the very first examples of the \textit{chamber modules} we define in Section \ref{sec:Chambermodules}. Furthermore, using the fact that prefundamental $\ell$-weights clearly generate the group $\mathfrak{r}$, we deduce from Remark \ref{rem:prefundMonGen} that all simple modules of $\mathcal{O}_{sh}$ can be realized as subquotients of tensor products of (positive and negative) prefundamental representations. 
\end{Rem}
\begin{Example}[{\cite[Proposition 2.6]{chari1990yangians}}]\label{ex:KR} Fix again $i\in I$ with $a\in \C$. Also let $k\in \mathbb{Z}_{>0}$. The simple module $$\mathsf{W}_{k,a}^{(i)}=L(\mathsf{Y}_{i,a-2(k-1)}\mathsf{Y}_{i,a-2(k-2)}\dots \mathsf{Y}_{i,a})=L(\mathsf{\Psi}_{i,a-2k}\mathsf{\Psi}_{i,a}^{-1})$$ of $\mathcal{O}_0$
is called a \textit{Kirillov-Reshetikhin module} (or \textit{fundamental representation} if $k=1$). It is finite-dimensional and can be realized explicitly, for $\mathfrak{g}=\mathfrak{sl}_2$ (and $i=1$), using the $Y_0$-action given on the vector space with basis $\{v_n\}_{n=0}^k$  by
$$\textstyle e_1(u)v_n=\frac{2(1-\delta_{n,0})}{u-a+2(n-1)}v_{n-1},\ \ f_1(u)v_n=\frac{2(k-n)(n+1)}{u-a+2n}v_{n+1}$$
and
$$ \textstyle h_1(u)v_n=\frac{(u-a-2)(u-a+2k)}{(u-a+2(n-1))(u-a+2n)}v_n$$
for $0\leq n\leq k$. We refer to \cite{chari2008beyond} and the references therein for a glimpse of the importance of the $\mathsf{W}_{k,a}^{(i)}$'s in the study of (unshifted) Yangians and quantum affine algebras.
\end{Example}

\subsection{$\ell$-characters and block decomposition}\label{subsec:ell_characters_and_blocks}
Consider the group morphism $\wt:\mathfrak{r}\rightarrow %
\mathfrak{h}^*
$ given by taking half of the leading coefficient of the underlying rational Laurent series, i.e.%
$$ \textstyle \wt(\psi) = \frac{1}{2}\sum_{i\in I} \psi_{i,1-\langle\mu,\alpha_i\rangle}\varpi_i\in \mathfrak{h}^*$$
for $\psi=(\psi_j(u))_{j\in I}\in \mathfrak{r}_{\mu}$,  where $ \textstyle \psi_j(u) = u^{\langle\mu,\alpha_j\rangle}+\sum_{r>-\langle\mu,\alpha_j\rangle}\psi_{j,r}u^{-r}$
as in \eqref{eq:lwt}. %
If $i\in I$ and $a\in \C$,
$$\textstyle \wt(\mathsf{\Psi}_{i,a})=-\frac{a}{2}\varpi_i\text{ and }\wt(\mathsf{Y}_{i,a})=\varpi_i.$$
Also, if $V$ lies in $\mathcal{O}_{sh}$, then, for each $\psi\in\mathfrak{r}$, the space $V_{\psi}$ is contained in the finite-dimensional vector space $V_{\wt(\psi)}$. Thus ($\mathcal{O}1$)--($\mathcal{O}2$) imply that, as vector spaces%
$$\textstyle  V\simeq \bigoplus_{\omega\in \mathfrak{h}^*}V_{\omega}\simeq \bigoplus_{\omega\in \mathfrak{h}^*}\bigoplus_{\substack{\psi\in \mathfrak{r};\, \wt(\psi)=\omega}} V_{\psi},$$
i.e.~$V$ decomposes as the direct sum of its $\ell$-weight spaces (which are all finite-dimensional). The idea of $\ell$-characters is to encode this decomposition in a generating series.

Let $\mathcal{E}_{\ell}$ be the set of functions $c:\mathfrak{r}\rightarrow \Z$ satisfying
\begin{enumerate}
\item $\{\psi\in \mathfrak{r}\,|\,\wt(\psi)=\omega \text{ and }c(\psi)\neq 0\}$ is a finite set for all $\omega\in \mathfrak{h}^*$, and
\item there exists $\omega_1,\dots,\omega_N\in \mathfrak{h}^*$ such that, for all $\psi\in \mathfrak{r}$,
$$ c(\psi)\neq 0 \Longrightarrow \wt(\psi)\leq\omega_k \text{ for some }1\leq k\leq N.$$
\end{enumerate}
For every $\psi\in \mathfrak{r}$, define a function $[\psi]\in \mathcal{E}_{\ell}$ via $[\psi](\psi')=\delta_{\psi,\psi'}$. Then, the natural ring structure on $\mathcal{E}_{\ell}$ coming from addition and convolution of functions satisfies $[\psi][\psi'] = [\psi\psi']$.

\begin{Def}[\cite{frenkel1999qcharacters,hernandez2024shifted}]\label{def:ell-chars}
Fix $V$ in $\mathcal{O}_{sh}$. The \textit{$\ell$-character} of $V$ is the element of $\mathcal{E}_{\ell}$ given by
$$\textstyle \chi_{\ell}(V) = \sum_{\psi\in \mathfrak{r}}\dim V_{\psi}[\psi].$$
\end{Def}

\begin{Theorem}\label{thm:K0inj} 
The map from $K_0(\mathcal{O}_{sh})$ to $\mathcal{E}_{\ell}$ given by taking $\ell$-characters of objects is an injective ring morphism. In particular, $K_0(\mathcal{O}_{sh})$ is a commutative (and associative) ring.
\end{Theorem}

Suppose that an object $V$ of $\mathcal{O}_{sh}$ has a unique $\ell$-weight $\psi$ whose weight $\wt(\psi)$ is maximal amongst the weights of $V$, i.e.
$$ \wt(\psi) -\omega \in Q_+$$
for all $\omega\in \mathfrak{h}^*$ such that $V_{\omega}\neq 0$. The \textit{normalized $\ell$-character} of $V$ is defined by $$ \widetilde{\chi}_{\ell}(V)=[\psi^{-1}]\chi_{\ell}(V).$$
Normalized $\ell$-characters of highest $\ell$-weight modules in $\mathcal{O}_{sh}$ admit remarkable expressions. To give details about these expressions, set, for $i\in I$ and $a\in \C$,
\begin{equation}
\textstyle \mathsf{A}_{i,a} = \mathsf{Y}_{i,a}\mathsf{Y}_{i,a+2}\prod_{j\sim i}\mathsf{Y}^{-1}_{j,a+1}%
=\tfrac{\sfPsi_{i,a-2}}{\sfPsi_{i,a+2}}\prod_{j\sim i}\tfrac{\sfPsi_{j,a+1}}{\sfPsi_{j,a-1}}
\in \mathfrak{r}_0.
\end{equation}
Then $\wt(\mathsf{A}_{i,a})=\alpha_i$ %
so $\mathsf{A}_{i,a}$ %
can be thought of
as an analogue of a simple root %
in $\mathfrak{r}$. Using~this analogy and writing $\mathcal{A}_+$ (resp.~$\mathcal{A}$) for the set of monomials (resp. Laurent monomials) in the variables $\{\mathsf{A}_{i,a}\}_{i\in I,a\in \C}$, we define a partial order $\preceq$ on $\mathfrak{r}$, called \textit{Nakajima's order}, via
\begin{equation}\label{eq:NakajimaOrderDef}
\psi\preceq \xi \iff \xi\psi^{-1}\in \mathcal{A}_+.
\end{equation}
(Clearly, $\psi\preceq \xi$ implies $\wt(\xi)\in \wt(\psi)+ Q_+$). We will need the following technical result.
\begin{Lemma} Fix $i,j\in I$. Then
\begin{equation}\label{eq:EqZhang}
\textstyle (u-z-c_{ij})h_i(u)e_j(z) = (u-z+c_{ij})e_j(z)h_i(u)-2c_{ij}e_j(u-c_{ij})h_i(u)
\end{equation}
and, similarly,
\begin{equation}\label{eq:EqZhang2}
\textstyle (u-z+c_{ij})h_i(u)f_j(z) = (u-z-c_{ij})f_j(z)h_i(u)+2c_{ij}f_j(u+c_{ij})h_i(u).
\end{equation}
\end{Lemma}
\begin{proof}
We follow \cite[Section 2.4]{gautam2016yangians}. Let $\Gamma(u) =(u-z)[h_i(u),e_j(z)]+[h_i(u),e_{j,1}]$. Then, %
\begin{align*}
\Gamma(u)&=\sum_{p\in \mathbb{Z}}\sum_{q\geq 1}[h_{i,p},e_{j,q}](u-z)u^{-p}z^{-q}+\sum_{p\in \mathbb{Z}}[h_{i,p},e_{j,1}]u^{-p}\\
&=\sum_{p\in \mathbb{Z}}\sum_{q\geq 1}[h_{i,p+1},e_{j,q}]u^{-p}z^{-q}-\sum_{p\in \mathbb{Z}}\sum_{q\geq 0}[h_{i,p},e_{j,q+1}]u^{-p}z^{-q}+\sum_{p\in \mathbb{Z}}[h_{i,p},e_{j,1}]u^{-p}\\
&=\sum_{p\in \mathbb{Z}}\sum_{q\geq 1}([h_{i,p+1},e_{j,q}]-[h_{i,p},e_{j,q+1}])u^{-p}z^{-q}=c_{ij}\sum_{p\in \mathbb{Z}}\sum_{q\geq 1} (h_{i,p}e_{j,q}+e_{j,q}h_{i,p})u^{-p}z^{-q}\\&=c_{ij}(h_i(u)e_j(z)+e_j(z)h_i(u)),
\end{align*}
where we used the defining relation \eqref{H,E} of $Y_{\mu}$. Hence, %
\begin{equation*}
(u-z-c_{ij})h_i(u)e_j(z)-(u-z+c_{ij})e_j(z)h_i(u) = -[h_i(u),e_{j,1}]
\end{equation*}
with the right-hand side independent of $z$. Now, by evaluating this equality at $z=u-c_{ij}$, we get
$[h_i(u),e_{j,1}]=2c_{i,j}e_j(u-c_{i,j})h_i(u)$
and it thus follows from the definition of $\Gamma(u)$~that
$$\Gamma(u)=(u-z)[h_i(u),e_j(z)]+2c_{i,j}e_j(u-c_{i,j})h_i(u)
=c_{ij}(h_i(u)e_j(z)+e_j(z)h_i(u)),$$
which gives \eqref{eq:EqZhang} after simplification. The proof of \eqref{eq:EqZhang2} is similar.
\end{proof}
Clearly, if $\omega\in \mathfrak{h}^*$ and $V$ is in $\mathcal{O}_{sh}$, then \eqref{H,E} implies that, for $i\in I$ and $q\in \Z_{>0}$,
$$ e_{i,q}V_{\omega} \subseteq V_{\omega+\alpha_i} \ \text{ and }\ f_{i,q}V_{\omega}\subseteq V_{\omega-
\alpha_i}.$$
The next theorem (which generalizes \cite[Proposition 3.8]{mukhin2014affinization} and \cite[Proposition 6]{zhang2020yangians}) gives similar containments for $\ell$-weight spaces using \eqref{eq:EqZhang}--\eqref{eq:EqZhang2}.
\begin{Theorem}\label{thm:MY} Fix a coweight $\mu$ of $\mathfrak{g}$ and $V$ in $\mathcal{O}_{\mu}$. Let $\psi\in \mathfrak{r}$, $j\in I$ and $n\in \mathbb{Z}_{>0}$. Then
\begin{equation}\label{eq:MY}
\textstyle e_{j,n}V_{\psi}\subseteq \bigoplus_{i\in I,a\in \mathbb{C}} V_{\psi \mathsf{A}_{i,a}}\text{ and }f_{j,n}V_{\psi}\subseteq\bigoplus_{i\in I,a\in \mathbb{C}} V_{\psi \mathsf{A}_{i,a}^{-1}}.
\end{equation}
\end{Theorem}
\begin{proof}
We follow the proof of \cite[Proposition 6]{zhang2020yangians}. Suppose $v\in V_{\psi}$ is such that $e_{j,n}v\neq 0$. Then $e_{j,n}v$ has a non-zero component in $V_{\xi}$ for some $\ell$-weight $\xi\in\mathfrak{r}$. Fix bases $\{v_k\}_{k=1}^r\subseteq V_{\psi}$ and $\{w_{\ell}\}_{\ell=1}^s\subseteq V_{\xi}$. Without loss of generality, there exist series
$$\{a_{i,k,k'}(u)\}_{i\in I, 1\leq k'<k\leq r}\text{ and }\{b_{i,\ell',\ell}(u)\}_{i\in I,1\leq\ell<\ell'\leq s}\text{ in }u^{\langle\mu,\alpha_i\rangle-1}\mathbb{C}[[u^{-1}]]$$
such that, for all $i$, $k$ and $\ell$,
\begin{enumerate}
\item $(h_i(u)-\psi_i(u))v_k=\sum_{1\leq k'<k} a_{i,k,k'}(u)v_{k'}$ and
\item $(h_i(u)-\xi_i(u))w_{\ell}=\sum_{\ell<\ell'\leq s} b_{i,\ell',\ell}(u)w_{\ell'}.$
\end{enumerate}
Also, there must exist $1\leq k\leq r$ such that $e_{j}(z)v_k$ has a non-zero component in $V_{\xi}[[z^{-1}]]$. Fix a minimal choice $K$ of such $k$ so that
$$ \text{pr}_{\xi}(e_{j}(z)v_K)\neq 0 \text{ and } \text{pr}_{\xi}(e_{j}(z)v_k) =0 \text{ for all }1\leq k<K,$$
with $\text{pr}_{\xi}:V\rightarrow V_{\xi}$ the canonical (vector space) projection. %
Then
\begin{equation}\label{eq:MY1}
\textstyle \text{pr}_{\xi}(e_{j}(z)v_K) = \sum_{\ell=1}^s \Lambda_{\ell}(z)w_{\ell}
\end{equation}
for some $\{\Lambda_{\ell}(z)\}_{\ell=1}^s$ in $\mathbb{C}[[z^{-1}]]$. \smallskip\par Let $L=\min\{1\leq\ell\leq s\,|\,\Lambda_{\ell}(z)\neq 0\}$ and set $\Lambda(z)=\Lambda_L(z)$. Observe that, for $i\in I$,
\begin{align*}
\text{pr}_{\xi}(e_j(z)h_i(u)v_K)&=\textstyle \psi_i(u)\text{pr}_{\xi}(e_j(z)v_K)+\sum_{1\leq k<K}a_{i,K,k}(u)\text{pr}_{\xi}(e_j(z)v_{k})\\
&= \textstyle \psi_i(u)\Lambda(z)w_L+\psi_i(u)\sum_{L<\ell\leq s}\Lambda_{\ell}(z)w_{\ell}\\
&\in \psi_i(u)\Lambda(z)w_L+u^{\langle \mu,\alpha_i\rangle }\mathsf{Z}[[z^{-1},u^{-1}]]
\end{align*}
where $\mathsf{Z}\subseteq V_{\xi}$ is the subspace generated by $\{w_{\ell}\}_{L<\ell\leq s}.$ \medskip\par
Apply %
\eqref{eq:EqZhang} on $v_K$ and use $\text{pr}_{\xi}$. By %
\eqref{eq:MY1} (and as $\text{pr}_{\xi}$ commutes with the action of $h_i(u)$), the left-hand side becomes $(u-z-c_{ij})\xi_i(u)\Lambda(z)w_L $ modulo $u^{\langle \mu,\alpha_i\rangle }\mathsf{Z}[[z^{-1},u^{-1}]]$. In contrast, the right-hand side gives, by the above computations,
$$ ((u-z+c_{ij})\Lambda(z) -2c_{ij}\Lambda(u-c_{ij}))\psi_i(u)w_L,$$
again modulo $u^{\langle \mu,\alpha_i\rangle }\mathsf{Z}[[z^{-1},u^{-1}]]$. Thus, taking the coefficient of $w_L$,
\begin{equation}\label{eq:ZhangEqA}
(u-z-c_{ij})\xi_i(u)\Lambda(z) = ((u-z+c_{ij})\Lambda(z) -2c_{ij}\Lambda(u-c_{ij}))\psi_i(u).
\end{equation}
Expand $\Lambda(z)=\sum_{k\geq 1}\lambda_k z^{-k}$ and let $m=\min\{k\geq 1\,|\,\lambda_k\neq 0\}$. Let us take the coefficient of $z^{-m}$ in \eqref{eq:ZhangEqA}. The second term of the right-hand side does not contribute and we get
$$ ((u-c_{ij})\lambda_m-\lambda_{m+1})\xi_i(u)=((u+c_{ij})\lambda_m-\lambda_{m+1})\psi_i(u)$$
so that
$$ \xi_i(u) = \psi_i(u)\frac{u+c_{ij}-\lambda_{m+1}\lambda_m^{-1}}{u-c_{ij}-\lambda_{m+1}\lambda_m^{-1}} = (\psi \mathsf{A}_{j,\lambda_{m+1}\lambda_m^{-1}})_i(u).$$
Therefore $\xi=\psi \mathsf{A}_{j,\lambda_{m+1}\lambda_m^{-1}}$ and the inclusion $e_{j,n}V_{\psi}\subseteq \bigoplus_{i\in I,a\in \mathbb{C}} V_{\psi \mathsf{A}_{i,a}}$ follows. The second inclusion is proven in an analogous fashion using the relation \eqref{eq:EqZhang2} instead of \eqref{eq:EqZhang}.
\end{proof}
In this paper as in \cite{frenkel2015baxter}, we omit the brackets $[\cdot]$ when writing a Laurent monomial in the $\mathsf{Y}$'s (or the $\mathsf{A}$'s) in the (normalized) $\ell$-character of an object of $\mathcal{O}_{sh}$. A justification for this convention comes from the following easy corollary of Theorem \ref{thm:MY} and Definition~\ref{def:highestlwtmod}, which gives the mentioned expressions for normalized $\ell$-characters of simple objects of $\mathcal{O}_{sh}$. (Note that a similar result was obtained in \cite[Proposition 5.8]{hernandez2024shifted} by a different method.)
\begin{Corollary}\label{cor:Aqcarnorm}  Fix a highest $\ell$-weight module $V$ in $\mathcal{O}_{sh}$. Then $\widetilde{\chi}_{\ell}(V)\in\mathbb{Z}_{\geq 0}[[\mathsf{A}_{i,a}^{-1}]]_{i\in I,a\in \C}$.
\end{Corollary}
We illustrate the above corollary and our convention in an example.
\begin{Example}\label{ex:lchar} Fix $i\in I$ and $a\in \C$. Then, by Example \ref{ex:Prefund}, $\chi_{\ell}(L(\mathsf{\Psi}_{i,a}))=[\mathsf{\Psi}_{i,a}]$. Also, if $\mathfrak{g}=\mathfrak{sl}_2$, Example \ref{ex:KR} gives, for $k\in \Z_{>0}$,
$$\textstyle \chi_{\ell}(\mathsf{W}_{k,a}^{(1)})=\mathsf{Y}_{1,a-2(k-1)}\dots \mathsf{Y}_{1,a}(1+\sum_{r=0}^{k-1} \mathsf{A}_{1,a}^{-1}\dots \mathsf{A}_{1,a-2r}^{-1}).$$
Hence, using Example \ref{ex:Prefund} again, we get that, in $\Z[[\mathsf{A}_{1,a}^{-1}]]_{a\in \C}$,
$$\textstyle \lim_{k\rightarrow\infty}\widetilde{\chi}_{\ell}(\mathsf{W}_{k,a}^{(1)})=1+\sum_{r\geq 0} \mathsf{A}_{1,a}^{-1}\dots \mathsf{A}_{1,a-2r}^{-1}=\widetilde{\chi}_{\ell}(L(\mathsf{\Psi}_{1,a}^{-1})).$$
This is not a coincidence as the next proposition shows.
\end{Example}
\begin{Proposition}[{\cite{hernandez2012asymptotic,zhang2020yangians}}]\label{prop:Prefundlim} For $i\in I$, $a\in \C$ and $\mathfrak{g}$ arbitrary, %
$$ \textstyle \smash{\lim_{k\rightarrow\infty}\widetilde{\chi}_{\ell}(\mathsf{W}_{k,a}^{(i)}) = \widetilde{\chi}_{\ell}(L(\mathsf{\Psi}_{i,a}^{-1}))}$$
as formal series in $\smash{\Z[[\mathsf{A}_{i,a}^{-1}]]_{i\in I,a\in \C}}$.
\end{Proposition}
In general, $\ell$-characters of Kirillov-Reshetikhin and prefundamental modules are difficult to compute. However, it so happens that the $\ell$-weights appearing in these $\ell$-characters can always be described in terms of the crystal $\mathcal{B}$ of Section \ref{sec:Crystal} using Remark \ref{rem:prefundMonGen} and the~result below (which is due to \cite{nakajima2003tanalogues} in the context of quantum affine algebras and can be applied to representations of Yangians because of the work of \cite{gautam2016yangians}). \medskip\par
Recall that the set $\cB_{\C}$ of Section \ref{sec:nonInt} contains all the (possibly non-integral) fundamental crystals $\cB(\varpi_i^{\vee},a)$ of Remark \ref{rem:nonIntCryst}.

\begin{Theorem}[{\cite[Theorem 3.3]{nakajima2003tanalogues}}]\label{thm:NakCrys} Choose $i\in I$ and $a\in \C$. The set of $\ell$-weights of the %
module $L(\mathsf{Y}_{i,-a})$ (omitting multiplicities) is the image of the  fundamental crystal  $\cB(\varpi_i^{\vee},a)\subseteq\mathcal{B}_{\mathbb{C}}$ under the group morphism $\mathcal{B}_{\C}\rightarrow \mathfrak{r}_0$ given by $y_{j,b}\mapsto \mathsf{Y}_{j,-b}$%
.
\end{Theorem}
\begin{Example} For $\mathfrak{g}=\mathfrak{sl}_4$, the fundamental crystal $\cB(\varpi_2^{\vee},0)$ is, as a set,
$$\textstyle \cB(\varpi_2^{\vee},0) = \{y_{2,0},\frac{y_{1,-1}y_{3,-1}}{y_{2,-2}},\frac{y_{1,-1}}{y_{3,-3}},\frac{y_{3,-1}}{y_{1,-3}},\frac{y_{2,-2}}{y_{1,-3}y_{3,-3}},\frac{1}{y_{2,-4}}\},$$
and the $\ell$-character of the fundamental representation $L(\mathsf{Y}_{2,0})$ is
$$\smash{ \mathsf{Y}_{2,0}+ \mathsf{Y}_{1,1} \mathsf{Y}_{3,1} \mathsf{Y}_{2,2}^{-1}+ \mathsf{Y}_{1,1} \mathsf{Y}_{3,3}^{-1}+ \mathsf{Y}_{3,3} \mathsf{Y}_{1,3}^{-1}+ \mathsf{Y}_{2,2} \mathsf{Y}_{1,3}^{-1} \mathsf{Y}_{3,3}^{-1}+ \mathsf{Y}_{2,4}^{-1}.}$$
\end{Example}\vspace*{-0.91mm}
We end this subsection by giving a block decomposition for (the finite-length subcategory of) $\mathcal{O}_{sh}$. While new in this generality, the decomposition we give is reminiscent of the one obtained using \textit{elliptic characters} in \cite{etingof2003elliptic,chari2005characters} for the category $\mathscr{C}$ of finite-dimensional type I  representations of unshifted quantum affine algebras. This decomposition is moreover closely related to the map $\awt:\cB\rightarrow\cB/\Gamma$ of Section \ref{sec:awt} and the truncations of Section \ref{section: tsy} (see Section \ref{sec:CorTSC} for this last point). We start by showing the result below.%
\begin{Theorem}\label{thm:Ext1Block} Fix a coweight $\mu$ of $\mathfrak{g}$ with $\psi,\xi\in \mathfrak{r}_{\mu}$. Suppose that there exists a non-trivial extension %
of $L(\psi)$ by $L(\xi)$. Then either $\psi\preceq \xi$ or $\xi\preceq \psi$.
\end{Theorem}
\begin{proof}
Fix a non-trivial extension $M$ of $L(\psi)$ by $L(\xi)$. Then, if $\psi\not\preceq \xi$, Corollary \ref{cor:Aqcarnorm} gives $\dim L(\xi)_{\psi m} = 0$ for $m\in \mathcal{A}_+$. Thus, as $\chi_{\ell}(M)=\chi_{\ell}(L(\psi))+\chi_{\ell}(L(\xi))$, by this same result
\begin{equation}\label{eq:NontrivialExt1}
 \smash{\dim M_{\psi}=\dim L(\psi)_{\psi}=1 \text{ and } \dim M_{\psi \mathsf{A}_{i,a}}=\dim L(\psi)_{\psi \mathsf{A}_{i,a}}=0}
\end{equation}
for all $i\in I$ and $a\in \mathbb{C}$. Choose $v\in M_{\psi}$. By \eqref{eq:NontrivialExt1} and Theorem \ref{thm:MY},
$h_j(u)v=\psi_j(u)v$ and $e_j(u)v=0$ for every $j\in I$. Also, $v$ does not belong to the submodule $L\subseteq M$ isomorphic to $L(\xi)$ and it follows (from non-triviality of the extension $M$) that $M=Y_{\mu}\cdot v$ is a module~of highest $\ell$-weight $\psi$. In particular, $\xi\preceq \psi$ by Corollary \ref{cor:Aqcarnorm}. This ends the proof.
\end{proof}\vspace*{-0.91mm}
We can now define the aforementioned ``blocks'' of $\mathcal{O}_{sh}$. Recall that $\mathcal{A}$ denotes the group of Laurent monomials in the $\mathsf{A}$'s and fix $\mu\in P^{\vee}$. Let also $\pi_{\mu}:\mathfrak{r}_{\mu}\twoheadrightarrow (\mathfrak{r}/\mathcal{A})_{\mu}$ be the canonical projection and take $\tau\in (\mathfrak{r}/\mathcal{A})_{\mu}$. We let $_{\tau}\mathcal{O}_{\mu}$ be the Serre subcategory of $\mathcal{O}_{\mu}$ generated by the simple objects $L(\psi)$ with highest $\ell$-weight $\psi\in \mathfrak{r}_{\mu}$ satisfying $\pi_{\mu}(\psi)=\tau$. \medskip\par
We call $_{\tau}\mathcal{O}_{\mu}$ a ``block'' of $\mathcal{O}_{\mu}$ because of the following result\footnote{We do not show here that the category ${}_{\tau}\mathcal{O}_{sh}$ cannot be decomposed further (as in \cite[Theorem~7.2]{chari2005characters}), but still refer to this category as a ``block'' by slight abuse of terminology. 
}.
\begin{Theorem}\label{thm:BlockDec} Fix $\tau_1,\tau_2\in (\mathfrak{r}/\mathcal{A})_{\mu}$. Let $V_1 \in \phantom{}_{\tau_1}\mathcal{O}_{\mu}$ and $V_2 \in \phantom{}_{\tau_2}\mathcal{O}_{\mu}$ be of finite length. If $\tau_1\neq \tau_2$, then
$$ \smash{\Ext^1_{\mathcal{O}_{\mu}}(V_1,V_2)} = 0.$$
In particular, the full subcategory of finite-length objects in $\mathcal{O}_{\mu}$ decomposes as $\bigoplus_{\tau\in (\mathfrak{r}/\mathcal{A})_{\mu}} {}_{\tau}\mathcal{O}_{\mu}$ (as abelian categories). Summing over all $\mu$'s gives a similar decomposition for $\mathcal{O}_{sh}$.
\end{Theorem}
\begin{proof}
The first statement follows from Theorem \ref{thm:Ext1Block} by the usual argument involving long exact sequences and double induction (on the lengths of $V_1$ and $V_2$). The second statement is clear (as any finite-length indecomposable object of $\mathcal{O}_{\mu}$ must lie in $_{\tau}\mathcal{O}_{\mu}$ for some $\tau$, again by Theorem \ref{thm:Ext1Block}).
\end{proof}
\begin{Rem}\label{rem:Blockawt} The decomposition given in Theorem \ref{thm:BlockDec} is intrinsically related to~the~map $\awt:\cB\twoheadrightarrow \cB/\Gamma$ of Section \ref{sec:awt}. Indeed, identify $\mathfrak{r}$ with the group $\cB_{\C}$ of Section \ref{sec:nonInt} %
using~the isomorphism $y_{i,a}\mapsto \mathsf{\Psi}_{i,a}$%
. Then the subgroup $\mathcal{A}\subseteq \mathfrak{r}$ corresponds to the weight-0 component $(\Gamma_{\C})_0$ of $\Gamma_{\C}\subseteq \cB_{\C}$ and the map
$$\pi_{\mu}:\mathfrak{r}_{\mu}\rightarrow \mathfrak{r}_{\mu}/\mathcal{A}%
,$$
underlying the decomposition of Theorem \ref{thm:BlockDec}, is naturally identified with the map
$$(\cB_{\C})_{\mu}\rightarrow (\cB_{\C})_{\mu}/(\Gamma_{\C})_{0}$$
obtained from $\awt$ (or its natural extension to $\cB_{\C}$) by restricting to the weight-$\mu$ component. One can hence think of Theorem \ref{thm:BlockDec} as giving a decomposition with respect to $\awt$.
\end{Rem}
\subsection{Truncations and crystals}\label{sec:GTweight}
Take $\lambda\in P_+^{\vee}$ and $\mathbf{R}\in \C^{\lambda}$. Take also $\mu\in P^{\vee}$ with $\mu\leq \lambda$ and write $\textstyle \la-\mu=\sum_{i\in I}m_i\alpha_i^{\vee}$ for $(m_i)_{i\in I}\in \mathbb{N}^I$. Recall the projection $\Phi_{\mu}^{\la}(\bR):Y_{\mu}\twoheadrightarrow Y_{\mu}^{\lambda}(\bR)$ of Theorem \ref{GKLO homomorphism} and define $\mathcal{O}_{\mu}^{\lambda}(\bR)$ as the image of $\mathcal{O}_{\mu}$ in the category of $Y_{\mu}^{\lambda}(\bR)$-modules,~
i.e.

\begin{Def}\label{def:OmuLambda} The category $\mathcal{O}_{\mu}^{\lambda}(\bR)$ is the full subcategory of $Y_{\mu}^{\lambda}(\bR)$-Mod consisting~of all the objects $V$ for which the pullback of $V$ by the map $\Phi_{\mu}^{\la}(\bR):Y_{\mu}\twoheadrightarrow Y_{\mu}^{\lambda}(\bR)$ lies in %
$\mathcal{O}_{\mu}$.
\end{Def}

The categories $\{\mathcal{O}_{\mu}^{\lambda}(\bR)\}_{\la,\mu,\bR}$ %
 satisfy properties similar to those satisfied by blocks of the usual BGG category $\mathcal{O}$ in Lie theory (which is not surprising as clever choices of $\mu,\la,\bR$~give categories equivalent to these blocks by \cite{webster2020quantum}%
 ). In particular, by results of \cite{kamnitzer2019category}, the above categories $\mathcal{O}_{\mu}^{\la}(\bR)$ have finitely-many simple objects and are in fact all equivalent to categories of finite-dimensional modules over finite-dimensional algebras. We give details about these algebras, called \textit{parity KLRW algebras}, in Section \ref{sec:KLR} (for integral $\bR$).\medskip\par

The following theorem, which was conjectured in \cite{kamnitzer2019highest} and proven in \cite{kamnitzer2019category}, characterizes which simple objects of $\mathcal{O}_{\mu}$ descend\footnote{Here ``$V$ descends to $Y_{\mu}^{\lambda}(\bR)$'' means that $(\Ker \Phi_{\mu}^{\lambda}(\bR))\cdot V=0$.} to $\mathcal{O}_{\mu}^{\lambda}(\bR)$ (for fixed parameters $\lambda$ and $\bR$). Recall the identification $\cB_{\C}\simeq \mathfrak{r}$ of Remark \ref{rem:Blockawt}. This induces an embedding of $\cB(\la,\bR)$ in $\mathfrak{r}$
even if $\bR$ is not integral (see Remark \ref{rem:nonIntCryst}).
\begin{Theorem}[{\cite[Corollary 5.22]{kamnitzer2019category}}] \label{th:descend}
Take $\psi\in \mathfrak{r}_{\mu}$. Then $L(\psi)$ descends to $Y_{\mu}^\lambda(\bR)$ if and only if $\psi$ belongs to $\cB(\lambda,\mathbf{R})_{\mu}$ (via the identification above).
\end{Theorem}
In particular, as every Laurent monomial in $\cB_{\C}$ lies in one of the infinitely-many (possibly non-integral) product monomial crystals $\cB(\la,\bR)$, we get%
:
\begin{Corollary}\label{coro:every_object_descends_to_some_trunc}
Let $V$ be a simple object in $\mathcal{O}_{\mu}$. Then there exists %
$\lambda\in P^{\vee}_+$ and $\bR\in \C^{\lambda}$~for which $V$ descends to %
$Y_{\mu}^{\lambda}(\bR)$ (and hence naturally lies in $\mathcal{O}_{\mu}^{\lambda}(\bR)$).
\end{Corollary}
\begin{Rem} The above corollary and the fact that only finitely-many simple modules in $\mathcal{O}_{sh}$ can descend to a truncation have been extended to shifted Yangians associated~to~non simply-laced simple Lie algebras in \cite[Theorem 8.4 and Theorem 9.3]{hernandez2024shifted}. Observe however that these generalizations use the alternative definition of truncation given in Remark \ref{rem:tsyideal}. For a generalization using our definition of truncation, see the recent work \cite{varagnolo2025representations}.
\end{Rem}
Fix $i\in I$ and $a\in \C$. Then Lemma \ref{th:lowest} gives
\begin{equation}\label{eq:monomial_one_in_crystal}
1=y_{i,a}y_{i,a}^{-1}\in \cB(\varpi_i^{\vee}+\varpi_{i^*}^{\vee},\bR)
\end{equation}
with $\bR=(R_j)_{j\in I}$ containing only the scalars $a$ and $a+h$ in $R_i$ and $R_{i^*}$ (resp.).~In~particular,
a given $\psi\in \mathfrak{r}$ always lies in infinitely-many distinct product monomial crystals and Theorem \ref{th:descend} implies that simple modules in $\mathcal{O}_{sh}$ all lie in infinitely-many distinct categories of the form $\mathcal{O}_{\mu}^{\la}(\bR)$. We will talk about this ``multiplicity phenomenon'' again in Section \ref{sec:TSC} after we tie truncations to the block decomposition of the last subsection. In the meantime,~recall that the results of \cite{kamnitzer2019highest} also allow the computation of the Gelfand--Kirillov dimension of the simple objects in $\mathcal{O}_{sh}$ via crystal combinatorics\footnote{For details on $\GKdim$, see \cite{krause2000growth}.}. Indeed, as $\cB$ is a normal crystal\footnote{We write everything below for monomials in $\cB$, but the result appearing in \cite{kamnitzer2022lie} can be generalized to $\cB_{\C}$ using crystals for a Lie algebra of the form $\fg^{\oplus k}$ as in Remark \ref{rem:nonIntCryst}. We avoid such technicalities here.}, every monomial $y\in \cB$ lies in a unique connected component of the form $B(\la')$ with $\la'\in P_+^{\vee}$.

\begin{Lemma}[{\cite[Proposition 9.18]{kamnitzer2022lie}}]\label{lem:GKdimSimples} Take $\la'\in P_+^{\vee}$ and $\mu \in P^{\vee}$ such that $\mu\leq \la'$. Fix also $\psi\in \cB_{\mu}\subseteq \mathfrak{r}_{\mu}$ lying in a connected component isomorphic to $B(\la')$. Then %
$$ \GKdim L(\psi) = \htt(\la'-\mu),$$
where $\htt:Q_+\rightarrow \mathbb{N}$ is the usual height map.%
\end{Lemma}

\subsection{Chamber modules}\label{sec:Chambermodules} Fix $i\in I$, $a\in \C$ and $w$ in the Weyl group $W$ of $\fg$. Then, there is one element in the $w\varpi_i^{\vee}$-weight space of $\cB(\varpi_i^{\vee},a)\simeq \cB(\varpi_i^{\vee})$. Hence, by Theorem~\ref{th:descend},~the category $\mathcal{O}_{\mu}^{\la}(a)$ (for $\la=\varpi_i^{\vee}$ and $\mu=w\varpi_i^{\vee}$) contains, up to isomorphism, a unique simple object, denoted by $L_{w\varpi_i^{\vee},a}$ and called the \textit{$w\varpi_i^{\vee}$-chamber module of spectral parameter $a$}.
\begin{Example}\label{ex:Inflsl3} Special cases %
of chamber modules are positive prefundamental representations $L_{\varpi_i^{\vee},a} = L(\mathsf{\Psi}_{i,a})$ %
and negative prefundamental representations $L_{w_0\varpi_i^{\vee},a}=L(\mathsf{\Psi}_{i^*,a-h}^{-1})$ (see Lemma \ref{th:lowest} for the notation used here). Chamber modules form however a much~larger family than prefundamental modules. For instance, the chamber module
$$\textstyle L_{s_i\varpi_i^{\vee},a}=L(\mathsf{\Psi}^{-1}_{i,a-2}\prod_{j\sim i}\mathsf{\Psi}_{j,a-1})$$
can be realized\footnote{A similar module was defined in the context of shifted quantum affine algebras in \cite[Example~5.2]{hernandez2023representations}.} (for any $\fg$,~here simply-laced) on the vector space with basis $\{v_n\}_{n\geq 0}$ via
$$\textstyle e_{i}(u)v_n=\frac{1-\delta_{n,0}}{u-a+2n}v_{n-1},\ \ f_i(u)v_n=\frac{2(n+1)}{u-a+2(n+1)}v_{n+1}, \ \ h_i(u)v_n=\frac{u-a}{(u-a+2n)(u-a+2(n+1))}v_n$$
with %
$ e_j(u)v_n=f_j(u)v_n=0$ for $j\neq i$ and
$$ h_j(u)v_n=\left\{
\begin{array}{ll}
(u-a+2n+1)v_n & \text{if }i\sim j,\\
0 & \text{if }i\neq j\text{ and }i\not\sim j.
\end{array}
\right.$$
Note that this module
can be thought of as an example of \textit{$\{i\}$-inflation~to~$\fg$}~as in \cite{pinetinflations}. %
\end{Example} %

The highest $\ell$-weight $\mathsf{\Psi}_{w\varpi_i^{\vee},a}$ of the chamber module $L_{w\varpi_i^{\vee},a}$ can be expressed using the braid group action on $\mathfrak{r}$ considered in \cite[Section 4.2]{friesen2025braid}. Indeed, denote by $\mathcal{Y}\subseteq \mathfrak{r}_0$ the group of Laurent monomials in the variables $\{\mathsf{Y}_{i,a}\}_{i\in I,a\in\C}$. Consider the group isomorphism
$$\vartheta: \mathfrak{r}\rightarrow \mathcal{Y}$$
given by $\vartheta(\mathsf{\Psi}_{i,a})=\mathsf{Y}_{i,-a}$. Then Theorem \ref{thm:NakCrys} and Theorem \ref{th:descend} imply that $\vartheta$ induces, for all $(i,a)\in I\times_2\Z$, a bijection between the highest $\ell$-weights of the category
$$\textstyle \mathcal{O}_{sh}^{\varpi_i^{\vee}}(a)=\bigoplus_{\mu \in P^{\vee}} \mathcal{O}_{\mu}^{\varpi_i^{\vee}}(a)$$
(containing the chamber modules $L_{w\varpi_i^{\vee},a}$ for $w\in W$) and the $\ell$-weights of the fundamental module $L(\mathsf{Y}_{i,-a})$ of $Y$. Now, take $w\in W$, $i\in I$ and $a\in \C$. Then, by \cite[Proposition 4.8]{friesen2025braid}, the unique $\ell$-weight $\mathsf{Y}_w$ of $L(\mathsf{Y}_{i,-a})$ with weight $\wt(\mathsf{Y}_w)=w\varpi_i$ is
\begin{equation}\label{eq:YwFWW}
\mathsf{Y}_w=T_w(\mathsf{Y}_{i,-a})
\end{equation}
where $T_w=\tau_w^{\mathscr{M}}$ is the braid group operator associated to $w$ in \cite{friesen2025braid}. It is furthermore easy to see that\footnote{Note that $T_w$ gives an automorphism of $\mathcal{Y}$ for all $w\in W$ by \eqref{eq:YwFWW}. Hence $\vartheta\circ T_{s_i}^{-1}:\mathcal{Y}\to\mathfrak{r}$ is well-defined.}
 $T_{s_i} \circ \vartheta = \vartheta \circ T_{s_i}^{-1}$ for all $i\in I$ so that
$$T_w \circ \vartheta = T_{s_{i_1}}\circ \dots \circ T_{s_{i_{\ell}}}\circ\vartheta = \vartheta\circ T_{s_{i_1}}^{-1}\circ \dots \circ T_{s_{i_{\ell}}}^{-1} = \vartheta \circ T_{w^{-1}}^{-1}$$
if $w=s_{i_1}\dots s_{i_\ell}$ is a reduced expression for $w$. Hence, by the above bijection, the extremal $\ell$-weight $\mathsf{Y}_{w}$ of $L(\mathsf{Y}_{i,-a})$ corresponds to the highest $\ell$-weight
$$ \mathsf{\Psi}_{w\varpi_i^{\vee},a}= \vartheta^{-1}(\mathsf{Y}_w) = \vartheta^{-1}\circ T_w(\mathsf{Y}_{i,-a}) = \vartheta^{-1}\circ T_w\circ \vartheta(\mathsf{\Psi}_{i,a})=T_{w^{-1}}^{-1}(\mathsf{\Psi}_{i,a})$$
of the chamber module $L_{w\varpi_i^{\vee},a}$.
\medskip\par
We record the above discussion in a proper lemma%
. We also include in this lemma~an~easy consequence of the results of \cite{varagnolo2025representations,hernandez2026borel} with the (clear) compatibility between the braid group actions~given~in %
\cite[(3.9)--(3.10)]{frenkel2024extended} and %
\cite[Corollary 4.5]{friesen2025braid}%
.
\begin{Lemma}\label{lem:highest_ell_weight_chamber_module_braid}
The highest $\ell$-weight of the module $L_{w\varpi_i^{\vee},a}$ is
$$\mathsf{\Psi}_{w\varpi_i^{\vee},a}=T_{w^{-1}}^{-1}(\mathsf{\Psi}_{i,a}),$$ where $T_{w^{-1}}:\mathfrak{r}\rightarrow\mathfrak{r}$ is the braid group operator associated to $w^{-1}$ in \cite{friesen2025braid}. In particular, the chamber modules $\{L_{w\varpi_i^{\vee},a}\}_{i\in I,a\in \C,w\in W}$ correspond precisely via \cite[Corollary 1.2.1]{varagnolo2025representations} and \cite[Theorem 4.2]{hernandez2026borel} to those appearing in \cite[Conjectures 4.8, 5.9, 6.8 and 6.11]{frenkel2024extended}.
\end{Lemma}
\begin{Rem} As stated in Section~\ref{sec:Intro}, we show all of the above conjectures in this paper.
\end{Rem}
One can deduce further facts about chamber modules than simply their highest $\ell$-weights. Indeed, by Lemma \ref{lem:GKdimSimples}, 
$$\smash{\GKdim(L_{w\varpi_i^{\vee},a})=\text{ht}(\varpi_i^{\vee}-w\varpi_i^{\vee}).}$$
Also, as we will show in an upcoming article \cite{otherpaper}, when $w\varpi_i^{\vee}\in P_+^{\vee}$ is \textit{almost dominant} (i.e.~$\langle w\varpi_i^{\vee},\alpha\rangle\geq -1$ for all positive roots $\alpha$ of $\fg$), the normalized $\ell$-character of $L_{w\varpi_i^{\vee},a}$ can~be expressed as a generating function of reverse plane partitions and the module $L_{w\varpi_i^{\vee},a}$ can~be realized as a \textit{limit of Demazure submodules in KR-modules}. This extends results of \cite{hernandez2012asymptotic,zhang2020yangians} for $w\varpi_i^{\vee}=-\varpi_{i^*}^{\vee}$.  %
\subsection{GT-weights}\label{subsec:GT_weights}
For a fixed $\nu=\sum_{i\in I}m_i\alpha_i\in Q_+^\vee$ (with $(m_i)_{i\in I}\in \Z_{\geq 0}^I$), let\footnote{A similar notion is considered in \cite{brundan2006shifted} (for dominantly-shifted Yangians and $\fg=\fgl_n$).}
$$\textstyle \Lambda_\nu=\prod_{i\in I} \C^{m_i}/\Sigma_{m_i}.$$
We refer to $\textstyle\Lambda:=\bigsqcup_{\nu\in Q_+^\vee} \Lambda_{\nu}$
as the set of \textit{GT-weights} (where \textit{GT} stands for \textit{Gelfand--Tsetlin}). Point-wise union of multisets provides a binary operation
\begin{equation*}
\cup:\Lambda_\nu \times \Lambda_{\nu'}\to \Lambda_{\nu+\nu'}
\end{equation*}
which, in turn, endows $\Lambda$ with the structure of an abelian monoid.
This gives a $Q_+^\vee$-graded ring structure to the \textit{GT-character ring}, that is the set of functions
\begin{equation*}
\textstyle \EGT:=\bigoplus_{\nu\in Q_+^\vee}\mathcal{E}_{\nu},\hspace{1em} \mathcal{E}_{\nu}:=\Z^{\Lambda_\nu}
\end{equation*}
where addition is defined point-wise and where product is defined via convolution, i.e.
\begin{equation}\label{eq:product_in_cEGT}
(f\ast g)(\bS)=\sum_{\bS^{(1)} \cup \bS^{(2)}=\bS} f(\bS_1)g(\bS_2).
\end{equation}
This product is well-defined since there are only finitely many ways to write a GT-weight~$\bS$ as the union of two GT-weights $\bS^{(1)}$ and $\bS^{(2)}$. %
\medskip\par
Fix $\lambda\in P_+^\vee$ together with %
$\bR=(R_i)_{i\in I}\in \C^\lambda$ and %
$\mu=\lambda-\nu\in P^\vee$. Recall that \eqref{eq: def of A gens}~defines elements $a_{i,r}\in Y_\mu$, where $i\in I$ and $1\leq r \leq m_i$. We use these elements $a_{i,r}$ to define, given $\bS=(S_i)_{i\in I}\in \Lambda_\nu$ and a $Y_\mu$-module $V$, the ``\textit{$\bS$-weight space of $V$}'' as %
\begin{align*}
W_\bS(V) = \{ v \in V \,|\, \exists p\in \mathbb{N}\text{ with } (a_{i,r} - (-1)^r e_r(S_i))^p v = 0  \text{ for all } i\in I \text{ and } 1\leq r\leq m_i\},
\end{align*}
with $e_r(S_i)$ the $r^{\text{th}}$-elementary symmetric function of the multiset $S_i$. Note that $(-1)^r e_r(S_i)$ is exactly the coefficient of $u^{m_i-r}$ in the polynomial $p_{S_i}(u)=\prod_{c\in S_i}(u-c)$.

\begin{Def}
Fix an object $V$ in $\O_{\mu}$. Then the \textit{GT-character} of $V$ (with respect~to~$\bR$) is the element of $\cE_{\nu}\subseteq \EGT$ given by
\begin{equation*}
\textstyle\chiGT^{\bR}(V)=\sum_{\bS\in \Lambda_\nu} \dim W_{\bS}(V) [\bS]
\end{equation*}
where $[\bS]\in \cE_{\nu}$ is the map given by $[\bS](\bS')=\delta_{\bS,\bS'}$ for all $\bS'\in \Lambda_{\nu}$.
\end{Def}
We can relate GT-characters to the notion of $\ell$-characters introduced previously. Define a map $\Psi_{\mathbf{R}}:\Lambda%
\rightarrow \mathfrak{r}$ via
$$ \textstyle (\Psi_{\mathbf{R}}(\bS))_i(u)=p_{R_i}(u)\frac{ \prod_{j \sim i}  p_{S_j}(u -1)}{p_{S_i}(u) p_{S_i}(u-2)}$$
and observe that %
for all $i\in I$,
\begin{align*}
\langle \mu,\alpha_i\rangle&= \textstyle \langle \lambda-\sum_{j\in I}m_j\alpha_j^{\vee},\alpha_i\rangle = \langle\lambda,\alpha_i\rangle-\sum_{j\in I}c_{ij}m_j\\
&=\textstyle \deg p_{R_i}(u)+\sum_{j\sim i}\deg p_{S_j}(u-1)-\deg p_{S_i}(u)-\deg p_{S_i}(u-2)
\end{align*}
Thus $\Psi_{\bR}(\bS)\in \mathfrak{r}_{\mu}$ for all $\bS\in \Lambda_{\nu}$ and $\Psi_{\bR}$ restricts to a map from $\Lambda_{\nu}$ to $\mathfrak{r}_{\mu}$. The next lemma follows %
from \cite[(2.8)]{gerasimov2005class} (after using the invertibility of the Cartan matrix of $\fg$).
\begin{Lemma}\label{le:injectivityPsibR}
The map $\Psi_{\bR}:\Lambda%
\to \mathfrak{r}%
$ is injective.
\end{Lemma}
Consider the map $[\Psi_{\bR}]:\EGT \to \Z^{\mathfrak{r}}$ given on the topological basis $\{[\bS]\}_{\bS\in \Lambda}\subseteq \EGT$ by
\begin{equation*}\label{eq:pushforward_Psi}
[\bS]\mapsto [\Psi_{\bR}(\bS)].
\end{equation*}
Then $[\Psi_{\bR}]$ is easily seen to be also injective.\newpage
\begin{Proposition}\label{prop:injectivity_chiGT}
If $V$ is an object of the category $\O_\mu^\lambda(\bR)$ of Definition \ref{def:OmuLambda}, then
\begin{equation}\label{eq:chiellchiGT}
\chi_\ell(V)=([\Psi_{\bR}]\circ \chiGT^\bR)(V).
\end{equation}
In particular, after restricting to $\O_{\mu}^{\la}(\bR)$, the GT-character map $\chiGT^\bR:K_0(\O_\mu^\lambda(\bR))\to \EGT$ becomes an injective (group) morphism.
\end{Proposition}

\begin{proof}
Observe first that the two subalgebras of $Y_{\mu}$ generated by $\{a_{i,r}\,|\,i\in I,r\in \Z_{\geq 0}\}$ and $\{h_{i,r}\,|\,i\in I,r>-\langle \mu,\alpha_i\rangle \}$ coincide (again by invertibility of the Cartan matrix of $\fg$). Also, the elements $a_{i,r}$ with $r>m_i$ all lie in the kernel of the map $\Phi_{\mu}^{\la}(\bR)$ of Theorem \ref{GKLO homomorphism}. Thus, given $V$ in $\O_\mu^\lambda(\bR)$ and $\psi\in \mathfrak{r}_{\mu}$ such that $V_\psi\neq 0$, there is necessarily $\bS \in \Lambda_{\nu}$~with~$\Psi_R(\bS)=\psi$ and $V_{\psi}=W_{\bS}(V)$. This shows \eqref{eq:chiellchiGT}. The other statement follows from Theorem \ref{thm:K0inj}.
\end{proof}

\begin{Example}\label{ex:sl2_GTchar_negative_prefund}
For $\fg=\fsl_2$, take $\lambda=\varpi_1^\vee$, $\mu=-\varpi_1^\vee$ and %
$a\in \C$. %
For $n\in \Z_{\geq 0}$, one has
\begin{equation*}
\Psi_{\bR}(\{a-2n\})=\tfrac{(u-a-2)}{(u-(a-2n))(u-2-(a-2n))}%
\end{equation*}
and it hence follows from Example \ref{ex:Prefund} and Proposition \ref{prop:injectivity_chiGT} that
\begin{equation*}
\chiGT^{\{a+2\}}(L(\sfPsi_{1,a}^{-1}))=[\{a\}]+[\{a-2\}]+[\{a-4\}]+\dots
\end{equation*}
\end{Example}

Now, take $\lambda_1,\lambda_2\in P_+^\vee$ with $\bR_1\in \C^{\la_1}$ and $\bR_2\in \C^{\la_2}$. Take also $\mu_1,\mu_2\in P^{\vee}$ such that
$$\nu_1=\la_1-\mu_1\text{ and } \nu_2=\la_2-\mu_2$$
both lie in $Q_+^{\vee}$. Then, clearly, for all $\bS_1\in \Lambda_{\nu_1}$ and $\bS_2\in \Lambda_{\nu_2}$,
$$\Psi_{\bR_1\cup\bR_2}(\bS_1\cup\bS_2)=\Psi_{\bR_1}(\bS_1)\Psi_{\bR_2}(\bS_2)$$
so that the maps $[\Psi_{\bR_1\cup\bR_2}],[\Psi_{\bR_1}],[\Psi_{\bR_2}]:\EGT\to \Z^{\mathfrak{r}}$ satisfy
$$[\Psi_{\bR_1\cup\bR_2}]([\bS_1\cup\bS_2])=[\Psi_{\bR_1\cup\bR_2}(\bS_1\cup\bS_2)]=[\Psi_{\bR_1}(\bS_1)\Psi_{\bR_2}(\bS_2)]=[\Psi_{\bR_1}]([\bS_1])\ast [\Psi_{\bR_2}]([\bS_2])$$
where $\ast$ is convolution in $\cE_{\ell}\subseteq \Z^{\mathfrak{r}}$. This motivates the following proposition.

\begin{Proposition}\label{prop:GTchar_is_multiplicative}
Taking GT-characters is multiplicative, that is, given two objects $V_1$ in $\O_{\mu_1}^{\lambda_1}(\bR_1)$ and $V_2$ in $\O_{\mu_2}^{\lambda_2}(\bR_2)$ (with $\la_1,\la_2,\bR_1,\bR_2,\mu_1,\mu_2$ as above),
\begin{equation*}
\chiGT^{\bR_1\cup\bR_2}(V_1\otimes V_2)=\chiGT^{\bR_1}(V_1)\ast \chiGT^{\bR_2}(V_2)
\end{equation*}
where $\ast$ is convolution in $\EGT$.
\end{Proposition}

\begin{proof}
First, for $\psi\in \fr$, by multiplicativity of $\ell$-characters (see Theorem \ref{thm:K0inj}),
\begin{equation}\label{eq:tensor_prod_ellweight_spaces}
\textstyle \dim (V_1\otimes V_2)_\psi=\sum_{\psi=\psi_1\psi_2} \dim (V_1)_{\psi_1}\cdot \dim (V_2)_{\psi_2}.
\end{equation}
Also, as explained in the proof of Proposition \ref{prop:injectivity_chiGT}, the $\ell$-weight spaces $(V_1)_{\psi_1}$ and $(V_2)_{\psi_2}$~are non-zero only if $\psi_1=\Psi_{\bR_1}(\bS_1)$ and $\psi_2=\Psi_{\bR_2}(\bS_2)$ for some $\bS_1\in \Lambda_{\nu_1}$ and $\bS_2\in \Lambda_{\nu_2}$. Thus, since $V_{\psi_1}=W_{\bS_1}(V_1)$ and $V_{\psi_2}=W_{\bS_2}(V_2)$ in this case, and since $\Psi_{\bR_1\cup\bR_2}(\bS_1\cup\bS_2)=\psi_1\psi_2$ by the above discussion (so that $(V_1\otimes V_2)_{\psi_1\psi_2}=W_{\bS_1\cup\bS_2}(V_1\otimes V_2)$), we have that \eqref{eq:tensor_prod_ellweight_spaces} gives
\begin{equation*}
\textstyle \dim W_{\bS}(V_1\otimes V_2)=\sum_{\bS=\bS_1\cup\bS_2}\dim W_{\bS_1}(V_1)\cdot \dim W_{\bS_2}(V_2).
\end{equation*}
This recovers the product defined in \eqref{eq:product_in_cEGT} and the proposition hence follows.
\end{proof}

\begin{Rem} By the above proof, if $V_1$ and $V_2$ lie in $\O_{\mu_1}^{\la_1}(\bR_1)$ and $\O_{\mu_2}^{\la_2}(\bR_2)$ respectively, then all $\ell$-weights of the $Y_{\mu_1+\mu_2}$-module $V_1\otimes V_2$ are of the form $\Psi_{\bR_1\cup\bR_2}(\bS_1\cup\bS_2)$ for some $\bS_1\in \Lambda_{\nu_1}$ and $\bS_2\in \Lambda_{\nu_2}$. In particular, the elements $\{a_{i,r}\,|\,i\in I,r>\langle \nu_1+\nu_2,\varpi_i\rangle\}$ defined~via the set of parameters $\bR_1\cup\bR_2$ act nilpotently on $V_1\otimes V_2$ and it can hence seem reasonable to expect (because of Remark \ref{rem:tsyideal}) that $V_1\otimes V_2$ descends to the truncation $Y_{\mu_1+\mu_2}^{\la_1+\la_2}(\bR_1\cup\bR_2)$. We will show that this is indeed the case in the upcoming Section \ref{sec:TSC}.
\end{Rem}
We finish this section with a technical, but useful lemma. Let $\mu$, $\la$ and $\bR$ be as in the beginning of this section and suppose that there is a surjection
\begin{equation*}
Y_\mu^{\lambda'}(\bR')\twoheadrightarrow Y_\mu^\lambda(\bR)
\end{equation*}
for some $\lambda'\in P^{\vee}_+$ and $\bR'\in \C^{\la'}$. Then $\mathcal{O}_{\mu}^{\la}(\bR)$ is %
a full subcategory of $\mathcal{O}_{\mu}^{\la'}(\bR')$~%
and the set %
of simple objects of the latter category contains the set %
of simple objects of the former, i.e.
$$\cB(\la,\bR)_{\mu}\subseteq \cB(\la',\bR')_{\mu}.$$
Now, suppose that $\bR$ and $\bR'$ are integral sets of parameters%
. Then Theorem \ref{th:TFAE} shows~that $y_{\bR'}\in y_{\bR}\Gamma_+$, i.e.~there exists an (integral) $\bT\in \Lambda$ such that
$$ y_{\bR'} = y_{\bR}z_{\bT}.$$
Take $V$ in $\mathcal{O}_{\mu}^{\la}(\bR)$. Then the following ensures that the GT-characters $\chiGT^{\bR}(V)$ and $\chiGT^{\bR'}(V)$ agree (up to convolution with the function $[\bT]\in \EGT$, which is independent of $V$).

\begin{Lemma}\label{lem:chiGT_and_inclusions}
With the above notation,
$$\chiGT^{\bR'}(V)=\chiGT^{\bR}(V)\ast [\bT],$$ where $\ast$ is convolution in $\EGT$.
\end{Lemma}

\begin{proof}
Since $y_{\bR'}=y_{\bR}z_{\bT}$, one has, for $i\in I$,
\begin{equation*}
\big(\Psi_{\bR'}(\bT)\big)_i(u)=p_{R'_i}(u)\dfrac{\prod_{j\sim i} p_{T_j}(u-1)}{p_{T_i}(u)p_{T_i}(u-2)}=p_{R_i}(u)=\big(\Psi_{\bR}(\emptyset)\big)_i(u)
\end{equation*}
and thus $\Psi_{\bR'}(\bT)=\Psi_{\bR}(\emptyset)$.
It hence follows that, for any $\bS\in \Lambda$,
\begin{equation*}
[ \Psi_{\bR'} ]([\bS]\ast [\bT])=[\Psi_{\emptyset}]([\bS])\cdot[\Psi_{\bR'}]([\bT])=[\Psi_{\emptyset}]([\bS])\cdot[\Psi_{\bR}]([\emptyset])=[\Psi_{\bR}]([\bS]).
\end{equation*}
Consequently, by Proposition \ref{prop:injectivity_chiGT},
\begin{equation}\label{eq:get_back_to_chiell_from_chiGT}
[\Psi_{\bR'}](\chiGT^{\bR}(V)\ast [\bT]) = ([\Psi_{\bR}]\circ \chiGT^{\bR})(V)=\chi_{\ell}(V)%
=[\Psi_{\bR'}](\chiGT^{\bR'}(V))
\end{equation}
and the injectivity of $[\Psi_{\bR'}]:\EGT\to \Z^{\mathfrak{r}}$ implies the lemma.
\end{proof} 
	
\section{Results on tensor products}\label{sec:TensorResults}
This section proves various new results about tensor products in the category $\mathcal{O}_{sh}$ of the preceding section. More precisely, we consider in Section \ref{sec:Auto} distinguished anti-involutions of shifted Yangians and use them to construct a canonical contravariant autofunctor of $\mathcal{O}_{sh}$ that we prove is compatible with tensor products. We then use this functor and associators of \cite{zhang2024theta} to extend results of \cite{hernandez2024shifted}, producing along the way \textit{normalized $R$-matrices} for many pairs of simple objects of $\mathcal{O}_{sh}$, and deduce that tensor products of simple modules in $\mathcal{O}_{sh}$ are generically irreducible. Finally, we show in Section \ref{sec:TSC} that the coproducts $\Delta_{\mu,\nu}$ of Section \ref{sec:Coproduct} give %
coproducts $\Delta_{\mu_1,\mu_2}^{\lambda_1,\lambda_2}(\bR_1,\bR_2):Y_{\mu_1+\mu_2}^{\lambda_1+\lambda_2}(\bR_1\cup\bR_2)\rightarrow Y_{\mu_1}^{\lambda_1}(\bR_1)\otimes Y_{\mu_2}^{\lambda_2}(\bR_2)$ for truncations. Byproducts of this construction are: (1) %
another criterion equivalent~to~those of Theorem \ref{th:TFAE}
and (2) 
a link between the blocks of Theorem \ref{thm:BlockDec} and truncations. 
\subsection{An anti-involution}\label{sec:Auto} Fix a coweight $\mu$ and consider the involutive anti-automorphism $\text{tr}_{\mu}$ of $Y_{\mu}$ that interchanges %
$e_i(u)$ with $f_i(u)$ and fixes $h_i(u)$ for each $i\in I$. We call this map the \textit{transpose}. It defines a duality $\mathbf{tr}_{\mu}$ of the category $\mathcal{O}_{\mu}$ in the usual way, i.e.~it sends an object $V=\bigoplus_{\omega\in \mathfrak{h}^*}V_{\omega}$ to its graded dual $V^{\star} = \bigoplus_{\omega\in \mathfrak{h}^*}V_{\omega}^*$ with action given implicitly by
$$ (x\cdot f)(v) = f(\text{tr}_{\mu}(x)v) $$ 
(for all $x\in Y_{\mu}$, $f\in V^{\star}$ and $v\in V$) and sends a morphism $\phi:V_1\rightarrow V_2$ to the map $V_2^{\star}\rightarrow V_1^{\star}$ given by precomposition by $\phi$.\medskip\par
Define a duality autofunctor of $\mathcal{O}_{sh}$ via $\mathbf{tr}=\bigoplus_{\mu\in P^{\vee}} \mathbf{tr}_{\mu}$. The lemma below is clear%
.
\begin{Lemma}\label{lem:trSimples} Fix $V$ in $\mathcal{O}_{sh}$. Then $\chi_{\ell}(\mathbf{tr}(V))=\chi_{\ell}(V)$. Thus $\mathbf{tr}(L(\psi))\simeq L(\psi)$ for $\psi\in \mathfrak{r}$. 
\end{Lemma}
\begin{Theorem}\label{thm:trTensor} Let $\mu,\nu\in P^{\vee}$. Then $$\smash{(\text{tr}_{\mu}\otimes \text{tr}_{\nu})\circ \Delta_{\nu,\mu}^{\text{op}}\circ\text{tr}_{\mu+\nu}=\Delta_{\mu,\nu}}$$ 
and there %
are 
thus 
natural isomorphisms $\mathbf{tr}(V_1\otimes V_2)\simeq \mathbf{tr}(V_2)\otimes \mathbf{tr}(V_1)$ for all %
$V_1,V_2$~in~$\mathcal{O}_{sh}$.
\end{Theorem}

\begin{proof} Denote $\nabla_{\mu,\nu}=(\text{tr}_{\mu}\otimes \text{tr}_{\nu})\circ \Delta_{\nu,\mu}^{\text{op}}\circ\text{tr}_{\mu+\nu}$. By Theorem \ref{thm:CoproductUnicity}, it suffices to show that:
\begin{itemize}
\item[(i)] $\nabla_{0,0}=\Delta_{0,0}$ is the coproduct of the Yangian $Y$,
\item[(ii)] $\nabla_{\mu,\nu}(e_{i,n})=e_{i,n}\otimes 1$ for $n<-\langle\mu,\alpha_i\rangle$ and $\mu,\nu$ antidominant,
\item[(iii)] $\nabla_{\mu,\nu}(f_{i,m})=1\otimes f_{i,m}$ for $m<-\langle\nu,\alpha_i\rangle$ and $\mu,\nu$ antidominant,
\item[(iv)] the diagram \vspace*{-1.1mm}
$$
\adjustbox{scale=0.8}{
\begin{tikzcd}[column sep = 4.5em, row sep = 1.25em]
Y_{\mu+\nu} \ar[r,"\nabla_{\mu,\nu}"] \ar[d,"\iota_{\mu+\nu,\zeta,\eta}",swap] & Y_{\mu}\otimes Y_{\nu}\ar[d,"\iota_{\mu,\zeta,0}\otimes\iota_{\nu,0,\eta}"] \\
Y_{\mu+\nu+\zeta+\eta} \ar[r,"\nabla_{\mu+\zeta,\nu+\eta}"] & Y_{\mu+\zeta} \otimes Y_{\nu+\eta}
\end{tikzcd}}
$$\noindent
commutes for $\zeta,\eta$ antidominant (where the $\iota$'s are the maps of Section \ref{sec:Algebras}), and
\item[(v)] the diagram \vspace*{-1.1mm}
$$
\adjustbox{scale=0.8}{
\begin{tikzcd}[column sep = 4.5em, row sep = 1em]
Y_{\mu+\rho+\nu} \ar[r,"\nabla_{\mu+\rho,\nu}"] \ar[d,"\nabla_{\mu,\rho+\nu}",swap] & Y_{\mu+\rho}\otimes Y_{\nu}\ar[d,"\nabla_{\mu,\rho}\otimes 1"] \\
Y_{\mu}\otimes Y_{\rho+\nu} \ar[r,"1\otimes \nabla_{\rho,\nu}"] & Y_{\mu}\otimes Y_{\rho}\otimes Y_{\nu}
\end{tikzcd}}
$$\noindent
commutes for $\rho$ antidominant.
\end{itemize}\medskip
Properties (ii)--(v) are easily shown using the associated properties of the coproduct $\Delta_{\mu,\nu}$. For instance, to prove (v), one can construct the diagram
\begin{equation*}
\begin{tikzpicture}[baseline= (a).base]
\node[scale=0.68] (a) at (0,0){
\begin{tikzcd}[column sep = 2.3em]
Y_{\mu+\rho+\nu} \ar[d,"\text{tr}_{\mu+\rho+\nu}",swap]\ar[rr,"\text{tr}_{\mu+\rho+\nu}"] &&[-7pt] Y_{\mu+\rho+\nu} \ar[rr,"\Delta_{\nu,\mu+\rho}"] &&[-7pt] Y_{\nu}\otimes Y_{\mu+\rho} \ar[rrr,"\text{flip}"] \ar[dd,"1\otimes \Delta_{\rho,\mu}"] &&&[-10pt] Y_{\mu+\rho}\otimes Y_{\nu} \ar[dd,"\Delta_{\rho,\mu}\otimes 1"]\ar[rr,"\text{tr}_{\mu+\rho}\otimes \text{tr}_{\nu}"] && Y_{\mu+\rho}\otimes Y_{\nu}\ar[d,"\text{tr}_{\mu+\rho}\otimes 1"]\\
Y_{\mu+\rho+\nu} \ar[d,"\Delta_{\rho+\nu,\mu}",swap] && Y_{\mu+\rho+\nu}\ar[-,double line with arrow={-,-}]{u} \ar[-,double line with arrow={-,-}]{ll} \ar[d,"\Delta_{\rho+\nu,\mu}",swap] && &&& && Y_{\mu+\rho}\otimes Y_{\nu}\ar[d,"\Delta_{\rho,\mu}\otimes 1"] \\
Y_{\rho+\nu}\otimes Y_{\mu} \ar[d,"\text{flip}",swap] && Y_{\rho+\nu}\otimes Y_{\mu} \ar[-,double line with arrow={-,-}]{ll} \ar[d,"\text{flip}",swap] \ar[rr,"\Delta_{\nu,\rho}\otimes 1"] && Y_{\nu}\otimes Y_{\rho}\otimes Y_{\mu} \ar[d,"\text{flip}_{12}\circ\text{flip}_{23}"]\ar[rrr,"\text{flip}_{23}\circ\text{flip}_{12}"] &&& Y_{\rho}\otimes Y_{\mu}\otimes Y_{\nu} \ar[rr,"1\otimes 1\otimes \text{tr}_{\nu}"]\ar[d,"\text{flip}_{12}"] && Y_{\rho}\otimes Y_{\mu}\otimes Y_{\nu}\ar[d,"\text{flip}_{12}"]\\
Y_{\mu}\otimes Y_{\rho+\nu} \ar[-,double line with arrow={-,-}]{rr} \ar[d,"\text{tr}_{\mu}\otimes \text{tr}_{\rho+\nu}",swap] && Y_{\mu}\otimes Y_{\rho+\nu} \ar[d,swap,"\text{tr}_{\mu}\otimes 1"]\ar[rr,"1\otimes \Delta_{\nu,\rho}"] && Y_{\mu}\otimes Y_{\nu} \otimes Y_{\rho}\ar[d,"\text{tr}_{\mu}\otimes 1\otimes 1"]\ar[rrr,"\text{flip}_{23}"] &&& Y_{\mu}\otimes Y_{\rho}\otimes Y_{\nu} \ar[d,"\text{tr}_{\mu}\otimes 1\otimes 1"]\ar[rr,"1\otimes 1\otimes \text{tr}_{\nu}"] && Y_{\mu}\otimes Y_{\rho}\otimes Y_{\nu}\ar[d,"\text{tr}_{\mu}\otimes \text{tr}_{\rho}\otimes 1"]\\
Y_{\mu}\otimes Y_{\rho+\nu} \ar[rr,"1\otimes \text{tr}_{\rho+\nu}"] && Y_{\mu}\otimes Y_{\rho+\nu}\ar[rr,"1\otimes \Delta_{\nu,\rho}"] && Y_{\mu}\otimes Y_{\nu}\otimes Y_{\rho} \ar[rrr,"\text{flip}_{23}"] &&& Y_{\mu}\otimes Y_{\rho}\otimes Y_{\nu}\ar[rr,"1\otimes\text{tr}_{\rho}\otimes\text{tr}_{\nu}"] && Y_{\mu}\otimes Y_{\rho}\otimes Y_{\nu}
\end{tikzcd}
};
\end{tikzpicture}
\end{equation*}
for which all the interior polygons are easily seen to commute. The big boundary diagram, which is  equivalent to the one in (v), thus commutes as well and it only remains to show~(i). To do this, recall from Section \ref{sec:Coproduct} that the coproduct $\Delta_{0,0}$ of $Y$ is uniquely determined by the conditions $\Delta_{0,0}(x)=x\otimes 1+1\otimes x$ for $x\in \mathfrak{g}\subseteq Y$ and
$$\textstyle \Delta_{0,0}(h_{i,2}) = h_{i,2}\otimes 1+1\otimes h_{i,2}+h_{i,1}\otimes h_{i,1} +\sum_{\gamma\in \Delta_+}(\alpha_i,\gamma)f_{\gamma}\otimes e_{\gamma}$$
where $\Delta_+$ is the set of positive roots of $\mathfrak{g}$ and with $f_{\gamma}\in \mathfrak{g}_{-\gamma}$ and $e_{\gamma}\in \mathfrak{g}_{\gamma}$ suitably normalized root vectors. Now, since $\text{tr}_0$ restricts to the Chevalley involution of $\mathfrak{g}$, 
$$\nabla_{0,0}(x)=x\otimes 1+1\otimes x$$
for $x\in \mathfrak{g}$ and 
$$\text{tr}_0(f_{\gamma})=a_{\gamma}e_{\gamma}$$ 
for all $\gamma\in \Delta_+$ where $\{a_{\gamma}\}_{\gamma\in \Delta^+}$ is some subset of $\mathbb{C}^{\times}$. Hence, $(\text{tr}_0\otimes \text{tr}_0)(f_{\gamma}\otimes e_{\gamma})=e_{\gamma}\otimes f_{\gamma}$ by involutivity and $\nabla_{0,0}(h_{i,2}) = (\text{tr}_0\otimes \text{tr}_0)\circ\Delta_{0,0}^{\text{op}}(h_{i,2}) = \Delta_{0,0}(h_{i,2})$. This ends the proof.
\end{proof}
\begin{Rem} The contravariance and exactness of $\mathbf{tr}$ with the above two results imply that the tensor products $L(\psi_1)\otimes L(\psi_2)$ and $L(\psi_2)\otimes L(\psi_1)\simeq \mathbf{tr}(L(\psi_1)\otimes L(\psi_2))$ have reciprocal structures (i.e.~they share the same composition factors, but these factors are organized in reciprocal ways in the two, typically non-isomorphic, tensor products). This gives another proof (simpler than the one using Theorem \ref{thm:K0inj}) that %
$K_0(\mathcal{O}_{sh})$ is commutative (and~one~can even see the existence of the functor $\mathbf{tr}$ as \textit{categorifying} this commutativity statement).
\end{Rem}
\subsection{$R$-matrices in $\mathcal{O}_{sh}$}\label{sec:TensO} In \cite{hernandez2024shifted}, the authors define morphisms $R_{V,W}:V\otimes W\rightarrow W\otimes V$, that they call \textit{R-matrices}, for particular tensor products of modules in $\mathcal{O}_{sh}$. More precisely, Hernandez and Zhang first show the following result:
\begin{Theorem}[{\cite[Theorem 4.8]{hernandez2024shifted}}]\label{thm:HZhighest} Let $V,W$ be objects of $\mathcal{O}_{sh}$. Then the tensor product $V\otimes W$ is of highest $\ell$-weight if either 
\begin{itemize}
\item[(i)] $V$ is a positive prefundamental representation and $W$ is of highest $\ell$-weight, or
\item[(ii)] $V$ is of highest $\ell$-weight 
and $W$ is a negative prefundamental representation. %
\end{itemize}\smallskip
Similarly, $V\otimes W$ is a co-highest $\ell$-weight module if either 
\begin{itemize}
\item[(iii)] $V$ is of co-highest $\ell$-weight and $W$ is a positive prefundamental representation, or
\item[(iv)] $V$ is a negative prefundamental representation and $W$ is of highest $\ell$-weight. 
\end{itemize}
\end{Theorem}\smallskip
Here, a $Y_{\mu}$-module $U$ of \textit{co-highest $\ell$-weight $\psi$} is a representation such that \smallskip
\begin{itemize}
\item[(1)] $U$ is \textit{top-graded}, that is $\dim U_{\wt(\psi)}=1$ with $U_{\omega} = 0$ unless $\omega\leq \wt(\psi)$, and
\item[(2)] $U_{\psi}$ is contained in every non-zero submodule of $U$. 
\end{itemize}\smallskip
Clearly, a $Y_{\mu}$-module $U$ is of highest $\ell$-weight $\psi$ if and only if %
its dual $\mathbf{tr}(U)$ is of co-highest $\ell$-weight $\psi$. Furthermore, in this situation, 
$$ \head(U)\simeq L(\psi)\simeq \soc(\mathbf{tr}(U))$$
and there is, up to scalar, a unique morphism $U\rightarrow\mathbf{tr}(U)$. An easy corollary of this is:
\begin{Corollary}[{\cite[Theorem 5.2]{hernandez2024shifted}}]\label{cor:RmatHZ} Fix $V,W$ simple objects in $\mathcal{O}_{sh}$ with highest $\ell$-weight vectors $v\in V$ and $w\in W$. Suppose that either
\begin{itemize}
\item[(i)] $V$ is a positive prefundamental representation, or
\item[(ii)] $W$ is a negative prefundamental representation.
\end{itemize}
Then there is a unique morphism $R_{V,W}:V\otimes W\rightarrow W\otimes V$ that sends $v\otimes w$ to $w\otimes v$.
\end{Corollary}\newpage
\begin{proof}
Consider case (i). Then $W$ is of highest $\ell$-weight (since it is simple) and Theorem~\ref{thm:HZhighest} implies that $V\otimes W$ is also a highest $\ell$-weight module. Furthermore, $\mathbf{tr}(V\otimes W)\simeq W\otimes V$ by Lemma \ref{lem:trSimples} and Theorem \ref{thm:trTensor} so that the above discussion ends the proof.
\end{proof}
Our goal in this subsection is to construct maps $R_{V,W}$ as above for more general pairs of modules $(V,W)$ in $\mathcal{O}_{sh}$. Call $\mathbf{P}=(\mathbf{P}_i(u))_{i\in I}\in \mathfrak{r}$ \textit{polynomial} if $\mathbf{P}_i(u)\in \C[u]$ for all $i\in I$. (In other words, $\mathbf{P}\in \mathfrak{r}$ is polynomial if it lies inside $\cB_{\C,+}$ via the association~$\mathfrak{r}\simeq\cB_{\C}$ of Remark \ref{rem:Blockawt}%
.) Consider the following generalization of Theorem \ref{thm:HZhighest}:
\begin{Theorem}\label{thm:HZgen} Fix %
$V$ in $\mathcal{O}_{sh}$. Fix also $\mathbf{P},\mathbf{Q}\in \mathfrak{r}$ polynomial and let $W=L(\frac{\mathbf{P}}{\mathbf{Q}})$. Then, 
\begin{enumerate}[label=(\roman*)]
\item\label{item:highest_ell_weight_criterion} $V\otimes W$ (resp.~$W\otimes V$) is of highest $\ell$-weight if $V\otimes L(\mathbf{P})$ (resp.~$L(\mathbf{Q}^{-1})\otimes V$) is%
,~and
\item $V\otimes W$ (resp.~$W\otimes V$) is of co-highest $\ell$-weight if $L(\mathbf{P})\otimes V$ (resp.~$V\otimes L(\mathbf{Q}^{-1})$) is%
.
\end{enumerate}
\end{Theorem}
\begin{proof}
The product $L(\mathbf{P})\otimes L(\mathbf{Q}^{-1})$ is of highest $\ell$-weight by Theorem \ref{thm:HZhighest}. There is hence~a surjective morphism
$$\textstyle L(\mathbf{P})\otimes L(\mathbf{Q}^{-1})\twoheadrightarrow L(\frac{\mathbf{P}}{\mathbf{Q}}) = W$$
which can be combined\footnote{Recall that the shifted coproduct $\Delta_{\mu,\nu}$ of Section \ref{sec:Coproduct} is not co-associative. The isomorphisms appearing in \eqref{eq:HZgen1}--\eqref{eq:HZgen2} are hence non-trivial (see the introduction of \cite{zhang2024theta} for details).} with \cite[Theorem 3.3 and Corollary 5.9]{zhang2024theta} to obtain surjections
\begin{equation}\label{eq:HZgen1}
(V\otimes L(\mathbf{P}))\otimes L(\mathbf{Q}^{-1}) \simeq V\otimes (L(\mathbf{P})\otimes L(\mathbf{Q}^{-1}))\twoheadrightarrow V\otimes W
\end{equation}
and 
\begin{equation}\label{eq:HZgen2}
L(\mathbf{P})\otimes(L(\mathbf{Q}^{-1})\otimes V) \simeq (L(\mathbf{P})\otimes L(\mathbf{Q}^{-1}))\otimes V\twoheadrightarrow W\otimes V.
\end{equation}
Now, assuming $V\otimes L(\mathbf{P})$ (resp.~$L(\mathbf{Q}^{-1})\otimes V$) is of highest $\ell$-weight, Theorem \ref{thm:HZhighest} implies that $(V\otimes L(\mathbf{P}))\otimes L(\mathbf{Q}^{-1})$ (resp.~$L(\mathbf{P})\otimes(L(\mathbf{Q}^{-1})\otimes V)$) is also of highest $\ell$-weight. Part (i) thus follows from either \eqref{eq:HZgen1} or \eqref{eq:HZgen2} since quotients of highest $\ell$-weight modules are always themselves of highest $\ell$-weight. Part (ii) also follows from part (i) using the duality $\mathbf{tr}$.
\end{proof}
\begin{Rem} Theorem \ref{thm:HZgen} reduces to Theorem \ref{thm:HZhighest} if either $\mathbf{P}$ or $\mathbf{Q}$ is the highest $\ell$-weight $\mathbbm{1}=(1)_{i\in I}$ of the trivial representation of $Y=Y_0$.
\end{Rem}%
The following corollary is an easy consequence of the fact that modules in $\mathcal{O}_{sh}$ are simple if and only if they are both of highest $\ell$-weight and of co-highest $\ell$-weight. (Recall also~that, by Theorem \ref{thm:K0inj}, a product $V\otimes W$ is simple if and only if the product $W\otimes V$ also is.)
\begin{Corollary}\label{cor:HZgenSimple} Fix $\mathbf{P}_1,\mathbf{P}_2,\mathbf{Q}_1,\mathbf{Q}_2\in \mathfrak{r}$ polynomial with $V=L(\frac{\mathbf{P}_1}{\mathbf{Q}_1})$ and $W=L(\frac{\mathbf{P}_2}{\mathbf{Q}_2})$.~Then
\begin{itemize}
\item[(i)] $V\otimes W$ is of highest $\ell$-weight if either $V\otimes L(\mathbf{P}_2)$ or $L(\mathbf{Q}_1^{-1})\otimes W$ is, and
\item[(ii)] $V\otimes W$ is of co-highest $\ell$-weight if either $L(\mathbf{P}_1)\otimes W$ or $V\otimes L(\mathbf{Q}_2^{-1})$ is.
\end{itemize}\smallskip
In particular, $V\otimes W$ is simple if both products $V\otimes L(\mathbf{P}_2)$ and $W\otimes L(\mathbf{P}_1)$ are simple (which is equivalent to them being both of highest $\ell$-weight). 
\end{Corollary}
Hence, to determine if a tensor product $V\otimes W$ (of GK-dimension $\GKdim V+\GKdim W$) is simple, it can be (and actually is, in most cases) enough to determine the simplicity of~two simpler tensor products (of respective GK-dimension $\GKdim V$ and $\GKdim W$). Moreover, as shown in \cite{hernandez2024shifted}, the question of whether a product $V\otimes L(\mathbf{P})$ is simple (with $V$ a simple module of $\mathcal{O}_{sh}$ and $\mathbf{P}\in \mathfrak{r}$ polynomial) is equivalent to applying a combinatorial criterion on the normalized $\ell$-character $\chi_{\ell}(V)$ of $V$. We recall this remarkable result below:
\begin{Theorem}[{\cite[Corollary 5.10]{hernandez2024shifted}}]\label{thm:HZcriterionTensSimpPol} Fix $V$ simple in $\mathcal{O}_{sh}$ with $\mathbf{P}\in \mathfrak{r}$ polynomial. Write %
$$ \mathbf{P}=\mathsf{\Psi}_{i_1,a_1}\dots\mathsf{\Psi}_{i_k,a_k}$$
for some $i_1,\dots,i_k\in I$ and $a_1,\dots,a_k\in \C$. Then the following statements are equivalent:\smallskip
\begin{itemize}
\item[(i)] $V\otimes L(\mathbf{P})$ is simple, and
\item[(ii)] the variable $\mathsf{A}_{i_r,a_r}^{-1}$ does not appear in the normalized $\ell$-character of $V$ for $1\leq r \leq k$.
\end{itemize}
\end{Theorem}
\begin{Rem}\label{rem:OtherCombCritDuality} In Section \ref{sec:duality}, we will use an algebra automorphism $D$ of $K_0(\O_{sh})$ (that generalizes a morphism given in \cite{hernandez2016cluster}) to obtain a similar combinatorial criterion~as~above, but for the simplicity of tensor products of the form $V\otimes L(\mathbf{Q}^{-1})$ (with $V$ a simple module of $\mathcal{O}_{sh}$ and $\mathbf{Q}\in \mathfrak{r}$ polynomial). We suspect that this other combinatorial criterion is new.
\end{Rem}

\begin{Example} It is natural to wonder whether the two distinct sufficient conditions given in part (i) of Corollary \ref{cor:HZgenSimple} are in fact equivalent. The answer is unfortunately no in general. Indeed, take $\fg = \mathfrak{sl}_6$ and consider $V=L(\frac{\mathbf{P'}}{\mathbf{Q}})$ with $W=L(\frac{\mathbf{P}}{\mathbf{Q}})$, where
$$\mathbf{P}=\mathsf{\Psi}_{3,1},\,\ \mathbf{Q}=\mathsf{\Psi}_{2,0}\mathsf{\Psi}_{4,0}\,\text{ and }\,\mathbf{P}'=\mathsf{\Psi}_{1,1}\mathsf{\Psi}_{3,1}\mathsf{\Psi}_{5,1}=\mathbf{P}\mathsf{\Psi}_{1,1}\mathsf{\Psi}_{5,1}.$$
We claim that $L(\mathbf{Q}^{-1})\otimes W\simeq L(\psi)$ for $\psi = \frac{\mathbf{P}}{\mathbf{Q}^2}$. Indeed, it is easy to show that 
$$\textstyle \frac{\mathbf{P}}{\mathbf{Q}}\in\cB(\varpi_3^{\vee},5)\,\text{ and }\,\mathbf{Q}^{-1}\in \cB(\varpi_2^{\vee}+\varpi_4^{\vee},(\{6\},\{6\}))$$ 
under the identification of Remark \ref{rem:Blockawt}. Hence $\psi\in \cB(\la,\bR)_{\mu}$ where 
$$\la = \varpi_2^{\vee}+\varpi_3^{\vee}+\varpi_4^{\vee},\ \,\mu = \varpi_3^{\vee}-2\varpi_2^{\vee}-2\varpi_4^{\vee}\,\text{ and }\,\bR=(\{6\},\{5\},\{6\}).$$
Moreover, the integral crystal $\cB(\la,\bR)$ is easily proven to be isomorphic to $\cB(\la)$ (as abstract $\fg^{\vee}$-crystals). Thus, since $\mu = s_3w_0\la$ lies in the Weyl orbit $W\la$, the category $\mathcal{O}^{\la}_{\mu}(\bR)$ contains only one simple object (up to isomorphism)%
. On the other hand, Theorem \ref{thm:truncoprod} (that~we~will prove in Section \ref{sec:TSC}) shows that 
$\mathcal{O}_{\mu}^{\la}(\bR)$ naturally contains the product 
$L(\mathbf{Q}^{-1})\otimes W$. Hence
$$ \chi_{\ell}(L(\mathbf{Q}^{-1})\otimes W)=\chi_{\ell}(L(\mathbf{Q}^{-1}))\chi_{\ell}(W) = N\chi_{\ell}(L(\psi))$$
for some $N\in \mathbb{Z}_{>0}$, but $N$ must equal 1 since %
$\dim(L(\mathbf{Q}^{-1})\otimes W)_{\psi}=1$. This shows our claim by Theorem \ref{thm:K0inj}. %
In particular, $V\otimes W$ is of highest $\ell$-weight by Corollary \ref{cor:HZgenSimple}. However, $V\otimes L(\mathbf{P})$ is not simple as follows from the $\ell$-character formulas to appear in \cite{otherpaper}. 
\end{Example}
\begin{Rem}\label{rem:Curtis} It is worth noting that the combinatorial criterion obtained by combining the last part of Corollary \ref{cor:HZgenSimple} with Theorem \ref{thm:HZcriterionTensSimpPol} greatly extends a known result for products of finite-dimensional simple modules over the Yangian $Y_0$ (see,~e.g.,~\cite[Theorem 1.4]{gautam2023poles}). 
\end{Rem}
An easy consequence of Corollary \ref{cor:HZgenSimple} is the following extension of Corollary \ref{cor:RmatHZ}. This  is the main result of this subsection.
\begin{Corollary}\label{cor:Rmatgen} Take $\psi_1,\psi_2\in \mathfrak{r}$ and fix highest $\ell$-weight vectors $v\in V$ and $w\in W$ where $V=L(\psi_1)$ and $W=L(\psi_2)$. Suppose that either:%
\begin{itemize}
\item[(i)] $V\otimes L(\mathbf{P})$ is simple for some $\mathbf{P},\mathbf{Q}\in \mathfrak{r}$ satisfying $\psi_2=\frac{\mathbf{P}}{\mathbf{Q}}$, or
\item[(ii)] $L(\mathbf{Q}^{-1})\otimes W$ is simple for some $\mathbf{P},\mathbf{Q}\in \mathfrak{r}$ satisfying $\psi_1=\frac{\mathbf{P}}{\mathbf{Q}}$.
\end{itemize}
Then there is a unique morphism $R_{V,W}:V\otimes W\rightarrow W\otimes V$ that sends $v\otimes w$ to $w\otimes v$.
\end{Corollary}\newpage
Following \cite{hernandez2024shifted}, we call \textit{$R$-matrices} the maps defined in the above corollary. We expect these maps to facilitate, in the future, the proof of numerous results about the category~$\mathcal{O}_{sh}$ (like the maps from Corollary \ref{cor:RmatHZ} did in \cite{hernandez2024shifted}) and thus think of Corollary \ref{cor:Rmatgen} as a result of independent interest from the principal results of this article. In particular, as mentioned in Section \ref{sec:Intro}, we expect our $R$-matrices to be part of a collection of \textit{renormalized $R$-matrices} (in the spirit of Kashiwara--Kim--Oh--Park, see, e.g., \cite{kashiwara2024monoidal}) which could potentially play a key role in the proof that $\mathcal{O}_{sh}$ is, as conjectured in
\cite{geiss2024representations}, a monoidal categorification~of a cluster algebra. We also expect that our $R$-matrices could help %
define \textit{monoidal Jantzen filtrations} (as introduced in \cite{fujita2024monoidal}) for the setting of modules over shifted Yangians. This would then give rise to a canonical deformation $\mathscr{K}_t$ of the Grothendieck ring $K_0(\mathcal{O}_{sh})$ (see \cite{paganelli2025quantum} for a related deformation coming from the perspective of cluster algebras, and the upcoming Remark \ref{rem:t_deformation_from_shuffles} for yet another deformation coming from parity KLRW-algebras).\medskip\par
We end this subsection by giving an explicit consequence of the above results, namely an extension to $\mathcal{O}_{sh}$ of (a weak version of) the main results of \cite{hernandez2010simple,hernandez2019cyclicity} (which~can~themselves be seen as the culmination of results shown in \cite{chari1991quantum,frenkel2001combinatorics,chari2002braid,hernandez2010cluster})\footnote{We were told by D. Hernandez that the main results of \cite{hernandez2010simple,hernandez2019cyclicity} need the additional hypothesis that all simple modules considered (except maybe one) are real. Note that this hypothesis is not needed~here.}.~We expect that this extension will be helpful for the proof that $\mathcal{O}_{sh}$ is a monoidal categorification of a cluster algebra, and will talk briefly about stronger versions of it after~we~introduce in Section \ref{sec:duality} the automorphism $D$ mentioned in Remark \ref{rem:OtherCombCritDuality}.
\medskip\par
We will assume the (widely believed) conjecture below, for which substantial evidence~was recently communicated to us by H. Zhang. %
\begin{Conjecture}\label{conj:Associators} 
Take $V_1,V_2,V_3$ in $\mathcal{O}_{sh}$. Then there exists an isomorphism 
$$V_1\otimes(V_2\otimes V_3)\simeq (V_1\otimes V_2)\otimes V_3.$$
Thus, $k$-fold tensor products in $\mathcal{O}_{sh}$ are independent of the choice of parenthesization, up~to (possibly non-canonical) isomorphism.
\end{Conjecture}
\begin{Corollary}\label{cor:BigTensProd} For $1\leq r\leq k$, fix $\mathbf{P}_r,\mathbf{Q}_r\in \mathfrak{r}$ polynomial and let $V_r=L(\frac{\mathbf{P}_r}{\mathbf{Q}_r})$. Suppose 
that Conjecture \ref{conj:Associators} holds. Then (omitting parentheses)
\begin{itemize}
\item[(i)] $V_1\otimes \dots \otimes V_k$ is of highest $\ell$-weight if $V_r\otimes L(\mathbf{P}_s)$ is simple for $1\leq r<s\leq k$, and
\item[(ii)] $V_1\otimes \dots \otimes V_k$ is of co-highest $\ell$-weight if $V_r\otimes L(\mathbf{P}_s)$ is simple for $1\leq s<r\leq k$.
\end{itemize}
In particular,
\begin{itemize}
\item[(iii)] $V_1\otimes \dots \otimes V_k$ is simple if $V_r\otimes L(\mathbf{P}_s)$ is simple for all $1\leq r,s\leq k$ with $r\neq s$.
\end{itemize}
\end{Corollary}
\begin{proof}
We prove (i) by induction on $k$ with the case $k=1$ being trivial. Assume thus $k\geq 2$. Then $V=V_1\otimes \dots\otimes V_{k-1}$ is of highest $\ell$-weight by the induction hypothesis. In~addition, %
$$\textstyle V_1\otimes \dots\otimes V_k\simeq V\otimes V_k$$
(by Conjecture \ref{conj:Associators}) and it suffices to show, because of Theorem \ref{thm:HZgen}, that $V\otimes L(\mathbf{P}_k)$ is of highest $\ell$-weight. For this, remark that, by our hypothesis (and Conjecture \ref{conj:Associators} again),
$$ \textstyle V\otimes L(\mathbf{P}_k) \simeq V'\otimes (V_{k-1}\otimes L(\mathbf{P}_k))\simeq V'\otimes L(\frac{\mathbf{P}_{k-1}\mathbf{P}_k}{\mathbf{Q}_{k-1}})=V_1\otimes \dots \otimes V_{k-2}\otimes L(\frac{\mathbf{P}_{k-1}\mathbf{P}_k}{\mathbf{Q}_{k-1}})$$
where $V'=V_1\otimes \dots \otimes V_{k-2}$ (which we choose to be the trivial representation of $Y$ if $k=2$).\newpage Using the induction hypothesis again, we see that it suffices to show that $V_r\otimes L(\mathbf{P}_{k-1}\mathbf{P}_k)$ is simple for all $1\leq r<k-1$. Fix hence such a $1\leq r< k-1$ and write 
$$ \textstyle \mathbf{P}_{k-1}=\prod_{(i,a)\in \mathscr{P}_1} \mathsf{\Psi}_{i,a} \,\text{ and }\, \mathbf{P}_{k}=\prod_{(j,b)\in \mathscr{P}_2} \mathsf{\Psi}_{j,b}$$
for multisets $\mathscr{P}_1,\mathscr{P}_2$ of elements of $I\times \C$. Let $\tilde{\chi}_{\ell}=\widetilde{\chi}_{\ell}(V_r)$. By Theorem \ref{thm:HZcriterionTensSimpPol}, 
\begin{align*}
V_r\otimes L(\mathbf{P}_{k-1}\mathbf{P}_k) \text{ is simple} &\iff \mathsf{A}_{i,a}^{-1}\text{ does not appear in }\tilde{\chi}_{\ell} \text{ for all }(i,a)\in \mathscr{P}_1\cup\mathscr{P}_2\\
&\iff \mathsf{A}_{i,a}^{-1}\text{ is not in }\tilde{\chi}_{\ell} \text{ for all }(i,a)\in \mathscr{P}_1 \text{ and all } (i,a)\in \mathscr{P}_2\\
&\iff \text{both }V_r\otimes L(\mathbf{P}_{k-1})\text{ and } V_r\otimes L(\mathbf{P}_{k}) \text{ are simple}
\end{align*}
(with the last statement true by hypothesis). This ends the proof of (i) and, simultaneously, the proof of the whole corollary as (ii) follows from (i) using the duality $\mathbf{tr}$ and (iii) is simply the combination of (i) and (ii).
\end{proof}
\subsection{Generic simplicity of products} 
\label{sec:genericsimplicity}
Write $\mathfrak{r}(0)$ for the group of \textit{integral $\ell$-weights}, that is the image of the integral crystal $\cB\subseteq \cB_{\C}$ via the correspondence $\cB_{\C}\simeq \mathfrak{r}$ of Remark \ref{rem:Blockawt}. This group gives rise to a monoidal full subcategory $\mathcal{O}_{sh}(0)$ of $\mathcal{O}_{sh}$ via the next definition:

\begin{Def}\label{def:intCatO} The category $\mathcal{O}_{sh}(0)$ is defined as the Serre subcategory of $\mathcal{O}_{sh}$ generated by the simple objects with highest $\ell$-weights in $\mathfrak{r}(0)$.
\end{Def}

Clearly, one could also define a monoidal Serre subcategory $\mathscr{O}%
\subseteq \mathcal{O}_{sh}(0)$ by imposing~the stronger condition that \textbf{all} the $\ell$-weights of \textbf{all} the objects of %
the subcategory lie inside~$\mathfrak{r}(0)$. However, as we will show below, this category $\mathscr{O}$ actually coincides with $\mathcal{O}_{sh}(0)$.
\smallskip\par
We will use the result below, shown in \cite{hernandez2024shifted} for all types except $E_8$, but generalizable~to this special type because of \cite[Theorem 1.9]{neguct2025category}\footnote{More precisely, \cite[Theorem 1.9]{neguct2025category} gives the $\ell$-characters of positive prefundamental representations of quantum affine Borel algebras for any finite type $\fg$. This and \cite[Corollary 1.2.1]{varagnolo2025representations} then~give~expressions for the characters %
of negative prefundamental modules over shifted Yangians which extend %
to~type~$E_8$~the~formulas %
in \cite[Theorem 3.16]{hernandez2024shifted}%
 . The extensions of \cite[Proposition 4.12 and Theorem 4.15]{hernandez2024shifted}~easily~follow.
}.
\begin{Proposition}[{\cite[Proposition 4.12 and Theorem 4.15]{hernandez2024shifted}}]\label{prop:WeylStandards} Fix $\mathbf{P},\mathbf{Q}\in \mathfrak{r}$ polynomial and let $V$ be a module of highest $\ell$-weight $\frac{\mathbf{P}}{\mathbf{Q}}$ in $\mathcal{O}_{sh}$. Then $V$ is a quotient of %
$L(\mathbf{P})\otimes L(\mathbf{Q}^{-1})$.
\end{Proposition}
The above proposition allows us to show the following refinement of Corollary \ref{cor:Aqcarnorm}.
\begin{Lemma}\label{lem:ellcharOint} Fix V in $\mathcal{O}_{sh}$ of highest $\ell$-weight $\psi\in\mathfrak{r}(0)$. Then 
$$\widetilde{\chi}_{\ell}(V)\in \mathbb{Z}_{\geq 0}[[\mathsf{A}_{i,a}^{-1}]]_{(i,a)\in I\times_2\Z}.$$ 
In particular, the categories $\mathcal{O}_{sh}(0)$ and $\mathscr{O}$ coincide.
\end{Lemma}
\begin{proof} The result follows from Theorem \ref{thm:NakCrys} if $V=L(\mathsf{Y}_{i,a})$ for some $(i,a)\in I\times_2\Z$,~i.e.~if $V$ is a (finite-dimensional) fundamental representation in $\mathcal{O}_{sh}(0)$. Hence, the result also~holds (by Remark \ref{rem:prefundMonGen} and multiplicativity of $\ell$-characters) for $V$ a Kirillov--Reshetikhin module of $\mathcal{O}_{sh}(0)$ and the case where $V=L(\mathsf{\Psi}_{i,a}^{-1})$ for some $(i,a)\in I\times_2\Z$ then follows directly~from Proposition \ref{prop:Prefundlim}. Finally, the result holds trivially for any simple module $L(\mathbf{P})$ with $\mathbf{P}\in \mathfrak{r}_0$ polynomial and the above discussion shows that it also holds for tensor products of negative prefundamental representations. One can thus finish the proof using Proposition \ref{prop:WeylStandards} (with Remark \ref{rem:prefundMonGen} and the multiplicativity of $\ell$-characters again).
\end{proof}
\begin{Rem} Fix $\mu\in P^{\vee}$ and notice that highest $\ell$-weight modules in $\mathcal{O}_{\mu}$ all belong~to~the category $\mathcal{O}_{\mu}^{\text{HZ}}$ of Remark \ref{rem:2categoriesO}. Indeed, for such a module $M$ with highest $\ell$-weight vector $m$, each weight vector $m'\in M$ has the form $m'=xm$ for some $x\in Y_{\mu}^-$ of $Q$-degree $\alpha\in-Q_+$ (where this $Q$-degree is defined by $\deg(e_{i,q})=-\deg(f_{i,q})=\alpha_i$ and $\deg(h_{i,p})=0$ for $i\in I$, $q\in \Z_{>0}$~and $p\in \Z$). In particular, using \eqref{H,F}, we get 
$$h_{i,1-\langle\mu,\alpha_i\rangle}m'=h_{i,1-\langle\mu,\alpha_i\rangle}xm = 2\langle \omega+\alpha, \alpha_i^{\vee}\rangle m'$$ 
for all $i\in I$, where $\omega\in \mathfrak{h}^*$ is such that $m\in M_{\omega}$. This explains why we can use Proposition \ref{prop:WeylStandards} here even though this result was proven using the (smaller) category $\mathcal{O}_{\mu}^{\text{HZ}}$ in \cite{hernandez2024shifted}.
\end{Rem}
An interesting consequence of Lemma \ref{lem:ellcharOint} is that $\mathcal{O}_{sh}(0)=\mathscr{O}$ is a monoidal subcategory of $\mathcal{O}_{sh}$. More precisely, combining Lemma \ref{lem:ellcharOint} with Remark \ref{rem:prefundMonGen} easily gives:
\begin{Corollary}\label{cor:tensor_product_decomp_in_Osh}
Fix $\psi,\psi'\in\mathfrak{r}(0)$. Then, in $K_0(\mathcal{O}_{sh})$,
\begin{equation}\label{eq:IntTensorDecomp}
\textstyle  [L(\psi)][L(\psi')]=[L(\psi\psi')]+\sum_{\substack{\xi\in \mathfrak{r}(0) \\ \xi\preceq \psi\psi'}} n_{\xi}[L(\xi)]
\end{equation}
for some $n_{\xi}$'s in $\Z_{\geq 0}$.
\end{Corollary}
\begin{Rem} The number of non-zero summands in the right-hand side of \eqref{eq:IntTensorDecomp} is finite~by \cite[Theorem 9.5]{hernandez2024shifted} (this also follows from the results shown in Section \ref{sec:TSC}).
\end{Rem}

We would now like to show that tensor products of simple objects in the category $\mathcal{O}_{sh}$~are generically irreducible. For this, take $\mu\in P^{\vee}$ with $b\in \C$ and recall from \cite{hernandez2024shifted} the spectral shift automorphism $\tau_b:Y_{\mu}\rightarrow Y_{\mu}$ defined by 
$$ e_i(u)\mapsto e_i(u-b),\ \  f_i(u)\mapsto f_i(u-b)\,\text{ and }\,h_i(u)\mapsto h_i(u-b).$$
We denote by $V(b)$ the pullback of a $Y_{\mu}$-module $V$ by $\tau_b$. Clearly, $V(b)$ lies in $\mathcal{O}_{\mu}$~if~$V$ does and $\chi_{\ell}(V(-b))=(\chi_{\ell}(V))\circ \tilde{\tau}_{b}$ with $\tilde{\tau}_b$ the group automorphism of $\mathfrak{r}$ given by $\mathsf{\Psi}_{i,a}\mapsto \mathsf{\Psi}_{i,a+b}$. Thus, identifying the maps $\tau_b$ associated to distinct $\mu$'s, we have that, for all $V_1,V_2$ in $\mathcal{O}_{sh}$,
\begin{equation}\label{eq:specshifttensor}
[(V_1\otimes V_2)(b)] = [V_1(b)][V_2(b)].
\end{equation}
Let $\mathfrak{r}(b)$ be the image of $\mathfrak{r}(0)$ under $\tilde{\tau}_b$. Define also, for each $S\subseteq \C$, a subgroup $\mathfrak{r}(S)\subseteq \mathfrak{r}$~using every possible finite product of elements in $\bigcup_{b\in S}\mathfrak{r}(b)$. We call two subsets $S,S'$ \textit{integrally disconnected} if $S\cap (S'+2\Z)=\emptyset$.
\begin{Rem} Note that $\mathfrak{r}(S)\cap\mathfrak{r}(S')=\{1\}$ if $S,S'\subseteq \C$ are integrally disconnected.
\end{Rem}
The following theorem is the central result of this subsection. We believe that it gives the first fully general answer to the question of generic simplicity of tensor products~for~simple modules over shifted Yangians although similar results were known in more restrictive contexts (such as finite-dimensional modules over quantum affine algebras \cite[Section~1]{hernandez2019cyclicity}).
\begin{Theorem}\label{thm:TensSimpInt} Fix integrally disconnected subsets $S,S'\subseteq \C$ with $\psi\in\mathfrak{r}(S)$ and $\xi\in \mathfrak{r}(S')$. Then $L(\psi)\otimes L(\xi)$ is simple. Thus, for $V$ and $W$ simple in $\mathcal{O}_{sh}$, the tensor product $V(a)\otimes W$ is irreducible for all but countably many $a\in\C$.
\end{Theorem}
\begin{proof}
Fix $a_1,\dots,a_r\in S$ such that $\psi=\psi_1\dots\psi_r$ for some $\psi_1\in \mathfrak{r}(a_1),\dots,\psi_r\in \mathfrak{r}(a_r)$.~Then $L(\psi)$ is a composition factor of $L(\psi_1)\otimes \dots L(\psi_r)$ by Remark \ref{rem:prefundMonGen}. Hence, by (the spectral shifts of) Lemma \ref{lem:ellcharOint} and multiplicativity of $\ell$-characters,
$$\tilde{\chi}_q(L(\psi))\in \mathbb{Z}[[\mathsf{A}_{i,a_1+c}^{-1},\dots,\mathsf{A}_{i,a_r+c}^{-1}]]_{(i,c)\in I\times_2 \mathbb{Z}}$$
and it follows easily from our hypothesis that the normalized $\ell$-character $\tilde{\chi}_q(L(\psi))$ contains no variables of the form $\mathsf{A}_{i,b+c}^{-1}$ with $b\in S'$ and $(i,c)\in I\times_2\Z$. In particular, all products~of the form $L(\psi)\otimes L(\mathbf{P})$ --- with  $\mathbf{P}\in \mathfrak{r}(S')$ polynomial --- are irreducible by Theorem \ref{thm:HZcriterionTensSimpPol} and Corollary \ref{cor:HZgenSimple} proves that the tensor product $L(\psi)\otimes L(\xi)$ is of highest $\ell$-weight.~An~analogous reasoning also shows that this product is of co-highest $\ell$-weight%
.
\end{proof}
\begin{Rem} Take $\psi,\xi\in \mathfrak{r}$ with $V=L(\psi)$ and $W=L(\xi)$. By the above theorem, for all but countably many $a\in \C$, there is (up to scalar) a unique isomorphism,
$$ V(a)\otimes W \simeq L(\tilde{\tau}_a(\psi)\xi) \simeq W\otimes V(a).$$ 
The category $\mathcal{O}_{sh}$ thus contains \textit{generic $R$-matrices} for its simple modules.%
\end{Rem}
For $S\subseteq \C$, let $\mathcal{O}_{sh}(S)$ be the Serre subcategory of $\mathcal{O}_{sh}$ generated by the simple modules $L(\psi)$ with $\psi\in \mathfrak{r}(S)$. Then $\mathcal{O}_{sh}(S)\subseteq \mathcal{O}_{sh}$ is a monoidal subcategory and contains the trivial module $L(\mathbbm{1})$ for all %
 $S \subseteq \C$. However, this common object $L(\mathbbm{1})$ is,~in~the sense of Corollary \ref{cor:ExtInt} below, the only place at which we can ``glue'' (i.e.~obtain non-trivial extensions from) objects coming from $\mathcal{O}_{sh}(S)$ and $\mathcal{O}_{sh}(S')$ 
with $S,S'\subseteq \C$ integrally disconnected subsets.

\begin{Lemma}\label{lem:ExtInt} Choose $\psi,\xi\in\mathfrak{r}$ with $\psi\in \mathfrak{r}(S)$ and $\xi\not\in \mathfrak{r}(S)$ for some $S\subseteq \C$. Suppose~$\psi\not\preceq \xi$. Then $\Ext^1_{\mathcal{O}_{sh}}(L(\psi),L(\xi))=\Ext^1_{\mathcal{O}_{sh}}(L(\xi),L(\psi))=0$.
\end{Lemma}

\begin{proof}
Assume $E$ is a non-trivial extension of $L(\psi)$ by $L(\xi)$. By the proof of Theorem \ref{thm:Ext1Block}, $E$ is a module of highest $\ell$-weight $\psi$ and it follows from Lemma \ref{lem:ellcharOint} that $\psi=\xi x$ for a monomial $x$ in the $\mathsf{A}_{i,a+b}$'s with $(i,a)\in I\times_2\Z$ and $b\in S$. This however implies
$\xi\in \mathfrak{r}(S)$, which contradicts our hypothesis and gives the equality $\Ext^1_{\mathcal{O}_{sh}}(L(\psi),L(\xi))=0$.~The~other equality follows from applying the functor $\mathbf{tr}$ of Section \ref{sec:Auto}.
\end{proof}
\begin{Corollary}\label{cor:ExtInt} Fix $S,S'\subseteq \C$ integrally disconnected with $V$ in $\mathcal{O}_{sh}(S)$ and $W$ in $\mathcal{O}_{sh}(S')$. Suppose $V$ and $W$ are both of finite-length. Then 
$\Ext^1_{\mathcal{O}_{sh}}(V,W)=0$ unless $V$ or $W$ contains the trivial representation $L(\mathbbm{1})$ as a composition factor.
\end{Corollary}
\begin{proof}
The result follows from Lemma \ref{lem:ExtInt} if $V$ and $W$ are non-trivial simple modules%
.~It can be proven by double induction on the lengths of $V$ and $W$ in the general case.
\end{proof}
\begin{Example}\label{ex:qWronskian} We illustrate how the previous results can help understand the structure~of certain modules in $\mathcal{O}_{sh}$. Take $V=L(\mathsf{\Psi_{1,0}})$ and $W=L(\mathsf{\Psi}_{1,0}^{-1})$ for $\fg=\mathfrak{sl}_2$. Then 
$$\mathscr{V}=V(1)\otimes V\otimes W\otimes W(1)$$
is of highest $\ell$-weight by Theorem \ref{thm:HZhighest}, and thus indecomposable. Moreover, in $K_0(\mathcal{O}_{sh})$,
\begin{align*} 
[\mathscr{V}]&=[V][W][V(1)][W(1)]=(1+[V(2)][W(-2)])\cdot(1+[V(3)][W(-1)])\\
&=1+[V(2)][W(-2)]+[V(3)][W(-1)]+[V(2)][W(-2)][V(3)][W(-1)]
\end{align*}
where we used the \textit{quantum Wronskian} of \cite[(1.3)]{bazhanov2011baxter} (which is also an example of the \textit{extended QQ-relations} we prove in this paper). Now, using Theorem \ref{thm:HZcriterionTensSimpPol} and Theorem \ref{thm:TensSimpInt} with the $\ell$-character formulas of Example \ref{ex:lchar}, it is easy to show that $V(2)\otimes W(-2) \simeq L(\psi)$, $V(3)\otimes W(-1)\simeq L(\psi')$ and
$$(V(2)\otimes W(-2))\otimes(V(3)\otimes W(-1)) \simeq L(\psi\psi')$$
for $\psi = \frac{\mathsf{\Psi}_{1,2}}{\mathsf{\Psi}_{1,-2}}=\mathsf{A}_{1,0}^{-1}\in \mathfrak{r}(0)$ and $\psi'=\frac{\mathsf{\Psi}_{1,3}}{\mathsf{\Psi}_{1,-1}}=\mathsf{A}_{1,1}^{-1}=\tilde{\tau}_1(\psi)\in \mathfrak{r}(1)$. \medskip\par Hence, by the above,
$$ [\mathscr{V}] = 1+[L(\psi)]+[L(\psi')]+[L(\psi\psi')],$$
and $V$ has 4 composition factors (all with multiplicity 1). Also, by Lemma \ref{lem:ExtInt}, as $\psi \not\preceq \psi'$,
$$ \Ext^1_{\mathcal{O}_{sh}}(L(\psi),L(\psi'))=\Ext^1_{\mathcal{O}_{sh}}(L(\psi'),L(\psi))=0$$
and using this same lemma with $\psi\not\preceq \psi\psi'=\psi \mathsf{A}_{1,1}^{-1}\not\in \mathfrak{r}(0)$ and $\psi'\not\preceq \psi\psi'%
\not\in \mathfrak{r}(1)$ gives
\begin{align*}
0&=\Ext^1_{\mathcal{O}_{sh}}(L(\psi),L(\psi\psi')) = \Ext^1_{\mathcal{O}_{sh}}(L(\psi\psi'),L(\psi)) \\&= \Ext^1_{\mathcal{O}_{sh}}(L(\psi'),L(\psi\psi'))=\Ext^1_{\mathcal{O}_{sh}}(L(\psi\psi'),L(\psi')).
\end{align*}
Hence, by indecomposability of $\mathscr{V}$, we must have a non-split short exact sequence 
$$ 0\rightarrow L(\psi)\oplus L(\psi')\oplus L(\psi\psi')\rightarrow \mathscr{V}\rightarrow L(\mathbbm{1})\rightarrow 0$$
and the structural diagram (or Loewy/Alperin diagram) of $\mathscr{V}$ is thus%
{
$$
\begin{tikzpicture}[scale = 1]
\node (A) at (0,0) {{\small $L(\psi)$}};
\node (B) at (1.9,0) {{\small $L(\psi')$}};
\node (C) at (4,0) {{\small $L(\psi\psi')$}};
\node (D) at (1.9,1.5) {{\small $L(\mathbbm{1})$}};
\draw[->] (D) -- (A);
\draw[->] (D) -- (B);
\draw[->] (D) -- (C);
\end{tikzpicture}
$$}
\end{Example}

\subsection{Truncated shifted coproducts}\label{sec:TSC}
The goal of this section is to show Theorem \ref{thm:truncoprod}, that is to prove that the shifted coproducts $\Delta_{\mu,\nu}:Y_{\mu+\nu}\to Y_{\mu}\otimes Y_{\nu}$ of Section \ref{sec:Coproduct} are compatible with truncations. For this, it will be essential to consider a parametric version of truncated shifted Yangians, i.e.~a family $Y_\mu^\la$ of algebras over $\C^\la$, which recovers the algebras $Y_\mu^\lambda(\bR)$ after specializ%
ing at %
$\bR \in \C^\la$.  We use a slight variation of the definition of \cite[\S B(ii)]{braverman2016coulomb}.\medskip\par 
Recall that, for $\la=\sum_{i\in I}\la_i\varpi_i^{\vee}\in P^{\vee}_+$,
$$\textstyle \C^\lambda = \prod_{i \in I}  \C^{\lambda_i} / \Sigma_{\lambda_i},$$ which is naturally an algebraic variety, isomorphic to an affine space. Its coordinate ring~$\coord^\lambda$ is a ring of partially symmetric polynomials, and %
is thus itself a polynomial ring in variables  
$\{R_{i,s}\}_{i\in I,1\leq s\leq \la_i}$, where $R_{i,s}$ extracts the $u^s$-coefficient of the polynomial $p_{R_i}(u)$ given in \eqref{eq:polyp} (or, equivalently, represents the elementary symmetric polynomial $(-1)^s e_s( R_i)$). For $i\in I$, we may think of the polynomial
\begin{equation}
	\label{eq:polypvar}
	u^{\lambda_i} + R_{i,1} u^{\lambda_i-1} + \ldots + R_{i,\lambda_i-1} u +  R_{i,\lambda_i} 
\end{equation}
as an element of $\coord^\lambda[u]$. Evaluation %
at a point $\bR\in \C^\lambda$ defines a map %
$\coord^\lambda[u] \rightarrow \C[u]$, and~the image of the above polynomial is exactly $p_{R_i}(u) \in \C[u]$. \medskip\par
The results of Section \ref{section: tsy} have parametric versions, where we replace all instances of the polynomial $p_{R_i}(u)$ by the polynomials in \eqref{eq:polypvar}.  In particular, we define in this way elements $a_{i,r} \in Y_\mu \otimes \coord^\lambda$ as in \eqref{eq: def of A gens} (for $\mu\in P^{\vee}$ with $\la\geq \mu$), and construct a ``GKLO homomorphism''
$$
\Phi_\mu^\lambda : Y_\mu \otimes \coord^\lambda \longrightarrow \widetilde{\mathscr{A}}_{\lambda-\mu} \otimes \coord^\lambda
$$
as in Theorem \ref{GKLO homomorphism}. The following definition gives the parametric version of truncation.
\begin{Def}
The \emph{parametric truncated shifted Yangian} $Y_\mu^\lambda $ is defined to be the image of $\Phi_\mu^\lambda$.  Note that $Y_\mu^\lambda$ is naturally an algebra over $\coord^\lambda$.   
\end{Def}

We will now take a brief detour into the theory of Coulomb branches for $3d$ $\mathcal{N}=4$ quiver gauge theories, as developed in the foundational work of Braverman--Finkelberg--Nakajima \cite{BFN1}.  Given a triple $(\mathbf{G}, \mathbf{N}; \mathbf{F})$ consisting of a reductive group $\mathbf{G}$, its representation $\mathbf{N}$, and a choice of so-called 
flavour symmetry group $\mathbf{F}$, \cite{braverman2016coulomb} defines the \emph{(flavour deformed) quantized Coulomb branch} $\mathcal{A}_\hbar(\mathbf{G},\mathbf{N}; \mathbf{F})$ which is %
an associative algebra over the equivariant cohomology ring
$$H_{\mathbf{F} \times \C^\times}^\bullet(pt) = H_{\mathbf{F}}^\bullet(pt) \otimes \C[\hbar].$$
In particular, specializing at $\hbar =2$ we obtain an algebra $\mathcal{A}_{\hbar =2}(\mathbf{G}, \mathbf{N}; \mathbf{F})$ which is  an algebra over $H_{\mathbf{F}}^\bullet(pt)$. The next result was proven for dominant $\mu$ in \cite[Corollary B.28]{braverman2016coulomb}, and extended to all $\mu$ in \cite[Theorem A]{weekes2019}. (Recall that we write $\lambda-\mu=\sum_{i\in I} m_i\alpha_i^\vee\in Q_+^\vee$.)
\begin{Theorem}\label{thm:BFN}
There is an isomorphism of filtered algebras between $Y_\mu^\la$ and the quantized Coulomb branch $\mathcal{A}_{\hbar =2}(\mathbf{G}, \mathbf{N}; \mathbf{F})$ for the quiver gauge theory given by
$$
\mathbf{G} = \prod_{i \in I} \operatorname{GL}(m_i), \ \  \mathbf{N} = \bigoplus_{\substack{i, j \in I, \\ i \rightarrow j}} \operatorname{Hom}(\C^{m_i}, \C^{m_j}) \oplus \bigoplus_{i \in I} \operatorname{Hom}(\C^{m_i}, \C^{\la_i})\,\text{ and }\,\mathbf{F} = \prod_{i\in I} \operatorname{GL}(\la_i).
$$
This %
also identifies the commutative subalgebras $\smash{\coord^\la \subseteq Y_\mu^\la}$ and $\smash{H_{\mathbf{F}}^\bullet(pt) \subseteq \mathcal{A}_{\hbar =2}(\mathbf{G}, \mathbf{N}; \mathbf{F})}$. 
\end{Theorem}
Since the Coulomb branch is constructed via the equivariant cohomology of an equivariantly formal space \cite[\S 2]{braverman2016coulomb}, the algebra $\mathcal{A}_{\hbar=2}(\mathbf{G}, \mathbf{N}; \mathbf{F})$ is free as a module over $H_{\mathbf{F}}^\bullet(pt)$.  We can hence deduce the following modest consequence of the Coulomb branch theory:
\begin{Corollary}\label{cor:CBfree}
The algebra $Y_\mu^\lambda$ is free as a module over $\coord^\lambda$. Also, for any $\bR \in \C^\lambda$ with associated evaluation map $\coord^\lambda \rightarrow \C$, there is an algebra isomorphism $Y_\mu^\lambda \otimes_{\coord^\lambda}\C \cong Y_\mu^\lambda(\bR)$.
\end{Corollary}
For $\la_1,\la_2\in P_+^{\vee}$ such that $\la=\la_1+\la_2$, there is a natural map $\C^{\la_1}\times \C^{\la_2}\to \C^{\la}$ given by 
componentwise union $(\bR_1, \bR_2) \mapsto \bR = \bR_1 \cup \bR_2$ of multisets.  This induces an injective~map of coordinate rings 
\begin{equation*}
	\fact^{\lambda_1, \lambda_2} : \coord^\la \longrightarrow \coord^{\la_1} \otimes \coord^{\la_2},
\end{equation*}
which is given explicitly by 
$$%
\textstyle \fact^{\la_1, \la_2}(R_{i,s}) = \sum_{a+b = s} R_{i,a} \otimes R_{i,b}.$$ 
Now, fix $\mu_1,\mu_2\in P^{\vee}$ such that $\mu=\mu_1+\mu_2$ and for which $\mu_1\leq \la_1$ and $\mu_2\leq \la_2$. Then taking the tensor product of $\fact^{\la_1, \la_2}$ with the coproduct $\Delta_{\mu_1, \mu_2} $ of Theorem \ref{thm:CoproductUnicity} and permuting middle tensor factors gives the composed map
\begin{equation}
\label{eq:paramcoprod}
Y_\mu \otimes \coord^\la \xrightarrow{\Delta_{\mu_1, \mu_2}\otimes \fact^{\lambda_1, \lambda_2}} Y_{\mu_1} \otimes Y_{\mu_2} \otimes \coord^{\lambda_1} \otimes \coord^{\lambda_2} \xrightarrow{1 \otimes (\text{flip}) \otimes 1}  (Y_{\mu_1} \otimes \coord^{\lambda_1}) \otimes  (Y_{\mu_2} \otimes \coord^{\lambda_2})
\end{equation}
By composing also with the surjections $\Phi_{\mu_1}^{\lambda_1}:Y_{\mu_1} \otimes \coord^{\la_1} \twoheadrightarrow Y_{\mu_1}^{\la_1}$ and $\Phi_{\mu_2}^{\lambda_2}:Y_{\mu_2} \otimes \coord^{\la_2} \twoheadrightarrow Y_{\mu_2}^{\la_2}$, we thus obtain a map 
\begin{equation}
\label{eq:paramcoprod2}
\mathsf{\Phi}%
:Y_\mu \otimes S^\la \longrightarrow Y_{\mu_1}^{\la_1} \otimes Y_{\mu_2}^{\la_2}
\end{equation}
Consider the map 
$$\mathsf{\Pi}(\bR)%
:Y_{\mu}^{\la} \simeq Y_{\mu}^{\la}\otimes_{S^{\la}}S^{\la}\twoheadrightarrow Y_\mu^\lambda(\bR)$$ 
given by Corollary \ref{cor:CBfree} as well as the corresponding maps $\mathsf{\Pi}_1(\bR_1):Y_{\mu_1}^{\la_1}\otimes S^{\la_1}\twoheadrightarrow Y_{\mu_1}^{\la_1}(\bR_1)$ and $\mathsf{\Pi}_2(\bR_2):Y_{\mu_2}^{\la_2}\otimes S^{\la_2}\twoheadrightarrow Y_{\mu_2}^{\la_2}(\bR_2)$.\medskip\par  We have the following ideals inside $Y_\mu \otimes \coord^\la$:
\begin{align*}
I_\mu^\la  &= \Ker%
\Phi_{\mu}^{\la}%
, \\
 I_\mu^\la(\bR)  &= \Ker\big(\mathsf{\Pi}(\bR)\circ\Phi_{\mu}^{\la}
 \big), \\
 J_{\mu_1,\mu_2}^{\la_1, \la_2}  & = \Ker%
 \mathsf{\Phi}
,  \\
J_{\mu_1,\mu_2}^{\la_1, \la_2}(\bR_1, \bR_2)  &= \Ker\big((\mathsf{\Pi}_{1}(\bR_1)\otimes \mathsf{\Pi}_{2}(\bR_2))\circ\mathsf{\Phi}%
\big).
\end{align*}
The following key result is an application of Corollary \ref{cor:CBfree}:
\begin{Lemma} Take $\la_1,\la_2,\la,\mu_1,\mu_2,\mu,\bR_1,\bR_2,\bR$ as above. Then
\label{lem:zariskidense}
\begin{enumerate}
\item[(i)] If $\cZ \subseteq \C^\lambda$ is a Zariski dense set of closed points, 
$$%
I_\mu^\lambda = \bigcap_{\bR \in \cZ} I_\mu^\lambda(\bR).$$
\item[(ii)] If $\cZ' \subseteq \C^{\lambda_1} \times \C^{\la_2}$ is a Zariski dense set of closed points, 
$$%
J_{\mu_1, \mu_2}^{\lambda_1, \la_2} = \bigcap_{(\bR_1, \bR_2) \in \cZ'} J_{\mu_1,\mu_2}^{\lambda_1, \la_2}(\bR_1, \bR_2).$$
\end{enumerate}
\end{Lemma}
\begin{proof}
Fix a basis $\mathbb{B}$ for the free $\coord^\lambda$--module $Y_\mu^\lambda$ and note that $\mathbb{B}$ induces a $\C$-basis for all the specializations $Y_\mu^\lambda(\bR)$.  In particular, if $x \in Y_\mu \otimes \coord^\lambda$ has expansion $\Phi_\mu^\lambda(x) = \sum_{b \in \mathbb{B}} x_b b \in Y_\mu^\lambda$ in $\mathbb{B}$ (for some $\{x_b\}_{b\in\mathbb{B}} \subseteq \coord^\lambda$), then the image of $x$ in $Y_\mu^\lambda(\bR)$ is its evaluation $\sum_{b \in \mathbb{B}} x_b(\bR) b$.  Therefore
$$
\begin{array}{ccl}
x \in I_\mu^\lambda & \ \ \ \Longleftrightarrow \ \ \ & x_b  = 0 \text{ for all } b \in \mathbb{B} \\
& \Longleftrightarrow & x_b(\bR) = 0 \text{ for all } b \in \mathbb{B} \text{ and all } \bR \in \mathbb{C}^{\la} \\
& \Longleftrightarrow & x_b(\bR) = 0 \text{ for all } b \in \mathbb{B} \text{ and all } \bR \in \mathcal{Z} \\
& \Longleftrightarrow & x \in I_\mu^\lambda(\bR) \text{ for all } \bR \in \mathcal{Z}.
\end{array}
$$
This completes the proof of (i). For (ii), observe that, as both $\coord^{\la}$ and $\coord^{\la_1}\otimes\coord^{\la_2}$ are partially symmetric polynomial algebras, $\coord^{\la_1} \otimes \coord^{\la_2}$ is free over $\coord^\la$ by the Chevalley--Shephard--Todd theorem. Hence, the free $\coord^{\la_1}\otimes \coord^{\la_2}$-module $Y_{\mu_1}^{\la_1}\otimes Y_{\mu_2}^{\la_2}$ is also free over $\coord^{\la}$ and one can~thus use the same strategy as in part (i) to prove (ii).
\end{proof}

\begin{Corollary}
\label{cor:reductioncoproductsgeneric}
The following are equivalent:
\begin{enumerate}
	\item[(i)] Theorem \ref{thm:truncoprod} holds for all pairs $(\bR_1, \bR_2) \in \C^{\la_1} \times \C^{\la_2}$.
	
	\item[(ii)] There exists a Zariski dense set $\cZ' \subseteq \C^{\la_1} \times \C^{\la_2}$ of closed points such that Theorem \ref{thm:truncoprod} holds for all $(\bR_1, \bR_2) \in \cZ'$.
	
	\item[(iii)] There exists a map making the following diagram commute %
	{$$
	\begin{tikzcd}
	Y_\mu\otimes S^{\lambda} \ar[r,"\text{\eqref{eq:paramcoprod}}"] \ar[d, two heads,"\Phi_{\mu}^{\la}",swap] & (Y_{\mu_1}\otimes \coord^{\lambda_1}) \otimes (Y_{\mu_2} \otimes S^{\lambda_2}) \ar[d,two heads,"\Phi_{\mu_1}^{\la_1}\otimes \Phi_{\mu_2}^{\la_2}"] \\
	Y_\mu^\lambda \ar[r,dashed] & Y_{\mu_1}^{\lambda_1}\otimes Y_{\mu_2}^{\lambda_2}
	\end{tikzcd}
	$$}
\end{enumerate}
\end{Corollary}
\begin{proof}
Clearly (i) implies (ii).  Specializing (iii) at $\bR = \bR_1 \cup \bR_2$ with $(\bR_1, \bR_2) \in \C^{\lambda_1} \times \C^{\lambda_2}$, we recover the commutative diagram from Theorem \ref{thm:truncoprod}. Thus (iii) implies (i).\medskip\par  

Assume that (ii) holds. Let $\cZ' \subseteq \C^{\la_1} \times \C^{\la_2}$ be Zariski dense, and let $\cZ \subseteq \C^\la$ be its image. Note that $\cZ$ is also Zariski dense. Our assumption (ii) implies that $I_\mu^\la(\bR ) \subseteq J_{\mu_1, \mu_2}^{\la_1, \la_2}(\bR_1, \bR_2)$ for all $(\bR_1, \bR_2) \in \cZ'$ and $\bR = \bR_1 \cup \bR_2 \in \cZ$.  Applying Lemma \ref{lem:zariskidense}, we thus get
$$
\textstyle I_\mu^\la = \bigcap_{\bR \in \cZ} I_\mu^\lambda(\bR) \subseteq  \bigcap_{(\bR_1, \bR_2) \in \cZ'} J_{\mu_1,\mu_2}^{\la_1, \la_2}(\bR_1, \bR_2) =  J_{\mu_1, \mu_2}^{\la_1, \la_2} 
$$
and (iii) holds as desired.
\end{proof}
Armed with this corollary, it only remains to establish Theorem \ref{thm:truncoprod} for some sufficiently generic set $\cZ'\subseteq \C^{\la_1}\times\C^{\la_2}$.  We will do so using the representation theory of $Y_\mu^\la(\bR)$.\medskip\par 
By \cite[Lemma 4.11]{kamnitzer2022lie} (see also Lemma \ref{lem:GKdimSimples}), every object $V$ in $\mathcal{O}_\mu^\la(\bR)$ satisfies
$$
\operatorname{GKdim} V \leq \operatorname{ht}(\la - \mu).
$$
Let $\mathcal{O}_\mu^\la(\bR)_{\topRM}$ be the quotient of $\mathcal{O}_\mu^\la(\bR)$ by the full (Serre) subcategory of objects for which $\operatorname{GKdim} V< \operatorname{ht}(\la - \mu)$ and let $k$ be the number of distinct integrality classes in $\bR$. By \cite[Proposition 9.21]{kamnitzer2022lie}, there is a categorical action of the Lie algebra $\mathfrak{g}_{\bR}:=(\mathfrak{g}^{\vee})^{\oplus k}$~on %
$\mathcal{O}_{sh}^{\la}(\bR)=\bigoplus_\mu \mathcal{O}_\mu^\la(\bR)$. This induces an action of $\fg_{\bR}$ on the quotient category~$\bigoplus_\mu \mathcal{O}_\mu^\la(\bR)_{\topRM}$, which in turn allows us to see the complexified Grothendieck group $\bigoplus_\mu K_\C\big(\mathcal{O}_\mu^\la(\bR)_{\topRM}\big)$~as~an irreducible quotient of the  $\fg_{\bR}$-module $K_0(\mathcal{O}_{sh}^{\la}(\bR))$. 
\begin{Theorem}\label{thm:faithfulmodules} Consider $\la,\mu$ as above with $\bR=(R_i)_{i\in I}\in \C^{\la}$. 
\begin{enumerate}
\item[(i)] Let $V$ be in $\mathcal{O}_\mu^\la(\bR)$.  Then 
\begin{center}
$V$ is a faithful $Y_\mu^\la (\bR)$-module $\Longleftrightarrow$ $\GKdim V =\operatorname{ht}(\la - \mu)$.
\end{center}
\item[(ii)] Suppose every element of 
$\bR$ lies in a distinct integrality class\footnote{In other words, $b -a+\overline{j}-\overline{i}
\not\in 2\Z$ for all $a\in R_i$ and $b\in R_j$ with $i,j\in I$.}. Then every non-zero %
module $V$ in $\O_\mu^\lambda(\bR)$ is faithful for $Y_\mu^\lambda(\bR)$ and satisfies%
$$
I_\mu^\lambda(\bR) = \operatorname{Ann}_{Y_\mu \otimes \coord^\lambda} (V).
$$
\end{enumerate}
\end{Theorem}
\begin{proof}
Let $V \in \mathcal{O}_\mu^\la(\bR)$ with $J = \operatorname{Ann}_{Y_\mu^\la(\bR)}(V)$.  On the one hand, as $Y_\mu^\la(\bR)$ is a domain,
$$
\operatorname{GKdim}\big( Y_\mu^\lambda(\bR) / J \big) \leq  \GKdim Y_\mu^\lambda(\bR),
$$
with equality if and only if $J = \{ 0\}$ by \cite[Satz 3.4]{borhokraft}. Also, by \cite[Lemma 4.12]{kamnitzer2022lie},%
$$
2 \GKdim V = \operatorname{GKdim}\big( Y_\mu^\lambda(\bR) / J\big)  
$$
and hence $V$ is faithful over $Y_{\mu}^{\la}(\bR)$ if and only if $2\GKdim V=\GKdim Y_{\mu}^{\la}(\bR)$. Finally,~the embedding 
$Y_\mu^\la(\bR) \subseteq \widetilde{\mathscr{A}}_{\la - \mu}$ gives
\begin{equation*}
\GKdim Y_\mu^\la(\bR) \leq \GKdim \widetilde{\mathscr{A}}_{\la-\mu} \leq 2\operatorname{ht}(\la - \mu)
\end{equation*}
and, since $Y_\mu^\la(\mathbf{R})$ quantizes the generalized affine Grassmannian slice $\overline{\mathcal{W}}{}_\mu^\la$, we get (see \cite[proof of Lemma 2.7]{braverman2016coulomb} and \cite[Lemma 6.5]{krause2000growth})
$$ \GKdim Y_\mu^\la(\bR) \geq \dim \overline{\mathcal{W}}_\mu^\la =  2\operatorname{ht}(\la - \mu).$$
This ends the proof of (i). For part (ii), take $\bR\in \C^{\la}$ containing $k = |\bR|$ distinct integrality classes and let $\fg_{\bR}:=(\fg^{\vee})^{\oplus k}$. By the results of \cite{kamnitzer2022lie}, $K_\C( \mathcal{O}_{sh}^\la(\mathbf{R}))$ is irreducible~as~a $\fg_{\bR}$-module since it is a non-zero submodule of a product of the form $ V(\varpi_{i_1})\otimes \dots \otimes V(\varpi_{i_k}),$
with $V(\varpi_{i_r})$ a fundamental representation for the $r^{\text{th}}$-copy of $\fg^{\vee}$ in the direct sum $\fg_{\bR}$. The surjective map of $\mathfrak{g}_{\bR}$-modules
$$
\textstyle  K_\C \big( \mathcal{O}_{sh}^\la(\mathbf{R}) \big) \twoheadrightarrow \bigoplus_{\mu} K_\C \big( \mathcal{O}_\mu^\la(\bR)_{\topRM} \big)
$$ 
is therefore also injective. Thus $\mathcal{O}_{\mu}^{\la}(\bR) = \mathcal{O}_{\mu}^{\la}(\bR)_{\text{top}}$ for all $\mu\in P^{\vee}$, that is any $V$ in $\mathcal{O}_{\mu}^{\la}(\bR)$ satisfies $\GKdim V=\operatorname{ht}(\la-\mu)$ and part (ii) follows from part (i).
\end{proof}
\begin{Rem}\label{rem:faithful_and_Btop}
For $\mathbf{R} \in \C^\la$, the faithful simple modules of $\mathcal{O}_\mu^\la(\mathbf{R})$ are easy to characterize using crystal combinatorics. Indeed, let again $k$ be the number of distinct integrality~classes in $\bR$ and let $\mathcal{B}(\la, \mathbf{R})_{\topRM}$ be the connected component of $\cB(\la,\mathbf{R})$ generated by the monomial $y_{\bR}$ given in \eqref{eq:yR}. Then $\mathcal{B}(\la, \mathbf{R})$ is a $\fg_{\bR}$-crystal (cf.~Remark \ref{rem:nonIntCryst}) and $\mathcal{B}(\la, \mathbf{R})_{\topRM}$ is the crystal of the simple $\fg_{\bR}$-module $\bigoplus_\mu K_\C \big( \mathcal{O}_\mu^\la(\bR)_{\topRM} \big)$ above. Using this, one shows that, for $\psi\in \mathfrak{r}_{\mu}$,
$$
L(\psi)  \text{ is faithful for } Y_\mu^\la(\mathbf{R}) \ \ \Longleftrightarrow  \ \ \psi \in \mathcal{B}(\la, \mathbf{R})_{\topRM}.
$$
In particular, (the non-integral version of) Lemma \ref{lem:ConvexHullWts} shows that there is always a faithful simple module in $\mathcal{O}_{\mu}^{\la}(\bR)$ when $\cB(\la,\bR)_{\mu}\neq \emptyset$ (i.e.~when the category $\mathcal{O}_{\mu}^{\la}(\bR)$ is non-empty). This generalizes \cite[Theorem B]{webster2020quantum}.
\end{Rem}

\begin{proof}[Proof of Theorem \ref{thm:truncoprod}]
Let 
$$\mathcal{Z}=\{\bR\in \C^{\la}\,|\,\text{every element of } \bR\text{ lies in a distinct integrality class}\}\subseteq \C^{\la},$$ 
and let $\mathcal{Z}' \subseteq \C^{\la_1} \times \C^{\la_2}$ denote the preimage of $\cZ$ under the natural map $\C^{\la_1} \times \C^{\la_2} \rightarrow \C^\la$.  Then both $\mathcal{Z}$ and $\mathcal{Z}'$ are Zariski dense.%

Fix $(\bR_1, \bR_2) \in \mathcal{Z}'$ with corresponding $\bR = \bR_1 \cup \bR_2 \in \mathcal{Z}$.  Take any $\psi \in \cB(\lambda_1, \bR_1)_{\mu_1}$ and $\xi \in \cB(\lambda_2, \bR_2)_{\mu_2}$. Then $\psi \xi \in \cB(\lambda, \bR)_\mu$ by the definition of product monomial crystals.~Also, as $(\bR_1, \bR_2) \in \mathcal{Z}'$, $\psi\in \mathfrak{r}(S)$ and $\xi\in \mathfrak{r}(S')$ for two integrally disconnected subsets $S,S'\subseteq\C$. Thus, Theorem \ref{thm:TensSimpInt} applies in this context and, as $Y_{\mu}$-modules,
$$
L(\psi \xi) \cong L(\psi) \otimes L(\xi).
$$ 
Now, by definition of $\cZ$, the sets of parameters $\mathbf{R}$, $\mathbf{R}_1$ and $\mathbf{R}_2$ all satisfy the hypothesis of part (2) of Theorem \ref{thm:faithfulmodules}. Thus $L(\psi\xi)$ and $L(\psi)\otimes L(\xi)$ are faithful modules~for~the~algebras $Y_{\mu}^{\la}(\bR)$ and $Y_{\mu_1}^{\la_1}(\bR_1) \otimes Y_{\mu_2}^{\la_2}(\bR_2)$, respectively, so that 
\begin{align*}
J_{\mu_1,\mu_2}^{\la_1,\la_2}(\bR_1,\bR_2)&=\Ker((\mathsf{\Pi}_1(\bR_1)\otimes\mathsf{\Pi}_2(\bR_2))\circ\mathsf{\Phi})\\
&= \{x\in Y_{\mu}\otimes \coord^{\la}\,|\,(\mathsf{\Pi}_1(\bR_1)\otimes\mathsf{\Pi}_2(\bR_2))(\mathsf{\Phi}(x)) \text{ acts trivially on } L(\psi)\otimes L(\xi)\}\\
&= \operatorname{Ann}_{Y_\mu \otimes \coord^\lambda} (L(\psi) \otimes L(\xi)) = \operatorname{Ann}_{Y_\mu \otimes \coord^\lambda} L( \psi\xi) = I_\mu^\lambda(\bR).
\end{align*}
Hence Theorem \ref{thm:truncoprod} holds for all $(\mathbf{R}_1, \mathbf{R}_2) \in \cZ'$ %
and using Corollary \ref{cor:reductioncoproductsgeneric} ends the proof.
\end{proof}
\begin{Rem}
\label{rem:injective} By \cite[Theorem 4.7 and Theorem 4.14]{kamnitzer2022hamiltonian}, the shifted coproduct $\Delta_{\mu_1,\mu_2}$ quantizes a multiplication map of schemes 
$$\mathcal{W}_{\mu_1} \times \mathcal{W}_{\mu_2} \rightarrow \mathcal{W}_{\mu},$$ 
while the defining surjection $Y_\mu \twoheadrightarrow Y_\mu^\la(\mathbf{R})$ quantizes the inclusion $\overline{\mathcal{W}}_\mu^\la \subseteq \mathcal{W}_\mu$ of a generalized affine Grassmannian slice. Hence, with $\bR=\bR_1\cup\bR_2$, the map 
\begin{equation}\label{eq:trunccoprod}
Y_{\mu}^{\la}(\bR)\to Y_{\mu_1}^{\la_1}(\mathbf{R}_1) \otimes  Y_{\mu_2}^{\la_2}(\mathbf{R}_2)
\end{equation}
quantizes the multiplication map $\overline{\mathcal{W}}_{\mu_1}^{\la_1} \times \overline{\mathcal{W}}_{\mu_2}^{\la_2}  \rightarrow \overline{\mathcal{W}}_{\mu}^{\la}$ of generalized %
slices given in \cite[\S 2(vii)]{braverman2016coulomb}, which was shown to be dominant in \cite[Proposition 5.7]{krylovperunov2021}. This implies~that~the quantized map \eqref{eq:trunccoprod} is injective (see, e.g., \cite[Lemma 6.8(1)]{kamnitzer2022hamiltonian}).  %
\end{Rem}
\begin{Rem}
\label{rem:truncoprodA}
Suppose that  $\fg$ is of type A. As recalled in Remark \ref{rem:tsyideal},
$$Y_\mu^\la(\bR) = Y_\mu / \langle a_{i,r}\,|\, i \in I\text{ and } r > m_i\rangle$$  
and, using the fact that the elements $a_{i,r}$ are related to quantum minors (see for example \cite[Section 5.3]{kamnitzer2019highest}), one can explicitly compute the elements $\Delta_{\mu_1,\mu_2}(a_{i,r})$ in order to prove, as was done for $\mu_1,\mu_2$ dominant in \cite[Proposition 4.1.13]{weekesthesis}, that the coproduct 
$\Delta_{\mu_1,\mu_2}:Y_{\mu}\to Y_{\mu_1}\otimes Y_{\mu_2}$
induces a map
$$ Y_{\mu_1+\mu_2}^{\la_1+\la_2}(\bR_1\cup \bR_2)\to Y_{\mu_1}^{\la_1}(\bR_1)\otimes Y_{\mu_2}^{\la_2}(\bR_2),$$
even when the hypotheses 
$$\cB(\la_1,\bR_1)_{\mu_1} \neq \emptyset \neq \cB(\la_2,\bR_2)_{\mu_2}$$ 
of Theorem \ref{thm:truncoprod} do not hold. This generalizes the recent result \cite[Theorem 21]{milot2025} obtained by a similar approach for $\fg=\mathfrak{sl}_2$ and general~$\mu\in P^{\vee}$.
\end{Rem}

\begin{Rem}
Fix $\psi_1,\psi_2\in \mathfrak{r}$. Let $\cB(\la_1,\bR_1)$ and $\cB(\la_2,\bR_2)$ be the unique minimal product monomial crystals that contain $\psi_1$ and $\psi_2$, respectively (cf.~Corollary \ref{cor:minimal-pmc}). Theorem \ref{thm:truncoprod} and Corollary \ref{cor:tensor_product_decomp_in_Osh} imply that the highest $\ell$-weights of the simple factors of the tensor product $L(\psi_1)\otimes L(\psi_2)$ must all lie in the subset $$\{\psi\in \cB(\la_1+\la_2, \bR_1 \cup \bR_2)\,|\,\psi\preceq \psi_1\psi_2\}\subseteq \cB(\la_1+\la_2, \bR_1 \cup \bR_2).$$
\end{Rem}
\subsection{Corollaries}\label{sec:CorTSC} Fix $\la_1,\la_2\in P^{\vee}_+$ along with sets of parameters $\bR_1\in \C^{\lambda_1}$ and $\bR_2\in \C^{\lambda_2}$. Recall that Theorem \ref{th:TFAE} gives three statements that are equivalent to the containment 
$$\cB(\lambda_1,\bR_1)\subseteq\cB(\lambda_2,\bR_2).$$
Below, using techniques of Section \ref{sec:TSC}, we establish yet another such equivalent statement.

\begin{Theorem}\label{thm:TFAEsurjection} Take $\mu\in P^{\vee}$ with $\cB(\la_1,\bR_1)_{\mu}\neq \emptyset$. Then $\cB(\la_1,\bR_1)_{\mu}\subseteq \cB(\la_2,\bR_2)_{\mu}$ if~and only if there is a surjective algebra map $Y_\mu^{\lambda_2}(\bR_2)\twoheadrightarrow Y_\mu^{\lambda_1}(\bR_1)$ %
completing the diagram%
:
{\small
$$\begin{tikzcd}[row sep = 1.25em]
Y_{\mu} \ar[d,"\Phi_{\mu}^{\la_2}(\bR_2)",swap,two heads] \ar[dr,"\Phi_{\mu}^{\la_1}(\bR_1)",two heads] \\
Y^{\la_2}_{\mu}(\bR_2) \ar[r, two heads, dashed] & Y_{\mu}^{\lambda_1}(\bR_1)
\end{tikzcd}$$}
\end{Theorem}

\begin{proof} The ``if'' direction follows from Theorem \ref{th:descend}. For the converse, use Remark \ref{rem:faithful_and_Btop}~to fix a faithful simple module $L(\psi)$ in $\mathcal{O}_\mu^{\lambda_1}(\bR_1)$. Then $\psi \in \cB(\la_1,\bR_1)_{\mu}\subseteq \cB(\la_2,\bR_2)_{\mu}$~and~$L(\psi)$ is also a simple object of $\mathcal{O}_{\mu}^{\la_2}(\bR_2)$. Hence, denoting by 
$$\rho_1:Y_{\mu}^{\la_1}(\bR_1)\to \End_{\C}(L(\psi)) \text{ and }\rho_2:Y_{\mu}^{\la_2}(\bR_2)\to \End_{\C}(L(\psi)) $$
the representation morphisms, we have that $\rho_1$ is injective and the diagram 
{\small \begin{equation*}
\begin{tikzcd}[row sep = 1.25em]
Y_\mu \arrow[r,"\Phi_\mu^{\lambda_2}(\bR_2)"]\arrow[d,"\Phi_\mu^{\lambda_1}(\bR_1)"'] & Y_\mu^{\lambda_2}(\bR_2) \arrow[d,"\rho_2"]\\
Y_\mu^{\lambda_1}(\bR_1) \arrow[r,hook,"\rho_1"] & \End_\C(L(\psi))
\end{tikzcd}
\end{equation*}}\noindent
commutes. In particular, for $x\in \Ker\Phi_\mu^{\lambda_2}(\bR_2)$, we get $\Phi_\mu^{\lambda_1}(\bR_1)(x)=0$ as $\rho_1$ is injective~with
\begin{equation*}
\smash{0=(\rho_2\circ\Phi_\mu^{\lambda_2}(\bR_2))(x)=(\rho_1\circ\Phi_\mu^{\lambda_1}(\bR_1))(x).}
\end{equation*}
This implies that $\Ker\Phi_\mu^{\lambda_2}(\bR_2)\subseteq\Ker\Phi_\mu^{\lambda_1}(\bR_1)$, proving the claim.
\end{proof}
\begin{Rem} A more immediate version of Theorem \ref{thm:TFAEsurjection}, using the alternative definition of truncation given in Remark \ref{rem:tsyideal}, can be found in \cite[Proposition 10.7]{hernandez2026crossroad}.
\end{Rem}
We finally want to end Section \ref{sec:TensorResults} by tying together the blocks appearing in Theorem~\ref{thm:BlockDec} and the ``truncated categories'' $\mathcal{O}_{\mu}^{\la}(\bR)$. For this, note that Theorem \ref{thm:truncoprod}, Theorem \ref{th:descend} and Proposition \ref{prop:WeylStandards} imply the following corollary (similar to \cite[Theorem 8.4]{hernandez2024shifted}):
\begin{Corollary}\label{cor:highestTrunc} All highest $\ell$-weight modules in $\mathcal{O}_{sh}$ 
descend to a truncation $Y_{\mu}^{\la}(\bR)$.
\end{Corollary}
Fix $\psi\in \mathfrak{r}$ and let $V$ be a module of co-highest $\ell$-weight $\psi$. Using the functor $\textbf{tr}$~of~Section \ref{sec:Auto} on Proposition \ref{prop:WeylStandards} shows that $V$ can be realized as a subquotient of a tensor product of two simple objects of $\mathcal{O}_{sh}$. Using Theorem \ref{thm:truncoprod} and Theorem \ref{th:descend} thus %
gives%
: %
\begin{Corollary}\label{cor:cohighestTrunc} All co-highest $\ell$-weight modules in $\mathcal{O}_{sh}$ 
descend to a truncation $Y_{\mu}^{\la}(\bR)$.
\end{Corollary}
Now, for $\psi,\xi\in \mathfrak{r}$ with $\psi\neq \xi$, Theorem \ref{thm:Ext1Block} shows that there can be non-trivial extensions of $L(\psi)$ by $L(\xi)$ in $\mathcal{O}_{sh}$ only if $\psi\preceq \xi$ or $\xi\preceq \psi$%
. Moreover, if $\psi\preceq \xi$, the proof of Theorem~\ref{thm:Ext1Block} shows that all the non-trivial extensions in $\Ext^1_{\mathcal{O}_{sh}}(L(\psi),L(\xi))$ are of highest $\ell$-weight, and, via $\mathbf{tr}$, we get that, dually, all non-trivial extensions of $L(\psi)$ by $L(\xi)$ in $\mathcal{O}_{sh}$ are of co-highest $\ell$-weight when $\xi\preceq \psi$. This gives the corollary below.
\begin{Corollary}\label{cor:ExttoTrunc} Fix $\psi,\xi\in \mathfrak{r}$ with $\psi\neq \xi$. Then, all non-trivial extensions of $L(\psi)$ by $L(\xi)$ in $\mathcal{O}_{sh}$ descend to a truncation $Y_{\mu}^{\lambda}(\mathbf{R})$.
\end{Corollary}
Recall from Remark \ref{rem:Blockawt} that the partial order $\preceq$ and the group of Laurent monomials~$\mathcal{A}$ are %
related to the map $\awt$ and the subgroup $\Gamma_{\C}\subseteq \cB_{\C}$ of Section \ref{sec:awt}. In particular, fixing $\psi,\xi\in \mathfrak{r}_{\mu}$ with $\psi \in \xi\mathcal{A}$, Lemma \ref{lem:InductiveSystemC} states that both simple modules $L(\psi)$ and $L(\xi)$ descend to a common truncation $Y_{\mu}^{\la}(\bR)$. There is hence in this case a natural inclusion of groups:
\begin{equation}\label{eq:ExtTruncorno}
\Ext^1_{\mathcal{O}_{\mu}^{\lambda}(\mathbf{R})}(L(\psi),L(\xi))\subseteq \Ext^1_{\mathcal{O}_{\mu}}(L(\psi),L(\xi)).
\end{equation}
One may wonder whether there is some special choice of $(\lambda,\mathbf{R})$ for which \eqref{eq:ExtTruncorno} is an equality. The answer is unfortunately no, as shown by the example below.
\begin{Example}\label{ex:SelfExt} Let $Y=Y(\mathfrak{sl}_2)$ and consider the $4$-dimensional $Y$-module $V$ with action
\begin{equation*}
[e(u)] = \adjustbox{scale=0.78}{$
\begin{pmatrix}
0 & \frac{2}{u} & 0 & 0\\
0 & 0 & 0 & 0\\
0 & \frac{4}{u^2} & 0 & \frac{2}{u}\\
0 & 0 & 0 & 0
\end{pmatrix}$}, \ [f(u)] = \adjustbox{scale=0.78}{$
\begin{pmatrix}
0 & 0 & 0 & 0\\
\frac{2}{u} & 0 & 0 & 0\\
0 & 0 & 0 & 0\\
\frac{4}{u^2} & 0 & \frac{2}{u} & 0
\end{pmatrix}
$}\, \text{ and }\, [h(u)] = \adjustbox{scale=0.78}{$
\begin{pmatrix}
\frac{u+2}{u} & 0 & 0 & 0\\
0 & \frac{u-2}{u} & 0 & 0\\
\frac{4}{u^2} & 0 & \frac{u+2}{u} & 0\\
0 & -\frac{4}{u^2} & 0 & \frac{u-2}{u}
\end{pmatrix}$}.
\end{equation*}
It is not hard to prove that $V$ is a well-defined $Y$-module, and, more precisely, that it is the unique (up to isomorphism\footnote{See \cite{grigorev2025representations} for a similar (unique) extension in the case of the quantum affine algebra $U_q(\widehat{\mathfrak{sl}}_2)$.%
}) non-trivial self-extension of the %
fundamental module~$L(\mathsf{Y}_{1,0})$ of Example \ref{ex:KR}. Also, as we will show in Appendix \ref{sec:Inflations}, $V$ does not descend to any truncation $\smash{Y_0^{\la}(\bR)}$ with $\la\in \smash{P_+^{\vee}}$ and $\bR\in \C^{\la}$ so that\vspace*{-0.2mm}
$$ {0=\Ext^1_{\mathcal{O}_{0}^{\lambda}(\mathbf{R})}(L(\mathsf{Y}_{1,0}),L(\mathsf{Y}_{1,0}))\not\simeq \Ext^1_{\mathcal{O}_{0}}(L(\mathsf{Y}_{1,0}),L(\mathsf{Y}_{1,0}))\simeq \C}$$
for all $\la$ and $\bR$. Hence \eqref{eq:ExtTruncorno} can never be an equality in this case. 
\end{Example}\newpage
Fortunately, as we show below, the problems arising in Example \ref{ex:SelfExt} do not arise (under a mild hypothesis) for extensions between non-isomorphic simple modules. Recall \eqref{eq:yR} and define, for $\la\in P_{+}^{\vee}$ and $\bR=(R_i)_{i\in I}\in \C^{\la}$,
$$ \textstyle \mathsf{\Psi}_{\bR} = \prod_{i\in I}\prod_{a\in R_i}\mathsf{\Psi}_{i,a}.$$
 \begin{Theorem}\label{thm:ExtTrunc}
Fix $\psi,\xi\in \mathfrak{r}_{\mu}$ with $\psi\neq\xi$ and for which 
$$0\neq \dim_{\C} \Ext^1_{\mathcal{O}_{sh}}(L(\psi),L(\xi))<\infty.$$ 
Then \eqref{eq:ExtTruncorno} is an isomorphism for some $\lambda$ and $\mathbf{R}$.
\end{Theorem}
\begin{proof}
Fix non-trivial extensions $E_1,\dots,E_N$ of $L(\psi)$ by $L(\xi)$ such that 
$$\Ext^1_{\mathcal{O}_{\mu}}(L(\psi),L(\xi))=\spa_{\C}(\{E_1,\dots,E_N\}).$$
By Corollary \ref{cor:ExttoTrunc}, each $E_k$ ($1\leq k\leq N$) passes to a truncation $Y_{\mu}^{\lambda_k}(\mathbf{R}_k)$ with $\la_k\in P^{\vee}_+$ and $\bR_k\in \C^{\la_k}$. Thus, $\psi \in \bigcap_{1\leq k\leq N}\cB(\la_k,\bR_k)$ and, by the non-integral version of Lemma \ref{th:qm2},
$$\textstyle \psi \in \bigcap_{1\leq k\leq N}\mathsf{\Psi_{\bR_k}}\mathcal{A}_+^{-1}$$
so that Lemma \ref{lem:InductiveSystemC} shows that there exists $\la\in P_+^{\vee}$ and $\bR\in \C^{\la}$ for which
$$\textstyle  \bigcup_{1\leq k\leq N}\cB(\la_k,\bR_k)\subseteq \cB(\la,\bR).$$
Hence, by Theorem \ref{thm:TFAEsurjection}, we have, for all $1\leq k\leq N$, a surjective algebra morphism 
$$ Y_{\mu}^{\la}(\bR)\twoheadrightarrow Y_{\mu}^{\la_k}(\bR_k),$$
and it follows that the extensions $E_1,\dots,E_N$ all lie in the category $\mathcal{O}_{\mu}^{\la}(\bR)$. Thus
$$\Ext^1_{\mathcal{O}_{\mu}}(L(\psi),L(\xi))=\spa_{\C}(\{E_1,\dots,E_N\})\subseteq \Ext^1_{\mathcal{O}_{\mu}^{\la}(\bR)}(L(\psi),L(\xi))$$
and \eqref{eq:ExtTruncorno} is an equality for our choice of $(\la,\bR)$. This ends the proof.
\end{proof}
To conclude, for different $\psi,\xi\in \mathfrak{r}$ satisfying the condition $\dim_{\C}\Ext^1_{\mathcal{O}_{sh}}(L(\psi),L(\xi))<\infty$ (which we believe always holds), all non-trivial extensions of $L(\psi)$ by $L(\xi)$ in $\mathcal{O}_{sh}$ (or in~one of the blocks ${}_{\tau}\mathcal{O}_{sh}$ of Theorem \ref{thm:BlockDec}) can be computed in a ``big enough'' truncation.

\section*{Interlude}\label{sec:interlude}

Starting next section, we make some significant notational changes.\medskip

\begin{itemize}
\item We replace the Lie algebra $\fg$ by its Langlands dual $\fg^\vee$ (and consider shifted Yangians for $\fg^{\vee}$). This replaces all weights by coweights, all roots by coroots, and vice-versa.

\item We replace the category $\mathcal{O}_{sh}$ by the full subcategory of finite-length objects in the integral category $\mathcal{O}_{sh}(0)$ of Definition \ref{def:intCatO}%
. Similarly, we replace $\mathfrak{r}$ by the group~$\mathfrak{r}(0)$ of integral $\ell$-weights, which is by definition isomorphic to the monomial crystal $\cB$, and replace $\mathcal{A}$ by $\mathcal{A}\cap\mathfrak{r}(0)$. Recall that %
$\mathcal{O}_{sh}(0)$
is monoidal by Corollary \ref{cor:tensor_product_decomp_in_Osh} (since we work under the assumption that Conjecture \ref{conj:Associators} holds).

\item All sets of parameters considered will be integral. The collection of (integral) sets of parameters~of level $\lambda\in P$
will be denoted $\Z^{\la}$.

\item We replace the set of GT-weights $\Lambda_{\nu}=\prod_{i\in I}\C^{m_i}/\Sigma_{m_i}$ (where $\nu=\sum_{i\in I}m_i\alpha_i\in Q_+$) by its integral counterpart, i.e.~by the subset 
$$\{\bS=(S_i)_{i\in I}\in \Lambda_{\nu}\,|\,S_i\subseteq \overline{i}+2\Z \text{ for each }i\in I\}.$$
\end{itemize}
In particular, since the category $\mathcal{O}_{sh}$ is (now) of finite-length, a natural $\Z$-basis
~for~the~ring $K_0(\mathcal{O}_{sh})$ is given by the set of classes of simple objects, i.e.~by the $[L(\psi)]$'s with $\psi\in \mathfrak{r}\simeq \cB$. In addition, by Theorem \ref{thm:TensSimpInt}, a simple object in the general (i.e.~non-integral) version~of~the category $\mathcal{O}_{sh}$ of the previous pages can always be factorized as a tensor product of spectral shifts of simple objects with integral highest $\ell$-weights. We hence do not lose any generality in considering the Grothendieck ring of our integral category $\mathcal{O}_{sh}$ instead of that of (the~full subcategory of finite-length objects in) the non-integral version of this category%
.
					
\section{KLR algebras and their variants}\label{sec:KLR}

We briefly recall the framework of KLR algebras due independently to \cite{khovanov2009diagrammatic,khovanov2011diagrammatic}~and \cite{rouquier20082}. These diagrammatic algebras provide a categorification of the negative-half~$U_q(\mathfrak{n}_-)$ of the quantum group $U_q(\mathfrak{g})$ (for $q$ a formal parameter associated to a grading functor) and their cyclotomic quotients categorify irreducible representations of $\mathfrak{g}$ by \cite{kang2012categorification} or \cite{webster2017knot}. We will not give proper definitions here and refer the reader to the above-mentioned papers for a detailed introduction. Throughout, we work with non-graded $\C$-algebras.

\subsection{KLR algebras} \label{subsec:KLR_algebras_def}
For $\nu\in Q_+$, let $R_\nu$ be the \textit{KLR algebra} for $\fg$ of height $\nu$. This~algebra comes equipped with %
 mutually orthogonal idempotents $\{e(\bi)\}_{\bi\in \Seq_{\nu}}%
 $, where (for $m=\rht(\nu)$)
\begin{equation*}
\textstyle\smash{\Seq_\nu=\{(i_1,\dots, i_m)\in I^m\;;\; \nu=\sum_{j=1}^m \alpha_{i_j}\},}
\end{equation*}
and is unital with unit %
$1_\nu=\sum_{\bi\in \Seq_\nu} e(\bi)$. For $\nu_1,\nu_2\in Q_+$ such that $\nu=\nu_1+\nu_2$, there is a non-unital morphism of algebras $\iota_{\nu_1,\nu_2}:R_{\nu_1}\otimes R_{\nu_2}\to R_\nu$ which produces an 
induction and a 
restriction functor
\begin{equation*}
\circ:R_{\nu_1}\fmod\times R_{\nu_2}\fmod\to R_{\nu}\fmod\,\text{ and }\,\Res_{\nu_1,\nu_2}:R_{\nu}\fmod\to R_{\nu_1}\fmod\times R_{\nu_2}\fmod,
\end{equation*}
where $\mathrm{fmod}$ stands for the category of finite-dimensional (left-)modules. 
In particular, for $i\in I$ such that $\nu-\alpha_i\in Q_+$, we get a functor
$$\smash{\cE_i^\nu:R_{\nu}\fmod\to R_{\nu-\alpha_i}\fmod }$$
given by the composition of $\Res_{\alpha_i,\nu-\alpha_i}$ with the pullback by the %
canonical algebra inclusion $R_{\nu-\alpha_i}\hookrightarrow  R_{\alpha_i}\otimes R_{\nu-\alpha_i}$. Extend this construction by setting $\cE_i^\nu=0$ for $\nu-\alpha_i\not\in Q_+$ and~let $R\fmod=\bigoplus_{\nu\in Q_+} R_\nu\fmod$. Then the functors
$$\textstyle \smash{\{\cE_i=\bigoplus_{\nu\in Q_+} \cE_i^\nu:R\fmod\to R\fmod\}_{i\in I}},$$
endow $R\fmod$ with a categorical left-action of $\mathfrak{n}$ which, on the ring\footnote{The multiplication in $K_{0}(R\fmod)$ %
comes from the induction functors. Also, to be precise, one should replace $R\fmod$ in the isomorphism with $\C[N]$ (and in Proposition \ref{thm:character_KLR} below) by its subcategory of \textit{nilpotent modules} (see \cite[Section 2.7]{lauda2011crystals} and \cite[Lemma 3.13]{kamnitzer2019category}). This will however not cause issues for us as all finite-dimensional modules over the cyclotomic quotients we study in the next subsections are nilpotent.} $K_{\C}(R\fmod)\simeq \C[N]$, matches the left-$\mathfrak{n}$-action by differential operators \cite{rouquier20082,khovanov2009diagrammatic}.  
\medskip\par
Let $\KLRchar_\nu=\Z^{\Seq_\nu}$ be the set of $\Z$-valued functions with domain $\Seq_{\nu}$. Endow 
$\textstyle \KLRchar=\bigoplus_{\nu\in Q_+}\KLRchar_\nu$
with a $Q_+$-graded ring structure where addition is defined point-wise and where the product 
$$\smash{\shuffle : \KLRchar_{\nu_1}\otimes \KLRchar_{\nu_2}\to \KLRchar_{\nu_1+\nu_2}}$$ 
is \textit{shuffle product}. More precisely, if $f\in \KLRchar_{\nu_1}$ and $g\in \KLRchar_{\nu_2}$, then, for $\bi\in \Seq_{\nu_1+\nu_2}$,
\begin{equation}\label{eq:def_of_shuffle_prod}
\smash{\textstyle(f\shuffle g)(\bi)=\sum_{\mathbf{j}_1,\,\mathbf{j}_2} f(\mathbf{j}_1)g(\mathbf{j}_2)}
\end{equation}
where the sum runs over all the ways we can write $\bi$ as a shuffle of $\mathbf{j}_1\in \Seq_{\nu_1}$ and $\mathbf{j}_2\in \Seq_{\nu_2}$ (i.e.~as an interleaving of such $\mathbf{j}_1$ and $\mathbf{j}_2$ that preserves the internal order of both sequences). 
\begin{Def} Fix an $R_{\nu}$-module $M$. The \textit{character} of $M$ is the element of $\KLRchar_{\nu}$ given by 
\begin{equation*}
\smash{\textstyle\ch(M)=\sum_{\bi\in \Seq_\nu} \dim(e(\bi)M) [\bi]},
\end{equation*} 
where $[i]\in \mathcal{S}_{\nu}$ is the map given by $[\bi](\bi')=\delta_{\bi,\bi'}$ for all $\bi'\in \Seq_{\nu}$.
\end{Def}
\begin{Proposition}[{\cite[Theorem~3.17]{khovanov2009diagrammatic}}]\label{thm:character_KLR}
The map $\ch:K_0(R\fmod)\to \KLRchar$ is an injective ring homomorphism. 
\end{Proposition}
\subsection{Cyclotomic KLR algebras} For every $\lambda=\sum_{i\in I}\la_{i}\varpi_i\in P_+$ and $\nu\in Q_+$, let $I^\lambda_\nu$ be the (two-sided) ideal of $R_\nu$ generated by elements of the form
\begin{equation*}
x^{\lambda_{i_m}}_me(\bi)=
\raisebox{-2.2em}{
\begin{tikzpicture}[every node/.style={font=\scriptsize}]
\draw (0,0.8) -- (0,0) node[below] {$i_1$};
\draw (0.5,0.8) -- (0.5,0) node[below] {$i_2$};
\node at (1,0.4) {$\dots$};
\draw (1.4,0.8) -- (1.4,0) node[below] {$i_{m-1}$};
\draw (1.9,0.8) -- (1.9,0) node[below] {$\ \ i_m$};
\fill (1.9,0.4) circle (0.07) node[right] {$\lambda_{i_m}$};
\end{tikzpicture}
}
\end{equation*}
where the notations are borrowed from \cite{khovanov2009diagrammatic} and $\bi=(i_1,\dots,i_m)$. Set also $R^\lambda_{\nu}=R_\nu/I_\nu^\lambda$.

The (finite-dimensional!) quotient $R^\lambda_{\nu}$ is called the \textit{cyclotomic KLR algebra} of level $\lambda$. It was introduced in \cite{khovanov2009diagrammatic}, where the following statement (now shown for~all~symmetrizable Kac--Moody algebras) was conjectured:

\begin{Theorem}[{\cite{kang2012categorification},\cite{webster2017knot}}]\label{thm:cyclotomic_categorification}
The functors $\{\cE_i=\bigoplus_{\nu\in Q_+} \cE_i^{\nu}$, $\cF_i=\bigoplus_{\nu\in Q_+} \cF_i^\nu\}_{i\in I}$ with
\begin{gather*}
\cF_i^\nu=R_{\nu+\alpha_i}^\lambda \otimes_{R_{\alpha_i}\otimes R_{\nu}} (R_{\alpha_i}\otimes (\trou)):R_{\nu}\fmod\to R_{\nu+\alpha_i}^{\la}\fmod\subseteq R_{\nu+\alpha_i}\fmod
\end{gather*}
induce endofunctors of $R^{\la}\fmod=\bigoplus_{\nu\in Q_+}R^{\la}_{\nu}\fmod$ which categorify the simple left-module $V(\lambda)$ of $\fg$ of highest-weight $\la$. Also, the $\fg$-equivariant isomorphism 
\begin{equation}\label{eq:CatKLR}
K_{\C}(R^{\la}\fmod)\simeq V(\la)
\end{equation}
coming from this categorification identifies the $\mu$-weight space $V(\lambda)_\mu$ with $K_\C(R^\lambda_{\lambda-\mu}\fmod)$, and the functors $\cE_i$ and $\cF_i$ with the Chevalley generators $e_i$ and $f_i$ of $\fg$ (respectively). 
\end{Theorem}

The natural map $K_0(R^{\la}\fmod)\hookrightarrow K_0(R\fmod)$ can be composed with the character~map of Theorem \ref{thm:character_KLR} to give an injective linear map $\ch:K_0(R^\lambda\fmod)\to \KLRchar$. We will need the following lemma, which is an easy consequence of the definition of the endofunctors $\{\mathcal{E}_i\}_{i\in I}$.

\begin{Lemma}\label{lem:chEiCycKLR} 
Fix $M$ in $R^{\la}_{\nu}\fmod$ with $i\in I$. Then $e(i,\bi)M=e(\bi)\cE_i(M)$ for all $\bi\in \Seq_{\nu-\alpha_i}$, where $e(i,\bi)$ is the idempotent associated to the sequence where $i$ has been added before $\bi$.
\end{Lemma}

In \cite{lusztig1990canonicalI}, Lusztig constructed the dual canonical basis of an arbitrary simple representation of $\fg$ using properties of the quantum group $U_q(\mathfrak{n}_-)$. Since then, Lusztig's construction has been related to KLR algebras by the result below, which we will use in Section \ref{sec:duality}.~Note that the assumption that $\fg$ is simply-laced is used here.

\begin{Theorem}[{\cite{varagnolo2011canonical}}]\label{thm:dual_canonical_KLR}
Under the isomorphism \eqref{eq:CatKLR}, the dual canonical basis of $V(\lambda)$ is identified with the basis given by classes of simple objects in $K_{\C}(R^{\la}\fmod)$.
\end{Theorem}

\subsection{KLRW algebras}\label{sec:def_KLRW}

In \cite{webster2017knot}, %
diagrammatic algebras generalizing the classical KLR algebras were introduced to categorify tensor products of irreducible representations. They depend on a list of dominant weights 
$$\ulambda=(\lambda^{(1)},\dots, \lambda^{(\ell)}),$$
together with an element of the positive root cone $\nu\in Q_+$. We denote the KLRW algebra associated to $\ulambda,\nu$ by $\smash{\tilde{T}{}^{\ulambda}_\nu}$. It is a unital algebra whose presentation is usually given in terms of ``braid-like'' decorated diagrams that contain both \textbf{black} strands, labelled by simple~roots (and allowed to cross), and \textcolor{red}{\textbf{red}} strands, which are not allowed to cross and~are labelled by the dominant weights of the list $\ulambda$ (in the order prescribed by that list). 
\medskip\par
The algebra $\smash{\tilde{T}{}^{\ulambda}_\nu}$ comes equipped, similarly to $R_{\nu}$, with mutually orthogonal idempotents $e(\bi,\kappa)%
$, labelled by pairs $(\bi,\kappa)$ where $\bi\in \Seq_{\nu}$ and where
$$\kappa:\{1,\dots,\ell\}\rightarrow \{0,\dots,m\}$$
is a weakly-increasing function. Moreover, $\smash{\tilde{T}{}^{\ulambda}_\nu}$ is unital with unit $1_{\nu}^{\ulambda}= \sum_{\bi,\kappa} e(\bi,\kappa)$, where the sum runs over all sequences and charges.
\begin{Rem}\label{rem:diagKLRW} In the above presentation of $\smash{\tilde{T}{}^{\ulambda}_\nu}$, the idempotents $e(\bi,\kappa)$ are represented by ``straight-line diagrams'' with the red strand labelled $\lambda^{(k)}$ positioned between the $\kappa(k)^\text{th}$ and $(\kappa(k)+1)^\text{th}$ black strands (i.e.~to the left/right of all black strands if $\kappa(k)=0$/$\kappa(k)=m$).
\end{Rem}

Let $J^{\ulambda}_\nu\subseteq \smash{\tilde{T}{}^{\ulambda}_\nu}$ be the (two-sided) ideal generated by all idempotents $e(\bi,\kappa)$ with $\kappa(\ell)<m$. Then the (finite-dimensional!) quotient 
$$T^{\ulambda}_\nu=\smash{\tilde{T}{}^{\ulambda}_\nu/J^{\ulambda}_\nu}$$ 
is called the \textit{cyclotomic KLRW algebra} (or also sometimes the ``\textit{steadied quotient}'') \cite{webster2017knot}. Similarly to the  KLR-case, endofunctors $\{\mathcal{E}_i\}_{i\in I}$ of the category 
\begin{equation}\label{eq:Tlambdafmod}
\textstyle T^{\ulambda}_{\textcolor{white}{;}}\fmod:=\bigoplus_{\nu\in Q_+} T^{\ulambda}_\nu\fmod
\end{equation}
can be constructed using non-unital morphisms of algebras
\begin{equation}\label{eq:non_unital_map_KLRW}
\iota^{\ulambda}_{\nu_1,\nu_2}:R_{\nu_1} \otimes \smash{\tilde{T}{}^{\ulambda}_{\nu_2}}\to \smash{\tilde{T}{}^{\ulambda}_{\nu_1+\nu_2}}
\end{equation}
(see \cite{webster2017knot}). One can also define endofunctors $\{\cF_i=\bigoplus_{\nu
} \cF_i^{\ulambda,\nu}\}_{i\in I}$ of this category via 
\begin{center}
$\cF_i^{\ulambda,\nu}:=T_{\nu+\alpha_i}^{\ulambda} \otimes_{R_{\alpha_i}\otimes \tilde{T}_{\nu}^{\ulambda}} (R_{\alpha_i}\otimes (\trou)).$
\end{center}
\begin{Theorem}[{\cite[Theorem~B]{webster2017knot}}]\label{thm:cyclotomic_categorification_tensorprod}
The category $T^{\ulambda}\fmod$ with the functors $\{\cE_i,\cF_i\}_{i\in  I}$ categorify the left $\fg$-module $V(\ulambda)=V(\la^{(1)})\otimes \dots \otimes V(\la^{(\ell)})$. In addition, the isomorphism 
\begin{equation}\label{eq:CatKLRW}
K_{\C}(T^{\ulambda}\fmod)\simeq V(\ulambda)
\end{equation}
coming from this categorification identifies the $\mu$-weight space $V(\ulambda)_\mu$ with $K_\C(T^{\ulambda}_{\lambda-\mu}\fmod)$, and the functors $\cE_i$ and $\cF_i$ with the Chevalley generators $e_i$ and $f_i$ of $\fg$. 
\end{Theorem}
\begin{Rem}\label{rem:KLRWtoCyclo} By \cite[Theorem~4.18]{webster2017knot}, for $\ulambda=(\lambda)$, we have $T^{(\lambda)}_\nu\simeq R^\lambda_\nu$ as algebras %
 and the above result recovers Theorem \ref{thm:cyclotomic_categorification}. 
\end{Rem}
Recall the construction of the standardization functors of \cite[Chapter~5]{webster2017knot} which are exact functors 
\begin{equation*}
\bbS^{\ulambda}:R^{\lambda^{(1)}}\fmod\times\dots \times R^{\lambda^{(\ell)}}\fmod \to T^{\ulambda}\fmod
\end{equation*} 
given via tensor product with a distinguished bimodule. In particular, recall that $\bbS^{\ulambda}$ sends tensor products of simple objects to \textit{standard modules} so that the corresponding map 
\begin{equation}\label{eq:SK0}
 [\bbS^{\ulambda}]:K_{0}(R^{\lambda^{(1)}}\fmod)\otimes\dots\otimes K_{0}(R^{\lambda^{(\ell)}}\fmod)\to K_0(T^{\ulambda}_{}\fmod)
\end{equation}
induces an isomorphism of $\fg$-modules that recovers \eqref{eq:CatKLRW} after using the isomorphisms given by \eqref{eq:CatKLR}. Consider %
the dual canonical basis given in  \cite[Chapter~24]{lusztig1993introduction} for tensor products of irreducible representations (see also Section \ref{subsec:effect_of_inv_on_basis}). Then %
we have (under our assumption that $\fg$ is finite-dimensional and simply-laced):

\begin{Theorem}[{\cite[Theorem~A]{webster2015canonical}}]\label{thm:dual_canonical_tensorprod}
The isomorphism \eqref{eq:CatKLRW} (equivalently, the map $[\bbS^{\ulambda}]$) identifies the dual canonical basis of $V(\ulambda)$ with the basis of simple classes in $K_{\C}(T^{\ulambda}_{}\fmod)$.
\end{Theorem}
\subsection{Characters for KLRW algebra modules}\label{subsec:character_KLRW_alg} We will use a notion of characters for finite-dimensional modules over KLRW algebras similar to the one appearing in Section~\ref{subsec:KLR_algebras_def}. Analogous techniques 
were used in \cite{silverthorne2024gelfand} to study Gelfand--Tsetlin modules for $U(\mathfrak{gl}_n)$.

Call \textbf{black} (or \textbf{\textcolor{red}{red}}) labels the elements of $I$ (or $\{\redbold{1},\dots,\redboldell\}$) in the alphabet $I\sqcup \{\redbold{1},\dots,\redboldell\}$. Also let $\Seq_\nu^{\ulambda}$ be the set of sequences in the alphabet $I\sqcup \{\redbold{1},\dots,\redboldell\}$ for which
\begin{itemize}
\item[(1)] the black labels sum to $\nu$, and
\item[(2)] the sequence reduces to $(\redbold{1},\dots,\redboldell)$ after forgetting the black labels.
\end{itemize}

Note that Remark \ref{rem:diagKLRW} gives a natural bijection between the idempotents $e(\bi,\kappa)\in \smash{\tilde{T}{}^{\ulambda}_\nu}$ and elements of $\Seq_{\nu}^{\ulambda}$. %
We abuse notation and write $(\bi,\kappa)$ for the sequence associated~to~$e(\bi,\kappa)$.
\begin{Def} Let $\smash{\KLRchar_\nu^{\ulambda}}=\Z^{\Seq_{\nu}^{\ulambda}}$ and fix a finite-dimensional $\smash{\tilde{T}{}^{\ulambda}_\nu}$-module $M$. The \textit{character} of $M$ is the element of $\smash{\KLRchar_\nu^{\ulambda}}$ defined by 
\begin{equation*}
\textstyle\ch(M)=\sum_{(\bi,\kappa)\in\Seq_\nu^{\ulambda}} \dim(e(\bi,\kappa)M) [(\bi,\kappa)]
\end{equation*}
where $[(\bi,\kappa)]\in \smash{\KLRchar_\nu^{\ulambda}}$ is the map defined by $[(\bi,\kappa)](\bi',\kappa')=\delta_{(\bi,\kappa),(\bi',\kappa')}$ for all $(\bi',\kappa')\in \Seq_{\nu}^{\ulambda}$.
\end{Def}
\begin{Proposition}\label{prop:ch_is_hom_KLRW}
The map $\ch:K_0(\smash{\tilde{T}{}^{\ulambda}_\nu}\fmod)\to \smash{\KLRchar_\nu^{\ulambda}}$ is an injective $\Z$-module morphism.
\end{Proposition}
\begin{proof}
This follows from the same techniques as used in the proof of \cite[Theorem~3.17]{khovanov2009diagrammatic} (i.e. from properties of %
the functors $\cE_i$, %
see \cite[Theorem~5.3.1]{kleshchev2005linear}, \cite[Section~5.5]{vazirani1999irreducible}).
\end{proof}
Let $\ulambda=(\lambda^{(1)},\dots, \lambda^{(\ell)})$ and $\umu=(\mu^{(1)},\dots, \mu^{(k)})$ be two lists of dominant weights and fix $w\in \Sigma_{\ell+k}$ a shuffle of the corresponding red strands, i.e. a permutation such that 
\begin{center}
$w(1)<\dots<w(\ell)\hspace{1em}\text{and}\hspace{1em}w(\ell+1)<\dots<w(\ell+k).$
\end{center}
Denote by $w(\ulambda,\umu)$ the list in $P_+$ obtained by applying $w$ to the concatenation of the lists~$\ulambda$ and $\umu$ (in this order). Also, let
\begin{equation}\label{eq:wShuffle}
\shuffle_w : \KLRchar_{\nu_1}^{\ulambda}\otimes\KLRchar_{\nu_2}^{\umu}\to \KLRchar_{\nu_1+\nu_2}^{w(\ulambda,\umu)}
\end{equation}
be the shuffle product as in \eqref{eq:def_of_shuffle_prod}, i.e. if $f\in \KLRchar_{\nu_1}^{\ulambda}$ and $g\in \KLRchar_{\nu_2}^{\umu}$, then, for $\bi\in \Seq_{\nu_1+\nu_2}^{w(\ulambda,\umu)}$,
\begin{equation*}
\textstyle(f\shuffle_w g)(\bi)=\sum_{\mathbf{j}_1,\,\mathbf{j}_2} f(\mathbf{j}_1)g(\mathbf{j}_2)
\end{equation*}
where the sum now runs over all the ways to write $\bi$ as a shuffle of $\bj_1\in \Seq_{\nu_1}^{\ulambda}$ and $\bj_2\in \Seq_{\nu_2}^{\umu}$. We say that the product $\shuffle_w$ is \textit{$w$-twisted} to keep in mind our non-canonical choice of $w$. In particular, for the identity element $w=e\in \Sigma_{\ell}$, we get a map of the form 
$$ \shuffle_e : \KLRchar^{(\lambda^{(1)})}_{\nu_1}\otimes\dots \otimes \KLRchar^{(\lambda^{(\ell)})}_{\nu_{\ell}}\to \KLRchar_{\nu_1+\dots+\nu_{\ell}}^{(\la^{(1)},\dots,\la^{(\ell)})}.$$
The next result follows by the same argument used in the proof of \cite[Proposition~2.18]{khovanov2009diagrammatic} (or via \cite[Proposition~3E.3(f)]{mathas2024cellularity}):
\begin{Proposition}\label{prop:standardization_and_shuffles}
Let $\ulambda=(\lambda^{(1)},\dots,\lambda^{(\ell)})$ be a list in $P_+$ and fix a sum $\nu=\nu_1+\dots+\nu_\ell$ of elements of $Q_+$. Use Remark \ref{rem:KLRWtoCyclo} and the above to write the diagram
\begin{equation*}
\adjustbox{scale=0.87}{
\begin{tikzcd}[scale=0.5, column sep = 4em, row sep = 1.75em]
K_0(T^{(\lambda^{(1)})}_{\nu_{1}}\fmod)\otimes\dots \otimes K_0(T^{(\lambda^{(\ell)})}_{\nu_{\ell}}\fmod) \arrow[r,"\ch\otimes\dots\otimes\ch"]\arrow[d,"{[}\bbS^\ulambda{]}"']& \KLRchar^{(\lambda^{(1)})}_{\nu_{1}}\otimes\dots \otimes \KLRchar^{(\lambda^{(\ell)})}_{\nu_{\ell}}\arrow[d,"\shuffle_e"]\\
K_0(T^{\ulambda}_{\nu}\fmod)\arrow[r,"\ch"] & \KLRchar^{\ulambda}_{\nu}
\end{tikzcd}}
\end{equation*}
where $e\in \Sigma_{\ell}$ is the identity element. Then the above diagram is commutative.
\end{Proposition}
\begin{Example}\label{ex:shuffle_product_KLRW_sl3}
Fix $\fg=\fsl_3$ and $\ulambda=(\varpi_1,\varpi_2$). Then the category $T^{(\varpi_1,\varpi_2)}_0\fmod$ categorifies the zero weight space of the tensor product $V(\varpi_1)\otimes V(\varpi_2)$ and thus has three simple objects. By Proposition \ref{prop:standardization_and_shuffles} and \cite[Section 8.1]{kleshchev2011representations}, we get
{\footnotesize
\begin{align*}
\ch(\bbS^{(\varpi_1,\varpi_2)}(L_{\varpi_1}, L_{-\varpi_1}))&=[({\color{red}\bf1})]\shuffle_e [(1,2,{\color{red}\bf 2})]=[({\color{red}\bf 1},1,2,{\color{red}\bf 2})]+[(1,{\color{red}\bf 1},2,{\color{red}\bf 2})]+[(1,2,{\color{red}\bf 1},{\color{red}\bf 2})]\\
\ch(\bbS^{(\varpi_1,\varpi_2)}(L_{\varpi_2-\varpi_1}, L_{\varpi_1-\varpi_2}))&=[(1,{\color{red}\bf 1})]\shuffle_e [(2,{\color{red}\bf 2})]=[(1,{\color{red}\bf 1},2,{\color{red}\bf 2})]+[(1,2,{\color{red}\bf 1},{\color{red}\bf 2})]+[(2,1,{\color{red}\bf 1},{\color{red}\bf 2})]\\
\ch(\bbS^{(\varpi_1,\varpi_2)}(L_{-\varpi_2}, L_{\varpi_2}))&=[(2,1,{\color{red}\bf 1})]\shuffle_e [({\color{red}\bf 2})]=[(2,1,{\color{red}\bf 1},{\color{red} \bf 2})]
\end{align*}}\noindent 
where the simple KLR modules appearing in the argument of the standardization functor $\bbS^{(\varpi_1,\varpi_2)}$ are the unique $1$-dimensional modules over the corresponding cyclotomic quotients. 
\end{Example}
\subsection{Parity KLRW algebras}\label{subsec:parity_KLRW_alg}
We now recall yet another variant of KLR algebras, called \textit{parity KLRW algebras}, which are idempotent sandwiches of KLRW algebras introduced in \cite{kamnitzer2019category} to help prove Theorem~\ref{th:descend}.\medskip\par
Fix a total order ``$<$'' on $I$%
. This induces a total order on $I\times_2\Z$ by reverse lexicographic order, i.e. we first order according to the integer and then break ties using the order on $I$.\medskip\par

Let $\lambda=\sum_{i\in I}\lambda_i\varpi_i\in P_+$ and fix an integral set of parameters $\bR=(R_i)_{i\in I}\in \Z^{\la}$ of level $\la$. Using the above total order, the data of $\bR$ produces an ordered list 
\begin{equation}\label{eq:ordered_list_bR}
(p_1,r_1)\leq\dots\leq (p_\ell,r_\ell)
\end{equation}
in $I\times_2\Z$, where $R_i=\{r_j\;|\; p_j=i\}$ and $\ell=\sum_{i\in I}\la_i$. Set
\begin{equation}\label{eq:varpibR}
\varpi_\bR=(\varpi_{p_1},\dots,\varpi_{p_\ell}).
\end{equation}
Fix $\nu\in Q_+$ with $\rht(\nu)=m$ and consider the KLRW algebra $\tilde{T}^{\bR}_\nu:=\tilde{T}^{\varpi_{\bR}}_\nu$ of Section \ref{sec:def_KLRW}. As the list $\varpi_{\bR}$ consists only of fundamental weights, one can think of sequences in $\Seq^{\varpi_{\bR}}_{\nu}$ as having letters in the alphabet $I\sqcup {\color{red}\bf I}$ given %
by two copies of the Dynkin diagram of $\fg$. One can thus talk about the \textit{parity of a strand}, it being the parity of the associated~Dynkin~node. We will also call \textit{longitude} the element of $\Z$ associated to a red strand using the data of $\bR$ (i.e.~the longitude of the $k^{\text{th}}$-red strand, associated to the fundamental weight $\varpi_{p_k}$, is $r_k$).\medskip\par
Let $(\bi,\kappa)\in \Seq_{\nu}^{\varpi_{\bR}}$ be a sequence. Write $(\bi,\kappa)=(j_1,\dots,j_{m+\ell})$ and fix $1\leq a\leq m+\ell-1$. Define 
\begin{equation*}
\delta(j_a,j_{a+1})=\begin{cases}
2 & \text{if }j_a,\;j_{a+1}\text{ have the same parity and }j_a\in I,\;j_{a+1}\in{\color{red}\bf I},\\
1 & \text{if }j_a,\;j_{a+1}\text{ have different parity},\\
0 &\text{otherwise}.
\end{cases}
\end{equation*}
We extend this definition by declaring
\begin{center}
$\delta(j_a,j_b)=\delta(j_b,j_a)$, $\delta(j_a,j_a)=0$, and $\delta(j_a,j_c)=\delta(j_a,j_b)+\delta(j_b,j_c)$ 
as soon as $a<b<c$.
\end{center}
We think of $\delta$ as measuring the distance between elements of the sequence $(\bi,\kappa)$.
\begin{Def}[{\cite[Definition~3.7]{kamnitzer2019category}}]\label{def:paritySeq}
A sequence $(\bi,\kappa)\in \Seq^{\varpi_{\bR}}_\nu$ is \textit{parity} if, for all pairs of red strands $j,j'$ of $(\bi,\kappa)$ having longitudes $r,r'$ respectively, %
$$\delta(j,j')\leq |r-r'|.$$
\end{Def}
The subset of parity sequences in $\Seq^{\varpi_{\bR}}_\nu$ will be denoted by $\Seq^{\bR}_\nu$. We call $e(\bi,\kappa)\in \tilde{T}^{\bR}_\nu$ a \textit{parity idempotent} if the sequence $(\bi,\kappa)$ is parity.

\begin{Example}\label{ex:many_sequences}
Take $\fg=\fsl_3$. %
Then, 
\begin{align*}
\Seq^{(\varpi_1,\varpi_2)}_{\alpha_1+\alpha_2}=\{
&(1,2,{\color{red}\bf 1},{\color{red}\bf 2}),
(2,1,{\color{red}\bf 1},{\color{red}\bf 2}),
(1,{\color{red}\bf 1},2,{\color{red}\bf 2}),
(2,{\color{red}\bf 1},1,{\color{red}\bf 2}),
({\color{red}\bf 1},1,2,{\color{red}\bf 2}),\\
&({\color{red}\bf 1},2,1,{\color{red}\bf 2}),
({\color{red}\bf 1},1,{\color{red}\bf 2},2),
({\color{red}\bf 1},2,{\color{red}\bf 2},1),
({\color{red}\bf 1},{\color{red}\bf 2},1,2),
({\color{red}\bf 1},{\color{red}\bf 2},2,1)\}.
\end{align*}
Choose $k\in \Z$ and consider the set of parameters $\bR=(\{a\},\{a+2k+1\})\in \C^{\varpi_1+\varpi_2}$ (as in Example \ref{ex:sl3_prod_mon_crystal}). Then the sequence $({\color{red}\bf 1},1,2,{\color{red}\bf 2})$ has distance $0+1+2=3$ between its~two~red strands and is hence parity if and only if $k\geq 1$. The same holds for the sequence $(1,{\color{red}\bf 1},2,{\color{red}\bf 2})$.
\end{Example}
Similarly to \eqref{eq:ordered_list_bR}, the data of a GT-weight $\bS=(S_i)_{i\in I}\in \Lambda_\nu$ (as in Section \ref{subsec:GT_weights}) produces an ordered list
\begin{equation}\label{eq:ordered_list_bS}
(i_1,s_1)\leq \dots \leq (i_m,s_m)
\end{equation}
where $S_i=\{s_k\;|\; i_k=i\}$. Notice that $(i_1,\dots, i_m)\in \Seq_\nu$.
\medskip\par
We order the set ${\color{red}\bf I}$ using the same order as on $I$ and then order $I\sqcup  {\color{red}\bf I}$ by considering \textcolor{red}{\bf red} vertices to be smaller than \textbf{black} ones. This makes $(I\sqcup {\color{red}\bf I})\times_2 \Z$ totally ordered via reverse lexicographic order and allows us to associate to each pair $(\bS,\bR)\in \Lambda_{\nu}\times \Z^{\la}$, a sequence
$$\bj_\bS=(j_1,\dots,j_{m+\ell})\in \Seq^{\varpi_{\bR}}_\nu $$ 
by combining the sequences associated to $\bR$ and $\bS$ via \eqref{eq:ordered_list_bR}--\eqref{eq:ordered_list_bS}. It is not hard to see that $\bj_\bS$ is always a parity sequence. Let $e(\bS):=e(\bj_\bS)\in \tilde{T}_\nu^\bR$ be the associated idempotent. The following lemma is an immediate consequence of the defining relations of KLRW algebras.
\begin{Lemma}[{\cite[Lemma~3.25]{kamnitzer2019category}}]
An idempotent $e(\bi,\kappa)\in \tilde{T}_\nu^\bR$ is parity if and only if there exists a GT-weight $\bS\in \Lambda_\nu$ such that $\tilde{T}_\nu^\bR e(\bi,\kappa)\simeq \tilde{T}_\nu^\bR e(\bS)$ as left $\tilde{T}_\nu^\bR$-modules.
\end{Lemma}
Let $e_{\bR}\in \tilde{T}^{\bR}_\nu$ be the sum of all parity idempotents of $\tilde{T}^{\bR}_\nu$. By abuse of notation, write $e_{\bR}\in T^{\bR}_\nu$ for the class of this sum in the cyclotomic quotient $T^{\bR}_\nu=\smash{\tilde{T}{}^{\ulambda}_\nu/J^{\ulambda}_\nu}$. 
\begin{Example}\label{ex:parity_sequences_sl3}
Following Example \ref{ex:many_sequences}, the sequences in $\Seq^{(\varpi_1,\varpi_2)}_{\alpha_1+\alpha_2}$ which give non-zero idempotents in the quotient $T^{\bR}_\nu$ can be seen to be (using the defining relations of $\smash{\tilde{T}{}^{\ulambda}_\nu}$) 
\begin{equation*}
(1,2,{\color{red}\bf 1},{\color{red}\bf 2}),
(2,1,{\color{red} \bf 1},{\color{red}\bf 2}),
(1,{\color{red}\bf 1},2,{\color{red}\bf 2}),
({\color{red}\bf 1},1,2,{\color{red}\bf 2}).
\end{equation*}
The first two are parity for $\bR=(\{a\},\{a+2k+1\})$, while the other are if and only if $k\geq 1$. 
\end{Example}
\begin{Def}[{\cite[Definition 3.7]{kamnitzer2019category}}] The \textit{parity KLRW algebra} is the idempotent sandwich algebra $P^\bR_\nu=e_{\bR}T^{\bR}_\nu e_{\bR}$.
\end{Def}
By construction, there is an idempotent truncation functor $$e_{\bR}(\trou):T^\bR_\nu\fmod\to P^\bR_\nu\fmod.$$
Also, as the %
map \eqref{eq:non_unital_map_KLRW} sends parity idempotents to parity idempotents, we have%
:

\begin{Theorem}[{\cite[Lemma~3.11 and Theorem~3.16]{kamnitzer2019category}}]\label{thm:parity_KLRW_categorifies_VR} The category 
$$\textstyle P^\bR\fmod=\bigoplus_{\nu\in Q_+} P^\bR_\nu\fmod$$ 
inherits, via the functor $e_{\bR}(\trou)$, a categorical (left-)action of $\fg$ from the category $T^{\bR}\fmod$. Moreover, the submodule
$$ K_\C(P^\bR\fmod)\subseteq K_{\C}(T^{\bR}\fmod)\simeq V(\varpi_{\bR})$$
has associated normal $\fg$-crystal given by the product monomial crystal $\cB(\lambda,\bR)$. 
\end{Theorem}
The sequences $(\bi,\kappa)\in \Seq^{\varpi_{\bR}}_\nu$ with $\kappa(1)=m=\rht(\nu)$ are clearly always parity. Taking the sum $e_{cyc}$ of these special parity idempotents instead of $e_{\bR}$ gives algebra isomorphisms
$$ R^{\la}_{\nu}\simeq e_{cyc}T_{\nu}^{\bR}e_{cyc}  \simeq e_{cyc}P_{\nu}^{\bR}e_{cyc},$$
where the first isomorphism is shown in \cite[Theorem~4.18~and~Proposition~5.31]{webster2017knot}. The next result follows from the relations of $T^{\bR}_{\nu}$ along with Theorem \ref{thm:charaterization_max_sing_crystals} and Definition \ref{def:paritySeq}.
\begin{Lemma}\label{coro:maximally_singular_categorifies_irrep} If $\bR$ is maximally singular, then $e_{cyc}=e_{\bR}$ in $T_{\nu}^{\bR}$. In particular, $R_{\nu}^{\la}\simeq P_{\nu}^{\bR}$. 
\end{Lemma}
\subsection{Characters for parity KLRW algebra modules}\label{subsec:char_parity_KLRW}
To define a notion of characters of parity KLRW algebra modules, %
let $\KLRchar^{\bR}_\nu=\Z^{\Seq^{\bR}_\nu}$ and define, for $M$ an object of $P^\bR_\nu\fmod$,
\begin{equation}\label{eq:charKLRWparity}
\textstyle\ch(M)=\sum_{(\bi,\kappa)\in \Seq^{\bR}_\nu} \dim (e(\bi,\kappa)M) [(\bi,\kappa)]\in \KLRchar^{\bR}_\nu.
\end{equation}
As $\smash{P^\bR_\nu\fmod}$ admits a categorical $\fg$-action, the reasoning used %
for Proposition \ref{prop:ch_is_hom_KLRW} shows:
\begin{Proposition}
The map $\ch:K_0(\smash{P^\bR_\nu\fmod})\to \smash{\KLRchar_\nu^{\bR}}$ is an injective $\Z$-module morphism.
\end{Proposition}
By construction, the following proposition holds:
\begin{Proposition}\label{prop:ch_parity_vs_klrw}
The diagram
\begin{equation*}
\adjustbox{scale=0.9}{
\begin{tikzcd}[scale=0.75]
K_0(T^\bR_\nu\fmod) \arrow[r,"\ch"]\arrow[d,"{[e_\bR(\trou)]}"']& \KLRchar^{\varpi_\bR}_\nu\arrow[d,two heads] \\
K_0(P^\bR_\nu\fmod) \arrow[r,"\ch"]& \KLRchar^{\bR}_\nu
\end{tikzcd}}
\end{equation*}
commutes, where the left-vertical map is the natural projection given by restriction to $\Seq^{\bR}_\nu$.
\end{Proposition}
For sets of parameters $\bR_1$ and $\bR_2$ of respective sizes $\ell_1$ and $\ell_2$, the union $\bR_1\cup \bR_2$ defines a unique element $w\in \Sigma_{\ell_1+\ell_2}$ of minimal length for which $w(\varpi_{\bR_1},\varpi_{\bR_2})=\varpi_{\bR_1\cup \bR_2}$. Hence, the $w$-twisted shuffle product of \eqref{eq:wShuffle} gives a map
$$ \shuffle_w:\KLRchar_{\nu_1}^{\varpi_{\bR_1}}\otimes \KLRchar_{\nu_2}^{\varpi_{\bR_2}} \rightarrow \KLRchar_{\nu_1+\nu_2}^{\varpi_{\bR_1\cup\bR_2}}. $$
\begin{Proposition}\label{prop:shuffle_sets_of_parameters}
There is a unique map $\shuffle_{\bR_1,\bR_2}:\KLRchar^{{\bR_1}}_{\nu_1}\otimes \KLRchar^{{\bR_2}}_{\nu_2}\to \KLRchar_{\nu_1+\nu_2}^{\bR_1\cup\bR_2}$ that makes 
\begin{equation*}
\adjustbox{scale=0.9}{
\begin{tikzcd}
\KLRchar^{\varpi_{\bR_1}}_{\nu_1}\otimes \KLRchar^{\varpi_{\bR_2}}_{\nu_2} \arrow[d,two heads]\arrow[r,"\shuffle_w"] & \KLRchar^{\varpi_{\bR_1\cup\bR_2}}_{\nu_1+\nu_2}\arrow[d,two heads]\\
\KLRchar^{{\bR_1}}_{\nu_1}\otimes \KLRchar^{{\bR_2}}_{\nu_2}\arrow[r,"\shuffle_{\bR_1,\bR_2}"] & \KLRchar^{{\bR_1\cup\bR_2}}_{\nu_1+\nu_2}
\end{tikzcd}}
\end{equation*}
commute (where both vertical arrows are the natural projections and $w$ is defined as above).
\end{Proposition}
\begin{proof}
Indeed, $w$-twisted shuffle products involving a sequence that is not parity cannot be parity (as the relative order of sequences is preserved by shuffling). This ends the proof.
\end{proof}\newpage
\begin{Example}
Continuing with Example \ref{ex:parity_sequences_sl3}, fix $\bR_1=(\emptyset,\{a+1\})$ and $\bR_2=(\{a\},\emptyset)$. Then $w=s_1$ and we see that (compare with Example \ref{ex:shuffle_product_KLRW_sl3})
\begin{align*}
[({\color{red}\bf 2})]\shuffle_{\bR_1,\bR_2}[(2,1,{\color{red}\bf 1})]&=[(2,1,{\color{red}\bf 1},{\color{red}\bf 2})], \\
[(2,{\color{red}\bf 2})\shuffle_{\bR_1,\bR_2}(1,{\color{red}\bf 1})]&=[(2,1,{\color{red}\bf 1},{\color{red}\bf 2})]+[(1,2,{\color{red}\bf 1},{\color{red}\bf 2})], \\
[(1,2,{\color{red}\bf 2})]\shuffle_{\bR_1,\bR_2}[({\color{red}\bf 1})]&=[(1,2,{\color{red}\bf 1},{\color{red}\bf 2})].
\end{align*}
\end{Example}

\subsection{Relation with truncated shifted Yangians}\label{subsec:KLRW_and_TSY} We now study an equivalence of categories which relates modules over parity KLRW algebras and truncated shifted Yangians. The main goal here is to explain the precise relationship between the notion of $\ell$-characters given in Section \ref{subsec:ell_characters_and_blocks} and our notion of characters for parity KLRW algebra modules.
\medskip\par
Given a multiset $S\in \C^{m}/\Sigma_m$, let $\sigma(S)$ be the size of the stabilizer in $\Sigma_m$ of a preimage of $S$ in $\C^m$. Namely, $ \textstyle \sigma(S)=\prod_{i=1}^k n_i!$ if $S$ consists of the elements $s_1,\dots,s_k$ with respective multiplicity $n_1,\dots, n_k$. Also let $\sigma(\bS)=\prod_{i\in I} \sigma(S_i)$ for $\bS=(S_i)_{i\in I}\in \Lambda_{\nu}$.
\begin{Theorem}[{\cite[Theorem~5.2]{kamnitzer2019category}}]\label{thm:equiv_of_cats_parity_and_Osh}
Let $\lambda\in P_+$, $\mu\in P$ and choose $\bR\in \Z^\lambda$ a set of parameters. Then, there is an equivalence of categories 
\begin{equation*}
\smash{\Theta_\bR:P^\bR_{\lambda-\mu}\fmod\to\O_\mu^\lambda(\bR)}
\end{equation*}
satisfying, for each finite-dimensional $P^\bR_{\lambda-\mu}$-module $M$,
\begin{equation}\label{eq:equality_of_S_weight_spaces}
\smash{\dim e(\bS)M=\sigma(\bS)\dim W_\bS(\Theta_\bR(M))}.
\end{equation}
\end{Theorem}
\begin{Rem} The above gives an equivalence $P^{\bR}\fmod\to\mathcal{O}_{sh}^{\la}(\bR)$ that we also note~$\Theta_{\bR}$.
\end{Rem}
\begin{Rem}
The factor $1/\sigma(\bS)$ appearing above accounts for the fact that the idempotents $e(\bS)$ are not primitive (see the discussion following \cite[(2.47)]{khovanov2009diagrammatic}). Remark that in \cite[(2.1.16)]{varagnolo2025representations}, the idempotents defined come from choosing indecomposable projectives in a product of nil-Hecke algebras (parametrized by $\bS\in \Lambda_{\nu}$). This removes the need for~the extra $1/\sigma(\bS)$ factor in their analogue of Theorem \ref{thm:equiv_of_cats_parity_and_Osh} (see also \cite[Section~4.2]{varagnolo2011canonical}).
\end{Rem}
\begin{Rem}\label{rem:OHZ_is_too_small}
The above result justifies the choice of category $\O$ alluded to in Remark~\ref{rem:2categoriesO}.  Indeed, to have an equivalence as in Theorem \ref{thm:equiv_of_cats_parity_and_Osh}, one must have, for $\fg=\mathfrak{sl}_2$,
$$\smash{\mathcal{O}_0^{\alpha_1}(\{0,0\})\simeq P_{\alpha_1}^{\{0,0\}}\fmod\simeq R^{\alpha_1}_{\alpha_1}\fmod \simeq \nicefrac{\C[x]}{(x^2)}\fmod}$$
as $P_{\alpha_1}^{\{0,0\}}\simeq R^{\alpha_1}_{\alpha_1}\simeq \nicefrac{\C[x]}{(x^2)}$ (by Lemma \ref{coro:maximally_singular_categorifies_irrep} and the definition of cyclotomic KLR algebras). However, considering only modules that decompose in non-generalized weight spaces for~$\mathcal{O}_0$ (i.e.~replacing this category by the full subcategory $\mathcal{O}_0^{\text{HZ}}$ defined in Remark \ref{rem:2categoriesO}) is easily seen to lead to a ``truncated category $\mathcal{O}$'' equivalent to the category $\C\fmod$ of finite-dimensional vector spaces for $\la=\alpha_1$ and $\bR=\{0,0\}$. Also, even if we stop caring about~parity KLRW algebras, 
the uniqueness property of categorifications in \cite{chuang2008derived} implies that truncations~of the form $\mathcal{O}_{\mu}^{\text{HZ}}$ cannot in general be endowed with a categorical $\fg$-action.
\end{Rem}

Inspired by \eqref{eq:equality_of_S_weight_spaces}, we define a $\Z$-linear map $\theta_\bR:\KLRchar^\bR_\nu \to \EGT$ via
\begin{align*}
[(\bi,\kappa)] \mapsto \sum_{\bS\in\Lambda_\nu,\; \bj_\bS=(\bi,\kappa)} \tfrac{1}{\sigma(\bS)}[\bS]
\end{align*}
where the sum runs over all GT-weights for which $e(\bS)=e(\bi,\kappa)$. As an immediate corollary of the definition of this map, we get the result below.
\begin{Lemma}\label{lem:commutative_diag_Theta_and_theta}
The following diagram commutes:
\begin{equation*}
\adjustbox{scale=0.9}{
\begin{tikzcd}
K_0(P^\bR_\nu\fmod) \arrow[r,"\ch"] \arrow[d,"{{[}\Theta_\bR}{]}"']& \KLRchar^{\bR}_\nu\arrow[d,"{\theta_\bR}"]\\
K_0(\O^\lambda_\mu(\bR))\arrow[r,"{\chiGT^\bR}"] & \EGT
\end{tikzcd}}
\end{equation*}
\end{Lemma}
\begin{Example}\label{ex:gt_character_neg_prefund_parity_example}
Take $\fg=\mathfrak{sl}_2$. Then the parity KLRW algebra $\smash{P^{\{a\}}_{-\varpi_1}}\simeq R^{\varpi_1}_{\alpha_1}$ has a unique 1-dimensional irreducible module $M$ (up to isomorphism) whose character is $\ch(M)=[(1,{\color{red}\bf 1})]$. Also, a GT-weight $\bS=(\{b\})\in \Lambda_{\alpha_1}$ satisfies $\bj_\bS=(1,{\color{red}\bf 1})$ if and only if $b=a-2s$ with $s\in \Z_{\geq 1}$. 
Thus, we deduce that, as expected by Example \ref{ex:sl2_GTchar_negative_prefund},
$$\theta_\bR([(1,{\color{red}\bf 1})])=[\{a-2\}]+[\{a-4\}]+[\{a-6\}]+\dots=\chiGT^{\{a\}}(L(\sfPsi_{a-2}^{-1})).$$
\end{Example}

While there is \textit{a priori} no monoidal structure\footnote{As mentioned in Section \ref{sec:Intro}, it would be interesting to define such a monoidal structure that makes the equivalences of Theorem \ref{thm:equiv_of_cats_parity_and_Osh} monoidal (with respect to the truncated shifted coproducts of Section \ref{sec:TSC}). } on the category of parity KLRW algebra modules, we can use the twisted shuffle products of Proposition \ref{prop:shuffle_sets_of_parameters} to multiply characters. We show that this multiplication agrees --- via the maps $\theta_{\bR}$ --- with the multiplication given in Section \ref{subsec:GT_weights} for %
$\EGT$ (and hence with tensor products in $\mathcal{O}_{sh}$ by Proposition \ref{prop:GTchar_is_multiplicative}).
\begin{Proposition}\label{prop:theta_vs_Theta_characters}
Fix $\lambda_1,\lambda_2\in P_+$ and choose sets of parameters $\bR_1\in \C^{\lambda_1}$ and $\bR_2\in \C^{\lambda_2}$. Fix also $\nu_1,\nu_2\in Q_+$ and take $M_1$ and $M_2$ in $P^{\bR_1}_{\nu_1}\fmod$ and $P^{\bR_2}_{\nu_2}\fmod$ (respectively). Then, %
\begin{equation*}
\theta_{\bR_1}(\ch(M_1))\ast \theta_{\bR_2}(\ch(M_2))=\theta_{\bR_1\cup\bR_2}\big(\ch(M_1)\shuffle_{\bR_1,\bR_2} \ch(M_2)\big)
\end{equation*}
\end{Proposition}
The $\fsl_2$-case of the above result boils down to the fact that, given a multiset $S\in \C^m/\Sigma_m$ and a partition $m=m_1+m_2$, 
\begin{equation*}
\sum_{S_1,S_2}\tfrac{\sigma(\bS)}{\sigma(S_1)\sigma(S_2)}={\binom{m}{m_1}}
\end{equation*}
where the sum runs over all pairs of multisets $(S_1,S_2)\in \C^{m_1}/\Sigma_{m_1}\times \C^{m_2}/\Sigma_{m_2}$ that satisfy $S_1\cup S_2=S$. %
Our proof extends this idea.

\begin{proof}[Proof of Proposition \ref{prop:theta_vs_Theta_characters}]
Fix $\bj_1 \in \Seq^{\bR_1}_{\nu_1}$ and $\bj_2 \in \Seq^{\bR_2}_{\nu_2}$. We will prove that 
\begin{equation*}
\theta_{\bR_1}([\bj_1])\ast \theta_{\bR_2}([\bj_2])=\theta_{\bR}\big(\,[\bj_1]\shuffle_{\bR_1,\bR_2} [\bj_2]\, \big)
\end{equation*}
which, by linearity, implies the proposition. Unpacking definitions, we see that it is enough to prove that 
$$
\sum_{\bS_1, \bS_2} \tfrac{1}{\sigma(\bS_1) \sigma(\bS_2)} [\bS_1 \cup \bS_2] = \sum_{\bj \in \Seq^{\bR}_{\nu}}  n(\bj_1, \bj_2; \bj)   \sum_{\bS\in\Lambda_\nu,\; \bj_\bS= \bj}  \tfrac{1}{\sigma(\bS)} [\bS] = \sum_{\bS\in \Lambda_{\nu}}\tfrac{1}{\sigma(\bS)}n(\bj_1,\bj_2;\bj_{\bS})[\bS]
$$ 
where the first sum runs over pairs $(\bS_1, \bS_2) $ such that $(\bj_{\bS_1},\bj_{\bS_2}) = (\bj_1,\bj_2)$, and with $ n(\bj_1, \bj_2; \bj) $ the number of shuffles of $ \bj_1 $ and $ \bj_2 $ which produce $ \bj $.\medskip\par
Choose $ \bS \in \Lambda_\nu $ and let $ \bj= \bj_\bS$.  Extracting the coefficient of $ [\bS] $ on the left and right sides above, we see that it suffices to show that
\begin{equation}
\label{eq:toShow}
\sum_{\bS_1, \bS_2} \tfrac{\sigma(\bS)}{\sigma(\bS_1) \sigma(\bS_2)} = n(\bj_1, \bj_2; \bj)
\end{equation}
where the sum on the left-hand side runs over all pairs $ (\bS_1, \bS_2)$ such that $(\bj_{\bS_1},\bj_{\bS_2}) = (\bj_1,\bj_2) $ and $ \bS_1 \cup \bS_2 = \bS $. Consider the ordered list $ (i_1, s_1) \le \cdots \le (i_m, s_m) $ associated to $\bS$ via \eqref{eq:ordered_list_bS} and let $m_1$, $m_2 $ be the number of black strands in $\bj_1$, $\bj_2 $ (resp.). Note that $m=m_1+m_2$.\medskip\par
Let $A_1 \sqcup A_2%
$ be an ordered  set partition of $ \{1,\dots, m\} $ into sets of respective sizes $m_1$ and $m_2$.  Consider the GT-weight $ \bS(A_1)=(\bS(A_1)_i)_{i\in I} \in \Lambda_{\nu_1} $ given by 
\begin{equation*}
\bS(A_1)_i = \{ s_a \,|\, a \in A_1 \text{ and } i_a = i \} ,
\end{equation*}
and define $ \bS(A_2)\in \Lambda_{\nu_2}$ similarly.  Then $ \bS(A_1) \cup \bS(A_2) = \bS $.  \medskip\par 
Finally, let $ X(\bj_1, \bj_2; \bS)  $ be the set of ordered set partitions $(A_1, A_2) $ of $\{1,\dots,m\}$ such that $ \bj_{\bS(A_1)} = \bj_1 $ and $ \bj_{\bS(A_2)} = \bj_2$. By construction, the set $ X(\bj_1, \bj_2; \bS) $ is in bijection with the shuffles of $ \bj_1 $ and $ \bj_2 $ which produce $\bj$, and thus has cardinality $ |X(\bj_1, \bj_2; \bS)| = n(\bj_1, \bj_2; \bj) $.

Let 
\begin{equation*}
\textstyle\Sigma(\bS) = \{ w \in \Sigma_m \mid i_{w(a)} = i_a \text{ and } s_{w(a)} = s_a \}\leq \prod_{i\in I} \Sigma_{m_i}
\end{equation*}
be the stabilizer of an element of the preimage of $\bS$ in $\prod_{i\in I}\C^{m_i}$. The group $ \Sigma(\bS) $ has size $ \sigma(\bS) $ %
and acts on $ X(\bj_1, \bj_2; \bS) $ through the usual action of $ \Sigma_m $ on set partitions. Furthermore, two set partitions $(A_1, A_2), (B_1, B_2) \in X(\bj_1, \bj_2; \bS) $ lie in a common $ \Sigma(\bS)$-orbit if and only if $ \bS(A_1) = \bS(B_1) $ and $ \bS(A_2) = \bS(B_2)$. Thus, the set of orbits $ X(\bj_1, \bj_2; \bS) / \Sigma(\bS) $ is in bijection with pairs $ (\bS_1, \bS_2)$ satisfying $(\bj_{\bS_1},\bj_{\bS_2}) = (\bj_1,\bj_2) $ and $ \bS_1 \cup \bS_2 = \bS $, i.e.~with the indexing set of the sum in \eqref{eq:toShow}. Hence, summing over orbits and noticing that the stabilizer in $\Sigma(\bS)$ of $(A_1,A_2)\in X(\bj_1,\bj_2;\bS)$ has size $ \sigma(\bS(A_1)) \sigma(\bS(A_2)) $, the orbit-stabilizer theorem gives
\begin{equation*}
\sum_{\bS_1, \bS_2}\tfrac{\sigma(\bS)}{\sigma(\bS_1)\sigma(\bS_2)}=\sum_{\bS_1, \bS_2}\tfrac{|\Sigma(\bS)|}{\sigma(\bS_1)\sigma(\bS_2)}=|X(\bj_1,\bj_2;\bS)|=n(\bj_1,\bj_2;\bj),
\end{equation*}
where the sums run over the same indices as in \eqref{eq:toShow}. This completes the proof.
\end{proof}

\begin{Rem}\label{rem:t_deformation_from_shuffles}
The work of \cite{webster2017knot} shows that the algebras $T^{\bR}_\nu$ have a natural $\Z$-grading that propagates to the algebras %
$P^{\bR}_\nu$ (as the idempotents $e(\bi,\kappa)$ have degree zero). Thus,~this grading can be transported to the categories $\O_{sh}^\lambda(\bR)$ via the equivalences of Theorem \ref{thm:equiv_of_cats_parity_and_Osh}. It is also well-known that the shuffle product defined above has a $\Z$-graded analogue. Using these two facts, it is possible to use Proposition \ref{prop:theta_vs_Theta_characters} to give a $t$-deformation of the product of the ring $\EGT$ of GT-characters, where $t$ is the parameter coming from the shift~functor. In future work, the authors and H.~Murata plan to investigate whether these $t$-deformations glue into a $t$-deformation of $K_0(\mathcal{O}_{sh})$ (see %
Section \ref{sec:glueing_of_g_action} for details on the gluing process).
\end{Rem}

\begin{Example}\label{ex:ParityKLRqWronsk}
Take $\fg=\fsl_2$ and let $L_{\varpi_1,a}$ and $L_{-\varpi_1,a}$ be the unique (up to isomorphism) $1$-dimensional simple objects of $P^{\{a\}}\fmod\simeq R^{\varpi_1}\fmod$ (for $a\in 2\Z$). Then 
$$\ch(L_{\varpi_1,a})=[(\redbold{1})]\text{ and }\ch(L_{-\varpi_1,a})=[(1,\redbold{1})]$$
with $\textstyle \theta_{\{a\}}\big(\,[(\redbold{1})]\,\big)=[\emptyset]$ and $\theta_{\{a\}}\big(\,[(1,\redbold{1})]\,\big)=\sum_{s\geq 1}\; [\{a-2s\}]$ (compare with Example \ref{ex:gt_character_neg_prefund_parity_example}).
Also, using Theorem \ref{th:descend}, we get
$$\Theta_{\{a\}}(L_{\varpi_1,a})=L(\mathsf{\Psi}_{1,a})\text{ and } \Theta_{\{a\}}(L_{-\varpi_1,a})=L(\mathsf{\Psi}_{1,a-2}^{-1}).$$
Now, choose $k\in \Z_{\geq 0}$ and let $\bR_1=\{a-2k\}$ with $\bR_2=\{a+2\}$. By definition, 
\begin{align*}
\ch(L_{\varpi_1,a-2k})\shuffle_{\bR_1,\bR_2}\ch(L_{-\varpi_1,a+2})&=[(\redbold{1},1,\redbold{1})]+[(1,\redbold{1},\redbold{1})],\\
&=[(\redbold{1},1,\redbold{1})]+\ch(L_{-\varpi_1,a-2k})\shuffle_{\bR_1,\bR_2}\ch(L_{\varpi_1,a+2})
\end{align*}
so that applying the above with Proposition \ref{prop:GTchar_is_multiplicative}, Lemma \ref{lem:commutative_diag_Theta_and_theta} and Proposition \ref{prop:theta_vs_Theta_characters} gives
\begin{align*}
\chiGT^{\bR}(L(\mathsf{\Psi}_{1,a-2k})\otimes L(\mathsf{\Psi}_{1,a}^{-1}))&=\theta_{\bR_1}(\ch(L_{\varpi_1,a-2k}))\ast\theta_{\bR_2}(\ch(L_{-\varpi_1,a+2}))
\\&= \theta_{\bR}(\,[(\redbold{1},1,\redbold{1})]\,\big)+\chiGT^{\bR}(L(\mathsf{\Psi}_{1,a-2(k+1)}^{-1})\otimes L(\mathsf{\Psi}_{1,a+2}))
\end{align*}
for $\bR=\bR_1\cup\bR_2$. On the other hand, using the map $[\Psi_{\bR}]%
$ of Section \ref{subsec:GT_weights}, we easily get
\begin{equation*}
\textstyle ([\Psi_{\bR}]\circ\theta_{\bR})\big(\,[(\redbold{1},1,\redbold{1})]\,\big)=\mathsf{Y}_{1,a-2(k-1)}\dots \mathsf{Y}_{1,a}(1+\sum_{r=0}^{k-1} \mathsf{A}_{1,a}^{-1}\dots \mathsf{A}_{1,a-2r}^{-1}) = \chi_{\ell}(\textsf{W}_{k,a}^{(1)})
\end{equation*}
(see Example \ref{ex:lchar}) and it thus follows from Proposition \ref{prop:injectivity_chiGT} that (setting $\chi_{\ell}(W_{0,a}^{(1)})=1$) 
$$ \chi_{\ell}(L(\mathsf{\Psi}_{1,a-2k}))\chi_{\ell}( L(\mathsf{\Psi}_{1,a}^{-1}))=\chi_{\ell}(\textsf{W}_{k,a}^{(1)})+\chi_{\ell}(L(\mathsf{\Psi}_{1,a-2(k+1)}^{-1}))\chi_{\ell}(L(\mathsf{\Psi}_{1,a+2}))$$
which corresponds to a well-known relation in $K_0(\mathcal{O}_{sh})$.
\end{Example}

Fix a set of parameters $\bR\in \Z^\lambda$ and order its elements  
$$(p_1,r_1)\leq \dots \leq (p_\ell,r_\ell)$$
as in \eqref{eq:ordered_list_bR}. Let $\nu_1,\dots,\nu_\ell\in Q_+$ and set $\nu=\nu_1+\dots+\nu_\ell$. %
Recall that the standardization functor $\bbS^{\varpi_{\bR}}$ and Remark \ref{rem:KLRWtoCyclo} give a morphism of Grothendieck groups
\begin{equation*}
[\bbS^{\varpi_\bR}]: K_0(T^{(\varpi_{p_1})}_{\nu_1}\fmod)\otimes \dots \otimes K_0(T^{(\varpi_{p_\ell})}_{\nu_\ell}\fmod)\to K_0(T^{\bR}_{\nu}).
\end{equation*}

\begin{Theorem}\label{thm:the_big_pentagone_commutes}
The following pentagon commutes:
\begin{equation*}
\adjustbox{scale = 0.85}{
\begin{tikzcd}[column sep=0em,row sep=0.5em]
& & {\bigotimes_{i=1}^\ell K_0(\O^{\varpi_{p_i}}_{\varpi_{p_i}-\nu_i}(r_i))}\arrow[drr,"\otimes"] & &\\
\bigotimes_{i=1}^\ell K_0(T^{(\varpi_{p_i})}_{\nu_i}\fmod)\arrow[rd,"{[\bbS^{\varpi_\bR}]}"']\arrow[rru,"{\otimes_{i=1}^\ell [\Theta_{r_i}]}"] & & & & K_0(\O^\lambda_{\lambda-\nu}(\bR)) \\
& {K_0(T^{\bR}_\nu\fmod)}\arrow[rr,"{[e_{\bR}(\trou)]}"'] & & {K_0(P^\bR_\nu\fmod)}\arrow[ru,"{[\Theta_\bR]}"'] &
\end{tikzcd}}
\end{equation*}
\end{Theorem}
\begin{proof}
Consider the diagram
\begin{equation}\label{eq:DiagPentEtendu}
\def\angle{60}
\adjustbox{scale = 0.85}{
\begin{tikzcd}[column sep = 4em,row sep=1.3em]
\bigotimes_{i=1}^\ell K_0(T^{(\varpi_{p_i})}_{\nu_i}\fmod) 
	\arrow[bend right=\angle,swap]{ddd}{\otimes_{i=1}^\ell [\Theta_{r_i}]}
	\arrow[r,"{[\bbS^{\varpi_{\bR}}]}"]
	\arrow[d,"\otimes_{i=1}^\ell \ch"]
& K_0(T^{\bR}_\nu\fmod) 
	\arrow[r,"{[e_\bR(\trou)]}"]
	\arrow[d,"\ch"]
& K_0(P^\bR_\nu\fmod) 
	\arrow[bend left=\angle]{ddd}{[\Theta_\bR]}
	\arrow[d,"\ch"']\\
\bigotimes_{i=1}^\ell\KLRchar_{\nu_i}^{(\varpi_{p_i})}
	\arrow[d,"{\otimes_{i=1}^\ell \theta_{r_i}}"] 
	\arrow[r,"\shuffle_e"]
& \KLRchar_{\nu}^{\varpi_\bR} 
	\arrow[r,two heads]
& \KLRchar_{\nu}^{\bR}
	\arrow[d,"\theta_\bR"']\\
\bigotimes_{i=1}^\ell \EGT\arrow[rr,"\ast"] 
& 
& \EGT\\
\bigotimes_{i=1}^\ell K_0(\O^{\varpi_{p_i}}_{\varpi_{p_i}-\nu_i}(r_i))
	\arrow[rr,"\otimes"]
	\arrow[u,"\otimes_{i=1}^\ell\chiGT^{r_i}"'] 
& 
& K_0(\O^\lambda_{\lambda-\nu}(\bR))
	\arrow[u,"\chiGT^\bR"]
\end{tikzcd}}
\end{equation}\noindent
The top left/right squares respectively commute by Proposition \ref{prop:standardization_and_shuffles} and Proposition \ref{prop:ch_parity_vs_klrw}. The middle/bottom rectangles also respectively commute by Proposition \ref{prop:theta_vs_Theta_characters} and Proposition \ref{prop:GTchar_is_multiplicative}. Finally, the leftmost/rightmost parts commute by Lemma \ref{lem:commutative_diag_Theta_and_theta}. Hence, given%
\begin{center}
$M_1$ in $T_{\nu_1}^{(\varpi_{p_1})}\fmod$, \dots, $M_{\ell}$ in $T_{\nu_{\ell}}^{(\varpi_{p_{\ell}})}\fmod$,
\end{center}
we get
$$ \chiGT^{\bR}\big(\,\Theta_{r_1}(M_1)\otimes\dots\otimes \Theta_{r_\ell}(M_\ell)\,\big)= \chiGT^{\bR}(\Theta_{\bR}(e_{\bR}\bbS^{\varpi_{\bR}}(M_1\otimes \dots\otimes M_{\ell})))$$
and the commutativity of the outside boundary of \eqref{eq:DiagPentEtendu} follows from Proposition \ref{prop:injectivity_chiGT}.
\end{proof}
We end this section with a final corollary related to dual canonical bases. Use the above notation with \eqref{eq:Tlambdafmod} and consider the isomorphism 
\begin{equation}\label{eq:IsoDualCanKLRW}
\smash{\textstyle V(\varpi_\bR) \cong \bigotimes_{j=1}^{\ell}K_{\C}(T^{(\varpi_{p_j})}%
\fmod)\simeq \bigotimes_{j = 1}^\ell K_\C(\O^{\varpi_{p_j}}_{sh}(r_j))}
\end{equation}
where the second map is $\otimes_{j=1}^{\ell}[\Theta_{r_j}]$ . We call \textit{dual canonical basis} of $\bigotimes_{j = 1}^\ell K_\C(\O^{\varpi_{p_j}}_{sh}(r_j))$~the image of  the dual canonical basis of the tensor product $V(\varpi_{\bR})$ under \eqref{eq:IsoDualCanKLRW}.
\begin{Corollary}\label{cor:multiplication_and_DC_basis}
Multiplication $\bigotimes_{j=1}^{\ell}K_{\C}(\mathcal{O}_{sh}^{\varpi_{p_j}}(r_j))\to K_{\C}(\mathcal{O}_{sh}^{\la}(\bR))$ takes elements of the dual canonical basis to simple classes or zero.
\end{Corollary}

\begin{proof} By Theorem \ref{thm:the_big_pentagone_commutes}, the composition of \eqref{eq:IsoDualCanKLRW} with the multiplication map is equal to
$$\textstyle V(\varpi_{\bR})\simeq \bigotimes_{j=1}^{\ell}K_{\C}(T^{(\varpi_{p_j})}%
\fmod)\xrightarrow{[\bbS^{\varpi_{\bR}}]}K_{\C}(T^{\varpi_{\bR}}%
\fmod)\xrightarrow{[\Theta_{\bR}]\circ[e_{\bR}(\trou)]}K_{\C}(\mathcal{O}_{sh}^{\la}(\bR)).$$
The result then follows from Theorem \ref{thm:dual_canonical_tensorprod} since equivalences preserve simplicity (and since idempotent truncation functors always send simple classes to simple classes or zero).
\end{proof}

\section{Gluing of the \texorpdfstring{$G$}{G}-action}\label{sec:glueing_of_g_action}
Using parity KLRW algebras, each category $\O_{sh}^\lambda(\bR)$ can be equipped with a categorical left $\fg$-action %
via \textit{transport de structure} along the equivalence given in Theorem \ref{thm:equiv_of_cats_parity_and_Osh}. Recall that, by Theorem \ref{thm:parity_KLRW_categorifies_VR}, the (left) $\fg$-module $K_{\C}(\O_{sh}^\lambda(\bR))$ is characterized by the fact that~its $\fg$-crystal is the product monomial crystal $\cB(\lambda,\bR)$. In particular, for $\bR$ maximally singular%
, $K_{\C}(\O_{sh}^\lambda(\bR))$ is isomorphic to the irreducible module $V(\lambda)$ (see also Lemma \ref{coro:maximally_singular_categorifies_irrep}) and the corresponding isomorphism 
$V(\la)\simeq K_{\C}(\mathcal{O}_{sh}^{\la}(\bR))$
maps the highest weight vector~$v_{\lambda}$~of~$V(\lambda)$ to the class of the $1$-dimensional module $\bigotimes_{i\in I}\bigotimes_{a\in \bR_i} L_{\varpi_i,a}$.

Fix (as in Section \ref{sec:Intro}) a connected and simply connected complex algebraic group $G$ with $\op{Lie}(G)\simeq\fg$. There is an equivalence of categories of complex representations $\mathrm{Rep}\,\fg\simeq \mathrm{Rep}\,G$
and the $\fg$-action on $K_{\C}(\mathcal{O}_{sh}^{\la}(\bR)$ can be integrated to a $G$-action. We use interchangeably both the above $G$ and $\fg$-actions in what follows and establish that %
the actions coming~from the %
various truncations glue to a $G$-action on the ring $K_\C(\O_{sh})$ (under the transition maps coming from inclusions of %
crystals or, as we explain below, Theorem \ref{thm:TFAEsurjection}).

Let $\lambda_1,\lambda_2\in P_+$ and write $\bR_1\leq \bR_2$ whenever sets of parameters $\bR_1\in \Z^{\lambda_1}$ and $\bR_2\in \Z^{\lambda_2}$ satisfy one of the equivalent conditions of Theorem \ref{th:TFAE}. Given $\bR_1\leq \bR_2$ and $\mu\in P$ such that $\cB(\la_1,\bR_1)_{\mu}\neq \emptyset$, we can use Theorem \ref{thm:TFAEsurjection} to obtain the commutative diagram
$$
\adjustbox{scale=0.8,center}{
\begin{tikzcd}[column sep=0.7em, row sep=1em]
& Y_{\mu} \ar[dl,"\Phi_{\mu}^{\la_2}(\bR_2)",swap,two heads] \ar[dr,"\Phi_{\mu}^{\la_1}(\bR_1)",two heads] &\\
Y^{\la_2}_{\mu}(\bR_2) \ar[rr, two heads] &  & Y_{\mu}^{\lambda_1}(\bR_1)
\end{tikzcd}}
$$
\noindent and pullback through the horizontal 
map gives a fully faithful functor $\O_{\mu}^{\lambda_1}(\bR_1)\hookrightarrow \O_{\mu}^{\lambda_2}(\bR_2)$,
which in turn induces an embedding 
of the corresponding Grothendieck groups.
Summing over $\mu$ (and using Theorem \ref{th:TFAE}), we get a $\Z$-module inclusion
\begin{equation}\label{eq:embedding_of_KOshlambda}
\smash{K_0(\O_{sh}^{\lambda_1}(\bR_1))\subseteq K_0(\O_{sh}^{\lambda_2}(\bR_2)).}
\end{equation}
We start by proving, using %
GT-characters, that the corresponding inclusion of $K_{\C}(\O_{sh}^{\lambda_1}(\bR_1))$ in $K_{\C}(\O_{sh}^{\lambda_2}(\bR_2))$ is $G$-equivariant, and then prove equivariance of multiplication in $K_{\C}(\mathcal{O}_{sh})$.

\subsection{Equivariance of inclusions}\label{subsec:equiv_of_inclusions}
Fix $\la\in P_+$ and $\bR\in \Z^{\la}$. Also fix $\nu\in Q_+$ and let~$i\in I$. Similarly to Lemma \ref{lem:chEiCycKLR}, the definition of the functor $\cE_i$ gives 
$e(i,\bj)M=e(\bj)\cE_i(M)$
for every $\bj\in \Seq^\bR_{\nu-\alpha_i}$ and each $P_{\nu}^{\bR}$-module $M$. The Yangian analogue of this is given by:
\begin{Lemma}\label{lem:GT_character_of_Ei_V}
Fix $V$ in $\O_{\mu}^\lambda(\bR)$ where $\mu=\la-\nu$. Then, %
for all %
$\bS\in \Lambda_{\nu-\alpha_i}$ and all $n\in \overline{i}+2\Z$ smaller than every element of $\bR$ and $\bS$,
\begin{equation}\label{eq:GTspaceEiO}
\dim W_{\bS}(\cE_i(V))=\dim W_{\bS\cup(n)_i}(V),
\end{equation}
where $(n)_i\in\Lambda_{\alpha_i}$ is the GT-weight containing only $n$.
\end{Lemma}
\begin{proof}
Take $M$ in $P_{\nu}^{\bR}\fmod$ with $\Theta_{\bR}(M)\simeq V$ and let $\bj\in \Seq_{\nu-\alpha_i}^{\bR}$ be such that $e(\bj)=e(\bS)$. Then $e(\bS\cup (n)_i)=e(i,\bj)$ and $\sigma(\bS)=\sigma(\bS\cup(n)_i)$ as $n$ is smaller than every element of $\bR$ and $\bS$. Thus, by \eqref{eq:equality_of_S_weight_spaces} and the discussion above, we have that 
\begin{align*}
\dim W_{\bS}(\cE_i (V))&=\dim W_{\bS}(\Theta_{\bR}(\cE_i (M)))=\tfrac{1}{\sigma(\bS)}\dim e(\bj)\cE_i(M)=\tfrac{1}{\sigma(\bS)}\dim e(i,\bj)M\\
&=\tfrac{1}{\sigma(\bS\cup(n)_i)}\dim e(\bS\cup(n)_i)M=\dim W_{\bS\cup(n)_i}(V)
\end{align*}
as desired.
\end{proof}
\begin{Rem} The class of $\cE_i(V)$ in $K_0(\O_{\mu+\alpha_i}^\lambda(\bR))$ is uniquely characterized by \eqref{eq:GTspaceEiO} (since the GT-character map %
is injective by Proposition \ref{prop:injectivity_chiGT}).
\end{Rem}
Recall that we want to prove the $G$-equivariance of the embeddings \eqref{eq:embedding_of_KOshlambda}.  We will need:%

\begin{Lemma}%
\label{lem:Bequiv_implies_Gequiv}
Fix %
finite-dimensional $G$-representations $V$ and $W$. Then every $B$-equivariant linear map from $V$ to $W$ is also $G$-equivariant. 
\end{Lemma}
\begin{proof}
This follows from %
$\Hom_{G}(V,W)\simeq (V^\ast\otimes W)^G\simeq (V^\ast\otimes W)^B\simeq \Hom_{B}(V,W)$, %
where the middle isomorphism holds as any (non-zero) $B$-invariant vector in a finite-dimensional $G$-module generates a copy of the trivial $G$-module, and is therefore $G$-invariant.
\end{proof}

\begin{Proposition}\label{prop:embedding_of_KOshs_is_equivariant}
The embeddings given in \eqref{eq:embedding_of_KOshlambda} are $G$-equivariant.
\end{Proposition}
\begin{proof}
Take $\la_1,\la_2,\bR_1,\bR_2$ as in \eqref{eq:embedding_of_KOshlambda}. By Lemma \ref{lem:Bequiv_implies_Gequiv}, it suffices to prove that the embedding is $B$-equivariant%
. Also, as $\fb=\text{Lie}(B)=\fh\oplus\fn$, it suffices to show that this map is equivariant for the action of $\fh$ and $\fn$, but it is clearly $\fh$-equivariant since it preserves weight spaces.~We thus only have to show $\mathfrak{n}$-equivariance.\medskip\par 
For $i \in I$, let $\cE_i^{\bR_1}$ and $\cE_i^{\bR_2}$ be the functors acting on the categories $\O_{sh}^{\lambda_1}(\bR_1)$ and $\O_{sh}^{\lambda_2}(\bR_2)$ (respectively). Take $V$ in the first category and recall from Lemma \ref{lem:chiGT_and_inclusions} that we must have
\begin{equation}\label{eq:equality_of_GTchars}
\chiGT^{\bR_2}(V)=\chiGT^{\bR_1}(V)\ast[\bT]
\end{equation}
for some $\bT\in\Lambda_{\lambda_2-\lambda_1}$. We claim that we also have
\begin{equation}\label{eq:equality_of_GTchars2}
\chiGT^{\bR_2}(\cE_i^{\bR_2}(V))=\chiGT^{\bR_1}(\cE_i^{\bR_1}(V))\ast[\bT].
\end{equation}
Let $\bS\in \Lambda_{\la_2-\mu-\alpha_i}$ and choose $n\in \overline{i}+2\Z$ smaller than every element of $\bR_1$, $\bR_2$, $\bS$ and $\bT$. By Lemma \ref{lem:GT_character_of_Ei_V},
$$\dim W^{\bR_2}_{\bS}(\cE_i^{\bR_2}(V))=\dim W^{\bR_2}_{\bS\cup (n)_i}(V),$$
where we use the notations $W_{\bS}^{\bR_1}$ and $W_{\bS}^{\bR_2}$ to distinguish between taking GT-weight spaces  with respect to $\bR_1$ and $\bR_2$. In addition, by \eqref{eq:equality_of_GTchars} (and as $n$ is smaller than the elements~of~$\bT$), there exists $\bS'\in \Lambda_{\lambda_1-\mu-\alpha_i}$ such that $\bS\cup(n)_i = \bS'\cup\bT\cup(n)_i$ (or, equivalently, $\bS=\bS'\cup\bT$). 
In particular, by applying again \eqref{eq:equality_of_GTchars} and Lemma \ref{lem:GT_character_of_Ei_V}, we get
$$ \dim W^{\bR_2}_{\bS\cup(n)_i}(V)= \dim W^{\bR_2}_{\bS'\cup\bT\cup(n)_i}(V)=\dim W^{\bR_1}_{\bS'\cup (n)_i}(V)=\dim W_{\bS'}(\cE_i^{\bR_1}(V)),$$
and the equality \eqref{eq:equality_of_GTchars2} follows easily. Finally, replacing $V$ by $\mathcal{E}_i^{\bR_1}(V)$ in \eqref{eq:equality_of_GTchars} gives
$$ \chiGT^{\bR_2}(\mathcal{E}_i^{\bR_1}(V))=\chiGT^{\bR_1}(\cE_i^{\bR_1}(V))\ast[\bT]=\chiGT^{\bR_2}(\cE_i^{\bR_2}(V))$$
and using Proposition \ref{prop:injectivity_chiGT} ends the proof. %
\end{proof}

Consider the complex torus $\aT$ whose weight lattice is given by the quotient $\cP=\cB/\Gamma$ of Section \ref{sec:Crystal}.
Explicitly, points of $\aT$ are arrays $(t_{i,a})_{(i,a)\in I\times_2\Z}\in (\C^\times)^{I\times_2 \Z}$ satisfying
\begin{equation*}
\textstyle t_{i,a} \cdot  t_{i,a+2} = \prod_{j\sim i} t_{j,a+1}
\end{equation*}
for all $(i,a)%
$. Recall from Lemma \ref{lem:height_function_identifies_groups} that fixing a height function $\xi$ gives an isomorphism $P\simeq\cP$, and hence an isomorphism $T\simeq \aT$ %
(given by $t\mapsto (t_{i,\xi_i+2s})_{i\in I,s\in \Z}=((c^s\varpi_i)(t))_{i\in I,s\in \Z}$ for $c$ the Coxeter element associated to $\xi$). Recall also from Theorem \ref{thm:BlockDec} and Remark~\ref{rem:Blockawt} that~the category $\O_{sh}$ has~a decomposition indexed by $\cP=\cB/\Gamma$, i.e.
$$\textstyle \smash{\O_{sh}\simeq \bigoplus_{\tau\in \cP} {}_\tau\O_{sh}}$$
where ${}_\tau\O_{sh}$ is the Serre full subcategory of $\O_{sh}$ generated by the simple objects $L(\psi)$ such that $\awt(\psi)=\tau$. We consider the right-action of the torus $\aT$ on the algebra $K_\C(\O_{sh})$ with $\tau$-weight space $K_\C({}_{\tau}\O_{sh})$ for all $\tau\in \cP$. More precisely, given $\tau\in \cP=\cB/\Gamma$ with preimage $\prod_{j=1}^N y_{i_j,a_j}^{n_j}$ in $\cB$, we define%
\begin{equation*}
\textstyle \smash{[V]\cdot (t_{i,a})_{(i,a)\in I\times_2\Z} = \big(\prod_{j=1}^N t_{i_j,a_j}^{n_j}\big)[V]}
\end{equation*}
for all $V$ in $\smash{{}_{\tau}\O_{sh}}$.
\begin{Lemma}\label{lem:subgroup_tau_Osh}
Fix $\tau\in \cP$. Then, in $K_0(\mathcal{O}_{sh})$,
$$\textstyle \bigcup_{\la,\bR} K_0(\mathcal{O}_{sh}^{\la}(\bR))= K_0({}_{\tau}\O_{sh})%
$$ 
where the union runs over all pairs of $\la\in P_+$ and $\bR\in \Z^{\la}$ for which $\awt(y_{\bR})= \tau$.
\end{Lemma}
\begin{proof}
Suppose $\la\in P_+$ and $\bR\in \Z^{\la}$ are such that $\awt(y_{\bR})= \tau$. Then, since $\awt$ is constant on $\cB(\la,\bR)$, we have $\awt(y)= \tau$ for all $y\in \cB(\la,\bR)$ and hence $K_0(\O_{sh}^{\la}(\bR))\subseteq K_0({}_{\tau}\O_{sh})$ by Theorem \ref{th:descend}. This shows the inclusion $\subseteq$.\par  
For the other inclusion, choose $V$ in ${}_\tau\mathcal{O}_{sh}$. By definition of ${}_\tau\mathcal{O}_{sh}$, we can write 
$$ [V] = [L(\psi_1)]+\dots+[L(\psi_k)]\in K_0({}_\tau\mathcal{O}_{sh})$$
where $\psi_1,\dots,\psi_k\in \mathfrak{r}$ satisfy $\awt(\psi_1)=\dots=\awt(\psi_k)=\tau$. Also, for $1\leq s\leq k$, Corollary \ref{coro:every_object_descends_to_some_trunc} implies that $\psi_s\in \cB(\la_s,\bR_s)$ for some $\la_s\in P_+$ and $\bR_s\in \Z^{\la_s}$. Clearly, using again the fact that $\awt$ is constant on each product monomial crystal, 
$$ \awt(y_{\bR_1})=\dots=\awt(y_{\bR_k})=\awt(\psi_1)=\dots=\awt(\psi_k)=\tau,$$
and using repeatedly Corollary \ref{co:directed} shows that there exists $\la\in P_+$ and $\bR\in \Z^{\la}$ such that 
$$ \cB(\la_1,\bR_1)\cup\dots\cup\cB(\la_k,\bR_k)\subseteq \cB(\la,\bR)$$
with $\awt(y_{\bR})=\tau$. In particular, 
$ [V] = [L(\psi_1)]+\dots+[L(\psi_k)]\in K_0(\mathcal{O}_{sh}^{\la}(\bR))$ as desired. This ends the proof.
\end{proof}

The following is the ``gluing result'' announced at the beginning of the section:

\begin{Theorem}\label{thm:construction_of_gaction_on_KOsh}
There exists a (unique) $G$-action on $K_\C(\O_{sh})$ commuting with the $\aT$-action and such that every inclusion $K_\C(\O_{sh}^\lambda(\bR))\hookrightarrow K_\C(\O_{sh})$ is $G$-equivariant.
\end{Theorem}

\begin{proof}
By Proposition \ref{prop:embedding_of_KOshs_is_equivariant} and Lemma \ref{lem:subgroup_tau_Osh}, for $\tau\in \cP$, the subspace $K_{\C}({}_{\tau}\O_{sh})$ is the direct limit of a system of subspaces of $K_{\C}(\mathcal{O}_{sh})$ with $G$-equivariant injective transition maps.~It hence inherits a (unique) $G$-action from the $K_{\C}(\O_{sh}^{\la}(\bR))$'s with $\awt(y_{\bR})=\tau$. This gives~a unique $G$-action on $K_{\C}(\O_{sh})\simeq \bigoplus_{\tau \in \cP}K_{\C}(_{\tau}\O_{sh})$ that commutes with the $A$-action.  
\end{proof}

\begin{Corollary}\label{cor:Invariants} 
The subalgebras of (left-)invariants in $K_{\C}(\mathcal{O}_{sh})$ are given by:
\begin{enumerate}
\item ${}^TK_{\C}(\mathcal{O}_{sh})=K_{\C}(\mathcal{O}_0)$, 
\item ${}^NK_{\C}(\mathcal{O}_{sh})=K_{\C}(\mathscr{C}_{sh})$ for $\mathscr{C}_{sh}\subseteq \mathcal{O}_{sh}$ the category of finite-dimensional objects,~and
\item ${}^GK_{\C}(\mathcal{O}_{sh})=K_{\C}(\mathscr{C}_0)$ for $\mathscr{C}_0=Y\fmod$ the category of finite-dimensional $Y$-modules.
\end{enumerate}
\end{Corollary}

\begin{proof}
It suffices to show the statements in $K_{\C}(\mathcal{O}_{sh}^{\la}(\bR))$ for some fixed $\la\in P_+$ and $\bR\in \Z^{\la}$. For this, note that (1) follows directly from the fact that $K_{\C}(\mathcal{O}_{\mu}^{\la}(\bR))$ is the $\mu$-weight space for the $T$-action on $K_{\C}(\mathcal{O}_{sh}^{\la}(\bR))$. Also, letting $\mathscr{C}_{sh}^{\la}(\bR)\subseteq \mathcal{O}_{sh}^{\la}(\bR)$ be the full subcategory~of finite-dimensional objects, we can use \cite[Theorem 5.21]{kamnitzer2019category} to deduce isomorphisms
\begin{equation*}
K_{\C}(\mathscr{C}_{sh}^{\la}(\bR))\simeq {}^{N}K_{\C}(P^\bR\fmod) \simeq {}^{N}K_{\C}(\mathcal{O}_{sh}^{\la}(\bR))
\end{equation*}
for which the composition from left to right comes from the inclusion of $\mathscr{C}_{sh}^{\la}(\bR)$ into $\mathcal{O}_{sh}^{\la}(\bR)$. This shows statement (2). The last statement is a combination of the first two.
\end{proof}
\subsection{Equivariance of multiplication}\label{sec:equivMult} We now want to show that the map
\begin{equation*}
K_\C(\O_{sh})\otimes K_\C(\O_{sh})\xrightarrow[]{\otimes} K_\C(\O_{sh})
\end{equation*}
induced by tensor product is a morphism of both $G$-representations and $\aT$-representations. We start by proving this at the level of truncations via Theorem \ref{thm:truncoprod}. For this, fix $\lambda_1,\lambda_2\in P_+$ with $\bR_1\in \C^{\lambda_1}$ and $\bR_2\in\C^{\lambda_2}$. Also let $\lambda=\lambda_1+\lambda_2$ and $\bR=\bR_1\cup\bR_2$.

\begin{Proposition}\label{prop:mult_is_equiv_truncation}
The multiplication map 
\begin{equation}\label{eq:mult_KC_level_of_truncations}
K_\C(\O_{sh}^{\lambda_1}(\bR_1))\otimes K_\C(\O_{sh}^{\lambda_2}(\bR_2))\xrightarrow[]{\otimes} K_\C(\O_{sh}^{\lambda}(\bR))
\end{equation}
is equivariant for both the left $G$-action and the right $\aT$-action.
\end{Proposition}

\begin{proof}
Observe that \eqref{eq:mult_KC_level_of_truncations} has the right codomain by Theorem \ref{thm:truncoprod} and is trivially $A$-equivariant as $\awt(y_{\bR})=\awt(y_{\bR_1}) \awt(y_{\bR_2})$ (see also the proof of Lemma \ref{lem:subgroup_tau_Osh}). To prove that it is~also $G$-equivariant, we use a similar technique as in the proof of Proposition \ref{prop:embedding_of_KOshs_is_equivariant} and only %
show $\fn$-equivariance (since $\fh$-equivariance is immediate). \medskip\par 
Take $\mu_1,\mu_2\in P$ with respective objects $V_1$ and $V_2$ in $\mathcal{O}_{\mu_1}^{\la_1}(\bR_1)$ and $\mathcal{O}_{\mu_2}^{\la_2}(\bR_2)$. Also let $i\in I$ and let $\nu_1=\la_1-\mu_1$ with $\nu_2=\la_2-\mu_2$ and $\nu=\nu_1+\nu_2$. Given $\bS\in \Lambda_{\nu-\alpha_i}$ with $n\in \overline{i}+2\Z$ smaller than every element of $\bS$ and $\bR$, Lemma \ref{lem:GT_character_of_Ei_V} gives
\begin{equation*}
\dim W_{\bS}^{\bR}\big(\cE_i(V_1\otimes V_2)\big)=\dim W_{\bS\cup (n)_i}^{\bR}(V_1\otimes V_2)
\end{equation*}
and Proposition \ref{prop:GTchar_is_multiplicative} allows us to write the right-hand side of this equation as
\begin{equation*}
\textstyle \dim W_{\bS\cup (n)_i}^{\bR}(V_1\otimes V_2)=\sum_{\bS_1\cup\bS_2=\bS\cup (n)_i} \dim W_{\bS_1}^{\bR_1}(V_1)\cdot \dim W_{\bS_2}^{\bR_2}(V_2)
\end{equation*}
where the sum runs over all ways of writing $\bS\cup (n)_i$ as the union of GT-weights $\bS_1\in \Lambda_{\nu_1}$ and $\bS_2\in \Lambda_{\nu_2}$. Clearly, one can split the above sum as
\begin{equation*}
\sum_{\bS_1'\cup\bS_2=\bS} \dim W_{\bS_1'\cup (n)_i}^{\bR_1}(V_1)\cdot \dim W_{\bS_2}^{\bR_2}(V_2)+\sum_{\bS_1\cup\bS_2'=\bS} \dim W_{\bS_1}^{\bR_1}(V_1)\cdot \dim W_{\bS_2'\cup(n)_i}^{\bR_2}(V_2)
\end{equation*}
where the first (resp.~second) sum runs over all ways of writing $\bS$ as the union of $\bS_1'\in \Lambda_{\nu_1-\alpha_i}$ (resp.~$\bS_1\in \Lambda_{\nu_1}$) and $\bS_2\in \Lambda_{\nu_2}$ (resp.~$\bS_2'\in \Lambda_{\nu_2-\alpha_i}$). In particular, using again Lemma \ref{lem:GT_character_of_Ei_V} with Proposition \ref{prop:GTchar_is_multiplicative} gives (since $n$ is smaller than every element of $\bS$, $\bR_1$ and $\bR_2$)
\begin{align*}
\sum_{\bS_1'\cup\bS_2=\bS} \dim W_{\bS_1'\cup (n)_i}^{\bR_1}(V_1)\cdot \dim W_{\bS_2}^{\bR_2}(V_2)&=
\sum_{\bS_1'\cup\bS_2=\bS} \dim W_{\bS_1'}^{\bR_1}(\cE_i(V_1))\cdot \dim W_{\bS_2}^{\bR_2}(V_2)\\
&=\dim W^{\bR}_{\bS}(\cE_i(V_1)\otimes V_2)
\end{align*}
and, similarly,
\begin{align*}
\sum_{\bS_1\cup\bS_2'=\bS} \dim W_{\bS_1}^{\bR_1}(V_1)\cdot \dim W_{\bS_2'\cup (n)_i}^{\bR_2}(V_2)
&=\dim W^{\bR}_{\bS}(V_1\otimes \cE_i(V_2)).
\end{align*}
Combining everything, we get that, for all $\bS\in \Lambda_{\nu-\alpha_i}$,
$$ \dim W_{\bS}^{\bR}\big(\cE_i(V_1\otimes V_2)\big)= \dim W^{\bR}_{\bS}(\cE_i(V_1)\otimes V_2)+\dim W^{\bR}_{\bS}(V_1\otimes \cE_i(V_2))$$
so that
$$\chiGT^\bR\big(\,\cE_i(V_1\otimes V_2)\,\big)=\chiGT^\bR\big(\,[\cE_i(V_1)\otimes V_2]+[V_1\otimes \cE_i(V_2)]\,\big).$$
The desired result then follows from applying Proposition \ref{prop:injectivity_chiGT}.
\end{proof}

\begin{Theorem}\label{thm:mult_in_KOsh_is_G_equiv}
Multiplication in $K_\C(\O_{sh})$ is equivariant for both the $G$ and $\aT$ actions.
\end{Theorem}

\begin{proof}
It is enough to prove the theorem for pairs of simple classes in $K_0(\mathcal{O}_{sh})$, but then the result follows from Proposition \ref{prop:mult_is_equiv_truncation} (thanks to Theorem \ref{th:descend}). This finishes the proof.
\end{proof}

\subsection{A generating set} We end this section by identifying generators for the ring $K_0(\O_{sh})$. To do so, we first show that multiplication is surjective at the level of truncations, which is very natural from the %
product monomial crystal picture. We use the notation of Section~\ref{sec:equivMult}.
\begin{Proposition}\label{prop:mult_is_surj_truncation}
The multiplication map \eqref{eq:mult_KC_level_of_truncations}
is surjective.
\end{Proposition}

\begin{proof}
Consider the $\Z$-bases 
\begin{equation*}
\big\{[L(\xi_1)]\otimes [L(\xi_2)]\;|\;\xi_1\in\cB(\lambda_1,\bR_1),\,\xi_2\in\cB(\lambda_2,\bR_2)\big\}\text{ and }\big\{[L(\psi)]\;|\;\psi\in\cB(\lambda,\bR)\big\}
\end{equation*}
of $K_0(\O_{sh}^{\lambda_1}(\bR_1))\otimes K_0(\O_{sh}^{\lambda_2}(\bR_2))$ and $K_0(\O_{sh}^{\lambda}(\bR))$ (respectively).
Remark that 
$$\cB(\la,\bR)=\cB(\la_1,\bR_1)\cB(\la_2,\bR_2)$$
and order the elements of $\cB(\lambda,\bR)$ according to Nakajima's partial order, i.e.
$$\cB(\lambda,\bR)=\{\psi_1,\dots,\psi_p\}$$ 
where $i\leq j$ whenever $\psi_i\succeq \psi_j$. For all $1\leq i\leq p$, choose $(\psi_{1,i},\psi_{2,i})\in\cB(\lambda_1,\bR_1)\times \cB(\lambda_2,\bR_2)$ such that $\psi_i=\psi_{1,i}\psi_{2,i}$. Then Corollary \ref{cor:tensor_product_decomp_in_Osh} shows that the integer ($p{\times}p$)-matrix $A=(a_{ij})$ with coefficients given by
\begin{equation*}
[L(\psi_{1,j})\otimes L(\psi_{2,j})]=\textstyle\sum_{i=1}^p a_{ij}[L(\psi_i)]
\end{equation*} 
is upper unitriangular, and thus invertible over $\Z$. This completes the proof.
\end{proof}
\begin{Theorem}\label{thm:generators_of_KOsh}
The ring $K_0(\O_{sh})$ is generated by the classes
\begin{equation*}
\smash{\{[L(\psi)]\,|\, \psi\in \cB(\varpi_i,a)\}_{(i,a)\in I\times_2\Z}.}
\end{equation*}
\end{Theorem}\vspace*{-2mm}

\begin{proof}
This follows easily from \eqref{eq:MonCrys}, Corollary \ref{coro:every_object_descends_to_some_trunc} and Proposition \ref{prop:mult_is_surj_truncation}.
\end{proof}

\begin{Rem}
The above generating set contains the classes of chamber modules defined in Section \ref{sec:Chambermodules}, together with the classes associated to non-extremal weight spaces (if any)%
. 
\end{Rem}

\section{\texorpdfstring{$G$}{G}-action and multiplication}\label{sec:description_of_mult}

The goal of this section is to use the actions of $G$ and $A$ on the algebra $K_{\C}(\mathcal{O}_{sh})$ to show new results on its multiplication. The main byproducts of this pursuit are various~inclusions of the coordinate ring of the base affine space $N_-\bbslash G$ in $K_0(\mathcal{O}_{sh})$ and proofs of conjectures of Frenkel--Hernandez and Geiss--Hernandez--Leclerc on extended $QQ$-systems (cf.~Section \ref{sec:Intro}). We start by relating distinct chamber modules through the $G$-action.
\subsection{Chamber modules and the $G$-action}\label{subsec:chamber_modules_and_G_action}
For every $i\in I$, there is a group morphism $\tau_i:\SL_2\to G$ integrating the Lie algebra map $\fsl_2\to \fg$ given by $e\mapsto e_i$, $f\mapsto f_i$ and $h\mapsto h_i$. This morphism produces a distinguished lift of the simple reflection $s_i\in W$ to $G$ via
\begin{equation}\label{eq:lift_of_simple_refl}
\dot{s}_i=\tau_i\big(\adjustbox{scale=0.76}{$\begin{pmatrix}
0 & -1\\
1 & 0
\end{pmatrix}$}
\big).
\end{equation}
Given $w\in W$ with a choice of reduced expression $\underline{w}=(s_{i_1},\dots,s_{i_n})$, consider $\dot{w}=\dot{s}_{i_1}\dots\dot{s}_{i_n}$. It is well-known that $\dot{w}$ is independent of the choice of reduced expression for $w$ and~that $\dot{w}T=w\in N(T)/T$. Also, the subgroup $\smash{\widetilde{W}}=\langle \dot{s}_i\rangle_{i\in I}\subseteq G$ is a 2-covering of the Weyl group (in the sense that the generators satisfy the braid relations with $(\dot{s}_i)^4=e$). Finally, given a $G$-module $V$ with $v\in V$ a vector of weight $\mu\in P$ satisfying $\langle \mu,\alpha_i^{\vee}\rangle=n\geq 0$, the lift \eqref{eq:lift_of_simple_refl} is such that $\smash{f^{(n)}_i} v=\dot{s}_i v$, where $\smash{f^{(n)}_i=\tfrac{1}{n!}f_i^n}$ is a divided power (see \cite[Section~3.8]{kac1990infinite} for details about this and the above statements).

Recall from Section \ref{sec:Chambermodules} that, given $(i,a)\in I\times_2\Z$ with a chamber weight $\gamma\in W\varpi_i$, the category $\O^{\varpi_i}_{\gamma}(a)$ has (up to isomorphism) a unique irreducible object $L_{\gamma,a}$ that we call the \textit{$\gamma$-chamber module of spectral parameter $a$}.

\begin{Lemma}\label{lem:Gaction_and_chamber_modules}
Take $(i,a)\in I\times_2\Z$ and $\gamma\in  W\varpi_i$. Also, fix $j\in I$ such that $\langle\gamma,\alpha_j^\vee\rangle=n>0$. 
Then, the class of the chamber module $L_{\gamma,a}$ in $K_0(\O^{\varpi_i}_{sh}(a))%
$ satisfies 
$[L_{s_j\gamma,a}]=\dot{s}_j[L_{\gamma,a}].$
\end{Lemma}
To prove this, we use \cite[Theorem 6.4]{chuang2008derived} which~shows that, if
$\textstyle \mathcal{C}=\bigoplus_{\mu\in P} \mathcal{C}_\mu$
is an abelian category endowed with a categorical $\mathfrak{sl}_2$-action such that $V=K_{\C}(\mathcal{C})$ is finite-dimensional, then the bounded derived categories of $\mathcal{C}_{\mu}$ and $\mathcal{C}_{s_1\mu}$ are equivalent for all~$\mu$~(via~an~explicit complex that lifts the isomorphism $\dot{s}_1:V_{\mu}\to V_{s_1\mu}$). Moreover, when $V_{\mu+\alpha_1}=0$, the above derived equivalence reduces to an honest equivalence $\mathcal{C}_{\mu}\simeq \mathcal{C}_{s_1\mu}$ and $\mathfrak{sl}_2$-restriction allows~us to deduce the following result (which could also be obtained using \cite[Theorem 5.9]{khovanov2012extended}):
\begin{Lemma}[\cite{chuang2008derived,khovanov2012extended}]\label{lem:ChuangRouquier} Let $\mathcal{C}=\bigoplus_{\mu\in P}\mathcal{C}_{\mu}$ be an abelian category with a categorical $\fg$-action. %
Fix $\mu\in P$ a weight of $V=K_{\C}(\mathcal{C})$ with $i\in I$. Suppose that $\dim V<\infty$ and that $\mu+\alpha_i$ is not a weight of $V$. Then the categories $\mathcal{C}_{\mu}$ and $\mathcal{C}_{s_i\mu}$ are equivalent. In addition,~if $\langle \mu,\alpha_i^{\vee}\rangle=n>0$, the %
equivalence is given by the ``divided power functors''
$$ \smash{\mathcal{E}_i^{(n)}:\mathcal{C}_{s_i\mu}\to\mathcal{C}_{\mu}\ \text{ and }\ \mathcal{F}_i^{(n)}:\mathcal{C}_{\mu}\to\mathcal{C}_{s_i\mu}}$$
 of \cite[Section 4.1.1]{rouquier20082} that categorify the action of the divided powers $\smash{e_i^{(n)}}$ and $\smash{f_i^{(n)}}$ on~%
$V$.
\end{Lemma}
\begin{proof}[Proof of Lemma \ref{lem:Gaction_and_chamber_modules}]
By Lemma \ref{lem:ChuangRouquier}, the functor 
$$\smash{\cF_j^{(n)}:\O^{\varpi_i}_{\gamma}(a)\to \O^{\varpi_i}_{s_j\gamma}(a)}$$
gives an equivalence of categories and thus sends the unique simple module $L_{\gamma,a}$ of $\O^{\varpi_i}_{\gamma}(a)$ to the unique simple module $L_{s_j\gamma,a}$ of $\O^{\varpi_i}_{s_j\gamma}(a)$. Passing to $K_0(\mathcal{O}_{sh}^{\varpi_i}(a))$, we deduce that
$$\smash{[L_{s_j\gamma,a}]=[\cF_j^{(n)}(L_{\gamma,a})]=f_j^{(n)}[L_{\gamma,a}]=\dot{s}_j[L_{\gamma,a}],} $$
as required. %
\end{proof}

\begin{Corollary}\label{cor:W_action_on_chamber_modules}
Let $w\in W$. Then, $\dot{w}[L_{\varpi_i,a}]=[L_{w\varpi_i,a}]$ in $K_0(\O_{sh}^{\varpi_i}(a))$.
\end{Corollary}

\begin{Rem}\label{rem:WactionEell} By Theorem \ref{thm:construction_of_gaction_on_KOsh}, the group $\widetilde{W}=\langle \dot{s}_i \rangle_{i\in I}\subseteq G$ acts on $K_\C(\O_{sh})$. Transporting this action to $\chi_{\ell}(K_{\C}(\O_{sh}))\subseteq \mathcal{E}_{\ell}$, we see, using Example \ref{ex:Inflsl3}, that
\begin{equation*}
\dot{s}_i\cdot [\sfPsi_{j,a}]%
=\begin{cases}
[\sfPsi_{j,a}] & \text{if }i\neq j,\\
[\sfPsi_{i,a-2}](1+\sfA_{i,a-2}^{-1}+\sfA_{i,a-2}^{-1}\sfA_{i,a-4}^{-1}+\dots) & \text{if }i=j.
\end{cases}
\end{equation*}
This recovers (a rational analogue of) the results of \cite[Proposition~3.17]{frenkel2024extended} and \cite[Section~7.3]{geiss2024representations} %
(see also \cite{frenkel2022weyl}), where a $\smash{\widetilde{W}}$-action on (a completion of) $\cE_\ell$ is also constructed.
\end{Rem}
\subsection{Height functions and the coordinate ring of the base affine space}\label{subsec:base_affine_space_and_Osh}
We now take a slight geometric tangent in preparation for the upcoming Sections \ref{sec:bott-samselson-background}--\ref{sec:bi-infinite-bott-samelson}. %
\medskip\par We adopt the convention that geometric objects are endowed with \textit{right-actions} and~that their ring of functions are equipped with the induced \textit{left-actions} (and vice-versa). We start by clarifying the impact of these conventions.\medskip\par
The group $G$ acts on itself by left and right translations%
. This induces a $(G,G)$-bimodule structure on the coordinate ring $\C[G]$ via
\begin{equation*}
(g_1 \cdot f\cdot g_2)(x)=f(g_2xg_1)
\end{equation*}
for $f\in \C[G]$ and $g_1,g_2,x\in G$. The Peter-Weyl theorem shows that, as $(G,G)$-bimodules, 
\begin{equation}\label{eq:iso_peter_weyl}
\textstyle\C[G]\simeq \bigoplus_{\lambda\in P_+} V(\lambda)^\ast \otimes V(\lambda)
\end{equation}
with $V(\lambda)$ the simple left $G$-module of highest weight $\lambda$ and where $V(\lambda)^\ast:=\Hom_\C(V(\lambda),\C)$ acquires the structure of a simple right $G$-module of highest weight $\lambda$. Note that the $(G,G)$-bimodule structure on
$V(\la)^{\ast}\otimes V(\la)$ which is compatible with \eqref{eq:iso_peter_weyl} is
\begin{equation*}
g_1 \cdot \langle v_1^\ast, (\trou)v_2\rangle \cdot g_2=\langle v_1^\ast g_2,g_1(\trou)v_2\rangle,
\end{equation*}
where $ \langle v_1^\ast, (\trou)v_2\rangle\in \C[G]$ is the matrix coefficient corresponding to $v_1^\ast\otimes v_2\in V(\la)^\ast\otimes V(\la)$.\medskip\par
Fix $\la\in P_+$ and let $v_\lambda^\ast\in V(\lambda)^\ast$ be the vector of weight $\la$ that pairs as 1 with the highest weight vector $v_{\la}\in V(\la)$ that we chose in Section \ref{sec:Intro}. Note that $v_{w\lambda}=\dot{w}v_\lambda$ is (up to scalar) the unique non-zero vector in $V(\la)_{w\la}$ and remark that $v_{w\lambda}^\ast:=v_\lambda^\ast (\dot{w})^{-1}$ satisfies
$$ \langle v_{w\la}^\ast,v_{w\la}\rangle=\langle v_{\la}^\ast(\dot{w})^{-1},\dot{w}v_{\la}\rangle=\langle v_{\la}^\ast,v_{\la}\rangle=1.$$
In addition, $v_{w_1\la}=v_{w_2\la}$ and $v_{w_1\la}^*=v_{w_2\la}^*$ whenever $w_1,w_2\in W$ satisfy $w_1\la=w_2\la$, and~the vectors $v_{w\la}$ and $v_{w\la}^*$ thus depend only on the chamber weight $w\la$ (and not the pair $(\la,w)$). 

Finally, recall Berenstein--Zelevinsky's \textit{generalized minors} from \cite{berenstein1997total,fomin199double}, which are the functions $\{\Delta_{w_1\varpi_i,w_2\varpi_i}\}_{i\in I,w_1,w_2\in W} \subseteq \C[G]$ defined by 
\begin{equation}\label{eq:Minors}
\Delta_{w_1\varpi_i,w_2\varpi_i}(g)=\langle v_{\varpi_i}^\ast,(\dot{w}_1)^{-1}g\dot{w}_2 v_{\varpi_i}\rangle=\langle v_{w_1\varpi_i}^\ast, gv_{w_2\varpi_i}\rangle.
\end{equation} 
\begin{Def} The \textit{base affine space} $N_-\bbslash G$ is the affine closure of the variety $N_-\backslash G$ or, equivalently, the affine scheme given by
$$\smash{N_-\bbslash G = \Spec(\C[G]^{N_-})},$$
where $\C[G]^{N_-}$ is the ring of invariants with respect to the right-action of $N_-$ on $\C[G]$.
\end{Def}
The ring of functions of the base affine space can be categorified using %
shifted Yangians. To explain this, observe that \eqref{eq:iso_peter_weyl} gives a $G$-equivariant algebra isomorphism
\begin{equation}\label{eq:PeterWeylBaseAff}
\textstyle  \C[N_-\bbslash G]=\C[G]^{N_-}\simeq \bigoplus_{\lambda\in P_+} V(\lambda)
\end{equation}
where the %
sum %
is endowed with the \textit{Cartan product}, i.e.~%
for $v_1\in V(\la_1)$ and $v_2\in V(\la_2)$,
$$ v_1\cdot v_2 = \pi_{\la_1,\la_2}(v_1\otimes v_2),$$
with $\pi_{\lambda_1,\lambda_2}$ the unique $G$-equivariant map $V(\lambda_1)\otimes V(\lambda_2)\twoheadrightarrow V(\lambda_1+\lambda_2)$ such that $$\pi_{\la_1,\la_2}(v_{\lambda_1}\otimes v_{\lambda_2})= v_{\lambda_1+\lambda_2}.$$\par %
Choose a height function $\xi$ and recall the corresponding isomorphisms $\aT\simeq T$ and $\cP\simeq P$ (see Section \ref{subsec:equiv_of_inclusions}). By Theorem \ref{thm:charaterization_max_sing_crystals}, Corollary \ref{coro:maximally_singular_categorifies_irrep} and Theorem \ref{thm:equiv_of_cats_parity_and_Osh}, given $\la\in P_+$ and $\bR_{\xi}\in \Z^{\la}$ the \textit{set of parameters given by $\xi$} (in the sense of Section \ref{subsec:maximally_sing_crystals}), the category %
\begin{equation*}
\O_{sh}^{\lambda}(\xi):=\O_{sh}^\lambda(\bR_\xi)
\end{equation*}
is such that $K_{\C}(\mathcal{O}_{sh}^{\la}(\xi))\simeq K_{\C}(R^{\la}\fmod)\simeq V(\la)$ (as $G$-modules).\medskip\par
We thus have an injective $G$-equivariant morphism
\begin{equation}\label{eq:InjCartanProd}
\textstyle \smash{\bigoplus_{\la\in P_+}}V(\la)\simeq \smash{\bigoplus_{\la\in P_+}\smash{K_{\C}(\mathcal{O}_{sh}^{\la}(\xi))}} \hookrightarrow K_{\C}(\O_{sh}),
\end{equation} 
which we now show is an algebra morphism. 
\begin{Theorem}\label{thm:subalgebra_iso_base_affine_space} For each height function $\xi$, the $G$-equivariant vector space isomorphism 
\begin{equation}\label{eq:CartanCyc2}
\textstyle \smash{\bigoplus_{\la\in P_+}\smash{K_{\C}(\mathcal{O}_{sh}^{\la}(\xi))}} \simeq \smash{\bigoplus_{\la\in P_+}V(\la)}
\end{equation}
is an algebra isomorphism (where multiplication for the right sum is Cartan product).~Thus, %
\eqref{eq:InjCartanProd} gives a $G$-equivariant algebra embedding of $\C[N_-\bbslash G]$ in $K_{\C}(\mathcal{O}_{sh})$.
\end{Theorem}
\begin{proof}
The only thing left to prove is the compatibility of \eqref{eq:CartanCyc2} with multiplication. For~this, fix $\la_1,\la_2\in P_+$ and %
use Theorem \ref{thm:mult_in_KOsh_is_G_equiv} to deduce that all the maps in the diagram
\begin{equation}\label{eq:multBaseAff} 
\adjustbox{scale=0.9}{
\begin{tikzcd}[column sep = 2em]
V(\la_1)\otimes V(\la_2)\arrow[d,"\simeq"] \arrow[r,"\pi_{\la_1,\la_2}"]& V(\la_1+\la_2)\\
K_{\C}(\mathcal{O}_{sh}^{\la_1}(\xi))\otimes K_{\C}(\mathcal{O}_{sh}^{\la_2}(\xi))\arrow[r,"\otimes"] & K_{\C}(\mathcal{O}_{sh}^{\la_1+\la_2}(\xi))\arrow[u,"\simeq"]
\end{tikzcd}}
\end{equation}
are $G$-equivariant. Also, by Theorem \ref{thm:dual_canonical_KLR}, the composition of the bottom horizontal~map~and the left vertical isomorphism sends $v_{\la_1}\otimes v_{\la_2}$ to the class of the tensor product of the unique simple modules of $\smash{\mathcal{O}^{\la_1}_{\la_1}(\xi)}$ and $\smash{\mathcal{O}^{\la_2}_{\la_2}(\xi)}$, which are clearly both %
$1$-dimensional (as $\smash{Y_{\la}^{\la}(\bR)}\simeq \C$ for all $\la\in P_+$ and $\bR\in \Z^{\la}$). In particular, the above tensor product is %
the unique simple module inside $\smash{\mathcal{O}^{\la_1+\la_2}_{\la_1+\la_2}(\xi)}$, and Theorem \ref{thm:dual_canonical_KLR} again (with the unicity property of $\pi_{\la_1,\la_2}$)~shows that \eqref{eq:multBaseAff} is commutative, which is precisely what we needed to show.
\end{proof}

\begin{Rem}\label{rem:AequivBaseAff} The category $\mathcal{O}_{sh}^{\la}(\xi)$ belongs to the block associated to $\la\in P\simeq \cP= \mathcal{B}/\Gamma$~in Theorem \ref{thm:BlockDec}, and thus \eqref{eq:InjCartanProd} is $A$-equivariant (for the usual left-action of $A\simeq T$ on $N_-\bbslash G$).
\end{Rem}

\begin{Rem}\label{rem:base_affine_space_map_factors}
By the above, choosing a height function $\xi$ for $\fg$ gives a $G$-equivariant map $\Spec K_\C(\O_{sh}) \rightarrow N_-\bbslash G$. We will show during the proof of Lemma \ref{lem:Zhat-torsor} (see also Theorem \ref{thm:GeoSubalgebrasK0}) that this map factors through the embedding of $N_-\backslash G$ in its affine closure $N_-\bbslash G$.
\end{Rem}

Let $(i,a)\in I\times_2\Z$ and fix a height function $\xi$ such that $\xi(i)=a$. Then, by Theorem~\ref{thm:dual_canonical_KLR}, the class in $K_0$ of the %
module $L_{\varpi_i,a}$ is sent, by \eqref{eq:CartanCyc2}, to the highest weight vector $v_{\varpi_i}\in V(\varpi_i)$ or, using \eqref{eq:PeterWeylBaseAff}, to the minor $\Delta_{\varpi_i,\varpi_i}\in \C[N_-\bbslash G]=\C[G]^{N_-}$ (cf.~\eqref{eq:Minors}). Hence, given $w\in W$, one can use Corollary \ref{cor:W_action_on_chamber_modules} to identify the class
$$ [L_{w\varpi_i,a}]=\dot{w}[L_{\varpi_i,a}]\in K_0(\mathcal{O}_{sh}^{\varpi_i}(\xi))$$
with the generalized minor 
$$\Delta_{w\varpi_i}:=\Delta_{\varpi_i,w\varpi_i}.$$
This implies a relation in the ring $\cE_\ell$ which was observed in \cite[Proposition~7.7]{geiss2024representations} (and which corresponds to a mutation relation in the cluster algebra studied therein).
\begin{Corollary}\label{cor:RelationBaseAff}
Fix $ i, j \in I $ with $ i \sim j $ and let $ w \in W $ be such that $ ws_i > w$ and $ws_j > w$. Finally, choose a spectral parameter $a\in \smash{\overline{i}+2\Z}$. Then, in $\smash{K_0(\mathcal{O}_{sh})}$,
$$
[L_{w s_i \varpi_i , a} ] [L_{w s_j \varpi_j, a+1}] = [L_{w \varpi_i, a}][L_{w s_i s_j \varpi_j, a+1}] + [L_{w s_j s_i \varpi_i, a}][L_{w \varpi_j, a+1}].
$$
\end{Corollary}

\begin{proof}
Fix a height function $\xi$ %
with $\xi(i)=a$ and $\xi(j)=a+1$. Then the result follows from the Plücker relations of %
\cite[Corollary~6.6]{berenstein1997total} after using Theorem \ref{thm:subalgebra_iso_base_affine_space} %
and the %
above.
\end{proof}

\begin{Corollary} \label{co:Rel12}
Fix $ i, j \in I $ with $ i \sim j $ and choose $ a \in \smash{\overline{i}+2\Z}$. Then
\begin{enumerate}
\item The multiplication $\smash{K_\C(\O^{\varpi_i}_{sh}(a))} \otimes \smash{K_\C(\O^{\varpi_i}_{sh}(a))} \rightarrow K_\C(\O_{sh}) $ annihilates all $G$-isotypic components other than $ V(2\varpi_i) $.
\item The multiplication $ \smash{K_\C(\O^{\varpi_i}_{sh}(a))}\otimes \smash{K_\C(\O^{\varpi_j}_{sh}(a+1))} \rightarrow K_\C(\O_{sh}) $ annihilates all $G$-isotypic components other than $ V(\varpi_i + \varpi_j) $.
\end{enumerate}
\end{Corollary}

\begin{proof}
Fix (again) a height function $ \xi $ with $(\xi(i),\xi(j)) = (a,a+1)$ and use Theorem \ref{thm:subalgebra_iso_base_affine_space}.
\end{proof}

\subsection{The extended $QQ$-system}\label{subsec:cat_O_of_QQ}
Fix $i\in I$ and note that, as $G$ is simply-laced,
\begin{equation}\label{eq:varpi_svarpi}
\textstyle\varpi_i+s_i\varpi_i=2\varpi_i-\alpha_i=\smash{\sum_{j\sim i}\varpi_j}.
\end{equation}
Consider the $G$-module $\smash{V(\varpi_i)^{\otimes 2}}$. Then, the space $\smash{(V(\varpi_i)^{\otimes 2})_{2\varpi_i-\alpha_i}}$ is isomorphic to
$$ (V(\varpi_i)_{\varpi_i}\otimes V(\varpi_i)_{s_i\varpi_i})\oplus (V(\varpi_i)_{s_i\varpi_i}\otimes  V(\varpi_i)_{\varpi_i})\simeq \C^2$$
and contains a unique (up to scaling) $N$-invariant vector which we can fix to be 
$$v_{\varpi_i}\wedge v_{s_i\varpi_i}=v_{\varpi_i}\otimes v_{s_i\varpi_i}-v_{s_i\varpi_i}\otimes v_{\varpi_i}\in \smash{\wedge^2 V(\varpi_i)}\subseteq V.$$
There is thus a unique $G$-equivariant map 
\begin{equation}\label{eq:Gequivariant_map_wedge}
\iota_1 : V(2\varpi_i-\alpha_i)\hookrightarrow V(\varpi_i)\otimes V(\varpi_i)
\end{equation}
with $\iota_1(v_{2\varpi_i-\alpha_i})= v_{\varpi_i}\otimes v_{s_i\varpi_i}-v_{s_i\varpi_i}\otimes v_{\varpi_i}$. %
There is also %
a unique $G$-equivariant map
\begin{equation}\label{eq:Gequivariant_map_from_identity}
\textstyle \iota_2:V(2\varpi_i-\alpha_i)\hookrightarrow \smash{\bigotimes_{j\sim i}V(\varpi_j)}
\end{equation}
satisfying $\iota_2(v_{2\varpi_i-\alpha_i})=\bigotimes_{j\sim i} v_{\varpi_j}$ (see~\eqref{eq:varpi_svarpi}).\medskip\par 
Now, fix $a\in \overline{i}+2\Z$ and consider the category $\O^{2\varpi_i}_{sh}(a,a+2)$ (where the set of parameters is $(a)_i\cup(a+2)_i$). Denote by $\cB(2\varpi_i,\{a,a+2\})$ the corresponding product monomial crystal and denote by $\cB(2\varpi_i-\alpha_i,\{a+1\})$ the crystal with set of parameters $\cup_{j\sim i} (a+1)_j$.

\begin{Lemma}\label{lem:description_of_crystals_QQ_system}
The two monomials contained in $\cB(2\varpi_i,\{a,a+2\})_{2\varpi_i-\alpha_i}$ are 
\begin{equation*}
\textstyle \prod_{j\sim i} y_{j,a+1}=y_{i,a}\tilde{f}_i (y_{i,a+2})\hspace{1em}\text{and}\hspace{1em}\tfrac{y_{i,a+2}}{y_{i,a-2}}\prod_{j\sim i}y_{j,a-1}=y_{i,a+2}\tilde{f}_i (y_{i,a}).
\end{equation*}
Moreover, there is an inclusion of sets $\cB(2\varpi_i-\alpha_i,\{a+1\})\subseteq\cB(2\varpi_i,\{a,a+2\})$, and thus an inclusion of categories
$$\smash{\O^{2\varpi_i-\alpha_i}_{sh}(a+1)}\subseteq\smash{\O^{2\varpi_i}_{sh}(a,a+2)}.$$
\end{Lemma}
\begin{proof}
The first statement is easy to check and the second follows  from Theorem \ref{th:TFAE} (see also the beginning of Section \ref{sec:glueing_of_g_action}).
\end{proof}
Combining Lemma \ref{lem:description_of_crystals_QQ_system} with \eqref{eq:Gequivariant_map_wedge}--\eqref{eq:Gequivariant_map_from_identity} and using the isomorphism $V(\varpi_i)\simeq \smash{K_{\C}(\O_{sh}^{\varpi_i}(a))}$ gives the heptagonal 
diagram of $G$-equivariant maps depicted below.
\begin{center}
\adjustbox{scale=0.9,center}{
\begin{tikzpicture}[
  baseline=(current bounding box.center),
  every node/.style={inner sep=2pt},
  arr/.style={->, shorten >=4pt, shorten <=4pt}
]
  \node (A1) at (0,2.3) {$V(2\varpi_i-\alpha_i)$};

  \node (A2) at (-3.6,1.5) {$ \bigotimes_{j\sim i} V(\varpi_j)$};
  \node (A7) at (3.6,1.5) {$V(\varpi_i)\otimes V(\varpi_i)$};

  \node (A3) at (-3.6,0.1) {$ \bigotimes_{j\sim i} K_\C(\O_{sh}^{\varpi_j}(a+1))$};
  \node (A6) at (3.6,0.1) {$K_\C(\O_{sh}^{\varpi_i}(a))\otimes K_\C(\O_{sh}^{\varpi_i}(a+2))$};

  \node (A4) at (-2.4,-1.15)
  {$K_\C(\O^{2\varpi_i-\alpha_i}_{sh}(a+1))$};

\node (A5) at ( 2.4,-1.15)
  {$K_\C(\O^{2\varpi_i}_{sh}(a,a+2))$};

  \draw[arr] (A1) -- (A2);
  \draw[arr] (A1) -- (A7);

  \draw[arr] (A2) -- (A3);
  \draw[arr] (A7) -- (A6);

  \draw[arr] (A3) -- (A4);
  \draw[arr] (A6) -- (A5);

  \draw[arr] (A4) -- (A5);
\end{tikzpicture}
}
\end{center}\smallskip%
Denote by $\mathsf{i}_1$ and $\mathsf{i}_2$ the $G$-equivariant maps $V(2\varpi_i-\alpha_i)\to \smash{K_\C(\O^{2\varpi_i}_{sh}(a,a+2))}$ defined via the right and left paths (respectively) in the previous diagram. 
\begin{Theorem}\label{thm:the_heptagon_commutes}
The above heptagon commutes, that is $\mathsf{i}_1=\mathsf{i}_2$.
\end{Theorem}
\begin{proof}
Since there is only one copy of $V(2\varpi_i-\alpha_i)$ in %
$V(\varpi_i)^{\otimes 2}$, by Schur's Lemma,
$$\Hom_G(V(2\varpi_i-\alpha_i),V(\varpi_i)^{\otimes 2})\simeq \C$$ 
and $\mathsf{i}_1$ must agree with $\mathsf{i}_2$ up to a scalar. To compute this scalar, we evaluate the two maps on the highest vector $v_{2\varpi_i-\alpha_i}$. Using \eqref{eq:Gequivariant_map_wedge}--\eqref{eq:Gequivariant_map_from_identity}, we easily get
\begin{equation}\label{eq:i1_on_hw_vector}
\mathsf{i}_1(v_{2\varpi_i-\alpha_i})= [L_{\varpi_i, a} \otimes L_{s_i \varpi_i, a+2}] - [L_{s_i \varpi_i, a} \otimes L_{\varpi_i, a+2}]
\end{equation}
and (since positive prefundamental modules are 1-dimensional)
\begin{equation}\label{eq:i2_on_hw_vector}
\textstyle \mathsf{i}_2(v_{2\varpi_i - \alpha_i})= [ \bigotimes_{j\sim i} L_{\varpi_j, a+1}]=[L(\prod_{j\sim i}\sfPsi_{j,a+1})].
\end{equation}
Clearly $\mathsf{i}_2\neq 0$ as $\mathsf{i}_2(v_{2\varpi_i - \alpha_i})$ is the class of a module. Hence, the above shows that $\mathsf{i}_1=c\,\mathsf{i}_2$ for some $c\in \C$ and it follows that
\begin{equation}\label{eq:what_is_coeff_c}
\textstyle[L_{\varpi_i, a} \otimes L_{s_i \varpi_i, a+2}] - [L_{s_i \varpi_i, a} \otimes L_{\varpi_i, a+2}]=c\,[L(\prod_{j\sim i}\sfPsi_{j,a+1})].
\end{equation}
Set 
$$\textstyle\psi=\sfPsi_{i,a+2}\sfPsi_{i,a-2}^{-1}\prod_{j\sim i} \sfPsi_{j,a-1}.$$ By Corollary \ref{cor:tensor_product_decomp_in_Osh} and Example \ref{ex:Inflsl3}, the %
product $V=L_{s_i\varpi_i,a}\otimes L_{\varpi_i,a+2}$  contains exactly one copy of $L(\psi)$ and all its other composition factors have highest $\ell$-weights strictly smaller than $\psi$ (under Nakajima's partial order $\preceq$).~However, by Lemma \ref{lem:description_of_crystals_QQ_system} (and our isomorphism $\cB\simeq \mathfrak{r}$) the only other possible highest $\ell$-weight for~a simple module in $\mathcal{O}_{2\varpi_i-\alpha_i}^{2\varpi_i}(a,a+2)$ is 
$$\textstyle \prod_{j\sim i}\sfPsi_{j,a+1} = \psi \mathsf{A}_{i,a} \succ\psi.$$
Hence $V$ must be isomorphic to the simple module $L(\psi)$ and \eqref{eq:what_is_coeff_c} can be written as
$$ \textstyle [L_{\varpi_i, a} \otimes L_{s_i \varpi_i, a+2}] = c[L(\prod_{j\sim i}\sfPsi_{j,a+1})]+[L(\psi)],$$
but Corollary \ref{cor:tensor_product_decomp_in_Osh} gives instead
$$ \textstyle [L_{\varpi_i, a} \otimes L_{s_i \varpi_i, a+2}] = [L(\prod_{j\sim i}\sfPsi_{j,a+1})]+m[L(\psi)],$$
for some $m\in \Z_{\geq 0}$. The only solution is clearly $m=c=1$ and thus $\mathsf{i}_1=\mathsf{i}_2$.
\end{proof}
Theorem \ref{thm:the_heptagon_commutes} allows us to show the %
rational analogue of \cite[Conjecture 6.11]{frenkel2024extended}:
\begin{Theorem}\label{thm:extended_QQ_holds}
Fix $i\in I$ and let $w\in W$ be such that $\ell(ws_i)>\ell(w)$. Then, %
in $K_0(\O_{sh})$,
\begin{equation*}
\textstyle [L_{w\varpi_i,a}][L_{ws_i\varpi_i,a+2}]-[L_{ws_i\varpi_i,a}][L_{w\varpi_i,a+2}]=\prod_{j\sim i} [L_{w\varpi_j,a+1}].
\end{equation*}
\end{Theorem}
\begin{proof}
Applying the maps $\mathsf{i}_1,\mathsf{i}_2$ on the extremal vector $v_{w(2\varpi_i-\alpha_i)}=\dot{w}v_{2\varpi_i-\alpha_i}\in V(2\varpi_i-\alpha_i)$ gives, by $G$-equivariance and \eqref{eq:i1_on_hw_vector}--\eqref{eq:i2_on_hw_vector},
\begin{equation*}
\mathsf{i}_1(v_{w(2\varpi_i-\alpha_i)})=\dot{w}\cdot \mathsf{i}_1(v_{2\varpi_i-\alpha_i})=\dot{w}\cdot \big(\,[L_{\varpi_i, a}] [L_{s_i \varpi_i, a+2}] - [L_{s_i \varpi_i, a} ][L_{\varpi_i, a+2}]\,\big)
\end{equation*}
and 
\begin{equation*}
\mathsf{i}_2(v_{w(2\varpi_i-\alpha_i)})=\dot{w}\cdot \mathsf{i}_2(v_{2\varpi_i-\alpha_i})=\dot{w}\cdot \big(\,\textstyle\prod_{j\sim i} [L_{\varpi_j, a+1}] \,\big).
\end{equation*}
By assumption, $\ell(ws_i)>\ell(w)$ and hence $\dot{w}\dot{s}_i$ is a lift of $ws_i\in W$ to $G$. Using Theorem~\ref{thm:mult_in_KOsh_is_G_equiv} and Corollary \ref{cor:W_action_on_chamber_modules} thus allows us to rewrite the above equations as
\begin{align*}
\mathsf{i}_1(v_{w(2\varpi_i-\alpha_i)})%
&=(\dot{w}[L_{\varpi_i, a}])(\dot{w}\dot{s}_i[L_{\varpi_i, a+2}]) - (\dot{w}\dot{s}_i[L_{\varpi_i, a} ])(\dot{w}[L_{\varpi_i, a+2}])\\
&=[L_{w\varpi_i, a}][L_{ws_i\varpi_i, a+2}]-[L_{ws_i\varpi_i, a} ][L_{w\varpi_i, a+2}]
\end{align*}
and 
\begin{equation*}
\mathsf{i}_2(v_{w(2\varpi_i-\alpha_i}))=\dot{w}\cdot (\textstyle\prod_{j\sim i} [L_{\varpi_j, a+1}])=\textstyle\prod_{j\sim i} (\dot{w}[L_{\varpi_j, a+1}])=\textstyle\prod_{j\sim i} [L_{w\varpi_j, a+1}].
\end{equation*}
The desired result then follows from Theorem \ref{thm:the_heptagon_commutes}.
\end{proof}

The relation appearing in Theorem \ref{thm:extended_QQ_holds} is (the rational analogue of) Frenkel--Hernandez's \textit{extended $QQ$-system}. More precisely, set
$$\mathfrak{A}=\Z[[A_{j,b}^{-1}]]_{(j,b)\in I\times_2\Z}.$$
Then, in \cite[Sections 3--5]{frenkel2024extended} (see also \cite[Section 7]{geiss2024representations}),~the authors 
construct elements
\begin{equation}\label{eq:Qvariables}
\{Q_{w\varpi_i,a}\}_{(i,a)\in I\times_2 \Z, w\in W}\subseteq \mathcal{E}_{\ell},
\end{equation}
that admit a factorization of the form 
\begin{equation}\label{eq:factorisationQvars}
Q_{w\varpi_i,a}=\sfPsi_{w\varpi_i,a}\Sigma_{w\varpi_i,a}
\end{equation}
with $\Sigma_{w\varpi_i,a}\in \mathfrak{A}$ a series of constant term 1 satisfying $\Sigma_{\varpi_i,a}=1$. Subsequently, they show the obvious trigonometric counterpart of:
\begin{Theorem}[{\cite{frenkel2024extended,geiss2024representations}}]\label{th:GHLQQ} Fix $(i,a)\in I\times_2\Z$ with $w\in W$ such that $\ell(ws_i)>\ell(w)$. %
Then the ``extended QQ-system''
\begin{equation}\label{eq:QQQvars}
\textstyle Q_{w \varpi_i, a} Q_{w s_i \varpi_i, a+2} - Q_{w s_i \varpi_i, a} Q_{w \varpi_i, a+2} = \prod_{j \sim i} Q_{w \varpi_j, a+1}
\end{equation}
holds in $\mathcal{E}_{\ell}$.
\end{Theorem}
Motivated by \eqref{eq:factorisationQvars}, Frenkel--Hernandez and Geiss--Hernandez--Leclerc moreover proposed:
\begin{Conjecture}[{\cite[Conjecture 6.8]{frenkel2024extended}}, {\cite[Conjecture 9.19]{geiss2024representations}}]\label{conj:FHQvars} Fix $(i,a)\in I\times_2\Z$ and $w\in W$. Then \vspace*{-0.75mm}
$$\smash{ Q_{w\varpi_i,a}=\chi_{\ell}(L_{w\varpi_i,a}).}$$
\end{Conjecture}\vspace*{-1mm}
Our goal for the rest of this subsection is to show the above conjecture. For this,~we~use~a unicity result for the solutions of \eqref{eq:QQQvars} that %
is somewhat implicit in \cite{frenkel2024extended,geiss2024representations}. However, preliminary work is needed before we can prove this result. \medskip\par
Fix $(i,a)\in I\times_2\Z$ and recall from Example \ref{ex:Inflsl3} (or the proof of Theorem \ref{thm:the_heptagon_commutes}) that 
$$ \frac{\sfPsi_{s_i\varpi_i,a}\sfPsi_{i,a+2}}{\sfPsi_{i,a}\sfPsi_{s_i\varpi_i,a+2}} = \mathsf{A}_{i,a}^{-1}%
$$
and
$$ \frac{\prod_{j\sim i}\sfPsi_{j,a+1}}{\sfPsi_{i,a}\sfPsi_{s_i\varpi_i,a+2}}=1$$
Recall also from Section \ref{sec:Chambermodules} the braid group operators $\{T_w\}_{w\in W}$ of \cite[Corollary~4.5]{friesen2025braid} (which are group automorphisms of $\mathfrak{r}$) and take $w\in W$ with $\ell(ws_i)>\ell(w)$. Then, Lemma \ref{lem:highest_ell_weight_chamber_module_braid} and the above computations give
\begin{align*}
T_{w^{-1}}^{-1}(\mathsf{A}_{i,a}^{-1})&=\frac{T_{w^{-1}}^{-1}\big(\,T_{s_i}^{-1}(\sfPsi_{i,a})\sfPsi_{i,a+2}\,\big)}{T_{w^{-1}}^{-1}\big(\,\sfPsi_{i,a}T_{s_i}^{-1}(\sfPsi_{i,a+2})\,\big)}\\&=\frac{(T_{s_i}T_{w^{-1}})^{-1}(\sfPsi_{i,a})T_{w^{-1}}^{-1}(\sfPsi_{i,a+2})}{T_{w^{-1}}^{-1}(\sfPsi_{i,a})(T_{s_i}T_{w^{-1}})^{-1}(\sfPsi_{i,a+2})}=\frac{\sfPsi_{ws_i\varpi_i,a}\sfPsi_{w\varpi_i,a+2}}{\sfPsi_{w\varpi_i,a}\sfPsi_{ws_i\varpi_i,a+2}}
\end{align*}
and, similarly,
$$ 1=\frac{\prod_{j\sim i}T_{w^{-1}}^{-1}(\sfPsi_{j,a+1})}{T_{w^{-1}}^{-1}(\sfPsi_{i,a})(T_{s_i}T_{w^{-1}})^{-1}(\sfPsi_{i,a+2})}=\frac{\prod_{j\sim i}\sfPsi_{w\varpi_j,a+1}}{\sfPsi_{w\varpi_i,a}\sfPsi_{ws_i\varpi_i,a+2}}.$$
Using \cite[Corollary 4.1]{friesen2025braid}, we easily get
$$ T_{s_k}(\mathsf{A}_{j,b})=\left\{
\begin{array}{ll}
\mathsf{A}_{j,b+2}^{-1} & \text{if }k=j,\\
\mathsf{A}_{j,b}\mathsf{A}_{k,b+1} & \text{if }k\sim j,\\
\mathsf{A}_{j,b} & \text{else}.
\end{array}
\right. $$
Thus $T_{w^{-1}}^{-1}(\mathsf{A}_{i,a}^{-1})\in \mathcal{A}$ since
\begin{equation}\label{eq:BraidAinv}
T_{s_k}^{-1}(\mathsf{A}_{j,b})=\left\{
\begin{array}{ll}
\mathsf{A}_{j,b-2}^{-1} & \text{if }k=j,\\
\mathsf{A}_{j,b}\mathsf{A}_{k,b-1} & \text{if }k\sim j,\\
\mathsf{A}_{j,b} & \text{else}.
\end{array}
\right. 
\end{equation}
Combining the above, we see that dividing \eqref{eq:QQQvars} by $\sfPsi_{w\varpi_i,a}\sfPsi_{ws_i\varpi_i,a+2}\in \mathfrak{r}$ gives the following ``renormalized extended $QQ$-system'' 
\begin{equation}\label{eq:renormalizedQQ}
\textstyle \Sigma_{w\varpi_i,a}\Sigma_{ws_i\varpi_i,a+2}-\varrho_{i,a,w}\Sigma_{ws_i\varpi_i,a}\Sigma_{w\varpi_i,a+2}=\prod_{j\sim i}\Sigma_{w\varpi_j,a+1}
\end{equation}
where $\varrho_{i,a,w}:=T_{w^{-1}}^{-1}(\mathsf{A}_{i,a}^{-1})\in \mathcal{A}$.
In particular, since the $\Sigma$'s belong to $\mathfrak{A}=\Z[[\mathsf{A}_{j,b}^{-1}]]_{(j,b)\in I\times_2\Z}$ and have constant term 1, \vspace*{-1mm}
$$\varrho_{i,a,w}\in \mathcal{A}_+^{-1}\subseteq \mathcal{A}$$
(where $\mathcal{A}_+$ is the set of monomials in the $\mathsf{A}$'s). Thus \eqref{eq:renormalizedQQ} is actually a relation in $\mathfrak{A}$. \medskip\par 
On the other hand, the normalized $\ell$-characters of chamber modules satisfy (by Theorem \ref{thm:extended_QQ_holds}) a similar relation, i.e.
\begin{equation}\label{eq:renormalizedQQ2}
\textstyle \tilde{\chi}_{\ell}(L_{w\varpi_i,a})\tilde{\chi}_{\ell}(L_{ws_i\varpi_i,a+2})-\varrho_{i,a,w}\tilde{\chi}_{\ell}(L_{ws_i\varpi_i,a})\tilde{\chi}_{\ell}(L_{w\varpi_i,a+2})=\prod_{j\sim i}\tilde{\chi}_{\ell}(L_{w\varpi_j,a+1}),
\end{equation}
again assuming $\ell(ws_i)>\ell(w)$. The following is adapted from \cite[Lemma 2.5]{frenkel2022weyl}. We will use the fact that $\Sigma_{\varpi_i,a}$ is invertible in $\mathfrak{A}$ (since it has constant term 1). 
\begin{Theorem}\label{thm:ConjFH} For $(i,a)\in I\times_2\Z$ and $w\in W$, 
\begin{equation}\label{eq:SigmaNormellChar}
\tilde{\chi}_{\ell}(L_{w\varpi_i,a})=\Sigma_{w\varpi_i,a}
\end{equation}
in $\mathfrak{A}$. In particular, Conjecture \ref{conj:FHQvars} holds.
\end{Theorem}
\begin{proof}
We show \eqref{eq:SigmaNormellChar} by induction on $\ell(w)$. Note that the case $\ell(w)=0$ is trivial since~both $\tilde{\chi}_{\ell}(L_{\varpi_i,a})$ and $\Sigma_{\varpi_i,a}$ are equal to 1. Suppose hence $\ell(w)\geq 1$ and choose $(i,a)\in I\times_2\Z$.~We consider the following two cases:
\begin{enumerate}
\item Assume $\ell(ws_i)>\ell(w)$. Then $w\varpi_i=w'\varpi_i$ for some $w'\in W$ with $\ell(w')<\ell(w)$ and
$$ \tilde{\chi}_{\ell}(L_{w\varpi_i,a})=\tilde{\chi}_{\ell}(L_{w'\varpi_i,a})=\Sigma_{w'\varpi_i,a}=\Sigma_{w\varpi_i,a}$$
by the induction hypothesis.
\item Assume $\ell(ws_i)<\ell(w)$. Then, $w'=ws_i$ satisfies $\ell(w')<\ell(w's_i)=\ell(w)$. Thus,~using the induction hypothesis and \eqref{eq:renormalizedQQ}--\eqref{eq:renormalizedQQ2} easily gives the equality
\begin{equation*}
\Sigma_{w'\varpi_i,a}(\Sigma_{w\varpi_i,a+2}-\tilde{\chi}_{\ell}(L_{w\varpi_i,a+2}))=\varrho_{i,a,w'}(\Sigma_{w\varpi_i,a}-\tilde{\chi}_{\ell}(L_{w\varpi_i,a}))\Sigma_{w'\varpi_i,a+2}.
\end{equation*}
Let $\tau:\mathfrak{A}\to\mathfrak{A}$ be the ring automorphism given by $\mathsf{A}_{j,b}^{-1}\mapsto \mathsf{A}_{j,b+2}^{-1}$ for $(j,b)\in I\times_2\Z$. Then 
the definition of the $\Sigma$'s in \cite{frenkel2024extended} gives directly
$$ \tau(\Sigma_{w\varpi_i,a})=\Sigma_{w\varpi_i,a+2}$$
and the above equality can be rewritten as
\begin{center}
$\tau(\chi)=\varrho_{i,a,w'}\chi$
\end{center}
where 
$$\chi=\Sigma_{w'\varpi_i,a}^{-1}(\Sigma_{w\varpi_i,a}-\tilde{\chi}_{\ell}(L_{w\varpi_i,a}))\in \mathfrak{A}.$$ 
We claim that $\chi=0$. Indeed, suppose otherwise and consider a monomial $M\in \mathcal{A}_+^{-1}$ of $\tau(\chi)$ with maximal weight $\wt(M)$ (for $\wt:\mathfrak{r}\to \mathfrak{h}$ the map given in Section~\ref{subsec:ell_characters_and_blocks}~and
where the order on $\mathfrak{h}$ is $\omega \leq \omega' \iff \omega'-\omega\in Q_+^{\vee}$). As $\tau(\chi)=\varrho_{i,a,w'}\chi$, there must exist a monomial $M'$ of $\chi$ for which $$\wt(M)=\wt(M')+\wt(\rho_{i,a,w'}).$$ 
However, $\wt(M') =\wt(\tau(M'))$ is also the weight of a monomial of $\tau(\chi)$ and 
$$ 0\neq\wt(\rho_{i,a,w'})\in \wt(\mathcal{A}_+^{-1}) \subseteq  -Q_+^{\vee}$$
contradicts the maximality of $\wt(M)$%
.~Thus, $\chi =0$ and \eqref{eq:SigmaNormellChar} follows.
\qedhere
\end{enumerate}
\end{proof}
\begin{Rem} To show that their $Q$-variables solve the extended $QQ$-system \eqref{eq:QQQvars}, Frenkel--Hernandez transfer in \cite{frenkel2024extended} known results about the usual $QQ$-system (i.e.~the case $w=e$) using an action of the $2$-covering $\smash{\widetilde{W}}$ of $W$ on a completion of the $\ell$-character ring $\mathcal{E}_{\ell}$ (see~also Remark \ref{rem:WactionEell}). Our method for proving Theorem \ref{thm:extended_QQ_holds} is very similar to this, but uses instead the $\smash{\widetilde{W}}$-action on $K_0(\mathcal{O}_{sh})$ of Theorem \ref{thm:construction_of_gaction_on_KOsh}. It would be interesting to relate the two actions.
\end{Rem}
\begin{Rem} We recall that the results shown here for representations of shifted Yangians imply the analogous results for representations of shifted quantum affine algebras~(by~\cite[Corollary 1.2.1]{varagnolo2025representations}). Also, one can use the above results in the setting of Borel quantum affine algebras using recent results \cite[Section~4]{hernandez2026borel}.
\end{Rem}
As explained in Section \ref{sec:Intro}, the $Q$-variables \eqref{eq:Qvariables} and the extended $QQ$-system \eqref{eq:QQQvars} respectively correspond to initial cluster variables and mutation relations for the cluster algebra studied in \cite{geiss2024representations}. Thus Theorem \ref{thm:extended_QQ_holds} and Theorem \ref{thm:ConjFH} can be seen as first steps toward proving \cite[Conjecture 9.16]{geiss2024representations} (which, in essence, states that the category $\mathcal{O}_{sh}$ gives a monoidal categorification of the above cluster algebra). %
\subsection{Lie-theoretic interpretation of the extended $QQ$-system
}\label{subsec:extendQQ_and_principal_block} By Lemma \ref{lemma:EqualityTrunc}, given $(i,a)\in I\times_2\Z$, the root embedding associated to $\alpha_i$ identifies the truncated shifted Yangian $\smash{Y^{2\varpi_i}_{2\varpi_i-\alpha_i}(\{a,a+2\})}$ (of type $\fg^{\vee}$) with the truncation 
$$\smash{Y_0^{2\varpi_i}(\{a,a+2\},\mathfrak{g}^{\vee}_{\{i\}})},$$
where the underlying Lie-algebra is $\fg_{\{i\}}^{\vee}\simeq \mathfrak{sl}_2$ (see Appendix \ref{sec:Inflations} for details on our notation). Also, by \cite[Example~4.6]{kamnitzer2014yangians}, the inclusion $U(\mathfrak{sl}_2)\subseteq Y_0(\mathfrak{sl}_2)$ induces an isomorphism
$$ \smash{Y_0^{2\varpi_i}(\{a,a+2\},\mathfrak{g}^{\vee}_{\{i\}})}\simeq U(\mathfrak{sl}_2)/\langle C\rangle$$
(where $C\in U(\mathfrak{sl}_2)$ is the usual Casimir element), and Corollary \ref{cor:identification_of_cat_O} together with a result~of Soergel \cite[Théorème 1]{soergel1986equivalence} imply that there is an equivalence of categories 
\begin{equation}\label{eq:EquivCatsl2we}
\smash{\mathcal{O}_{2\varpi_i-\alpha_i}^{2\varpi_i}(\{a,a+2\})}\simeq \mathcal{O}_{\chi_0}(\mathfrak{sl}_2)
\end{equation}
where $\mathcal{O}_{\chi_0}(\mathfrak{sl}_2)$ is the principal block of the BGG category $\mathcal{O}$ of $\mathfrak{sl}_2$. \medskip\par
We will need the following elementary %
result%
:
\begin{Lemma}\label{lem:weightsVvarpii}
Fix $w\in W$ and $j\in I$. Assume $\langle w(2\varpi_i-\alpha_i),\alpha_j^{\vee}\rangle > 0$. Then $w(2\varpi_i-\alpha_i)+\alpha_j$ is not a weight of the $G$-module $V(\varpi_i)^{\otimes 2}$.
\end{Lemma}
\begin{proof}
Let $\beta=w^{-1}(\alpha_j)\in \Delta$. By hypothesis, 
$$\textstyle \langle\sum_{k\sim i}\varpi_k,\beta^{\vee}\rangle=\langle 2\varpi_i-\alpha_i,\beta^{\vee}\rangle=\langle w(2\varpi_i-\alpha_i),\alpha_j^{\vee} \rangle  >0,$$
and $\beta$ must satisfy $\beta\geq \alpha_k$ for some $k\sim i$. In particular, $\omega =2\varpi_i-\alpha_i+\beta$ is not a weight~of $V(\varpi_i)^{\otimes 2}$ as this would give $\omega\leq 2\varpi_i$, that is, equivalently, $\beta\leq \alpha_i$. This ends the proof.
\end{proof}
Using Lemma \ref{lem:weightsVvarpii} and the results of Section \ref{subsec:chamber_modules_and_G_action}, we extend the above discussion to:
\begin{Theorem}\label{thm:QQ_equiv_bloc_principal}
For all $w\in W$%
, there is an equivalence of categories 
$$\mathcal{O}_{w(2\varpi_i-\alpha_i)}^{2\varpi_i}(\{a,a+2\})\simeq \mathcal{O}_{\chi_0}(\mathfrak{sl}_2).$$
\end{Theorem}
\begin{proof}
Let $\la=2\varpi_i-\alpha_i\in P_+$. We proceed by induction on $r_w=\rht(\la-w(\la))\geq 0$ with~the case $r_w=0$ (i.e.~$w=e$) shown above. Suppose hence $r_w>0$ (so $w(\la)%
\not\in P_+$) and fix $j \in I$ for which
$ m=\langle w(\la),\alpha_j^{\vee}\rangle <0$. Then $w'=s_jw$ satisfies %
$r_{w'}=r_w+m<r_w$ %
and
$$\smash{\langle w'(\la),\alpha_j^{\vee}\rangle} = -m >0.$$
In particular, the induction hypothesis with Lemma \ref{lem:ChuangRouquier} and Lemma \ref{lem:weightsVvarpii} gives equivalences
$\smash{\mathcal{O}_{\chi_0}(\mathfrak{sl}_2)}\simeq \smash{\mathcal{O}_{w'(2\varpi_i-\alpha_i)}^{2\varpi_i}(\{a,a+2\})}\simeq \smash{\mathcal{O}_{w(2\varpi_i-\alpha_i)}^{2\varpi_i}(\{a,a+2\})},$
and we can end the proof.
\end{proof}
We use Theorem \ref{thm:QQ_equiv_bloc_principal} to give a Lie-theoretic interpretation of the extended $QQ$-system. For this, recall from the proof of Theorem \ref{thm:the_heptagon_commutes} that the tensor products
\begin{center}
$\bigotimes_{j\sim i} L_{\varpi_j,a+1}\text{ and }\textstyle L_{s_i\varpi_i,a}\otimes L_{\varpi_i,a+2}$
\end{center}
are simple in the category $\mathcal{O}_{2\varpi_i-\alpha_i}^{2\varpi_i}(\{a,a+2\})$ where the $QQ$-system for $w=e$ takes place. In addition, checking dimensions, we see that the equivalence of categories \eqref{eq:EquivCatsl2we} sends these products to the simple objects of respective highest weight $0$ and $-\alpha_1$ in $\mathcal{O}_{\chi_0}(\mathfrak{sl}_2)$.~Thus,~for $w\in W$ such that $\ell(ws_i)>\ell(w)$, since
$$ \textstyle [L_{ws_i\varpi_i,a}][L_{w\varpi_i,a+2}] = \dot{w}[L_{s_i\varpi_i,a}][L_{\varpi_i,a+2}] \text{ and }\prod_{j\sim i}[L_{w\varpi_j,a+1}] = \dot{w}[L_{\varpi_j,a+1}]$$
(by Theorem \ref{thm:mult_in_KOsh_is_G_equiv} and Corollary \ref{cor:W_action_on_chamber_modules}), and since the isomorphism 
$$ \smash{K_{\C}(\mathcal{O}_{2\varpi_i-\alpha_i}^{2\varpi_i}(\{a,a+2\}))}\simeq \smash{K_{\C}(\mathcal{O}_{w(2\varpi_i-\alpha_i)}^{2\varpi_i}(\{a,a+2\}))}$$
(coming from the proof of Theorem \ref{thm:QQ_equiv_bloc_principal}) is given by the action of $\dot{w}$, we directly get:
\begin{Proposition} Fix $w\in W$ with $\ell(ws_i)>\ell(w)$. Then the equivalence of Theorem \ref{thm:QQ_equiv_bloc_principal} identifies the products
$\bigotimes_{j\sim i}L_{w\varpi_j,a+1}$ and $\textstyle L_{ws_i\varpi_i,a}\otimes L_{w\varpi_i,a+2}$ with the simple modules of highest weight $0$ and $-\alpha_1$ in $\mathcal{O}_{\chi_0}(\mathfrak{sl}_2)$ (resp.). In particular, these products are simple.
\end{Proposition}
\begin{Corollary}\label{cor:extQQissl2decomp} Fix $w\in W$ as above. Then the equivalence of Theorem \ref{thm:QQ_equiv_bloc_principal} identifies~the extended $QQ$-system of Theorem \ref{thm:extended_QQ_holds} with the relation 
$$ [\Delta(0)]=1+[\Delta(-\alpha_1)]$$
giving the decomposition of the class of the Verma module $\Delta(0)$ with highest weight $0$ as the sum of the class of the trivial $\mathfrak{sl}_2$-module and the class of the simple Verma module~$\Delta(-\alpha_1)$ of highest weight $-\alpha_1$.
\end{Corollary}
By the above results, the extended $QQ$-system for $w$ can be written as the relation
\begin{equation}\label{eq:extQQsimples}
\textstyle [V_w]=[L(\sfPsi_{ws_i\varpi_i,a}\sfPsi_{w\varpi_i,a+2})]+[L(\prod_{j\sim i}\sfPsi_{w\varpi_j,a+1})],
\end{equation}
which states that $V_w:=L_{w\varpi_i,a}\otimes L_{ws_i\varpi_i,a+2}$ can be constructed out of the two composition factors $L(\sfPsi_{ws_i\varpi_i,a}\sfPsi_{w\varpi_i,a+2})\simeq L_{ws_i\varpi_i,a}\otimes L_{w\varpi_i,a+2}$ and $L(\prod_{j\sim i}\sfPsi_{w\varpi_j,a+1})\simeq \bigotimes_{j\sim i}L_{w\varpi_j,a+1}$ (with multiplicity 1). Also, for $w=e$, Theorem \ref{thm:HZhighest} (i) shows that $V_e$ is of highest $\ell$-weight and \eqref{eq:extQQsimples} thus comes from the non-split short exact sequence
\begin{equation}\label{eq:SECQQcaswe}
\textstyle 0\rightarrow L(\sfPsi_{s_i\varpi_i,a}\sfPsi_{\varpi_i,a+2})\to V_e \to L(\prod_{j\sim i}\sfPsi_{\varpi_j,a+1})\to 0
\end{equation}
\begin{Conjecture}\label{con:sec_qq_systems}
Choose %
$w\in W$ such that $\ell(ws_i)>\ell(w)$. Then $V_w=L_{w\varpi_i,a}\otimes L_{ws_i\varpi_i,a+2}$ is of highest $\ell$-weight and there is a non-split short exact sequence
\begin{equation}\label{eq:SECQQ}
\textstyle 0 \to L_{ws_i\varpi_i,a}\otimes L_{w\varpi_i,a+2} \to  V_w \to \bigotimes_{j\sim i} L_{w\varpi_j,a+1}\to  0
\end{equation}
that categorifies the extended QQ-system.
\end{Conjecture}
\begin{Rem} Suppose Conjecture \ref{con:sec_qq_systems} holds for $w$. Then the reasoning of the beginning of Section \ref{sec:TensO} shows that the unique map (up to scaling)
\begin{center}
$ R_w:L_{w\varpi_i,a}\otimes L_{ws_i\varpi_i,a+2} \to L_{ws_i\varpi_i,a+2}\otimes L_{w\varpi_i,a},$
\end{center}
satisfies 
\begin{center}
$\op{Im} R_w\simeq \bigotimes_{j\sim i}L_{w\varpi_j,a+1}$ and $\Ker R_w\simeq L_{ws_i\varpi_i,a}\otimes L_{w\varpi_i,a+2}$.
\end{center}
The sequence \eqref{eq:SECQQ} thus ``comes from an $R$-matrix'', in the sense that it is the first row of
$$ 
\adjustbox{scale=0.85}{
\begin{tikzcd}[row sep = 1.65em]
0\arrow[r] & \Ker R\arrow[-,double line with arrow={-,-}]{d}\arrow[r] & L_{w\varpi_i,a}\otimes L_{ws_i\varpi_i,a+2}\arrow[r]\arrow[d,"R", shift left] & \op{Im} R \arrow[-,double line with arrow={-,-}]{d}\arrow[r]& 0\\
0 & \op{Im} \mathbf{tr}(R) \arrow[l] & L_{ws_i\varpi_i,a+2}\otimes L_{w\varpi_i,a}\arrow[u,"\mathbf{tr}(R)",shift left]\arrow[l] & \Ker \mathbf{tr}(R) \arrow[l] & \arrow[l] 0
\end{tikzcd}}
$$
with $\mathbf{tr}$ the duality of Section \ref{sec:Auto}. We expect this to be useful in proving \cite[Conjecture 9.16]{geiss2024representations}, exactly as how %
$R$-matrices helped prove important results in the work of  Kang--Kashiwara--Kim--Oh--Park (see, e.g., \cite{kang2018monoidal,kashiwara2024monoidal}) and Cautis--Williams \cite{cautis2019cluster} on monoidal categorifications %
arising from KLR algebras or geometry.
\end{Rem}
We will show the above conjecture for shifted Yangians of type $A$ %
in an upcoming paper \cite{otherpaper} using Corollary \ref{cor:HZgenSimple} and $\ell$-character formulae for chamber modules (see also~Theorem \ref{thm:HZcriterionTensSimpPol} and Proposition \ref{prop:CombCritNeg}). In the meantime, we illustrate why it holds in an example. 
\begin{Example} Fix $\fg=\mathfrak{sl}_4$, $i=2$ and $a=3$. Then the tensor products to consider are:\vspace*{-1.3mm}
\begin{multicols}{3}
\begin{enumerate}[label=(\roman*)]
\item\label{eq:one} $L(\sfPsi_{2,3})\otimes L(\frac{\sfPsi_{1,4}\sfPsi_{3,4}}{\sfPsi_{2,3}})$
\item $L(\sfPsi_{2,3})\otimes L(\frac{\sfPsi_{1,4}}{\sfPsi_{3,2}})$
\item $L(\sfPsi_{2,3})\otimes L(\frac{\sfPsi_{3,4}}{\sfPsi_{1,2}})$
\item\label{eq:four} $L(\sfPsi_{2,3})\otimes L(\frac{\sfPsi_{2,3}}{\sfPsi_{1,2}\sfPsi_{3,2}})$
\item\label{eq:five} $L(\frac{\sfPsi_{1,2}\sfPsi_{3,2}}{\sfPsi_{2,1}})\otimes L(\frac{\sfPsi_{1,4}}{\sfPsi_{3,2}})$
\item\label{eq:six} $L(\frac{\sfPsi_{1,2}\sfPsi_{3,2}}{\sfPsi_{2,1}})\otimes L(\frac{\sfPsi_{3,4}}{\sfPsi_{1,2}})$
\item $L(\frac{\sfPsi_{1,2}}{\sfPsi_{3,0}})\otimes L(\frac{\sfPsi_{2,3}}{\sfPsi_{1,2}\sfPsi_{3,2}})$
\item\label{eq:eight} $L(\frac{\sfPsi_{3,2}}{\sfPsi_{1,0}})\otimes L(\frac{\sfPsi_{2,3}}{\sfPsi_{1,2}\sfPsi_{3,2}})$
\item\label{eq:nine} $L(\frac{\sfPsi_{1,2}\sfPsi_{3,2}}{\sfPsi_{2,1}})\otimes L(\sfPsi_{2,1}^{-1})$
\item $L(\frac{\sfPsi_{1,2}}{\sfPsi_{3,0}})\otimes  L(\sfPsi_{2,1}^{-1})$
\item $L(\frac{\sfPsi_{3,2}}{\sfPsi_{1,0}})\otimes L(\sfPsi_{2,1}^{-1})$
\item\label{eq:twelve} $L(\frac{\sfPsi_{2,1}}{\sfPsi_{1,0}\sfPsi_{3,0}})\otimes  L(\sfPsi_{2,1}^{-1})$
\end{enumerate}
\end{multicols}\vspace*{-2.75mm}
Clearly, %
\ref{eq:one}--\ref{eq:four} and \ref{eq:nine}--\ref{eq:twelve} are of highest $\ell$-weight by Theorem~\ref{thm:HZhighest}. To prove that~\ref{eq:five} is of highest $\ell$-weight, %
we 
show that $\smash{L(\sfPsi_{2,1}^{-1})\otimes L(\frac{\sfPsi_{1,4}}{\sfPsi_{3,2}})}$ is simple (cf.~Corollary~\ref{cor:HZgenSimple}). However, 
$$ \psi:=\smash{\tfrac{\sfPsi_{1,4}}{\sfPsi_{2,1}\sfPsi_{3,2}}}\in \cB(2\varpi_2,\{5,5\})_{\varpi_1-\varpi_2-\varpi_3}$$
and $\cB(2\varpi_2,\{5,5\})\simeq B(2\varpi_2)$ by Theorem \ref{thm:charaterization_max_sing_crystals}%
. Therefore, the category $\smash{\mathcal{O}_{\varpi_1-\varpi_2-\varpi_3}^{2\varpi_2}}(\{5,5\})$
contains only $\dim V(2\varpi_2)_{\varpi_1-\varpi_2-\varpi_3}=1$ simple object by Theorem \ref{th:descend}, and the fact that the $\ell$-weight space $(L(\sfPsi_{2,1}^{-1})\otimes L(\tfrac{\sfPsi_{1,4}}{\sfPsi_{3,2}}))_{\psi}$ has dimension 1 implies 
$$ \smash{L(\sfPsi_{2,1}^{-1})\otimes L(\tfrac{\sfPsi_{1,4}}{\sfPsi_{3,2}})}\simeq L(\psi),$$
as desired. The proof that \ref{eq:six}--\ref{eq:eight} are of highest $\ell$-weight is totally analogous.
\end{Example}
Finally, we end this subsection by observing that the analysis done for the extended $QQ$-system of Theorem \ref{thm:extended_QQ_holds} can be carried out equally well for the relations of Corollary \ref{cor:RelationBaseAff} (that come from the inclusions $\C[N_-\bbslash G]\subseteq K_{\C}(\mathcal{O}_{sh})$). In this case, the category~$\mathcal{O}_{\chi_0}(\mathfrak{sl}_2)$ is replaced by a singular block of a parabolic category $\mathcal{O}$ for $\fg=\mathfrak{sl}_3$ and Corollary \ref{cor:RelationBaseAff} can be seen as giving the ``semisimple decomposition'' of the (singular) parabolic Verma~module of highest weight $-\varpi_1$ into its (two) simple composition factors. %
\subsection{Chamber modules are real prime simple modules}\label{sec:application_real_and_prime} Recall that a simple object~$V$ of $\mathcal{O}_{sh}$ is said to be \textit{real} if $V^{\otimes 2}$ is simple, and \textit{prime} if, for every factorization $V\simeq V_1\otimes V_2$, either $V_1$ or $V_2$ is the trivial representation $L(\mathbbm{1})$. In this section, we show that all chamber modules are real and prime%
.\medskip\par 

Fix $\lambda\in P_+$, %
$\bR\in \Z^\lambda$ and $w\in W$. Since the $w\la$-weight space of $\textstyle \smash{\cB(\lambda,\bR)\subseteq \bigotimes_{i\in I} B(\varpi_i)^{\otimes \lambda_i}}$ contains one element, Theorem \ref{th:descend} shows that $\smash{\mathcal{O}_{w\la}^{\la}(\bR)}$ has a unique simple object $L_{w\la,\bR}$ (up to isomorphism). We call such a module $L_{w\la,\bR}$ \textit{extremal}. Clearly, chamber modules~are extremal modules associated to fundamental weights.
\begin{Proposition}\label{prop:extremalmodulesReal}
Fix $\lambda\in P_+$ and $w\in W$. Fix also sets of parameters $\bR_1,\dots, \bR_n\in \Z^\lambda$. Then the product $L_{w\lambda,\bR_1}\otimes \dots \otimes L_{w\lambda,\bR_n}$ is simple. In particular, extremal modules are real.
\end{Proposition}
\begin{proof}
Fix $\bR =\bR_1\cup\dots \cup\bR_n$ and note that $V=L_{w\lambda,\bR_1}\otimes \dots \otimes L_{w\lambda,\bR_n}$ belongs to $\smash{\O_{w(n\lambda)}^{n\lambda}(\bR)}$ by Theorem \ref{thm:truncoprod}. Moreover, by the above paragraph,%
$$ [V] = k[\smash{L_{w(n\la),\bR}}] $$
for some $k\in \Z_{\geq 1}$, and clearly $k=1$ because of Corollary \ref{cor:tensor_product_decomp_in_Osh}.
\end{proof}
\begin{Rem}
The above generalizes the well-known fact that tensor products of negative (or positive) prefundamental modules are always simple (see, e.g., \cite[Corollary 4.9]{hernandez2024shifted}). Also, for $\la \in P_+$ and $\bR\in \Z^{\la}$, the above result implies that the module $L_{w\lambda,\bR}$ is isomorphic to a tensor product of chamber modules (given by decomposing $\bR$ in 1-element multisets).
\end{Rem}
Now, fix $(i,a)\in I\times_2\Z$. The following is the main result of this subsection:
\begin{Theorem}\label{thm:chamber_modules_are_prime}
For all $w\in W$, the chamber module $L_{w\varpi_i,a}$ is real and prime.
\end{Theorem}
\begin{proof}
Using Proposition \ref{prop:extremalmodulesReal}, we see that it is enough to show that $L_{w\varpi_i,a}$ is prime, which we do by induction on $r_w=\rht(\varpi_i-w\varpi_i)$. Hence, we first assume that $L_{\varpi_i,a}\simeq V_1\otimes V_2$~for $V_1,V_2$ in $\mathcal{O}_{sh}$ and note that, in this case, $\dim V_1=\dim V_2=1$. %
Thus, by \cite[Lemma~4.1]{hernandez2024shifted},
\begin{equation*}
\smash{V_1\simeq L(\mathbf{P}_1) \ \text{ and }\ V_2\simeq L(\mathbf{P}_2)}
\end{equation*}
for some $\mathbf{P}_1,\mathbf{P}_2\in \mathfrak{r}$ polynomial, and $\sfPsi_{i,a}=\mathbf{P}_1\mathbf{P}_2$ since
\begin{equation*}
\smash{L(\mathsf{\Psi}_{i,a})=L_{w\varpi_i,a}\simeq V_1\otimes V_2 \simeq L(\mathbf{P}_1\mathbf{P}_2)}.
\end{equation*}
In particular, either $\mathbf{P}_1=\mathbbm{1}$ or $\mathbf{P}_2=\mathbbm{1}$, which ends the proof for the case $r_w=0$ (i.e.~$w=e$).
Now, assume $r_w>0$ and fix $w'=s_jw$ where $j\in I$ is chosen such that $m=-\langle w\varpi_i,\alpha_j^{\vee}\rangle>0$. Then a similar proof as the one used for Lemma \ref{lem:weightsVvarpii} shows that $w'\varpi_i+\alpha_j$ is not a weight of $V(\varpi_i)$. Hence Lemma \ref{lem:ChuangRouquier} implies that the divided power functor
$$ \smash{\mathcal{E}_j^{(m)}}:\mathcal{O}^{\varpi_i}_{w\varpi_i}(a)\to \mathcal{O}_{w'\varpi_i}^{\varpi_i}(a)$$
is an equivalence of categories, quasi-inverse to $\smash{\mathcal{F}_j^{(m)}}$. Moreover, by Corollary \ref{cor:W_action_on_chamber_modules}, 
$$ \smash{L_{w'\varpi_i,a}\simeq \smash{\mathcal{E}_j^{(m)}\mathcal{F}_j^{(m)}}(L_{w'\varpi_i,a})\simeq  \smash{\mathcal{E}_j^{(m)}}(L_{w\varpi_i,a})}.$$
Suppose as above $L_{w\varpi_i,a}\simeq V_1\otimes V_2$ for $V_1,V_2$ in $\mathcal{O}_{sh}$. Then we can assume that $V_1$ and $V_2$ are simple, so that Corollary \ref{coro:every_object_descends_to_some_trunc} implies that they respectively lie in
$$ \smash{\mathcal{O}_{sh}^{\la_1}(\bR_1)\ \text{ and }\ \mathcal{O}_{sh}^{\la_2}(\bR_2)}$$
for some  $\lambda_1,\lambda_2\in P_+$ with $\bR_1\in \Z^{\lambda_1}$ and $\bR_2\in \Z^{\lambda_2}$. Using Theorem \ref{thm:mult_in_KOsh_is_G_equiv}, we thus easily get
\begin{align*}
[L_{w'\varpi_i,a}]%
=[\smash{\mathcal{E}_j^{(m)}}(L_{w\varpi_i,a})]%
=e_j^{(m)}[V_1\otimes V_2]%
=\textstyle \sum_{n=0}^m[\smash{\mathcal{E}_j^{(n)}}(V_1)\otimes \smash{\mathcal{E}_j^{(m-n)}}(V_2)],
\end{align*}
and the simplicity of $L_{w'\varpi_i,a}$ allows us to reduce the sum of classes above to a~single term. Consequently, there exists $0\leq n\leq m$ such that 
$$L_{w'\varpi_i,a}\simeq \smash{\mathcal{E}_j^{(n)}}(V_1)\otimes \smash{\mathcal{E}_j^{(m-n)}}(V_2),$$
but $L_{w'\varpi_i,a}$ is prime by the induction hypothesis (as $r_{w'}=r_w-m<r_w$). Hence either
$$\smash{\mathcal{E}_j^{(n)}}(V_1)\simeq L(\mathbbm{1})\ \text{ or }\ \smash{\mathcal{E}_j^{(m-n)}}(V_2)\simeq L(\mathbbm{1}).$$
This clearly implies $V_1\simeq L(\mathbbm{1})$ or $V_2\simeq L(\mathbbm{1})$ (since the only simple $G$-submodule of $K_0(\mathcal{O}_{sh})$ that contains $[L(\mathbbm{1})]$ is the trivial representation).
\end{proof}

\begin{Rem}\label{rem:realprimecluster} The initial cluster variables of the cluster algebra of \cite{geiss2024representations} all correspond to classes of chamber modules via the identification of Theorem \ref{thm:ConjFH}. The fact that chamber modules are real hence proves that \cite[Conjecture 9.16]{geiss2024representations} holds (at least) for the initial cluster variables whereas the fact that they are prime is instead compatible with the original definition of monoidal categorifications of cluster algebras given in\footnote{This definition has now been relaxed a little, see for example \cite[Definition 4.3]{hernandez2025symmetries}.} \cite[Definition 2.1]{hernandez2010cluster}.
\end{Rem}

\section{Background on Bott--Samelson varieties} \label{sec:bott-samselson-background}

Throughout Sections \ref{sec:bott-samselson-background} and \ref{sec:bi-infinite-bott-samelson}, all schemes and pro-varieties are over $\C$. Our convention, as in Section  \ref{subsec:base_affine_space_and_Osh}, is that the geometric objects we study are equipped with \emph{right} $G$-actions, and their rings of functions (or more generally, spaces of sections of equivariant line bundles) are equipped with the induced \emph{left} $G$-action.
 Hence we will adopt the convention of considering $V(\lambda)^*$ as a \emph{right} $G$-module. Moreover, $v_{w\lambda}^*$ will denote the vector of weight $w\lambda$ in the right module $V(\lambda)^*$, normalized so that its pairing with $v_{w\lambda}$ is $1$. 

\subsection{Generalized flag varieties}

Given a subset 
$J\subseteq I$, let $\Fl_J$ be the $G$-orbit of $([v_{\varpi_i}^*])_{i\in J}$ in $\prod_{i\in J} \PP(V(\varpi_i)^*)$. The stabilizer of $([v_{\varpi_i}^*])_{i\in J}$ is the parabolic subgroup $P_{\bar{J}}^-$ \cite[Chapter 10.6.2]{procesi2007liegroups} where $\bar{J}=I\setminus J$, so $\op{Fl}_J$ is isomorphic to the generalized partial flag variety $P_{\bar{J}}^-\backslash G$. Moreover, since $P_{\bar{J}}^-\backslash G$ is complete, $\op{Fl}_J$ is a closed subscheme of $\prod_{i\in J} \PP(V(\varpi_i)^*)$. This is typically called the \emph{Plücker embedding} of $P_{\bar{J}}^-\backslash G$. We will denote the full flag variety $\Fl_I$ by $\op{Fl}$ and we will denote the generalized partial flag variety $\Fl_{\{i\}}$ corresponding to a maximal parabolic by $\op{Gr(i)}$ (as these are the usual Grassmannians in type $A$). \medskip\par

By a theorem of Kostant (see \cite[Chapter 10.6.6]{procesi2007liegroups} for a proof), as a closed subscheme of $\mathbb{P}(V(\varpi_i)^*)$, the generalized Grassmannian $\op{Gr}(i)$ is defined by the quadratic equations
\[x\otimes x \in \mathbb{P}(V(2\varpi_i)^*) \subseteq \mathbb{P}(V(\varpi_i)^*\otimes V(\varpi_i)^*).\] 
In other words, for every $\C$-algebra $R$, the $R$-points of $\op{Gr}(i)$ are rank $1$ direct summands\footnote{Or equivalently, rank $1$ projective quotients of $V(\varpi_i)_R$.} of $V(\varpi_i)^*_R := V(\varpi_i)^* \otimes_\C R$ whose tensor square is contained in $V(2\varpi_i)^*_R \subseteq V(\varpi_i)^*_R\otimes_R V(\varpi_i)^*_R$.\medskip\par

More generally, we can easily adapt the proof of Kostant's theorem to obtain the following description of the equations defining $\Fl_J$ inside $\prod_{i\in J}\PP(V(\varpi_i)^*)$\footnote{This result is folklore, but we include a proof for the convenience of the reader.}.

\begin{Theorem} \label{thm:quad-equations}
    For every $\C$-algebra $R$, an $R$-point $(x_i)_{i\in J}$ of $\prod_{i\in J} \PP(V(\varpi_i)^\ast)$ lies in $\Fl_J$ if and only if it satisfies, for every $i,j \in J$, the quadratic equation
    \[ x_i\otimes x_j \in \mathbb{P}(V(\varpi_i+\varpi_j)^*) \subseteq \mathbb{P}(V(\varpi_i)^*\otimes V(\varpi_j)^*).\]
\end{Theorem}

\begin{proof}
    By the Borel--Weil theorem, the multi-homogeneous coordinate ring of \[  \textstyle\Fl_J \subseteq \prod_{i\in J} \PP(V(\varpi_i)^*)\] is given by
  \[   \textstyle\bigoplus_{\lambda \in P_+^J} V(\lambda) \]
    where $P_+^J=\bigoplus_{i\in J} \N\varpi_i$ and the multiplication maps $ V(\lambda)  \otimes  V(\mu) \twoheadrightarrow V(\lambda+\mu)$ are projections onto highest weight isotypic components. This ring is clearly generated in degree 1 (i.e.~$\Fl_J$ is projectively normal inside $\prod_{i\in J} \PP(V(\varpi_i)^*)$), so we have a surjection
    \[  \textstyle p: \smash{\bigotimes_{i\in J} \op{Sym}^\bullet V(\varpi_i) \twoheadrightarrow \bigoplus_{\lambda \in P_+^J} V(\lambda)}\]
    of $P_+^J$-graded rings. We claim that the kernel of $p$ is generated in degree 2 as an ideal of the polynomial ring $\bigotimes_{i\in J} \op{Sym}^\bullet V(\varpi_i)$, i.e.~that it is generated by the kernels of~the~maps $V(\varpi_i) \otimes V(\varpi_j) \twoheadrightarrow V(\varpi_i+\varpi_j)$ for $i,j \in J$ with $i \neq j$, and $\op{Sym^2} V(\varpi_i) \twoheadrightarrow V(2\varpi_i)$ for $i \in J$ (which would imply that $\Fl_J$ is defined, as a closed subscheme of $\prod_{i\in J}\PP(V(\varpi_i)^*)$, by the corresponding quadratic equations, as wanted).%
    \medskip\par
    Let $\mathsf{I}$ be the ideal of the polynomial ring $A = \bigotimes_{i\in J} \op{Sym}^\bullet V(\varpi_i)$ generated by the %
    above quadratic elements, and let $B=A/\mathsf{I}$. We need to prove that the surjection from $B$ to $\bigoplus_{\lambda \in P_+^J} V(\lambda)$ is an isomorphism. Since this surjection is $G$-equivariant and respects $P_+^J$-grading, it is enough to check that the degree $\lambda$ part of $B$ is an irreducible right $G$-module with highest weight $\lambda$ for every $\lambda\in P_+^J$. This can be done by considering the action of the Casimir element $C\in (U\fg)^\fg$, which acts on $V(\lambda)$ by the scalar $C(\lambda)=(\lambda+\rho, \lambda+\rho)-(\rho, \rho)$.\medskip\par

    As $C$ can be written as a linear combination of products of pairs of elements of $\fg$,~we~can write the iterated coproducts $\Delta^{(n)}(C)\in (U\fg)^{\otimes n}$ in terms of $\Delta(C)$ and $C$ themselves,~i.e. %
    \[ \textstyle {\Delta^{(n)}(C) = \sum_{1 \le k < \ell \le n} \Delta(C)_{k,\ell} - (n-2)\sum_{1\leq k \leq n} C_k},\]
    with $\Delta(C)_{k,\ell}$ the image of $\Delta(C)\in (U\fg)^{\otimes 2}$ in $(U\fg)^{\otimes n}$ under the inclusion $(U\fg)^{\otimes 2} \hookrightarrow (U\fg)^{\otimes n}$ using the $k^\text{th}$ and $\ell^\text{th}$ tensor factors, and similarly for $C_k$.\medskip\par

    This allows us to express the action of $C$ on higher degree elements of $B$ in terms of the action on quadratic and linear elements. Moreover, by definition of $B$, the degree $\lambda$ part of $B$ is isomorphic to $V(\lambda)$, on which $C$ acts by the scalar $C(\lambda)$, whenever $\lambda$ is a sum of at most two fundamental weights. Therefore, given a sequence $a_1, \ldots, a_n$ of linear elements of $B$ of degrees $\varpi_{i_1}, \ldots, \varpi_{i_n}$, we have
    \[\textstyle  {C \cdot (a_1a_2\cdots a_n)= \Big(\,\sum_{1 \le k < \ell \le n} C(\varpi_{i_k}+\varpi_{i_\ell}) - (n-2) \sum_{1\leq k\leq n} C(\varpi_{i_k})\,\Big) a_1a_2\cdots a_n.} \]
    A straightforward calculation shows that the scalar on the right-hand side is equal to $C(\lambda)$, where $\lambda=\sum_{k=1}^n \varpi_{i_k}$, so $C$ acts on the degree $\lambda$ part of $B$ by the scalar $C(\lambda)$. Since a priori all irreducible representations that can appear in the degree $\lambda$ part of $B$ have highest weight $\le \lambda$ (as it is a quotient of $V(\varpi_{i_1})\otimes \cdots \otimes V(\varpi_{i_n})$), it follows that this degree $\lambda$ part is in fact irreducible of highest weight $\lambda$ by \cite[Lemma 10.6.6.2]{procesi2007liegroups}.
\end{proof}

Given two points $x, y$ of $\Gr(i), \Gr(j)$ respectively, we say that $x$ is \emph{incident} to $y$, denoted $x \sim y$, if they satisfy the condition
\[ x\otimes y \in \mathbb{P}(V(\varpi_i+\varpi_j)^*) \subseteq \mathbb{P}(V(\varpi_i)^*\otimes V(\varpi_j)^*)\]
(this definition applies to $R$-points for every $\C$-algebra $R$, and in fact defines a closed subscheme of $\Gr(i)\times \Gr(j)$). Then Theorem \ref{thm:quad-equations} is equivalent to saying that~$\op{Fl}_J$~parametrizes tuples $(x_i)_{i\in J} \in \prod_{i\in J}\op{Gr}(i)$ satisfying the incidence equations $x_i \sim x_j$ for all $i\neq j$ in $J$.

\begin{Example}\label{ex:SLnincidence}
    If $G=\op{SL}_n$, then $\op{Gr}(i)$ is the Grassmannian of $i$-dimensional subspaces of $\C^n$, and the incidence conditions correspond simply to inclusion of subspaces. Then $\Fl$ precisely parametrizes tuples of subspaces $(W_1, \ldots, W_{n-1})$ of $\C^n$ with $\dim W_i = i$, such that $W_i$ is contained in $W_j$ for every $i<j$.
\end{Example}\newpage

The above example suggests that not all incidence relations are necessary to define $\Fl_J$~in general, as some of them follow automatically from the others (for instance, we only need to impose consecutive inclusions $W_i \subseteq W_{i+1}$ in Example \ref{ex:SLnincidence}). This is indeed the case, as the following result of %
 \cite{magyar1998schubert} shows. Note that Magyar only works at the level of $\C$-points, but the proof can easily be adapted (as we do here) to show the scheme-theoretic~statement.
\begin{Lemma}
    Assume that $i,j,k$ are three vertices of $I$ such that $j$ lies on the path between $i$ and $k$. Then if $x_i, x_j,x_k$ are $R$-points of $\Gr(i), \Gr(j), \Gr(k)$ respectively (for some $\C$-algebra $R$) such that $x_i \sim x_j$ and $x_j \sim x_k$, then $x_i \sim x_k$.
\end{Lemma} 

\begin{proof}
For any $i\in I$, denote by $\widehat{P}_i^-$ the maximal parabolic subgroup $\smash{P_{I\setminus \{i\}}^-}$. It's not hard to see using the Bruhat decomposition that, at the level of $\C$-points, we have $\widehat{P}_j^- \subseteq \widehat{P}_i^- \widehat{P}_k^-$. As both sides are actually reduced subschemes of $G$, this inclusion in fact holds at a scheme-theoretic level, which means that any $R$-point of $\widehat{P}_j$ can be written as the product of an $R$-point of $\widehat{P}_i$ and an $R$-point of $\widehat{P}_k$, at the cost of replacing $R$ by an fppf extension of itself.\medskip\par

Now write $x_i= [v_{\varpi_i}^*] \cdot g_i$, $x_j=[v_{\varpi_j}^*]\cdot g_j$ and $x_k= [v_{\varpi_k}^*]\cdot g_k$, for some $g_i, g_j, g_k \in G(R)$ (which is possible by the definition of $\Gr(i)$ as a $G$-orbit, at the cost of passing to an fppf extension of $R$). Since $x_i \sim x_j$, the pair $(x_i, x_j)$ is a point of $\Fl_{\{i,j\}}$, which implies that $g_ig_j^{-1}\in \widehat{P}_i^-\widehat{P}_j^-$. Similarly, we have $g_jg_k^{-1}\in \widehat{P}_j^- \widehat{P}_k^-$, so $g_ig_k^{-1} \in \widehat{P}_i^-\widehat{P}_j^-\widehat{P}_k^-=\widehat{P}_i^-\widehat{P}_k^-$ (the last equality following from the previous paragraph). Therefore, after passing to an fppf extension and multiplying $g_i$ (resp. $g_k$) on the left by a suitable $R$-point of $\widehat{P}_i^-$ (resp. $\widehat{P}_k^-$), we can assume that $g_ig_k^{-1}=1$, so $(x_i, x_k)= ([v_{\varpi_i}^*], [v_{\varpi_k}^*])\cdot g$ for some $g\in G(R)$. The incidence condition $x_i\sim x_k$ follows immediately.
\end{proof}

\begin{Corollary} \label{cor:adj-incidence-eqs}
The full flag variety $\Fl=\Fl_I$ is the closed subscheme of $\prod_{i\in I} \Gr{(i)}$ defined by the incidence equations $x_i \sim x_j$ for all \emph{adjacent} vertices $i,j$ in $I$.
\end{Corollary}

\begin{Rem}
    While the proof of Theorem \ref{thm:quad-equations} goes by proving the corresponding result for the multi-homogeneous coordinate ring, Corollary \ref{cor:adj-incidence-eqs} is of a somewhat different flavor. It does not mean that the defining ideal of the multi-homogeneous coordinate ring is generated by the kernels of $\op{Sym^2} V(\varpi_i) \twoheadrightarrow V(2\varpi_i)$ and of $V(\varpi_i) \otimes V(\varpi_j) \twoheadrightarrow V(\varpi_i+\varpi_j)$ for $i$ adjacent to $j$ (in fact, this must be false simply for degree reasons).
\end{Rem}
\subsection{Bott--Samelson and free Bott--Samelson varieties}

In this section, we recall important facts about Bott--Samelson varieties. We will use a realization of Bott--Samelson varieties as subschemes of products of $\op{Gr}(i)$'s defined by incidence relations. This description (at a set-theoretic rather than scheme-theoretic level) first appears in \cite{magyar1998schubert} to the best of our knowledge.\medskip\par 

It will be convenient for us to think of Bott--Samelson varieties as associated primarily to combinatorial gadgets introduced by Viennot called \emph{heaps} \cite{viennot1986heaps}, which are in bijection with words in the alphabet $I$ up to commutation relations.

\begin{Def} \label{def:heap}
    A \textit{heap} is a finite partially ordered set $(H, \le)$ along with a map $c: H \to I$ satisfying
    \begin{enumerate}
        \item For $i,j\in I$  that are either adjacent or equal, $c^{-1}(\{i,j\}) \subseteq H$ is totally ordered. 
        \label{item:heap-cond1}
        \item For every covering relation $h \lessdot h'$ in $H$, $c(h)$ and $c(h')$ are either adjacent or equal. \label{item:heap-cond2}
    \end{enumerate}
\end{Def}

Given a heap $(H,\le, c)$, we will denote by $H_i$ the subset $c^{-1}(i)\subseteq H$, so that $H=\bigsqcup_{i\in I} H_i$. We say that a heap $H$ is \emph{alternating} if $H_i$ and $H_j$ are always perfectly interlaced, i.e.~for every edge $i\sim j$ of the Dynkin diagram and every $h, h'\in H_i$ with $h<h'$, there exists $h''\in H_j$ such that $h<h''<h'$. All the heaps that we will consider later in Section \ref{sec:bi-infinite-bott-samelson} will be alternating.

\begin{Def}
    Given two heaps $(H_1,\le_1, c_1)$ and $(H_2,\le_2,c_2)$, the product $H_1 \odot H_2$ is the heap with underlying set $H_1 \sqcup H_2$, whose map to $I$ extends $c_1:H_1\to I$ and $c_2:H_2\to I$, and whose partial order is defined as the transitive closure of the relation defined by $h \le h'$ if either
    \begin{itemize}
        \item $h,h'\in H_1$ and $h\le_1 h'$, or
        \item $h,h'\in H_2$ and $h\le_2 h'$, or
        \item $h\in H_1$, $h'\in H_2$ and $c_1(h)$, $c_2(h')$ are either adjacent or equal.
    \end{itemize}
    Given a word $\ul{i}=(i_1, \ldots, i_r) \in I^r$, let
    \[ H(\ul{i}) = H(i_1) \odot  \cdots \odot H(i_r),\]
    where $H(i)$ is the heap with a unique element $*$ and $c(*)=i$.
\end{Def}

\begin{Lemma} \cite[Lemmas 3.2 and 3.3]{viennot1986heaps}
    The map $\ul{i} \mapsto  H(\ul{i})$ defines a bijection between words in $I$ up to commutation relations (i.e. the exchange of $i_k$ and $i_{k+1}$ in $\ul{i}$ if $i_k \not\sim i_{k+1}$) and heaps (up to isomorphism). 
    Also, given a heap $H$, the set of words $\ul{i}$ such that $H\cong H(\ul{i})$ is in bijection with linear extensions of $H$, i.e.~enumerations $h_1, \ldots, h_r$ of the elements of $H$ such that $h_k \le h_\ell$ implies $k\le \ell$, where a linear extension $h_1, \ldots, h_r$ corresponds to the word $(c(h_1), \ldots, c(h_r))$.
\end{Lemma}

One way to represent heaps visually is via the ``abacus model'' \cite{kleshchev2010homogeneous,dranowski2024heaps}, where the elements on $H$ are represented by beads on an abacus with rods indexed by $I$, where the rods are placed in such a way that the beads on two adjacent rods are forced to interlace (i.e.~the distance between two adjacent rods is smaller than the diameter of the beads). Any pairs of beads touching each other correspond to a covering relation in the poset $H$. Given a word $\ul{i}=(i_1, \ldots,i_r)$, the corresponding heap $H(\ul{i})$ is obtained by dropping a bead on rods $i_1$, then one on rod $i_2$, and so on until the last bead on rod $i_r$. Figure \ref{fig:abacus}, for example, shows a heap of type $A_3$ using this representation. 
\definecolor{bead_color}{gray}{0.8}
\definecolor{rod_color}{gray}{0}
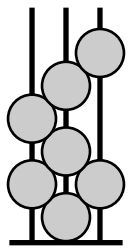
\begin{figure}[ht]
    \centering
    \begin{tikzpicture}[line width=2pt,scale=0.75]
\def\r{0.4} %
\def\hspace{0.6} %
\def\h{6*\d+3*\r+5*\vspace} %
\def\vspace{sqrt(4*\r^2-\hspace^2} %
\def\d{0.05} %

\foreach \i in {0,1,2} {
	\draw[rod_color] (\i*\hspace,0) -- (\i*\hspace,\h);
}

\draw[rod_color] (-\r,0) -- (2*\hspace+\r,0);

\draw (\hspace,\d+\r) circle (\r);
\fill[bead_color] (\hspace,\d+\r) circle (\r);
\draw (0,2*\d+\r+\vspace) circle (\r);
\fill[bead_color] (0,2*\d+\r+\vspace) circle (\r);
\draw (2*\hspace,2*\d+\r+\vspace) circle (\r);
\fill[bead_color] (2*\hspace,2*\d+\r+\vspace) circle (\r);
\draw (\hspace,3*\d+\r+2*\vspace) circle (\r);
\fill[bead_color] (\hspace,3*\d+\r+2*\vspace) circle (\r);
\draw (0,4*\d+\r+3*\vspace) circle (\r);
\fill[bead_color] (0,4*\d+\r+3*\vspace) circle (\r);
\draw (\hspace,5*\d+\r+4*\vspace) circle (\r);
\fill[bead_color] (\hspace,5*\d+\r+4*\vspace) circle (\r);
\draw (2*\hspace,6*\d+\r+5*\vspace) circle (\r);
\fill[bead_color] (2*\hspace,6*\d+\r+5*\vspace) circle (\r);

\end{tikzpicture}
    \caption{The heap $H(2,1,3,2,1,2,3)=H(2,3,1,2,1,2,3)$. This heap is not alternating, as there is no bead on the $3^{\text{rd}}$ rod between the highest~beads of the $2^{\text{nd}}$ rod (but adding a bead there would yield an alternating heap).}
    \label{fig:abacus}
\end{figure}
\begin{Def} \label{def:Gamma_H}
    Given a finite linearly ordered set $S$, let $I(S)$ denote the set of decompositions $S=L\sqcup U$ of $S$ into a lower set $L$ and an upper set $U$, ordered by 
    $$(L,U) \le (L',U') \Leftrightarrow L \subseteq L' \Leftrightarrow U' \subseteq U.$$ 
    This yields a finite linearly ordered set with one more element than $S$.
    Given %
    $s\in S$, we denote by $s_-$ the element of $I(S)$ corresponding to the decomposition 
    \[ S = \{ s'\in S:s'<s\} \sqcup \{ s'\in S:s'\ge s\}\]
    and by $s_+$ the element of $I(S)$ corresponding to the decomposition 
    \[ S = \{ s'\in S:s'\le s\} \sqcup \{ s'\in S:s'> s\}.\]
    
    Given a heap $H$, we associate an $I$-coloured graph $\Gamma_H$ whose set of vertices is $\bigsqcup_{i\in I} I(H_i)$. There is an edge between $(L,U)\in I(H_i)$ and $(L',U') \in I(H_j)$ if and only if $i\sim j$ and $L\sqcup L'$ is a lower set of $H_i\sqcup H_j$ (equivalently $U\sqcup U'$ is an upper set of $H_i\sqcup H_j$). We will also denote by $c:\Gamma_H\to I$ the map which sends $I(H_i)$ to $i$.
\end{Def}

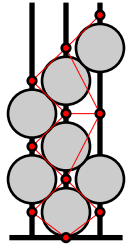
\begin{figure}[ht]
    \centering
    \begin{tikzpicture}[line width=2pt,scale=0.75]

\def\r{0.4} %
\def\hspace{0.6} %
\def\h{6*\d+3*\r+5*\vspace} %
\def\vspace{sqrt(4*\r^2-\hspace^2} %
\def\d{0.05} %
\def\vertexR{0.075}

\foreach \i in {0,1,2} {
  \draw[rod_color] (\i*\hspace,0) -- (\i*\hspace,\h);
}

\draw[rod_color] (-\r,0) -- (2*\hspace+\r,0);

\draw (\hspace,\d+\r) circle (\r);
\fill[bead_color] (\hspace,\d+\r) circle (\r);
\draw (0,2*\d+\r+\vspace) circle (\r);
\fill[bead_color] (0,2*\d+\r+\vspace) circle (\r);
\draw (2*\hspace,2*\d+\r+\vspace) circle (\r);
\fill[bead_color] (2*\hspace,2*\d+\r+\vspace) circle (\r);
\draw (\hspace,3*\d+\r+2*\vspace) circle (\r);
\fill[bead_color] (\hspace,3*\d+\r+2*\vspace) circle (\r);
\draw (0,4*\d+\r+3*\vspace) circle (\r);
\fill[bead_color] (0,4*\d+\r+3*\vspace) circle (\r);
\draw (\hspace,5*\d+\r+4*\vspace) circle (\r);
\fill[bead_color] (\hspace,5*\d+\r+4*\vspace) circle (\r);
\draw (2*\hspace,6*\d+\r+5*\vspace) circle (\r);
\fill[bead_color] (2*\hspace,6*\d+\r+5*\vspace) circle (\r);

\draw[red,thin] (\hspace,0) -- (0,\d+\r);
\draw[red,thin] (\hspace,0) -- (2*\hspace,\d+\r);
\draw[red,thin] (2*\hspace,\d+\r) -- (\hspace,2*\d+\r+\vspace);
\draw[red,thin] (0,\d+\r) -- (\hspace,2*\d+\r+\vspace);
\draw[red,thin] (\hspace,2*\d+\r+\vspace) -- (0,3*\d+\r+2*\vspace);
\draw[red,thin] (\hspace,2*\d+\r+\vspace) -- (2*\hspace,4*\d+\r+3*\vspace);
\draw[red,thin] (0,3*\d+\r+2*\vspace) -- (\hspace,4*\d+\r+3*\vspace);
\draw[red,thin] (2*\hspace,4*\d+\r+3*\vspace) -- (\hspace,4*\d+\r+3*\vspace);
\draw[red,thin] (\hspace,4*\d+\r+3*\vspace) -- (0,5*\d+\r+4*\vspace);
\draw[red,thin] (0,5*\d+\r+4*\vspace) -- (\hspace,6*\d+\r+5*\vspace);
\draw[red,thin] (2*\hspace,4*\d+\r+3*\vspace) -- (\hspace,6*\d+\r+5*\vspace);
\draw[red,thin] (\hspace,6*\d+\r+5*\vspace) -- (2*\hspace,7*\d+\r+6*\vspace);

\filldraw[fill=red,draw=black,line width=1pt] (\hspace,0) circle (\vertexR);
\filldraw[fill=red,draw=black,line width=1pt] (0,\d+\r) circle (\vertexR);
\filldraw[fill=red,draw=black,line width=1pt] (2*\hspace,\d+\r) circle (\vertexR);
\filldraw[fill=red,draw=black,line width=1pt] (\hspace,2*\d+\r+\vspace) circle (\vertexR);
\filldraw[fill=red,draw=black,line width=1pt] (0,3*\d+\r+2*\vspace) circle (\vertexR);
\filldraw[fill=red,draw=black,line width=1pt] (\hspace,4*\d+\r+3*\vspace) circle (\vertexR); 
\filldraw[fill=red,draw=black,line width=1pt] (2*\hspace,4*\d+\r+3*\vspace) circle (\vertexR); 
\filldraw[fill=red,draw=black,line width=1pt] (0,5*\d+\r+4*\vspace) circle (\vertexR); 
\filldraw[fill=red,draw=black,line width=1pt] (\hspace,6*\d+\r+5*\vspace) circle (\vertexR); 
\filldraw[fill=red,draw=black,line width=1pt] (2*\hspace,7*\d+\r+6*\vspace) circle (\vertexR);

\end{tikzpicture}
\caption{The graph $\Gamma_{H}$ for $H=H(2,1,3,2,1,2,3)$.}
\label{fig:Gamma_H}
\end{figure}

For example, Figure \ref{fig:Gamma_H} shows in red the graph $\Gamma_H$ for the heap from Figure \ref{fig:Gamma_H}. A vertex in $I(H_i)$ correspond to a decomposition $(L,U)$ is represented by a dot on rod $i$ above all beads of $L$ and below all beads of $U$.

\vspace*{-2.75mm}

\begin{Rem} \label{rem:alt-heaps}
    For alternating heaps, $\Gamma_H$ is almost the same as the Hasse diagram of $H$. Specifically, the map $h\mapsto h_+$ gives an isomorphism between the Hasse diagram of $H$ and the induced subgraph of $\Gamma_H$ obtained by deleting the minimal element $p_i^{\min}=(\varnothing, H_i)$ of $I(H_i)$ for every $i$. Moreover, $p_i^{\min}$ is connected to $h_+$ if and only if $i\sim c(h)$ and there are no elements in $H_i$ smaller than $h$.
    Dually, $h\mapsto h_-$ gives an isomorphism between the Hasse diagram of $H$ and the induced subgraph of $\Gamma_H$ obtained by deleting the maximal element of $I(H_i)$ for every $i$. Outside the alternating case, however, the graph $\Gamma_H$ can look quite different from the Hasse diagram of $H$.
\end{Rem}
\begin{Def}
    Given a heap $H$, the \textit{Bott--Samelson variety} $Z_H$ is the subscheme of 
 \begin{equation}\label{eq:ambiantprodBS}
 \textstyle\prod_{p\in \Gamma_H} \op{Gr}(c(p))
 \end{equation}
defined by the equations
    \begin{enumerate}
        \item \label{item:BS-incidence-eq} $x_p \sim x_q$ for every edge $(p,q)$ of the graph $\Gamma_H$, 
        \item \label{item:BS-boundary-eq} $x_{p_i^{\min}} = [v_{\varpi_i}^*]$, for every $i\in I$, where $p_i^{\min}=(\varnothing, H_i)$ is the minimal element of $I(H_i)$.
    \end{enumerate}
    Here $x_p \in \op{Gr}(c(p))$ is the component indexed by $p\in \Gamma_H$ of a tuple in the above product.
    
    Similarly, the \textit{free Bott--Samelson variety} $\cZ_H$ is the subscheme of \eqref{eq:ambiantprodBS} defined only by the incidence equations (\ref{item:BS-incidence-eq}), without the boundary conditions (\ref{item:BS-boundary-eq}).
    
    For $p\in \Gamma_H$, we denote by $\pi_p: Z_H\to \op{Gr}(c(p))$ the projection to the factor indexed by $p$. We also denote $\pi_p: \cZ_H\to \op{Gr}(c(p))$ the corresponding projection for the free version.
    
  From the abacus perspective, in $ Z_H $, we begin with the standard flag and every time we drop a bead on the $ i^{\text{th}}$ rod, we change the $ \Gr(i) $ component of our flag.  In $ \cZ_H $, we begin with an arbitrary flag and proceed in the same way.
    
\end{Def}

Note that there are obvious right actions of $B_-$ on $Z_H$ and $G$ on $\cZ_H$ induced by the action of $G$ on $\op{Gr}(i)$. It is easy to see from Corollary \ref{cor:adj-incidence-eqs} that the map $Z_H \times^{B_-} G\to \cZ_H$ induced by the action of $G$ is an isomorphism. \medskip\par

Note also that, for every heap $H$ and every $i\in I$, the coloured graph $\Gamma_{H(i)\odot H}$ can be identified with the graph obtained from $\Gamma_H$ by adjoining a minimal $i$-coloured element $p_i^{\min}$ that is connected to the minimal $j$-coloured elements of $\Gamma_H$ for every $j\sim i$. It then also follows easily from Corollary \ref{cor:adj-incidence-eqs} that, for every heap $H$, the map 
$$ Z_H \times^{B_-} P_i^- \to Z_{H(i)\odot H}$$ 
sending $((x_p),g)$ to the point $(x_p')_{p\in \Gamma_{H(i)\odot H} = \{p_i^{\min}\} \sqcup \Gamma_H} \in Z_{H(i)\odot H}$
of coordinates $x_p'=x_p \cdot g$ (for $p\in \Gamma_H$ and $x'_{p_i^{\min}}=[v_{\varpi_i}^*]$) is an isomorphism.\medskip\par 
    
As a consequence, we obtain, for any word $\ul{i}$, an isomorphism
\begin{equation} \label{eq:bott-samelson-iso}
    Z_{\ul{i}}:=Z_{H(\ul{i})} \cong B_-\backslash P_{i_r}^- \times^{B_-} \cdots \times^{B_-} P_{i_2}^- \times^{B_-} P_{i_1}^-,
\end{equation} 
which is (up to our conventions about right actions) the usual definition of Bott--Samelson varieties. In particular, this shows that $Z_H$ is a smooth variety of dimension $\#H$ (and $\cZ_H$ is a smooth variety of dimension $\#H + \dim (B_-\backslash G)$).\medskip\par
Denote by $\Gamma_H^+\subseteq \Gamma_H$ the subset $\bigsqcup_i (I(H_i) \setminus\{p_i^{\min}\})$ of vertices corresponding to $\op{Gr}(i)$ factors that are not frozen to $[v_{\varpi_i}^*]$ in the definition of $Z_H$. For every $p\in \Gamma_H^+$, there is a line bundle $\mathcal{O}_p(1)$ on $Z_H$ obtained by pulling back the antitautological line bundle $\O(1)$ under \[\pi_p: Z_H \to \Gr(c(p)) \hookrightarrow \mathbb{P}(V(\varpi_{c(p)})^*).\] It is not hard to see that these line bundles correspond to the $\O(1)$-basis of \cite[Section 3.1]{lauritzen2004line} (under the isomorphism (\ref{eq:bott-samelson-iso})). In particular, their classes form a basis for the Picard group of $Z_H$. Similarly, the line bundles $\mathcal{O}_p(1)$ on $\cZ_H$ defined similarly (but for any $p\in \Gamma_H$, possibly $p=p_i^{\min}$) form a basis for the Picard group of $\cZ_H$. Given $\ul{m}=(m_p)\in \Z^{\Gamma_H^+}$ (resp. $\Z^{\Gamma_H}$), we will denote by $\cL_{\ul{m}}$ the line bundle $\bigotimes_{p\in \Gamma_H^+} \O_p(m_p)$ (resp. $\bigotimes_{p\in \Gamma_H} \O_p(m_p)$) on $Z_H$ (resp. $\cZ_H$), where $\O_p(m):=\O_p(1)^{\otimes m}$. Note that the line bundles $\O_p(1)$ (hence also all $\cL_{\ul{m}}$) on $Z_H$ come naturally equipped with the structure of $B_-$-equivariant line bundles. Similarly, the line bundles $\cL_{\ul{m}}$ on $\cZ_H$ are naturally $G$-equivariant line bundles.\medskip\par
For any $h \in H$, let $H\setminus \{h\}$ denote the heap with underlying set $H\setminus \{h\}$, equipped with the finest partial order that is coarser than the original order of $H$ while satisfying the heap condition (\ref{item:heap-cond2}) from Definition \ref{def:heap} (equivalently, $h'\le h''$ if and only if there exists a chain $h'=h_0\le h_1\le \cdots \le h_t=h''$ of elements of $H\setminus \{h\}$ such that $c(h_k)$ is adjacent or equal to $c(h_{k+1})$ for every $0\le k < t$). The corresponding graph $\Gamma_{H\setminus \{h\}}$ is obtained from $\Gamma_H$ by merging the two vertices $h_-$ and $h_+$.\medskip\par

Then, for all $h \in H$, we have a closed embedding $Z_{H\setminus \{h\}} \hookrightarrow Z_{H}$ whose image is~the~locus given by the equation $\pi_{h_-}(x) = \pi_{h_+}(x)$. This defines a prime divisor $D_h$ in $Z_H$. The union \[ \textstyle \bigcup_{h\in H} D_h\] is a simple normal crossing divisor. We denote its complement by $Z^\circ_H$, called the open Bott--Samelson cell. If $H\cong H({\ul{i}})$ with $\ul{i}$ a \emph{reduced} word for $w\in W$, then the map $Z_H \to \Fl$ sending $x$ to $(\pi_{p_i^{\max}}(x))_{i\in I} \in \Fl$ (where $p_i^{\max}=(H_i,\varnothing)$ is the maximal element of $I(H_i)$) induces an isomorphism from $Z_H^\circ$ to the Schubert cell indexed by $w$ in $\Fl$ (i.e. the $B_-$-orbit of $([v_{w\varpi_i}^*])_i=([v_{\varpi_i}^*])_i \cdot \dot{w}^{-1}$) \cite[Chapter II.13.6]{jantzen2003representations}.\medskip\par

We can also define divisors $D_h$ in $\cZ_H$ indexed by $h\in H$ in the same way. Their union is still a simple normal crossing divisor and its complement is denoted $\cZ_H^\circ$. If $H\cong H(\ul{i})$ with $\ul{i}$ a \emph{reduced} word for $w\in W$, then the map $\cZ_H \to \Fl\times \Fl$ sending $x$ to $((\pi_{p_i^{\min}}(x))_{i}, (\pi_{p_i^{\max}}(x))_{i}) \in \Fl\times \Fl$ induces an isomorphism from $\cZ_H^\circ$ to the Schubert cell indexed by $w$ in $\Fl\times \Fl$ (i.e. the $G$-orbit of $(([v_{\varpi_i}^*])_i, ([v_{w\varpi_i}^*])_i)=(([v_{\varpi_i}^*])_i, ([v_{\varpi_i}^*])_i \cdot \dot{w}^{-1})$).

\begin{Rem} \label{rem:R-points-open-cell}
Let $ (L_p)_{p \in \Gamma_H} $ be an $ R$-point of $ \cZ_H$, where $ R $ is a $ \C$-algebra and $ L_p $~is~a~rank $1$ direct summand of $V(\varpi_{c(p)})^*_R$.  This point lies in the open subscheme $\cZ_H^\circ$ if and only if $L_{h_-}\cap L_{h_+}=\{0\}$ and $L_{h_-}\oplus L_{h_+}$ is a rank $2$ direct summand of $V(\varpi_{c(h)})^*_R$ for every $h\in H$.
 \end{Rem}

\begin{Lemma} \label{lem:BS-divisor-expansion}
   For any heap $H$ and any $h\in H_i$, the expansion of the line bundle $\O(D_h)$ on $\cZ_H$ in the $\O(1)$ basis is given by
   \[ \textstyle \O(D_h) \cong \O_{h_-}(1) \otimes \O_{h_+}(1) \otimes \bigotimes_{j\sim i} \O_{h_j}(-1),\]
   with $h_-, h_+\in I(H_i)$ as in Definition \ref{def:Gamma_H} and where $h_j\in I(H_j)$ corresponds to %
   \[H_j = \{h'\in H_j:h'< h\} \sqcup \{h'\in H_j:h'> h\}.\]
   The same expansion holds on $Z_H$, with the convention that $\O_p(1)$ is the trivial line bundle if $p$ is the minimal element of $I(H_i)$.
\end{Lemma}

\begin{proof}
    Pick a connected induced subgraph $T\subseteq \Gamma_H$ such that
    \begin{itemize}
        \item $T\cap I(H_j)$ is a singleton for any $j\neq i$.
        \item $T\cap I(H_j) = \{h_j\}$ for $j\sim i$.
        \item $T\cap I(H_i) = \{h_-,h_+\}$.
    \end{itemize}
    Then $T$ can be identified with $\Gamma_{H(i)}$ as an $I$-coloured graph, so we have a projection map $\cZ_H \to \cZ_{H(i)}=\cZ_{i}$  given by remembering the coordinates $x_t$ for $t\in T$. This map pullbacks $\O_t(1)$ to $\O_t(1)$ for all $t\in T$ and pullbacks the unique boundary divisor of $\cZ_{i}$ to $D_h$. Thus, without loss of generality we can assume $H\cong H(i)$ (i.e. $h$ is the unique element~of~$H$). 
    
    The composition 
    $$\textstyle \bigwedge^2V(\varpi_i)^* \twoheadrightarrow V(2\varpi_i-\alpha_i)^*\hookrightarrow \bigotimes_{j\sim i} V(\varpi_j)^*$$ (the first map being the dual of \eqref{eq:Gequivariant_map_wedge} and the second being the right inverse of the dual of \eqref{eq:Gequivariant_map_from_identity}) induces a map of sheaves
    \[ \textstyle \O_{h_-}(-1) \otimes \O_{h_+}(-1) \to \O_{\cZ_H} \otimes \bigotimes_{j\sim i} V(\varpi_j)^*.\]
    Moreover, we claim that the image of this map is contained in the subbundle \[\textstyle \bigotimes_{j \sim i} \O_{h_j}(-1) \subseteq \O_{\cZ_H} \otimes \bigotimes_{j\sim i} V(\varpi_j)^*.\]
    To prove the claim, note that it's enough to check it on the dense open set $\cZ_H^\circ$. Since $G$ acts transitively on $\cZ_H^\circ$ and the map above is $G$-equivariant, it is in fact enough to check it at %
    the point given by $x_{h_-} = [v_{s_i\varpi_i}^*] = [v_{\varpi_i-\alpha_i}^*]$ and $x_p=[v_{\varpi_{c(p)}}^*]$ for every other $p\in \Gamma_H$, in which case the claim is clear by comparing weights (using \eqref{eq:varpi_svarpi}).
    
    We therefore have a map
    \[ \textstyle \O_{h_-}(-1) \otimes \O_{h_+}(-1) \to \bigotimes_{j \sim i} \O_{h_j}(-1), \]
    or equivalently a section of 
    \[ \textstyle \O_{h_-}(1) \otimes \O_{h_+}(1) \otimes \bigotimes_{j \sim i} \O_{h_j}(-1), \]
    which is nonvanishing on $\cZ_H^\circ$ (again by transitivity of the $G$-action) but vanishes to order $1$ on the boundary divisor $D_h$. This yields the desired isomorphism. The statement about $Z_H$ is an immediate consequence of the one about $\cZ_H$.
\end{proof}

\subsection{Global sections of line bundles on (free) Bott--Samelson varieties}

It is proven in \cite{lauritzen2004line} that the line bundle $\cL_{\ul{m}}$ on $Z_H$ is globally generated if and only if $m_p \ge 0$ for all $p$. Moreover, a description of the space of global sections, as a $B_-$-module, is given in this globally generated case in \cite{lakshmibai2002standard}: it is dual to a generalized Demazure module, which is a submodule of $\bigotimes_{p\in \Gamma_H^+} V(\varpi_{c(p)})^*$ for which an explicit character formula is known. \medskip\par

To state this formula, we need the Demazure operators $\Lambda_i$ for $i\in I$, which are operators on the group ring $\Z [P]\simeq \Z[e^{\la}\,|\,\la\in P]$ given by
\[ \Lambda_i (f) = \frac{e^{\alpha_i}f -  s_i(f)}{e^{\alpha_i}-1}.\]
These operators satisfy braid relations, so we can also define operators $\Lambda_w$ for any $w\in W$ by $\Lambda_w = \Lambda_{i_1} \circ  \cdots \circ \Lambda_{i_r}$ if $w=s_{i_1}\cdots s_{i_r}$ is a reduced expression for $w$.

\begin{Theorem}[\cite{lakshmibai2002standard}] \label{thm:demazure-character-for-BS}
    Fix a heap $H$ with $h_1,\ldots, h_r$ a linear extension of $H$ corresponding to an isomorphism $H\cong H(\ul{i})$. If $\ul{m}\in \mathbb{N}^{{\Gamma_H^+}}$, the character of $H^0(Z_H, \cL_{\ul{m}} )$ (for the left action of $T\subseteq B_-$) is given by
    \[ \Lambda_{i_1}(e^{m_{h_{1+}}\varpi_{i_1}} \Lambda_{i_2}(e^{m_{h_{2+}}\varpi_{i_2}} \cdots \Lambda_{i_r}(e^{m_{h_{r+}}\varpi_{i_r}}) \cdots )).\]
\end{Theorem}
One consequence of the above character formula proven in \cite{lakshmibai2002standard} is a projective normality result for $Z_H$. To state this result, consider, for $\ul{m}\in \mathbb{N}^{\Gamma_H^+}$, the map
\begin{align*}
Z_H &\hookrightarrow \textstyle \prod_{p\in \Gamma_H^+} \Gr(c(p)) \hookrightarrow \prod_{p\in \Gamma_H^+} \mathbb{P}(V(\varpi_{c(p)})^*)\\
& \to\textstyle  \prod_{p\in \Gamma_H^+} \mathbb{P}(V(\varpi_{c(p)})^*)^{m_p} \hookrightarrow \mathbb{P}\Big( \bigotimes_{p\in \Gamma_H^+} \left(V(\varpi_{c(p)})^*\right)^{\otimes m_p}\Big),
\end{align*}
where the third map is a product of diagonal maps and the last one is the Segre embedding. The pullback of $\mathcal{O}(1)$ under this map is $\cL_{\ul{m}}$, so we have an induced map
\begin{equation} \label{eq:BS-proj-embedding}
   \textstyle \bigotimes_{p\in \Gamma_H^+} V(\varpi_{c(p)})^{\otimes m_p} \longrightarrow H^0(Z_H, \cL_{\ul{m}}).
\end{equation}

\begin{Theorem}[\cite{lakshmibai2002standard}] \label{thm:BS-proj-normality}
    The map (\ref{eq:BS-proj-embedding})
    is surjective.
\end{Theorem}

\begin{Corollary} \label{cor:BS-mult-surj}
    Let $H$ be a heap and let $\ul{m}\in \mathbb{N}^{\Gamma_H^+}$. The multiplication map
    \[ \textstyle \bigotimes_{p\in \Gamma_H^+} H^0(Z_H, \mathcal{O}_p(1))^{\otimes m_p} \to H^0(Z_H, \cL_{\ul{m}})\]
    is surjective.
\end{Corollary}

\begin{proof}

We have a commutative diagram
$$\adjustbox{scale=0.9}{
\begin{tikzcd}[column sep=1em,row sep=0.75em]
\bigotimes_{p\in \Gamma_H^+} (V(\varpi_{c(p)}))^{\otimes m_p} \arrow[dd, two heads] \arrow[rr, Rightarrow, no head] &  & \bigotimes_{p\in \Gamma_H^+} V(\varpi_{c(p)})^{\otimes m_p} \arrow[dd, two heads] \\
                                                         &  &                          \\
\bigotimes_{p\in \Gamma_H^+} H^0(Z_H, \mathcal{O}_p(1))^{\otimes m_p} \arrow[rr]                                            &  & H^0(Z_H, \cL_{\ul{m}})                      
\end{tikzcd}}$$
where the two vertical arrows are surjective by the above theorem. Thus, the bottom arrow must be surjective as well.
\end{proof}

While the results above give a good understanding of global sections of line bundles on Bott--Samelson varieties $Z_H$, we are most interested for the purposes of this paper~in~studying line bundles on $\cZ_H$. We can however reduce questions about the latter to questions about the former with the help of the following %
special case of \cite[Theorem 2.20]{fujita2020flag}:

\begin{Theorem} \label{thm:BS-to-GBS-global-sections-iso}
    Let $H$ be a heap and $\ul{i}$ be a reduced word for the longest element $w_0 \in W$. The map \[Z_{H(\ul{i}) \odot H}\to \cZ_H\]
    induced by the obvious embedding of $\Gamma_H$ in $\Gamma_{H(\ul{i})\odot H}$ induces a $B_-$-equivariant isomorphism
    \[ H^0(\cZ_H, \cL_{\ul{m}}) \overset{\sim}{\longrightarrow} H^0(Z_{H(\ul{i})\odot H}, \cL_{\ul{m}})\]
    for every $\ul{m}\in \Z^{\Gamma_H}$.
\end{Theorem}

\begin{Corollary} \label{cor:GBS-mult-surj}
    Let $\ul{m}\in \mathbb{N}^{\Gamma_H}$. The multiplication map
    \[ \textstyle \bigotimes_{p\in \Gamma_H} H^0(\cZ_H, \mathcal{O}_p(1))^{\otimes m_p} \to H^0(\cZ_H, \cL_{\ul{m}})\]
    is surjective.
\end{Corollary}

\begin{proof}
    This is immediate from the combination of Corollary \ref{cor:BS-mult-surj} and Theorem \ref{thm:BS-to-GBS-global-sections-iso}
\end{proof}
	
\section{The bi-infinite Bott--Samelson pro-variety and its open cell}\label{sec:bi-infinite-bott-samelson}

In this section, we introduce the bi-infinite Bott--Samelson pro-variety $\cZ_{\infty}$, which informally parametrizes tuples 
\[ (x_{i,a})\in \prod_{(i, a) \in I\times _2 \Z}\Gr(i)\]
such that $x_{i,a} \sim x_{j,a+1}$ for all $(i,a)\in I\times_2\Z$ and $j\sim i$.

\subsection{An inverse limit of Bott--Samelson varieties}
The precise definition of $\cZ_{\infty}$ is:
\begin{Def} \label{def:Z}
    For two height functions $\xi, \xi'$ with $\xi\le \xi'$ (i.e.~with $\xi_i\le \xi_i'$ for all $i\in I$), let $\cZ_{\xi, \xi'}$ be the closed subscheme of \[ \prod_{\substack{(i, a) \in I\times _2 \Z \\ \xi_i \le a \le \xi'_i}}\Gr(i)\]
    defined by the equations $x_{i,a} \sim x_{j,a+1}$
    whenever $i\sim j$, $\xi_i \le a \le \xi'_i$ and $\xi_j \le a+1 \le \xi'_j$, where $x_{i,a}\in \op{Gr(i)}$ denotes the $(i,a)$-component of a tuple in the above product. We denote by $\pi_{i,a}:\cZ_{\xi,\xi'}\to \Gr(i)$ the projection onto the factor indexed by $(i,a)$ and by $\pi_{\xi''}: \cZ_{\xi,\xi'}\to \Fl$ the map given, for any height function $\xi''$ with $\xi\le \xi''\le \xi'$, by
    \[  x \mapsto (\pi_{i,\xi''_i}(x))_{i\in I} \in\cZ_{\xi'',\xi''}= \Fl\subseteq\prod_{i\in I} \op{Gr}(i).\]
There are obvious projections $\smash{\pi^{\xi_1,\xi_1'}_{\xi_2,\xi_2'}}: \cZ_{\xi_1, \xi_1'} \to \cZ_{\xi_2,\xi'_2}$ whenever $\xi_1 \le \xi_2\le \xi_2'\le \xi_1'$, so the varieties $\cZ_{\xi,\xi'}$ form an inverse system. We define the pro-variety (see Appendix \ref{app:pro_varieties})
    \[ \cZ_{\infty} = \lim_{\leftarrow} \cZ_{\xi, \xi'}\]
    as the formal inverse limit of this inverse system. We have maps $\pi_{i,a}:\cZ_\infty \to \op{Gr(i)}$ and $\pi_\xi: \cZ_\infty \to \Fl$ for every $(i,a)\in I\times_2 \Z$ and every height function $\xi$ induced from the maps with the same name defined above.
\end{Def}

The variety $\cZ_{\xi,\xi'}$ is a free Bott--Samelson variety $\cZ_{H(\xi,\xi')}$ for the heap $H(\xi,\xi')$, for which the vertex set is
$$\{(i,b)\in I\times^{\mathrm{op}}_2 \Z\,|\, \xi_i< b < \xi_i'\},$$ 
the partial order is the one defined in Section \ref{sec:Notation}, and the map to $I$ is given by projection onto the first factor. Indeed, the graph $\Gamma_{H(\xi,\xi')}$~can be identified with the set 
$$\{(i,a)\in I\times_2 \Z\,|\,\xi_i \le a \le \xi_i'\}$$ 
by matching a pair $(i,a)$ with%
\[ H(\xi,\xi')_i = \{(i,b):\xi_i <b<a\} \sqcup \{(i,b):a<b<\xi_i'\}\]
in such a way that the edges of $\Gamma_{H(\xi,\xi')}$ match the incidence conditions %
defining $\cZ_{\xi,\xi'}$.

\begin{Rem} \label{rem:Hxi-alternating}
    The heaps $H(\xi,\xi')$ are clearly all alternating. Conversely, it is not hard to see that every alternating heap is isomorphic to $H(\xi,\xi')$ for some %
    height functions $\xi \le \xi'$.
\end{Rem}
For every $(i, b) \in I\timesop_2 \Z$ with $\xi_i < b < \xi'_i$, we have a prime divisor $D_{i,b}$ on $\cZ_{\xi,\xi'}$ defined by the equation $\pi_{i,b-1}(x)=\pi_{i,b+1}(x)$. These divisors are all compatible with the projections $\pi^{\xi_1,\xi_1'}_{\xi_2,\xi_2'}$, so they define divisors $D_{i,b}$ on $\cZ_\infty$ indexed by $I\timesop_2 \Z$.
\begin{Def}
    Given %
    height functions $\xi, \xi'$ with $\xi\le \xi'$, let $\cZ^\circ_{\xi, \xi'} \subseteq \cZ_{\xi, \xi'}$ be the complement of the simple normal crossing divisor
    \[ \bigcup_{\substack{(i, b) \in I\timesop_2 \Z \\ \xi_i < b < \xi'_i}} D_{i,b}.\] 
    Also, define \[{\cZ_\infty^\circ= \lim_\leftarrow \cZ^\circ_{\xi, \xi'},}\]
    with transition maps given by the restrictions of the transition maps $\smash{\pi^{\xi_1,\xi_1'}_{\xi_2,\xi_2'}}$ for $\cZ_\infty$.
\end{Def}

\begin{Rem}
    The transition maps in the definition of $\cZ_\infty^\circ$ are affine morphisms (in fact %
     affine space bundles), so the limit $\cZ_\infty^\circ$ is not only a pro-variety but can also be regarded as an honest scheme by Remark \ref{rem:affine-transition-inverse-limits}. This does not hold for $\cZ_\infty$.
\end{Rem}

\begin{Example}\label{ex:example_Zinft_SL5}
	The pro-variety $\cZ_\infty$ for $G=\SL_5$ parametrizes fillings of the bi-infinite diagram shown in Figure \ref{fig:Zinfty-A4-example}, where each $W_{i,a}$ is an $i$-dimensional subspace of $\C^5$ and every edge represents an inclusion. The open cell $\cZ_\infty^\circ$ parametrizes those fillings for which $W_{i,a+2}$ always differs from $W_{i,a}$. When $\xi=(-1,-2,-1,-2)$ and $\xi'=(1,2,3,2)$, the variety $\cZ_{\xi,\xi'}$ parametrizes fillings of the (finite) red part of the diagram from Figure \ref{fig:Zinfty-A4-example}.
\end{Example}
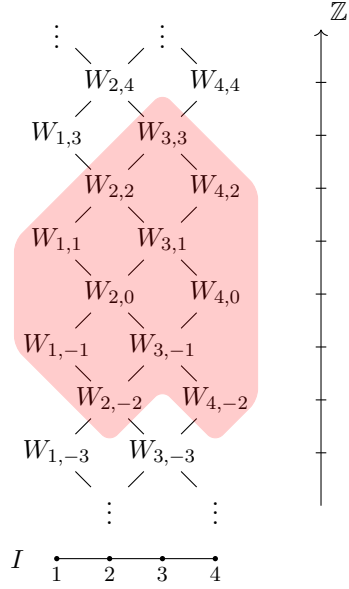
\begin{figure}[ht] 
\begin{tikzpicture}[scale=0.7,every node/.style={scale=0.9}]
\def\rad{0.5}
\def\xnudge{0.8}
\def\ynudge{0.8}

\foreach \i in {2} {
	\foreach \x in {-4,-2,0,2,4} {
			\ifnum\x>-4
				\draw (\i,\x) -- (\i+1,\x-1);
			\else
			\fi
			\draw (\i,\x) -- (\i+1,\x+1);
	}
}

\foreach \i in {1,3} {
	\foreach \x in {-3,-1,1,3,5} {
		\ifnum\x<5
			\draw (\i,\x) -- (\i+1,\x+1);
		\else
		\fi
		\draw (\i,\x) -- (\i+1,\x-1);

	}
}

\foreach \i in {1,3} {
	\foreach \x in {-3,-1,1,3,5} {
		\fill[fill=white] (\i,\x) circle (\rad);
		\ifnum\x<5
			\node at (\i,\x) {$W_{\i,\x}$};
		\else
		\fi
	}
}

\foreach \i in {2,4} {
	\foreach \x in {-4,-2,0,2,4} {
		\fill[fill=white] (\i,\x) circle (\rad);
		\ifnum\x>-4
			\node at (\i,\x) {$W_{\i,\x}$};
		\else
		\fi
	}
}

\filldraw[rounded corners,red,opacity=0.2] 	(1-\xnudge,-1) -- (2,-2-\ynudge) -- (3,-1-\ynudge) -- (4,-2-\ynudge) -- 
		(4+\xnudge,-2) -- (4+\xnudge,2) -- (3,3+\ynudge) -- (1,1+\ynudge) -- (1-\xnudge,1) -- cycle;
		
\node at (2,-4) {$\vdots$};
\node at (4,-4) {$\vdots$};
\node at (1,5) {$\vdots$};
\node at (3,5) {$\vdots$};

\def\hI{-5}
\def\sizeV{0.05}

\draw (4,\hI) -- (1,\hI);
\foreach \i in {1,2,3,4} {
	\fill (\i,\hI) circle (\sizeV) node[below] {\footnotesize $\i$};
}

\node at (1-0.75,\hI) {$I$};
\draw[->] (6,-4)--(6,5) node [above right] {$\Z$};
\foreach \n in {-3,...,4} {
	\draw (5.9,\n)-- (6.1,\n);
}
\end{tikzpicture}
\caption{The incidence diagram for $\cZ_\infty$ in type $A_4$.}
\label{fig:Zinfty-A4-example}
\end{figure}\vspace*{-3mm}

\begin{Lemma} \label{lem:wedge-eq-tensor}
    For every point $(x_{i,a})_{i,a\in I\times_2 \Z}$ of $\cZ_\infty^\circ$ and every $(i,b)\in I\timesop_2\Z$, the wedge $x_{i,b-1}\wedge x_{i,b+1}$ lies in the highest isotypic component $V(2\varpi_i-\alpha_i)^* \subseteq \bigwedge\nolimits^{\!2} V(\varpi_i)^*$, the tensor product $\bigotimes_{j\sim i} x_{j,b}$ lies in the highest isotypic component $V(2\varpi_i-\alpha_i)^* \subseteq \bigotimes_{j\sim i} V(\varpi_j)^*$, and %
    \[\textstyle x_{i,b-1}\wedge x_{i,b+1} = \bigotimes_{j\sim i} x_{j,b} \]
    in $\mathbb{P}(V(2\varpi_i-\alpha_i)^*)$.
\end{Lemma}
\begin{proof}
    Pick two height functions $\xi\le \xi'$ that differ only in position $i$, such that $\xi_i=b-1$, $\xi_i'=b+1$, and therefore $\xi_j=\xi'_j=b$ for every $j\sim i$. Then this is really a statement about the finite-dimensional variety $\cZ_{\xi,\xi'}^\circ$. The group $G$ acts transitively on $\cZ_{\xi,\xi'}^\circ$ (by the isomorphism with the Schubert cell in $\Fl \times \Fl$), so we may assume without loss of generality that $x_{i,b-1}=[v_{\varpi_i}^*]$, $x_{i,b+1}=[v_{s_i\varpi_i}^*]$ and $x_{j,b}=[v_{\varpi_j}^*]$ for $j\sim i$. But in that case the lemma is clear just by comparing weights (using \eqref{eq:varpi_svarpi}).
\end{proof}
\begin{Rem} For any height function $\xi$, the heap $H(\xi, \xi+2)$ is isomorphic to $H(\ul{i})$,~where $\ul{i}$ is any reduced word for the Coxeter element $c$ associated to $\xi$. Thus $\cZ_{\xi, \xi+2}^\circ$ is isomorphic to the Schubert cell associated to $c$ in $\Fl\times \Fl$. It follows that, if we fix a height function $\xi$, we obtain an isomorphism $x\mapsto (\pi_{\xi+2s}(x))_{s\in \Z}$ between $\cZ^\circ_\infty$ and the pro-variety parametrizing infinite sequences of points in $\Fl$, any two consecutive ones having relative position $c$. Interestingly, this description is quite reminiscent of the notion of $q$-opers~\cite{frenkel2024qopers}.
\end{Rem}

\begin{Rem} \label{rem:Z-and-bands}
	The Hasse diagram of the heap $H(\xi, \xi^*+h)$ (see Section \ref{subsubsec:height_functions} for notation) is isomorphic to the Auslander-Reiten quiver of the path algebra whose quiver corresponds to the Dynkin diagram of $\fg$, together with the orientation corresponding to $\xi$ (this is clear from the description of the Auslander-Reiten quiver given in \cite[Section 2.2]{ringel1980tame}).~It~is~well-known that the corresponding commutation-class of words consists of reduced words for $w_0$ \cite{bedard1999commutation}. It follows that $\pi_\xi \times \pi_{\xi^*+h} : \cZ_{\xi,\xi^\ast+h}^\circ \to \Fl \times \Fl$ induces an isomorphism between $\cZ_{\xi,{\xi^*+h}}^\circ$ and the open Schubert cell in $\Fl \times \Fl$. The group $G$ acts transitively on this open cell and the stabilizer of the point $(([v_{\varpi_i^*}])_i, ([v_{w_0\varpi_i}^*])_i)$ is $T$, so there is an isomorphism~from $T\backslash G$ to that cell sending $[g]$ to $(([v_{\varpi_i}^*])_i, ([v_{w_0\varpi_i}^*])_i) \cdot g$.
\medskip\par
    Furthermore, a point of $\cZ_{\xi, \xi^*+h}^\circ$ corresponding to a pair $(p_1, p_1')\in \Fl \times \Fl$ and a point of $\cZ_{\xi+2, \xi^*+h+2}^\circ$ corresponding to a pair $(p_2, p_2') \in \Fl \times \Fl$ agree on $\cZ_{\xi+2, \xi^*+h}^\circ$ if and only if
    \begin{itemize}
        \item $p_1, p_2$ are in relative position $c$,
        \item $p_2, p_1'$ are in relative position $c^{-1}w_0$ and 
        \item $p_1', p_2'$ are in relative position $c^*:=w_0cw_0$.
    \end{itemize}\smallskip
    If $(p_1, p_1') = (([v_{\varpi_i}^*])_i , ([v_{w_0\varpi_i}^*])_i)\cdot g_1$ and $(p_2, p_2') = (([v_{\varpi_i}^*])_i, ([v_{w_0\varpi_i}^*])_i)\cdot g_2$, then these three conditions are equivalent to
    \begin{align*}
        g_1g_2^{-1} &\in B_-cB_- \cap  BcB_- \cap BcB = Bc \cap cB_-.
    \end{align*}
    It follows from this observation that for $m \in \N$, $\cZ_{\xi, \xi^*+h+2m}^\circ$ can be identified with the variety that parametrizes sequences $([g_s])_{s=0}^m \in (T\backslash G)^{m+1}$ of elements $[g_s]\in T\backslash G$ satisfying
    \begin{equation} \label{eq:bands-condition}
        g_s g_{s+1}^{-1} \in Bc \cap cB_-.
    \end{equation}
    In the limit, we get that $\cZ_\infty^\circ$ parametrizes sequences $([g_s])_{s\in \Z} \in (T\backslash G)^{\Z}$ satisfying (\ref{eq:bands-condition}). \medskip\par
    
    In \cite{francone2025cluster}, Francone and Leclerc defined a scheme $B(G,c)$, called the \emph{scheme of bands}, which parametrizes infinite sequences $(g_s)_{s\in \Z} \in G^{\Z}$ satisfying $g_s g_{s+1}^{-1} \in N\dot{c} \cap \dot{c}N_-$. There is a free left $T$-action on $B(G,c)$ given by $t\cdot (g_s)_{s\in \Z} = (c^{-s}tc^s g_s)_{s\in \Z}$. It follows from the discussion above that a choice of height function $\xi$ with corresponding Coxeter element $c$ induces a $G$-equivariant isomorphism between $\cZ_\infty^\circ$ and the quotient $T\backslash B(G,c)$. Explicitly, this isomorphism sends $(g_s)_{s \in \Z}$ to the point $(x_{i,a}) \in \cZ_\infty^\circ$ determined by $x_{i,\xi_i+2s} = [v_{\varpi_i}^*] \cdot g_s$.
\end{Rem}
\subsection{The Cox ring of $\cZ_\infty^\circ$ and the universal torus bundle}\label{subsec:cox_ring_of_Zinftycirc}

Since $\cZ_{\xi,\xi'}$ is a free Bott--Samelson variety, its Picard group is a free abelian group with basis given by the classes of $\O_{i,a}(1)$ for $(i,a)\in I\times_2 \Z$, $\xi_i \le a\le \xi_i'$, where $\O_{i,a}(1)$ is the pullback of $\O(1)$ via
\[ \cZ_\infty \overset{\pi_{i,a}}{\longrightarrow} \Gr(i) \hookrightarrow \mathbb{P}(V(\varpi_i)^*).\]
Taking the limit, we get that the Picard group of $\cZ_\infty$ is free abelian with basis given by the classes $[\O_{i,a}(1)]$ for all $(i,a)\in I\times_2 \Z$. We identify $ \Pic(\cZ_\infty) $ with $ \cB $ by sending $ [\O_{i,a}(1)] $ to the monomial $ y_{i,a} $.\medskip\par

As the complement of all divisors $D_{i,b}$ for $\xi_i< b < \xi_i'$, the open cell $\cZ_{\xi,\xi'}^\circ$ has Picard group given by the quotient of $\op{Pic}(\cZ_{\xi,\xi'})$ by the classes of $\O(D_{i,b})$ for $\xi_i< b < \xi_i'$ (this is clear from the isomorphism of the Picard groups of $\cZ_{\xi,\xi'}$ and $\cZ_{\xi,\xi'}^\circ$ with their divisor class group). By Lemma \ref{lem:BS-divisor-expansion}, this class is given by $[\O_{i,b-1}(1)]+[\O_{i,b+1}(1)]-\sum_{j\sim i} [\O_{j,b}(1)]$. In the limit, we obtain that $\op{Pic}(\cZ_{\infty}^\circ)$ is the quotient of the free abelian group on the classes $[\O_{i,a}(1)]$ for $(i,a)\in I\times_2 \Z$ by the subgroup generated by
\[ \textstyle [\O_{i,b-1}(1)]+[\O_{i,b+1}(1)]-\sum_{j\sim i} [\O_{j,b}(1)]\]
for $(i,b)\in I\timesop_2 \Z$. We can (and will) therefore identify $\op{Pic}(\cZ_\infty^\circ)$ with the finitely generated free abelian group $\cP = \cB/\Gamma$, by sending the class $[\O_{i,a}(1)]$ to $\tau_{i,a}:=\awt(y_{i,a})$. 

\begin{Rem} \label{rem:growth-of-pic-and-sections}
    It's clear from the discussion above (or the fact that the projection maps are affine space bundles) that all $\pi^{\xi_1,\xi_1'}_{\xi_2,\xi_2'}: \cZ_{\xi_1,\xi_1'}^\circ \to \cZ_{\xi_2,\xi_2'}^\circ$ induce isomorphisms on Picard groups. Hence the Picard group of $\cZ_{\xi,\xi}^\circ$ is independent of $\xi, \xi'$, though the space of sections of a given line bundle does grow as $\xi$ and $\xi'$ get further apart.\medskip\par
    
    The inverse system $(\cZ_{\xi,\xi'})$, on the other hand, has the opposite behaviour. The Picard group grows as a function of $\xi, \xi'$, but the space of sections is independent of $\xi, \xi'$, in the sense that the projections $\pi^{\xi_1,\xi_1'}_{\xi_2,\xi_2'}: \cZ_{\xi_1,\xi_1'} \to \cZ_{\xi_2,\xi_2'}$ induce isomorphisms\vspace*{-1mm}
    \[ {H^0(\cZ_{\xi_2,\xi_2'}, \mathcal{L}) \overset{\sim}{\to} H^0(\cZ_{\xi_1,\xi_1'}, \pi^{\xi_1,\xi_1'*}_{\xi_2,\xi_2'}\mathcal{L})}\]
    for every line bundle $\mathcal{L}$ on $\cZ_{\xi_2,\xi_2'}$.
    This follows from the projection formula \cite[Exercise II.5.1.(d)]{hartshorne1977ag} together with the fact that these projections are proper with connected fibres, hence pushforward the structure sheaf to the structure sheaf \cite[\href{https://stacks.math.columbia.edu/tag/0AY8}{Tag 0AY8}]{stacks-project}.
\end{Rem}

We now fix for every $i\in I$ a choice of identification of the isotypic components of highest weight $2\varpi_i-\alpha_i$ in $\bigwedge\nolimits^{\!2} V(\varpi_i)^*$ and $\bigotimes_{j\sim i} V(\varpi_j)^*$. This choice is unique up to scalar, but we can fix it by requiring that $v_{\varpi_i}^*\wedge v_{s_i \varpi_i}^*$ goes to $\otimes_{j\sim i} v_{\varpi_j}^*$ (dually to the identification in Section \ref{subsec:cat_O_of_QQ}). As can be seen from the proof of Lemma \ref{lem:BS-divisor-expansion}, this identification provides us with a choice of trivialization of \[\textstyle \O_{i,b-1}(1) \otimes \O_{i,b+1}(1) \otimes \bigotimes_{j\sim i} \O_{j,b}(-1)\]
on $\cZ_\infty^\circ$ for every $(i,b)\in I\timesop_2 \Z$. Now, given a class $\tau = \sum_{i,a} m_{i,a}\tau_{i,a} \in \cP = \op{Pic}(\cZ_\infty^\circ)$, we consider the line bundle $\cL_\tau=\bigotimes_{i,a} \O_{i,a}(m_{i,a})$ which has isomorphism class $\tau$. Note that two different expansions of $\tau$ as linear combinations of classes $\tau_{i,a}$ yield \emph{canonically} isomorphic line bundles $\cL_\tau$ (due to our choice of trivialization of $\O_{i,b-1}(1) \otimes \O_{i,b+1}(1) \otimes \bigotimes_{j\sim i} \O_{j,b}(-1)$), so it makes sense to denote any of these line bundles by $\cL_\tau$. There are obvious isomorphisms $\cL_{\tau} \otimes \cL_{\tau'} \overset{\sim}{\to} \cL_{\tau+\tau'}$ for every $\tau, \tau'\in \cP$ that satisfy an associativity condition.
\begin{Def} \label{def:cox-ring}
    The Cox ring of $\cZ_\infty^\circ$ is the $\cP$-graded ring
    \[\textstyle \cR = \bigoplus_{\tau\in \cP} H^0(\cZ_\infty^\circ, \cL_{\tau}),\]
    with multiplication given by tensor product of sections.
\end{Def}

Note that we have $\cR =\displaystyle\lim_\to \cR_{\xi, \xi'}$, where
\[ \textstyle \cR_{\xi, \xi'} =  \bigoplus_{\tau\in \cP} H^0(\cZ_{\xi, \xi'}^\circ, \cL_{\tau}).\]
Here, note that we can make sense of $\cL_{\tau}$ as a line bundle over $\cZ_{\xi, \xi'}^\circ$ since $(\tau_{i,\xi_i})_{i\in I}$ spans $\cP$ for any height function $\xi$.\medskip\par

For any smooth variety $X$ whose Picard group is free of finite rank and with only constant invertible functions, there is a universal principal torus bundle $\hat{X} \to X$, characterized by the property that taking associated line bundles yields an isomorphism from the character lattice of the torus acting on $\hat{X}$ to $\op{Pic}(X)$. The bundle $\hat{X}$ can be constructed as the relative spectrum of the Cox sheaf of $X$, and 
the ring of functions of $\hat{X}$ is given by the Cox ring of $X$ \cite[Chapter 1.6]{arzhantsev2015coxrings}. For example, $\widehat{\Gr}(i)$ is the complement of the origin in the affine cone of $\op{Gr}(i)\subseteq \mathbb{P}(V(\varpi_i)^*)$ (whose $R$-points are elements $x\in V(\varpi_i)_R^*$ such that $Rx$ is a rank $1$ direct summand of $ V(\varpi_i)_R^*$ and $x\otimes x \in V(2\varpi_i)_R^*$). We now give an explicit presentation of this universal torus bundle for the bi-infinite Bott--Samelson cell $\cZ_\infty^\circ$.

\begin{Def} \label{def:Zhat}
    Given height functions $\xi\le \xi'$, let $\widehat{\cZ}^\circ_{\xi,\xi'}$ be the closed subscheme of \[ \prod_{\substack{(i, a) \in I\times _2 \Z \\ \xi_i \le a \le \xi'_i}}\widehat{\Gr}(i)\]
    parametrizing arrays $(x_{i,a})$ satisfying the equations $x_{i,a} \sim x_{j,a+1}$, i.e.
    \begin{equation*}
        x_{i,a}\otimes x_{j,a+1} \in V(\varpi_i+\varpi_j)^* \subseteq V(\varpi_i)^*\otimes V(\varpi_j)^*,
    \end{equation*} 
    whenever $\xi_i \le a \le \xi'_i$, $\xi_j \le a+1 \le \xi'_j$, and $i\sim j$, as well as
    \begin{equation*}
       \textstyle x_{i,b-1} \wedge x_{i,b+1} = \bigotimes_{j\sim i} x_{j,b}
    \end{equation*}
    (as elements of $V(2\varpi_i-\alpha_i)^*$) for every $(i,b)\in I\timesop_2\Z$ with $\xi_i < b < \xi'_i$.\medskip\par 

    Let \[\widehat{\cZ}_\infty^\circ= \displaystyle\lim_{\leftarrow}\widehat{\cZ}^\circ_{\xi,\xi'}.\]
\end{Def}

Note that we use our fixed identification of the isotypic components of highest weight $2\varpi_i-\alpha_i$ in $\bigwedge\nolimits^{\!2} V(\varpi_i)^*$ and $\bigotimes_{j\sim i} V(\varpi_j)^*$ for every $i\in I$ to make sense of the last condition. Furthermore, by Lemma \ref{lem:wedge-eq-tensor}, this last condition $x_{i,b-1} \wedge x_{i,b+1} = \otimes_{j\sim i} x_{j,b}$ is automatically satisfied up to a scalar; we only require that this scalar is $1$.\medskip\par
Recall the definition of the torus $A$ from Section \ref{subsec:equiv_of_inclusions}. By definition, the character lattice of $A$ is $\cP$ and
points of $A$ are arrays $(t_{i,a})\in \mathbb{G}_m^{I\times_2 \Z}$ satisfying
\[\textstyle t_{i,b-1} \cdot  t_{i,b+1} = \prod_{j\sim i} t_{j,b}\]
for each $(i,b)\in I\timesop_2 \Z$. There is an action of $A$ on $\widehat{\cZ}_\infty^\circ$ via $(t_{i,a})_{i,a}\cdot (x_{i,a})_{i,a}= (t_{i,a}x_{i,a})_{i,a}$.\medskip\par

Notice that by construction, the scheme $\widehat{\cZ}^\circ_{\xi,\xi}$ is isomorphic to $N_{-}\backslash G$ and the projection $\widehat{\cZ}^\circ_{\xi,\xi}\to \cZ^\circ_{\xi,\xi}$ is identified with the principal $T$-bundle $N_-\backslash G\to B_-\backslash G$. This generalizes as:

\begin{Lemma} \label{lem:universal-torus-bundle}
    For any $\xi \le \xi'$, the projection $\widehat{\cZ}^\circ_{\xi,\xi'} \to \cZ^\circ_{\xi,\xi'}$ is a principal $A$-bundle. Moreover, for each $\xi_1 \le \xi_2\le \xi_2'\le \xi_1'$, the square
    \begin{equation} \label{eq:universal-bundle-cart-square}
    \adjustbox{scale=0.9}{
    \begin{tikzcd}[row sep = 1em]
	{\widehat{\cZ}^\circ_{\xi_1,\xi_1'}} & {\widehat{\cZ}^\circ_{\xi_2,\xi_2'}} \\
	{\cZ^\circ_{\xi_1,\xi_1'}} & {\cZ^\circ_{\xi_2,\xi_2'}}
	\arrow[from=1-1, to=1-2]
	\arrow[from=1-1, to=2-1]
	\arrow[from=1-2, to=2-2]
	\arrow[from=2-1, to=2-2]
    \end{tikzcd}
    }
    \end{equation}
    is Cartesian. Finally, for every $\tau \in \cP$, the associated line bundle on $\cZ^\circ_{\xi,\xi'}$ obtained from~the weight $\tau \in \cP =\mathbb{X}^\bullet(A)$ and the principal $A$-bundle $\widehat{\cZ}^\circ_{\xi,\xi'}$ is $\cL_{\tau}$.\medskip\par
    
    In the limit, the projection $\widehat{\cZ}^\circ_{\infty}\to \cZ^\circ_{\infty}$ is a principal $A$-bundle for which the associated bundle to a weight $\tau\in \cP$ is $\cL_{\tau}$.
\end{Lemma}

\begin{proof}
    We start by proving that (\ref{eq:universal-bundle-cart-square}) is Cartesian. It is enough to prove this in the case where either $\xi_1=\xi_2$ and $\xi_1'$ differs from $\xi_2'$ only in a single index $i\in I$, or $\xi_1'=\xi_2'$ and $\xi_1$ differs from $\xi_2$ only 
    at $i\in I$.
    Let us assume that we are in the first case (the argument~is~the same in the second case), and let $b\in \Z$ be such that $\xi_{1,i}'=b+1$ and $\xi_{2,i}'=b-1$.\medskip\par
    
    Let $R$ be a $\C$-algebra and suppose that we are given an $R$-point $(x_{i,a})$ of $\widehat{\cZ}_{\xi_2,\xi_2'}^\circ$, as well as an $R$-point of $\Gr(i)$ (given by a rank $1$ direct summand $L\subseteq V(\varpi_i)_R^*$), satisfying the incidence equations $L\sim Rx_{j,b}$ for every $j\sim i$ and such that $L\oplus Rx_{i,b-1}$ is a rank $2$ direct summand of $V(\varpi_i)_R^*$ (so that $((Rx_{i,a}), L)$ is an $R$-point of $\cZ^\circ_{\xi_1,\xi_1'}$). Then, by Lemma \ref{lem:wedge-eq-tensor}, we have an equality $x_{i,b-1}\wedge L=R \cdot \otimes_{j\sim i} x_{j,b}$ of rank $1$ direct summands of $V(2\varpi_i-\alpha_i)^*_R$. This means that there is a unique generator $x_{i,b+1}\in L$ that satisfies $x_{i,b-1}\wedge x_{i,b+1}=R \cdot \otimes_{j\sim i} x_{j,b}$. In other words, there is a unique way to complete $(x_{i,a})$ to a point of $\widehat{\cZ}_{\xi_1,\xi_1'}^\circ$ whose image in $\cZ_{\xi_1,\xi_1'}^\circ$ is given by $((Rx_{i,a}), L)$, meaning that (\ref{eq:universal-bundle-cart-square})~is~Cartesian.\medskip\par

    It is clear that $\widehat{\cZ}_{\xi,\xi}^\circ\to \cZ_{\xi,\xi}^\circ = \Fl$ is a principal $A$-bundle for any height function $\xi$. Hence $\widehat{\cZ}_{\xi,\xi'}^\circ\to \cZ_{\xi,\xi'}^\circ$ is a principal $A$-bundle for any $\xi\le \xi'$ by the Cartesian square (\ref{eq:universal-bundle-cart-square}).\medskip\par

    Finally, whenever $\xi_i\le a \le \xi_i'$, the associated bundle to the weight $\tau_{i,a}\in \cP$ relative to the principal $A$-bundle $\widehat{\cZ}^\circ_{\xi,\xi'}\to \cZ^\circ_{\xi,\xi'}$ is canonically isomorphic to $\O_{i,a}(1)$. This is straightforward to check when $\xi=\xi'$, and the general case follows again from the Cartesian square (\ref{eq:universal-bundle-cart-square}). Thus, for each $\tau\in \cP$, the associated bundle is canonically isomorphic to $\cL_{\tau}$.
\end{proof}

\begin{Corollary}\label{cor:iso_R_CZ}
	For any pair of height functions $\xi\leq \xi'$, there is a canonical isomorphism $\cR_{\xi,\xi'} \cong \C[\widehat{\cZ}^\circ_{\xi,\xi'}]$ which, upon taking the limit, yields an isomorphism $\cR \cong \C[\widehat{\cZ}^\circ_\infty]$.
\end{Corollary}
\begin{proof}
Let $X=\widehat{\cZ}^\circ_{\xi,\xi'}$ for notational convenience, and let 
$$\O_{X} =  \textstyle \bigoplus_{\tau\in \cP} \O_{X,\tau}$$
be the weight decomposition of the structure sheaf of $X$ induced by the $A$-action. The associated line bundle to $\tau\in \cP$ is given by the pushforward $p_* \O_{X,\tau}$, where $p:X\to \cZ^\circ_{\xi,\xi'}$ is the projection. By Lemma \ref{lem:universal-torus-bundle}, this associated bundle can be identified with $\cL_{\tau}$,~and taking global sections before summing over all $\tau$ yields the first isomorphism. The second isomorphism follows as the above identifications are compatible with the maps $\pi^{\xi_1,\xi_1'}_{\xi_2,\xi_2'}$.
\end{proof}

\begin{Rem} \label{rem:Zhat-and-bands}
    Recall from Remark \ref{rem:Z-and-bands} that a choice of height function $\xi$ induces an isomorphism between $\cZ_\infty^\circ$ and $T\backslash B(G,c)$, where $B(G,c)$ is the scheme of bands of \cite{francone2025cluster}. Moreover, the principal $T$-bundle $B(G,c)\to T\backslash B(G,c)$ has the property that taking associated line bundles yield an isomorphism between the weight lattice of $T$ and the Picard group of $T\backslash B(G,c)$. Specifically, it's not hard to see that the associated bundle for the weight $\varpi_i$ is isomorphic to $\O_{i,\xi_i}(1)$ under the isomorphism $T\backslash B(G,c) \cong \cZ_\infty^\circ$. This isomorphism of weight lattices induces an isomorphism between $T$ and $A$, and the isomorphism of associated bundles implies that the principal torus bundles $B(G,c)$ and $\widehat{\cZ}_\infty^\circ$ are themselves isomorphic. Therefore, our universal torus bundle $\widehat{\cZ}_\infty^\circ$ gives a realization of Francone and Leclerc's scheme of bands that is independent of the choice of Coxeter element $c$.  In fact, the isomorphism $ B(G, c) \rightarrow \widehat{\cZ}_\infty^\circ$ is explicitly given by $ (g_s) \mapsto (x_{i,a}) $ where $ x_{i,\xi_i + 2s} = v_{\varpi_i}^* \cdot g_s$.
\end{Rem}

\begin{Lemma} \label{lem:Zhat-torsor}
	For every height function $\xi$, the variety $\widehat{\cZ}_{\xi,\xi^*+h}^\circ$ is a $G$-torsor.
\end{Lemma}

\begin{proof}
	Compare the Cartesian squares
	\begin{equation*}
	\adjustbox{scale=0.94}{
		\begin{tikzcd}[row sep = 1.25em]
			{\widehat{\cZ}^\circ_{\xi,\xi^*+h}} & {\widehat{\cZ}^\circ_{\xi,\xi}} \\
			{\cZ^\circ_{\xi,\xi^*+h}} & {\cZ^\circ_{\xi,\xi}}
			\arrow[from=1-1, to=1-2]
			\arrow[from=1-1, to=2-1]
			\arrow[from=1-2, to=2-2]
			\arrow[from=2-1, to=2-2]
		\end{tikzcd}}
		\hspace{1em}
		\text{and}
		\hspace{1em}
		\adjustbox{scale=0.94}{
		\begin{tikzcd}[row sep = 1.25em]
		{G} & {N_-\backslash G} \\
		{T\backslash G} & {B_-\backslash G}
		\arrow[from=1-1, to=1-2]
		\arrow[from=1-1, to=2-1]
		\arrow[from=1-2, to=2-2]
		\arrow[from=2-1, to=2-2]
	\end{tikzcd}}
	\end{equation*}
We have $G$-equivariant isomorphisms $\cZ^\circ_{\xi,\xi}\cong B_-\backslash G$, $\cZ^\circ_{\xi,\xi^*+h}\cong T \backslash G$ and $\widehat{\cZ}^ \circ_{\xi,\xi}\cong N_- \backslash G$~that are
compatible with the maps in both diagrams. Hence %
$\widehat{\cZ}^\circ_{\xi,\xi^*+h}\cong G$ as $G$-varieties.
\end{proof}

\begin{Rem}\label{rem:embeddings_of_CG} 
We can trivialize the $G$-torsor $\widehat{\cZ}_{\xi,\xi^*+h}^\circ$ by picking the base point $p_0=(x_{i,a})$ given by $x_{i, \xi_i + 2s} = v_{\varpi_i}^* \cdot \dot{c}^{-s}$ (which is the image of the point $(g_s=\dot{c}^{-s})_{s\in \Z}$ of $B(G,c)$ under the map $B(G,c)\overset{\sim}{\to} \widehat{\cZ}_\infty^\circ \twoheadrightarrow \widehat{\cZ}_{\xi,\xi^*+h}^\circ$). The resulting isomorphism $G \overset{\sim}{\to} \widehat{\cZ}_{\xi,\xi^*+h}^\circ$ fits into the following commutative diagram
\[\adjustbox{scale=0.94}{ \begin{tikzcd}[row sep = 1.25em]
	{B(G,c)} & {\widehat{\cZ}_{\infty}^\circ} \\
	G & {\widehat{\cZ}_{\xi,\xi^*+h}^\circ}
	\arrow["\sim", from=1-1, to=1-2]
	\arrow[two heads, from=1-1, to=2-1]
	\arrow[two heads, from=1-2, to=2-2]
	\arrow["\sim", from=2-1, to=2-2]
\end{tikzcd}}\]
in which the left vertical map simply takes $(g_s)_{s\in \Z}\in B(G,c)$ to $g_0$.
\end{Rem}

It follows from Lemma \ref{lem:Zhat-torsor} that $\widehat{\cZ}_{\xi,\xi^*+h}^\circ$ is an affine scheme. Hence $\widehat{\cZ}_{\xi,\xi'}^\circ$, for $\xi' \ge \xi^*+h$, and $\widehat{\cZ}_{\infty}^\circ$ itself, are also affine since they are affine bundles over $\widehat{\cZ}_{\xi,\xi^*+h}^\circ$. This is not obvious from Definition \ref{def:Zhat}, which presents $\widehat{\cZ}_{\xi,\xi'}^\circ$ as a locally closed subscheme of $\prod_{i,a} V(\varpi_i)^*$, hence a priori only as a quasi-affine scheme. It turns out however that $\widehat{\cZ}_\infty^\circ$ is actually a \emph{closed} subscheme of $\prod_{(i,a)\in I\times_2\Z} V(\varpi_i)^*$, as we now explain.

For $i\in I$, let  $\langle\cdot, \cdot\rangle : V(\varpi_i)^*\times V(\varpi_{i^*})^* \to \C$ be the unique $G$-invariant pairing~for which 
$$\langle v_{\varpi_i}^\ast,v_{-\varpi_i}^\ast\rangle=1.$$ For $(i,a) \in I \times_2 \Z$ and $k \in \Z$, let $f_{k,a}^{(i)} \in \C[\widehat{\mathcal{Z}}^\circ_{\infty}]$ be the function given by
\begin{equation}\label{eq:def_G_invariant_functions}
f_{k,a}^{(i)}(x)=\langle x_{i,a-2k}, x_{i^*, a+h} \rangle.
\end{equation}
By construction, this function is $G$-invariant.

\begin{Lemma}\label{lem:constant_function_1}
For every $(i,a)\in I\times_2 \Z$, the function $\smash{f_{0,a}^{(i)}}$ is constant and equal to $1$.
\end{Lemma}

\begin{proof}
Choose a height function $\xi$ such that $\xi_i=a$, thus making $\smash{f_{0,a}^{(i)}}$ a $G$-invariant function on $\widehat{\cZ}^\circ_{\xi,\xi^*+h}$. Consider the isomorphism $\widehat{\cZ}^\circ_{\xi,\xi^*+h}\simeq G$ of Remark \ref{rem:embeddings_of_CG}. 
For $x\in %
\widehat{\mathcal{Z}}^\circ_{\xi,\xi^*+h}$,~let $g\in G$ be such that $p_0g=x$, with $p_0$ the base point specified by $\xi$. Using the combinatorics of Section \ref{subsubsec:height_functions}, we get
\begin{equation*}
f_{0,a}^{(i)}(x)=\langle x_{i,\xi_i}, x_{i^*,\xi_i+h} \rangle=\langle x_{i,\xi_i}, x_{i^*,\xi_{i^\ast}+2m_{i^\ast}} \rangle=\langle v_{\varpi_i}^\ast,v_{\varpi_{i^\ast}} \dot{c}^{-m_{i^*}}\rangle=\langle v_{\varpi_i}^\ast,v_{-\varpi_{i}} \rangle=1,
\end{equation*}
where the second to last equality follows from \cite[Lemma 5.9]{francone2025cluster}.
\end{proof}

\begin{Rem}\label{rem:GinvariantsGeom} 
Fix $(i,a) \in I \times_2 \Z$ and $k\in \Z_{\geq 1}$. 
After picking a height function $\xi$, the map $f_{k,a}^{(i)}$ coincides, under the isomorphism of Remark \ref{rem:Zhat-and-bands}, with the function $\theta_{i,k}^{(s)}\in \C[B(G,c)]$
of \cite[Section 7.1]{francone2025cluster}, where $a=\xi_i+2(s+k)$. In particular, by \cite[Corollary 6.4]{francone2025cluster},~the~set 
$$\{f_{1,a}^{(i)}\,|\,(i,a)\in I\times_2\Z\}$$
freely generates the ring of $G$-invariants ${}^G\C[\widehat{\mathcal{Z}}_{\infty}^{\circ}]\simeq {}^G \cR$.
\end{Rem}

Using Lemma \ref{lem:constant_function_1}, we see that $\widehat{\cZ}_\infty^\circ$ is actually a closed subscheme of $\prod_{(i,a)\in I\times_2\Z}V(\varpi_i)^*$. Indeed,
the open condition $x_{i,a}\ne 0$ in the definition of $\widehat{\cZ}_\infty^\circ$ (which is implicit in Definition \ref{def:Zhat} but appears in the definition~of~$\widehat{\op{Gr}}(i)$) can be replaced by the closed condition
\begin{equation*}
    \langle x_{i,a}, x_{i^*,a+h}\rangle = 1.
\end{equation*}
We saw that this equation is automatically satisfied for every point of $\widehat{\cZ}_\infty^\circ$, and conversely the equation above clearly implies nonvanishing of $x_{i,a}$. For easy reference, we summarize the presentation of $\widehat{\cZ}_\infty^\circ$ that follows from the discussion above as follows:

\begin{Proposition} \label{prop:presentation of Zhat}
    $\widehat{\cZ}_\infty^\circ$ is the closed subscheme of 
    $\prod_{(i,a)\in I\times_2\Z} V(\varpi_i)^*$
defined by the following equations for every $(i,a)\in I\times_2 \Z$:
    \begin{equation} \label{eq:Zhat-Gr-equation}
        x_{i,a}\otimes x_{i,a} \in V(2\varpi_i)^* \subseteq V(\varpi_i)^*\otimes V(\varpi_i)^*,
    \end{equation}
    \begin{equation} \label{eq:Zhat-incidence-equation}
        x_{i,a}\otimes x_{j,a+1} \in V(\varpi_i+\varpi_j)^* \subseteq V(\varpi_i)^*\otimes V(\varpi_j)^*
    \end{equation} 
    whenever $i\sim j$,
    \begin{equation} \label{eq:Zhat-wedge-tensor-equation}
       \textstyle x_{i,a} \wedge x_{i,a+2} = \bigotimes_{j\sim i} x_{j,a+1},
    \end{equation}
    and 
    \begin{equation} \label{eq:Zhat-pairing-equation}
        \langle x_{i,a}, x_{i^*,a+h}\rangle = 1.
    \end{equation}
\end{Proposition}

\begin{Corollary} \label{cor:presentation-of-R}
    The ring $\cR= \C[\widehat{\cZ}_\infty^\circ]$ is the quotient of the infinite polynomial ring
    \[ \bigotimes_{(i,a)\in I\times_2\Z} \op{Sym}^\bullet V(\varpi_i,a) \]
    (where $V(\varpi_i,a)$ is just a copy of $V(\varpi_i)$ indexed by an integer $a$) by the ideal generated by:\smallskip
    \begin{enumerate}
        \item\label{item:rel1} The kernel of $\op{Sym}^2V(\varpi_i,a) \twoheadrightarrow V(2\varpi_i)$ for every $(i,a)\in I\times_2\Z$,
        \item\label{item:rel2} The kernel of $V(\varpi_i,a) \otimes V(\varpi_j,a+1) \twoheadrightarrow V(\varpi_i+\varpi_j)$ for all $(i,a)\in I\times_2\Z$~and~$j\sim i$,
        \item\label{item:rel3} $\iota_1(f)-\iota_2(f)$ for every $(i,b)\in I\timesop_2\Z$ and $f\in V(2\varpi_i-\alpha_i)$, where 
		\begin{align*}
		\iota_1: V(2\varpi_i-\alpha_i)&\hookrightarrow  V(\varpi_i, b-1) \otimes V(\varpi_i, b+1),\\
		\iota_2: V(2\varpi_i-\alpha_i) &\hookrightarrow \otimes_{j\sim i} V(\varpi_j,b)
		\end{align*}		        
        are defined in \eqref{eq:Gequivariant_map_wedge} and \eqref{eq:Gequivariant_map_from_identity}, and
        \item\label{item:rel4} $\Omega_{i,a}-1$ for each $(i,a)\in I\times_2 \Z$, with $\Omega_{i,a}\in V(\varpi_i,a)\otimes V(\varpi_{i^*},a+h)$ the $G$-invariant element corresponding to the normalized bilinear pairing $\langle \cdot, \cdot\rangle$ between $V(\varpi_i)^*$ and $V(\varpi_{i^*})^*$ (i.e.~the generator of the trivial representation in $V(\varpi_i,a)\otimes V(\varpi_{i^*},a+h)$).
    \end{enumerate}
\end{Corollary}

\begin{Rem}
To make sense of \eqref{eq:Zhat-wedge-tensor-equation}, note that by Lemma \ref{lem:wedge-eq-tensor}, the equations \eqref{eq:Zhat-Gr-equation} and \eqref{eq:Zhat-incidence-equation} imply that both sides of \eqref{eq:Zhat-wedge-tensor-equation} lie in the isotypic component $ V(2\varpi_i - \alpha_i)^* $.  This also explains why condition (\ref{item:rel3}) in Corollary \ref{cor:presentation-of-R} is equivalent to \eqref{eq:Zhat-wedge-tensor-equation}.
\end{Rem}

\subsection{A filtration of $\cR$}\label{subsec:filt_of_R}

We now give a filtration of the graded ring $\cR$ indexed by sets of parameters $\bR\in \Z^\lambda$,
analogous to the filtration of $K_0(\O_{sh})$ by $K_0(\O^\lambda_{sh}(\bR))$.

\begin{Def} \label{def:filtration-of-R}
    For $\lambda\in P_+$ and $\bR\in \Z^\lambda$,
    let $\cL_{\lambda, \bR}$ be the line bundle on $\cZ_\infty$ defined by
    \[ \textstyle \cL_{\lambda, \bR} := \bigotimes_{i\in I} \bigotimes_{a\in R_i} \O_{i,a}(1). \]
    Note that $\cL_{\lambda, \bR}|_{\cZ_\infty^\circ}= \cL_{\tau_\bR}$, where $\tau_\bR=\awt(y_{\bR})\in \cP$. 
    Let 
    \begin{equation}\label{eq:VlambdaRAntoine}
     V(\lambda, \bR) := H^0(\cZ_\infty, \cL_{\lambda, \bR}) \subseteq H^0(\cZ_\infty^\circ, \cL_{\tau_\bR}) \subseteq \cR.
     \end{equation}
\end{Def}
Given two pairs $(\lambda_1, \bR_1)$, $(\lambda_2, \bR_2)$, we clearly have
\[ V(\lambda_1, \bR_1) \cdot V(\lambda_2, \bR_2) \subseteq V(\lambda_1 + \lambda_2, \bR_1 \cup \bR_2).\]
In fact, it turns out that the reverse inclusion also holds.

\begin{Lemma} \label{le:surjectR}
    The multiplication map
\begin{center}
$V(\lambda_1, \bR_1) \otimes V(\lambda_2, \bR_2) \to V(\lambda_1 + \lambda_2, \bR_1 \cup \bR_2)$
\end{center}
    is surjective.
\end{Lemma}

\begin{proof}
    It's enough to show that, for every pair $(\lambda, \bR)$, the multiplication map
    \[\textstyle \bigotimes_{i\in I} \bigotimes_{a\in R_i} V(\varpi_i,a) \to V(\lambda, \bR)\]
    is surjective. Since 
    \[ \smash{V(\lambda, \bR) = H^0(\cZ_\infty, \cL_{\lambda, \bR}) =\lim_{\to} H^0(\cZ_{\xi,\xi'}, \cL_{\lambda, \bR}),}\]
    it is enough to show the surjectivity of 
    \[\textstyle \bigotimes_{i\in I} \bigotimes_{a\in R_i} H^0(\cZ_{\xi,\xi'}, \O_{i,a}(1)) \to H^0(\cZ_{\xi,\xi'}, \cL_{\lambda, \bR})\]
    for $\xi, \xi'$ such that $\bR$ is supported between $\xi$ and $\xi'$, which follows from Corollary \ref{cor:GBS-mult-surj}.
\end{proof}

We now prove that $\{V(\lambda, \bR)\}_{\lambda, \bR}$ really is a filtration of $\cR$ as a graded ring, in the sense that every homogeneous element of $\cR$ lies in some $V(\lambda, \bR)$. To do so, we first need: %
\begin{Lemma} \label{lem:Oia-tensor-Oi*a+h}
    For every $(i,a)\in I\times_2 \Z$, the class of the line bundle $\O_{i,a}(1)\otimes \O_{i^*,a+h}(1)$ in $\op{Pic}(\cZ_\infty)$ is a nonnegative linear combination of the classes of the boundary divisors $D_{j,b}$, for $(i,a+1) \le (j,b) \le (i^*, a+h-1)$
\end{Lemma}
\begin{proof} 
    This follows easily from Lemma \ref{th:lowest} (by identifying $\Pic(\cZ_\infty)$ with $\cB$, which matches the classes of the boundary divisors $D_{j,b}$ to the generators $z_{i,a}$ of $\Gamma_+$), but we can also give a simple geometric proof as follows.\medskip\par
    
    The nondegenerate $G$-invariant pairing $\langle\cdot, \cdot\rangle : V(\varpi_i)^*\times V(\varpi_{i^*})^* \to \C$ induces a section of $\O_{i,a}(1)\otimes \O_{i^*,a+h}(1)$ on $\cZ_\infty$. We claim that $s$ is nonvanishing on $\cZ_\infty^\circ$. To see this, pick a height function $\xi$ such that $\xi_i=a$. Then $s$ can be thought of as a section of $\O_{i,a}(1)\otimes \O_{i^*,a+h}(1)$ on $\cZ_{\xi, \xi^*+h}$. But by Remark \ref{rem:Z-and-bands}, $\cZ_{\xi, \xi^*+h}^\circ$ can be identified with the big Schubert cell in $\Fl\times \Fl$, on which it's clear that the section of $\O(\varpi_i)\boxtimes \O(\varpi_{i^*})$ induced by the bilinear pairing $V(\varpi_i)^*\times V(\varpi_{i^*})^* \to \C$ is nonvanishing.\medskip\par

    If $m_{i,b}\in \N$ is the order of vanishing of $s$ at the boundary divisor $D_{i,b}$, then $s$ therefore induces an isomorphism between $\O_{i,a}(1)\otimes \O_{i^*,a+h}(1)$ and $\bigotimes_{i,b} \O(m_{i,b}D_{i,b})$.\medskip\par
    
    Taking the height function $\xi_j = a+d(i,j)$ shows that only boundary divisors indexed by $(j,b)\ge (i,a+1)$ can appear with a nonzero coefficient. Taking the height function $\xi_j = a-d(i,j)$ shows that only boundary divisors indexed by $(j,b)\le (i^*,a+h-1)$ can appear with a nonzero coefficient.
\end{proof}

\begin{Proposition} \label{th:cRgamma}
    For every $\tau \in \cP$, we have
    \[ \cR_\tau = \bigcup_{\substack{\lambda\in P_+,\bR\in \Z^\lambda\\\awt(y_{\bR}) = \tau}} V(\lambda, \bR),\]
    where $\cR_\tau=H^0(\cZ_\infty^\circ, \cL_{\tau})$ is the degree $\tau$ part of $\cR$.
\end{Proposition}

\begin{proof}
    Let $s\in R_\tau= H^0(\cZ_\infty^\circ, \cL_{\tau})$. Then $s\in H^0(\cZ_{\xi,\xi'}^\circ, \cL_{\tau})$ for some height functions $\xi, \xi'$. Pick a decomposition \[\tau=\sum_{i\in I}\sum_{\substack{\xi_i \le a \le \xi_i' \\ a\equiv_2 i}} m_{i,a}\tau_{i,a}\]
    of $\tau$ as a linear combination of classes $\tau_{i,a}$, which yields a line bundle $\cL_{\ul{m}}$ on $\cZ_{\xi,\xi'}$ extending the line bundle $\cL_{\tau}$ on $\cZ_{\xi,\xi'}^\circ$. Then $s$ can be thought of as a meromorphic section of $\cL_{\ul{m}}$. If $s$ has a pole of order $k$ at one of the boundary divisors $D_{i,b}$, we can get rid of this pole by increasing $m_{i,b-1}$ and $m_{i,b+1}$ by $k$ while decreasing $m_{j,b}$ by $k$ for $j\sim i$, which correspond to twisting $\cL_{\ul{m}}$ by $\O(kD_{i,b})$. We can therefore assume, after changing $\ul{m}$, that $s$ has no pole at any of the boundary divisors, and is therefore a regular section of $\cL_{\ul{m}}$.\medskip\par

    If $m_{i,a}\ge 0$ for all $i,a$, then we are done: we have $\cL_{\ul{m}}=\cL_{\lambda, \bR}$ for some $(\lambda, \bR)$, and then $s\in V(\lambda, \bR)$ since it extends to a regular section of $\cL_{\lambda, \bR}$ on $\cZ_{\infty}$. If $m_{i,a}=-k<0$ for some $i,a$, then let $\ul{m}'$ be obtained from $\ul{m}$ by increasing both $m_{i,a}$ and $m_{i^*,a+h}$ by $k$ (which might require replacing $\xi'$ by $\xi'^*+h$ so that $(i^*,a+h)$ is still in the correct range). By Lemma \ref{lem:Oia-tensor-Oi*a+h}, $\cL_{\ul{m}'}$ is also an extension of $\cL_{\tau}$ and we have an inclusion $\cL_{\ul{m}}\subseteq\cL_{\ul{m}'}$ as extensions of $\cL_{\tau}$ to $\cZ_{\xi,\xi'}$, so $s$ is still a regular section on $\cL_{\ul{m}'}$. We can repeat this to get rid of all negative entries in $\ul{m}$, and then we are done.
\end{proof}

The last result of this section is a description of $V(\lambda, \bR)$ as a representation of $G$.

\begin{Theorem} \label{th:charcR}
    The character of $V(\lambda, \bR)$, as a left $G$-module, is equal to the character of the product monomial crystal $\cB(\lambda, \bR)$.
\end{Theorem}
\begin{proof}
    To prove this, we will compare the Demazure character formula for Bott--Samelson varieties given in Theorem \ref{thm:demazure-character-for-BS} with the Demazure character formula for product monomial crystals from \cite{gibson2021demazure}. Fix two height functions $\xi, \xi'$ such that $R_i$ is supported on the interval $(\xi_i,\xi_i']$ for all $i$. Let $(i_1, b_1), (i_2,b_2), \ldots, (i_r, b_r)$ be an enumeration of the pairs 
    $(i,b)\in I\timesop_2\Z$  with $\xi_i< b < \xi_i'$,
    compatible with the partial order on 
    $I\timesop_2\Z$. 
    Then $H(\xi, \xi')=H(i_1, \ldots, i_r)$, so by Theorems \ref{thm:demazure-character-for-BS} and \ref{thm:BS-to-GBS-global-sections-iso} we have
    \begin{align}
        \op{ch}(V(\lambda, \bR))
        &= \op{ch}H^0(\cZ_{\infty}, \cL_{\lambda, \bR}) = \op{ch}H^0(\cZ_{\xi,\xi'}, \cL_{\lambda, \bR}) \notag\\
        &=  \Lambda_{w_0}(\Lambda_{i_1}(e^{m_{i_1,b_1+1}\varpi_{i_1}} \Lambda_{i_2}(e^{m_{i_2,b_2+1}\varpi_{i_2}} \cdots \Lambda_{i_r}(e^{m_{i_r,b_r+1}\varpi_{i_r}}) \cdots ))), \label{eq:demazure-character-for-RlambdaR}
    \end{align} 
(where $m_{i,a}$ is the multiplicity of $a$ in the multiset $R_i$) and the second equality follows from Remark \ref{rem:growth-of-pic-and-sections}.

    We now compute the character of $\cB(\lambda, \bR)$ via the method given in \cite{gibson2021demazure}. This method involves auxiliary truncations $\cB(\lambda, \bR, J)\subseteq \cB(\lambda, \bR)$ depending on an upward-closed set $J\subseteq I\times_2 \Z$, consisting of the monomials in $\cB(\lambda, \bR)$ whose support lies inside $J$. Define a sequence of triples $(\lambda_0, \bR_0, J_0), \ldots, (\lambda_r, \bR_r, J_r)$ by the following rules:
    \begin{itemize}
        \item $\lambda_r=0$, $\bR_r=\varnothing$, and $J_r$ is the upward-closed set with boundary $\xi'$, i.e.~
        $$J_r = \{(i,a)\in I\times_2 \Z : a\ge \xi'_i\}.$$
        \item For $k\le r$, $\lambda_{k-1}=\lambda_{k}+m_{i_k,b_k+1} \varpi_{i_k}$, $\bR_{k-1}$ is obtained from $\bR_{k}$ by adjoining $b_k+1$ to $R_{i_k}$ with multiplicity $m_{i_k,b_k+1}$, and $J_k = J_{k+1} \sqcup \{(i_k, b_k-1)\}$.
    \end{itemize}
    Note that $(\lambda_0, \bR_0)=(\lambda, \bR)$ and $\bR_k$ is supported in $J_k$ for all $k$. By \cite[Theorem 5.9]{gibson2021demazure},%
    \[\op{ch} \cB(\lambda_r, \bR_r, J_r) = 1\]
    and
    \[ \op{ch} \cB(\lambda_{k-1}, \bR_{k-1}, J_{k-1}) = \Lambda_{i_k}(e^{m_{i_k,b_k+1} \varpi_{i_k}} \op{ch} \cB(\lambda_{k}, \bR_{k}, J_{k}))\]
    for $k\le r$. Therefore, we have
    \[ \op{ch} \cB(\lambda_0, \bR_0, J_0) = \Lambda_{i_1}(e^{m_{i_1,b_1+1}\varpi_{i_1}} \Lambda_{i_2}(e^{m_{i_2,b_2+1}\varpi_{i_2}} \cdots \Lambda_{i_r}(e^{m_{i_r,b_r+1}\varpi_{i_r}}) \cdots )).\]
    Moreover, by \cite[Corollary 5.14]{gibson2021demazure}, we have
    \[ \op{ch} \cB(\lambda, \bR) = \op{ch} \cB(\lambda_0, \bR_0) = \Lambda_{w_0} (\op{ch} \cB(\lambda_0, \bR_0, J_0)).\]
    Comparing with (\ref{eq:demazure-character-for-RlambdaR}) yields the desired equality of characters. 
\end{proof}

We are now in a position to fulfil a promise made in the proof of Theorem \ref{thm:charaterization_max_sing_crystals}.

\begin{Corollary}
Fix $\lambda\in P_+$ together with a height function $\xi$. Let $\bR\in \Z^\lambda$ be the set of parameters given by $\xi$. Then, $\bR$ is maximally singular. 
\end{Corollary}

\begin{proof}
As $\cZ_{\xi,\xi}\simeq \cZ_{\xi,\xi}^\circ\simeq B_- \backslash G$, the representation $V(\lambda,\bR)$ coincides with the representation $V(\lambda)$ by the Borel-Weil Theorem. Since the previous theorem shows that $V(\lambda,\bR)$ has the same character as $\cB(\lambda,\bR)$, the statement follows.
\end{proof}

\begin{Rem}\label{rem:generalized_schur_modules}
Take $ G = \SL_m$ %
.  In \cite[Theorem 6.23]{gibson2021demazure}, Gibson shows that $ \cB(\la, \bR) $ is the crystal of a generalized Schur module.  On the other hand, in \cite{magyar1998borel}, Magyar gives a Borel-Weil construction for the generalized Schur modules.  Thus, Theorem \ref{th:charcR} can be regarded as a generalization of Magyar's result from type A to any simply-laced group.
\end{Rem}

\section{Isomorphism between \texorpdfstring{$K_\C(\O_{sh})$}{KC(Osh)} and \texorpdfstring{$\cR$}{R}}\label{sec:iso_between_KOsh_and_R}
In this section, we prove \cite[Conjecture 6.5]{hernandez2025jordan} by constructing an algebra isomorphism between $\cR$ and $K_{\C}(\mathcal{O}_{sh})$. We then use this isomorphism~to describe particular subalgebras of $\mathcal{R}$ and deduce an alternative presentation for $K_{\C}(\mathcal{O}_{sh})$ in type A.
\subsection{The morphism}
By Corollary \ref{cor:presentation-of-R}, the algebra $\cR$ is generated by the %
fundamental modules $V(\varpi_i,a)$ for $ (i,a) \in I \times_2 \Z $, with %
explicit relations. In particular, for $(i,a) \in I\times_2\Z$ and $\gamma=w\varpi_i\in W\varpi_i$%
, we have a generalized minor $\Delta_{\gamma,a}\in \cR$ which evaluates on %
$x\in \widehat{\cZ}_\infty^\circ$~as
\begin{equation}\label{eq:Deltagammaa}
\Delta_{\gamma,a}(x)%
=\langle x_{i,a},v_{\gamma}\rangle
,
\end{equation}
with $\langle.,.\rangle:V(\varpi_i)^\ast \otimes V(\varpi_i)\to \C$ the evaluation pairing, $v_{\gamma}:=\dot{w} v_{\varpi_i}\in V(\varpi_i)$, and $x_{i,a}$ the $(i,a)$-component of $x$%
.
\begin{Rem}\label{rem:what_are_generalized_minors}
Take a height function $\xi$ and let $c\in W$ be the associated Coxeter element. For $i\in I$, write %
 $m_i=\tfrac{1}{2}(\xi_i^*-\xi_i+h)\in \Z$ (see Section \ref{subsubsec:height_functions}).
The identification $G\simeq \smash{\widehat{\cZ}^\circ_{\xi,\xi^*+h}}$ associated~to~$\xi$ (see Remark \ref{rem:embeddings_of_CG}) is given by 
$\smash{g\mapsto x=(x_{i,\xi_i+2s})_{i\in I,0\leq s\leq m_i}}%
$ with
$$x_{i,\xi_i+2s}=v_{\varpi_i}^\ast \dot{c}^{-s}g\in V(\varpi_i)^\ast.$$
In addition, if $x\in\smash{\widehat{\cZ}_\infty^\circ}$ is sent to $g$ by the composition $\smash{\widehat{\cZ}_\infty^\circ}\twoheadrightarrow \smash{\widehat{\cZ}^\circ_{\xi,\xi^*+h}}\simeq G$, then, for $i\in I$, $0\leq s\leq m_i$ and $\gamma=w\varpi_i\in W\varpi_i$, we have\footnote{Note that $v_{\varpi_i}^*\dot{c}^{-s}=v_{\varpi_i}^*\dot{(c^{-s})}$ for $0\leq s\leq m_i$ because of \cite[Lemma 5.9]{francone2025cluster}.
}
\begin{equation*}
\Delta_{\gamma,\xi_i+2s}(x)=\langle x_{i,\xi_i+2s},v_{\gamma}\rangle=\langle v_{\varpi_i}^\ast \dot{c}^{-s}g,v_{w\varpi_i}\rangle
=\Delta_{c^s\varpi_i,w\varpi_i}(g)
\end{equation*}
according to the notation of Section \ref{subsec:base_affine_space_and_Osh}. This explains why we call the functions $\Delta_{\gamma,a}\in \cR$ ``generalized minors''. 
\end{Rem}

Recall the identification $V(\varpi_i,a)\simeq K_\C(\O^{\varpi_i}_{sh}(a))$ of Section \ref{sec:glueing_of_g_action}, which takes $v_{\gamma}\in V(\varpi_i)$~(or, equivalently, the generalized minor $\Delta_{\gamma,a}$) to the class of the chamber module  $L_{\gamma,a}$ (see the paragraph following Remark \ref{rem:base_affine_space_map_factors}). We use this discussion to prove the following result:

\begin{Theorem}\label{thm:Omega_is_surj}
Identifying $V(\varpi_i,a)$ and $ K_\C(\O^{\varpi_i}_{sh}(a))$ for all $(i,a)\in I\times_2\Z$ extends to a surjective $(G\times A)$-equivariant morphism of $\C$-algebras
\begin{equation*}
\Omega:\cR \rightarrow K_\C(\O_{sh}).
\end{equation*}  
\end{Theorem}
\begin{proof}
To show that $\Omega$ is an algebra morphism, we must check the four relations of Corollary \ref{cor:presentation-of-R}.  Relations \eqref{item:rel1}--\eqref{item:rel2} follow from Corollary \ref{co:Rel12} whereas \eqref{item:rel3} follows from Theorem \ref{thm:the_heptagon_commutes}. For \eqref{item:rel4}, let $ (i,a) \in I \times_2 \Z $. For $\bR=(a)_i\cup(a+h)_{i^\ast}$, \eqref{eq:monomial_one_in_crystal} gives
\begin{equation*}
1\in \cB(\varpi_i + \varpi_{i^*},\bR)
\end{equation*}
and the trivial representation $L(\mathbbm{1})$ is thus an object of $\O^{\varpi_i + \varpi_{i^*}}_{sh}(\bR)$. Hence, 
$$1=[L(\mathbbm{1})]\in K_0(\O^{\varpi_i + \varpi_{i^*}}_{sh}(\bR)),$$ 
and \eqref{item:rel4} easily follows%
. Surjectivity also follows from Theorem \ref{thm:generators_of_KOsh} and equivariance is clear since $\Omega$ is equivariant on generators (see, e.g., Remark \ref{rem:AequivBaseAff} and \eqref{eq:VlambdaRAntoine}).
\end{proof}

We will now give two separate proofs of the injectivity of $\Omega$.

\subsection{Cluster algebraic proof of the injectivity}
Fix a height function $\xi$ with Coxeter element $c=c_\xi$ and let $m_i=\frac{1}{2}(\xi_i^*-\xi_i+h)$ for $i\in I$%
. By \cite[Theorem 1.3]{francone2025cluster},~%
the~complexification~of the cluster algebra of \cite{geiss2024representations} is isomorphic to the coordinate ring $\C[B(G,c)]$ of the scheme of $(G,c)$-bands %
of Remark \ref{rem:Z-and-bands}.
Moreover, using Francone--Leclerc's notation, the initial~seed (associated to $\xi$) of Geiss--Hernandez--Leclerc's cluster algebra consists entirely of functions  
$\{u_{i,r}\}_{(i,r)\in I\times_2\Z}$ of 
 the form \cite[Definition 5.1]{francone2025cluster}
$$ u_{i,r}={\Delta^{(s)}_{c^{k}\varpi_i,\tilde{c}^{\ell}\varpi_i}}$$
for some $s,\ell\in \Z$ with $0\leq k\leq m_i$, where $\tilde{c}=w_0c^{-1}w_0$, and where $\smash{\Delta_{c^k\varpi_i,\tilde{c}^{\ell}\varpi_i}^{(s)}}$ is~the~generalized minor $\Delta_{c^k\varpi_i,\tilde{c}^{\ell}\varpi_i}$ for the $\smash{s^{\text{th}}}$-copy of $\C[G]$ in $\C[B(G,c)]$. Also,
$${\Delta^{(s)}_{c^{k}\varpi_i,\tilde{c}^{\ell}\varpi_i}}={\Delta^{(s+k)}_{\varpi_i,\tilde{c}^{\ell}\varpi_i}}$$
for $s,\ell\in \Z$ with $0\leq k\leq m_i$ by \cite[Proposition~5.11]{francone2025cluster}, and it hence follows from Remark \ref{rem:what_are_generalized_minors} that the cluster variables in the above initial seed can all be associated to generalized minors of the form \eqref{eq:Deltagammaa} via the isomorphism
$$ \C[B(G,c)]\simeq \C[\smash{\widehat{\mathcal{Z}}_{\infty}^{\circ}}]\simeq \cR$$
of Remark \ref{rem:Zhat-and-bands}. This can be summarized as follows:

\begin{Lemma}\label{lem:initialSeed} There exist functions $\{\gamma_i:\overline{i}+2\Z\to W\varpi_i\}_{i\in I}$ and $\{a_i:\overline{i}+2\Z\to \overline{i}+2\Z\}_{i\in I}$ such that  the initial seed for the above cluster algebra structure on $\cR$ is
$$ \{ \Delta_{\gamma_i(r), a_i(r)}\,|\, (i,r) \in I \times_2 \Z \}.$$ 
 \end{Lemma}
In particular, the Laurent phenomenon \cite[Theorem 3.1]{fomin2002cluster} gives:
\begin{Corollary}
The algebra $\cR$ is contained in $\C\otimes_{\Z}\mathscr{L}$ with 
$$ \mathscr{L}:=\Z[\Delta_{\gamma_i(r), a_i(r)}^{\pm 1}\,|\, (i,r) \in I \times_2 \Z ].$$
\end{Corollary}
 
Consider the ring morphism $\mathscr{L}\to \cE_{\ell}$ given by 
$${\Delta_{\gamma_i(r), a_i(r)}}\mapsto{Q_{\gamma_i(r),a_i(r)}}$$ 
where $Q_{\gamma_i(r),a_i(r)}$ is the $Q$-variable of \eqref{eq:Qvariables} (which is invertible by \eqref{eq:factorisationQvars}). Then this morphism is injective by \cite[Proposition 8.1]{geiss2024representations} and we are thus led to the diagram
\begin{equation} \label{eq:diamond}
\adjustbox{scale=0.9}{
\begin{tikzcd}
\cR \ar[d,"\Omega"',two heads] \ar[r,hook]  &  \C\otimes_{\Z}\mathscr{L}  \ar[d,hook] \\
K_\C(\O_{sh}) \ar[r,"\chi_\ell",hook] & \C\otimes_\Z\cE_\ell
\end{tikzcd}}
\end{equation}
We are now in a position to prove:

\begin{Theorem}\label{thm:omega_is_inj_cluster_side}
The map $ \Omega :\cR \twoheadrightarrow K_\C(\O_{sh})$ is injective.
\end{Theorem}
\begin{proof}
Fix $(i,a)\in I\times_2\Z$ with $\gamma\in W\varpi_i$. Then the definition of $\Omega$ and Theorem \ref{thm:ConjFH} give
$$ (\chi_{\ell}\circ \Omega)(\Delta_{\gamma,a}) =  \chi_{\ell}(L_{\gamma,a})=Q_{\gamma,a},$$
showing that the diagram \eqref{eq:diamond} commutes (at least) for the initial cluster variables given in Lemma \ref{lem:initialSeed}. However, since the top horizontal arrow is an injection of $\cR$ in the ring~$\C\otimes_{\Z}\mathscr{L}$ of Laurent polynomials in these initial cluster variables, it follows that \eqref{eq:diamond} commutes for any element of $\cR$. Thus, the theorem follows from the fact that both the upper horizontal and the right vertical arrows are injective.
\end{proof}
The above proof %
also allows us to establish \cite[Conjecture 6.5]{hernandez2025jordan}. Indeed, denote~by~$\AGHL$ Geiss--Hernandez--Leclerc's cluster algebra (over $\Z$) and let $\Oshtot$ be the category of Definition \ref{def:intCatO}, i.e.~the integral category $\mathcal{O}_{sh}$ considered since Section \ref{sec:KLR}, but for which the objects are not necessarily of finite-length. Clearly, the ring $K_0({\Oshtot})$ can be thought of as a completion of $K_0(\mathcal{O}_{sh})$ where the basis of classes of simple objects becomes a ``topological basis'' (see, e.g.,~\cite{hernandez2024shifted}). Finally, let $\AGHLcomp$ be the completion of the cluster algebra $\AGHL$ studied in \cite{geiss2024representations}, i.e.~$\AGHLcomp$ is the closure of $\AGHL$ in $\cE_{\ell}$ under the embedding $\AGHL\hookrightarrow \mathscr{L}\hookrightarrow \cE_{\ell}$.
\begin{Theorem}\label{thm:conjecture_HZ} The injective morphism $\AGHL\hookrightarrow \cE_{\ell}$ induces an isomorphism $\AGHL\simeq K_0(\mathcal{O}_{sh})$. Thus, \cite[Conjecture 6.5]{hernandez2025jordan} holds.
\end{Theorem} 
We will need the result below, which easily follows from the definition of $\AGHLcomp$ and the fact that simple classes give a $\C$-basis (resp.~a topological $\Z$-basis) for $K_{\C}(\mathcal{O}_{sh})$ (resp.~$K_0(\Oshtot)$).
\begin{Lemma}\label{lem:IntersectionsHZ} With the above notation, 
$$(\C\otimes_{\Z}\AGHL)\cap \AGHLcomp=\AGHL\ \text{ and }\ K_{\C}(\mathcal{O}_{sh})\cap K_0(\Oshtot)=K_0(\mathcal{O}_{sh}).$$
\end{Lemma}
\begin{proof}[Proof of Theorem \ref{thm:conjecture_HZ}] Unpacking what was done in this section, we see that the diagram
\begin{equation}\label{eq:diag1HZ}
\adjustbox{scale=0.88}{\begin{tikzcd}
\C\otimes_{\Z}\AGHL\ar[r,hook]\ar[d,"\simeq"] & \C\otimes_{\Z}\mathscr{L}\ar[r,hook] & \mathbb{C}\otimes_{\mathbb{Z}} \cE_{\ell}\\ \C[B(G,c)] \ar[r,"\simeq",swap] & \cR\ar[u,hook]\ar[r,"\simeq"',"\Omega"] & K_{\C}(\mathcal{O}_{sh})\ar[u,"\chi_{\ell}", hook]\
\end{tikzcd}}
\end{equation}
commutes. Moreover, by \cite[Proposition 9.14 and Theorem 9.15]{geiss2024representations}, the $\ell$-character~map $\chi_{\ell}:K_0(\mathcal{O}_{sh})\to \cE_{\ell}$ and the injective morphism $\AGHL\hookrightarrow \cE_{\ell}$ extend to ring isomorphisms
$$K_0(\Oshtot)\simeq \cE_{\ell} \text{ and }\AGHLcomp\simeq \cE_{\ell},$$
which induces an algebra isomorphism $\C\otimes_{\Z}\AGHLcomp\simeq K_{\C}(\Oshtot)$%
. Denote by $\phi$ this isomorphism and write $\tilde{\Omega}$ for the isomorphism $\C\otimes_{\Z}\AGHL\simeq K_{\C}(\mathcal{O}_{sh})$ appearing in \eqref{eq:diag1HZ}. Using~the~commutativity of \eqref{eq:diag1HZ}, it is easy to show that the diagram 
$$ 
\adjustbox{scale=0.88}{\begin{tikzcd}[column sep = 1.5em, row sep = 1.25em] 
\C\otimes_{\Z}\AGHL\arrow[rr,"\simeq"',"\tilde{\Omega}"]\ar[ddr,out=180, in=180, looseness=2, hook']\arrow[d,hook] &  & K_{\C}(\mathcal{O}_{sh})\ar[d,hook]\ar[ddl,out=0,in=0, looseness=2, hook]\\
\C\otimes_{\Z} \AGHLcomp\arrow[rr,"\simeq"',"\phi"]\ar[dr,"\simeq"] & & K_{\C}(\Oshtot)\ar[dl,"\simeq"]\\
& \C\otimes_{\Z}\cE_{\ell} &
\end{tikzcd}
}
$$
commutes. In particular, $\phi(\C\otimes_{\Z}\AGHL)=K_{\C}(\mathcal{O}_{sh})$, and, by construction, $\phi(\AGHLcomp)=K_{0}(\Oshtot)$. Hence, by Lemma \ref{lem:IntersectionsHZ}, the isomorphism $\phi$ restricts to an isomorphism 
$$\AGHL= (\C\otimes_{\Z}\AGHL)\cap \AGHLcomp \simeq K_{\C}(\mathcal{O}_{sh})\cap K_0(\Oshtot)=K_0(\mathcal{O}_{sh}),$$
which makes (again by construction) the diagram 
$$
\adjustbox{scale=0.88}{\begin{tikzcd}[row sep = 0.8em, column sep = 4em]
\AGHL\ar[dd,hook]\ar[dr,hook]\ar[rr,"\simeq"] & & K_0(\mathcal{O}_{sh})\ar[dd,hook]\ar[dl,hook',"\chi_{\ell}",swap]\\
& \cE_{\ell} & \\
\AGHLcomp\ar[rr,"\simeq"]\ar[ur,"\simeq"]& & K_0(\Oshtot)\ar[ul,"\simeq",swap]
\end{tikzcd}
}
$$
commutative. 
\end{proof}

\subsection{Geometric proof of the injectivity}
Fix $\la\in P_+$ with $\bR\in \Z^{\la}$ and consider 
$$V(\lambda, \bR) = H^0(\cZ_\infty, \cL_{\lambda, \bR})$$
where we follow the notation of Section \ref{subsec:filt_of_R}. Also, let $\Omega^{\lambda, \bR} : V(\lambda, \bR) \rightarrow  K_{\C}(\O_{sh})$
be the restriction of $\Omega:\cR\to K_{\C}(\mathcal{O}_{sh})$ to $V(\lambda, \bR)$. 

\begin{Theorem} \label{thm:injectivity-geo}
For all $\lambda$ and $\bR$, the map $\Omega^{\lambda,\bR}$ gives a $G$-equivariant isomorphism 
\begin{equation*}
V(\la,\bR)\cong K_\C(\O^\lambda_{sh}(\bR)). 
\end{equation*}
Consequently, $\Omega$ is injective.
\end{Theorem}
\begin{proof}
We begin by proving the first statement. For this, note that the elements of $ V(\lambda, \bR) $ can all be realized as products of elements from fundamental representations (as shown in the proof of Lemma \ref{le:surjectR}). Thus $\Omega$ indeed maps $ V(\lambda, \bR)$ to $ K_\C(\O^\lambda_{sh}(\bR))$. Also, $\Omega^{\la,\bR}$~is $G$-equivariant and surjective by Theorems \ref{prop:mult_is_equiv_truncation} and \ref{prop:mult_is_surj_truncation}, and, finally, Theorems \ref{thm:parity_KLRW_categorifies_VR}~and~\ref{th:charcR} show that $V(\la,\bR)$ has the same character %
as $K_{\C}(\O^\lambda_{sh}(\bR))$ as a left $G$-module. Hence $ \Omega^{\lambda, \bR} $ is indeed an isomorphism.

Now, we show that $ \Omega$ is injective.  Since $ \Omega $ is $ \aT $-equivariant, it suffices to check injectivity on each $ \aT $-weight space.  However, for $\tau\in \cP$,  Proposition \ref{th:cRgamma} gives
\begin{equation*}
\cR_\tau = \textstyle\bigcup_{\awt(y_{\bR})=\tau} V(\lambda, \bR),
\end{equation*}
and the injectivity of each $ \Omega^{\la, \bR} $ implies the injectivity of $ \Omega $ restricted to $ \cR_\tau$.  
\end{proof}

\subsection{Finite-dimensional modules}

Recall from Corollary \ref{cor:Invariants} the categories $ \scrC_{sh} $~and~$ \scrC_0 $ of finite-dimensional modules for shifted and unshifted Yangians.  Combining Theorem \ref{thm:Omega_is_surj} with this corollary, we obtain:

\begin{Corollary} \label{cor:Omega_restrict_NG}
The isomorphism $ \Omega$ restricts to algebra isomorphisms
$$
{}^N \cR \cong K_{\C}(\mathscr{C}_{sh}) \ \text { and }\ {}^G \cR \cong K_{\C}(\scrC_0)
$$
\end{Corollary}\vspace*{-1mm}
These isomorphisms are consistent with the two isomorphisms 
$$ {}^N\C[B(G,c)]\simeq K_{\C}(\scrC_{sh}) \text{ and } {}^G\C[B(G,c)]\simeq K_{\C}(\scrC_{0})$$
obtained in \cite[Section 8]{francone2025cluster}%
. The following result shows that our isomorphism is the same as the one given in \cite[Prop 8.1]{francone2025cluster}.  Recall the functions $\smash{f_{k,a}^{(i)}} \in {}^G \cR $ from \eqref{eq:def_G_invariant_functions}.

\begin{Corollary} 
For $ k \ge 1$ and $ (i,a) \in I \times_2 \Z$, the isomorphism $\Omega$ maps $f_{k,a}^{(i)} $ to the class $[\smash{\mathsf{W}_{k,a}^{(i)}}]\in K_0(\scrC_0)$ of the KR-module of Example \ref{ex:KR}.
\end{Corollary}

\begin{proof}
By definition, we have that $f^{(i)}_{k,a} \in {}^GV(\varpi_i + \varpi_{i^*},\bR)$ with $\bR = (a-2k)_i \cup (a + h)_{i^*}$. Hence, by Theorem \ref{thm:injectivity-geo} and the proof of Corollary \ref{cor:Invariants}, %
$\smash{\Omega(f_{k,a}^{(i)})}$ lies in 
$${}^GK_{\C}(\mathcal{O}_{sh}^{\varpi_i+\varpi_{i^*}}(\bR))= K_{\C}(\mathscr{C}^{\varpi_i+\varpi_{i^*}}_0(\bR)),$$
with $\smash{\mathscr{C}^{\varpi_i+\varpi_{i^*}}_0(\bR)\subseteq \mathcal{O}_0^{\varpi_i+\varpi_{i^*}}(\bR)}$ the full subcategory of finite-dimensional objects. Also, 
$$ \mathsf{Y}_{i,a-2(k-1)}\dots \mathsf{Y}_{i,a+2}\mathsf{Y}_{i,a}=\tfrac{\sfPsi_{i,a-2k}}{\sfPsi_{i,a}}\in \cB(\varpi_i+\varpi_{i^*},\bR)_0,$$
and thus $\smash{\mathsf{W}_{k,a}^{(i)}=L(\tfrac{\sfPsi_{i,a-2k}}{\sfPsi_{i,a}})}$ lies in $\mathscr{C}_0^{\varpi_i+\varpi_{i^*}}(\bR)$. On the other hand, 
$$ \dim K_{\C}(\mathscr{C}_0^{\varpi_i+\varpi_{i^*}}(\bR))=\dim \,{}^GV(\varpi_i+\varpi_{i^*},\bR)\leq \dim\,{}^G(V(\varpi_i)\otimes V(\varpi_{i^*})) = 1,$$
and therefore 
$\Omega(f_{k,a}^{(i)})=m{[\mathsf{W}_{k,a}^{(i)}]}$
for some $m\in \C$. Now, by %
Corollary \ref{cor:tensor_product_decomp_in_Osh},
$$ \textstyle [L(\sfPsi_{i,a-2k})][L(\sfPsi_{i,a}^{-1})]=[\mathsf{W}_{k,a}^{(i)}]+\sum_{\zeta\prec  %
\mathsf{Y}_{i,a-2(k-1)\dots \mathsf{Y}_{i,a+2}\mathsf{Y}_{i,a}}}n_{\zeta}[L(\zeta)]$$
for some $n_{\zeta}$'s in $\mathbb{Z}_{\geq 0}$, but, by definition of the pairings $\langle.,.\rangle:V(\varpi_i)^*\otimes V(\varpi_{i^*})^*\to \C$ given in Section \ref{sec:bi-infinite-bott-samelson}, the composition
$$ V(\varpi_i,a-2k)\otimes V(\varpi_{i^*},a+h)\xlongrightarrow{\text{mult}} V(\varpi_i+\varpi_{i^*},\bR)
\twoheadrightarrow {}^GV(\varpi_i+\varpi_{i^*},\bR)$$
sends $\Delta_{\varpi_i,a-2k}\otimes \Delta_{w_0\varpi_{i^*},a+h}$ to $\smash{f_{k,a}^{(i)}}$. This shows that $m=1$ since $\Delta_{\varpi_i,a-2k}$ and $\Delta_{w_0\varpi_{i^*},a+h}$ respectively correspond to $[L(\sfPsi_{i,a-2k})]$ and $[L(\sfPsi_{i,a}^{-1})]$ via the isomorphism $\Omega$.
\end{proof}
\subsection{Application to global sections of line bundles on $\cZ_\infty$} \label{sec:application-global-sections}
Recall that the Picard group of $\cZ_\infty$ is the free abelian group on the classes of the line bundles $\{\O_{i,a}(1)\}_{(i,a)\in I\times_2\Z}$. 
This group can thus be identified with the crystal $\cB$, by associating (as in Definition \ref{def:filtration-of-R}), the monomial $\textstyle m=\prod_{i,a} y_{i,a}^{m_{i,a}}\in \cB$
to the line bundle 
$$\textstyle  \smash{\cL_m := \bigotimes_{i,a} \O_{i,a}(m_{i,a}).}$$
One of the key ingredients in the geometric proof of the injectivity of %
$\Omega: \cR \overset{\sim}{\to} K_{\C}(\O_{sh})$ was a good understanding of the space of global sections of globally generated line bundles on $\cZ_\infty$ (i.e.~the $\cL_m$'s with $m\in \cB_+$), in the form of the Demazure character formula~(Theorem \ref{thm:demazure-character-for-BS}) or its product monomial crystal interpretation (Theorem \ref{th:charcR}). In this section, we use our isomorphism $\Omega$ and the combinatorics of product monomial crystals to deduce~a description of $H^0(\cZ_\infty, \cL_{m})$ for any Laurent monomial $m\in \cB$.\medskip\par
Noticing that $\cL_m|_{\cZ_\infty^\circ}= \cL_{\awt(m)}$, we can generalize Definition \ref{def:filtration-of-R} and define
\[ \smash{V(m) := H^0(\cZ_\infty, \cL_m) \subseteq H^0(\cZ_\infty^\circ, \cL_{\awt(m)}) \subseteq \cR.}\]
By construction, $V(m)=V(\la,\bR)$
if $m=y_\bR\in \cB_+
$ for $\la\in P_+$ and $\bR\in \Z^{\la}$.

It will be convenient to consider the partial order $\leq$ on $\cB$ given by $m \leq m'$ if $m' \in m\Gamma_+$ (note that the restriction of this order to $\cB_+$ is exactly relation (4) of Theorem \ref{th:TFAE}, %
 which gives rise to the order on sets of parameters considered in Section \ref{sec:glueing_of_g_action}). Also, $m'\leq m$ implies $V(m') \subseteq V(m)$ since then $\cL_{m'} = \cL_{m}(-D)$, where $D$ is a non-negative linear combination of boundary divisors $D_{i,a}$ (this was already used implicitly in the proof of Proposition \ref{th:cRgamma}).

It turns out that we can reduce the %
study of the $V(m)$'s to the case where $m\in \cB_+$ (for which we already have a fairly good understanding thanks to Theorems \ref{th:charcR} and \ref{thm:injectivity-geo}):
\begin{Theorem} \label{thm:global-sections-negative-line-bundles}
 For $m\in \cB$, %
 $V(m)$ is spanned by %
 the
 $V(m')$~with $m'\in \cB_+$ and $m' \leq m$.
\end{Theorem}
\begin{proof}
Fix $m\in \cB$ and let
\begin{equation*}
V'(m) = \sum_{\substack{m' \in \cB_+\\ \,m' \leq m}}V(m')\hspace{1em}\text{and}\hspace{1em}V''(m) = \bigcap_{\substack{m'' \in \cB_+ \\ m'' \geq m}} V(m'').
\end{equation*}
Then %
$V'(m) \subseteq V(m) \subseteq V''(m)$, so it suffices to prove that  the inclusion $V'(m)\subseteq V''(m)$ is an equality. For this, define
\begin{equation*}
\cB'(m) = \bigcup_{\substack{\lambda\in P_+, \bR\in \Z^\lambda \\ y_\bR \leq m}} \cB(\lambda, \bR)  \hspace{1em}\text{and}\hspace{1em}\cB''(m) = \bigcap_{\substack{\lambda\in P_+, \bR\in \Z^\lambda \\ y_\bR \geq m}} \cB(\lambda, \bR).
\end{equation*}
By Theorems \ref{th:descend} and \ref{thm:injectivity-geo} (with the identification $\Omega: \cR \overset{\sim}{\to} K_{\C}(\O_{sh})$), $V'(m)$ and $V''(m)$ have bases given by the classes of simple objects of $\O_{sh}$ with highest $\ell$-weight in $\cB'(m)$ and  $\cB''(m)$, respectively. We are therefore reduced to showing that the inclusion $\cB'(m)\subseteq \cB''(m)$ (coming combinatorially from Theorem \ref{th:TFAE}) is an equality. \medskip\par %
Fix $p\in \cB''(m)$. Then Corollary \ref{cor:minimal-pmc} gives $(\la_p,\bR_p)$ such that $p\in \cB(\la_p,\bR_p)\subseteq \cB''(m)$,~and hence $y_{\bR_p}\in \cB''(m)\cap \cB_+$. In particular, $y_{\bR_p} \in m\Gamma$ and it is enough to prove that $y_{\bR_p}\leq m$, as this would imply $p\in \cB(\la_p,\bR_p)\subseteq \cB'(m)$.\medskip\par 
Suppose $y_{\bR_p}\not\leq m$. Then there exists $(i,a)\in I\times_2\Z$ such that %
$y_{\bR_p}m^{-1}\in \Gamma$ has a positive exponent associated to the variable $z_{i,a}$. Fix $q\in \cB_+\cap m\Gamma_+$ (%
see Lemma \ref{th:qm}) and let $k\geq 0$ be the exponent of $z_{i,a}$ in the expansion~of $qm^{-1}\in \Gamma_+$. Consider
$$ m''=q \smash{\big( y_{i,a}y_{i^*,a-h}\cdot y_{i,a+2}y_{i^*,a+h+2}\cdot z_{i,a}^{-1} \big)^k.}$$
Then $m''$ lies in $\cB_+$ as $y_{i,a}y_{i,a+2}z_{i,a}^{-1}=\prod_{j\sim i}y_{j,a+1}$ does. Also, by Lemmas \ref{th:lowest}, \ref{th:qm2}~and~\ref{lem:Oia-tensor-Oi*a+h}, both $y_{i,a+2}y_{i^*,a+h+2}$ and $y_{i^*,a-h}y_{i,a}$ are monomials of $\Gamma_+$ in which the variable $z_{i,a}$ does not appear. This same variable hence cannot appear in 
$$m''m^{-1}\in \smash{(qm^{-1}z_{i,a}^{-k})}\Gamma_+\subseteq \Gamma_+$$
by definition of $k$. However, by the above, $m''$ appears in the indexing set of the intersection defining $V''(m)$, and hence $y_{\bR_p}\leq m''$%
, but then $m''m^{-1}\in y_{\bR_p}m^{-1}\Gamma_+\cap \Gamma_+$ contradicts the fact that $z_{i,a}$ has a positive exponent in the expansion of $y_{\bR_p}m^{-1}$. This ends the proof.
\end{proof}

\begin{Rem} \label{rem:support-preceq}
    If $\xi\le \xi'$ is a pair of height functions such that $m$ is supported between $\xi$ and $\xi'$ (i.e.~$m$ is a Laurent monomial in the $y_{i,a}$'s with $(i,a)\in I\times_2\Z$ and $\xi_i\leq a\leq \xi_i'$), then any monomial $m'$ as in %
    Theorem \ref{thm:global-sections-negative-line-bundles} is also supported between $\xi$ and $\xi'$. Indeed, following the reasoning used for example in \cite[Section 5.2.4]{hernandez2010cluster}, fix $i\in I$ so that, if
$$ \textstyle m' = m \prod_{(i,a)\in I\times_2\Z} z_{i,a}^{-c_{i,a}}\in \cB_+,$$ 
and $a$ is maximal (resp.~minimal) amongst all pairs with $c_{i,a}>0$, then %
$a< \xi_i'$ (resp.~$a\geq \xi_i$) as else the exponent of $y_{i,a+2}$ (resp.~$y_{i,a}$) in $m'$ would be $-c_{i,a}<0$, contradicting $m'\in \cB_+$. A consequence of this is that there are only finitely many $m'$ as in Theorem~\ref{thm:global-sections-negative-line-bundles}. Indeed, 
\begin{enumerate}
\item $\wt(m')$ is a dominant weight bounded above by $\wt(m)$, and
\item for each dominant weight $\la$, there are only finitely many elements of $\cB_+$ supported between $\xi$ and $\xi'$ with weight $\la$.
\end{enumerate}
\end{Rem}

\begin{Example}
  Take $\fg=\mathfrak{sl}_4$, and 
  $$\smash{m=\tfrac{y_{1,-1}y_{1,1}y_{3,-1}y_{3,1}}{y_{2,0}}.}$$
Then there are three $m'\in \cB_+$ with $m'\leq m$, namely $y_{2,0}$, $y_{1,-1}y_{1,1}$ and $y_{3,-1}y_{3,1}$. However,
$$ \smash{V(\varpi_2,0) = V(2\varpi_1, \{-1,1\})\cap V(2\varpi_3, \{-1,1\})},$$
and thus $V(m)$ is the finite sum%
  \[ \smash{V(m) = V(2\varpi_1, \{-1,1\}) + V(2\varpi_3, \{-1,1\}).} \]
\end{Example}

\begin{Rem} \label{rem:negative-line-bundles}
Given height functions $\xi\leq \xi'$, Theorem \ref{thm:global-sections-negative-line-bundles} implies that, for every line bundle $\cL_{m}$ on $\cZ_{\xi, \xi'}$ outside the globally generated cone (i.e.~for which $m\not\in \cB_+$), the space $H^0(\cZ_{\xi,\xi'},\cL_{m})$ is generated by its subspaces $H^0(\cZ_{\xi,\xi'},\cL_{m}(-D_{i,b}))$ of sections vanishing on the boundary divisors $\{D_{i,b}\,|\,(i,b)\in H(\xi, \xi')\}$. From this perspective, it is natural~to~ask~if the corresponding statement holds for an arbitrary free Bott--Samelson variety $\cZ_H$ (and/or for the non-free version $Z_H$). Indeed, this statement holds for $\cZ_H$ whenever $H$ is alternating (because of Remark \ref{rem:Hxi-alternating} and Theorem \ref{thm:global-sections-negative-line-bundles}), but it is reasonable to hope that a more direct geometric proof exists and would work for any $H$.
\end{Rem}

We record one significant consequence of Theorem \ref{thm:global-sections-negative-line-bundles}.  Choose height functions $\xi\le \xi'$ and recall the subalgebra
    \[ \textstyle \smash{\cR_{\xi,\xi'}=\bigoplus_{\tau\in \cP} H^0(\cZ_{\xi, \xi'}^\circ, \cL_{\tau}) \subseteq \cR}.\]

\begin{Corollary} \label{cor:Rxi-xi'-gen}
     The $V(\varpi_i,a)$'s with $(i,a)\in I\times_2\Z$ and $\xi_i\le a \le \xi_i'$ generate $\cR_{\xi,\xi'}$. 
\end{Corollary}

\begin{proof}
    Let $\cR'_{\xi,\xi'}$ be the subalgebra of $\cR$ generated by the $V(\varpi_i,a)$'s for which $(i,a)\in I\times_2\Z$ satisfies $\xi_i\leq a\leq \xi_i'$. For such a pair $(i,a)$, we have
    \begin{equation*}
    V(\varpi_i,a)=H^0(\cZ_{\infty},\O_{i,a}(1))\simeq H^0(\cZ_{\xi,\xi'},\O_{i,a}(1))\subseteq H^0(\cZ_{\xi,\xi'}^\circ,\cL_{\tau_{i,a}})\subseteq\cR_{\xi,\xi'},
    \end{equation*}
which implies $\cR'_{\xi,\xi'} \subseteq \cR_{\xi,\xi'}$. Conversely, following the proof of Proposition \ref{th:cRgamma} gives
$$ \textstyle \smash{H^0(\cZ_{\xi, \xi'}^\circ, \cL_{\tau})}=\bigcup_m V(m),$$
where the union runs over the monomials $m\in \cB$ that satisfy $\awt(m)=\tau$ and are supported between $\xi$ and $\xi'$. Thus, as $V(m)$ is contained in $\cR'_{\xi,\xi'}$ if $m\in \cB_+$ by Lemma \ref{le:surjectR}, it follows from Theorem \ref{thm:global-sections-negative-line-bundles} and Remark \ref{rem:support-preceq} that $V(m)$ is contained in $\cR'_{\xi,\xi'}$ for any $m$.
\end{proof}

\begin{Rem}
When $\xi$ and $\xi'$ are sufficiently far apart (explicitly, $\xi'\ge \xi+2h-2$),~one~can give a much simpler proof of Corollary \ref{cor:Rxi-xi'-gen} using that, for each $m\in \cB$ supported between $\xi$ and $\xi'$, there exists $m'' \in \cB_+$, also supported between $\xi$ and $\xi'$, that satisfies  $m''\geq m$ (this follows from $y_{i,a}y_{i^*,a+h}\in \Gamma_+$ and the fact that, if $\xi_i \le a \le \xi_i'$, then either $\xi_{i^*} \le a-h \le \xi'_{i^*}$ or $\xi_{i^*} \le a+h \le \xi'_{i^*}$). This simpler proof is also related to the observation that $\widehat{\cZ}_{\xi,\xi'}^\circ$ is, in the ``far apart case'', a closed subscheme of 
\[\prod_{\substack{(i, a) \in I\times _2 \Z \\ \xi_i \le a \le \xi'_i}} V(\varpi_i)^*\]
(as one has ``enough space'' between $\xi$ and $\xi'$ for relations of the form (\ref{eq:Zhat-pairing-equation})). This however fails if $\xi$ and $\xi'$ are not far apart enough. For example, $\widehat{\mathcal{Z}}_{\xi,\xi}^{\circ}\simeq N_-\backslash G$ is not an affine~scheme.
\end{Rem}

\subsection{Subalgebras of $K_\C(\O_{sh})$}\label{subsec:subalgebras_of_KOsh}
Fix a pair of height functions $\xi \leq \xi'$ as above and consider the Serre subcategory $\mathcal{O}_{sh}(\xi,\xi')$ of $\mathcal{O}_{sh}$ generated by the simple objects lying in $\mathcal{O}_{sh}^{\la}(\bR)$ for some $\la\in P_+$ and $\bR\in \Z^{\la}$, with $y_{\bR}$ supported between $\xi$ and $\xi'$. Then $\mathcal{O}_{sh}(\xi,\xi')\subseteq \mathcal{O}_{sh}$~is~a monoidal subcategory by Theorem \ref{thm:truncoprod}. We thus have a ring inclusion
\begin{equation}\label{eq:Incxixi'}
K_{0}(\mathcal{O}_{sh}(\xi,\xi'))\hookrightarrow K_0(\mathcal{O}_{sh}),
\end{equation}
and Proposition \ref{prop:mult_is_surj_truncation} allows us to identify $ K_{0}(\mathcal{O}_{sh}(\xi,\xi'))$ with the subring of $K_0(\mathcal{O}_{sh})$~generated by the subgroups $K_0(\mathcal{O}^{\varpi_i}_{sh}(a))$ for which $\xi_i\leq a \leq \xi'_i$. 

We also consider the category $ \scrC_{sh}(\xi, \xi') $ of finite-dimensional objects of $ \mathcal{O}_{sh}(\xi,\xi') $ and its unshifted version $ \scrC_0(\xi, \xi') \subseteq \scrC_0$. 
Applying Corollaries \ref{cor:Omega_restrict_NG}  and \ref{cor:Rxi-xi'-gen}, we conclude:

\begin{Proposition}\label{prop:IncXiXi'}
For height functions $\xi \leq \xi'$, the map $ \Omega $ gives algebra isomorphisms
$$
\cR_{\xi, \xi'} \cong K_\C(\O_{sh}(\xi,\xi')) \quad {}^N \cR_{\xi, \xi'} \cong K_\C(\scrC_{sh}(\xi, \xi')) \quad {}^G \cR_{\xi, \xi'} \cong K_\C(\scrC_0(\xi, \xi'))
$$
\end{Proposition}\vspace*{-1mm}

We highlight interesting cases of the above that come from special choices of $ \xi$ and $\xi'$.  

\begin{Theorem}\label{thm:GeoSubalgebrasK0}
Let $ \xi $ be a height function. The map $ \Omega $ restricts to $G$-equivariant algebra isomorphisms
\begin{enumerate}
\item $ \C[N_- \backslash G] \cong K_\C(\O_{sh}(\xi, \xi)), $
\item 
$\C[G]  \cong K_\C(\O_{sh}(\xi, \xi^* + h)), \text{ and}  $
\item
$ \C[ G \times_{N_- \backslash G} G] \cong K_\C(\O_{sh}(\xi, \xi + 2h)),$
\end{enumerate}
where, in the fiber product $G\times_{N_- \backslash G} G$, the left map $G\to N_-\backslash G$ is 
$g\mapsto N_-(\dot{w_0})^{-1}g$,  
and the right map is the tautological one.
\end{Theorem}

\begin{proof}
Recall the algebra isomorphism $\cR_{\xi,\xi'}\simeq \C[\hat{\cZ}_{\xi,\xi'}^{\circ}]$ of Corollary \ref{cor:iso_R_CZ}. Then (1) follows from  combining Proposition \ref{prop:IncXiXi'} with the canonical isomorphism $\hat{\cZ}_{\xi,\xi}^{\circ}\simeq N_-\backslash G$. Similarly, (2) follows from using Proposition \ref{prop:IncXiXi'} and the isomorphism $\hat{\cZ}_{\xi,\xi+h}^{\circ}\simeq G$ of Remark \ref{rem:embeddings_of_CG}.  %
Finally, for (3), we note that, by the definition of the Bott-Samelson varieties, for any triple $ \xi \le \xi' \le \xi''$ we have a $(G\times A)$-equivariant isomorphism (note that $\cZ_{\xi', \xi'}=\cZ_{\xi',\xi'}^{\circ}=\op{Fl}$)
$$
\cZ_{\xi, \xi''}^\circ\simeq  \cZ_{\xi, \xi'}^\circ \times_{\cZ_{\xi', \xi'}} \cZ_{\xi', \xi''}^\circ,
$$ 
where $G$ (resp.~$A$) acts diagonally on the right (resp.~left) of $\cZ_{\xi, \xi'}^\circ \times_{\cZ_{\xi', \xi'}} \cZ_{\xi', \xi''}^\circ$.
Applying~this to $(\xi',\xi'')=(\xi^*+h,\xi+2h)$ and passing to the universal $ A$-bundles, we get
$$ 
\widehat{\cZ}_{\xi, \xi+2h}^\circ \simeq G \times_{N_- \backslash G} G,
$$
and thus $ \cR_{\xi, \xi + 2h} \simeq \C[G \times_{N_- \backslash G} G ]$. Note that the map $\hat{\cZ}_{\xi,\xi^*+h}^{\circ}\to \hat{\cZ}_{\xi^*+h,\xi^*+h}^{\circ}$ corresponds~to $g\mapsto (\dot{w_0})^{-1}g$ because of Remark \ref{rem:what_are_generalized_minors}.
\end{proof}

Taking $ G$ invariants in (3) above and using the known fact that the map 
$$ N_- \rightarrow (G \times_{N_- \backslash G} G)/G $$ 
given by $ n \mapsto ((\dot{w_0})^{-1},n) $ is an isomorphism, we obtain:

\begin{Corollary} \label{cor:Ccoeur}
Let $ \xi $ be a height function.  The map $ \Omega $ restricts to an isomorphism
$$
\C[ N_-]  \cong K_\C(\scrC_0(\xi, \xi + 2h)).
$$
\end{Corollary}
\begin{Rem} Choose a height function $\xi$, and denote by $Q$ the orientation of the Dynkin diagram $I$ determined by $\xi$ (see Section \ref{subsubsec:height_functions}). In \cite{hernandez2011quantum}, the authors consider the Serre subcategory $\mathcal{C}_Q\subseteq\mathscr{C}_0$ generated by the simple objects whose highest $\ell$-weight~is~a~monomial in the $\sfY_{i,a}$'s with $\xi_i+2\leq a\leq \xi_{i^*}+h$. Also, Hernandez--Leclerc show, using cluster algebras, that there is an algebra isomorphism
$$ K_{\C}(\mathcal{C}_Q)\simeq \C[N].$$
Hence, the subcategories $\mathcal{C}_Q$ and $\mathscr{C}_0(\xi,\xi+2h)$ of $\mathscr{C}_0$ have isomorphic Grothendieck rings. In addition, using 
$$ \sfY_{i,a}=\tfrac{\sfPsi_{i,a-2}}{\sfPsi_{i,a}}\in \cB(\varpi_i+\varpi_{i^*}, (a-2)_i\cup(a+h)_{i^*}),$$
it is not hard to prove that $\mathcal{C}_Q\subseteq \mathscr{C}_0(\xi,\xi+2h)$. We expect this inclusion to be an equality. We also expect the category $\mathscr{C}_0(\xi,\xi')$ associated to a general pair of height functions $\xi\leq \xi'$ to be the same as the Serre subcategory of $\mathscr{C}_0$ generated by the simple objects whose highest $\ell$-weight is a monomial in the  $\sfY_{i,a}$'s with $\xi_i+2\leq a\leq \xi_{i^*}'-h$. This would enable us to see the categories $\mathscr{C}_0(\xi,\xi')$ as generalizations of the categories denoted by $\{\mathcal{C}_{\ell}\}_{\ell\geq 0}$ in \cite{hernandez2010cluster}.
\end{Rem}
\subsection{An alternative presentation in type $ A $}\label{subsec:alternate_presentation_type_A}
We now specialize to $\fg=\mathfrak{sl}_n$ and identify $V(\varpi_i)$ with $\wedge^i \C^n = \wedge^i V(\varpi_1)$. This also specifies a highest weight vector for $V(\varpi_i)$, namely 
$v_{\varpi_i}=e_1\wedge\dots \wedge e_i$, where $\{e_1,\dots, e_n\}$ is the standard basis of $\C^n$%
. \medskip\par
Following Francone--Leclerc \cite{francone2025cluster}, we call $\SL_n$-\emph{band} a complex $(\Z{\times}n)$-matrix $$b=(b_{rv})_{r\in \Z, v\in \{1,\dots, n\}}$$
such that every contiguous $(n{\times}n)$-submatrix
\begin{equation}\label{eq:Contiguous}
b(s):=(b_{rv})_{r\in \{s,\dots, s+n-1\},v\in \{1,\dots, n\}} 
\end{equation}
has determinant 1. We write $ B(\SL_n) $ for the scheme of $\SL_n$-bands. By \cite[Section~3.1]{francone2025cluster},
sending a $\SL_n$-band $b=(b_{rv})_{r\in \Z, v\in \{1,\dots, n\}}$ to the sequence $(b(s))_{s\in \Z}$ identifies  $B(\SL_n)$ with the scheme of bands $B(\SL_n, c_{st})$ associated to the standard Coxeter element $c_{st}=s_1\dots s_{n-1}$ of $\SL_n$ (see Remark \ref{rem:Z-and-bands}).  \medskip\par
Let $\xi_{st}$ be the \textit{standard height function}, that is $(\xi_{st})_i=i$ for all $i\in \{1,\dots,n-1\}$.~Then $c_{\xi_{st}}=c_{st}$ and Remark \ref{rem:Zhat-and-bands} gives an isomorphism $\widehat \cZ_\infty^\circ  \simeq B(\SL_n,c_{st})\cong B(\SL_n)$. 

\begin{Proposition}\label{prop:identification_type_A}
With the notation above, we have the following.
\begin{enumerate}
\item\label{item:wedges_of_x_ones} Let $x\in \widehat \cZ_\infty^\circ $.  Then, for all $ (i, a) \in I \times_2 \Z $, the $(i,a)$-component of $x$ is  
$$ x_{i,a} = x_{1, a-i+1} \wedge x_{1, a-i +3} \wedge \cdots \wedge x_{1, a+i-3} \wedge x_{1, a+i-1}. $$
As a consequence, a point $x$ as above is uniquely determined by the array $(x_{1,a})_{a\in 2\Z+1}$.  
\item\label{item:array_of_droites} The scheme $\cZ_\infty^\circ$ parametrizes bi-infinite sequences $(l_{a})_{a\in \Z}\in (\PP^{n-1})^\Z$ of lines in $\C^n$ such that every contiguous $n$-tuple generates $\C^n$.
\item\label{item:iso_with_SLn_bands} The isomorphism $ B(\SL_n)\simeq \widehat{\cZ}_\infty^\circ$ induced by $\xi_{st}$ is given by $b%
\mapsto x%
$ where, for $r\in \Z$,
\begin{equation*}
x_{1,2r+1} = \textstyle \sum_{k=1}^n b_{rk} e_k^\ast =(b_{rv})_{v=1,\dots,n}.
\end{equation*}
In other words, $x_{1, 2r+1} $ is the $r^{\text{th}}$-row of the matrix $b$.
\end{enumerate}
\end{Proposition}

\begin{proof}
The first statement follows easily from \eqref{eq:Zhat-incidence-equation}--\eqref{eq:Zhat-pairing-equation} and induction on $i$. The second statement is immediate from the first.
Finally, the third statement is also clear since, given $b\in B(\SL_n)$, by Remark \ref{rem:Zhat-and-bands},
\begin{equation*}
x_{1,2r+1}=\textstyle v_{\varpi_1}^\ast b(r)=e_1^*b(r) =\sum_{k=1}^n b_{rk} e_k^\ast,
\end{equation*}
where $b(r)\in\SL_n$ is the $\smash{r^{\text{th}}}$-contiguous submatrix of $b$ (as in \eqref{eq:Contiguous}).
\end{proof}
As can be seen using the above proposition, for $w\in S_n$, the isomorphism $B(\SL_n)\simeq \widehat{\cZ}_\infty^\circ$ identifies $\Delta_{w\varpi_i,2r+i}\in \C[\widehat{\cZ}_{\infty}^{\circ}]$ with the function which maps $b\in B(\SL_n)$ to the minor of the matrix $b(r)$ given by the rows labelled $\{1,\dots, i\}$ and the columns labelled $\{w(1),\dots,w(i)\}$. As a special case, the matrix coordinate function $b\mapsto b_{rv}$ on $B(\SL_n)$ is identified with the minor $\Delta_{\varepsilon_v, 2r+1}$ on $\widehat{\cZ}_{\infty}^{\circ}$, where $\varepsilon_v \in S_n\varpi_1$ is the usual standard weight, i.e. $\varepsilon_1=\varpi_1$, and
$$ \varepsilon_v = s_{v-1}\dots s_1\varpi_1
$$
for $1< v\leq n$. 
Hence, since the isomorphism $\Omega:\cR\to K_{\C}(\mathcal{O}_{sh})$ sends $\Delta_{\varepsilon_v,2r+1}$ to the class of the chamber module $L_{\varepsilon_v,2r+1}$, the definition of $B(\SL_n)$ and the above isomorphism 
$$\C[B(\SL_n)]\simeq \C[\widehat{\cZ}_\infty^\circ]\simeq \cR$$ 
give the following alternative presentation of $K_{\C}(\mathcal{O}_{sh})$:
\begin{Corollary}
The algebra $K_\C(\O_{sh}) $ is generated by the classes $\{ [L_{\varepsilon_v, 2r + 1}] \}_{r \in \Z, v\in \{1,\dots,n\}} $ subject to the relations
\begin{equation} \label{eq:det}
\smash{\textstyle\sum_{w \in S_n} (-1)^{\ell(w)} \prod_{v=1}^n [L_{\varepsilon_{v}, 2r - 1 + 2w(v)}] = 1}
\end{equation}
for all $r \in \Z $.
\end{Corollary}
Fix $r \in \Z $ %
and recall that $ L_{\varepsilon_v,2r+1} $ is a module over %
$ Y^{\varpi_1}_{\varepsilon_v}(2r+1) $. Hence, since
$
\textstyle \sum_{j=1}^n \varepsilon_j %
= 0 %
,$ 
we deduce that the relation \eqref{eq:det} takes place in $K_\C(\O^{n \varpi_1}_0(\bR(r)))$, where
\begin{equation*}
\bR(r) = \{2r +1,2r+3, \dots, 2r +(2n-1) \}.
\end{equation*}
We will need the following (very) special case of \cite[Theorem~4.3]{webster2020quantum}.
\begin{Proposition}\label{prop:qMV}
The natural map $ U \mathfrak{sl}_n \rightarrow Y_0$ (cf.~Section \ref{sec:Coproduct}) induces an isomorphism 
\begin{equation*}
U \mathfrak{sl}_n / \langle Z_+ \rangle \cong Y^{n \varpi_1}_0(\bR(r)),
\end{equation*}
where $ Z_+ $ is the positive degree part of the centre of $ U \mathfrak{sl}_n$.  
\end{Proposition}

Let $\mathcal{O}_{\chi_0}$ be the principal block of the BGG category $\mathcal{O}$ for $\mathfrak{sl}_n$. For $w\in S_n$, denote~by~$\Delta_w$ the Verma module in $\mathcal{O}_{\chi_0}$ of highest weight %
 $w\cdot 0=w\rho-\rho$, where $\rho=\sum_{i=1}^{n-1}\varpi_i$ is the~Weyl vector. Then the BGG resolution (or equivalently the Weyl character formula) shows that
\begin{equation}\label{eq:BGG}
\textstyle \sum_{w \in S_n} (-1)^{\ell(w)} [\Delta_w] = 1
\end{equation}
(where the term on the right corresponds to the trivial $U\mathfrak{sl}_n$-module). It is hence natural to wonder if \eqref{eq:det} is just the image of \eqref{eq:BGG} under the isomorphism of Proposition \ref{prop:qMV}.~We show that this is indeed the case (tacitly assuming Conjecture \ref{conj:Associators}).

\begin{Theorem}\label{thm:verma_is_tensor_prod}
Fix $w\in W$. Then the map of Proposition \ref{prop:qMV} induces an identification
\begin{equation*}
\Delta_{w^{-1}}\cong L_{\varepsilon_1, 2r-1+2w(1)} \otimes L_{\varepsilon_2, 2r-1+2w(2)} \otimes \cdots \otimes L_{\varepsilon_n, 2r-1+2w(n)},
\end{equation*}
and thus \eqref{eq:det} is naturally identified with \eqref{eq:BGG}. In particular, \eqref{eq:det} is categorified by the long exact sequence coming from the BGG resolution in $\mathcal{O}_{\chi_0}$ of the trivial $\mathfrak{sl}_n$-module.
\end{Theorem}
\begin{proof}
Let $a\in 2\Z+1$ and choose $k\in \{1,\dots,n\}$. Then it follows easily from Lemma \ref{lem:highest_ell_weight_chamber_module_braid} and \cite[Corollary 4.5]{friesen2025braid} that the highest $\ell$-weight of the chamber module $L_{\varepsilon_k,a}$ is
\begin{equation}\label{eq:highest_ell_weight_of_Lepsilon}
T^{-1}_{(s_{k-1}\dots s_1)^{-1}}(\sfPsi_{1,a})=(T^{-1}_{s_{k-1}}\circ\dots\circ  T^{-1}_{s_1})(\sfPsi_{1,a})=\tfrac{\sfPsi_{k,a-k+1}}{\sfPsi_{k-1,a-k}},
\end{equation}
where we use the convention that $\sfPsi_{0,a-1}=\sfPsi_{n,a-n+1}=1$. In addition, by Lemma \ref{lemma:EqualityTrunc}~(and since $\varepsilon_k=\varpi_1-\alpha_1-\dots-\alpha_{k-1}$), for $k\geq 2$, 
$$ Y_{\varepsilon_k}^{\varpi_1}(a) \simeq Y_{-\varpi_{k-1}}^{\varpi_1}(a,\fg_{\{1,\dots,k-1\}}),$$
where the algebra on the right is a truncated shifted Yangian for $\fg_{\{1,\dots,k-1\}}\simeq \mathfrak{sl}_k$. Therefore, using Corollary \ref{cor:identification_of_cat_O}, we deduce that the chamber module $L_{\varepsilon_k,a}$ is a $\{1,\dots,k-1\}$-inflation of the negative prefundamental representation of $\mathcal{O}_{-\varpi_{k-1}}^{\varpi_1}(a,\fg_{\{1,\dots,k-1\}})$. (Note that this is very specific to type A and $\varpi_1$.) In particular, its normalized character (in the usual sense, cf.~\cite[Section 3.3]{hernandez2024shifted}) is easily seen to be (see Corollary \ref{cor:identification_of_cat_O} and \cite[Theorem 3.16]{hernandez2024shifted})
\begin{equation}\label{eq:charChamberA}
\textstyle \tilde{\chi}(L_{\varepsilon_k,a})=\prod_{i=1}^{k-1}(1-e^{-(\alpha_{i}+\dots+\alpha_{k-1})})^{-1},
\end{equation}
and there are thus no variables of the form $\mathsf{A}_{j,b}^{-1}$ with $j\not\in \{1,\dots,k-1\}$ in the normalized $\ell$-character of $L_{\varepsilon_k,a}$. We abbreviate this last property by declaring that $L_{\varepsilon_k,a}$ is \textit{concentrated on $\{1,\dots,k-1\}\subseteq I$}. (Note that $L_{\varepsilon_1,a}=L_{\varpi_1,a}$ has a trivial normalized $\ell$-character, which we abbreviate by declaring that this module is \textit{concentrated on $\emptyset\subseteq I$.}) \medskip\par
Now, choose $a_1,\dots, a_n\in 2\Z+1$ and consider the module
\begin{equation*}
V=L_{\varepsilon_1, a_1} \otimes L_{\varepsilon_2, a_2}\otimes \dots \otimes L_{\varepsilon_n,a_n}.
\end{equation*}
Fix $1\leq k< s\leq n$. Then the product $L_{\varepsilon_k,a_k}\otimes L(\sfPsi_{s,a_s-s+1})$ is simple by Theorem \ref{thm:HZcriterionTensSimpPol} since $L_{\varepsilon_k,a_k}$ is concentrated on the subset $\{1,\dots,k-1\}$ (in which $s$ does not lie). Thus,~applying Corollary \ref{cor:BigTensProd} shows that $V$ is of highest $\ell$-weight $\psi=\prod_{k=1}^n\sfPsi_{\varepsilon_k,a_k}$. Moreover, by \eqref{eq:highest_ell_weight_of_Lepsilon},
\begin{align*}\textstyle \wt(\psi)%
&=\textstyle \frac{1}{2}\sum_{k=1}^{n-1}(a_{k+1}-a_k)\varpi_k-\rho
\end{align*}
and, by \eqref{eq:highest_ell_weight_of_Lepsilon}--\eqref{eq:charChamberA}, the normalized character of $V$ is
\begin{equation}\label{eq:NormCharVA}
\textstyle \tilde{\chi}(V) = \prod_{k=1}^n\prod_{i=1}^{k-1}(1-e^{-(\alpha_{i}+\dots+\alpha_{k-1})})^{-1}=\prod_{\alpha\in \Delta_+}(1-e^{-\alpha})^{-1}.
\end{equation} 
Finally, suppose that $a_k=2r+(2w(k)-1)$ for all $1\leq k\leq n$. Since, for each $k$,
$$ \langle w^{-1}\rho,\alpha_k^{\vee}\rangle = \langle \rho,w\alpha_k^{\vee}\rangle = w(k+1)-w(k)$$
(as one easily deduces from %
type A combinatorics), the above expression for $\wt(\psi)$ gives 
$$ \textstyle \wt(\psi) = \sum_{k=1}^{n-1}(w(k+1)-w(k))\varpi_k-\rho=\sum_{k=1}^{n-1}\langle w^{-1}\rho, \alpha_k^{\vee}\rangle\varpi_k-\rho = w^{-1}\rho-\rho%
$$
In particular, when the $a_k$'s are as above, the $U\mathfrak{sl}_n$-module associated to $V$ (via Proposition \ref{prop:qMV}) is a highest weight module for which the normalized character \eqref{eq:NormCharVA} and the highest weight $\wt(\psi)=w^{-1}\cdot 0$ are exactly those of %
$\Delta_{w^{-1}}$%
. Thus
$V\simeq \Delta_{w^{-1}}$ as $U\mathfrak{sl}_n$-modules%
.
\end{proof}
The results of this section can be summarized as follows:
\begin{Corollary}\label{cor:alternate_pres_type_A}
When $\fg=\mathfrak{sl}_n$, the algebra $K_{\C}(\O_{sh})$ is generated by the classes~of~chamber modules associated to %
$\varpi_1$, with only the relations coming from the BGG resolution.
\end{Corollary}

We find it remarkable that all the relations in $ K_0(\O_{sh}) $ can in this case be deduced from a single standard fact about the representation theory of $ \mathfrak{sl}_n $.

\section{Extension of Hernandez--Leclerc's duality}\label{sec:duality}
This section defines an algebra involution $D$ of $K_0(\mathcal{O}_{sh})$ which extends the rational analogue of the isomorphism given in \cite[Section 5.3]{hernandez2016cluster} for quantum affine Borel algebras~and in \cite[Theorem 8.7]{hernandez2023representations} for shifted quantum affine algebras. We show that the involution~$D$ sends classes of simple modules to classes of simple modules using the general framework of \textit{cobased modules} (see Definition \ref{def:cobasedModule}) and Lusztig's work  \cite{lusztig1993introduction} on canonical bases.~This answers (the obvious counterpart of) a question asked in \cite[Appendix A]{pinet2024functor}.

\subsection{An involution of \texorpdfstring{$G$}{G}}
Let $ \invomega : G \rightarrow G $ be the involutive automorphism of $ G $ which integrates the Lie algebra involution given by
\begin{equation}\label{eq:involution_of_G}
e_i \mapsto f_i ,\quad f_i \mapsto e_i, \quad h_i \mapsto -h_i.
\end{equation}
For any $G$-representation $ V $, let ${}^\invomega V$ denote the pullback of $ V $ by $\omega$.  Clearly,  ${}^\invomega V(\la)\simeq V(\la^*)$, and there exists a unique $G$-equivariant isomorphism 
\begin{equation}\label{eq:def_of_invomega}
\invomega_\lambda :\smash{ {}^\invomega V(\lambda) }\rightarrow V(\lambda^*) 
\end{equation}
such that $ \invomega_\lambda(v_\lambda) = v_{\lambda^*}^{\low} $, where $v_{\lambda^*}^{\low}:=\dot{w_0} v_{\lambda^*}$ for $w_0$ the longest element of $W$. In addition, $\omega_{\la^*}\circ \omega_{\la}=\mathrm{id}_{V(\la)}$ since $\omega(\dot{w_0})=\dot{w_0}^{-1}$.

For $ \lambda, \mu\in P_+$, consider the unique $G$-equivariant maps
\begin{align}\label{eq:lambda_mu_injection_projection}
\iota_{\lambda,\mu} : V(\lambda + \mu) \hookrightarrow V(\lambda) \otimes V(\mu) \ \text{ and }\ \pi_{\lambda,\mu}:V(\lambda) \otimes V(\mu)\twoheadrightarrow V(\lambda+\mu)
\end{align}
such that $\iota_{\lambda,\mu}(v_{\lambda+\mu})=v_\lambda\otimes v_\mu$ and $\pi_{\lambda,\mu}(v_\lambda\otimes v_\mu)=v_{\lambda+\mu}$. The result below is immediate.

\begin{Proposition} \label{th:iotaomega}
Fix $ \lambda, \mu\in P_+$. Then the following diagrams commute:
\begin{equation*}
\adjustbox{scale=0.9}{
\begin{tikzcd}[row sep =1.75em]
{}^{\invomega}V(\lambda + \mu) \ar[r, "\iota_{\lambda,\mu}"] \ar[d,"\invomega_{\lambda + \mu}"'] & {}^\invomega V(\lambda) \otimes {}^\invomega V(\mu) \ar[d,"\invomega_\lambda \otimes \invomega_\mu"] \\
V(\lambda^* + \mu^*) \ar[r, "\iota_{\lambda^*, \mu^*}"] & V(\lambda^*) \otimes V(\mu^*)
\end{tikzcd}}
\ \text{ and }\ 
\adjustbox{scale=0.9}{ \begin{tikzcd}[row sep =1.75em]
{}^\invomega V(\lambda) \otimes {}^\invomega V(\mu) \ar[r, "\pi_{\lambda,\mu}"] \ar[d,"\invomega_\lambda \otimes \invomega_\mu"'] & {}^\invomega V(\lambda + \mu) \ar[d,"\invomega_{\lambda + \mu}"] \\
V(\lambda^*) \otimes V(\mu^*) \ar[r, "\pi_{\lambda^*, \mu^*}"] & V(\lambda^* + \mu^*)
\end{tikzcd}}
\end{equation*}
\end{Proposition}\newpage

Recall, for $ i \in I $, the $G$-equivariant map 
\begin{center}
$\iota_{1,i}:=\iota_{1}: V(2\varpi_i - \alpha_i) \rightarrow V(\varpi_i) \otimes V(\varpi_i)$
\end{center}
defined in \eqref{eq:Gequivariant_map_wedge} via
\begin{equation*}
\iota_{1,i}(v_{2\varpi_i - \alpha_i})=v_{\varpi_i}\otimes v_{s_i\varpi_i}-v_{s_i\varpi_i}\otimes v_{\varpi_i}%
=v_{\varpi_i}\wedge f_iv_{\varpi_i}
\end{equation*}
(where we use the fact that $v_{s_i\varpi_i}=f_iv_{\varpi_i}$). Similarly, denote by 
\begin{equation}\label{eq:def_of_inverseiota}
\inverseiota_{1,i}: V(\varpi_i) \otimes V(\varpi_i)\to V(2\varpi_i - \alpha_i)
\end{equation}
the unique $G$-equivariant left-inverse of $\iota_{1,i}$.

\begin{Proposition}\label{th:iota2omega}
Fix $i\in I$ and write $\omega_i=\omega_{\varpi_i}:{}^\omega V(\varpi_i)\to V(\varpi_{i^*})$. The diagrams
\begin{equation*}
\adjustbox{scale=0.9}{
\begin{tikzcd}[row sep =1.75em]
{}^\invomega V(2 \varpi_i - \alpha_i) \ar[r,"\iota_{1,i}"] \ar[d,"\invomega_{2 \varpi_i - \alpha_i}"'] & {}^\invomega V(\varpi_i) \otimes{}^\invomega V(\varpi_i) \ar[d,"\invomega_{i} \otimes \invomega_{i}"] \\
V(2 \varpi_{i^*} - \alpha_{i^*}) \ar[r,"-\iota_{1,i^\ast}"] & V(\varpi_{i^*}) \otimes V(\varpi_{i^*})
\end{tikzcd}} \ \text{ and }\ 
\adjustbox{scale=0.9}{
\begin{tikzcd}[row sep =1.75em]
{}^\invomega V(\varpi_i) \otimes {}^\invomega V(\varpi_i) \arrow[d,"\invomega_i\otimes\invomega_i"']\arrow[r,"\inverseiota_{1,i}"] & 
{}^\invomega V(2 \varpi_i - \alpha_i) \arrow[d,"\invomega_{2 \varpi_i - \alpha_i}"]\\ 
V(\varpi_{i^*}) \otimes V(\varpi_{i^*}) \arrow[r,"-\inverseiota_{1,i^\ast}"]&
V(2 \varpi_{i^*} - \alpha_{i^*})
\end{tikzcd}
}
\end{equation*}
commute (notice the minus sign on the lower horizontal arrows).
\end{Proposition}

\begin{proof}
The commutativity of the first diagram easily follows from the observation that
$$ ((\invomega_{i} \otimes \invomega_{i}) \circ\iota_{1,i})(v_{2 \varpi_i - \alpha_i}) = (\invomega_{i} \otimes \invomega_{i})(v_{\varpi_i} \wedge  f_i v_{\varpi_i}) = v_{\varpi_{i^*}}^{\low} \wedge e_iv_{\varpi_{i^*}}^{\low},$$
whereas, by $G$-equivariance of $\iota_{1,i^*}$ and since $\dot{w}_0f_{i^*}=-e_i\dot{w}_0$,
\begin{align*}
(-\iota_{1,i^\ast}\circ \invomega_{2 \varpi_i - \alpha_i})(v_{2\varpi_i - \alpha_i})  &= -\iota_{1,i^\ast}(\dot{w}_0 v_{2 \varpi_{i^*} - \alpha_{i^*}})=-\dot{w}_0(v_{\varpi_{i^*}}\wedge f_{i^*}v_{\varpi_{i^*}})\\&=-\dot{w}_0v_{\varpi_{i^*}}\wedge \dot{w}_0f_{i^*}v_{\varpi_{i^*}}=v^{\low}_{\varpi_{i^*}}\wedge e_i v^{\low}_{\varpi_{i^*}}.
\end{align*}
For the second diagram, observe that it suffices to check the commutativity on the image of $\iota_{1,i}:{}^\invomega V(2\varpi_i-\alpha_i)\to {}^\invomega V(\varpi_i)\otimes {}^\invomega V(\varpi_i)$ since its $G$-equivariant left-inverse $\pi_{1,i}$ annihilates~all other simple composition factors of ${}^\invomega V(\varpi_i)\otimes {}^\invomega V(\varpi_i)$ (and $\omega_i\otimes \omega_i$ sends $\op{Im} \iota_{1,i}$ to $\op{Im} \iota_{1,i^*}$%
). This verification however directly follows from the commutativity of the first diagram.
\end{proof}

\subsection{Compatibility with the dual canonical basis of tensor products}\label{subsec:effect_of_inv_on_basis}

Let $q$ be an indeterminate and let $\bfU= U_q(\fg)$ be the quantum group associated to $ \fg $, which is an algebra over $\C(q)$ with generators $\{e_i,f_i,k_i^{\pm 1}\}_{i\in I}$.  Recall that the coproduct $\Delta : \mathbf U \rightarrow \mathbf U \otimes \mathbf U $ is~the algebra morphism given by
\begin{equation}\label{eq:def_of_coproduct_Uqg}
\def\xspace{0.75}
\Delta(e_i)=e_i \otimes 1 + k_i \otimes e_i, \hspace{\xspace em} \Delta(f_i)= f_i \otimes k_i^{-1} + 1 \otimes f_i\ \text{ and }\  \Delta(k_i) = k_i \otimes k_i.
\end{equation}
Consider the \textit{bar-involution} $\overline{\phantom{a}}%
$ of $\C(q)$, which is the $\C$-linear involution %
satisfying $\overline{q}=q^{-1}$. Also, let $\sigma$ be the \textit{bar-involution} of $\bfU$, i.e.~the $\C$-linear algebra involution %
defined by
\begin{equation*}
\def\xspace{0.75}
\sigma(e_i)=e_i,\hspace{\xspace em}\sigma(f_i)=f_i,\hspace{\xspace em}\sigma(k_i^{\pm1})=k_i^{\mp1}\ \text{ and }\ \sigma(q^{\pm 1})=\overline{(q^{\pm 1})}=q^{\mp 1}.
\end{equation*} 
Denote by $\omega$ the $\C(q)$-linear involution of $\bfU$ given by
\begin{equation*}
\def\xspace{0.75}
\invomega(e_i)=f_i, \hspace{\xspace em} \invomega(f_i)=e_i\ \text{ and }\ \invomega(k_i^{\pm1})=k_i^{\mp1}.
\end{equation*}
Clearly, $\omega$ quantizes %
 \eqref{eq:involution_of_G}. Finally, let $\tau$ be the $ \C$-linear anti-involution of $\bfU$ defined by
\begin{equation*}
\def\xspace{0.75}
\tau(e_i)=f_i, \hspace{\xspace em} \tau(f_i)=e_i,\hspace{\xspace em} \tau(k_i^{\pm1})=k_i^{\mp1}\ \text{ and }\ \tau(q^{\pm 1})=\overline{(q^{\pm 1})}=q^{\mp 1}.
\end{equation*}

\newcommand{\dual}{\vee}

We use the anti-involution $\tau$ to define what we %
call 
\textit{contragredient duals}. 
\begin{Def} Fix a $\bfU$-module $V$. The \textit{contragredient dual} $V^\dual$ is the $\C$-vector space $\Hom_{\C(q)}(V,\C(q))$ equipped with the $\bfU$-action uniquely determined by
\begin{equation*}
\langle a \cdot \alpha, v \rangle = \langle \alpha, \tau(a) \cdot v \rangle 
\end{equation*}
for $ a \in \bfU$, $v \in V$ and $\alpha \in V^\dual $, where $ \langle .,. \rangle $ is the evaluation pairing between $V^\dual$ and $V$. 
\end{Def}
\begin{Rem} The evaluation pairing is $\C(q)$-linear in the second variable. We say that it is \textit{$\C(q)$-anti-linear} in the first variable  as $\langle x \cdot \alpha, v \rangle=\overline{x}\langle \alpha,v\rangle$ for $x\in \C(q)$, $v \in V$~and~$\alpha \in V^\dual$.
\end{Rem}

From \eqref{eq:def_of_coproduct_Uqg}, we see that $\invomega$ and $\tau$ satisfy
\begin{equation*}
\smash{ (\invomega \otimes \invomega) \circ \Delta = \Delta^{\mathrm{op}}\circ \invomega\quad\text{and}\quad(\tau \otimes \tau) \circ \Delta = \Delta^{\mathrm{op}}\circ \tau},
\end{equation*}
where $\Delta^{\mathrm{op}}$ is the opposite coproduct. Thus, given finite-dimensional $\bfU$-modules $V$ and~$W$, the swap of tensor factors $P:V\otimes W\to W\otimes V$ induces an isomorphism
\begin{equation*}
\smash{{}^{\omega}(V\otimes W)\simeq {}^{\omega}W\otimes {}^{\omega}V},
\end{equation*}
where ${}^{\omega}V$ denotes as before the pullback of the %
module $V$ by $\omega$. Similarly, the swap gives %
rise to a canonical isomorphism of $\bfU$-modules
\begin{equation}\label{eq:contra_and_op}
\smash{(V \otimes W)^\dual \cong W^\dual \otimes V^\dual}.
\end{equation}
\begin{Rem}
In what follows, we will often use isomorphisms of the form \eqref{eq:contra_and_op} implicitly. In particular, we will often write 
$\langle \beta\otimes \alpha,v\otimes w\rangle =  \langle \alpha,v\rangle \langle \beta,w\rangle$ when $\alpha\in V^{\dual}$, $\beta\in W^{\dual}$,~$v\in V$ and $w\in W$. This gives, for $x,y\in \bfU$,
\begin{equation}\label{eq:tensorevaluationpairing}
\langle (x\cdot \beta)\otimes (y\cdot \alpha),v\otimes w\rangle = \langle y\cdot \alpha,v\rangle\langle x\cdot \beta,w\rangle = \langle\beta\otimes \alpha, (\tau(y)\cdot v)\otimes (\tau(x)\cdot w)\rangle.
\end{equation}
\end{Rem}

Recall that %
the category of finite-dimensional $\mathbf U$-modules (of type I) is semisimple, with simple objects indexed by %
$P_+$, just as for finite-dimensional $G$-modules. Moreover, specialization of $q$ to 1 is compatible with taking weight spaces of representations and identifies the simple $\bfU$-module of highest weight $\la$ with $V(\la)$. We will hence also denote by $\{V(\la)\}_{\la\in P_+}$ the irreducible $\bfU$-modules (which are now viewed as $\C(q)$-vector spaces).

For $\la\in P_+$, let $(.,.)_{\la}$ be the \textit{Jantzen-Shapovalov form} associated to $V(\la)$, i.e.~the unique $\C$-bilinear inner product on $V(\la)$ that is $\C(q)$-linear~in~the second factor and satisfies
\begin{equation*}
(v_\lambda,v_\lambda)_{\la}=1\hspace*{0.75em}\text{and}\hspace*{0.75em}(a v_1, v_2)_\lambda = (v_1, \tau(a) v_2)_\lambda
\end{equation*}
for all $v_1,v_2\in V(\la)$ and $ a \in \bfU $.  The forms $(.,.)_{\la}$ induce $ \bfU $-module isomorphisms
\begin{equation}\label{eq:shap_induces_iso}
\begin{aligned}
V(\la) &\cong V(\la)^\dual,\\
v&\mapsto (v, -)_\lambda,
\end{aligned}
\end{equation}
leading in turn (with \eqref{eq:contra_and_op}) to isomorphisms %
\begin{equation}\label{eq:shap_induces_iso_tensor_prod}
(V(\la_1) \otimes \cdots \otimes V(\la_n))^\dual \cong V(\la_n)^\dual \otimes \cdots \otimes V(\la_1)^\dual \cong V(\la_n) \otimes \dots \otimes V(\la_1)
\end{equation}
for all $\la_1,\dots,\la_n\in P_+$.

Equip each simple $\mathbf{U}$-module $V(\lambda)$ with its dual canonical basis $B(\lambda)$ (which is a weight basis of $V(\lambda)$, see \cite{lusztig1990canonicalI} along with Theorem \ref{thm:dual_canonical_KLR}) %
and fix an isomorphism 
\begin{equation}\label{eq:def_of_invomega2}
\smash{\invomega_\lambda : {}^\invomega V(\lambda) \rightarrow V(\lambda^*)}
\end{equation}
such that 
\begin{equation*}
\smash{\invomega_\lambda(v_\lambda) = v_{\lambda^*}^{\low}}, 
\end{equation*}
with $\smash{v_{\lambda^*}^{\low}}$ the unique lowest weight vector of $V(\lambda^{\ast})$ contained in $B(\lambda^\ast)$ (see Remark~\ref{Rem:bothmapsomegalambdaareequal}). Extend this to lists of dominant weights $\ul{\la}= (\lambda_1, \dots, \lambda_n)$ using the $\bfU$-module isomorphisms
\begin{equation}\label{eq:omegabR}
\begin{aligned}
\invomega_{\ul{\lambda}} : {}^{\invomega}(V(\lambda_1) \otimes \cdots \otimes V(\lambda_n)) &\rightarrow V(\lambda_n^*) \otimes \cdots\otimes V(\lambda_1^*),  \\
v_1 \otimes \cdots \otimes v_n &\mapsto \smash{\invomega_{\lambda_n}(v_n) \otimes \cdots \otimes \invomega_{\lambda_1}(v_1)}.
\end{aligned}
\end{equation}
Finally, denote by $B(\ul{\lambda})$ the dual canonical basis of the product $V(\ul{\la})=V(\la_1)\otimes \dots \otimes V(\la_n)$ (see \eqref{eq:def_of_dc_basis_tensorprod}). Our aim is to show the result below, proven for $n=1$ in \cite[Prop 21.1.2]{lusztig1993introduction}.

\begin{Theorem} \label{th:basistobasis}
Fix $\ul{\la}$ as above. The map $ \invomega_{\ul{\la}} $ carries the dual canonical basis $B(\ul{\la}) $ to the dual canonical basis $ B(\la_n^*,\dots,\la_1^*)$ of $V(\la_n^*)\otimes \dots \otimes V(\la_1^*)$. 
\end{Theorem}

We will need Lusztig's formalism of \textit{based modules} \cite[Section~27.1.2]{lusztig1993introduction}, or, more~precisely the dual formalism (which we formulate slightly differently than Khovanov \cite{khovanov1997graphical}).%

\begin{Def}\label{def:cobasedModule}
Given a finite-dimensional $\bfU$-module $V$ and a $\C(q)$-basis $B$ of $V$, let $B^\dual $ be the basis~of~$V^\dual$ dual to $B$ with respect to the evaluation pairing $\langle .,.\rangle$. The pair $(V,B)$~is a \emph{cobased module} if $ (V^\dual,B^\dual) $ is a based module in the sense of \cite[Section 27.1.2]{lusztig1993introduction}.  
\end{Def}

We will not recall the definition of based modules here, but we will use the following~consequences of this definition.  First, if $ (V, B) $ is a cobased module, then $B$ is a weight basis~of $V$. For the second consequence, let 
$\smash{\sigma_B:V\to V}$ be the \textit{bar-involution} of $V$, i.e.~the $\C$-linear map given (for $x\in \C(q)$ and $b\in B$) by 
\begin{equation*}
\smash{\sigma_B(xb)=\overline{x}b}.
\end{equation*}
Then, by the results of \cite[Section~27.1.2(c)]{lusztig1993introduction}, $\sigma_B$ is compatible 
with the bar-involution $\sigma$ of $\bfU$ %
in the sense that
\begin{equation}\label{eq:sigmaCobased}
\sigma_B(av)=\sigma(a)\sigma_B(v)
\end{equation}
for
$a\in \bfU$ and $v\in V$. %
\par The maps $\sigma_B$ and $\sigma_{B^\dual}$ are adjoint in the sense below:

\begin{Def}\label{def:adjoint} Given a $\C(q)$-anti-linear map $T:V\to W$, the \textit{adjoint} $T^\dual: W^\dual\to V^\dual$~of $T$ is the unique $\C(q)$-anti-linear map satisfying 
$${\langle \beta,T(v)\rangle=\overline{\langle T^\dual(\beta),v\rangle}}$$ 
for all $v\in V$ and $\beta\in W^\dual$.
\end{Def}

The dual canonical bases $\{B(\la)\}_{\la\in P_+}$ provide the primary example of cobased modules. 
\begin{Theorem}[{\cite[Section~27.1.4]{lusztig1993introduction}}]\label{thm:dc_basis_is_cobased}
For $\la\in P_+$, $(V(\lambda),B(\lambda))$ is a cobased module.%
\end{Theorem}

\begin{Example}\label{ex:dc_basis_type_A1}
If $\fg=\fsl_2$, endow $V(\varpi_1)%
$ and $V(2\varpi_1)%
$ with their canonical bases $\{v_1',v_{-1}'\}$ and $\{w_2',w_0',w_{-2}'\}$ (see \cite[Section 1.4]{frenkel1997canonical}). The $\bfU$-action on these bases yields diagrams
\begin{equation*}
\adjustbox{scale=1,center}{
\begin{tikzcd}
v_1'\arrow[loop below,"q"] \arrow[r,bend left,"{1}"{pos=0.55}]& v_{-1}'\arrow[l,bend left,"{1}"{pos=0.45}] \arrow[loop below,"q^{-1}"]
\end{tikzcd}
\hspace{2em}
\begin{tikzcd}
w_2'\arrow[loop below,"q^2"] \arrow[r,bend left,"{1}"]& w_0'\arrow[loop below,"1"]\arrow[r,bend left,"{q+q^{-1}}"{pos=0.58}]\arrow[l,bend left,"{q+q^{-1}}"]& w_{-2}'\arrow[l,bend left,"{1}"{pos=0.45}] \arrow[loop below,"q^{-2}"]
\end{tikzcd}}
\end{equation*}
where the action of the generator $e_1$ (resp.~$f_1$) of $\bfU$ points left (resp.~right), and the loops indicate the action of $k_1$. The arrows are decorated with the scalars appearing in the action.\par 
Now, let $B(\varpi_1)^{\vee}$ and $B(2\varpi_1)^{\vee}$ be the bases obtained via \eqref{eq:shap_induces_iso}, namely
$$ B(\varpi_1)^{\vee}=
\{v_1^{\vee}%
,v_{-1}^\dual%
\}\subseteq  V(\varpi_1)^{\vee} \text{ and } B(2\varpi_1)^{\vee}=\{w_2^{\dual}%
,w_0^{\dual}%
,w_{-2}^{\dual}%
\}\subseteq  V(2\varpi_1)^{\vee},$$
where $v_i^{\vee}=(v_i',.)_{\varpi_1}$ and $w_j^{\vee}=(w_j',.)_{2\varpi_1}$ for $i\in \{-1, 1\}$ and $j\in \{-2,0,2\}$. \par The dual canonical basis  $B(2\varpi_1)=\{w_2,w_0,w_{-2}\}$ is the basis that is dual to $B(2\varpi_1)^{\vee}$~via the evaluation pairing, i.e.~$\langle w_i^{\vee},w_j\rangle = \delta_{i,j}$ for $i,j\in \{-2,0,2\}$. It is easy to see that $w_2=w_2'$ and $w_{-2}=w_{-2}'$ using the properties of the form $(.,.)_{2\varpi_1}$. Moreover, these properties and
$$ \langle w_0^{\vee},w_0'\rangle=(w_0',w_0')_{2\varpi_1}= (f_1w_2',w_0')_{2\varpi_1}=(w_{2}',e_1w_0')_{2\varpi_1}=(w_2',(q+q^{-1})w_2')_{2\varpi_1}=q+q^{-1}$$ 
show that $w_0=(q+q^{-1})^{-1}w_0'$. Repeating the same analysis for $B(\varpi_1)$ gives the diagrams
\begin{equation*}
\adjustbox{scale=1,center}{
\begin{tikzcd}
v_1\arrow[loop below,"q"] \arrow[r,bend left,"{1}"]& v_{-1}\arrow[l,bend left,"{1}"{pos=0.45}] \arrow[loop below,"q^{-1}"]
\end{tikzcd}
\hspace{2em}
\begin{tikzcd}
w_2\arrow[loop below,"q^2"] \arrow[r,bend left,"{q+q^{-1}}"]& w_0\arrow[loop below,"1"]\arrow[r,bend left,"{1}"]\arrow[l,bend left,"{1}"]& w_{-2}\arrow[l,bend left,"{q+q^{-1}}"{pos=0.43}] \arrow[loop below,"q^{-2}"]
\end{tikzcd}}
\end{equation*}
which describe the $\bfU$-action on $B(\varpi_1)$ and $B(2\varpi_1)$. This ends our example.
\end{Example}

Given two cobased modules $(V, B)$ and $(W,C)$, the pair $(V\otimes W,B\otimes C)$ is not a cobased module in general. However, there exists a basis $B\hsop C\subseteq V\otimes W$ such that $(V\otimes W,B\hsop C)$ is a cobased module. To describe this basis, recall the partial order $\geq_{B\times C}$  on the set $B\times C$ given in \cite[Section~27.3.1]{lusztig1993introduction}, namely, for $(b,c), (b',c')\in B\times C$ with 
\begin{equation*}
b\in V_\mu, \quad b'\in V_{\mu'},\quad c\in W_\nu, \quad c'\in W_{\nu'}, 
\end{equation*}
we say that $(b,c)\geq_{B\times C} (b',c')$ if and only if
\begin{equation}\label{eq:def_of_Lusztigs_partial_order}
\mu\geq \mu' %
\text{ and }%
\mu+\nu=\mu'+\nu'.
\end{equation} 
Consider also the quasi-$R$-matrix $\Theta\in \bfU\hat{\otimes}\bfU$ \cite[ch.~4]{lusztig1993introduction} and define a $\C$-linear map
\begin{equation}\label{eq:def_of_inv_Psi}
\begin{aligned}
\Psi:V\otimes W &\to V\otimes W,\\
v\otimes w&\mapsto \Theta(\sigma_{B}(v)\otimes \sigma_{C}(w))
\end{aligned}
\end{equation}
(which depends on both $B$ and $C$ despite the notation). By \eqref{eq:sigmaCobased} along with standard~facts about $\Theta$, the map $\Psi$ is an involution.\par
We will prove the following dual version of \cite[Theorem~27.3.2]{lusztig1993introduction}.

\begin{Theorem}\label{thm:cobased_properties}
Let $(V,B)$ and $(W,C)$ be cobased modules. For all $(b,c)\in B\times C$, there exists a unique vector $b\smallheart c\in V\otimes W$ such that
\begin{enumerate}
\item\label{item:prop1cobased} $\Psi(b\smallheart c)=b\smallheart c$ and 
\item\label{item:prop2cobased} $b\otimes c-b\smallheart c\in q \op{span}_{\Z[q]}\{b'\otimes c' \,|\,(b',c')\in B\times C\text{ and }  (b,c)\geq_{B\times C} (b',c')\}$.
\end{enumerate}
Moreover, the basis $B\hsop C=\{b\smallheart c\,|\,(b,c)\in B\times C\}$ %
makes $(V\otimes W,B\hsop C)$ a cobased module.
\end{Theorem}

Let $P:\bfU\hat{\otimes}\bfU\to \bfU\hat{\otimes}\bfU$ be the swap %
and set $\sigma^{\otimes 2}:=\sigma\otimes\sigma$ with $\tau^{\otimes 2}:=\tau\otimes \tau$. Clearly,~all three of $P$, $\sigma^{\otimes 2}$ and $\tau^{\otimes 2}$ commute. To prove the above theorem, we need:

\begin{Lemma}\label{lem:quasi_R_matrix_sigmatauP}
The quasi-$R$-matrix satisfies $(\sigma^{\otimes 2}\circ\tau^{\otimes 2}\circ P)(\Theta)=\Theta$.
\end{Lemma}

\begin{proof}
Recall the $Q$-grading on $\bfU$ defined by $\deg(e_i)=-\deg(f_i)=\alpha_i$ and $\deg(k_i)=0$ for $i\in I$. Write $\bfU^+$ (resp. $\bfU^-$) for the subalgebra of $\bfU$ generated by the $e_i$'s (resp. the $f_i$'s) and let $\bfU^{\pm}_{\pm \nu}$ be the subspace of $\bfU^{\pm}$ of homogeneous elements of degree $\pm \nu$ (with $\nu\in Q_+$).\par By \cite[Theorem 4.1.2(a)]{lusztig1993introduction}, the quasi-$R$-matrix $\Theta$ of $\bfU\hat{\otimes}\bfU$ is uniquely characterized by the following two properties:
\begin{enumerate}
\setlength{\itemsep}{1pt}
\item[(1)] $\Theta=1\otimes 1+\sum_{\nu\in Q_+\backslash\{0\}}\Theta_{\nu}$ where $\Theta_{\nu}\in \bfU_{-\nu}^-\otimes \bfU_{\nu}^+$ for all $0\neq \nu\in Q_+$, and
\item[(2)] $\Delta(x)\cdot \Theta=\Theta\cdot (\sigma^{\otimes 2}\circ \Delta\circ \sigma)(x)$ for all $x\in \bfU$.
\end{enumerate}
Let $\tilde{\Theta}:=(\sigma^{\otimes 2}\circ \tau^{\otimes 2}\circ P)(\Theta)$. Using (1) above and the definition of $P$, $\sigma$ and $\tau$, we get
$$\textstyle \tilde{\Theta}\in 1\otimes 1+\sum_{\nu\in Q_+\backslash\{0\}}\bfU_{-\nu}^-\otimes \bfU_{\nu}^+,$$
and thus $\tilde{\Theta}$ also satisfies (1). Moreover, an easy computation 
shows that, for $x\in \bfU$,
\begin{align*}
\Delta(x)\cdot \tilde{\Theta}%
=(\sigma^{\otimes 2}\circ \tau^{\otimes 2}\circ P)\big(\Theta\cdot (\sigma^{\otimes 2}\circ \Delta\circ \sigma)((\sigma\circ \tau)(x))\big).
\end{align*}
Hence, applying (2) above, we find that
\begin{align*}
\Delta(x)\cdot\tilde{\Theta}&= (\sigma^{\otimes 2}\circ \tau^{\otimes 2}\circ P)(\Delta((\sigma\circ \tau)(x))\cdot \Theta)%
=\tilde{\Theta}\cdot (\sigma^{\otimes 2}\circ\Delta\circ \sigma)(x),
\end{align*}
and thus $\tilde{\Theta}$ also satisfies (2). This implies $\Theta=\tilde{\Theta}$ by uniqueness. 
\end{proof}

\begin{proof}[Proof of Theorem \ref{thm:cobased_properties}]
Write $B=\{b_1,\dots, b_k\}$ and $C=\{c_1,\dots, c_\ell\}$, and denote by
\begin{equation*}
B^\dual=\{b_1^\dual,\dots, b_k^\dual\}\subset V^\dual \text{ and } C^\dual=\{c_1^\dual,\dots, c_\ell^\dual\}\subset W^\dual
\end{equation*}
the corresponding dual bases. By definition, $(V^{\vee},B^{\vee})$ and $(W^{\vee},C^{\vee})$ are based modules.\par 

Let $\Psi'$ be the endomorphism of $W^{\dual}\otimes V^{\dual}%
$ defined by \eqref{eq:def_of_inv_Psi}. On the one hand, by Lemma \ref{lem:quasi_R_matrix_sigmatauP} and \eqref{eq:tensorevaluationpairing}, we have that, for $v\in V$, $w\in W$, $\alpha\in V^\dual$ and $\beta\in W^\dual$,
\begin{align*}
\langle \Psi'(\beta \otimes \alpha), v \otimes w \rangle %
&= \langle \sigma_{C^{\dual}}(\beta)\otimes \sigma_{B^\dual}(\alpha), (\tau^{\otimes 2}\circ P)(\Theta)\cdot(v\otimes w)\rangle\\
&=\overline{\langle \beta\otimes \alpha,(\sigma_B\otimes \sigma_C)\big((\tau^{\otimes 2}\circ P)(\Theta)\cdot (v\otimes w)\big)\rangle}\\
&= \overline{\langle \beta\otimes \alpha, \Theta\cdot(\sigma_B(v)\otimes \sigma_C(w))\rangle}=\overline{\langle \beta\otimes \alpha, \Psi(v\otimes w)\rangle},
\end{align*}
i.e.~$\Psi'$ is the adjoint of $\Psi$ in the sense of Definition \ref{def:adjoint} (when seen as an endomorphism of $(V\otimes W)^{\dual}$ via \eqref{eq:contra_and_op}). On the other hand, using \cite[Theorem 27.3.2]{lusztig1993introduction}, we deduce~that,~for every $1\leq i\leq k$ and $1\leq j\leq \ell$, there exists a unique vector $c_{j}^\dual\diamond b_{i}^\dual\in W^\dual\otimes V^\dual$ such that
\begin{enumerate}
\item[(1')] $\Psi'(c_{j}^\dual\diamond b_{i}^\dual)=c_{j}^\dual\diamond b_{i}^\dual$ with
\item[(2')] $c_{j}^\dual\otimes b_{i}^\dual- c_{j}^\dual\diamond b_{i}^\dual\in q^{-1} \op{span}_{\Z[q^{-1}]}\{c^\dual\otimes b^\dual \;|\; (c_{j}^\dual,b_{i}^\dual)\geq_{C^{\dual}\times B^{\dual}} (c^\dual,b^\dual)\}$,
\end{enumerate}
and the set 
$${C^\dual\dsop B^\dual=\{c_{j}^\dual\diamond b_{i}^\dual\}_{1\leq i\leq k,1\leq j\leq \ell}}$$ 
is a basis %
for which $(W^{\vee}\otimes V^{\vee},C^{\vee}\dsop B^{\vee})$ is a based module. %
Let 
$${B\hsop C:=\{b_i\smallheart c_j\}_{1\leq i\leq k,1\leq j\leq \ell}\subset V\otimes W}$$ 
be the corresponding dual basis, that is (using \eqref{eq:tensorevaluationpairing})
$$ {\langle c_{j'}^{\vee}\diamond b_{i'}^{\vee},b_i\smallheart c_j\rangle =\delta_{i,i'}\delta_{j,j'}}$$
for all $1\leq i,i'\leq k$ and $1\leq j,j'\leq \ell$. By definition, $(V\otimes W,B\hsop C)$ is a cobased module. Moreover, for $1\leq i,i'\leq k$ and $1\leq j,j'\leq \ell$, by (1') above,
\begin{equation*}
\langle c_{j'}^\dual\diamond b_{i'}^\dual,\Psi(b_{i}\smallheart c_j)\rangle=\overline{\langle \Psi'(c_{j'}^\dual\diamond b_{i'}^\dual),b_{i}\smallheart c_j\rangle}=\overline{\langle c_{j'}^\dual\diamond b_{i'}^\dual,b_{i}\smallheart c_j\rangle}=\delta_{i,i'}\delta_{j,j'},
\end{equation*} 
and thus $\Psi(b_{i}\smallheart c_j)=b_{i}\smallheart c_j$. This shows property (\ref{item:prop1cobased}) of the statement of the theorem.\par 
For property (\ref{item:prop2cobased}), note that, by (2'), the change-of-basis matrix from $C^\dual\dsop B^\dual$ to $C^\dual\otimes B^\dual$ is unitriangular with respect to the partial order $\geq_{C^{\vee}\times B^{\vee}}$ of \eqref{eq:def_of_Lusztigs_partial_order} and that the coefficients above the diagonal are in $q^{-1}\Z[q^{-1}]$. Taking dual bases, it follows that the change-of-basis matrix from $B\otimes C$ to $B\hsop C$ is unitriangular with respect to the partial order $\geq_{B\times C}$ and that the coefficients above the diagonal belong to $\overline{q^{-1}\Z[q^{-1}]}=q\Z[q]$. The same must~be~true for the inverse change-of-basis matrix %
by the usual inversion formula for unitriangular matrices. Hence, (\ref{item:prop2cobased}) holds for $B\hsop C$. Finally, uniqueness %
follows from \cite[Theorem~27.3.2(1)]{lusztig1993introduction}.
\end{proof}

\begin{Rem}\label{rem:heartAssociative} 
Given three cobased modules $(V_1,B_1)$, $(V_2,B_2)$ and $(V_3,B_3)$, the bases $B_1\hsop (B_2\hsop B_3)$ and $(B_1\hsop B_2)\hsop B_3$ are equal by \cite[Section~27.3.6]{lusztig1993introduction}.
\end{Rem}

For a product $V(\ul{\lambda})=V(\la_1)\otimes \dots\otimes V(\la_n)$ with $\la_1,\dots,\la_n\in P_+$,
the dual canonical~basis $B(\ul{\lambda})$ is defined as (see Remark \ref{rem:heartAssociative}) 
\begin{equation}\label{eq:def_of_dc_basis_tensorprod}
B(\ul{\lambda}):=B(\lambda_1)\hsop \dots \hsop B(\lambda_n).
\end{equation}
By Theorem \ref{thm:cobased_properties}, this basis makes $(V(\ul{\lambda}),B(\ul{\lambda}))$ into a cobased module.

\begin{Example}
Building on Example \ref{ex:dc_basis_type_A1}, the dual canonical basis $B(\varpi_1)\hsop B(2\varpi_1)$  of $V(\varpi_1)\otimes V(2\varpi_1)$ is obtained by computing the basis which is dual under \eqref{eq:shap_induces_iso_tensor_prod} to the
basis $B(2\varpi_1)^\dual\dsop B(\varpi_1)^\dual$ of $V(2\varpi_1)^\dual\otimes V(\varpi_1)^\dual%
$. Using %
 \cite[Proposition~1.7]{frenkel1997canonical}, we see that the latter basis can be described as
\begin{gather*}
w_2^\dual\diamond v_{-1}^\dual=w_2^\dual\otimes v_{-1}^\dual+ q^{-1}w_0^\dual\otimes v_1^\dual ,\hspace{2em} w_{0}^\dual\diamond v_{-1}^\dual=w_0^\dual\otimes v_{-1}^\dual+q^{-2}w_{-2}^\dual\otimes v_1^\dual, \\
w_2^\dual\diamond v_1^\dual =w_2^\dual\otimes v_1^\dual,\hspace{0.6em} w_0^\dual\diamond v_1^\dual =w_0^\dual\otimes v_1^\dual ,\hspace{0.6em} w_{-2}^\dual\diamond v_1^\dual =w_{-2}^\dual\otimes v_1^\dual ,\hspace{0.6em} w_{-2}^\dual\diamond v_{-1}^\dual =w_{-2}^\dual\otimes v_{-1}^\dual.
\end{gather*}
For $x,y\in \C(q)$ and $v=x v_1\otimes w_0 + y v_{-1}\otimes w_2$, a direct computation gives
$$ \langle w_2^{\vee}\diamond v_{-1}^{\vee},v\rangle=\langle w_2^{\vee}\otimes v_{-1}^{\vee},v\rangle+q\langle w_0^{\vee}\otimes v_1^{\vee},v\rangle=y+qx \hspace*{0.75em} \text{and}\hspace*{0.75em}  \langle w_0^{\vee}\diamond v_{1}^{\vee},v\rangle=x.$$
Consequently, the vectors
\begin{equation*}
v_{-1}\smallheart w_2=v_{-1}\otimes w_{2}\text{ and }v_{1}\smallheart w_0=v_1\otimes w_0-q v_{-1}\otimes w_2
\end{equation*}
give a basis for the weight space $(V(\varpi_1)\otimes V(2\varpi_1))_{\varpi_1}$ which is dual to $\{w_2^\dual\diamond v_{-1}^\dual,w_0^\dual\diamond v_1^\dual\}$.
Using similar computations, we get that the dual canonical basis $B(\varpi_1)\hsop B(2\varpi_1)$ is 
\begin{gather*}
v_{1} \smallheart w_{0} =v_{1} \otimes w_{0} -q v_{-1} \otimes w_{2}  ,\hspace{2em} v_{1} \smallheart w_{-2} =v_{1} \otimes w_{-2} -q^2 v_{-1} \otimes w_{0} , \\
v_1 \smallheart w_2  =v_1 \otimes w_2 ,\hspace{0.75em} v_{-1} \smallheart w_2  =v_{-1} \otimes w_2  ,\hspace{0.75em} v_{-1} \smallheart w_{0}  =v_{-1} \otimes w_{0}  ,\hspace{0.75em} v_{-1} \smallheart w_{-2}  =v_{-1} \otimes w_{-2} .
\end{gather*}
This is compatible with Theorem \ref{thm:dual_canonical_tensorprod}. Indeed, the unique simple modules $L_{-1}$ and $L_{2}$~of~the cyclotomic KLRW-algebras $T^{(\varpi_1)}_{-\varpi_1}$ and $T^{(2\varpi_1)}_{2\varpi_1}$ satisfy
\begin{equation*}
\ch(L_{-1}) = [(1,\redbold{1})]\hspace*{0.75em}\text{and}\hspace*{0.75em}\ch(L_{2})=[(\redbold{1})],
\end{equation*}
and the $T^{(\varpi_1,2\varpi_1)}_{\varpi_1}$-module $\bbS^{\varpi_1,2\varpi_1}(L_{-1},L_{2})$ must therefore be simple, compatibly with the equality $v_{-1}\smallheart w_2=v_{-1}\otimes w_{2}$, since, by Proposition \ref{prop:standardization_and_shuffles},
\begin{equation*}
\ch(\bbS^{\varpi_1,2\varpi_1}(L_{-1}, L_2))=\ch(L_{-1})\shuffle_e \ch(L_2)=[(1,\redbold{1},\redbold{2})].
\end{equation*}
On the other hand, the unique simple modules $L_{1}$ and $L_{0}$ of $T_{\varpi_1}^{(\varpi_1)}$ and $T_{0}^{(2\varpi_1)}$ satisfy
\begin{equation*}
\ch(L_{1}) = [(\redbold{1})]\quad\text{and}\quad\ch(L_{0})=[(1,\redbold{1})],
\end{equation*}
and we recover the specialization at $q=1$ of the equality $v_{1}\smallheart w_0=v_1\otimes w_0-qv_{-1}\otimes w_2$ via
\begin{equation*}
\ch(\bbS^{\varpi_1,2\varpi_1}(L_{1}, L_0))=[(\redbold{1},1,\redbold{2})]+[(1,\redbold{1},\redbold{2})]=[(\redbold{1},1,\redbold{2})]+\ch(\bbS^{\varpi_1,2\varpi_1}(L_{-1},L_2)).
\end{equation*}
\end{Example}

Fix cobased modules $(V, B)$, $(W,C)$, $(\widetilde V, \widetilde B)$ and $(\widetilde W, \widetilde C) $. Suppose there exist $\bfU$-module isomorphisms  $\varphi_V:{}^{\omega}V\to\widetilde{V}$ and $\varphi_W:{}^{\omega}W\to\widetilde{W}$ with $\varphi_V(B)=\widetilde{B}$ and $\varphi_W(C)=\widetilde{C}$. Note that, in this case, $\varphi_V$ (resp.~$\varphi_W$) intertwines $\sigma_{B}$ and $\sigma_{\widetilde{B}}$ (resp.~$\sigma_{C}$ and $\sigma_{\widetilde{C}}$).

Consider the $\bfU$-module isomorphism 
\begin{equation}\label{eq:varphiOmega}
\begin{aligned} 
\varphi :  {}^{\invomega} (V \otimes W)  &\to \widetilde W \otimes \widetilde V, \\
v \otimes w &\mapsto \varphi_W(w) \otimes \varphi_V(v)=((\varphi_W\otimes \varphi_V)\circ P)(v\otimes w)
\end{aligned}
\end{equation}
and let $\Psi\in \End_\C(V\otimes W)$ and $\tilde{\Psi}\in \End_\C(\widetilde{W}\otimes\widetilde{V})$ be the involutions defined in \eqref{eq:def_of_inv_Psi}.

\begin{Lemma}\label{lem:varphi_intertwines}
The map $\varphi$ intertwines $\Psi$ and $\tilde{\Psi}$, that is the following diagram commutes
\begin{equation*}
\adjustbox{scale=0.85}{
\begin{tikzcd}
{}^{\omega}(V\otimes W) \ar[r,"\Psi"]\ar[d,"\varphi"]& {}^{\omega}(V\otimes W)\ar[d,"\varphi"]\\
\widetilde W \otimes \widetilde V \ar[r,"\tilde{\Psi}"]& \widetilde W \otimes \widetilde V
\end{tikzcd}}
\end{equation*}
\end{Lemma}

\begin{proof}
Fix $v\in V$, $w\in W$. Since $\varphi_V$ (resp.~$\varphi_W$) intertwines $\sigma_{B}$ and $\sigma_{\widetilde{B}}$ (resp.~$\sigma_{C}$ and $\sigma_{\widetilde{C}}$),
\begin{equation*}
(\tilde\Psi \circ \varphi)(v \otimes w)=\Theta\cdot \big((\varphi_W\otimes \varphi_V)\circ (\sigma_{C}\otimes\sigma_{B})\circ P\big) (v\otimes w).
\end{equation*}
Moreover, since the $\varphi_V$ and $\varphi_W$ are isomorphisms of $\bfU$-modules, the above reduces to 
\begin{equation}\label{eq:PsiVarphi1}
(\tilde\Psi \circ \varphi)(v \otimes w)=(\varphi_W\otimes \varphi_V)\big((\omega\otimes\omega)(\Theta)\cdot \big((\sigma_{C}\otimes\sigma_{B})\circ P\big)(v\otimes w)\big),
\end{equation}
where $\cdot$ now denotes the $\bfU$-action on $W\otimes V$. By \cite[Theorem 4.1.2]{lusztig1993introduction}, we have 
$$(\omega\otimes\omega)(\Theta)=P(\Theta).$$ 
Hence \eqref{eq:PsiVarphi1} can be rewritten as
\begin{align*} 
(\tilde\Psi \circ \varphi)(v \otimes w)&=(\varphi_W\otimes \varphi_V)\big(P(\Theta)\cdot ((\sigma_{C}\otimes\sigma_{B})\circ P)(v\otimes w)\big)\\
&=\big((\varphi_W\otimes\varphi_V)\circ P\big)\big(\Theta\cdot (\sigma_{B}\otimes \sigma_{C})(v\otimes w)\big),
\end{align*}
where $\cdot$ now denotes the $\bfU$-action on $V\otimes W$. Using the definition of $\Psi$ and $\varphi$ thus gives
\begin{align*}
 (\tilde\Psi \circ \varphi)(v \otimes w)&=\varphi(\Theta\cdot(\sigma_{B}\otimes \sigma_{C})(v\otimes w))=(\varphi\circ \Psi)(v\otimes w),
\end{align*}
as desired.
\end{proof}

The following is immediate:

\begin{Lemma}\label{lem:bijection_posets}
The map from $B\times C$ to  $\widetilde{C}\times \widetilde{B}$ given by
\begin{align*} 
\smash{(b,c)\mapsto (\varphi_W(c),\varphi_V(b))}
\end{align*}
is a bijection of posets. 
\end{Lemma}

\begin{Lemma} \label{le:based}
With the above setup, $ \varphi(B \hsop C) = \widetilde C \hsop \widetilde B $.
\end{Lemma}

\begin{proof}
By Lemma \ref{lem:varphi_intertwines}, if $b\smallheart c\in B\hsop C$, then 
\begin{equation*}
(\widetilde{\Psi}\circ \varphi)(b\smallheart c)=(\varphi\circ \Psi)(b\smallheart c) =\varphi(b\smallheart c),
\end{equation*}
and $\varphi(B\hsop C)$ therefore satisfies the first property characterizing $\widetilde{C}\hsop\widetilde{B}$ in Theorem \ref{thm:cobased_properties}. Now, fix $\tilde{b}\in \widetilde{B}$ and  $\tilde{c}\in \widetilde{C}$. There are $b\in B$ and $c\in C$ such that $\varphi_V(b)=\tilde{b}$ and $\varphi_W(c)=\tilde{c}$. In addition, the element $b\smallheart c\in B\hsop C$ satisfies, by Theorem \ref{thm:cobased_properties},
\begin{equation*}
b\otimes c-b\smallheart c\in q \op{span}_{\Z[q]}\{b'\otimes c'\;|\; (b,c)\geq_{B\times C} (b',c')\}.
\end{equation*}\newpage
By Lemma \ref{lem:bijection_posets}, applying $\varphi$ gives that $\tilde{c}\otimes \tilde{b}-\varphi(b\smallheart c)=\varphi(b\otimes c-b\smallheart c)$ belongs to
\begin{equation*}
q\op{span}_{\Z[q]}\{\varphi_W(c')\otimes \varphi_V(b')\;|\; (b,c)\geq_{B\times C} (b',c')\}=q\op{span}_{\Z[q]}\{\tilde{c}'\otimes \tilde{b}'\;|\; (\tilde{c},\tilde{b})\geq_{\widetilde{C}\times \widetilde{B}} (\tilde{c}',\tilde{b}')\}.
\end{equation*}
Hence $\varphi(B\hsop C)$ satisfies both properties characterizing $\widetilde{C}\hsop \widetilde{B}$ in Theorem \ref{thm:cobased_properties}, and thus equals $\widetilde{C}\hsop \widetilde{B}$ by uniqueness of this basis.
\end{proof}

We can finally prove the main result of this subsection. Recall the list $\ul{\la}=(\la_1,\dots,\la_n)$.

\begin{proof}[Proof of Theorem \ref{th:basistobasis}]
We use induction on $n$, with the case $ n = 1$ proven in \cite[Prop 21.1.2]{lusztig1993introduction}. Suppose $n>1$ and write $\ul{\la}[2,n]=(\la_2,\dots,\la_n)$ with $\ul{\la}^*[n,2]=(\la_n^*,\dots,\la_2^*)$. Equip
\begin{equation*}
V=V(\la_1),\hspace*{0.5em}W=V(\ul{\la}[2,n]),\hspace*{0.5em}\widetilde{V}=V(\la^*_1)\ \text{ and }\ \widetilde{W}=V(\ul{\la}^*[n,2])
\end{equation*}
with their respective dual canonical bases
\begin{equation*}
B=B(\la_1),\hspace*{0.5em} C=B(\ul{\la}[2,n]),\hspace*{0.5em} \widetilde{B}=B(\la_1^*)\ \text{ and }\ \widetilde{C}=B(\ul{\la}^*[n,2]),
\end{equation*}
and let
\begin{equation*}
\varphi_V = \invomega_{\la_1}:{}^{\omega}V\to \widetilde{V}\ \text{ with }\ \varphi_W=\invomega_{\ul{\la}[2,n]}:{}^{\omega}W\to \widetilde{W}.
\end{equation*}
By the induction hypothesis, we have $\varphi_V(B)=\widetilde{B}$ and $\varphi_W(C)=\widetilde{C}$. Also, the map $\varphi$ given in \eqref{eq:varphiOmega} is canonically identified with $\omega_{\ul{\la}}$ and, by Lemma \ref{le:based},
\begin{equation*}
\omega_{\ul{\la}}(B\hsop C)=\varphi(B \hsop C) =  \widetilde C \hsop \widetilde B,
\end{equation*}
so that the result follows from the inductive definition of the dual canonical basis for tensor products (and Remark \ref{rem:heartAssociative}).
\end{proof}

\begin{Rem}\label{Rem:bothmapsomegalambdaareequal} 
The isomorphism $\omega_{\la}:{}^{\omega}V(\la)\to V(\la^*)$ given in \eqref{eq:def_of_invomega2} reduces precisely to \eqref{eq:def_of_invomega} after specializing $q$ to $1$. Indeed, it suffices to show that $\dot{w}_0v_{\lambda}$ lies in the dual canonical basis $B(\la)\subseteq V(\la)$, but this can be easily deduced from Theorem \ref{thm:dual_canonical_KLR} and Corollary \ref{cor:W_action_on_chamber_modules}.
\end{Rem}

\subsection{An involution on the bi-infinite Bott--Samelson variety}\label{subsec:invBottSam}
For each $i\in I$, the $\C$-vector space isomorphism $\invomega_{i}:V(\varpi_i)\simeq {}^{\omega}V(\varpi_i)\to V(\varpi_{i^\ast})$ given in \eqref{eq:def_of_invomega} can be dualized to produce a $\C$-vector space isomorphism $\invomega_{i}^\ast:V(\varpi_{i^\ast})^\ast\to V(\varpi_{i})^\ast$. Define a map 
\begin{align*}
\mathsf{F}:%
\prod_{(i,a) \in I \times_2 \Z } V(\varpi_i)^\ast&\to %
\prod_{(i,a) \in I \times_2 \Z } V(\varpi_i)^\ast%
\end{align*}
by declaring the $(i,a)$-component of $y=\mathsf{F}(x)$ to be related to the $(i^*,-a+h)$-component~of $x$ via
\begin{equation}\label{eq:componentsFbiinf}
y_{i,a} = \invomega_{i}^\ast(x_{i^*,-a+h})
\end{equation}
for each $(i,a)\in I\times_2\Z$. Clearly, $\mathsf{F}$ defines an involutive automorphism of $\prod_{(i,a) \in I \times_2 \Z } V(\varpi_i)^\ast$. Let $F$ be the restriction of $\mathsf{F}$ to the closed subscheme $\smash{\widehat \cZ_\infty^\circ}$ (see Proposition \ref{prop:presentation of Zhat}).
\begin{Theorem}\label{thm:F_is_involution_of_bi_infinite_bs}
The map $F$ defines an involutive automorphism of $ \widehat \cZ_\infty^\circ $.
\end{Theorem}
\begin{proof}
Choose $x \in \widehat \cZ_\infty^\circ$ and set $y=F(x)$. We need to check that the equations \eqref{eq:Zhat-Gr-equation}--\eqref{eq:Zhat-pairing-equation} hold for $y$ given the same equations for $x$ and \eqref{eq:componentsFbiinf}. Fix hence $i,j\in I$ and let  $\pi_{i,j}:=\pi_{\varpi_i,\varpi_j}$ be the $G$-equivariant map of \eqref{eq:lambda_mu_injection_projection}. Proposition \ref{th:iotaomega} implies
\begin{equation}\label{eq:Feq1}
 \pi_{i,j}^*\circ \invomega_{\varpi_i+\varpi_j}^*= (\invomega_{i}^*\otimes \invomega_{j}^*)\circ \pi_{i^\ast,j^\ast}^*.
\end{equation}
Fix $a\in\overline{i}+2\Z$. By \eqref{eq:Zhat-Gr-equation}, the element $ x_{i^*,-a+h} \otimes x_{i^*,-a+h} $ lies in $\op{Im}\pi^\ast_{i^\ast,i^\ast}$. Thus, by \eqref{eq:Feq1}, 
$$ y_{i,a} \otimes y_{i,a} =  (\invomega_{i}^\ast \otimes \invomega_{i}^\ast)( x_{i^*,-a+h} \otimes x_{i^*,-a+h} ) $$
lies in $\op{Im}\pi^\ast_{i,i}$ and \eqref{eq:Zhat-Gr-equation} holds for $y$.  Equation \eqref{eq:Zhat-incidence-equation} is handled in the same way. For \eqref{eq:Zhat-wedge-tensor-equation}, consider the two $G$-equivariant maps
\begin{equation*}
\inverseiota_{1,i}: V(\varpi_i)  \otimes V(\varpi_i)\to V(2 \varpi_i - \alpha_i) \ \text{ and }\ \pi_{2,i} : \textstyle\bigotimes_{j \sim i} V(\varpi_j) \to V(2 \varpi_i - \alpha_i)
\end{equation*}
defined respectively in \eqref{eq:def_of_inverseiota} and \eqref{eq:lambda_mu_injection_projection}. By Propositions \ref{th:iotaomega} and \ref{th:iota2omega},
\begin{equation*}
- (\invomega_{i}^* \otimes \invomega_{i}^*)\circ \inverseiota_{1,i^\ast}^* = \inverseiota_{1,i}^*\circ \invomega_{2 \varpi_i - \alpha_i}^*
\ \text{ and }\ 
\textstyle \pi_{2,i}^*\circ \invomega_{2 \varpi_i - \alpha_i}^* =  (\bigotimes_{j \sim i} \invomega_{j}^*) \circ \pi_{2,i^\ast}^*
\end{equation*}
respectively. Moreover, since \eqref{eq:Zhat-wedge-tensor-equation} holds for $x$, there exists $v\in V(2\varpi_{i^\ast}-\alpha_{i^\ast})^\ast$ such that
\begin{equation*}
\textstyle \pi_{1,i^\ast}^\ast(v)=x_{i^*,-a  +h-2} \wedge x_{i^*, -a+h} \ \text{ and }\ \pi_{2,i^\ast}^\ast(v)=\bigotimes_{j\sim i} x_{j^*, -a +h-1}.
\end{equation*}
Consequently, applying $\pi_{1,i}^\ast$ to the vector $\invomega_{2\varpi_i-\alpha_i}^\ast(v)\in V(2\varpi_i-\alpha_i)^\ast$ gives
\begin{align*}
(\pi_{1,i}^\ast\circ \invomega_{2\varpi_i-\alpha_i}^\ast)(v)&=-((\invomega_i^\ast\otimes \invomega_i^\ast)\circ \pi_{1,i^\ast}^\ast)(v)\\
&=-(\invomega_i^\ast\otimes \invomega_i^\ast)(x_{i^*,-a +h-2} \wedge x_{i^*, -a+h})\\
&=-y_{i,a+2}\wedge y_{i,a}=y_{i,a}\wedge y_{i,a+2},
\end{align*}
while applying $\pi_{2,i}^\ast$ to the same vector yields
\begin{align*}
(\pi_{2,i}^\ast\circ \invomega_{2\varpi_i-\alpha_i}^\ast)(v)&=\textstyle ((\bigotimes_{j\sim i} \invomega_{j}^\ast)\circ \pi_{2,i^\ast}^\ast)(v)\\
&=\textstyle \bigotimes_{j\sim i} \invomega_{j}^\ast(x_{j^*, -a +h-1})=\bigotimes_{j\sim i}\, y_{j,a+1},
\end{align*}
proving that \eqref{eq:Zhat-wedge-tensor-equation} holds for $y$. Finally, as explained in the discussion preceding Proposition \ref{prop:presentation of Zhat}, the fact that \eqref{eq:Zhat-pairing-equation} holds for $x$ (and the invertibility of the $\omega_i^*$'s) gives 
$$y_{i,a}=\omega_i^*(x_{i^*,-a+h}) \neq 0$$ 
for all $(i,a)\in I\times_2 \Z$, which is equivalent to saying that \eqref{eq:Zhat-pairing-equation} holds for $y$.
\end{proof}

Let $D:\C[\widehat \cZ_\infty^\circ]\to \C[\widehat \cZ_\infty^\circ]$ be the involution of $\C$-algebras induced from $F$ and recall that, for every $(i,a)\in I\times_2\Z$, there is a subspace $V(\varpi_i,a)\subset \cR\simeq \C[\widehat \cZ_\infty^\circ]$ coming from sections~of the line bundle $\O_{i,a}(1)$ over $\cZ_\infty$. By definition of $F$, we have:

\begin{Corollary}\label{cor:restriction_of_D_to_Via}
For $ (i,a) \in I \times_2 \Z $, the restriction of the map $D$ to $V(\varpi_i,a)$ has values~in $V(\varpi_{i^\ast},-a+h)$ and coincides with the linear map $\invomega_{i}:V(\varpi_i)\to V(\varpi_{i^\ast})$ of \eqref{eq:def_of_invomega}.
\end{Corollary}

\begin{Example}\label{ex:involution_HL_SLn}
Consider $\fg=\mathfrak{sl}_2$. By Proposition \ref{prop:identification_type_A}, $\widehat \cZ_\infty^\circ$ can be identified with the scheme of $(\Z{\times}2)$-matrices $b=(b_{r,v})_{r \in \Z, v \in \{1,2\}}$ for which every contiguous $(2{\times}2)$-submatrix has determinant $1$. From our choice of isomorphism $V(\varpi_1)\simeq \C^2$, we deduce that the linear map $\omega_1:V(\varpi_1)\simeq \C^2\to \C^2$ swaps the weight vectors $e_1$ and $e_2$. Hence, $F$ is given (under the above identification) by $ F(b) =c= (c_{r,v})_{r\in \Z,v\in\{1,2\}}$, where 
\begin{equation*}
c_{r,1} = b_{-r,2} \ \text{ and }\  c_{r,2} = b_{-r,1}.
\end{equation*}
In other words, we flip our matrices vertically and horizontally.

For $ G = \SL_n $, the situation is slightly more complicated and we need to replace each~$ b_{rv} $~in a matrix of $B(\SL_n)$ by the corresponding minor in the expansion of the determinant~of~the contiguous $(n{\times}n)$-submatrix with bottom row $b_{r1}, \dots, b_{rn}$. After this replacement, we flip the $(\Z{\times}n)$-matrix vertically.
\end{Example}

The involution $F$ that was given in Theorem \ref{thm:F_is_involution_of_bi_infinite_bs} in fact comes from an involution~on the bi-infinite Bott--Samelson variety $\cZ_{\infty}$. Indeed, the map $\omega^\ast_i:V(\varpi_{i^*})^*\to V(\varpi_i)^{*}$ induces~a morphism of homogeneous coordinate rings, and thus a map of projective varieties
\begin{equation*}
(\trou)^\perp:\Gr(i)\to \Gr(i^\ast).
\end{equation*}
This gives rise to an involution $f$ of $\cZ_{\infty}$ defined by 
$$ (x_{i,a})_{(i,a)\in I\times_2\Z} \mapsto (y_{i,a}:=x_{i^\ast,-a+h}^\perp)_{(i,a)\in I\times_2 \Z}$$ 
which recovers the involution $F$ after restricting to the open cell $\cZ_{\infty}^\circ$ and using the identifications $f^\ast \O_{i,a}(1)\simeq\O_{i^\ast,-a+h}(1)$ given by $\omega^\ast_i$ (that can be shown to be compatible with the choices of trivializations coming from Lemma \ref{lem:BS-divisor-expansion} using the proof of Theorem \ref{thm:F_is_involution_of_bi_infinite_bs}). 

\begin{Example}
When $G=\SL_n$, there is a simple description of $f$. Using the perspective illustrated in Example \ref{ex:example_Zinft_SL5}, we see that $f$ sends an array of subspaces $(W_{i,a}\subset \C^n)_{(i,a)\in I\times_2\Z}$ to the array of perpendicular spaces 
$$(W'_{i,a}:=W_{n-i,-a+h}^\perp\subseteq \C^n)_{(i,a)\in I\times_2\Z}$$ 
with respect to the symmetric bilinear form $(.,.)$ on $\C^n$ satisfying $(e_i,e_j)=(-1)^{i}\delta_{ij}$. 
\end{Example}
\subsection{An involution on \texorpdfstring{$K_\C(\O_{sh})$}{KC(Osh)}}
Under the isomorphism $\Omega:\cR\to K_\C(\O_{sh})$ of Theorem \ref{thm:Omega_is_surj}, the involution $ D$ of $\cR\simeq \C[\widehat \cZ_\infty^\circ] $ gives rise to an algebra involution 
$$ D:K_\C(\O_{sh})\to K_{\C}(\O_{sh}).$$
Also, by Corollary \ref{cor:restriction_of_D_to_Via}, this map $D$ restricts, for all $(i,a)\in I\times_2\Z$, to an isomorphism
\begin{equation}\label{eq:DVvarpiia}
K_\C(\O_{sh}^{\varpi_i}(a))\simeq V(\varpi_i,a) \to K_\C(\O_{sh}^{\varpi_{i^*}}(-a+h))\simeq V(\varpi_{i^\ast},-a+h)
\end{equation}
which can be identified with the map $\omega_i:V(\varpi_i)\to V(\varpi_{i^*})$ of \eqref{eq:def_of_invomega} (or \eqref{eq:def_of_invomega2}). We have: 
\begin{Theorem}\label{thm:DLpsi_is_still_the_class_of_a_simple}
The map $D$ sends classes of simple objects to classes of simple objects.
\end{Theorem}

\begin{proof}
Take $\psi \in \fr$. Then $L(\psi)$ lies in $\O_{sh}^\lambda(\bR)$ for some $ \lambda\in P_+$ and $\bR=(R_i)_{i\in I}\in \Z^\lambda$. Let $\varpi_{\bR}=(\varpi_{p_1},\dots,\varpi_{p_{\ell}})$ be the list of fundamental weights associated to $\bR$ in \eqref{eq:varpibR} and recall that this list depends on a choice of total order ``$<$'' on $I$. Denote by $\bR^*=(R^*_i)_{i\in I}\in \Z^{\la^*}$~the set of parameters given by
\begin{equation}\label{eq:def_of_bRast}
R_i^\ast=\{-r+h\;|\;r\in R_{i^\ast}\}
\end{equation}
for $i\in I$, and consider the corresponding list $\varpi_{\bR^*}$ where we now use for $I$ the total order ``$<^{\mathrm{*}}$'' given by $i<^* j$ if and only if $j^* < i^*$. It is easy to see that $\varpi_{\bR^*}=(\varpi_{p_{\ell}^*},\dots,\varpi_{p_1^*})$.

On the other hand, as the algebra involution $D$ restricts to isomorphisms of the form~\eqref{eq:DVvarpiia}, we have the commutative diagram
\begin{equation*}
\adjustbox{scale=1}{
\begin{tikzcd}
V(\varpi_\bR) \cong K_{\C}(\mathcal{O}_{sh}^{\varpi_{p_1}}(r_1))\otimes\dots \otimes K_{\C}(\mathcal{O}_{sh}^{\varpi_{p_{\ell}}}(r_{\ell})) \ar[r,"m_{\bR}"] \ar[d,"\mathsf{rev}\circ(D\otimes \dots\otimes D)"'] & K_\C(\O_{sh}^\lambda(\bR)) \ar[d,"D"] \\
V(\varpi_{\bR^\ast}) \cong 	K_{\C}(\mathcal{O}_{sh}^{\varpi_{\smash{p_{\ell}^*}}}(-r_{\ell}+h))\otimes \dots\otimes K_{\C}(\mathcal{O}_{sh}^{\varpi_{\smash{p^*_1}}}(-r_1+h)) \ar[r,"m_{\bR^\ast}"] & K_\C(\O_{sh}^{\lambda^\ast}(\bR^\ast))
\end{tikzcd}}
\end{equation*}
where $\mathsf{rev}$ reverses the tensor product, and where $m_{\bR}$ and $m_{\bR^\ast}$ are the multiplication maps. Clearly, by Corollary \ref{cor:restriction_of_D_to_Via}, the left vertical map coincides with the map $\invomega_{\varpi_{\bR}}$ of \eqref{eq:omegabR}, and thus sends by Theorem \ref{th:basistobasis} the dual canonical basis $B(\varpi_{\bR})$ of the product $V(\varpi_{\bR})$ to the dual canonical basis $B(\varpi_{\bR^*})$ of $V(\varpi_{\bR^\ast})$. \medskip\par 
Now, by Corollary \ref{cor:multiplication_and_DC_basis}, the map $m_{\bR}$ sends vectors of $B(\varpi_{\bR})$ to classes of simple modules or $0$. Hence, using Proposition \ref{prop:mult_is_surj_truncation}, we deduce that there must exist an element $b\in B(\varpi_{\bR})$ with $m_{\bR}(b)=[L(\psi)]$, but then the commutativity of the diagram gives
$$ D([L(\psi)])=(D\circ m_{\bR})(b)=%
m_{\bR^*}(\omega_{\varpi_{\bR}}(b))\in m_{\bR^*}(B(\varpi_{\bR^*}))$$
(where we used the identification between the left vertical map and $\omega_{\varpi_{\bR}}$). In particular,~by Corollary \ref{cor:multiplication_and_DC_basis} again%
, the element $D([L(\psi)])$ is either the class of a simple module or~$0$,~but~it cannot be $0$ since $D$ is an isomorphism.
\end{proof}
A direct consequence of the above is the following corollary:
\begin{Corollary}\label{cor:dsetinv} The involution $D$ of $K_{\C}(\mathcal{O}_{sh})$ induces an algebra involution of $K_0(\mathcal{O}_{sh})$. Also, there exists an involution $d$ of the set $\mathfrak{r}$ such that
$D([L(\psi)])=[L(d(\psi))]%
$
for all $\psi\in \mathfrak{r}$.
\end{Corollary}
Our goal for the rest of this subsection is to understand the properties of the involution $d$ of Corollary \ref{cor:dsetinv}. First, we give the image of highest $\ell$-weights of chamber modules~via~$d$.
\begin{Lemma}\label{lem:image_of_psi_ia}
Fix $ (i,a) \in I \times_2 \Z $ and $w\in W$. Then  
$$ \smash{d(\sfPsi_{w\varpi_i,a})=\sfPsi_{ww_0\varpi_{i^\ast},-a+h}.}$$ 
In particular, 
$ d(\sfPsi_{i,a})=\sfPsi_{i,-a}^{-1}.$
\end{Lemma}
\begin{proof}
The restriction of $D$ to $V(\varpi_i,a)$ coincides with the morphism $\omega_i:{}^\omega V(\varpi_i)\to V(\varpi_{i^*})$.  Since $ \omega $ acts by $ -1 $ on the Cartan subalgebra, we see that $ \omega_{i} $ takes a vector of weight $ w \varpi_i $ to a vector of weight $ -w\varpi_i = ww_0\varpi_{i^\ast}$.  Thus we deduce from \eqref{eq:DVvarpiia} that
$$D([L_{w\varpi_i,a}])\in K_0(\mathcal{O}_{ww_0\varpi_{i^*}}^{\varpi_{i^*}}(-a+h)),$$
and the first statement follows from Theorem \ref{thm:DLpsi_is_still_the_class_of_a_simple} as $L_{ww_0\varpi_{i}^*,-a+h}$ is the unique irreducible object of $\smash{\mathcal{O}_{ww_0\varpi_{i^*}}^{\varpi_{i^*}}}(-a+h)$. The second statement is clear (see Example  \ref{ex:Inflsl3}).
\end{proof}

We now investigate the multiplicativity of $d$. We first easily deduce:

\begin{Lemma}\label{lem:multiplicative_when_product_simple}
Fix $\psi_1,\psi_2\in \fr$ for which the tensor product $L(\psi_1)\otimes L(\psi_2)$ is simple. Then $$d(\psi_1\psi_2)=d(\psi_1)d(\psi_2).$$
\end{Lemma}
\begin{proof}
Let $\psi=\psi_1\psi_2$. Then $L(\psi_1)\otimes L(\psi_2)\simeq L(\psi)$ and, since $D$ is an algebra isomorphism,
\begin{align*}
[L(d(\psi))]&=D([L(\psi)])=D([L(\psi_1)][L(\psi_2)])%
=[L(d(\psi_1))][L(d(\psi_2))].
\end{align*}
Hence $L(d(\psi_1))\otimes L(d(\psi_2))\simeq L(d(\psi))$, and thus $d(\psi_1)d(\psi_2)=d(\psi)$.
\end{proof}

In particular, $d(\psi^2)=d(\psi)^2$ if $\psi\in \fr$ is the highest $\ell$-weight of a real simple object~of~$\mathcal{O}_{sh}$
(see Theorem \ref{thm:chamber_modules_are_prime}). On the other extreme, the next result shows that $d$ is also multiplicative on $\ell$-weights %
that lie in generic product monomial crystals.
\begin{Theorem}\label{thm:multdgeneric}
Fix $\lambda_1,\lambda_2\in P_+$ with $\bR_1\in \Z^{\lambda_1}$ and $\bR_2\in \Z^{\lambda_2}$ such that the multiplication 
\begin{equation}\label{eq:multCrist12}
\cB(\lambda_1,\bR_1)\times \cB(\lambda_2,\bR_2)\to \cB(\lambda_1+\lambda_2,\bR_1\cup\bR_2)
\end{equation}
is bijective. Then $d(\psi_1\psi_2)=d(\psi_1)d(\psi_2)$ for all $\psi_1\in \cB(\lambda_1,\bR_1)$ and $\psi_2\in\cB(\lambda_2,\bR_2)$. 
\end{Theorem}
We will need the following direct consequence of the proof of Theorem \ref{thm:DLpsi_is_still_the_class_of_a_simple}:
\begin{Lemma}\label{lem:RestricDlambda} Fix $\la\in P_+$ with $\bR\in \Z^{\la}$. Then $D$ restricts to an isomorphism 
$$ \smash{K_{\C}(\mathcal{O}_{sh}^{\la}(\bR))\simeq K_{\C}(\mathcal{O}_{sh}^{\la^*}(\bR^*))},$$
where $\smash{\bR^*\in \Z^{\la^*}}$ is defined as in \eqref{eq:def_of_bRast}.
\end{Lemma}
\begin{proof}[Proof of Theorem \ref{thm:multdgeneric}]
Let $\la=\la_1+\la_2$ with $\bR=\bR_1\cup\bR_2$ and denote by
$$\bR_1^\ast\in \Z^{\la^*_1},\hspace*{0.75em} \bR_2^\ast\in \Z^{\la_2^*}\ \text{ and }\ \bR^\ast\in \Z^{\la^*}$$
the sets of parameters constructed using $\bR_1$, $\bR_2$ and $\bR$ as in \eqref{eq:def_of_bRast}. By Theorem \ref{th:descend} and Proposition \ref{prop:mult_is_surj_truncation}, our hypothesis about the bijectivity of \eqref{eq:multCrist12} is equivalent to the fact that multiplication induces a group isomorphism
$$ \smash{K_{0}(\mathcal{O}_{sh}^{\la_1}(\bR_1))\otimes K_{0}(\mathcal{O}_{sh}^{\la_2}(\bR_2))\simeq K_{0}(\mathcal{O}^{\la}_{sh}(\bR))},$$
but using $D$ with Lemma \ref{lem:RestricDlambda} shows that multiplication also gives a group isomorphism
$$\smash{ K_{0}(\mathcal{O}_{sh}^{\la_1^*}(\bR_1^*))\otimes K_{0}(\mathcal{O}_{sh}^{\la_2^*}(\bR_2^*))\simeq K_{0}(\mathcal{O}^{\la^*}_{sh}(\bR^*))},$$
and it follows easily (using for instance Corollary \ref{cor:tensor_product_decomp_in_Osh}) that the multiplication
$$\cB(\lambda_1^\ast,\bR_1^\ast)\times \cB(\lambda_2^\ast,\bR_2^\ast)\to \cB(\lambda^\ast,\bR^\ast)$$
is bijective. Fix enumerations $\{\psi_1,\dots,\psi_p\}= \cB(\la,\bR)$ and $\{\psi_1',\dots,\psi_p'\}= \cB(\la^*,\bR^*)$ that are compatible with Nakajima's partial order (as in the proof of Proposition \ref{prop:mult_is_surj_truncation}). Moreover, for each $k\in \{1,\dots, p\}$, denote by
\begin{equation*}
(\zeta_{1,k},\zeta_{2,k})\in \cB(\lambda_1,\bR_1)\times \cB(\lambda_2,\bR_2)\ \text{ and }\ (\zeta_{1,k}',\zeta_{2,k}')\in \cB(\lambda_1^\ast,\bR_1^\ast)\times \cB(\lambda_2^\ast,\bR_2^\ast)
\end{equation*}
the pairs such that $\psi_k=\zeta_{1,k}\zeta_{2,k}$ and $\psi_k'=\zeta_{1,k}'\zeta_{2,k}'$. Consider the commutative diagram
\begin{equation}\label{eq:DiagMultGenD}
\adjustbox{scale=1}{
\begin{tikzcd}[column sep = 4em, row sep = 2.75em]
K_0(\O^{\lambda_1}_{sh}(\bR_1))\otimes_\Z K_0(\O^{\lambda_2}_{sh}(\bR_2)) \arrow[r,"D\otimes D"]\arrow[d,"\text{mult}"'] &
K_0(\O^{\lambda_1^\ast}_{sh}(\bR_1^\ast))\otimes_\Z K_0(\O^{\lambda_2^\ast}_{sh}(\bR_2^\ast)) \arrow[d,"\text{mult}"] \\
K_0(\O^\lambda_{sh}(\bR))\arrow[r,"D"] & K_0(\O^{\lambda^\ast}_{sh}(\bR^\ast))
\end{tikzcd}}
\end{equation}
All arrows in the above diagram are isomorphisms. Also, our choice of enumeration~for $\cB(\la,\bR)$ and $\cB(\la^*,\bR^*)$ give ordered bases 
$$\smash{\{[L(\psi_1)],\dots,[L(\psi_p)]\}\subseteq K_0(\O_{sh}^{\la}(\bR))\ \text{ and }\ \{[L(\psi'_1)],\dots,[L(\psi'_p)]\}\subseteq K_0(\O_{sh}^{\la^*}(\bR^*))}.$$
Similarly, the enumerations $\{(\zeta_{1,1},\zeta_{2,1}),\dots,(\zeta_{1,p},\zeta_{2,p})\}%
$ and $\{(\zeta'_{1,1},\zeta'_{2,1}),\dots,(\zeta'_{1,p},\zeta'_{2,p})\}$ give ordered bases for the other $\Z$-modules appearing in \eqref{eq:DiagMultGenD}.

With all these $\Z$-bases in mind, we denote by $A$ and $B$ the matrices corresponding to the left and right vertical arrows in \eqref{eq:DiagMultGenD}, and write $P$ and $Q$ for the matrices associated~to~the maps $D$ and $D\otimes D$ (resp.). \par 
We claim that $P=Q$. Indeed, by commutativity of \eqref{eq:DiagMultGenD}, we have %
$$PA=BQ,$$
where $P,Q$ are permutation matrices by Theorem \ref{thm:DLpsi_is_still_the_class_of_a_simple}, and where $A,B$ are upper unitriangular by Corollary \ref{cor:tensor_product_decomp_in_Osh} (see also the proof of Proposition \ref{prop:mult_is_surj_truncation}). This implies our~claim.

We are now ready to prove the desired statement. Indeed, for $k\in \{1,\dots,p\}$, by definition of the matrix $Q=(q_{ij})_{1\leq i,j\leq p}$, 
\begin{equation*}
(D\otimes D)\big(\,[L(\zeta_{1,k})]\otimes [L(\zeta_{2,k})]\,\big)=\sum_{1\leq i\leq p} q_{ik} [L(\zeta_{1,i}')]\otimes [L(\zeta_{2,i}')]=[L(\zeta_{1,j}')]\otimes [L(\zeta_{2,j}')]
\end{equation*}
for some $j\in\{1,\dots,p\}$, i.e. $q_{ik}=\delta_{i,j}$. In addition, by definition of $P=(p_{ij})_{1\leq i,j\leq p}$,
\begin{equation*}
D([L(\psi_k)])=\sum_{1\leq i\leq p} p_{ik} [L(\psi_i')]=[L(\psi_{j'}')]
\end{equation*}
for some $j'\in \{1,\dots,p\}$, i.e. $p_{ik}=\delta_{i,j'}$. Hence, since $P=Q$, we must have $j=j'$, and~thus, by combining the two equations, we get 
\begin{align*} 
d(\zeta_{1,k})d(\zeta_{2,k}) = \zeta'_{1,j}\zeta'_{2,j}=\psi_j'=d(\psi_k)=d(\zeta_{1,k}\zeta_{2,k}),
\end{align*}
as desired. This concludes the proof.
\end{proof}

Although the above technique only applies when the multiplication map \eqref{eq:multCrist12} is bijective, we expect its conclusion to remain valid more generally. We therefore propose:

\begin{Conjecture}\label{conj:d_preserves_Nak}
The involution $d$ of $\mathfrak{r}$ is a group automorphism.
\end{Conjecture}

This conjecture is actually equivalent to a seemingly weaker statement:

\begin{Theorem} 
Conjecture \ref{conj:d_preserves_Nak} is equivalent to the condition
\begin{equation}\label{eq:dNakOrder}
\psi\preceq \psi' \Longrightarrow d(\psi)\preceq d(\psi')
\end{equation}
for all $\psi,\psi'\in \mathfrak{r}$ (where $\preceq$ is Nakajima's partial order).
\end{Theorem}

\begin{proof}
Assuming first Conjecture \ref{conj:d_preserves_Nak}, condition \eqref{eq:dNakOrder} is easily seen to hold since
\begin{align*} 
\textstyle d(\mathsf{A}_{i,a})&= \textstyle\frac{d(\sfPsi_{i,a-2})}{d(\sfPsi_{i,a+2})}\prod_{j\sim i}\frac{d(\sfPsi_{j,a+1})}{d(\sfPsi_{j,a-1})}=\textstyle\frac{\sfPsi_{i,-a-2}}{\sfPsi_{i,-a+2}}\prod_{j\sim i}\frac{\sfPsi_{j,-a+1}}{\sfPsi_{j,-a-1}}=\mathsf{A}_{i,-a}.
\end{align*}
Conversely, if this condition holds, fix $\psi_1,\psi_2\in \fr$.  Corollary \ref{cor:tensor_product_decomp_in_Osh} gives 
\begin{align*}
[L(d(\psi_1))\otimes L(d(\psi_2))]&=[L(d(\psi_1)d(\psi_2))]+\textstyle \sum_{\zeta\prec d(\psi_1)d(\psi_2)} n_\zeta [L(\zeta)]
\end{align*}
for some $n_{\zeta}$'s in $\mathbb{Z}_{\geq 0}$, but decomposing rather $[L(\psi_1)\otimes L(\psi_2)]$ and using $D$ shows instead
\begin{align*}
[L(d(\psi_1))\otimes L(d(\psi_2))]&=D([L(\psi_1)\otimes L(\psi_2)])=[L(d(\psi_1\psi_2))]+\textstyle \sum_{\vartheta\prec d(\psi_1\psi_2)} n_\vartheta [L(\vartheta)]
\end{align*}
for some $n_{\vartheta}$'s in $\Z_{\geq 0}$. This clearly implies $d(\psi_1\psi_2)=d(\psi_1)d(\psi_2)$ using \eqref{eq:dNakOrder}.
\end{proof}

In an upcoming paper \cite{otherpaper}, we will prove Conjecture \ref{conj:d_preserves_Nak} by establishing a compatibility result between the partial order introduced in \eqref{eq:def_of_Lusztigs_partial_order} and Nakajima's partial order (see Lemma \ref{lem:bijection_posets}). We also provide evidence for Conjecture \ref{conj:d_preserves_Nak} in Section \ref{sec:HLD}.

We now conclude this subsection with results which were foreshadowed in Section \ref{sec:TensO}.

Let $W=L(\frac{\mathbf{P}}{\mathbf{Q}})$ for $\mathbf{P},\mathbf{Q}\in \mathfrak{r}$ polynomial and fix a simple object $V$ in $\mathcal{O}_{sh}$. Then Corollary \ref{cor:HZgenSimple} shows that the product $V\otimes W$ is of highest $\ell$-weight if either one of the products
\begin{center}
$V\otimes L(\mathbf{P})\ \text{ or }\ L(\mathbf{Q}^{-1})\otimes V$
\end{center}
is simple. Moreover, by Theorem \ref{thm:HZcriterionTensSimpPol}, the condition that $V\otimes L(\mathbf{P})$ is simple is equivalent to a combinatorial condition on the normalized $\ell$-character of $V$ and is thus somewhat ``easy to use'' in practice. However, the product $L(\mathbf{Q}^{-1})\otimes V$ always satisfies
$$ \GKdim (L(\mathbf{Q}^{-1})\otimes V)\geq \GKdim (V\otimes W),$$
and is hence harder to study. One can nevertheless use the involution $D$ introduced above to again reduce the question of whether such a product is simple to a combinatorial~condition. We give this condition below along with the dual version of Corollary~\ref{cor:BigTensProd}~%
mentioned\footnote{Notice that this result and Proposition \ref{prop:CombCritNeg} also work for the non-integral version of the category~$\mathcal{O}_{sh}$ studied in the first sections of this paper since Theorem \ref{thm:DLpsi_is_still_the_class_of_a_simple} can be easily extended to this setting.} in Section \ref{sec:TensorResults}. (Note that the latter dual version can be proven exactly as Corollary \ref{cor:BigTensProd}).
\begin{Proposition}\label{prop:CombCritNeg}
Fix $V=L(\psi)$ simple in $\mathcal{O}_{sh}$ with $\mathbf{Q}\in \mathfrak{r}$ polynomial. Write
$$ \mathbf{Q}=\mathsf{\Psi}_{i_1,a_1}\dots\mathsf{\Psi}_{i_k,a_k}$$
for some $(i_1,a_1),\dots,(i_k,a_k)\in I\times_2\Z$. Then the following statements are equivalent:\smallskip
\begin{enumerate}
\item $L(\mathbf{Q}^{-1})\otimes V$ is simple, and
\item the normalized $\ell$-character of $L(d(\psi))$ does not involve the variables $\{A_{i_r,-a_r}^{-1}\}_{1\leq r\leq k}$.
\end{enumerate}
\end{Proposition}
\begin{proof}
Note that $d(\mathbf{Q}^{-1}) = \mathsf{\Psi}_{i_1,-a_1}\dots \mathsf{\Psi}_{i_k,-a_k}$ by Lemma \ref{lem:multiplicative_when_product_simple} and Proposition \ref{prop:extremalmodulesReal}. Thus, by Theorem %
\ref{thm:DLpsi_is_still_the_class_of_a_simple} and Corollary \ref{cor:tensor_product_decomp_in_Osh},
\begin{align*}
L(\mathbf{Q}^{-1})\otimes V \text{ is simple}%
&\iff [L(\mathbf{Q}^{-1})][V]=[L(\psi\mathbf{Q}^{-1})]\text{ in }K_0(\mathcal{O}_{sh})\\
&\iff %
[L(d(\mathbf{Q}^{-1}))][L(d(\psi))]=[L(d(\psi\mathbf{Q}^{-1}))]
\text{ in }K_0(\mathcal{O}_{sh})\\
&\iff L(d(\psi))\otimes L(d(\mathbf{Q}^{-1})) \text{ is simple},
\end{align*}
and using Theorem \ref{thm:HZcriterionTensSimpPol} ends the proof.
\end{proof}
\begin{Corollary}
For $1\leq r\leq k$, fix $\mathbf{P}_r,\mathbf{Q}_r\in \mathfrak{r}$ polynomial and let $V_r=L(\frac{\mathbf{P}_r}{\mathbf{Q}_r})$. Suppose 
that Conjecture \ref{conj:Associators} holds. Then (omitting parentheses)
\begin{enumerate}
\item[(i)] $V_1\otimes \dots \otimes V_k$ is of highest $\ell$-weight if $L(\mathbf{Q}_r^{-1})\otimes V_s$ is simple for $1\leq r<s\leq k$, and
\item[(ii)] $V_1\otimes \dots \otimes V_k$ is of co-highest $\ell$-weight if $L(\mathbf{Q}_r^{-1})\otimes V_s$ is simple for $1\leq s<r\leq k$.
\end{enumerate}
In particular,
\begin{enumerate}\setcounter{enumi}{2}
\item[(iii)] $V_1\otimes \dots \otimes V_k$ is simple if $ L(\mathbf{Q}_r^{-1})\otimes V_s$ is simple for all $1\leq r,s\leq k$ with $r\neq s$.
\end{enumerate}
\end{Corollary}
\subsection{Relation to Hernandez--Leclerc's duality}\label{sec:HLD} As mentioned in the beginning of~this section, the involution $D$ of $K_0(\mathcal{O}_{sh})$ constructed here is closely related to the eponymous map in \cite{hernandez2016cluster}. More precisely, define $\mathcal{O}^{+}_{sh}$ (resp.~$\mathcal{O}^{-}_{sh}$) as the Serre subcategory of~$\mathcal{O}_{sh}$~for which simple objects have highest $\ell$-weights in the submonoid of~$\mathfrak{r}$ generated by the~$\mathsf{\Psi}_{i,a}$'s (resp.~the $\smash{\mathsf{\Psi}_{i,a}^{-1}}$'s) and the $\mathsf{Y}_{i,a}$'s. Also, denote by $\smash{K_0^+}\subseteq \smash{K_0(\mathcal{O}_{sh}^+)}$ (resp.~$\smash{K_0^-}\subseteq \smash{K_0(\mathcal{O}_{sh}^-)}$)~the subalgebra generated by the classes in $K_0$ of the $L^+_{i,a}$'s (resp.~the $L^-_{i,a}$'s) and the $L(\mathsf{Y}_{i,a})$'s. Using the results of \cite[Section 8.6]{hernandez2023representations} together with \cite[Corollary 1.2.1]{varagnolo2025representations}, we~can adapt \cite[Proposition 5.13 and Theorem 5.17]{hernandez2016cluster} to our context and obtain the following:
\begin{Theorem}[\cite{hernandez2016cluster}] The assignment $[L_{i,a}^+]\mapsto [L_{i,-a}^-]$ induces a unique algebra isomorphism $K_0^+\to K_0^-$, which in turn extends uniquely to an algebra isomorphism 
\begin{equation}\label{eq:HLD} 
\smash{D_{\normalfont\text{HL}}: K_0(\mathcal{O}^+_{sh})\to K_0(\mathcal{O}_{sh}^-).}
\end{equation}
\end{Theorem}
In particular, by Lemma \ref{lem:image_of_psi_ia}, our involution $D$ of $K_0(\mathcal{O}_{sh})$ restricts to the isomorphism \eqref{eq:HLD} (which we call \textit{Hernandez--Leclerc's duality}). This justifies the title of this section. 
\begin{Rem} The above can seem surprising, but actually follows quite easily from the \textit{Baxter $TQ$-relations} of \cite{frenkel2015baxter}. We briefly sketch this approach here for $\fg=\mathfrak{sl}_2$ and~refer~to \cite[Section 5.3]{hernandez2016cluster} for the general case. By Example \ref{ex:lchar}, we have 
$$ \smash{\chi_{\ell}(L(\mathsf{Y}_{1,a}))= \mathsf{Y}_{1,a}+\mathsf{Y}_{1,a}\mathsf{A}_{1,a}^{-1} = \tfrac{\sfPsi_{1,a-2}}{\sfPsi_{1,a}}+\tfrac{\sfPsi_{1,a+2}}{\sfPsi_{1,a}},}$$
which implies that, in $K_0(\mathcal{O}_{sh}^+)$, 
\begin{equation}\label{eq:TQsl2}
{[L_{1,a}^+][L(\mathsf{Y}_{1,a})]=[L_{1,a-2}^+]+[L_{1,a+2}^+].}
\end{equation}
Also, the results of \cite{frenkel2015baxter} give $[L_{1,a-2}^-][L(\mathsf{Y}_{1,a})]=[L_{1,a-4}^-]+[L_{1,a}^-]$
in $K_0(\mathcal{O}_{sh}^-)$, which, after applying $D$, implies
$$ {D([L_{1,a-2}^-][L(\mathsf{Y}_{1,a})]) = [L_{1,2-a}^+]D([L(\mathsf{Y}_{1,a})])=[L^+_{1,4-a}]+[L^+_{1,-a}]=D([L_{1,a-4}^-]+[L_{1,a}^-])}.$$
Using \eqref{eq:TQsl2} with the map $\chi_{\ell}$ thus gives $\sfPsi_{1,2-a}\chi_{\ell}(D([L(Y_{1,a})]))=\sfPsi_{1,2-a}\chi_{\ell}(L(Y_{1,2-a}))$ from which it follows that $D([L(\mathsf{Y}_{1,a})])=[L(\mathsf{Y}_{1,2-a})]$ in $K_0(\mathcal{O}_{sh})$ (i.e.~the image via $D$ of $L(\mathsf{Y}_{1,a})$ is uniquely determined by the image of the positive or negative prefundamental classes).
\end{Rem}

Combining  Theorem \ref{thm:DLpsi_is_still_the_class_of_a_simple} with the observation that our involution $D$ extends Hernandez--Leclerc's duality $D_{\text{HL}}$ answers a question of \cite[Appendix A]{pinet2024functor}. Moreover, the above fact allows us to deduce the following properties of the involution $D$ from the results of \cite{hernandez2016cluster}. Note that these properties give good evidence toward the validity of Conjecture \ref{conj:d_preserves_Nak}.

\begin{Proposition} The involution $d$ of the set $\mathfrak{r}$ restricts to monoid isomorphisms
$$ \smash{d:\mathfrak{r}^{\pm}\to \mathfrak{r}^{\mp},}$$
with $\mathfrak{r}^+,\mathfrak{r}^-\subseteq \mathfrak{r}$ the submonoids underlying the definition of the categories $\mathcal{O}_{sh}^+,\mathcal{O}_{sh}^-\subseteq \mathcal{O}_{sh}$. Also, $d(\mathsf{Y}_{i,a})=\mathsf{Y}_{i,{2-a}}$ for all $(i,a)\in I\times_2\Z$, and $D$ restricts to an involution of $K_0(Y\fmod)$ where $Y\fmod$ is the category of finite-dimensional modules over $Y=Y_0$.
\end{Proposition}

\begin{Rem} The fact that $D$ restricts to an involution of $K_0(Y\fmod)\subseteq K_0(\mathcal{O}_{sh})$ can be seen geometrically. Indeed, the involution $F$ of $\smash{\widehat{\mathcal{Z}}_{\infty}^\circ}$ given in Section \ref{subsec:invBottSam} is easily~seen~to send $G$-invariant elements to $G$-invariant elements (since $\varpi_i^*$ gives a $G$-equivariant isomorphism from $V(\varpi_{i^*})^*$ to ${}^\omega (V(\varpi_i)^*)$). Therefore, the involution $D$ of $\smash{\C[\widehat{\mathcal{Z}}_{\infty}^{\circ}]\simeq \cR\simeq K_{\C}(\mathcal{O}_{sh})}$ induces an involution of the $G$-invariant subalgebra (see Corollary \ref{cor:Invariants})
$$\smash{{}^{G}K_{\C}(\mathcal{O}_{sh})}\simeq \smash{K_{\C}(Y\fmod)}.$$
Note also that $\mathcal{O}^+_{sh}$ is equal to the full subcategory $\mathscr{C}_{sh}\subseteq \mathcal{O}_{sh}$ of finite-dimensional modules \cite{hernandez2024shifted}, and \eqref{eq:HLD} can thus be obtained from $D$ by restricting the domain to $N$-invariants.
\end{Rem}
Finally, we ask a question that is inspired by the work of the fourth author \cite{pinet2024functor}:
\begin{Question} Is there an involutive autofunctor $\mathscr{D}$ of $\mathcal{O}_{sh}$ that categorifies %
$D$%
?%
\end{Question}

\newpage
\appendix
\section{Pro-varieties}\label{app:pro_varieties}

The algebro-geometric objects that 
we consider in Section \ref{sec:bi-infinite-bott-samelson} are inverse limits of finite-dimensional algebraic varieties. Such %
limits do not always exist in the category of schemes, so we need to consider them as formal inverse limits in the pro-completion of the category of varieties. In this appendix, we set up some basic definitions regarding these pro-varieties as well as coherent sheaves on them.

Let $\kk$ be a field. By a variety over $\kk$, we mean a scheme of finite type over $\kk$. Let $\Var_{\kk}$ denote the category of varieties over $\kk$.

\begin{Def}
    The category of pro-varieties over $\kk$ is the pro-completion $\Pro(\Var_{\kk})$ of the category of $\kk$-varieties. Specifically, the objects are formal inverse limits \[\displaystyle\lim_{\substack{\leftarrow \\ i\in \cI}} X_i\] associated to diagrams $\cI \to \Var$ indexed by a cofiltered category $\cI$, and the morphisms~are %
    \[ \Hom_{\Pro(\Var_\kk)}\Bigg(\displaystyle\lim_{\substack{\leftarrow \\ i\in \cI}} X_i, \displaystyle\lim_{\substack{\leftarrow \\ j\in \cJ}} Y_j\Bigg) = \lim_{\substack{\leftarrow \\ j\in \cJ}} {\lim_{\substack{\rightarrow \\ i\in \cI}}} \Hom_{\Var_\kk}(X_i, Y_j)\]
\end{Def}

There is a fully faithful embedding of $\Var_\kk$ into $\Pro(\Var_\kk)$ given by considering trivial diagrams with a single object, and every object $X$ of $\Var_\kk$ is cocompact in $\Pro(\Var_\kk)$, meaning that $\Hom(-,X): \Pro(\Var_\kk)^{\op{op}} \to \op{Set}$ sends cofiltered limits in $\Pro(\Var_\kk)$ to filtered colimits in $\op{Set}$. Moreover, the pro-completion satisfies the following universal property: for every category $\mathcal{C}$ that admits cofiltered limits and every fully faithful functor $F:\Var_\kk\to \mathcal{C}$ that lands in cocompact objects, there is a unique (up to unique isomorphism) extension~of $F$ to a cofiltered limit-preserving functor $F': \Pro(\Var_\kk)\to \mathcal{C}$. Specifically, $F'$ maps~a~formal inverse limit $\displaystyle\lim_{\leftarrow} X_i$ %
to the inverse limit in $\mathcal{C}$ of the diagram $(F(X_i))_{i\in \cI}$.

\begin{Def}
Let \[X=\displaystyle\lim_{\substack{\leftarrow \\ i\in \cI}} X_i\] be a pro-variety (i.e.~an object of $\Pro(\Var_{\kk})$). The category of coherent sheaves $\op{Coh}(X)$~on $X$ is the direct limit of categories $\displaystyle\lim_{\rightarrow} \op{Coh}(X_i)$. 
    Concretely, this means that:
    \begin{itemize}
        \item objects of $\op{Coh}(X)$ are all of the form $\pi_i^*\cF$, where $\cF$ is a coherent sheaf on some $X_i$ and $\pi_i$ is the projection $X\to X_i$;
        \item given $\cF\in \Oh(X_i)$ and $\cG \in \op{Coh}(X_{i'})$, morphisms between $\pi_i^*\cF$ and $\pi_{i'}^* \cG$ are all of the form $\pi^*_j\varphi$, where $\varphi$ is a morphism between the pullbacks of $\cF$ and $\cG$ to $X_j$, for some $j\in \cI$ mapping to both $i$ and $i'$; and
        \item morphisms $\pi_j^*\varphi$ and $\pi_{j'}^*\psi$ are identified if and only if their pullbacks to $X_k$ agree~for some $k\in \cI$ mapping to both $j$ and $j'$. 
    \end{itemize}
\end{Def}

In fact, $\Oh$ is a $2$-functor from $\Pro(\Var_{\kk})^{\op{op}}$ to the bicategory of categories, meaning~that a morphism $f:X\to Y$ of pro-varieties always induces a pullback %
$f^*: \Oh(Y)\to \Oh(X)$, compatibly with composition of morphisms. Thus, isomorphisms of pro-varieties give equivalences of categories of coherent sheaves. 

By a line bundle on a pro-variety $X=\displaystyle\lim_{\leftarrow} X_i$, we mean a coherent sheaf of the form~$\pi_i^* \mathcal{L}$, where $\mathcal{L}$ is a line bundle on $X_i$. If $\Pic(X)$ denotes the group of isomorphism classes of line bundles on $X$, then it follows immediately from the definitions that 
$$\Pic(\displaystyle\lim_{\leftarrow} X_i)=\displaystyle\lim_{\rightarrow} \Pic(X_i).$$

\begin{Def}
Let \[X=\displaystyle\lim_{\substack{\leftarrow \\ i\in \cI}} X_i\] be a pro-variety and $\cF$ a coherent sheaf on $X_i$ for some $i$. The space of global sections of the coherent sheaf $\pi_i^*\cF$ on $X$ is defined as
    \[ H^0(X, \pi_i^*\cF)  = \lim_{\substack{\rightarrow \\ (j,\alpha)%
    }} H^0(X_j, \pi_{\alpha}^*\cF),\]
where the direct limit is over all pairs $(j,\alpha)\in \cI/i\times \Hom_{\cI}(j,i)$, and where $\pi_{\alpha}:X_j\to X_i$~is the transition map corresponding to $\alpha$.
\end{Def}

In other words, every global section is defined at some finite level, and two global sections agree on the inverse limit if and only if they agree at some finite level. Again, this definition is functorial, in the sense that a morphism $\cF\to \cG$ of coherent sheaves on $X$ induces a map $H^0(X, \cF)\to H^0(X,\cG)$ and that a morphism $f:X\to Y$ of pro-varieties induces a map $H^0(Y,\cF)\to H^0(X, f^*\cF)$ for every $\cF\in \Oh(Y)$, compatibly with compositions.

\begin{Rem} \label{rem:affine-transition-inverse-limits}
    Let $\op{Sch}_\kk$ be the category of quasi-compact quasi-separated schemes over~$\kk$. By \cite[\href{https://stacks.math.columbia.edu/tag/01YX}{Tag 01YX}]{stacks-project}, $\op{Sch}_\kk$ admits limits of cofiltered diagrams with \emph{affine} transition~morphisms\footnote{The reference we give assumes that the indexing category $\cI$ is (the dual of) a directed set rather than a cofiltered category. This restriction does not come with any loss of generality%
    , thanks to \cite[\href{https://stacks.math.columbia.edu/tag/0032}{Tag 0032}]{stacks-project}.}. Moreover, if $X$ is a variety, then $\Hom(-,X)$ sends such limits to colimits in $\op{Set}$ by \cite[\href{https://stacks.math.columbia.edu/tag/01ZC}{Tag 01ZC}]{stacks-project}.
    
    If we let $\Pro(\Var_\kk)_{\text{aff.}}$ be the full subcategory of $\Pro(\Var_\kk)$ consisting of the pro-varieties obtained as the limit of a diagram with affine transition morphisms, it follows that we have a fully faithful embedding $\Pro(\Var_\kk)_{\text{aff.}}\hookrightarrow \op{Sch}_{\kk}$ sending a formal inverse limit $\displaystyle\lim_{\leftarrow} X_i$ to the actual inverse limit in $\op{Sch}_\kk$. We can therefore think of such pro-varieties as actual schemes. Furthermore, by \cite[\href{https://stacks.math.columbia.edu/tag/01ZR}{Tag 01ZR}, \href{https://stacks.math.columbia.edu/tag/0B8W}{Tag 0B8W}, \href{https://stacks.math.columbia.edu/tag/01Z0}{Tag 01Z0}]{stacks-project}, our notions of coherent sheaves, line bundles and global sections agree, in this situation, with the usual notions for schemes.

    Without the assumption on affine transition morphisms, however, pro-varieties generally cannot be regarded as schemes. See for example \cite[\href{https://stacks.math.columbia.edu/tag/078E}{Tag 078E}]{stacks-project}, which explains why the product of infinitely many copies of $\PP^1$ (which is the same as our bi-infinite Bott--Samelson variety in type $A_1$) does not exist in the category of schemes.
\end{Rem}
				
\section{Inflations}\label{sec:Inflations}
Inflations were introduced in \cite{pinetinflations} as distinguished preimages for canonical restriction functors arising in the study of representations of shifted quantum affine algebras.~We~define their rational analogue, for representations of shifted Yangians, here and use truncations~to prove counterparts of conjectures/results of \cite{pinetinflations}. We use the notation of Sections \ref{sec:Crystal}--\ref{sec:TensorResults}.
\subsection{Definition}
Take $J\subseteq I$ and let $\fg_J\subseteq \fg$ be the Lie subalgebra of $\fg$ with Cartan matrix $(C_{i,j})_{i,j\in J}$. Set also
$\textstyle P_J^{\vee} = \bigoplus_{j\in J} \Z \varpi_j^{\vee}$ and let $\res_J:P^{\vee}\twoheadrightarrow P^{\vee}_J$ be the projection defined by 
$$ \res_J(\varpi_i^{\vee})=\left\{\begin{array}{ll}
\varpi_i^{\vee} & \text{if }i\in J,\\
0 & \text{else.}	
\end{array}\right. $$
For $\nu\in P_J^{\vee}$, denote by $Y_{\nu}(\fg_J)$ the shifted Yangian associated to $\nu$ and $\fg_J$, that is~the~algebra with generators $\{e_{j,q},f_{j,q},h_{j,p}\,|\,j\in J,\,q\in \Z_{>0},\,p\in \Z\}$ and relations \eqref{H,H}--\eqref{symF}. Then $Y_{\nu}(\fg_J)$ admits a category $\mathcal{O}_{\nu}(\fg_J)$ (as in Section \ref{sec:Osh}) with a triangular decomposition (as in \eqref{eq:TriangularDec}). Moreover, given $\mu\in P^{\vee}$ such that $\res_J(\mu)=\nu$, there is a canonical algebra map\footnote{In this section, notations and concepts which are not explicitly associated with the Lie algebra $\fg_J$ refer to the Lie algebra $\fg$ (and thus correspond to the notations used elsewhere in the paper).}
$$\iota_{\mu}:Y_{\nu}(\fg_J)\to Y_{\mu}$$
that sends generators of $Y_{\nu}(\fg_J)$ to generators with the same label in $Y_{\mu}$.

\begin{Lemma} 
Let $\mu\in P^{\vee}$ with $\nu=\res_J(\mu)$. Then the map $\iota_{\mu}:Y_{\nu}(\fg_J)\rightarrow Y_{\mu}$ is injective.
\end{Lemma}

\begin{proof}
This is shown using triangular decompositions %
as in \cite[Proposition 2.4]{pinetinflations}.
\end{proof}

We denote by $\res_J^{\mu}$ the functor from $\mathcal{O}_{\mu}$ to $\mathcal{O}_{\nu}(\fg_J)$ given by pullback with respect to $\iota_{\mu}$. The following definition is an adaptation of \cite[Definition 3.1]{pinetinflations}.
\begin{Def} Fix an %
object $W$ in $\mathcal{O}_{\nu}(\fg_J)$ with a coweight $\mu\in P^{\vee}$ such that $\res_J(\mu)=\nu$. Then, a \textit{$J$-inflation of $W$ of coweight $\mu$ to $\fg$} (or just \textit{inflation of $W$ to $\fg$} if $\mu$ and $J$~are~clear from the context) is an object $V$ of $\mathcal{O}_{\mu}$ satisfying the conditions
\begin{itemize}
\item[(i)] $\res_J^{\mu}(V)\simeq W$ as $Y_{\nu}(\fg_J)$-modules, and
\item[(ii)] $e_i(u)V=f_i(u)V=0$ whenever $i\not\in J$.
\end{itemize}
\end{Def} 
As explained in \cite[Section 4.2]{pinetinflations}, the notion of inflation allows one to reduce problems in $\mathcal{O}_{sh}$ to analogous (typically easier) problems in the category $\mathcal{O}_{sh}(\fg_J)=\bigoplus_{\nu\in P^{\vee}_J}\mathcal{O}_{\nu}(\fg_J)$. In particular, inflations are compatible with the notion of \textit{real} and \textit{prime modules}, and are expected to have natural applications in the setting of monoidal categorifications of cluster algebras. A natural (and very basic) question about these is however:

\begin{Question}\label{question:Infl} 
Do all simple modules in $\mathcal{O}_{sh}(\fg_J)$ admit an inflation to $\fg$?
\end{Question}

This question was answered positively for shifted quantum affine algebras of type A--B--G (or for finite-dimensional irreducible modules of arbitrary shifted quantum affine algebras) in \cite[Corollary 3.29 and Theorem 3.40]{pinetinflations} using: 
\begin{itemize}
\item[(1)] the compatibility between inflations and \textit{fusion/tensor products} with 
\item[(2)] a technical study of $\ell$-characters of negative prefundamental representations.
\end{itemize}
The purpose of this section is to answer Question \ref{question:Infl} in full generality, therefore proving~the rational analogue of \cite[Conjecture 1.4]{pinetinflations}. We use for this a different approach from the one used in \cite{pinetinflations}, and rely on notable equalities involving truncations for shifted Yangians of the form $Y_{\mu}$ and $Y_{\nu}(\fg_J)$ (as in \cite[Section 5.2]{kamnitzer2018reducedness}). We think that our approach could also be adapted to the context of shifted quantum affine algebras, but leave this adaptation for further work (see \cite[Section 4.1]{pinetinflations} for details%
).
\subsection{Existence} Fix $\la,\nu \in P_J^{\vee}$ such that $\la\in P_+^{\vee}$ and set 
$$ \textstyle \la-\nu = \sum_{j\in J} m_j\beta_j^{\vee},$$
with $\{\beta_j^{\vee}=\res_J(\alpha_j^{\vee})\}_{j\in J}$ the simple coroots of $\fg_J$. Set also $m_i=0$ for $i\not\in J$ and let 
$$\textstyle \mu = \la-\sum_{j\in J}m_j\alpha_j^{\vee} = \nu + \sum_{j\in J}m_j\sum_{\substack{i\not\in J;\,i\sim j}}\varpi_i^{\vee}.$$
Finally, fix
$\bR \in \C^{\la}$ %
  and consider the diagram %
{\footnotesize
\begin{center}
$\begin{tikzcd}[column sep = 5em, row sep = 1.5em]
Y_{\mu%
} \ar[r,"\Phi_{\mu%
}^{\la%
}(\bR%
)"] & \widetilde{\mathscr{A}}_{\la-\mu}=:\mathscr{A}\\ 
Y_{\nu}(\fg_J) \ar[u, "\iota_{\mu}", hook] \ar[ur,"\Phi_{\nu}^{\la}(\bR{,}\fg_J)",swap]
\end{tikzcd}$
\end{center}}\noindent
with $\Phi_{\mu%
}^{\la%
}(\bR%
)$ and $\Phi_{\nu}^{\la}(\bR,\fg_J)$ the algebra morphisms of Theorem \ref{GKLO homomorphism} (and where we used the fact that the algebra $\mathscr{A}$ of Definition \ref{def:DiffOpAlg} depends only on the gauge parameters $(m_i)_{i\in I}$,~and not on the underlying Lie algebra). %
Note that, for $i\not\in J$,
$$ (\Phi_{\mu}^{\la}(\bR))(e_i(u))=(\Phi_{\mu}^{\la}(\bR))(f_i(u))=0$$
with
\begin{equation}\label{eq:GKLOA} 
(\Phi_{\mu}^{\la}(\bR))(a_i(u))%
=1.
\end{equation}
Moreover, denoting by $\{a'_j(u)=1+\sum_{q\geq 1}a'_{j,q}u^{-q}\}_{j\in J}$ the currents of $Y_{\nu}(\fg_J)$ defined in~\eqref{eq: def of A gens}, we have, for $j\in J$, 
\begin{equation}\label{eq:GKLOAJ}
(\Phi_{\mu}^{\la}(\bR))(a_j(u))=u^{-m_j}W_j(u)=(\Phi_{\nu}^{\la}(\bR,\fg_J))(a_j'(u))
\end{equation}
so that 
$$\textstyle (\Phi_{\mu}^{\la}(\bR)\circ \iota_{\mu})(h_j(u)) %
= p_{R_j}(u)\frac{\prod_{k\in J;\,k\sim j}W_k(u-1)}{W_j(u)W_j(u-2)} = (\Phi_{\nu}^{\la}(\bR,\fg_J))(h_j(u)).$$
The lemma below easily follows from the above computations\footnote{Note that $\iota_{\mu}(a_j'(u))\neq a_j(u)$ if $j\in J$, but $(\Phi_{\mu}^{\la}(\bR)\circ \iota_{\mu})(a_j'(u))=(\Phi_{\nu}^{\la}(\bR,\fg_J))(a_j'(u))=(\Phi_{\mu}^{\la}(\bR))(a_j(u))$.}. 
\begin{Lemma}\label{lemma:EqualityTrunc} The above diagram commutes, that is $\Phi_{\mu%
}^{\la%
}(\bR%
)\circ \iota_{\mu}=\Phi_{\nu}^{\la}(\bR,\fg_J)$. Furthermore, 
\begin{center}
$Y_{\mu%
}^{\la%
}(\bR%
)=\op{Im}\Phi_{\mu}^{\la}(\bR)$ and $Y_{\nu}^{\la}(\bR,\fg_J)=\op{Im}\Phi_{\nu}^{\la}(\bR,\fg_J)$ 
\end{center}
coincide (as subalgebras of $\mathscr{A}$).
\end{Lemma}
Take $V_J$ in $\mathcal{O}_{\nu}^{\la}(\bR,\fg_J)$. Then the pullback $V$ of $V_J$ by $\Phi_{\mu}^{\la}(\bR):Y_{\mu}\twoheadrightarrow Y_{\mu}^{\la}(\bR)=Y_{\nu}^{\la}(\bR,\fg_J)$ is clearly a $Y_{\mu}$-module~for which $\res_J^{\mu}(V)\simeq V_J$ (as $Y_{\nu}(\fg_J)$-modules) and such that
$$ e_i(u)V=f_i(u)V=(a_i(u)-1)V=0$$
for $i\not\in J$. To show that $V$ is an inflation of $V_J$, it thus suffices to prove that $V$ belongs to~$\mathcal{O}_{\mu}$ (i.e.~that the equality $Y_{\mu}^{\la}(\bR)=Y_{\nu}^{\la}(\bR,\fg_J)$ identifies $\mathcal{O}_{\mu}^{\la}(\bR)$ with $\mathcal{O}_{\nu}^{\la}(\bR,\fg_J)$). For~this goal, note first that using \eqref{eq: def of A gens} easily gives, for all $i\in I$,
\begin{equation}\label{eq:HA1}
\textstyle h_{i,1-\langle \mu,\alpha_i\rangle} = \sum_{k\in I}c_{ik}(m_k-a_{k,1})-\sum_{c\in R_i}c,
\end{equation}
with, reciprocally 
\begin{equation}\label{eq:AH1}
\textstyle a_{i,1} = m_i-\sum_{k\in I}(C^{-1})_{ik}(h_{k,1-\langle \mu,\alpha_k\rangle}+\sum_{c\in R_k}c),
\end{equation}
where $C^{-1}$ is the inverse of the Cartan matrix of $\fg$.
Let now $\mathfrak{h}_J$ be the Cartan subalgebra of $\fg_J$. Then, given a weight $\omega\in \mathfrak{h}_J^*$ of the $Y_{\nu}(\fg_J)$-module $V_J$, the weight-space
\begin{align*}
(V_J)_{\omega} &= \{v\in V_J\,|\, \exists p\in \mathbb{N} \text{ such that } (h_{j,1-\langle \mu,\alpha_j\rangle}-2\langle\beta_j^{\vee},\omega\rangle)^pv=0 \text{ for all } j\in J\}\subseteq V_J
\end{align*}
can be equivalently described as a simultaneous generalized eigenspace for the action of~the elements $a_{j,1}'\in Y_{\nu}(\fg_J)$ (with $j\in J$) or, because of \eqref{eq:GKLOA}--\eqref{eq:GKLOAJ}, as a simultaneous generalized eigenspace in the pullback $V=(\Phi_{\mu}^{\la}(\bR))^*(V_J)$ for the action of the elements $a_{i,1}\in Y_{\mu}$ (with $i\in I$). Thus, putting everything together and tracking eigenvalues using \eqref{eq:HA1}--\eqref{eq:AH1}, we get 
$$ (V_J)_{\omega} = V_{\iota(\omega)}=\{v\in V\,|\, \exists p\in \mathbb{N} \text{ such that } (h_{i,1-\langle \mu,\alpha_i\rangle}-2\langle\alpha_i^{\vee},\iota(\omega)\rangle)^pv=0 \text{ for all } i\in I\}$$
where the weight $\iota(\omega)\in \mathfrak{h}^*$ is defined by 
\begin{equation}\label{eq:iotaInfl}
\textstyle \iota(\omega) = \infl_J(\omega-\sum_{j\in J}\sum_{c\in R_j}c\varpi_j)+\sum_{j\in J}\sum_{c\in R_j}c\varpi_j,
\end{equation}
with $\infl_J:\mathfrak{h}_J^*\to\mathfrak{h}^*$ the linear map given by $\infl_J(\beta_j)=\alpha_j$ for $j\in J$. In particular, since $V_J$ lies in $\mathcal{O}_{\nu}^{\la}(\bR,\fg_J)\subseteq \mathcal{O}_{\nu}(\fg_J)$, we get a decomposition 
$$\textstyle V = V_J\simeq \bigoplus_{\omega\in \mathfrak{h}_J^*}(V_J)_{\omega} = \bigoplus_{\omega\in\iota(\mathfrak{h}_J^*)}V_{\omega}$$
where the summands are all finite-dimensional vector spaces. Hence conditions $(\mathcal{O}1)$--$(\mathcal{O}2)$ of Definition \ref{def:Omu} hold for $V$. For condition $(\mathcal{O}3)$, note that, for $\omega,\omega'\in \mathfrak{h}_J^*$,
$$ \textstyle \omega -\omega'\in Q_{J,+}:=\sum_{j\in J}\Z_{\geq 0}\beta_j\ \Longrightarrow\ \iota(\omega)-\iota(\omega')=\infl_J(\omega-\omega')\in \infl_J(Q_{J,+})\subseteq Q_+$$
and therefore condition $(\mathcal{O}3)$ for $V$ follows easily from the same condition for $V_J$. This~(with Theorem \ref{thm:truncoprod}, Corollary \ref{cor:highestTrunc} and Corollary \ref{cor:cohighestTrunc}), finishes the proof of the following result:
\begin{Corollary}\label{cor:all_objects_inflatable}
All objects in the subcategory $\mathcal{O}_{\nu}^{\la}(\bR,\fg_J)\subseteq \mathcal{O}_{\nu}(\fg_J)$ admit inflations to $\fg$. In particular, all tensor products of highest $\ell$-weight and co-highest $\ell$-weight modules inside $\mathcal{O}_{sh}$ have such inflations~and the rational analogue of \cite[Conjecture 1.4]{pinetinflations} holds.
\end{Corollary}

In fact, a more precise consequence of the above discussion (and of \eqref{eq:GKLOA}--\eqref{eq:GKLOAJ}) is:

\begin{Corollary}\label{cor:identification_of_cat_O}
The equality $Y_{\mu}^{\la}(\bR)=Y_{\nu}^{\la}(\bR,\fg_J)$ identifies $\mathcal{O}_{\nu}^{\la}(\bR,\fg_J)$ with $\mathcal{O}_{\mu}^{\la}(\bR)$ and is compatible with $GT$-characters (in the sense that the $GT$-character of an object~of~$\mathcal{O}_{\nu}^{\la}(\bR,\fg_J)$ coincides with the $GT$-character of the corresponding object of $\mathcal{O}_{\mu}^{\la}(\bR)$ as functions in $\mathcal{E}_{\la-\mu}$). This equality also preserves weights/weight-spaces up to the map $\iota:\mathfrak{h}^*_J\to\mathfrak{h}^*$ given in \eqref{eq:iotaInfl}.
\end{Corollary}
\begin{Rem} Fix $\la'\in P_{I\backslash J}^{\vee}\cap P_+^{\vee}$ and $\bR'\in \C^{\la'}$. Then, %
as subalgebras of $\mathscr{A}=\widetilde{\mathscr{A}}_{\la-\mu}$,
\begin{equation}\label{eq:ShiftTrunc}
Y_{\mu}^{\la}(\bR) = Y_{\mu+\la'}^{\la+\la'}(\bR\cup\bR').
\end{equation}
Moreover, as is straightforward to verify, the equality \eqref{eq:ShiftTrunc}~fits in a diagram of the form\smallskip
{\footnotesize
\begin{equation*}
\begin{tikzcd}[row sep = 1.5em]
Y_{\mu+\la'}\ar[r]\ar[d, two heads, "\Phi_{\mu+\la'}^{\la+\la'}(\bR\cup\bR')",swap] & Y_{\mu}\ar[d, two heads, "\Phi_{\mu}^{\la}(\bR)"]\\
Y_{\mu+\la'}^{\la+\la'}(\bR\cup\bR')\ar[r,equal] & Y_{\mu}^{\la}(\bR)
\end{tikzcd}
\end{equation*}}\noindent
where the top arrow is one of the algebra maps defined in \cite[(4.25)]{hernandez2024shifted} (which generalize the shift morphisms of Section \ref{sec:Algebras}). In particular, $ \mathcal{O}_{\mu}^{\la}(\bR)$ and $\mathcal{O}_{\mu+\la'}^{\la+\la'}(\bR\cup\bR')$ are naturally identified, and, for $V_J$ in $\mathcal{O}_{\nu}^{\la}(\bR,\fg_J)$, we can construct not one, but infinitely-many distinct inflations of $V_J$ to $\fg$ using equalities of truncations (as in Lemma \ref{lemma:EqualityTrunc} and \eqref{eq:ShiftTrunc}).
\end{Rem}
\begin{Rem} By Theorem \ref{thm:BFN} and Corollary \ref{cor:CBfree}, we have algebra isomorphisms
$$Y_{\mu}^{\la}(\bR)\simeq \mathcal{A}_{\hbar=2}(\mathbf{G}, \mathbf{N}; \mathbf{F})\otimes_{S^{\la}}\C \simeq Y_{\nu}^{\la}(\bR,\fg_J),$$ 
where $\mathcal{A}_{\hbar=2}(\mathbf{G}, \mathbf{N}; \mathbf{F})\otimes_{S^{\la}}\C$ is the (specialized) Coulomb branch algebra for the triple
$$
\mathbf{G} = \prod_{i \in I} \operatorname{GL}(m_i), \ \  \mathbf{N} = \bigoplus_{\substack{i, j \in I, \\ i \rightarrow j}} \operatorname{Hom}(\C^{m_i}, \C^{m_j}) \oplus \bigoplus_{i \in I} \operatorname{Hom}(\C^{m_i}, \C^{\la_i})\,\text{ and }\,\mathbf{F} = \prod_{i\in I} \operatorname{GL}(\la_i).
$$
These isomorphisms can be seen as giving a ``Coulomb branch explanation'' for Lemma~\ref{lemma:EqualityTrunc}.
\end{Rem}
\subsection{A notable example}
It is natural to ask if all inflations %
can be constructed as above (i.e.~using equalities of truncations as in Lemma \ref{lemma:EqualityTrunc} and \eqref{eq:ShiftTrunc}). This is unfortunately not~the case. Indeed, choose $i\in I$ and let $J=\{i\}\subseteq I$. Then $\fg_J\simeq \mathfrak{sl}_2$ and one can understand the self-extension $V$ of Example \ref{ex:SelfExt} as a $Y_0(\fg_J)$-module. \medskip\par
The following result about $V$ was mentioned in Section \ref{sec:CorTSC}.
\begin{Proposition} Fix $n\in 2\mathbb{N}$ and $\bR\in \C^{n}/\Sigma_n$. Then $V$ does not descend to $Y_0^{n\varpi_1^{\vee}}(\bR,\fg_J)$.
\end{Proposition}
\begin{proof} Assume the contrary and let $\la=n\varpi_1^{\vee}$. Then $a_{r}:=a_{1,r}$ acts trivially on $V$ if $r>\nicefrac{n}{2}$ by Remark \ref{rem:tsyideal}. Moreover, as the current $\tilde{a}(u)=u^{\nicefrac{n}{2}}(1+\sum_{r\geq 1}a_ru^{-r})$ commutes~with~$h(u)$, the structure of $V$ implies that the matrix $[\tilde{a}(u)]$ representing the action of this current on the $\C$-basis of Example~\ref{ex:SelfExt} has the form
$$ [\tilde{a}(u)] = \adjustbox{scale=0.75}{$ \begin{pmatrix}
p_1(u) & 0 & 0 & 0\\
0 & p_2(u) & 0 & 0 \\
p_3(u) & 0 & p_1(u) & 0\\
0 & p_4(u) & 0 & p_2(u)
\end{pmatrix}$} $$
for polynomials $p_1(u),p_2(u),p_3(u),p_4(u)\in \C[u]$, with $p_1(u)$ and $p_2(u)$ monic of degree $\nicefrac{n}{2}$. On the other hand, by \eqref{eq: def of A gens}, the matrix $[\tilde{a}(u)]$ must satisfy the relation 
\begin{equation}\label{eq:SelfExtEq1}
[h(u)][\tilde{a}(u)][\tilde{a}(u-2)] = p(u) \op{Id}
\end{equation}
where $\op{Id}$ is the $(4{\times}4)$-identity matrix and where, as usual,
$$\textstyle p(u) = \prod_{c\in \bR}(u-c). $$
In particular, taking the $(3{,}1)$-entry of \eqref{eq:SelfExtEq1} and dividing by $\smash{\frac{u+2}{u}p_1(u)p_1(u-2)\neq 0}$ gives
\begin{equation*}
\textstyle \frac{4}{u(u+2)}+f(u)+f(u-2)=0
\end{equation*}
where $\smash{\textstyle f(u) = \frac{p_3(u)}{p_1(u)} \in \C(u)}$. Put differently,
\begin{equation}\label{eq:SelfExtEq2}
\textstyle f(u)+f(u-2)=-\frac{4}{u(u+2)} = \frac{2}{u+2}-\frac{2}{u}.
\end{equation}
We claim %
that the rational 
 function $f(u)$ is regular at $u=0$ and 
$u=-4$. Indeed,~if~$u=-4$ was a pole of $f(u)$, then it would also be a pole of $f(u-2)$ as the RHS of \eqref{eq:SelfExtEq2} is regular~at this point. Equivalently, $f(u)$ would be %
singular at $u=-6$, and repeating the above would prove that $f(u)$ has poles at  $u=-2m$ for all $m\in \mathbb{Z}_{>2}$, contradicting the rationality of this function. Similarly, if $f(u)$ had a pole at $u=0$, then $u=2$ would be a pole of $f(u-2)$,~but not of the RHS of \eqref{eq:SelfExtEq2}, and the above would imply that $f(u)$ is singular at $u=2m$ for~all $m\in \mathbb{N}$, again contradicting rationality. This ends the proof of our claim.%
\medskip\par
Using the above (now proven) claim and taking residues at $u=-2$ in \eqref{eq:SelfExtEq2} gives
$$\textstyle \Res_{-2}(f(u)) = \Res_{-2}(f(u)+f(u-2)) = \Res_{-2}(\frac{2}{u+2}-\frac{2}{u}) = 2,$$
but taking residues of the same equation at $u=0$ instead gives
$$\textstyle \Res_{-2}(f(u)) = \Res_{0}(f(u-2))= \Res_{0}(f(u)+f(u-2)) = \Res_{0}(\frac{2}{u+2}-\frac{2}{u}) = -2.$$
This contradiction shows that \eqref{eq:SelfExtEq2} has no rational solutions and ends the proof.
\end{proof}
Thus, inflations of $V$ to $\fg$ cannot be obtained using equalities of truncations as in Lemma \ref{lemma:EqualityTrunc} and \eqref{eq:ShiftTrunc}. Nevertheless, such inflations can still be constructed in other ways as shown in the lemma below (that can be proven easily by direct computation).
\begin{Lemma} The self-extension $V$ admits an inflation to $\fg$ where the matrices $[e_i(u)]$, $[f_i(u)]$ and $[h_i(u)]$ are those of Example \ref{ex:SelfExt}, and where
$$ [h_j(u)]  =\adjustbox{scale=0.75}{$\begin{pmatrix} 
u-1 & 0 & 0 & 0\\
0 & u+1 & 0 & 0\\
-2 & 0 & u-1 & 0\\
0 & -2 & 0 & u+1
\end{pmatrix}$}$$
for $j\sim i$, with $[h_j(u)]=\op{Id}$ the identity matrix if $i\not\sim j$ and $i\neq j$.
\end{Lemma}
We conclude this appendix with the following conjecture, which is intrinsically related~to \cite[Conjecture 4.4 and Corollary 4.6]{pinetinflations}.
\begin{Conjecture} Fix $J\subseteq I$ and let $V_J$ be an object of $\mathcal{O}_{\nu}^{\la}(\bR,\fg_J)$ for some $\la$, $\nu$~and~$\bR$. Then all inflations of $V_J$ to $\fg$ can be constructed using %
Lemma   \ref{lemma:EqualityTrunc} and \eqref{eq:ShiftTrunc}.
\end{Conjecture}
					
\newpage
\printbibliography

@book {arzhantsev2015coxrings,
    AUTHOR = {Arzhantsev, Ivan and Derenthal, Ulrich and Hausen, J\"urgen and Laface, Antonio},
     TITLE = {Cox rings},
    SERIES = {Cambridge Studies in Advanced Mathematics},
    VOLUME = {144},
 PUBLISHER = {Cambridge University Press, Cambridge},
      YEAR = {2015},
     PAGES = {viii+530},
      ISBN = {978-1-107-02462-5},
   MRCLASS = {14Cxx (14Jxx 14Lxx)},
  MRNUMBER = {3307753},
MRREVIEWER = {Alexandr\ V.\ Pukhlikov},
}

@article {bazhanov2011baxter,
    AUTHOR = {Bazhanov, Vladimir V. and Frassek, Rouven and Łukowski, Tomasz and Meneghelli, Carlo and Staudacher, Matthias},
     TITLE = {Baxter {$\mathbf{Q}$}-operators and representations of {Y}angians},
   JOURNAL = {Nuclear Phys. B},
  FJOURNAL = {Nuclear Physics. B. Theoretical, Phenomenological, and Experimental High Energy Physics. Quantum Field Theory and Statistical Systems},
    VOLUME = {850},
      YEAR = {2011},
    NUMBER = {1},
     PAGES = {148--174},
      ISSN = {0550-3213,1873-1562},
   MRCLASS = {82B23 (16T20 16T25)},
  MRNUMBER = {2803592},
MRREVIEWER = {Yuri\ Kozitsky},
       DOI = {10.1016/j.nuclphysb.2011.04.006},
       URL = {https://doi.org/10.1016/j.nuclphysb.2011.04.006},
}

@article {bedard1999commutation,
    AUTHOR = {B{\'e}dard, Robert},
     TITLE = {On commutation classes of reduced words in {W}eyl groups},
   JOURNAL = {European J. Combin.},
  FJOURNAL = {European Journal of Combinatorics},
    VOLUME = {20},
      YEAR = {1999},
    NUMBER = {6},
     PAGES = {483--505},
      ISSN = {0195-6698,1095-9971},
   MRCLASS = {05E15 (20F55)},
  MRNUMBER = {1703595},
       DOI = {10.1006/eujc.1999.0296},
       URL = {https://doi.org/10.1006/eujc.1999.0296},
}

@article{Besson2026,
  title={Tangent spaces of spherical Schubert varieties and counterexamples to the reducedness conjecture},
  author={Besson, Marc and Hong, Jiuzu and Yu, Huanhuan},
  journal={arXiv preprint arXiv:2603.17273},
  year = {2026}
}

@article {BFN1,
    AUTHOR = {Braverman, Alexander and Finkelberg, Michael and Nakajima, Hiraku},
     TITLE = {Towards a mathematical definition of {C}oulomb branches of 3-dimensional {$\mathcal{N}=4$} gauge theories, {II}},
   JOURNAL = {Adv. Theor. Math. Phys.},
  FJOURNAL = {Advances in Theoretical and Mathematical Physics},
    VOLUME = {22},
      YEAR = {2018},
    NUMBER = {5},
     PAGES = {1071--1147},
      ISSN = {1095-0761,1095-0753},
   MRCLASS = {57R57 (14J33 14N35 16G20 17B67 81T13)},
  MRNUMBER = {3952347},
MRREVIEWER = {Dave\ Auckly},
       DOI = {10.4310/ATMP.2018.v22.n5.a1},
       URL = {https://doi-org.ezproxy.usherbrooke.ca/10.4310/ATMP.2018.v22.n5.a1},
}

@article {braverman2016coulomb,
    AUTHOR = {Braverman, Alexander and Finkelberg, Michael and Nakajima, Hiraku},
     TITLE = {Coulomb branches of {$3d$} {$\mathcal{N}=4$} quiver gauge theories and slices in the affine {G}rassmannian},
      NOTE = {With two appendices by Braverman, Finkelberg, Joel Kamnitzer, Ryosuke Kodera, Nakajima, Ben Webster and Alex Weekes},
   JOURNAL = {Adv. Theor. Math. Phys.},
  FJOURNAL = {Advances in Theoretical and Mathematical Physics},
    VOLUME = {23},
      YEAR = {2019},
    NUMBER = {1},
     PAGES = {75--166},
      ISSN = {1095-0761,1095-0753},
   MRCLASS = {57R57 (14M15 16G20 17B81 81T13)},
  MRNUMBER = {4020310},
       DOI = {10.4310/ATMP.2019.v23.n1.a3},
       URL = {https://doi-org.ezproxy.usherbrooke.ca/10.4310/ATMP.2019.v23.n1.a3},
}

@article{brundan2006shifted,
    AUTHOR = {Brundan, Jonathan and Kleshchev, Alexander},
     TITLE = {Shifted {Y}angians and finite {$W$}-algebras},
   JOURNAL = {Adv. Math.},
  FJOURNAL = {Advances in Mathematics},
    VOLUME = {200},
      YEAR = {2006},
    NUMBER = {1},
     PAGES = {136--195},
      ISSN = {0001-8708,1090-2082},
   MRCLASS = {17B37},
  MRNUMBER = {2199632},
MRREVIEWER = {Chengming\ Bai},
       DOI = {10.1016/j.aim.2004.11.004},
       URL = {https://doi.org/10.1016/j.aim.2004.11.004},
}

@article {borhokraft,
    AUTHOR = {Borho, Walter and Kraft, Hanspeter},
     TITLE = {\"Uber die {G}elfand-{K}irillov-{D}imension},
   JOURNAL = {Math. Ann.},
  FJOURNAL = {Mathematische Annalen},
    VOLUME = {220},
      YEAR = {1976},
    NUMBER = {1},
     PAGES = {1--24},
      ISSN = {0025-5831,1432-1807},
   MRCLASS = {17B35},
  MRNUMBER = {412240},
MRREVIEWER = {Helmut\ Strade},
       DOI = {10.1007/BF01354525},
       URL = {https://doi-org.ezproxy.usherbrooke.ca/10.1007/BF01354525},
}

@article{bourbaki1975lie,
    AUTHOR = {Bourbaki, Nicolas},
     TITLE = {\'{E}l\'{e}ments de math\'{e}matique. {F}asc. {XXXVIII}:
              {G}roupes et alg\`ebres de {L}ie. {C}hapitre {VII}:
              {S}ous-alg\`ebres de {C}artan, \'{e}l\'{e}ments r\'{e}guliers.
              {C}hapitre {VIII}: {A}lg\`ebres de {L}ie semi-simples
              d\'{e}ploy\'{e}es},
    SERIES = {Actualit\'{e}s Scientifiques et Industrielles [Current Scientific and Industrial Topics], No. 1364},
 PUBLISHER = {Hermann, Paris},
      YEAR = {1975},
     PAGES = {271},
   MRCLASS = {17BXX (22E65)},
  MRNUMBER = {453824},
MRREVIEWER = {James\ E.\ Humphreys},
}

@article {berenstein1997total,
    AUTHOR = {Berenstein, Arkady and Zelevinsky, Andrei},
     TITLE = {Total positivity in {S}chubert varieties},
   JOURNAL = {Comment. Math. Helv.},
  FJOURNAL = {Commentarii Mathematici Helvetici},
    VOLUME = {72},
      YEAR = {1997},
    NUMBER = {1},
     PAGES = {128--166},
      ISSN = {0010-2571},
   MRCLASS = {14M15 (22E46)},
  MRNUMBER = {1456321},
MRREVIEWER = {Kailash C. Misra},
       DOI = {10.1007/PL00000363},
       URL = {https://doi.org/10.1007/PL00000363},
}

@incollection {chari2008beyond,
    AUTHOR = {Chari, Vyjayanthi and Hernandez, David},
     TITLE = {Beyond {K}irillov-{R}eshetikhin modules},
 BOOKTITLE = {Quantum affine algebras, extended affine {L}ie algebras, and their applications},
    SERIES = {Contemp. Math.},
    VOLUME = {506},
     PAGES = {49--81},
 PUBLISHER = {Amer. Math. Soc., Providence, RI},
      YEAR = {2010},
      ISBN = {978-0-8218-4507-3},
   MRCLASS = {17B37 (17B10 17B67)},
  MRNUMBER = {2642561},
MRREVIEWER = {Jan\ E.\ Grabowski},
       DOI = {10.1090/conm/506/09935},
       URL = {https://doi.org/10.1090/conm/506/09935},
}

@article{chari2002braid,
    AUTHOR = {Chari, Vyjayanthi},
     TITLE = {Braid group actions and tensor products},
   JOURNAL = {Int. Math. Res. Not.},
  FJOURNAL = {International Mathematics Research Notices},
      YEAR = {2002},
    NUMBER = {7},
     PAGES = {357--382},
      ISSN = {1073-7928,1687-0247},
   MRCLASS = {17B37 (20F36)},
  MRNUMBER = {1883181},
MRREVIEWER = {Rinat\ Kedem},
       DOI = {10.1155/S107379280210612X},
       URL = {https://doi.org/10.1155/S107379280210612X},
}

@article{chari2005characters,
    AUTHOR = {Chari, Vyjayanthi and Moura, Adriano A.},
     TITLE = {Characters and blocks for finite-dimensional representations of quantum affine algebras},
   JOURNAL = {Int. Math. Res. Not.},
  FJOURNAL = {International Mathematics Research Notices},
      YEAR = {2005},
    NUMBER = {5},
     PAGES = {257--298},
      ISSN = {1073-7928,1687-0247},
   MRCLASS = {17B67 (17B10 17B37)},
  MRNUMBER = {2130797},
MRREVIEWER = {Jacob\ Greenstein},
       DOI = {10.1155/IMRN.2005.257},
       URL = {https://doi.org/10.1155/IMRN.2005.257},
}

@article{chari1990yangians,
    AUTHOR = {Chari, Vyjayanthi and Pressley, Andrew},
     TITLE = {Yangians and {$R$}-matrices},
   JOURNAL = {Enseign. Math. (2)},
  FJOURNAL = {L'Enseignement Math\'{e}matique. Revue Internationale. 2e S\'{e}rie},
    VOLUME = {36},
      YEAR = {1990},
    NUMBER = {3-4},
     PAGES = {267--302},
      ISSN = {0013-8584},
   MRCLASS = {17B37 (81R50)},
  MRNUMBER = {1096420},
MRREVIEWER = {Vladimir\ A.\ Stukopin},
}

@article{chari1991quantum,
    AUTHOR = {Chari, Vyjayanthi and Pressley, Andrew},
     TITLE = {Quantum affine algebras},
   JOURNAL = {Comm. Math. Phys.},
  FJOURNAL = {Communications in Mathematical Physics},
    VOLUME = {142},
      YEAR = {1991},
    NUMBER = {2},
     PAGES = {261--283},
      ISSN = {0010-3616,1432-0916},
   MRCLASS = {17B37},
  MRNUMBER = {1137064},
MRREVIEWER = {Vladimir\ A.\ Stukopin},
       URL = {http://projecteuclid.org/euclid.cmp/1104248585},
}

@article {chuang2008derived,
    AUTHOR = {Chuang, Joseph and Rouquier, Rapha\"{e}l},
     TITLE = {Derived equivalences for symmetric groups and {$\mathfrak{sl}_2$}-categorification},
   JOURNAL = {Ann. of Math. (2)},
  FJOURNAL = {Annals of Mathematics. Second Series},
    VOLUME = {167},
      YEAR = {2008},
    NUMBER = {1},
     PAGES = {245--298},
      ISSN = {0003-486X,1939-8980},
   MRCLASS = {20C08 (17B10 18E15 20C33)},
  MRNUMBER = {2373155},
MRREVIEWER = {Bogdan\ Ion},
       DOI = {10.4007/annals.2008.167.245},
       URL = {https://doi.org/10.4007/annals.2008.167.245},
}

@article {cautis2019cluster,
    AUTHOR = {Cautis, Sabin and Williams, Harold},
     TITLE = {Cluster theory of the coherent {S}atake category},
   JOURNAL = {J. Amer. Math. Soc.},
  FJOURNAL = {Journal of the American Mathematical Society},
    VOLUME = {32},
      YEAR = {2019},
    NUMBER = {3},
     PAGES = {709--778},
      ISSN = {0894-0347,1088-6834},
   MRCLASS = {22E67 (13F60 14F05)},
  MRNUMBER = {3981987},
MRREVIEWER = {Truong\ Le\ Hoang},
       DOI = {10.1090/jams/918},
       URL = {https://doi.org/10.1090/jams/918},
}

@article{dumanski2025k,
  title={K-theoretic Hikita conjecture for quiver gauge theories},
  author={Dumanski, Ilya and Krylov, Vasily},
  journal={arXiv preprint arXiv:2509.06226},
  year={2025}
}

@article {dranowski2024heaps,
    AUTHOR = {Dranowski, Anne and Elek, Bal\'azs and Kamnitzer, Joel and Morton-Ferguson, Calder},
     TITLE = {Heaps, crystals, and preprojective algebra modules},
   JOURNAL = {Selecta Math. (N.S.)},
  FJOURNAL = {Selecta Mathematica. New Series},
    VOLUME = {30},
      YEAR = {2024},
    NUMBER = {5},
     PAGES = {Paper No. 94, 47},
      ISSN = {1022-1824,1420-9020},
   MRCLASS = {17B10 (14M15 16G20 17B37 22E57)},
  MRNUMBER = {4812519},
       DOI = {10.1007/s00029-024-00978-8},
       URL = {https://doi.org/10.1007/s00029-024-00978-8},
}

@article {drinfeld1985hopf,
    AUTHOR = {Drinfeld, Vladimir G.},
     TITLE = {Hopf algebras and the quantum {Y}ang-{B}axter equation},
   JOURNAL = {Dokl. Akad. Nauk SSSR},
  FJOURNAL = {Doklady Akademii Nauk SSSR},
    VOLUME = {283},
      YEAR = {1985},
    NUMBER = {5},
     PAGES = {1060--1064},
      ISSN = {0002-3264},
   MRCLASS = {58F07 (17B20 82A15)},
  MRNUMBER = {802128},
MRREVIEWER = {Alexander\ A.\ Pankov},
}

@article{drinfeld1987new,
    AUTHOR = {Drinfeld, Vladimir G.},
     TITLE = {A new realization of {Y}angians and of quantum affine algebras},
   JOURNAL = {Dokl. Akad. Nauk SSSR},
  FJOURNAL = {Doklady Akademii Nauk SSSR},
    VOLUME = {296},
      YEAR = {1987},
    NUMBER = {1},
     PAGES = {13--17},
      ISSN = {0002-3264},
   MRCLASS = {17B65 (16A24 17B45 58F07 81E99)},
  MRNUMBER = {914215},
MRREVIEWER = {J.\ S.\ Joel},
}

@article {elek2021bott,
    AUTHOR = {Elek, Bal\'{a}zs and Lu, Jiang-Hua},
     TITLE = {Bott-{S}amelson varieties and {P}oisson {O}re extensions},
   JOURNAL = {Int. Math. Res. Not. IMRN},
  FJOURNAL = {International Mathematics Research Notices. IMRN},
      YEAR = {2021},
    NUMBER = {14},
     PAGES = {10745--10797},
      ISSN = {1073-7928,1687-0247},
   MRCLASS = {22E46 (14M15 17B63 20F55)},
  MRNUMBER = {4285734},
MRREVIEWER = {William\ M.\ McGovern},
       DOI = {10.1093/imrn/rnz127},
       URL = {https://doi.org/10.1093/imrn/rnz127},
}

@article{etingof2003elliptic,
    AUTHOR = {Etingof, Pavel I. and Moura, Adriano A.},
     TITLE = {Elliptic central characters and blocks of finite dimensional representations of quantum affine algebras},
   JOURNAL = {Represent. Theory},
  FJOURNAL = {Representation Theory. An Electronic Journal of the American Mathematical Society},
    VOLUME = {7},
      YEAR = {2003},
     PAGES = {346--373},
      ISSN = {1088-4165},
   MRCLASS = {17B37 (20G42)},
  MRNUMBER = {2017062},
MRREVIEWER = {David\ Hernandez},
       DOI = {10.1090/S1088-4165-03-00201-2},
       URL = {https://doi.org/10.1090/S1088-4165-03-00201-2},
}

@article{frenkel2015baxter,
    AUTHOR = {Frenkel, Edward and Hernandez, David},
     TITLE = {Baxter's relations and spectra of quantum integrable models},
   JOURNAL = {Duke Math. J.},
  FJOURNAL = {Duke Mathematical Journal},
    VOLUME = {164},
      YEAR = {2015},
    NUMBER = {12},
     PAGES = {2407--2460},
      ISSN = {0012-7094,1547-7398},
   MRCLASS = {17B67 (17B10 17B37 82B23)},
  MRNUMBER = {3397389},
MRREVIEWER = {Kailash\ C.\ Misra},
       DOI = {10.1215/00127094-3146282},
       URL = {https://doi.org/10.1215/00127094-3146282},
}

@article{frenkel2024extended,
    AUTHOR = {Frenkel, Edward and Hernandez, David},
     TITLE = {Extended {B}axter relations and {QQ}-systems for quantum affine algebras},
   JOURNAL = {Comm. Math. Phys.},
  FJOURNAL = {Communications in Mathematical Physics},
    VOLUME = {405},
      YEAR = {2024},
    NUMBER = {8},
     PAGES = {Paper No. 190, 42},
      ISSN = {0010-3616,1432-0916},
   MRCLASS = {17B67 (17B37 18M05 81R50)},
  MRNUMBER = {4779559},
MRREVIEWER = {Zhaobing\ Fan},
       DOI = {10.1007/s00220-024-05051-1},
       URL = {https://doi.org/10.1007/s00220-024-05051-1},
}

@article {frenkel2022weyl,
    AUTHOR = {Frenkel, Edward and Hernandez, David},
     TITLE = {Weyl group symmetry of {$q$}-characters},
   JOURNAL = {Selecta Math. (N.S.)},
  FJOURNAL = {Selecta Mathematica. New Series},
    VOLUME = {31},
      YEAR = {2025},
    NUMBER = {4},
     PAGES = {Paper No. 72, 36},
      ISSN = {1022-1824,1420-9020},
   MRCLASS = {17B37 (17B10 17B67 81R50 82B23)},
  MRNUMBER = {4939703},
       DOI = {10.1007/s00029-025-01072-3},
       URL = {https://doi.org/10.1007/s00029-025-01072-3},
}

@article {fujita2024monoidal,
    AUTHOR = {Fujita, Ryo and Hernandez, David},
     TITLE = {Monoidal {J}antzen filtrations},
   JOURNAL = {Adv. Math.},
  FJOURNAL = {Advances in Mathematics},
    VOLUME = {495},
      YEAR = {2026},
     PAGES = {Paper No. 110963},
      ISSN = {0001-8708,1090-2082},
   MRCLASS = {17B37 (16T25 17B10 17B67 20G42 81R50)},
  MRNUMBER = {5061762},
       DOI = {10.1016/j.aim.2026.110963},
       URL = {https://doi.org/10.1016/j.aim.2026.110963},
}

@article{finkelberg2018comultiplication,
    AUTHOR = {Finkelberg, Michael and Kamnitzer, Joel and Pham, Khoa and Rybnikov, Leonid and Weekes, Alex},
     TITLE = {Comultiplication for shifted {Y}angians and quantum open {T}oda lattice},
   JOURNAL = {Adv. Math.},
  FJOURNAL = {Advances in Mathematics},
    VOLUME = {327},
      YEAR = {2018},
     PAGES = {349--389},
      ISSN = {0001-8708,1090-2082},
   MRCLASS = {17B37 (81T25)},
  MRNUMBER = {3761996},
MRREVIEWER = {Huafeng\ Zhang},
       DOI = {10.1016/j.aim.2017.06.018},
       URL = {https://doi.org/10.1016/j.aim.2017.06.018},
}

@article {frenkel1997canonical,
    AUTHOR = {Frenkel, Igor B. and Khovanov, Mikhail G.},
     TITLE = {Canonical bases in tensor products and graphical calculus for {$U_q(\mathfrak{sl}_2)$}},
   JOURNAL = {Duke Math. J.},
  FJOURNAL = {Duke Mathematical Journal},
    VOLUME = {87},
      YEAR = {1997},
    NUMBER = {3},
     PAGES = {409--480},
      ISSN = {0012-7094,1547-7398},
   MRCLASS = {17B37 (17B10)},
  MRNUMBER = {1446615},
MRREVIEWER = {Kailash\ C.\ Misra},
       DOI = {10.1215/S0012-7094-97-08715-9},
       URL = {https://doi.org/10.1215/S0012-7094-97-08715-9},
}

@article{francone2025cluster,
  title={Cluster structures on schemes of bands},
  author={Francone, Luca and Leclerc, Bernard},
  journal={arXiv preprint arXiv:2504.14012},
  year={2025}
}

@article {fujita2020flag,
    AUTHOR = {Fujita, Naoki and Lee, Eunjeong and Suh, Dong Youp},
     TITLE = {Algebraic and geometric properties of flag {B}ott-{S}amelson varieties and applications to representations},
   JOURNAL = {Pacific J. Math.},
  FJOURNAL = {Pacific Journal of Mathematics},
    VOLUME = {309},
      YEAR = {2020},
    NUMBER = {1},
     PAGES = {145--194},
      ISSN = {0030-8730,1945-5844},
   MRCLASS = {05E10 (14M15 57S25)},
  MRNUMBER = {4202007},
       DOI = {10.2140/pjm.2020.309.145},
       URL = {https://doi.org/10.2140/pjm.2020.309.145},
}

@article{frenkel2001combinatorics,
    AUTHOR = {Frenkel, Edward and Mukhin, Evgeny},
     TITLE = {Combinatorics of {$q$}-characters of finite-dimensional representations of quantum affine algebras},
   JOURNAL = {Comm. Math. Phys.},
  FJOURNAL = {Communications in Mathematical Physics},
    VOLUME = {216},
      YEAR = {2001},
    NUMBER = {1},
     PAGES = {23--57},
      ISSN = {0010-3616,1432-0916},
   MRCLASS = {17B37 (17B10 81R50)},
  MRNUMBER = {1810773},
MRREVIEWER = {Alexei\ P.\ Isaev},
       DOI = {10.1007/s002200000323},
       URL = {https://doi.org/10.1007/s002200000323},
}

@incollection{frenkel1999qcharacters,
    AUTHOR = {Frenkel, Edward and Reshetikhin, Nicolai},
     TITLE = {The {$q$}-characters of representations of quantum affine algebras and deformations of {$W$}-algebras},
 BOOKTITLE = {Recent developments in quantum affine algebras and related topics ({R}aleigh, {NC}, 1998)},
    SERIES = {Contemp. Math.},
    VOLUME = {248},
     PAGES = {163--205},
 PUBLISHER = {Amer. Math. Soc., Providence, RI},
      YEAR = {1999},
   MRCLASS = {17B37 (05E15 17B68)},
  MRNUMBER = {1745260},
       DOI = {10.1090/conm/248/03823},
       URL = {https://doi.org/10.1090/conm/248/03823},
}

@article {frenkel2024qopers,
    AUTHOR = {Frenkel, Edward and Koroteev, Peter and Sage, Daniel S. and Zeitlin, Anton M.},
     TITLE = {{$q$}-opers, {$QQ$}-systems, and {B}ethe ansatz},
   JOURNAL = {J. Eur. Math. Soc. (JEMS)},
  FJOURNAL = {Journal of the European Mathematical Society (JEMS)},
    VOLUME = {26},
      YEAR = {2024},
    NUMBER = {1},
     PAGES = {355--405},
      ISSN = {1435-9855,1435-9863},
   MRCLASS = {82B23 (14D24 17B80 22E10 81R50)},
  MRNUMBER = {4705654},
MRREVIEWER = {Jun\ Pei},
       DOI = {10.4171/jems/1268},
       URL = {https://doi.org/10.4171/jems/1268},
}

@incollection {finkelberg2019multiplicative,
    AUTHOR = {Finkelberg, Michael and Tsymbaliuk, Alexander},
     TITLE = {Multiplicative slices, relativistic {T}oda and shifted quantum affine algebras},
 BOOKTITLE = {Representations and nilpotent orbits of {L}ie algebraic systems},
    SERIES = {Progr. Math.},
    VOLUME = {330},
     PAGES = {133--304},
 PUBLISHER = {Birkh\"{a}user/Springer, Cham},
      YEAR = {2019},
   MRCLASS = {17B37 (81R10 81T13)},
  MRNUMBER = {3971731},
MRREVIEWER = {Kyungyong\ Lee},
       DOI = {10.1007/978-3-030-23531-4\{_}6}

@article{hernandez2021quantum,
    AUTHOR = {Fujita, Ryo and Hernandez, David and Oh, Se-jin and Oya, Hironori},
     TITLE = {Isomorphisms among quantum {G}rothendieck rings and propagation of positivity},
   JOURNAL = {J. Reine Angew. Math.},
  FJOURNAL = {Journal f\"{u}r die Reine und Angewandte Mathematik. [Crelle's Journal]},
    VOLUME = {785},
      YEAR = {2022},
     PAGES = {117--185},
      ISSN = {0075-4102,1435-5345},
   MRCLASS = {17B67 (17B37)},
  MRNUMBER = {4402493},
MRREVIEWER = {Euiyong\ Park},
       DOI = {10.1515/crelle-2021-0088},
       URL = {https://doi.org/10.1515/crelle-2021-0088},
}

@article{friesen2025braid,
    AUTHOR = {Friesen, Noah and Weekes, Alex and Wendlandt, Curtis},
     TITLE = {Braid group actions, {B}axter polynomials, and affine quantum groups},
   JOURNAL = {Trans. Amer. Math. Soc.},
  FJOURNAL = {Transactions of the American Mathematical Society},
    VOLUME = {378},
      YEAR = {2025},
    NUMBER = {2},
     PAGES = {1329--1372},
      ISSN = {0002-9947,1088-6850},
   MRCLASS = {17B37 (17B10)},
  MRNUMBER = {4850442},
MRREVIEWER = {Run-Qiang\ Jian},
       DOI = {10.1090/tran/9279},
       URL = {https://doi.org/10.1090/tran/9279},
}

@article{fomin2002cluster,
    AUTHOR = {Fomin, Sergey and Zelevinsky, Andrei},
     TITLE = {Cluster algebras. {I}. {F}oundations},
   JOURNAL = {J. Amer. Math. Soc.},
  FJOURNAL = {Journal of the American Mathematical Society},
    VOLUME = {15},
      YEAR = {2002},
    NUMBER = {2},
     PAGES = {497--529},
      ISSN = {0894-0347,1088-6834},
   MRCLASS = {16S99 (14M99 17B99)},
  MRNUMBER = {1887642},
MRREVIEWER = {Eric\ N.\ Sommers},
       DOI = {10.1090/S0894-0347-01-00385-X},
       URL = {https://doi.org/10.1090/S0894-0347-01-00385-X},
}

@article {fomin199double,
    AUTHOR = {Fomin, Sergey and Zelevinsky, Andrei},
     TITLE = {Double {B}ruhat cells and total positivity},
   JOURNAL = {J. Amer. Math. Soc.},
  FJOURNAL = {Journal of the American Mathematical Society},
    VOLUME = {12},
      YEAR = {1999},
    NUMBER = {2},
     PAGES = {335--380},
      ISSN = {0894-0347,1088-6834},
   MRCLASS = {20G20 (15A23)},
  MRNUMBER = {1652878},
       DOI = {10.1090/S0894-0347-99-00295-7},
       URL = {https://doi.org/10.1090/S0894-0347-99-00295-7},
}

@article {gerasimov2005class,
    AUTHOR = {Gerasimov, Anton A. and Kharchev, Sergei M. and Lebedev, Dmitry R. and Oblezin, Sergey V.},
     TITLE = {On a class of representations of the {Y}angian and moduli space of monopoles},
   JOURNAL = {Comm. Math. Phys.},
  FJOURNAL = {Communications in Mathematical Physics},
    VOLUME = {260},
      YEAR = {2005},
    NUMBER = {3},
     PAGES = {511--525},
      ISSN = {0010-3616},
   MRCLASS = {53C07 (17B37 53D17 53D30)},
  MRNUMBER = {2182434},
MRREVIEWER = {Olivier G. Schiffmann},
       DOI = {10.1007/s00220-005-1417-3},
       URL = {https://doi.org/10.1007/s00220-005-1417-3},
}

@article{geiss2024representations,
    AUTHOR = {Geiss, Christof and Hernandez, David and Leclerc, Bernard},
     TITLE = {Representations of shifted quantum affine algebras and cluster algebras {I}: {T}he simply laced case},
   JOURNAL = {Proc. Lond. Math. Soc. (3)},
  FJOURNAL = {Proceedings of the London Mathematical Society. Third Series},
    VOLUME = {129},
      YEAR = {2024},
    NUMBER = {3},
     PAGES = {Paper No. e12630, 75},
      ISSN = {0024-6115,1460-244X},
   MRCLASS = {17B67 (13F60 17B10 17B37 35J25 82B23)},
  MRNUMBER = {4793282},
MRREVIEWER = {Bing\ Duan},
       DOI = {10.1112/plms.12630},
       URL = {https://doi.org/10.1112/plms.12630},
}

@article {gibson2021demazure,
    AUTHOR = {Gibson, Joel},
     TITLE = {A {D}emazure character formula for the product monomial crystal},
   JOURNAL = {Algebr. Comb.},
  FJOURNAL = {Algebraic Combinatorics},
    VOLUME = {4},
      YEAR = {2021},
    NUMBER = {2},
     PAGES = {301--327},
   MRCLASS = {17B37 (05E10)},
  MRNUMBER = {4244375},
MRREVIEWER = {Jacinta Torres},
       DOI = {10.5802/alco.156},
       URL = {https://doi.org/10.5802/alco.156},
}

@article{grigorev2025representations,
  title={On representations of quantum affine {$\mathfrak{sl}_2$}},
  author={Grigorev, Andrei and Mukhin, Evgeny},
  journal={arXiv preprint arXiv:2505.11605},
  year={2025}
}

@article{guay2018coproduct,
    AUTHOR = {Guay, Nicolas and Nakajima, Hiraku and Wendlandt, Curtis},
     TITLE = {Coproduct for {Y}angians of affine {K}ac-{M}oody algebras},
   JOURNAL = {Adv. Math.},
  FJOURNAL = {Advances in Mathematics},
    VOLUME = {338},
      YEAR = {2018},
     PAGES = {865--911},
      ISSN = {0001-8708,1090-2082},
   MRCLASS = {17B37 (17B67)},
  MRNUMBER = {3861718},
MRREVIEWER = {Volodymyr\ Mazorchuk},
       DOI = {10.1016/j.aim.2018.09.013},
       URL = {https://doi.org/10.1016/j.aim.2018.09.013},
}

@article{gautam2016yangians,
    AUTHOR = {Gautam, Sachin and Toledano Laredo, Valerio},
     TITLE = {Yangians, quantum loop algebras, and abelian difference equations},
   JOURNAL = {J. Amer. Math. Soc.},
  FJOURNAL = {Journal of the American Mathematical Society},
    VOLUME = {29},
      YEAR = {2016},
    NUMBER = {3},
     PAGES = {775--824},
      ISSN = {0894-0347,1088-6834},
   MRCLASS = {17B67 (17B37)},
  MRNUMBER = {3486172},
MRREVIEWER = {Rutwig\ Campoamor-Stursberg},
       DOI = {10.1090/jams/851},
       URL = {https://doi.org/10.1090/jams/851},
}

@article{gautam2023poles,
    AUTHOR = {Gautam, Sachin and Wendlandt, Curtis},
     TITLE = {Poles of finite-dimensional representations of {Y}angians},
   JOURNAL = {Selecta Math. (N.S.)},
  FJOURNAL = {Selecta Mathematica. New Series},
    VOLUME = {29},
      YEAR = {2023},
    NUMBER = {1},
     PAGES = {Paper No. 13, 68},
      ISSN = {1022-1824,1420-9020},
   MRCLASS = {17B37 (81R05)},
  MRNUMBER = {4520258},
MRREVIEWER = {Cristian\ Vay},
       DOI = {10.1007/s00029-022-00813-y},
       URL = {https://doi.org/10.1007/s00029-022-00813-y},
}

@book {hartshorne1977ag,
    AUTHOR = {Hartshorne, Robin},
     TITLE = {Algebraic geometry},
    SERIES = {Graduate Texts in Mathematics},
    VOLUME = {52},
 PUBLISHER = {Springer-Verlag, New York-Heidelberg},
      YEAR = {1977},
     PAGES = {xvi+496},
      ISBN = {0-387-90244-9},
   MRCLASS = {14-01},
  MRNUMBER = {463157},
MRREVIEWER = {Robert\ Speiser},
}

@article{hernandez2010simple,
    AUTHOR = {Hernandez, David},
     TITLE = {Simple tensor products},
   JOURNAL = {Invent. Math.},
  FJOURNAL = {Inventiones Mathematicae},
    VOLUME = {181},
      YEAR = {2010},
    NUMBER = {3},
     PAGES = {649--675},
      ISSN = {0020-9910,1432-1297},
   MRCLASS = {17B37 (17B10)},
  MRNUMBER = {2660455},
MRREVIEWER = {Jonathan\ Brundan},
       DOI = {10.1007/s00222-010-0256-9},
       URL = {https://doi.org/10.1007/s00222-010-0256-9},
}

@article{hernandez2011quantum,
    AUTHOR = {Hernandez, David and Leclerc, Bernard},
     TITLE = {Quantum {G}rothendieck rings and derived {H}all algebras},
   JOURNAL = {J. Reine Angew. Math.},
  FJOURNAL = {Journal f\"{u}r die Reine und Angewandte Mathematik. [Crelle's
              Journal]},
    VOLUME = {701},
      YEAR = {2015},
     PAGES = {77--126},
      ISSN = {0075-4102,1435-5345},
   MRCLASS = {17B67 (17B10 17B37)},
  MRNUMBER = {3331727},
MRREVIEWER = {Youjun\ Tan},
       DOI = {10.1515/crelle-2013-0020},
       URL = {https://doi.org/10.1515/crelle-2013-0020},
}

@article{hernandez2019cyclicity,
    AUTHOR = {Hernandez, David},
     TITLE = {Cyclicity and {$R$}-matrices},
   JOURNAL = {Selecta Math. (N.S.)},
  FJOURNAL = {Selecta Mathematica. New Series},
    VOLUME = {25},
      YEAR = {2019},
    NUMBER = {2},
     PAGES = {Paper No. 19, 24},
      ISSN = {1022-1824,1420-9020},
   MRCLASS = {17B37 (17B10 81R50)},
  MRNUMBER = {3916087},
MRREVIEWER = {J\"{o}rg\ Feldvoss},
       DOI = {10.1007/s00029-019-0465-z},
       URL = {https://doi.org/10.1007/s00029-019-0465-z},
}

@article{hernandez2023representations,
    AUTHOR = {Hernandez, David},
     TITLE = {Representations of shifted quantum affine algebras},
   JOURNAL = {Int. Math. Res. Not. IMRN},
  FJOURNAL = {International Mathematics Research Notices. IMRN},
      YEAR = {2023},
    NUMBER = {13},
     PAGES = {11035--11126},
      ISSN = {1073-7928,1687-0247},
   MRCLASS = {17B37 (13F60)},
  MRNUMBER = {4609778},
MRREVIEWER = {Xin\ Tang},
       DOI = {10.1093/imrn/rnac149},
       URL = {https://doi.org/10.1093/imrn/rnac149},
}

@article{hernandez2025symmetries,
  title={Symmetries of Grothendieck rings in representation theory},
  author={Hernandez, David},
  journal={arXiv preprint arXiv:2501.03024},
  year={2025}
}

@article{hernandez2026crossroad,
  title={Representations and characters of quantum affine algebras at the crossroads between cluster categorification and quantum integrable models},
  author={Hernandez, David},
  journal={International Congress of Mathematicians 2026},
  year={2026},
  pages = {532--551}
}

@article{hernandez2012asymptotic,
    AUTHOR = {Hernandez, David and Jimbo, Michio},
     TITLE = {Asymptotic representations and {D}rinfeld rational fractions},
   JOURNAL = {Compos. Math.},
  FJOURNAL = {Compositio Mathematica},
    VOLUME = {148},
      YEAR = {2012},
    NUMBER = {5},
     PAGES = {1593--1623},
      ISSN = {0010-437X,1570-5846},
   MRCLASS = {17B37 (17B10)},
  MRNUMBER = {2982441},
MRREVIEWER = {Olga\ Bershtein},
       DOI = {10.1112/S0010437X12000267},
       URL = {https://doi.org/10.1112/S0010437X12000267},
}

@article{hernandez2010cluster,
    AUTHOR = {Hernandez, David and Leclerc, Bernard},
     TITLE = {Cluster algebras and quantum affine algebras},
   JOURNAL = {Duke Math. J.},
  FJOURNAL = {Duke Mathematical Journal},
    VOLUME = {154},
      YEAR = {2010},
    NUMBER = {2},
     PAGES = {265--341},
      ISSN = {0012-7094,1547-7398},
   MRCLASS = {17B37 (13F60)},
  MRNUMBER = {2682185},
MRREVIEWER = {Jan\ E.\ Grabowski},
       DOI = {10.1215/00127094-2010-040},
       URL = {https://doi.org/10.1215/00127094-2010-040},
}

@article{hernandez2016cluster,
    AUTHOR = {Hernandez, David and Leclerc, Bernard},
     TITLE = {Cluster algebras and category {$\mathcal{O}$} for representations of {B}orel subalgebras of quantum affine algebras},
   JOURNAL = {Algebra Number Theory},
  FJOURNAL = {Algebra \& Number Theory},
    VOLUME = {10},
      YEAR = {2016},
    NUMBER = {9},
     PAGES = {2015--2052},
      ISSN = {1937-0652,1944-7833},
   MRCLASS = {17B67 (13F60 17B10 17B37 82B23)},
  MRNUMBER = {3576119},
MRREVIEWER = {Jan\ E.\ Grabowski},
       DOI = {10.2140/ant.2016.10.2015},
       URL = {https://doi.org/10.2140/ant.2016.10.2015},
}

@article{hernandez2026borel,
  title={Borel and shifted category {$\mathcal{O}$}},
  author={Hernandez, David and Neguţ, Andrei},
  journal={arXiv preprint arXiv:2603.00928},
  year={2026}
}

@article {hernandez2024shifted,
    AUTHOR = {Hernandez, David and Zhang, Huafeng},
     TITLE = {Shifted {Y}angians and polynomial {$R$}-matrices},
   JOURNAL = {Publ. Res. Inst. Math. Sci.},
  FJOURNAL = {Publications of the Research Institute for Mathematical Sciences},
    VOLUME = {60},
      YEAR = {2024},
    NUMBER = {1},
     PAGES = {1--69},
      ISSN = {0034-5318},
   MRCLASS = {20G42 (16T25 81R50)},
  MRNUMBER = {4803333},
       DOI = {10.4171/prims/60-1-1},
       URL = {https://doi.org/10.4171/prims/60-1-1},
}

@article{hernandez2025jordan,
  title={Jordan-Hölder property for shifted quantum affine algebras},
  author={Hernandez, David and Zhang, Huafeng},
  journal={arXiv preprint arXiv:2501.16859},
  year={2025}
}

@book {jantzen2003representations,
    AUTHOR = {Jantzen, Jens Carsten},
     TITLE = {Representations of algebraic groups},
    SERIES = {Mathematical Surveys and Monographs},
    VOLUME = {107},
   EDITION = {Second},
 PUBLISHER = {American Mathematical Society, Providence, RI},
      YEAR = {2003},
     PAGES = {xiv+576},
      ISBN = {0-8218-3527-0},
   MRCLASS = {20G05 (17B10)},
  MRNUMBER = {2015057},
}

@article{jang2025unipotent,
    AUTHOR = {Jang, Il-Seung and Kwon, Jae-Hoon and Park, Euiyong},
     TITLE = {Unipotent quantum coordinate ring and cominuscule prefundamental representations},
   JOURNAL = {J. Algebra},
  FJOURNAL = {Journal of Algebra},
    VOLUME = {673},
      YEAR = {2025},
     PAGES = {260--303},
      ISSN = {0021-8693,1090-266X},
   MRCLASS = {17B37 (05E10 22E46)},
  MRNUMBER = {4878804},
MRREVIEWER = {Iwan\ Praton},
       DOI = {10.1016/j.jalgebra.2025.01.031},
       URL = {https://doi.org/10.1016/j.jalgebra.2025.01.031},
}

@book {kac1990infinite,
    AUTHOR = {Kac, Victor G.},
     TITLE = {Infinite-dimensional {L}ie algebras},
   EDITION = {Third},
 PUBLISHER = {Cambridge University Press, Cambridge},
      YEAR = {1990},
     PAGES = {xxii+400},
   MRCLASS = {17B65 (17B67 17B68 58F07)},
  MRNUMBER = {1104219},
       DOI = {10.1017/CBO9780511626234},
       URL = {https://doi.org/10.1017/CBO9780511626234},
}

@article{otherpaper,
  title={Almost-dominant chamber modules and reverse plane partitions},
  author={Kalmykov, Artem and Kamnitzer, Joel and Leroux-Lapierre, Alexis and Pinet, Théo and Weekes, Alex},
  note = {To appear},
  year={2026}
}

@article{kamnitzer2005mirkovic,
    AUTHOR = {Kamnitzer, Joel},
     TITLE = {Mirkovi\'{c}-{V}ilonen cycles and polytopes},
   JOURNAL = {Ann. of Math. (2)},
  FJOURNAL = {Annals of Mathematics. Second Series},
    VOLUME = {171},
      YEAR = {2010},
    NUMBER = {1},
     PAGES = {245--294},
      ISSN = {0003-486X,1939-8980},
   MRCLASS = {20G05 (14M15)},
  MRNUMBER = {2630039},
MRREVIEWER = {Peter\ Fiebig},
       DOI = {10.4007/annals.2010.171.245},
       URL = {https://doi.org/10.4007/annals.2010.171.245},
}

@article{kamnitzer2014yangians,
    AUTHOR = {Kamnitzer, Joel and Webster, Ben and Weekes, Alex and Yacobi, Oded},
     TITLE = {Yangians and quantizations of slices in the affine {G}rassmannian},
   JOURNAL = {Algebra Number Theory},
  FJOURNAL = {Algebra \& Number Theory},
    VOLUME = {8},
      YEAR = {2014},
    NUMBER = {4},
     PAGES = {857--893},
      ISSN = {1937-0652},
   MRCLASS = {17B37 (14D24 14M15 20G15 53D55)},
  MRNUMBER = {3248988},
MRREVIEWER = {Christian Ohn},
       DOI = {10.2140/ant.2014.8.857},
       URL = {https://doi.org/10.2140/ant.2014.8.857},
}

@article{kamnitzer2018reducedness,
    AUTHOR = {Kamnitzer, Joel and Muthiah, Dinakar and Weekes, Alex and Yacobi, Oded},
     TITLE = {Reducedness of affine {G}rassmannian slices in type {A}},
   JOURNAL = {Proc. Amer. Math. Soc.},
  FJOURNAL = {Proceedings of the American Mathematical Society},
    VOLUME = {146},
      YEAR = {2018},
    NUMBER = {2},
     PAGES = {861--874},
      ISSN = {0002-9939,1088-6826},
   MRCLASS = {22E67 (14M15)},
  MRNUMBER = {3731717},
MRREVIEWER = {Andrea\ Pinamonti},
       DOI = {10.1090/proc/13850},
       URL = {https://doi.org/10.1090/proc/13850},
}

@article{kamnitzer2019highest,
    AUTHOR = {Kamnitzer, Joel and Tingley, Peter and Webster, Ben and Weekes, Alex and Yacobi, Oded},
     TITLE = {Highest weights for truncated shifted {Y}angians and product monomial crystals},
   JOURNAL = {J. Comb. Algebra},
  FJOURNAL = {Journal of Combinatorial Algebra},
    VOLUME = {3},
      YEAR = {2019},
    NUMBER = {3},
     PAGES = {237--303},
      ISSN = {2415-6302},
   MRCLASS = {16G20 (05E10 16T99 17B10)},
  MRNUMBER = {4011667},
MRREVIEWER = {Aleksandr Panov},
       DOI = {10.4171/JCA/32},
       URL = {https://doi.org/10.4171/JCA/32},
}

@article{kamnitzer2019category,
    AUTHOR = {Kamnitzer, Joel and Tingley, Peter and Webster, Ben and  Weekes, Alex and Yacobi, Oded},
     TITLE = {On category {$\mathcal O$} for affine {G}rassmannian slices and categorified tensor products},
   JOURNAL = {Proc. Lond. Math. Soc. (3)},
  FJOURNAL = {Proceedings of the London Mathematical Society. Third Series},
    VOLUME = {119},
      YEAR = {2019},
    NUMBER = {5},
     PAGES = {1179--1233},
      ISSN = {0024-6115},
   MRCLASS = {14M15 (17B37 20G42)},
  MRNUMBER = {3968721},
MRREVIEWER = {Huafeng Zhang},
       DOI = {10.1112/plms.12254},
       URL = {https://doi.org/10.1112/plms.12254},
}

@article{kamnitzer2022lie,
    AUTHOR = {Kamnitzer, Joel and Webster, Ben and Weekes, Alex and Yacobi, Oded},
     TITLE = {Lie algebra actions on module categories for truncated shifted yangians},
   JOURNAL = {Forum Math. Sigma},
  FJOURNAL = {Forum of Mathematics. Sigma},
    VOLUME = {12},
      YEAR = {2024},
     PAGES = {Paper No. e18, 69},
   MRCLASS = {20G05 (14M15 17B67)},
  MRNUMBER = {4699878},
MRREVIEWER = {Linliang Song},
       DOI = {10.1017/fms.2024.3},
       URL = {https://doi.org/10.1017/fms.2024.3},
}

@article{kang2018monoidal,
    AUTHOR = {Kang, Seok-Jin and Kashiwara, Masaki and Kim, Myungho and Oh, Se-jin},
     TITLE = {Monoidal categorification of cluster algebras},
   JOURNAL = {J. Amer. Math. Soc.},
  FJOURNAL = {Journal of the American Mathematical Society},
    VOLUME = {31},
      YEAR = {2018},
    NUMBER = {2},
     PAGES = {349--426},
      ISSN = {0894-0347,1088-6834},
   MRCLASS = {13F60 (16Gxx 17B37 18D10 81R50)},
  MRNUMBER = {3758148},
MRREVIEWER = {Fan\ Qin},
       DOI = {10.1090/jams/895},
       URL = {https://doi.org/10.1090/jams/895},
}

@incollection {kashiwara2003realizations,
    AUTHOR = {Kashiwara, Masaki},
     TITLE = {Realizations of crystals},
 BOOKTITLE = {Combinatorial and geometric representation theory ({S}eoul, 2001)},
    SERIES = {Contemp. Math.},
    VOLUME = {325},
     PAGES = {133--139},
 PUBLISHER = {Amer. Math. Soc., Providence, RI},
      YEAR = {2003},
   MRCLASS = {17B37 (17B10)},
  MRNUMBER = {1988989},
MRREVIEWER = {Jae-Hoon Kwon},
       DOI = {10.1090/conm/325/05668},
       URL = {https://doi.org/10.1090/conm/325/05668},
}

@article{kashiwara2024monoidal,
    AUTHOR = {Kashiwara, Masaki and Kim, Myungho and Oh, Se-jin and Park, Euiyong},
     TITLE = {Monoidal categorification and quantum affine algebras {II}},
   JOURNAL = {Invent. Math.},
  FJOURNAL = {Inventiones Mathematicae},
    VOLUME = {236},
      YEAR = {2024},
    NUMBER = {2},
     PAGES = {837--924},
      ISSN = {0020-9910,1432-1297},
   MRCLASS = {17B37 (13F60 18M15)},
  MRNUMBER = {4728243},
MRREVIEWER = {Farrokh\ Razavinia},
       DOI = {10.1007/s00222-024-01249-1},
       URL = {https://doi.org/10.1007/s00222-024-01249-1},
}

@article{khovanov2012extended,
    AUTHOR = {Khovanov, Mikhail and Lauda, Aaron D. and Mackaay, Marco and Sto\v{s}i\'{c}, Marko},
     TITLE = {Extended graphical calculus for categorified quantum {${\mathfrak{sl}}(2)$}},
   JOURNAL = {Mem. Amer. Math. Soc.},
  FJOURNAL = {Memoirs of the American Mathematical Society},
    VOLUME = {219},
      YEAR = {2012},
    NUMBER = {1029},
     PAGES = {vi+87},
      ISSN = {0065-9266,1947-6221},
      ISBN = {978-0-8218-8977-0},
   MRCLASS = {81R50 (05Axx 18D05)},
  MRNUMBER = {2963085},
MRREVIEWER = {Fan\ Xu},
       DOI = {10.1090/S0065-9266-2012-00665-4},
       URL = {https://doi.org/10.1090/S0065-9266-2012-00665-4},
}

@thesis {khovanov1997graphical,
    AUTHOR = {Khovanov, Mikhail},
     TITLE = {Graphical calculus, canonical bases and {K}azhdan-{L}usztig theory},
      NOTE = {Thesis (Ph.D.)--Yale University},
 PUBLISHER = {ProQuest LLC, Ann Arbor, MI},
      YEAR = {1997},
     PAGES = {103},
      ISBN = {978-0591-43629-7},
   MRCLASS = {Thesis},
  MRNUMBER = {2695927},
       URL ={http://gateway.proquest.com/openurl?url_ver=Z39.88-2004&rft_val_fmt=info:ofi/fmt:kev:mtx:dissertation&res_dat=xri:pqdiss&rft_dat=xri:pqdiss:9733946},
}

@article {kang2012categorification,
    AUTHOR = {Kang, Seok-Jin and Kashiwara, Masaki},
     TITLE = {Categorification of highest weight modules via {K}hovanov-{L}auda-{R}ouquier algebras},
   JOURNAL = {Invent. Math.},
  FJOURNAL = {Inventiones Mathematicae},
    VOLUME = {190},
      YEAR = {2012},
    NUMBER = {3},
     PAGES = {699--742},
      ISSN = {0020-9910,1432-1297},
   MRCLASS = {17B67 (17B37)},
  MRNUMBER = {2995184},
MRREVIEWER = {Volodymyr\ Mazorchuk},
       DOI = {10.1007/s00222-012-0388-1},
       URL = {https://doi.org/10.1007/s00222-012-0388-1},
}

@book{krause2000growth,
    AUTHOR = {Krause, G\"{u}nter R. and Lenagan, Thomas H.},
     TITLE = {Growth of algebras and {G}elfand-{K}irillov dimension},
    SERIES = {Graduate Studies in Mathematics},
    VOLUME = {22},
   EDITION = {Revised},
 PUBLISHER = {American Mathematical Society, Providence, RI},
      YEAR = {2000},
     PAGES = {x+212},
      ISBN = {0-8218-0859-1},
   MRCLASS = {16P90},
  MRNUMBER = {1721834},
MRREVIEWER = {Martha\ K.\ Smith},
       DOI = {10.1090/gsm/022},
       URL = {https://doi.org/10.1090/gsm/022},
}

@article{khovanov2009diagrammatic,
    AUTHOR = {Khovanov, Mikhail and Lauda, Aaron D.},
     TITLE = {A diagrammatic approach to categorification of quantum groups {I}},
   JOURNAL = {Represent. Theory},
  FJOURNAL = {Representation Theory. An Electronic Journal of the American Mathematical Society},
    VOLUME = {13},
      YEAR = {2009},
     PAGES = {309--347},
   MRCLASS = {17B37},
  MRNUMBER = {2525917},
MRREVIEWER = {Fan Xu},
       DOI = {10.1090/S1088-4165-09-00346-X},
       URL = {https://doi.org/10.1090/S1088-4165-09-00346-X},
}

@article{khovanov2011diagrammatic,
    AUTHOR = {Khovanov, Mikhail and Lauda, Aaron D.},
     TITLE = {A diagrammatic approach to categorification of quantum groups {II}},
   JOURNAL = {Trans. Amer. Math. Soc.},
  FJOURNAL = {Transactions of the American Mathematical Society},
    VOLUME = {363},
      YEAR = {2011},
    NUMBER = {5},
     PAGES = {2685--2700},
      ISSN = {0002-9947},
   MRCLASS = {17B37 (16T20)},
  MRNUMBER = {2763732},
MRREVIEWER = {Volodymyr Mazorchuk},
       DOI = {10.1090/S0002-9947-2010-05210-9},
       URL = {https://doi.org/10.1090/S0002-9947-2010-05210-9},
}

@book {kleshchev2005linear,
    AUTHOR = {Kleshchev, Alexander},
     TITLE = {Linear and projective representations of symmetric groups},
    SERIES = {Cambridge Tracts in Mathematics},
    VOLUME = {163},
 PUBLISHER = {Cambridge University Press, Cambridge},
      YEAR = {2005},
     PAGES = {xiv+277},
      ISBN = {0-521-83703-0},
   MRCLASS = {20C30 (20C08)},
  MRNUMBER = {2165457},
MRREVIEWER = {Christine\ Bessenrodt},
       DOI = {10.1017/CBO9780511542800},
       URL = {https://doi.org/10.1017/CBO9780511542800},
}

@article {knight1995spectra,
    AUTHOR = {Knight, Harold},
     TITLE = {Spectra of tensor products of finite-dimensional representations of {Y}angians},
   JOURNAL = {J. Algebra},
  FJOURNAL = {Journal of Algebra},
    VOLUME = {174},
      YEAR = {1995},
    NUMBER = {1},
     PAGES = {187--196},
      ISSN = {0021-8693,1090-266X},
   MRCLASS = {17B37 (81R50)},
  MRNUMBER = {1332866},
MRREVIEWER = {Preeti\ Parashar},
       DOI = {10.1006/jabr.1995.1123},
       URL = {https://doi.org/10.1006/jabr.1995.1123},
}

@article{krylovperunov2021,
    AUTHOR = {Krylov, Vasily and Perunov, Ivan},
     TITLE = {Almost dominant generalized slices and convolution diagrams over them},
   JOURNAL = {Adv. Math.},
  FJOURNAL = {Advances in Mathematics},
    VOLUME = {392},
      YEAR = {2021},
     PAGES = {Paper No. 108034, 45},
      ISSN = {0001-8708},
   MRCLASS = {14M15 (20G05)},
  MRNUMBER = {4316674},
MRREVIEWER = {Ryan David Kinser},
       DOI = {10.1016/j.aim.2021.108034},
       URL = {https://doi.org/10.1016/j.aim.2021.108034},
}

@article{kamnitzer2022hamiltonian,
    AUTHOR = {Kamnitzer, Joel and Pham, Khoa and Weekes, Alex},
     TITLE = {Hamiltonian reduction for affine {G}rassmannian slices and truncated shifted {Y}angians},
   JOURNAL = {Adv. Math.},
  FJOURNAL = {Advances in Mathematics},
    VOLUME = {399},
      YEAR = {2022},
     PAGES = {Paper No. 108281, 52},
      ISSN = {0001-8708,1090-2082},
   MRCLASS = {14M15 (17B10)},
  MRNUMBER = {4385132},
MRREVIEWER = {Felipe\ Zald\'{\i}var},
       DOI = {10.1016/j.aim.2022.108281},
       URL = {https://doi.org/10.1016/j.aim.2022.108281},
}

@article {kleshchev2010homogeneous,
    AUTHOR = {Kleshchev, Alexander and Ram, Arun},
     TITLE = {Homogeneous representations of {K}hovanov-{L}auda algebras},
   JOURNAL = {J. Eur. Math. Soc. (JEMS)},
  FJOURNAL = {Journal of the European Mathematical Society (JEMS)},
    VOLUME = {12},
      YEAR = {2010},
    NUMBER = {5},
     PAGES = {1293--1306},
      ISSN = {1435-9855},
   MRCLASS = {20C08 (17B67)},
  MRNUMBER = {2677617},
MRREVIEWER = {Andrew Mathas},
       DOI = {10.4171/JEMS/230},
       URL = {https://doi.org/10.4171/JEMS/230},
}

@article {kleshchev2011representations,
    AUTHOR = {Kleshchev, Alexander and Ram, Arun},
     TITLE = {Representations of {K}hovanov-{L}auda-{R}ouquier algebras and combinatorics of {L}yndon words},
   JOURNAL = {Math. Ann.},
  FJOURNAL = {Mathematische Annalen},
    VOLUME = {349},
      YEAR = {2011},
    NUMBER = {4},
     PAGES = {943--975},
      ISSN = {0025-5831,1432-1807},
   MRCLASS = {16S99 (16G10 20C08)},
  MRNUMBER = {2777040},
MRREVIEWER = {Selvaraj\ Chelliah},
       DOI = {10.1007/s00208-010-0543-1},
       URL = {https://doi.org/10.1007/s00208-010-0543-1},
}

@article{Alexis,
  title={Category {$\mathcal{O}$} and asymptotic characters},
  author={Leroux-Lapierre, Alexis},
  journal={arXiv preprint arXiv:2507.16215},
  year={2025}
}

@article {lakshmibai2002standard,
    AUTHOR = {Lakshmibai, Venkatramani and Littelmann, Peter and Magyar, Peter},
     TITLE = {Standard monomial theory for {B}ott-{S}amelson varieties},
   JOURNAL = {Compositio Math.},
  FJOURNAL = {Compositio Mathematica},
    VOLUME = {130},
      YEAR = {2002},
    NUMBER = {3},
     PAGES = {293--318},
      ISSN = {0010-437X,1570-5846},
   MRCLASS = {14M15 (14L30)},
  MRNUMBER = {1887117},
MRREVIEWER = {E.\ A.\ Tevel\"ev},
       DOI = {10.1023/A:1014396129323},
       URL = {https://doi.org/10.1023/A:1014396129323},
}

@article {lauritzen2004line,
    AUTHOR = {Lauritzen, Niels and Thomsen, Jesper Funch},
     TITLE = {Line bundles on {B}ott-{S}amelson varieties},
   JOURNAL = {J. Algebraic Geom.},
  FJOURNAL = {Journal of Algebraic Geometry},
    VOLUME = {13},
      YEAR = {2004},
    NUMBER = {3},
     PAGES = {461--473},
      ISSN = {1056-3911,1534-7486},
   MRCLASS = {14M17 (14C22 14L35)},
  MRNUMBER = {2047677},
MRREVIEWER = {Vikram\ B.\ Mehta},
       DOI = {10.1090/S1056-3911-03-00358-8},
       URL = {https://doi.org/10.1090/S1056-3911-03-00358-8},
}

@article {lusztig1990canonicalI,
    AUTHOR = {Lusztig, George},
     TITLE = {Canonical bases arising from quantized enveloping algebras},
   JOURNAL = {J. Amer. Math. Soc.},
  FJOURNAL = {Journal of the American Mathematical Society},
    VOLUME = {3},
      YEAR = {1990},
    NUMBER = {2},
     PAGES = {447--498},
      ISSN = {0894-0347,1088-6834},
   MRCLASS = {17B35 (16A64)},
  MRNUMBER = {1035415},
MRREVIEWER = {Ya.\ S.\ So\u{\i}bel\cprime man},
       DOI = {10.2307/1990961},
       URL = {https://doi.org/10.2307/1990961},
}

@book {lusztig1993introduction,
    AUTHOR = {Lusztig, George},
     TITLE = {Introduction to quantum groups},
    SERIES = {Progress in Mathematics},
    VOLUME = {110},
 PUBLISHER = {Birkh\"{a}user Boston, Inc., Boston, MA},
      YEAR = {1993},
     PAGES = {xii+341},
      ISBN = {0-8176-3712-5},
   MRCLASS = {17B37 (16W30 17-02 17B35 81R50)},
  MRNUMBER = {1227098},
MRREVIEWER = {Jie\ Du},
}

@article{lauda2011crystals,
    AUTHOR = {Lauda, Aaron D. and Vazirani, Monica},
     TITLE = {Crystals from categorified quantum groups},
   JOURNAL = {Adv. Math.},
  FJOURNAL = {Advances in Mathematics},
    VOLUME = {228},
      YEAR = {2011},
    NUMBER = {2},
     PAGES = {803--861},
      ISSN = {0001-8708,1090-2082},
   MRCLASS = {17B37},
  MRNUMBER = {2822211},
MRREVIEWER = {Peter\ W.\ Tingley},
       DOI = {10.1016/j.aim.2011.06.009},
       URL = {https://doi.org/10.1016/j.aim.2011.06.009},
}

@article {magyar1998borel,
    AUTHOR = {Magyar, Peter},
     TITLE = {Borel-{W}eil theorem for configuration varieties and {S}chur modules},
   JOURNAL = {Adv. Math.},
  FJOURNAL = {Advances in Mathematics},
    VOLUME = {134},
      YEAR = {1998},
    NUMBER = {2},
     PAGES = {328--366},
      ISSN = {0001-8708,1090-2082},
   MRCLASS = {14M17 (05E10 14F17 14N15)},
  MRNUMBER = {1617793},
MRREVIEWER = {Laurent\ Manivel},
       DOI = {10.1006/aima.1997.1700},
       URL = {https://doi.org/10.1006/aima.1997.1700},
}

@article {magyar1998schubert,
    AUTHOR = {Magyar, Peter},
     TITLE = {Schubert polynomials and {B}ott-{S}amelson varieties},
   JOURNAL = {Comment. Math. Helv.},
  FJOURNAL = {Commentarii Mathematici Helvetici},
    VOLUME = {73},
      YEAR = {1998},
    NUMBER = {4},
     PAGES = {603--636},
      ISSN = {0010-2571,1420-8946},
   MRCLASS = {14M15 (05E15 16G20 20G05)},
  MRNUMBER = {1639896},
MRREVIEWER = {Witold\ Kra\'skiewicz},
       DOI = {10.1007/s000140050071},
       URL = {https://doi.org/10.1007/s000140050071},
}

@article{milot2025,
	title={Compatibility between truncation and coproducts for quantum affine algebra and {Y}angian of {$\mathfrak{sl}_2(\C)$}},
	author = {Milot, {J\'er\^ome}},
	journal = {arXiv preprint arXiv:2506.10544},
	year = 2025
}

@article {mathas2024cellularity,
    AUTHOR = {Mathas, Andrew and Tubbenhauer, Daniel},
     TITLE = {Cellularity and subdivision of {KLR} and weighted {KLRW} algebras},
   JOURNAL = {Math. Ann.},
  FJOURNAL = {Mathematische Annalen},
    VOLUME = {389},
      YEAR = {2024},
    NUMBER = {3},
     PAGES = {3043--3122},
      ISSN = {0025-5831,1432-1807},
   MRCLASS = {20C08 (18M30 18N25 20G43)},
  MRNUMBER = {4753081},
MRREVIEWER = {Lei\ Shi},
       DOI = {10.1007/s00208-023-02660-4},
       URL = {https://doi.org/10.1007/s00208-023-02660-4},
}

@article{mukhin2014affinization,
    AUTHOR = {Mukhin, Evgeny and Young, Charles A. S.},
     TITLE = {Affinization of category {$\mathcal{O}$} for quantum groups},
   JOURNAL = {Trans. Amer. Math. Soc.},
  FJOURNAL = {Transactions of the American Mathematical Society},
    VOLUME = {366},
      YEAR = {2014},
    NUMBER = {9},
     PAGES = {4815--4847},
      ISSN = {0002-9947,1088-6850},
   MRCLASS = {17B37},
  MRNUMBER = {3217701},
MRREVIEWER = {Christian\ Ohn},
       DOI = {10.1090/S0002-9947-2014-06039-X},
       URL = {https://doi.org/10.1090/S0002-9947-2014-06039-X},
}

@article{nakajima2003tanalogues,
    AUTHOR = {Nakajima, Hiraku},
     TITLE = {{$t$}-analogs of {$q$}-characters of quantum affine algebras of type {$A_n,D_n$}},
 BOOKTITLE = {Combinatorial and geometric representation theory ({S}eoul, 2001)},
    SERIES = {Contemp. Math.},
    VOLUME = {325},
     PAGES = {141--160},
 PUBLISHER = {Amer. Math. Soc., Providence, RI},
      YEAR = {2003},
   MRCLASS = {17B37 (81R50 82B23)},
  MRNUMBER = {1988990},
MRREVIEWER = {\c{S}erban Raianu},
       DOI = {10.1090/conm/325/05669},
       URL = {https://doi.org/10.1090/conm/325/05669},
}

@article{neguct2025category,
  title={Category {$\mathcal{O}$} for quantum loop algebras},
  author={Negu{\c{t}}, Andrei},
  journal={arXiv preprint arXiv:2501.00724},
  year={2025}
}

@article{paganelli2025quantum,
    AUTHOR = {Paganelli, Francesca},
     TITLE = {Quantum cluster algebras and representations of shifted quantum affine algebras},
   JOURNAL = {Math. Z.},
  FJOURNAL = {Mathematische Zeitschrift},
    VOLUME = {313},
      YEAR = {2026},
    NUMBER = {2},
     PAGES = {Paper No. 34},
      ISSN = {0025-5874,1432-1823},
   MRCLASS = {17B37 (13F60 16G20 17B67 19M05)},
  MRNUMBER = {5076886},
       DOI = {10.1007/s00209-026-04038-z},
       URL = {https://doi.org/10.1007/s00209-026-04038-z},
}

@article{pressley1991fundamental,
    AUTHOR = {Chari, Vyjayanthi and Pressley, Andrew},
     TITLE = {Fundamental representations of {Y}angians and singularities of {$R$}-matrices},
   JOURNAL = {J. Reine Angew. Math.},
  FJOURNAL = {Journal f\"{u}r die Reine und Angewandte Mathematik. [Crelle's Journal]},
    VOLUME = {417},
      YEAR = {1991},
     PAGES = {87--128},
      ISSN = {0075-4102,1435-5345},
   MRCLASS = {17B37},
  MRNUMBER = {1103907},
MRREVIEWER = {Vladimir\ A.\ Stukopin},
       DOI = {10.1515/crll.1991.417.87},
       URL = {https://doi.org/10.1515/crll.1991.417.87},
}

@article{pinet2024functor,
    AUTHOR = {Pinet, Th\'{e}o},
     TITLE = {A functor for constructing {$R$}-matrices in the category {$\mathcal{O}$} of {B}orel quantum loop algebras},
   JOURNAL = {J. Lond. Math. Soc. (2)},
  FJOURNAL = {Journal of the London Mathematical Society. Second Series},
    VOLUME = {109},
      YEAR = {2024},
    NUMBER = {1},
     PAGES = {Paper No. e12815, 46},
      ISSN = {0024-6107,1469-7750},
   MRCLASS = {17B37 (13F60 16D90 18M15 81R10)},
  MRNUMBER = {4680197},
MRREVIEWER = {Euiyong\ Park},
       DOI = {10.1112/jlms.12815},
       URL = {https://doi.org/10.1112/jlms.12815},
}

@article{pinetinflations,
    AUTHOR = {Pinet, Th\'{e}o},
     TITLE = {Inflations for representations of shifted quantum affine algebras},
   JOURNAL = {Adv. Math.},
  FJOURNAL = {Advances in Mathematics},
    VOLUME = {462},
      YEAR = {2025},
     PAGES = {Paper No. 110093, 51},
      ISSN = {0001-8708,1090-2082},
   MRCLASS = {17B37 (13F60 16T20 17B10 17B67)},
  MRNUMBER = {4843776},
MRREVIEWER = {Run-Qiang\ Jian},
       DOI = {10.1016/j.aim.2024.110093},
       URL = {https://doi.org/10.1016/j.aim.2024.110093},
}

@book {procesi2007liegroups,
    AUTHOR = {Procesi, Claudio},
     TITLE = {Lie groups},
    SERIES = {Universitext},
      NOTE = {An approach through invariants and representations},
 PUBLISHER = {Springer, New York},
      YEAR = {2007},
     PAGES = {xxiv+596},
   MRCLASS = {22E10 (05E10 05E15 14L24 14M15 17B10 20G05 22-02)},
  MRNUMBER = {2265844},
MRREVIEWER = {James\ E.\ Humphreys},
}

@incollection {ringel1980tame,
    AUTHOR = {Ringel, Claus Michael},
     TITLE = {On algorithms for solving vector space problems. {II}. {T}ame algebras},
 BOOKTITLE = {Representation theory, {I} ({P}roc. {W}orkshop, {C}arleton {U}niv., {O}ttawa, {O}nt., 1979)},
    SERIES = {Lecture Notes in Math.},
    VOLUME = {831},
     PAGES = {137--287},
 PUBLISHER = {Springer, Berlin},
      YEAR = {1980},
      ISBN = {3-540-10263-9},
   MRCLASS = {16A64 (15A21 16-02)},
  MRNUMBER = {607143},
MRREVIEWER = {Sheila\ Brenner},
}

@article{rouquier20082,
  title   = {2-{K}ac-{M}oody algebras},
  author  = {Rouquier, Rapha{\"e}l},
  journal = {arXiv preprint arXiv:0812.5023},
  year    = {2008}
}

@article {soergel1986equivalence,
    AUTHOR = {Soergel, Wolfgang},
     TITLE = {\'{E}quivalences de certaines cat\'{e}gories de $\mathfrak{g}$-modules},
   JOURNAL = {C. R. Acad. Sci. Paris S\'{e}r. I Math.},
  FJOURNAL = {Comptes Rendus des S\'{e}ances de l'Acad\'{e}mie des Sciences. S\'{e}rie I. Math\'{e}matique},
    VOLUME = {303},
      YEAR = {1986},
    NUMBER = {15},
     PAGES = {725--728},
      ISSN = {0249-6291},
   MRCLASS = {17B10 (22E47)},
  MRNUMBER = {872544},
MRREVIEWER = {James\ E.\ Humphreys},
}

@misc{stacks-project,
  author       = {The {Stacks project authors}},
  title        = {The Stacks project},
  howpublished = {\url{https://stacks.math.columbia.edu}},
  year         = {2025},
}

@article {silverthorne2024gelfand,
    AUTHOR = {Silverthorne, Turner and Webster, Ben},
     TITLE = {Gelfand-{T}setlin modules: canonicity and calculations},
   JOURNAL = {Algebr. Represent. Theory},
  FJOURNAL = {Algebras and Representation Theory},
    VOLUME = {27},
      YEAR = {2024},
    NUMBER = {2},
     PAGES = {1405--1455},
      ISSN = {1386-923X,1572-9079},
   MRCLASS = {17B10},
  MRNUMBER = {4741528},
MRREVIEWER = {Stefano\ Capparelli},
       DOI = {10.1007/s10468-024-10264-y},
       URL = {https://doi.org/10.1007/s10468-024-10264-y},
}

@book{vazirani1999irreducible,
    AUTHOR = {Vazirani, Monica Joy},
     TITLE = {Irreducible modules over the affine {H}ecke algebra: {A} strong multiplicity one result},
      NOTE = {Thesis (Ph.D.)--University of California, Berkeley},
 PUBLISHER = {ProQuest LLC, Ann Arbor, MI},
      YEAR = {1999},
     PAGES = {54},
      ISBN = {978-0599-31645-4},
   MRCLASS = {Thesis},
  MRNUMBER = {2699147},
       URL = {http://gateway.proquest.com/openurl?url_ver=Z39.88-2004&rft_val_fmt=info:ofi/fmt:kev:mtx:dissertation&res_dat=xri:pqdiss&rft_dat=xri:pqdiss:9931428},
}

@incollection {viennot1986heaps,
    AUTHOR = {Viennot, G\'erard Xavier},
     TITLE = {Heaps of pieces. {I}. {B}asic definitions and combinatorial lemmas},
 BOOKTITLE = {Combinatoire \'enum\'erative ({M}ontreal, {Q}ue., 1985)},
    SERIES = {Lecture Notes in Math.},
    VOLUME = {1234},
     PAGES = {321--350},
 PUBLISHER = {Springer, Berlin},
      YEAR = {1986},
      ISBN = {3-540-17207-6},
   MRCLASS = {05A99 (20M10)},
  MRNUMBER = {927773},
MRREVIEWER = {Dominique\ Perrin},
       DOI = {10.1007/BFb0072524},
       URL = {https://doi.org/10.1007/BFb0072524},
}

@article {varagnolo2011canonical,
    AUTHOR = {Varagnolo, Michela and Vasserot, Eric},
     TITLE = {Canonical bases and {KLR}-algebras},
   JOURNAL = {J. Reine Angew. Math.},
  FJOURNAL = {Journal f\"{u}r die Reine und Angewandte Mathematik. [Crelle's Journal]},
    VOLUME = {659},
      YEAR = {2011},
     PAGES = {67--100},
      ISSN = {0075-4102,1435-5345},
   MRCLASS = {17B37 (16T20)},
  MRNUMBER = {2837011},
MRREVIEWER = {Nicolas\ Jacon},
       DOI = {10.1515/CRELLE.2011.068},
       URL = {https://doi.org/10.1515/CRELLE.2011.068},
}

@article{varagnolo2025representations,
  title={Representations of shifted affine quantum groups and Coulomb branches},
  author={Varagnolo, Michela and Vasserot, Eric},
  journal={arXiv preprint arXiv:2503.06262},
  year={2025},
}

@article {webster2015canonical,
    AUTHOR = {Webster, Ben},
     TITLE = {Canonical bases and higher representation theory},
   JOURNAL = {Compos. Math.},
  FJOURNAL = {Compositio Mathematica},
    VOLUME = {151},
      YEAR = {2015},
    NUMBER = {1},
     PAGES = {121--166},
      ISSN = {0010-437X,1570-5846},
   MRCLASS = {17B37 (17B10 18D05)},
  MRNUMBER = {3305310},
MRREVIEWER = {Kevin\ D.\ Coulembier},
       DOI = {10.1112/S0010437X1400760X},
       URL = {https://doi.org/10.1112/S0010437X1400760X},
}

@article {webster2017knot,
    AUTHOR = {Webster, Ben},
     TITLE = {Knot invariants and higher representation theory},
   JOURNAL = {Mem. Amer. Math. Soc.},
  FJOURNAL = {Memoirs of the American Mathematical Society},
    VOLUME = {250},
      YEAR = {2017},
    NUMBER = {1191},
     PAGES = {v+141},
      ISSN = {0065-9266,1947-6221},
   MRCLASS = {57M27 (17B10 18D05 57M25)},
  MRNUMBER = {3709726},
MRREVIEWER = {Stefan\ K.\ Friedl},
       DOI = {10.1090/memo/1191},
       URL = {https://doi.org/10.1090/memo/1191},
}

@book {weekesthesis,
    AUTHOR = {Weekes, Alex},
     TITLE = {Highest {W}eights for {T}runcated {S}hifted {Y}angians},
      NOTE = {Thesis (Ph.D.) -- University of Toronto (Canada)},
 PUBLISHER = {ProQuest LLC, Ann Arbor, MI},
      YEAR = {2016},
     PAGES = {102},
      ISBN = {978-1369-85412-1},
   MRCLASS = {99-05},
  MRNUMBER = {3697595},
       URL = {http://gateway.proquest.com.ezproxy.usherbrooke.ca/openurl?url_ver=Z39.88-2004&rft_val_fmt=info:ofi/fmt:kev:mtx:dissertation&res_dat=xri:pqm&rft_dat=xri:pqdiss:10190626},
}

@article{weekes2019,
  title   = {Generators for {C}oulomb branches of quiver gauge theories},
  author  = {Weekes, Alex},
  journal = {arXiv preprint arXiv:1903.07734},
  year    = {2019},
}

@article{webster2020quantum,
    AUTHOR = {Webster, Ben and Weekes, Alex and Yacobi, Oded},
     TITLE = {A quantum {M}irkovi\'{c}-{V}ybornov isomorphism},
   JOURNAL = {Represent. Theory},
  FJOURNAL = {Representation Theory. An Electronic Journal of the American Mathematical Society},
    VOLUME = {24},
      YEAR = {2020},
     PAGES = {38--84},
      ISSN = {1088-4165},
   MRCLASS = {17B37 (16S80 20C99)},
  MRNUMBER = {4052554},
MRREVIEWER = {Aleksandr\ Panov},
       DOI = {10.1090/ert/536},
       URL = {https://doi.org/10.1090/ert/536},
}

@article{zhang2020yangians,
    AUTHOR = {Zhang, Huafeng},
     TITLE = {Yangians and {B}axter's relations},
   JOURNAL = {Lett. Math. Phys.},
  FJOURNAL = {Letters in Mathematical Physics},
    VOLUME = {110},
      YEAR = {2020},
    NUMBER = {8},
     PAGES = {2113--2141},
      ISSN = {0377-9017,1573-0530},
   MRCLASS = {17B37 (17B10 17B80)},
  MRNUMBER = {4126875},
       DOI = {10.1007/s11005-020-01285-x},
       URL = {https://doi.org/10.1007/s11005-020-01285-x},
}

@article{zhang2024theta,
    AUTHOR = {Zhang, Huafeng},
     TITLE = {Theta series for quantum loop algebras and {Y}angians},
   JOURNAL = {Comm. Math. Phys.},
  FJOURNAL = {Communications in Mathematical Physics},
    VOLUME = {405},
      YEAR = {2024},
    NUMBER = {10},
     PAGES = {Paper No. 230, 68},
      ISSN = {0010-3616,1432-0916},
   MRCLASS = {20G42 (16T20 81R50)},
  MRNUMBER = {4797741},
MRREVIEWER = {Cristian\ Vay},
       DOI = {10.1007/s00220-024-05110-7},
       URL = {https://doi.org/10.1007/s00220-024-05110-7},
}

\end{document}